\documentclass{article}

\usepackage{amsmath}
\usepackage{amsfonts}
\usepackage{amssymb}

\newtheorem{theorem}{Theorem}[section]

\newtheorem{definition}{Definition}[section]
\newtheorem{lemma}{Lemma}[section]
\newtheorem{corollary}{Corollary}[section]

\newtheorem{remark}{Remark}[section]

\newtheorem{example}{Example}[section]

\begin{document}

\author{V.S. Mikhalevich \and A.M. Gupal \and V.I. Norkin }

\title{METHODS OF NONCONVEX OPTIMIZATION}


\date{Nauka, Moscow, 1987.}
\maketitle
Mikhalevich V. S., Gupal A. M., Norkin V. I. Methods of Nonconvex Optimization.  Moscow, Nauka (Phys.-Mat. Litt.), 1987. 280 p. (Translated from Russian by V.I. Norkin, 2024)

In applied problems it is often necessary to deal with nondifferentiable functions and non-convex regions. Development of optimization methods for such problems is a great challenge. The book deals with finite-dimensional nonconvex nonsmooth optimization problems and numerical methods for their solution. 

For specialists in applied mathematics, computer science, economics, and cybernetics.

Bibliography: 188 items.

Reviewer: Doctor of Physical and Mathematical Sciences V.V. Fedorov

Book review: J Glob Optim 1, 205–206 (1991). https://doi.org/10.1007/BF00119992

\bigskip\bigskip
\noindent
English version, June, 2024,\\ translated by Vladimir I. Norkin, \\e-mail: 
\\v.norkin@kpi.ua
\\vladimir.norkin@gmail.com

\newpage
\textbf{CONTENT}

\smallskip
{\bf Preface $\ldots$} \pageref{FOREWORD}

\smallskip
{\bf Chapter 1. Elements of the theory of nonconvex analysis $\ldots$} \pageref{Ch.1}

$\S$ 1. Generalized differentiable functions $\ldots$ \pageref{Sec.1}

$\S$ 2. Smoothing of functions and basic properties of locally Lipschitz functions $\ldots$ \pageref{Sec.2}

$\S$ 3. Necessary optimality conditions $\ldots$ \pageref{Sec.3}

{\bf \smallskip
Chapter 2. Minimization of Lipschitz functions without computing gradients $\ldots$} \pageref{Ch.2}

$\S$ 4. Finite-difference method of minimization of Lipschitz functions $\ldots$ \pageref{Sec.4}

$\S$ 5. Construction of finite differences for the maximum function $\ldots$ \pageref{Sec.5}

$\S$ 6. Random finite-difference directions $\ldots$ \pageref{Sec.6}

$\S$ 7. Efficiency of finite-difference methods $\ldots$ \pageref{Sec.7}

{\bf\smallskip
Chapter 3. Generalized gradient descent methods  $\ldots$} \pageref{Ch.3}

$\S$ 8. Conditions for convergence of iterative algorithms for nonlinear programming  $\ldots$ \pageref{Sec.8}

\begin{footnotesize}
\hspace{0.3cm}1.Definitions...\pageref{Sec.8.1}. 2.Necessary and sufficient conditions for the convergence of algorithms...\pageref{Sec.8.2}. 3.Sufficient conditions for convergence...\pageref{Sec.8.3}.
\end{footnotesize}

$\S$ 9. The generalized gradient decent method $\ldots$ \pageref{Sec.9}

$\S$ 10. The generalized gradient method in a conditional optimization problem  $\ldots$ \pageref{Sec.10}

\begin{footnotesize}
\hspace{0.3cm}1.One inequality-constraint problem...\pageref{Sec.10.1}.  2.Multiple inequality constraints problem...\pageref{Sec.10.2}. 3.Problems with linear constraints...\pageref{Sec.10.3}
\end{footnotesize}

$\S$ 11. Construction of relaxation methods for nonconvex nonsmooth optimization  $\ldots$ \pageref{Sec.11}

$\S$ 12. On multiextremal optimization  $\ldots$ \pageref{Sec.12}

\begin{footnotesize}
\hspace{0.3cm}1.Combined algorithms...\pageref{Sec.12.1}. 2.The approximation method...\pageref{Sec.12.2}. 3.The smoothing method...\pageref{Sec.12.3}.
\end{footnotesize}

{\bf\smallskip
Chapter 4. Nonmonotonic methods with averaging descent directions $\ldots$} \pageref{Ch.4}

$\S$ 13. Methods with averaging of descent  directions  $\ldots$ \pageref{Sec.13}

$\S$ 14. Methods of averaged generalized gradients  $\ldots$ \pageref{Sec.14}

\begin{footnotesize}
\hspace{0.3cm}1.Stability of the generalized gradient method...\pageref{Sec.14.1}. 2.The method of averaged gradients...\pageref{Sec.14.2}. 3.Adjustment of parameters in the method of averaged gradients...\pageref{Sec.14.3}. 4.The heavy ball (Polyak) method...\pageref{Sec.14.4}. 5.The gully step (Nesterov) method...\pageref{Sec.14.5}.
\end{footnotesize}

{\bf\smallskip
Chapter 5. Solving extremal problems with Lipschitz functions under constraints $\ldots$} \pageref{Ch.5}

$\S$ 15. The conditional gradient method $\ldots$ \pageref{Sec.15}

$\S$ 16. The reduced gradient method  $\ldots$ \pageref{Sec.16}

$\S$ 17. The method of feasible directions  $\ldots$ \pageref{Sec.17}

$\S$ 18. Penalty function methods $\ldots$ \pageref{Sec.18}

\begin{footnotesize}
\hspace{0.3cm}1.Finite-difference methods for nonsmooth penalties...\pageref{Sec.18.1}. 2.Global features of non-smooth penalties...\pageref{Sec.18.2}.
\end{footnotesize}

$\S$ 19. Finite-difference method for minimizing Lipschitz functions under
constraints  $\ldots$ \pageref{Sec.19}

$\S$ 20. The finite-difference Errow-Hurwitz method with the trajectory averaging $\ldots$ \pageref{Sec.20}

{\bf\smallskip
Chapter 6. Random Lipschitz and random generalized differentiable
functions  $\ldots$} \pageref{Ch.6}

$\S$ 21. Measurable multivalued mappings  $\ldots$ \pageref{Sec.21}

\begin{footnotesize}
\hspace{0.3cm}1.Multivalued mappings...\pageref{Sec.21.1}. 2.Measurable multivalued mappings...\pageref{Sec.21.2}.
\end{footnotesize}

$\S$ 22. Random Lipschitz functions  $\ldots$ \pageref{Sec.22}

$\S$ 23. Random generalized differentiable functions and the calculus of
stochastic generalized gradients  $\ldots$ \pageref{Sec.23}

\begin{footnotesize}
\hspace{0.3cm}1.Random generalized differentiable functions...\pageref{Sec.23.1}. 2.Calculus of stochastic generalized gradients...\pageref{Sec.23.2}.
\end{footnotesize}

{\bf\smallskip
Chapter 7. Solving stochastic extremal problems  $\ldots$} \pageref{Ch.7}

$\S$ 24. Methods of averaged stochastic gradients  $\ldots$ \pageref{Sec.24}

\begin{footnotesize}
\hspace{0.3cm}1.The method of averaged stochastic gradients...\pageref{Sec.24.1}. 2.Method of stochastic generalized gradients...\pageref{Sec.24.2}. 3.Stochastic heavy ball method...\pageref{Sec.24.3}. 4.Stochastic gully step method...\pageref{Sec.24.4}.
\end{footnotesize}

$\S$ 25. Finite difference methods in stochastic programming  $\ldots$ \pageref{Sec.25}

\begin{footnotesize}
\hspace{0.3cm}1.The method of stochastic approximation...\pageref{Sec.25.1}. 2.A game stochastic problem...\pageref{Sec.25.2}. 3.Convergence rate...\pageref{Sec.25.3}.
\end{footnotesize}

$\S$ 26. The averaging operation  $\ldots$ \pageref{Sec.26}

$\S$ 27. Stochastic optimization based on the averaging operation  $\ldots$ \pageref{Sec.27}

\begin{footnotesize}
\hspace{0.3cm}1.Methods with averaging of descent directions...\pageref{Sec.27.1}. 2.Stochastic conditional gradient methods...\pageref{Sec.27.2}. 3.Convergence rate of the stochastic conditional gradient method...\pageref{Sec.27.3}. 4.Examples...\pageref{Sec.27.4}. 5.Stochastic reduced gradient method...\pageref{Sec.27.5}. 6.Stochastic methods of possible directions...\pageref{Sec.27.6}. 7.Stochastic methods for minimizing Lipschitz functions...\pageref{Sec.27.7}.
\end{footnotesize}

$\S$ 28. The stochastic finite-difference Arrow-Hurwicz method  $\ldots$ \pageref{Sec.28}

{\bf\smallskip
Bibliographic notes  $\ldots$} \pageref{BibNotes}

{\bf\smallskip
Bibliography  $\ldots$} \pageref{biblio}

\newpage

\begin{flushright}
\textbf{FOREWORD}
\end{flushright}
\label{FOREWORD}

This book is devoted to finite-dimensional problems of non-convex non-smooth optimization and numerical methods for their solution. The complexity theory of extremal problems states that nonconvex problems as an object of study are very complex, and the complexity of obtaining a guaranteed solution increases exponentially with increasing dimensions. However, the great practical importance of these problems convinces researchers to develop numerical methods despite such discouraging warnings about their complexity.

The problem of nonconvexity  is studied in the book on two main models of nonconvex dependencies: these are the so-called generalized differentiable functions and locally Lipschitz functions. The study of non-smooth functions is not a tribute to the current fashion. Non-smooth functions naturally arise in various applications. In addition, they often appear in the theory of extremal problems itself due to the operations of taking the maximum and minimum, decomposition techniques, exact non-smooth penalties, and duality.

The considered models of nonconvexity are quite general and cover the majority of practically important optimization problems; they clearly show all the difficulties of non-convex optimization. The method of studying the generalized differentiable functions is that for these functions a generalization of the concept of gradient is introduced, a calculus is constructed, and various properties of nonconvex problems are studied in terms of generalized gradients. As for numerical methods, it is possible to extend the theory and algorithms of subgradient descent of convex optimization to problems with generalized differentiable functions.

Methods for solving Lipschitz problems are characterized by the fact that the original functions are approximated by smoothed ones and iterative minimization procedures are applied to them. With this approach, it is possible to approximate the gradients of smoothed functions by stochastic finite differences and thus to construct methods without calculating gradients. A similar approach can be justified in Lipschitz stochastic programming. In this case, various generalizations of the classical stochastic approximation method are obtained for solving constrained Lipschitz stochastic programming problems.

To date, the mathematical theory of convex optimization has acquired complete features: a theory of complexity of convex programming problems has been developed and efficient numerical methods have been constructed. The theory of non-convex optimization has been developed to a much less extent, however, it can be seen that individual components of the theory (such as necessary optimality conditions) are advanced quite far, and the remaining components need further development. First of all, this concerns the development of non-convex optimization methods that overcome non-smoothness, non-convexity, and ravine functions, as well as taking into account the multi-extremality and possible stochastic nature of these problems.

In optimization theory, diverse classes of non-convex non-smooth functions were introduced and studied. These are functions: maximum, pseudo-convex, quasi-convex, weakly convex, semi-convex, locally convex, semi-smooth, quasi-differentiable, locally Lipschitz, discontinuous, etc. In Ch. 1 a new class of non-differentiable functions is also studied, generalized differentiable functions. The question arises: why is this another class of functions introduced?

Generalized differentiable functions are a direct generalization of continuously differentiable and convex functions. For them, pseudogradients are constructed similar to subgradients of convex functions, and a differential expansion similar to that of differentiable functions takes place. At the same time, this class of functions is quite wide: it includes continuously differentiable, convex and concave, weakly convex and weakly concave, semi-smooth functions. In addition, the class of generalized differentiable functions is closed under the operations of taking maximum, minimum, and superposition and taking the mathematical expectation. There are simple formulas (calculus) for calculating pseudogradients of complex generalized differentiable functions, similar to formulas for calculating ordinary gradients of complex differentiable functions, while for other classes of non-smooth functions, the construction of such a calculus encounters significant difficulties.

For generalized differentiable functions, the necessary optimality conditions are easily derived, and an analog of the mean value theorem takes place. And, finally, the most important point is that in order to minimize generalized differentiable functions, by analogy with the convex case, it is possible to construct numerical methods that use the pseudogradients of these functions.

The papers [145, 146] introduce the concepts of “generalized gradient” and “generalized directional derivative” for a fairly wide class of locally Lipschitz functions. It turns out that the class of locally Lipschitz functions is so wide that the formulas for calculating the generalized gradients of complex functions (for example, for the maximum, sum, and superposition functions) are not applicable. This means that generalized gradient descent methods cannot be widely used to minimize Lipschitz functions.

Ch. 2 shows that stochastic finite difference procedures lead to success, which do not require the calculation of generalized gradients, while ordinary finite difference approximations do not work. The main idea of choosing the direction of descent is that the final difference is calculated not at the current point but at a close, randomly chosen point, and the vector of mathematical expectation of the constructed finite differences is close to the set of generalized gradients. It is shown that the constructed finite difference methods find stationary points of Lipschitz functions. From the point of view of the number of calculations, it is important to note that for a maximum of functions, the finite difference is built on the function on which the maximum is achieved. In the case when one calculation of the goal function takes significant time, the construction of finite differences using random directions leads to a significant reduction in the number of calculations of the goal function.

An analysis of the efficiency of finite difference methods on convex programming problems shows that the replacement of subgradients by the introduced finite differences, roughly speaking, leads to a situation where subgradients are calculated with some random noise.

Among the methods of nonsmooth optimization, the method of generalized gradient descent [136] is attractive due to its simplicity, economy, certain optimality for convex problems of high dimension [73], and convergence in the nonconvex cases [4, 35, 80, 88].

In Ch. 3 the generalized gradient method is extended to problems of local minimization of generalized differentiable functions under constraints. In nonsmooth convex optimization, relaxation $\epsilon$-subgradient methods have become widespread. Ch. 3 offers somewhat different relaxation methods suitable for minimizing non-convex non-smooth functions.

Generalized gradient methods make it possible to find only local minima of functions, so a separate section is devoted to some approaches that make it possible to overcome the multi-extremality of problems. In particular, schemes of combined algorithms are studied, in which the developed local algorithms are used to approximately find local minima, and global algorithms are used to get out of zones of attraction of local extremums. We propose an algorithm for multivariate global search under constraints, generalizing the method [96]. The global properties of the function smoothing procedure widely used in this book are investigated. Note also that the heavy ball and gully step methods considered in Ch. 4 have certain global properties.

It has long been observed that (generalized) gradient methods converge slowly in gully situations. One way to overcome this drawback is to average the gradients (or corresponding finite differences) obtained at previous iterations in order to find the direction along the gully. In this connection, Ch. 4 studies non-monotonic gradient averaging (finite difference) methods for minimizing non-convex nonsmooth functions under constraints, which retain many of the advantages of gradient methods but already have anti-gully properties. It turns out that the appropriately generalized well-known heavy ball [100] and gully step [14, 74] methods fall into this class. Averaged gradient methods can also be interpreted as procedures for minimization of a sequence of some averaged functions.

Ch. 5 investigates numerical methods for minimizing Lipschitz functions under constraints. The methods of conditional and reduced gradients are studied. For problems with nonlinear constraints, methods of possible directions and penalty functions are constructed. To minimize Lipschitz functions, a very simple method of the penalty function type is proposed, in which the penalty coefficient is determined by an analytical expression.

An extensive field of application of the ideas of non-smooth optimization is stochastic programming. The stochastic problem differs from a deterministic non-linear programming problem in that the values of functions and their gradients in a general case cannot be calculated exactly in a finite time and, thus, are known only with some random error.

The first methods of stochastic programming go back to the methods of stochastic approximation. To solve problems of convex stochastic programming, methods of stochastic quasi-gradients were proposed [41]; important improvements were made in [73].

Significant difficulties arise in the way of solving nonconvex stochastic problems. In order to substantiate the stochastic analogues of gradient minimization procedures, first of all, it is necessary to show the commutability of the signs of generalized differentiation and taking the mathematical expectation. This fact is shown in Ch. 6 for generalized differentiable functions. On the basis of the apparatus of measurable multivalued mappings, a calculus of random pseudogradients of complex generalized differentiable functions is constructed, which is completely analogous to the calculus of pseudogradients of deterministic functions. Therefore, for these functions, the methods of stochastic generalized gradients can be substantiated.

For Lipschitz functions, the situation is different: the gradient procedures cannot be justified, since the generalized gradient of the expectation function is only included in the expectation of the (random) generalized gradient of the function under the sign of the expectation. In this regard, the developed finite-difference methods of Lipschitz stochastic programming, which largely generalize the Kiefer-Wolfowitz stochastic approximation methods, are of great importance.

In Ch. 7 we study stochastic analogs of the previously considered methods for solving deterministic problems.

\newpage
\begin{flushright}
\textbf{CHAPTER 1}
\label{Ch.1}

\textbf{ELEMENTS OF THE THEORY OF NONCONVEX ANALYSIS}
 \underline{\hspace{15cm}}

\end{flushright}
\bigskip\bigskip\bigskip\bigskip\bigskip\bigskip

\section*{$\S$ 1. Generalized differentiable functions}
\label{Sec.1}
\setcounter{section}{1}
\setcounter{definition}{0}
\setcounter{equation}{0}
\setcounter{theorem}{0}
\setcounter{lemma}{0}
\setcounter{remark}{0}
\setcounter{corollary}{0}
\setcounter{example}{0}
\numberwithin{equation}{section}

In this section, we introduce a class of nonconvex nonsmooth functions, which are here called generalized differentiable. These functions have many good optimization properties. They are a natural generalization of continuously differentiable and convex functions, for them some pseudogradients are introduced, similar to subgradients of convex functions, and a differential expansion involving pseudogradients is valid.

Generalized differentiable functions satisfy the local Lipschitz condition. They, generally,  have no directional derivatives but are continuously differentiable almost everywhere. Pseudogradients of the functions are not uniquely determined. For one function, there is a whole family of pseudo-gradient mappings, the inclusion minimal of which coincides with the Clarke [145] subdifferential mapping of this function. At the same time, any pseudogradient mapping defines a generalized differentiable function up to a constant.

Continuously differentiable, convex and concave, weakly convex and weakly concave [91], semi-smooth [167] functions are also generalized differentiable, and their gradients, subgradients, quasigradients, and generalized gradients are also pseudogradients.

The class of generalized differentiable functions is closed under the finite operations of taking the maximum, minimum, and superposition, and also, as will be shown in $\S$ 23, under taking the mathematical expectation. For complex generalized differentiable functions, there is a pseudo-gradient calculus.

For generalized differentiable functions, the necessary conditions for an extremum are written out ($\S$ 3), and the generalized Lagrange theorem on finite increments is valid. Finally, as will be shown in subsequent chapters, generalized gradient methods, including stochastic ones, can be developed to minimize generalized differentiable functions.

In this section, without comment, some general information from convex analysis, the theory of multivalued mappings ($\S$ 21), and locally Lipschitz functions ($\S$ 2) is used.
\begin{definition}
\label{def:1.1}
A function $f: E_n \rightarrow E_1$ is said to be generalized differentiable at a point $x \in E_n$, if in some neighborhood of the point $x$ there is defined a set-valued mapping $G_f$ upper semicontinuous in $x$ such that its values $G_f(y )\left(y \in E_n\right)$ are non-empty bounded convex closed sets, and in a neighborhood of the point $x$ the following expansion holds true
\begin{equation} \tag{1.1}
    f(y)=f(x)+\left<g, y-x\right>+o(x, y, g),
\end{equation}
where $g \in G_f(y)$, $\left<\cdot,\cdot\right>$ is the scalar product of two vectors, and the residual function $o(x, y, g)$ satisfies the condition
\begin{equation} \tag{1.2}
\lim _{k \rightarrow \infty}\left(o\left(x, y^k, g^k\right) /\left\|y^k-x\right\|\right)=0
\end{equation}
for any sequences $y^k \rightarrow x\left(y^k \neq x\right)$, $g^k \in G_f\left(y^k\right)$ or, equivalently,
\begin{equation} \tag{1.3}
\underset{\varepsilon \rightarrow+0}{\lim}\;
\underset{\{y \mid\|y-x\| \le \varepsilon, y \neq x\}} {\sup} \;
\underset{g \in G_f(y)}{\sup } \frac{|o(x, y, g)|}{\|y-x\|}=0 .
\end{equation}

A function is called generalized differentiable in a domain, if it is generalized differentiable at every point in the domain. The vectors $g \in G_j(y)$ will be called pseudo-gradients (generalized gradients) of the function $f$ at the point $y$.
\end{definition}

If condition (1.3) is satisfied for the function $o(x, y, g)$, then it is obvious that (1.2) also holds. The converse is easily proven by contradiction.
\begin{definition}
\label{def:1.2}
A function $f: E_n \rightarrow E_1$ is called locally Lipschitzian if for any compact set $K \subset E_n$ there exists a Lipschitz constant $L_K$ such that
$$
\left|f\left(x^{\prime}\right)-f\left(x^{\prime \prime}\right)\right| \le L_K\left\|x^{\prime}-x^{\prime \prime}\right\| \quad \forall x^{\prime}, x^{\prime \prime} \in K .
$$
\end{definition}
\begin{theorem}
\label{th:1.1}
Generalized differentiable functions are locally Lipschitz and hence continuous.
\end{theorem}

{\it P r o o f.} Let $f(x)$ be a generalized differentiable function, $K$ be a compact set in $E_n$, and let $co K$ be its convex hull. For $y, x \in K$ we represent
\begin{eqnarray}
f(y)-f(x)&=&\left<g, y-x\right>+o(x, y, g) \nonumber \\
& =&\left(\left<g, \frac{y-x}{\|y-x\|}\right>+\frac{o(x, y, g)}{\|y-x\|}\right)\| y-x\|,\nonumber
\end{eqnarray}
where $g \in G_f(y)$. By virtue of condition (1.2), for any $\varepsilon$ there exists $\delta(x, \varepsilon)$ such that
$$
|o(x, y, g)| /\|y-x\| \le \varepsilon \quad \text { with } \quad\|y-x\| \le \delta(x, \varepsilon), \quad g \in G_f(y) .
$$
Since the mapping $G_f$ is upper semicontinuous and its values $G_f(x)$ are compact sets, it is locally bounded; so there is a number
$$
\Gamma_K=\sup \left\{\|g\| \mid g \in G_f(y), y \in \operatorname{co} K\right\}<+\infty .
$$
Thus, for each $x \in \mbox{co} K$, a $\delta$-neighborhood of the point $x$ is found such that
$$
|f(y)-f(x)| \le\left(\Gamma_K+\varepsilon\right)\|y-x\| \quad \text { with } \quad\|y-x\| \le \delta, \quad y \in \operatorname{co} K .
$$
Let now $x^{\prime}, x^{\prime \prime} \in K$; consider the segment
$$
M=\left\{x=x^{\prime}+t\left(x^{\prime \prime}-x^{\prime}\right), t \in[0,1]\right\} , \quad M \subset \operatorname{co} K .
$$
The above $\delta$-neighborhoods of points $x \in M$ cover $M$; therefore, there exists a finite cover with centers $x^1=x^{\prime}, x^2, \ldots, x^k=x^{\prime \prime}$, and, respectively, $t_1=0<t_2<\ldots <t_k=1$. It is easy to see that for neighboring points $x^j$ and $x^{j+1}$ in all cases
$$
\left|f\left(x^j\right)-f\left(x^{j+1}\right)\right| \le\left(\Gamma_K+\varepsilon\right)\left\|x^i-x^{i+1}\right\|, 
$$
from where
$$
\left|f\left(x^{\prime}\right)-f\left(x^{\prime \prime}\right)\right| \le \sum_{i=1}^{k-1}\left|f\left(x^j\right)-f\left(x^{i+1}\right)\right| \le\left(\Gamma_K+\varepsilon\right)\left\|x^{\prime}-x^{\prime \prime}\right\| .
$$
The theorem is proved.
\begin{theorem}
\label{th:1.2}
Continuously differentiable functions are generalized differentiable, and their gradients can be taken as pseudogradients.
\end{theorem}

{\it P r o o f.} Let a function $f(y) \;\left(y \in E_n\right)$ be continuously differentiable, $g(y)$ be its gradient. Then at each point $x \in E_n$ the expansion holds
\begin{eqnarray} 
f(y)&=&f(x)+\left<g(x), y-x\right>+o(\| y -x \|) \nonumber \\ 
& =&f(x)+\left<g( y), y-x\right>+o(x, y, g),\nonumber
\end{eqnarray} 
where the remainder
$$
o(x, y, g)=\left<g(x)-g(y), y-x\right>+o(\|y-x\|)
$$
satisfies the required smallness conditions (1.2) and (1.3).
\begin{definition}
\label{df:1.3}
[88]. A continuous function $f: E_n \rightarrow E_1$ is said to be weakly convex if for any $x$ there exists a non-empty set $G(x)$ of quasi-gradients $g$ such that for all $g \in G(x)$ the inequality holds
\begin{equation} \tag{1.4}
f(y)\geqslant f(x)+\left<g, y-x\right>+r(x, y),
\end{equation}
where $r(x, y) /\|y-x\| \rightarrow 0, $ when $\| y-x \mid \rightarrow 0, $ and $x, y$ belong to a closed bounded set.
\end{definition}

Obviously, convex functions are also weakly convex; in this case $r(x, y) \equiv 0$. It is known [91] that the set of quasi-gradients $G(x)$ is non-empty, bounded, convex and closed, the set-valued mapping $G: x \rightarrow G(x)$ is upper semicontinuous. Weakly convex functions have directional derivatives.
\begin{theorem}
\label{th:1.3}
Convex and weakly convex functions are 
generalized differentiable, and their subgradients and quasigradients
can be taken as pseudogradients of these functions.
\end{theorem}

{\it P r o o f.} Here the role of $G_f$ is played by $G$. Let us prove the expansion (1.1) in a neighborhood of a point $x$. By the definition of a weakly convex function, inequality (1.4) holds. Let's write this inequality for points $y$ :
\begin{equation} \tag{1.5}
f(x) \geqslant f(y)+\left<g(y), x-y\right>+r(y, x).
\end{equation}
We introduce the notation
$$
o(x, y, g)=f(y)-f(x)-\left<g, y-x\right>.
$$
From (1.4) and (1.5) it follows
$$
\left<g(x)-g(y), y-x\right>+r(x, y) \le o(x, y, g(y)) \le-r(y, x),
$$
where $g(x) \in G(x)$, $g(y) \in G(y)$. We choose $g(x)$ such that
$$
\|g(x)-g(y)\|=\rho(g(y), \;\;\;G(x)) \equiv \text { inf}_g\{\|g(y)-g\|:\, g \in G(x)\} .
$$
Then
\begin{equation} \tag{1.6}
-\rho(g(y), G(x))\|y-x\|+r(x, y) \le o(x, y, g(y)) \le-r(y, x).
\end{equation}
Due to the upper semicontinuity of the mapping $G$ and the definition of a weakly convex function, $\rho(g(y), G(x)) \rightarrow 0$,
$$
r(x, y) /\|y-x\| \rightarrow 0, \quad r(y, x) /\|y-x\| \rightarrow 0
$$
for $\|y-x\| \rightarrow 0$ uniformly in $g(y) \in G(y)$. Then the required property $o(x, y, g)$ follows from (1.6).
\begin{definition}
\label{df:1.4}
[167]. A function $f: E_n \rightarrow E_1$ is called semi-smooth, if:

1) it is locally Lipschitz;

2) for any $x \in E_n$, any sequences $y^k \rightarrow x\left(\left(y^k-\right.\right.$ $\left.x) /\left\|y^ k-x\right\| \rightarrow d\right)$ and $g^k \in \partial f\left(y^k\right)$ there is a limit of scalar products $\lim _{k \rightarrow \infty}\left<g^k,d\right>$.
\end{definition}

Here $\partial f(y)$ is the Clarke subdifferential (set of generalized gradients) of the locally Lipschitz function $f(y)$. The subdifferential $\partial f$ is studied in detail in \S 2.

Obviously, for the same direction $d$, all inner products $\left<g^k, d\right>$ in Definition 1.4 have the same limit. It is easy to show, using condition 2) and the mean value Theorem 2.9, that semi-smooth functions have derivatives $f^{\prime}(x ; d)$ in the direction $d$ and
$$
\lim _{k \rightarrow \infty}\left<g^k, d\right>=f^{\prime}(x ; d)
$$
(see Theorem 1.14).
\begin{theorem}
\label{th:1.4}
Semi-smooth functions are generalized differentiable, and their generalized Clarke gradients can be taken as
pseudo-gradients.
\end{theorem}

{\it P r o o f.} Here the role of $G_f$ is played by a closed-valued upper semicontinuous mapping $\partial f$ (see $\S$ 2). Let us prove the expansion (1.1). Let's assume the opposite. Then there are such $x \in E_n$, $\varepsilon>0$, $y^k \rightarrow x$ $\left(\left(y^k-x\right) /\left\|y^k-x\right\| \rightarrow d\right)$, $g^k \in \partial f\left(y^k\right)$, that
$$
\left|\frac{f\left(y^k\right)-f(x)}{\left\|y^k-x\right\|}-\left<g^k, \frac{y^k-x} {\left\|y^k-x\right\|}\right>\right| \geq \varepsilon>0 .
$$
Using Theorem 2.9 on finite increments for locally Lipschitz functions, we rewrite the previous inequality as
$$
\left|\left<\overline{g}^{k}, \frac{y^{k}-x}{\left\|y^{k}-x\right\|}\right>-\left<g^{k}, \frac{y^{k}-x}{\left\|y^{k}-x\right\|}\right>\right| \geq \varepsilon>0,
$$
where $\bar{g}^{k} \in \partial f\left(z^{k}\right), \quad z^{k}=x+\theta^{k}\left(y^{k}-x\right), \quad 0<\theta^{k}<1$. Then
\[
 \left(z^{k}-x\right) /\left\|z^{k}-x\right\|=\left(y^{k}-x\right) /\left\|y^{k}-x\right\| \rightarrow d, 
\]
\begin{equation}
\left|\left<\bar{g}^{k}, \frac{z^{k}-x}{\left\|z^{k}-x\right\|}\right>-\left<g^{k}, \frac{y^{k}-x}{\left\|y^{k}-x\right\|}\right>\right| \geq \varepsilon>0 \;. \tag{1.7}
\end{equation}
Both scalar products in (1.7) have, by virtue of the definition of a semismooth function, the same limit as $k \rightarrow+\infty$, whence we get a contradiction: $0 \geqslant \varepsilon>0$. The theorem has been proven.

Thus, the class of generalized differentiable functions contains continuously differentiable, convex, and semismooth functions. Let us now show that it is closed with respect to the finite operations of maximization and superposition; hence it will follow that concave and weakly concave functions [91], as well as the minimum function of a finite number of generalized differentiable functions, are generalized differentiable.

\begin{theorem}
\label{th:1.5}
Let $f_{i}(y) \quad\left(y \in E_{n}, i=1, \ldots, m\right)$ are generalized differentiable functions, $G_{f_{i}}(y)$ are their pseudogradient maps. Then the maximum function
$$
f(y)=\max \left(f_{1}(y), \ldots, f_{m}(y)\right)
$$
is generalized differentiable and
\begin{equation}
G_{f}(y)=\operatorname{co}\left\{g \in E_{n} \mid g \in G_{f_{l}}(y),\; f_{i}(y)=f(y)\right\} .\tag{1.8}
\end{equation}
\end{theorem}

{\it P r o o f.} It suffices to consider a maximum of two functions. The mapping $G_{f}$ is convex-valued. Let us show that it is closed and locally bounded; hence its upper semicontinuity will follow. Closedness means that if $y^{k} \rightarrow y, \quad g^{k} \in G_{f}\left(y^{k}\right)$ and $g^{k} \rightarrow g$ it follows $g \in G_{f}(y)$. Let us represent
$$
g^{k}=\sum_{i=1}^{m} \lambda_{i}^{k} g_{l}^{k}, \quad \lambda_{i}^{k} \geqslant 0, \quad \sum_{i=1}^{m} \lambda_{i}^{k}=1,
$$
where $\lambda_{i}^{k}=0$, if
$$
f_{i}\left(y^{k}\right)<f\left(y^{k}\right)
$$
Due to the boundedness of $\lambda_{i}^{k}$, closedness and local boundedness of mappings $G_{f_{i}}$ and the number $m$ is finite, there are subsequences
$$
\lambda_{i}^{k_{s}} \rightarrow \lambda_{i}, \quad g_{i}^{k_{s}} \rightarrow g_{i}, \quad s=0,1, \ldots
$$
and
$$
\begin{gathered}
\lambda_{i} \geq 0, \quad \sum_{i=1}^{m} \lambda_{i}=1, \quad g_{i} \in G_{f_{i}}(y), \\
g=\lim _{s \rightarrow \infty} g^{k_{s}}=\sum_{i=1}^{m} \lambda_{i} g_{i} .
\end{gathered}
$$
If for some $i$
$$
f_{i}(y)<f(y),
$$
then this also holds in some neighborhood of the point $y$, so for all sufficiently large $k_{s}$ holds $\lambda_{i}^{k_{s}}=0$; hence,
$$
\lambda_{i}=\lim _{s \rightarrow \infty} \lambda_{i}^{k_{s}}=0.
$$
Thus, $g \in G_{f}(y)$, and closedness of $G_{f}$ is proven. Each of mapings $G_{f_{i}}$ is compact-valued and upper semicontinuous and hence locally bounded. Obviously, then $G_{f}$ is also locally bounded. And from the closedness and local boundedness it follows the upper semicontinuity of the mapping $G_{f}$.

Now it is necessary to establish the validity of expansion (1.1), i.e., it is necessary to show that for any $\varepsilon>0$ there is such $\delta(\varepsilon)>0$ that
\begin{equation}
-\varepsilon \le \frac{f(y)-f(x)-\left<g, y-x\right>}{\|y-x\|} \le \varepsilon \tag{1.9}
\end{equation}
for all $\{y \mid\|y-x\| \le \delta, y \neq x\}$ and all $g \in G_{f}(y)$. Due to the generalized differentiability of $f_{1}(y)$ and $f_{2}(y)$, for any given $\varepsilon>0$ there is $\delta_{1}(\varepsilon)>0$ and $\delta_{2}(\varepsilon)>0$ such that
\begin{equation}
-\varepsilon \le \frac{f_{1}(y)-f_{1}(x)-\left<g^{1}, y-x\right>}{\|y-x\|} \le \varepsilon \tag{1.10}
\end{equation}
for all $\left\{y \mid\|y-x\| \le \delta_{1}(\varepsilon), y \neq x\right\}$ and all $g^{1} \in G_{f_{1}}(y)$ and that
\begin{equation}
-\varepsilon \le \frac{f_{2}(y)-f_{2}(x)-\left<g^{2}, y-x\right>}{\|y-x\|} \le \varepsilon \tag{1.11}
\end{equation}
for all $\left\{y \mid\|y-x\| \le \delta_{2}(\varepsilon), y \neq x\right\}$ and all $g^{2} \in G_{f_{2}}(y)$.

If
$$
f_{1}(x) \neq f_{2}(x),
$$
for example, $f_{1}(x)>f_{2}(x)$, then the same inequality persists in some $\delta_{3}$-neighbourhood of the point $x$, and there $f(y) \equiv f_{1}(y)$. For the given $\varepsilon$ let's take
$$
\delta(\varepsilon)=\min \left(\delta_{1}(\varepsilon), \delta_{3}(\varepsilon)\right),
$$
then the required inequality (1.9) follows from (1.10).

Let now
$$
f_{1}(x)=f_{2}(x)=f(x).
$$
For an arbitrary $\varepsilon$ we take
$$
\delta(\varepsilon)=\min \left(\delta_{1}(\varepsilon), \delta_{2}(\varepsilon)\right)
$$
and verify inequality (1.9) for all $\{y \mid\|y-x\| \le \delta$, $y \neq x\}$ and all $g \in G_{f}(y)$, where $G_{f}(y)$ is constructed according to (1.8).

Let $y$ be such that $\|y-x\| \le \delta$ and, for example, $f_{1}(y)>f_{2}(y)$. Then
$$
f(y)=f_{1}(y), \quad g=g^{1} \in G_{f_{1}}(y),
$$
and (1.9) is valid due to (1.10).

Consider $y$ such that $\|y-x\| \le \delta$ and $f(y)=f_{1}(y)=f_{2}(y)$. Then, according to (1.8),
$$
g=\lambda_{1} g^{1}+\lambda_{2} g^{2}
$$
where $\lambda_{1}+\lambda_{2}=1,$ $g^{1} \in G_{f_{1}}(y),$ $g^{2} \in G_{f_{2}}(y)$. Substituting these $y,$ $g^{1},$ $g^{2}$ in (1.10) and (1.11), multiplying (1.10) and (1.11) by $\lambda_{1}$ and $\lambda_{2}$, and adding, we obtain (1.9). The theorem has been proven.
\begin{theorem}
\label{th:1.6}
Let $f_{0}(z)\left(z \in E_{m}\right),$ $f_{i}(y) \left(y \in E_{n}, i=1, \ldots, m
\right)$ be generalized differentiable functions. Then the compound function
$$
f(y)=f_{0}(z(y))=f_{0}\left(f_{1}(y), f_{2}(y), \ldots, f_{m}(y)\right)
$$
is also generalized differentiable and
\setcounter{equation}{11}
\begin{eqnarray}
G_{f}(y)&=&\operatorname{co}\left\{g \in E_{n} \mid g=\left[g^{1} \ldots g^{m}\right] g^{0},\right. \;  \nonumber\\
& &\left.\;\;\;\;\;g^{0} \in G_{f_{0}}(z(y)),\; g^{i} \in G_{f_{l}}(y), \;i=1, \ldots, m\right\},
\label{eqn:1.12}
\end{eqnarray}
where $\left[g^{1} \ldots g^{m}\right]$ is $n \times m$-matrix, composed of vectors (columns $g^{1}, \ldots, g^{m}),$ $z(y)=\left(f_{1}(y), \ldots, f_{m}(y)\right)$.
\end{theorem}

Formula (1.12) is a generalization of the chain rule for differentiating complex functions in mathematical analysis. Thus, the class of generalized differentiable functions forms a linear space and is closed under superpositions.

{\it P r o o f of Theorem 1.6.} The mapping $G_{f}$ is convex-valued. Let us show that it is locally bounded and closed; hence its upper semicontinuity will follow. Note that the mappings $G_{f_{i}} (i=1,2, \ldots, m)$ are locally bounded since they are compact-valued and upper semicontinuous. This is equivalent to their limitedness in every compact. Mapping $G_{f_{0}}(z(y))$ is locally bounded as a superposition of a continuous mapping $z(y)$ and a bounded in every compact  mapping $G_{f_{0}}(z)$. Then it is obvious from (1.12) that $G_{f}$ is also locally bounded. For $G_{f}$ to be closed, one must show that from $y^{k} \rightarrow y,$ $g^{k} \in G_{f}\left(y^{k}\right)$ and $g^{k} \rightarrow g$ it follows $g \in G_{f}(y)$. By virtue of the Carathéodory theorem, $g^{k}$ can be represented as a convex combination
$$
\begin{aligned}
& g^{k}=\sum_{j=1}^{n+1} \lambda_{j k}\left[\left(g^{1}\right)^{j k} \ldots\left(g^{m}\right)^{j k}\right]\left(g^{0}\right)^{j k}, \quad \lambda_{j k} \geqslant 0, \quad \sum_{j=1}^{n+1} \lambda_{j k}=1, \\
& \left(g^{i}\right)^{j k} \in G_{f_{i}}\left(y^{k}\right), \quad i=1, \ldots, m, \quad j=1, \ldots, n+1, \\
& \left(g^{0}\right)^{j k} \in G_{f_{0}}\left(z\left(y^{k}\right)\right), \quad j=1, \ldots, n+1,
\end{aligned}
$$
Sequences $\left\{\left(g^{i}\right)^{j k}\right\}_{k=1}^{\infty} \quad(i=0,1, \ldots, m, j=1, \ldots, n+1)$ are bounded due to the local boundedness of the mappings $G_{f_{i}}(i=1, \ldots$ $\ldots, m)$ and $G_{f_{0}}(z(y))$. 
Sequences $\left\{\lambda_{j k}\right\}_{k=1}^{\infty}(j=1, \ldots, n+1)$ are also bounded. Since $m$ and $n$ are finite, there exists a subsequence of indices $k_{s}$ such that
$$
\lambda_{j k_{s}} \rightarrow \lambda_{j}, \quad\left(g^{i}\right)^{j k_{s}} \rightarrow\left(g^{i}\right)^{j}, \quad i=0, \ldots, m, \quad j=1, \ldots, n+1
$$
wherein
$$
\lambda_{j} \geqslant 0, \quad \sum_{j=1}^{n+1} \lambda_{j}=1, \quad\left(g^{i}\right)^{j} \in G_{f_{i}}(y), \quad i=1, \ldots, m, \quad j=1, \ldots, n+1,
$$
and also
$$
\left(g^{0}\right)^{j} \in G_{f_{0}}(z(y))
$$
due to the closedness of mappings $G_{f_{i}}(y) \quad(i=1, \ldots, m)$ and $G_{f_{0}}(z(y))$. That's why
$$
g=\lim _{s \rightarrow \infty} g^{k_{s}}=\sum_{j=1}^{n+1} \lambda_{j}\left[\left(g^{1}\right)^{i} \ldots\left(g^{m}\right)^{i}\right]\left(g^{0}\right)^{j} \in G_{f}(y)
$$

The closedness of $G_{f}$ is proved. So $G_{f}$ is locally bounded and closed; hence it is upper semicontinuous.

Let us now prove the expansion (1.1). Let us represent
\begin{eqnarray}
f(y)&=&f_{0}(z(y))=f_{0}\left(f_{1}(y), \ldots, f_{m}(y)\right)= \nonumber\\
& =&f_{0}(z(x))+\left<g^{0}(z(y)), z(y)-z(x)\right>+o_{0}\left(z(x), z(y), g^{0}(z(y))\right) \nonumber\\
& =&f(x)+\sum_{i=1}^{m} g_{(i)}^{0}(z(y))\left(f_{i}(y)-f_{i}(x)\right)+o_{0}\left(z(x), z(y), g^{0}(z(y))\right) \nonumber\\
& =&f(x)+\sum_{i=1}^{m} g_{(i)}^{0}(z(y))\left(\left<g^{i}(y), y-x\right>+o_{i}\left(x, y, g^{i}\right)\right) \nonumber\\
& &+o_{0}\left(z(x), z(y), g^{0}(z(y))\right)=f(x)+\left<\left[g^{1} \ldots g^{m}\right] g^{0}(z(y)), y-x\right> \nonumber\\
&&+\sum_{i=1}^{m} g_{(i)}^{0}(z(y)) o_{i}\left(x, y, g^{i}\right)+o_{0}\left(z(x), z(y), g^{0}(z(y))\right),\label{eqn:1.13}
\end{eqnarray}
where
$$
\begin{aligned}
g^{0}(z(y))= & \left(g_{(1)}^{0}(z(y)), \ldots, g_{(m)}^{0}(z(y))\right) \in G_{f_{0}}(z(y)), \\
& g^{i}(y) \in G_{f_{i}}(y), \quad i=1, \ldots, m, \\
\end{aligned}
$$
$o_{i}\left(x, y, g^{i}\right) /\|y-x\| \rightarrow 0 \text { uniformly across } y \rightarrow x \text { and } g^{i} \in G_{f_{i}}(y).$

Representation (1.13) implies the required expansion (1.1)-(1.3). Indeed, we fix some $\delta$-neighbourhood of the point $x$. Because $G_{f_{0}}(z(y))$ is bounded in every compact set, then there exists a constant $c<\infty$ such that for every $\{y \mid\|y-x\| \le \delta\}$ and all $g^{0} \in G_{f_{0}}(z(y))$ there will be $\|g\| \le c$. by virtue of Theorem 1.1, mapping $z(y)$ is locally Lipschitz; so for every $\{y \mid\|y-x\| \le \delta\}$ we have
$$
\|z(y)-z(x)\| /\|y-x\| \le L<\infty .
$$
Let now $g \in G_{f}(y),\|y-x\| \le \delta$. Denote
\begin{equation}
o(x, y, g)=f(y)-f(x)-\left<g, y-x\right>. \label{eqn:1.14}
\end{equation}
By virtue of the Carathéodory theorem, the vector $g$ can be represented as a convex combination:
$$
\begin{gathered}
g=\sum_{j=1}^{n+1} \lambda_{j}\left[\left(g^{1}\right)^{j} \ldots\left(g^{m}\right)^{j}\right]\left(g^{0}\right)^{j}, \quad \lambda_{j} \geqslant 0, \quad \sum_{j=1}^{n+1} \lambda_{j}=1, \\
\left(g^{i}\right)^{j} \in G_{f_{j}}(y), \quad i=1, \ldots, m, \quad j=1, \ldots, n+1 ; \quad\left(g^{0}\right)^{j} \in G_{f_{0}}(z(y)) .
\end{gathered}
$$
Writing representation (1.13) for each $j$ and summing them with weights $\lambda_{j}$, for function (1.14), we obtain
\[
o(x, y, g)=\sum_{j=1}^{n+1} \lambda_{j} \sum_{i=1}^{m}\left(g_{(i)}^{0}\right)^{j} o_{i}\left(x, y,\left(g^{i}\right)^{j}\right)+\sum_{j=1}^{n+1} \lambda_{j} o_{0}\left(z(x), z(y),\left(g^{0}\right)^{j}\right), 
\]
\begin{eqnarray}
\frac{|o(x, y, g)|}{\|y-x\|} \le c \sum_{j=1}^{n+1} \lambda_{j} \sum_{i=1}^{m} \frac{\left|o_{i}\left(x, y,\left(g^{i}\right)^{j}\right)\right|}{\|y-x\|}+\sum_{j=1}^{n+1} \lambda_{j} \frac{\left|o_{0}\left(z(x), z(y),\left(g^{0}\right)^{j}\right)\right|}{\| y-x\|} \nonumber\\
\le c \sum_{j=1}^{n+1} \lambda_{j} \sum_{i=1}^{m} \frac{\left|o_{i}\left(x, y,\left(g^{i}\right)^{j}\right)\right|}{\|y-x\|}+L \sum_{j=1}^{n+1} \lambda_{j} \frac{\left|o_{0}\left(z(x), z(y),\left(g^{0}\right)^{j}\right)\right|}{\|z(y)-z(x)\|} .\;\;\;\label{eqn:1.15}
\end{eqnarray}
Using the obtained estimate for $|o(x, y, g)| /\|y-x\|$, let us check the fulfillment of condition (1.3). Let an arbitrary $\varepsilon>0$ be given. Due to the generalized differentiability of $f_{i}(y)(i=1, \ldots, m)$ there exists such $\delta_{1} \le \delta$
that
\begin{equation}
\frac{\left|o_{i}\left(x, y, g^{i}\right)\right|}{\|y-x\|} \le \frac{\varepsilon}{2 c m} \quad \text { at } \quad\|y-x\| \le \delta_{1}, \quad g^{i} \in G_{f_{i}}(y) .\label{eqn:1.16}
\end{equation}
By virtue of the generalized differentiability of $f_{0}(z)$, there exists such $\varepsilon_{1}$ that
\begin{equation}
\frac{\left|o_{0}\left(z(x), z, g^{0}\right)\right|}{\|z-z(x)\|} \le \frac{\varepsilon}{2 L} \quad \text { at } \quad\|z-z(x)\| \le \varepsilon_{1}, \quad g^{0} \in G_{f_{0}}(z) . \label{eqn:1.17}
\end{equation}
Since $z(y)$ is continuous, for a given $\varepsilon_{1}$ there is such $\delta_{2} \le \delta_{1}$ that
$$
\|z(y)-z(x)\| \le \varepsilon_{1} \quad \text { at } \quad\|y-x\| \le \delta_{2} .
$$
Thus, at $\|y-x\| \le \delta_{2}$ and $g \in G_{f}(y)$, from (1.15)-(1.17) it follows that
$$
|o(x, y, g)| / \| y-x \| \le \varepsilon .
$$
Condition (1.3) is verified, and the theorem is proved.
\begin{corollary}
\label{cor:1.1}
Concave and weakly concave [91] functions are generalized differentiable, and their subgradients and quasigradients can be taken as pseudogradients of these functions.
\end{corollary}
\begin{corollary}
\label{cor:1.2}
The minimum function of a finite number of generalized differentiable functions
$$
f(y)=\min \left\{f_{1}(y), \ldots, f_{m}(y)\right\}
$$
is generalized differentiable, and
\begin{equation}
G_{f}(y)=\operatorname{co}\left\{g \in E_{n} \mid g \in G_{f_{i}}(y), f_{i}(y)=f(y)\right\} .\label{eqn:1.18}
\end{equation}
\end{corollary}
\begin{corollary}
\label{cor:1.3}
Modulus function of a generalized differentiable function $f(x)$
$$
h(x)=|f(x)|=\max \{f(x),-f(x)\}
$$
is generalized differentiable.
\end{corollary}
\begin{corollary}
\label{cor:1.4}
The Lagrange function for a mathematical programming problem with generalized differentiable functions is generalized differentiable with respect to the set of direct and dual variables.
\end{corollary}

According to Theorem 1.1, generalized differentiable functions are locally Lipschitz. The following examples show that generalized differentiable functions generally have no directional derivatives, but not all Lipschitz functions (and even not all differentiable functions) are generalized differentiable.
\begin{example}
\label{exmp:1.1}
Consider a nontrivial example of a generalized differentiable function. Let us show that the function
$$
f(x)=\left\{\begin{array}{l}
(2 k-\ln |\ln | x||)|x|, \quad e^{-e^{2 k}} \le|x| \le e^{-e^{2 k-1}}, \\
(-2 k+\ln |\ln | x ||)|x|, \quad e^{-e^{2 k+1}} \le|x| \le e^{-e^{2 k}} \\
0, \quad x=0,
\end{array}\right.
$$
where $k=0,1, \ldots, \quad x \in E_{1},|x| \le e^{-1 / e}$, is generalized differentiable, but has no directional derivatives at point $x=0$.

Let us verify that the piecewise given function $f(x)$ is continuous. Let us introduce notation
$$
f^{+}(x)=\lim _{y \rightarrow x+0} f(y), \quad f^{-}(x)=\lim _{y \rightarrow x-0} f(y).
$$
It is necessary to check the continuity of $f(x)$ at the points $x=e^{-e^{n}} \quad(n=$ $=0,1, \ldots)$. At $n=2 k$ we have
$$
f^{+}\left(e^{-e^{2 k}}\right)=f^{-}\left(e^{-e^{2 k}}\right)=0.
$$
At $n=2 k+1$ we have
$$
f^{+}\left(e^{-e^{2 k+1}}\right)=f^{-}\left(e^{-e^{2 k+1}}\right)=e^{-e^{2 k+1}} \text {. }
$$
It remains to check the point $x=0$. It holds $0 \le f(x) \le|x|;$ that's why
$$
\lim _{x \rightarrow 0} f(x)=0.
$$
So $f(x)$ is continuous.

It is clear that $f(x)$ is piecewise continuously differentiable, and $(x>0)$
$$
f^{\prime}(x)=\left\{\begin{array}{cc}
2 k-\ln |\ln x|+\frac{1}{|\ln x|}, & e^{-e^{2 k}}<x<e^{-e^{2 k-1}}, \\
-2 k+\ln |\ln x|-\frac{1}{|\ln x|}, & e^{-e^{2 k+1}}<x<e^{-e^{2 k}}.
\end{array}\right.
$$

Denote $(x>0)$
$$
f_{+}^{\prime}(x)=\lim _{y \rightarrow x+0} f^{\prime}(y), \quad f_{-}^{\prime}(x)=\lim _{y \rightarrow x - 0} f^{\prime}(y)
$$
At $x=e^{-e^{2 k}}$ we have
$$
f_{+}^{\prime}(x)=e^{-2 k}=-f_{-}^{\prime}(x)
$$
At $x=e^{-e^{2 k+1}}$ we have
$$
f_{+}^{\prime}(x)=1-e^{-2 k-1}, \quad f_{-}^{\prime}(x)=1+e^{-2 k-1} .
$$

Let's define
$$
G_{f}(x)= \begin{cases}f^{\prime}(x), & x \neq e^{-e^{n}}, \quad n=0,1, \ldots, \\ \operatorname{co}\left\{f_{+}^{\prime}(x), f_{-}^{\prime}(x)\right\}, & x=e^{-e^{n}}, \quad n=0,1, \ldots, \\ {[0,1],} & x=0\end{cases}
$$

Closedness, convexity, boundedness of $G_{f}(x)$ and upper semicontinuity of $G_{f}$ is obvious. Differential representation (1.1) for $f(x)$, obviously takes place at break points $x=e^{-e^{n}} ;$ at these points $f(x)$ also has derivatives with respect to directions. It remains to show that representation (1.1) also holds at the point $x=0$, but there is no directional derivative. Function
$$
\left|\frac{o(x, y, g)}{y}\right|=\left|\frac{f(y)-f(0)}{y}-g(y)\right|=\left|\frac{f(y)}{y}-g(y)\right| \le \frac{1}{|\ln y|}
$$
tends to zero at $y \rightarrow+0$ uniformly across $g \in \operatorname{co}\left\{f_{+}^{\prime}(y), f_{-}^{\prime}(y)\right\}$. Therefore, representation (1.1) takes place at the point $x=0$. Consider the relation
$$
f(y) / y=(f(y)-f(0)) / y
$$
and show that it has no limit for $y \rightarrow 0$. Indeed,
$$
\lim _{k \rightarrow \infty} f\left(e^{-e^{2 k}}\right) / e^{-e^{2 k}}=0, \quad \lim _{k \rightarrow \infty} f\left(e^{-e^{2 k+1}}\right) / e^{-e^{2 k+1}}=1 .
$$

Therefore, the function $f(x)$ has no directional derivatives at the point ${x}=0$.
\end{example}
\begin{example}
\label{exmp:1.2}
Differentiable at each point function
$$
f(y)=y^{2} \sin (1 / y)
$$
is not generalized differentiable at zero. Note that it is locally Lipschitz but not continuously differentiable at zero.

Indeed, the function
$$
\frac{o\left(0, y, f^{\prime}(y)\right)}{y}=\cos \frac{1}{y}-2 y \sin \frac{1}{y}
$$
has no limit when $y \rightarrow 0$.
\end{example}
\begin{example}
\label{exmp:1.3}
Locally Lipschitz maximum function
$$
f(y)=\max _{k}\left(\frac{1}{2 \cdot 3^{k}}-\left|x-\frac{1}{3^{k}}\right|\right), \quad k=0,1, \ldots, \infty,
$$
is not generalized differentiable at zero. This example shows that the countable operation of taking maximum over generalized differentiable functions leads out of the class of generalized differentiable functions.

Indeed, $o\left(0,3^{-k},+1\right) / 3^{-k}$ is equal to -1 and does not tend to zero as $k \rightarrow \infty$.
\end{example}
\begin{theorem}
\label{th:1.7}
For generalized differentiable functions $f(x)$ ($\left.x \in E_{n}\right)$ the generalized Lagrange theorem on finite increments is valid, namely
$$f(y)-f(x)=\left<g, y-x\right>, \quad g \in G_{f}(x+\theta(y-x)), \quad 0<\theta<1.$$
\end{theorem}

{\it P r o o f.} Consider the function
$$
\varphi(t)=f(x+t(y-x))-t(f(y)-f(x))-f(x) .
$$
Generalized gradients of the function  $\varphi(t)$, according to (1.12), are
\begin{equation}
G_{\varphi}(t)=\left\{\left<g, y-x\right>-f(y)+f(x) \mid g \in G_{f}(x+t(y-x))\right\} .
\label{eqn:1.19}
\end{equation}
Note that $\varphi(0)=\varphi(1)=0$. Therefore, there is always an internal extreme point $t=\theta \in(0,1)$ of function $\varphi(t)$. According to the Theorem 3.1, $0 \in G_{\varphi}(\theta)$ regardless of whether it is the maximum or minimum extreme point. Then it follows from (1.19) that there exists a pseudogradient $g \in G_{f}(x+\theta(y-x))$ such  that
$$
\left<g, y-x\right>-f(y)+f(x)=0.
$$
The theorem has been proven.

Now consider the question: Is a pseudogradient mapping of a generalized differentiable function uniquely defined by Definition 1.1? We will show that for the same generalized differentiable function there exists a whole family of pseudogradient mappings satisfying Definition 1.1. However, among them there is a minimal inclusion mapping; it coincides with the Clarke subdifferential mapping of this function. In view of this, we should say that formulas (1.8), (1.12), (1.18) for calculating pseudo-gradients of complex functions specify only some, not necessarily minimal, pseudo-gradient mapping. That is why these formulas have equality and not inclusion in contrast to the corresponding formulas (2.9), (2.13), (2.17) for locally Lipschitz functions. However, the ambiguity of the definition of pseudo-gradient mappings is not an obstacle to their application both in theoretical problems and in construction of numerical methods for minimization of generalized differentiable functions. All pseudo-gradient mappings
are completely equivalent, since they only need to satisfy Definition 1.1.
\begin{theorem}
\label{th:1.8}
{If some pseudo-degenerate mapping of a generalized differentiable function is arbitrarily modified at a finite number of points such that it remains compact-valued, convex, and closed, then the modified mapping remains pseudo-degenerate one for that function.}
\end{theorem}

{\it P r o o f}. 
Indeed, the modified mapping is still locally bounded, closed, and hence semicontinuous from above. At the points of modification, the expansion (1.1) is preserved since the pseudogradients of these points essentially
are not used in (1.1). The decomposition (1.1) written down at other points
remains valid but may be for smaller neighborhoods of these points.
\begin{theorem}
\label{th:1.9}
{Let $f(x)$ and $h(x)$ be generalized differentiable functions. Then along with $G_f(x)$ the sets}
\[ G_f^{\prime}(x) =
  \begin{cases}
    G_f(x)\text{,}       & \quad f(x) \neq h(x),\\
    co\{G_f(x)\text{,} G_h(x)\} \text{,}  & \quad f(x) = h(x),
  \end{cases}
\]
\textit{ define another pseudogradient mapping of function $f(x)$.}
\end{theorem}

This theorem is proven similarly to Theorem 1.5. Here the extension of $G_f$ to $G^\prime_f$ turns out to be possible, since it does not violate the defining expansion (1.1) and other requirements for the pseudogradient mapping.

\begin{theorem}
\label{th:1.10}
The minimum inclusion pseudogradient mapping of a generalized differentiable function coincides with the subdifferential (generalized gradient) mapping according to Clarke 
(see $\S$ 2) of this function.
\end{theorem}

{\it P r o o f}. Recall that, according to Theorem 1.1, the generalized differentiable function $f(x)$ is locally Lipschitz, so there exists Clarke subdifferential $\partial f(x)$, which is the set of vectors $g\in E_n$ such that for any direction $d\in E_n$ there holds
\[\left<g,d\right>\leq \lim\limits_{\delta \to 0 } \sup_{\substack{ ||y-x|| \leq \delta \\
0 < \lambda \leq \delta}} [f(y+\lambda d)]\lambda^{-1}.
\]

For $g\in \partial f(x)$ by virtue of Theorem 1.7 and the semicontinuity from above of $G_f$, we can write
\[ \left<g,d\right>\leq \lim\limits_{\delta \to 0 } \sup_{\substack{ ||y-x|| \leq \delta \\
0 < \lambda \leq \delta}} \left<g(y+\bar{\lambda} d), d\right> \leq \sup_{g \in G_f(x)} \left<g, d\right>,
\]
where 
\[
g(y+\bar{\lambda} d)\in G_f(y+\bar{\lambda}d),\;\;\; 0 < \bar{\lambda}<\lambda\leq \delta.
\]

Due to convexity and closedness of $G_f(x)$ and arbitrariness of $d$, 
it follows that
\[
\partial f(x)\subset G_f(x).
\]
Then $\partial f$ together with $G_f$ satisfies expansion (1.1). In addition,
$\partial f(x)$ ia a nonempty bounded convex closed set, and
the mapping $x \to \partial f(x)$ is semicontinuous from above (cf. $\S$ 2). Hence,
$\partial f$ is the minimal inclusion pseudogradient mapping for
$f(x)$. The theorem is proved.

The following statements establish the relation between
ordinary and generalized differentiability of functions.
\begin{theorem}
\label{th:1.11}
{Let $f$ be a generalized differentiable function, $G_f$ is its pseudogradient mapping. If at point $x$ the set $G_f(x)$ consists of a single point, then the function $f$ is differentiable at $x$ (by Fr\'echet and by Gateaux).}
\end{theorem}

{\it P r o o f}. Denote $g(x) = G_f(x)$. The expansion (1.1) for $f$ is valid, and we represent it in the following form:
\[
f(y)=f(x)+\left<g(x), y-x\right>+\left<g-g(x), y-x\right>+o(x, y, g),
\]
where $g\in G_f(y)$. Note that the summand
\[
r(x,y)\equiv \left<g-g(x), y-x\right>+o(x,y,g)
\]
is in fact independent of $g$.
Thus,
\[
f(y) = f(x)+\left<g(x), y-x\right> + r(x,y).
\]

The estimation holds true
\[
\frac{|r(x,y)|}{\|y-x\|}\leq \sup_{g\in G_f(y)} \| g(x) -g \| + \sup_{g\in G_f(y)} \frac{|o(x,y,g)|}{\|y-x\|}, 
\]
whence, at a fixed $x$, by virtue of semicontinuity from above of $G_f$ and
condition (1.2), we get
\[
\lim_{y \to x} \frac{r(x,y)}{\|y-x\|} =0,
\]
i.e., the function $f (x)$ is differentiable at the point $x$ by Fr\'echet, and hence
also by Gateaux. The theorem is proved.
\begin{lemma}
\label{lem:1.1}
{Let $f(x) (x \in E_1)$ be a generalized differentiable function of one variable. Let us introduce a gradient discontinuity function}
\[
d(x) = \sup_{g\in G_f(x)} g-\inf_{g\in G_f(x)} g,
\]
\textit{then for any $\varepsilon>0$ and any bounded set $K \subset E_1$
the set of gradient discontinuity points $\{x\in K|d(x)\geq \varepsilon > 0\}$ is finite.}
\end{lemma}

{\it P r o o f.} Suppose the contrary. Then there exists
such an infinite sequence of points $\{x^k\}$ that $x^k\in K$ and $d(x^k) \geq \varepsilon > 0$. Let us choose a subsequence $x^{k_s} \to x (s= 0, 1,...)$. Define
\[
g^{k_s}_1 = \sup_{g\in G_f(x^{k_s})} g,\;\;\; g^{k_s}_2 = \inf_{g\in G_f(x^{k_s})} g.
\]
By virtue of the generalized differentiability of $f(x)$, we have
\[
f(x^{k_s}) = f(x) + \left<g^{k_s}_1, x^{k_s}-x\right>+o(x,x^{k_s},g^{k_s}_1),
\]
\[
f(x^{k_s}) = f(x) + \left<g^{k_s}_2, x^{k_s}-x\right>+o(x,x^{k_s},g^{k_s}_2),
\]
whereby
\[
\lim_{s\to \infty} \frac{o(x,x^{k_s},g^{k_s}_1)}{\|x^{k_s}-x\|} = 
\lim_{s\to \infty} \frac{o(x,x^{k_s},g^{k_s}_2)}{\|x^{k_s}-x\|} = 0.
\]

From this we get 
\[
\left<g^{k_s}_1-g^{k_s}_2, x^{k_s}-x\right> = o(x,x^{k_s},g^{k_s}_2) -o(x,x^{k_s},g^{k_s}_1)
\]
and we come to the contradiction:
\[
0<\varepsilon \leq g^{k_s}_1-g^{k_s}_2 \leq \frac{o(x,x^{k_s},g^{k_s}_1)}{||x^{k_s}-x||} + \frac{o(x,x^{k_s},g^{k_s}_2)}{||x^{k_s}-x||} \to 0,
\]
when $s\to \infty$. The lemma is proved.
\begin{corollary}
\label{cor:1.5}
{For a generalized differentiable function
$f(x) (x\in E_1)$ of one variable, the set of points of ambiguity
(discontinuity) of its pseudogradient mapping $G_f(x)$ is no more than countable.}
\end{corollary}

Indeed, the set $A$ of discontinuity points of a pseudogradient can be
represented as a countable union of finite point sets
\[ 
A = \bigcup^{\infty}_{i,j=1} \left \{ x\in E_i \middle| |x|\leq i, d(x)\geq \frac{1}{j} \right \}  
\]
\begin{definition}
\label{df:1.5}
{The diameter of the set $G \subset E_n $ in the direction $l\in E_n$ is the value}
\[
d_G(l) = \sup_{g\in G}\left<g,l\right> - \inf_{g\in G}\left<g,l\right>.
\]
\end{definition}

It is easy to see that the function $d_G(l)$ is convex in $l$.
\begin{definition}
\label{df:1.6}
The diameter of the set $G \subset E_n $ is the value
\[
 d_G = \sup \left \{\|g_1-g_2\| \text{ }|\text{ } g_1,g_2\in G \right \}.
\]
\end{definition}

It cab be seen that
\[
d_G = \sup\left\{d_g(l)\text{ }\middle|\text{ } \|l\|\leq 1\right\}.
\]
\begin{lemma}
\label{lem:1.2}
{For the set $G\subset E_n$ to consist of a single point, it is necessary and sufficient that for some system of $n$ linearly independent directions $\{l_i\}$ it holds $d_G(l_i) =0 (i=1,...,n).$}
\end{lemma}

{\it P r o o f}. The necessity is obvious. Let us prove the sufficiency.
It is clear that
\[
d_G(-l_i)=0, i=1,...,n.
\]
For any direction ${l\in E_n}$, the representation holds
\[
l= \sum^{n}_{i=1}\alpha_i l_i = \sum^{n}_{i=1} |\alpha_i|l_i \text{ sign }  \alpha_i,
\]
where
\[
 \text{sign }  \alpha = \begin{cases}
    1,       & \quad \alpha > 0, \\
    0,  & \quad \alpha = 0, \\
    -1, & \quad \alpha < 0. \\
  \end{cases}
\]
Due to the convexity of the function $d_G(l)$, we have
\[
d_G(l)\leq \sum^{n}_{i=1}|\alpha_i| d_G(l_i\text{ sign } \alpha_i)=0.
\]
Hence,
\[
d_G=sup\left\{d_G(l)\,\middle|\, \|l\|\leq 1 \right\} =0,
\]
i.e., $G$ consists of a single point. The lemma is proved.
\begin{theorem}
\label{th:1.12}
{The set of pseudogradients $G_f(x)$ of the generalized differentiable function $f(x)\,(x\in E_n)$ almost everywhere in $E_n$ consists of one point, and hence the generalized differentiable function is almost everywhere differentiable (by Fr\'echet and by Gateaux), and its gradient
is continuous on the set where it is defined.}
\footnote{Mapping $g(x)$ is called continuous on a set $X$, if for any 
$x\in X$ and $\{x^k\}\rightarrow x$ $(x^k\in X)$ it holds $g(x^k)\rightarrow g(x)$.}
\end{theorem}

{\it P r o o f}. Let the vectors $l_i\:(i=1,...,n)$ be the unit vectors of the coordinate axes in $E_n$. Denote by $d(x)$ the diameter of the set $G_f(x)$ and by $d(x,l)$ its diameter in the direction $l\in E_n$. By virtue of Lemma 1.2, the set
\[
A = \{x\in E_n \,|\, d(x)>0\}
\]
of multivaluedness of the pseudogradient mapping $G_f (x)$ 
can be represented in the form
\[
A=\bigcup^{n}_{i=1}A_i,\; A_i=\{x\in E_n \,|\,d(x,l_i)>0\}.
\]

Define the sets
\[
A_{ij}=\{x\in E_n\,|\,d(x,l_i)\geq1/j\},\; j=1,2,... .
\]

Obviously, $A=\bigcup_{i,j}A_{ij}$. It is easy to show that the function $d(\cdot,l_i)$, is semicontinuous from above; therefore the sets $A_{ij}$ are closed and hence Lebesgue measurable in $E_n$. Let us introduce the sets
\[
K_m=\{x\in E_n\,|\max_{1\leq s \leq n}|x_{(s)}|\leq m\},\quad A_{ijm}=A_{ij}\cap K_m
\]
They are also closed and measurable. We have $A=\cup_{i,j,m}A_{ijm}$.

    Let us show that the Lebesgue measure of the sets $A_{ijm}$ is zero. Let us define
a hyperplane
\[
H_i=\{x\in E_n \, |\,\left<x,l_i\right>=0\}.
\]
For points
\[
x=(x_{(1)},...\,,\,x_{(i-1)}, 0, x_{(x+1)},...\,,\,x_{(n)})\in H_i,
\]
denote $\Pi_i(x)$ the lines passing through $x$,
\[
\Pi_i(x)=\{x+\lambda l_i\,|\,-\infty<\lambda<\infty\}.
\]
The intersections of $H_i$ and $\Pi_i(x)$ with the cube $K_m$ are denoted by $H_{im}$ and $\Pi_{im}(x)$, respectively. Consider the sets
\[
A_{ijm}(x)=A_{ijm}\cap\Pi_{im}(x).
\]
It is easy to see that they are sets of discontinuity points of the pseudogradient for a generalized differentiable function
\[
\phi(\lambda)=f(x+\lambda l_i),\quad |\lambda|\leq m.
\]
According to Lemma 1.1, the sets $A_{ijm}(x)$ are finite or empty; hence their one-dimensional
Lebesgue measure is zero,
\[
\mu_i A_{ijm}(x)=0.
\]
By Fubini's theorem, the $n$-dimensional Lebesgue measure is
\[
\mu_n A_{ijm}=\int_{H_im} \mu_1 A_{ijm}(x)\mu_{n-1} \mathrm{d}x=0.
\]
The result is
\[
\mu_n A=\sum_{i,j,m}\mu_nA_{ijm}=0.
\]
Thus, almost everywhere in $E_n$ the mapping $G_f(x)$ is single-valued and continuous (i.e. semicontinuous from above and from below). By Theorem 1.11, it follows that generalized differentiable functions are differentiable almost everywhere in $E_n$ (by Fr\'echet and by Gateaux), and their gradients are continuous wherever they are defined. The theorem is proved.

Now let us prove for generalized differentiable functions the generalized Newton-Leibniz formula. From this it will follow that a generalized differentiable function is defined by any of its pseudogradient mappings up to a constant.
\begin{theorem}
\label{th:1.13}
{Let $f(x)\;(x\in E_n)$ be a generalized differentiable
function, $g(x)$ be an arbitrary single-valued Borel-measurable section 
of the pseudo-gradient mapping $G_f(x)$, i.e.,}
\[
g(x)\in G_f(x).
\]
\textit{Then the integral formula takes place:}
\[
f(x+h)-f(x)=\int_0^1\left<g(x+th),h\right>\mathrm{d}t,
\]
\textit{where the integral is understood in the Lebegues sense, $(\cdot,\cdot)$ denotes the scalar product.}
\end{theorem}

{\it P r o o f.} Consider the function
\[
\phi(x,h)=\int_0^1f(x+th)\mathrm{d}t
\]
and calculate its derivative $\phi_h^\prime(x,h)$ in the direction of $h$ in two ways. On the one hand,
\[\begin{array}{l}
\phi_h^\prime(x,h)=\lim_{\alpha \to +0}\frac{1}{\alpha} [\phi(x+\alpha h,h)-\phi(x,h)]  \\
=\lim_{\alpha \to +0}\frac{1}{\alpha}\left[\int_{\alpha}^{1+\alpha}f(x+t^{\prime}h)\mathrm{d}t^\prime-\int_{0}^{1}f(x+th)\mathrm{d}t\right]=f(x+h)-f(x).
\end{array}\]
On the other hand,
\begin{equation}
\phi_h^\prime(x,h)=\lim_{\alpha \to +0}\int_0^1\frac{1}{\alpha}[f(x+th+\alpha h) -f(x+th)]\mathrm{d}t.\quad \label{eqn:1.20}
\end{equation}
Here it is obvious that the integrand function is modulo bounded 
by the Lipschitz constant at $t\in [0,1]$ uniformly over $\alpha \in [0,1]$. Due to the generalized differentiability of $f(x)$, the representation holds
\[
\frac{1}{\alpha}[f(x+th+\alpha h)-f(x+th)]=
\]
\[
= \left<g(x+th+\alpha h),h\right>+\frac{1}{\alpha}o\left(x+th,x+th+\alpha h, g(x+th+\alpha h)\right),
\]
where
\[
\lim_{\alpha \to +0}\frac{1}{\alpha} o\left(x+th,x+th+\alpha h,g\right)=0.
\]
Function
\[
\psi (t)=f(x+th),\quad t\in[0,1],
\]
according to Theorem 1.6, is generalized differentiable; its pseudogradient mapping is given by the formula
\[
G_{\psi}(t)=\{(g,h)\, |\, g\in G_f(x+th)\}
\]
According to Corollary 1.5, the mapping $G_{\psi}(t)$ almost everywhere on $[0, 1]$
is single-valued and hence continuous. It follows that almost
for all $t\in[0,1]$
\[
\lim_{\alpha \to +0}\left<g(x+th+\alpha h),h\right>=\left<g(x+th),h\right>.
\]
Now, according to Lebesgue's theorem, the limit in (1.20) can be brought under the sign of the
integral. Thus, we obtain the relation
\[
f(x+h)-f(x)=\phi_h^{\prime}=\int_0^1\left<g(x+th),h\right>\mathrm{d}t.
\]
The theorem is proved.
\begin{lemma}
\label{lem:1.3}
{Let $f_1(x), f_2(x)$ be generalized differentiable functions and $G_{f_1}(x), G_{f_2}(x)$ be their pseudogradient mappings, and}
\[
G_{f_1}(x)\cap G_{f_2}(x) \neq \emptyset
\]
\textit{for all $x\in X$($X$ is a convex set). Then}
\[
f_1(x)-f_2(x)=const.
\]
\end{lemma}

{\it P r o o f}. According to Theorem 1.6, the function
\[
f(x)=f_1(x)-f_2(x)
\]
is generalized differentiable, and its pseudogradient
is given by the formula
\[
G_f(x)=\{g\in E_n\, |\, g=g_1-g_2,\,g_1\in G_{f_1}(x),\, g_2\in G_{f_2}(x)\}.
\]
Since
\[
G_{f_1}(x)\cap G_{f_2}(x)\neq \emptyset,
\]
then it follows that $0\in G_f(x)$ for $x\in X$. Let us show that in this case $f(x) = const$ under $x\in X$ . Indeed, by Theorem 1.13, for any
$x,\, y\in X$ we have
\[
f(y)=f(x)+\int_0^1\left<g(x+t(y-x)),\,y-x\right>\mathrm{d}t,
\]
where
\[
g(x+t(y-x))\in G_f(x+t(y-x)).
\]
Assuming
\[
g(x+t(y-x))\equiv 0,
\]
we get $f(y) = f(x)$. The lemma is proved.

\begin{corollary}
\label{cor:1.6}
{A generalized differentiable function is defined by any of its pseudogradient mapping up to a constant.}
\end{corollary}
\begin{corollary}
\label{cor:1.7}
{A semi-smooth function is defined by its subdifferential up to a constant, since it is generalized differentiable and its subdifferential can be taken as a pseudodifferential mapping.}
\end{corollary}

Let us conclude the paragraph by clarifying the relation between generalized differentiable and semi-smooth functions (see Definition 1.4).
\begin{definition}
\label{df:1.7}
{ A function is called generalized differentiable 
in the narrow sense in a domain, if it is generalized differentiable in a domain and has a derivative in any direction at each point of the domain.}
\end{definition}
\begin{theorem}
\label{th:1.14}
{Classes of the generalized differentiable in the narrow sense and the semi-smooth functions coincide.}
\end{theorem}

{\it P r o o f}. Let $f(x)$ be a semi-smooth function. According to Theorem 1.4, it is generalized differentiable.
Let us show that $f(x)$ has derivatives along the directions,
\[
f^{\prime}(x;d)=\lim_{\lambda_k \to +0}[f(x+\lambda_k d)-f(x)]/\lambda_k.
\]
By virtue of Theorem 2.9, it is possible to represent
\[
[f(x+\lambda_k d)-f(x)]/\lambda_k=\left<\bar{g}^{k},d\right>,\quad \bar{g}^{k}\in \partial f(x+\theta_k\lambda_kd),
\]
\[
0<\theta_k<1.
\]
Let's set
\[
y^k=x+\theta_k\lambda_kd,\;\;\;(y^k-x)/\|y^k-x\|=d/\|d\|.
\]
By the definition of the semi-smooth functions, there exists $\lim_{k\to\infty}(\bar{g}^{k},d)$, and hence there exists a directional derivative $f^{\prime}(x;d)$. So the semi-smooth function is generalized differentiable in the narrow sense.

Now let $f(x)$ be generalized differentiable and has directional derivatives at each point. Let us show that it is semi-smooth. According to Theorem 1.1, $f(x)$ is locally Lipschitz and, according to Theorem 1.10, $\partial f(x)$ is its pseudogradient mapping. Let
\[
y^k\to x,\;\;\; g^k\in\partial f(y^k),
\]
where
\[
(y^k-x)/\|y^k-x\|\to d.
\]
A valid expansion is
\[
f(y^k)=f(x)+\left<g^k,y^k-x\right>+o(x,y^k,g^k),
\]
whence
\[
\lim_{k\to\infty}\left<g^k, \frac{y^k-x}{||y^k-x||}\right>=\lim_{k\to\infty}(g^k,d)=\lim_{k\to\infty}\frac{f(y^k)-f(x)}{||y^k-x||}.
\]
The last limit exists and is equal to the derivative of the function $f$ in the direction $d$. Indeed, let us set
\[
\lambda_k=(y^k-x,d)
\]
and represent
\[
y^k=x+\lambda_kd+r_k.
\]
It is not difficult to check that
\[
\lambda_k\to 0,\;\;\; r^k/\lambda_k\to 0, \;\;\; \|y^k-x\|/\lambda_k\to 1.
\]
Then
\[\lim_{k\to\infty}\frac{f(y^k)-f(x)}{\|y^k-x\|}=\]
\[
=\lim_{k\to\infty} \left[
\frac{f(x+\lambda_kd)-f(x)}{\lambda_k}\frac{\lambda_k}{\|y^k-x\|}+
\frac{f(x+\lambda_kd+r_k)-f(x+\lambda_kd)}{\lambda_k}\frac{\lambda_k}{\|y^k-x\|}\right]=
\]
\[
=\lim_{k\to\infty}\frac{f(x+\lambda_kd)-f(x)}{\lambda_k}.
\]
We have shown that the limit of scalar products $\lim_{k\to\infty}(g^k,d)$ always exists; hence, $f(x)$ is a semi-smooth function. The theorem is proved. 

So, as classes, semi-smooth and generalized differentiable functions (in the narrow sense) coincide. The difference between them lies in the way the generalized gradients are defined. For the semi-smooth functions we use ready generalized Clarke gradients, or, viewed from Theorem 1.10, minimal in a certain sense generalized gradients. However, as can be seen from the results of this paragraph and will be confirmed later, it is useful to consider not only minimal, but also broader generalized gradient mappings. For instance, this allows one to construct a pseudo-gradient calculus, which does not exist for minimal generalized Clarke gradients. Furthermore, the transition from semi-smooth to generalized differentiable functions allows one to relax the differentiability requirements on the functions in question without any detriment. Finally, an important point is the transition to decomposition (1.1), from which, in fact, all the results follow.

\section*{$\S$ 2. Smoothing of functions and basic properties of locally Lipschitz functions}
\label{Sec.2}
\setcounter{section}{2}
\setcounter{definition}{0}
\setcounter{equation}{0}
\setcounter{theorem}{0}
\setcounter{lemma}{0}
\setcounter{remark}{0}
\setcounter{corollary}{0}
\setcounter{example}{0}

A successful approach to studying locally Lipschitz functions and developing numerical methods for solving them involves considering limiting extremal problems of the following type.
There is a sequence of functions $f(x,\alpha)$ that converges to $f(x)$ as $\alpha\rightarrow 0.$ The goal is to find the minimum of $f(x)$ in the domain $X$ using only the functions in the sequence $f(x,\alpha)$.
The operation of passing to the limit can spoil some of the good properties of the functions 
$f(x,\alpha),$ which is typical for approximation problems, when a "bad" function for some reasons is approximated by a sequence of "good" ones.
In such cases, an interesting problem arises of optimizing the limiting function $f(x)$ using information about the members of the sequence $f(x,\alpha)$.

In this paragraph, the study of a number of properties of locally Lipschitz functions is carried out on the basis of investigating the properties of smoothed functions $f(x,\alpha)$, converging to $f(x)$ as  $\alpha \to 0$. Artificial averaging of the function $f(x)$ allows in many cases to smooth out and reduce its oscillatory character. In practice, linear averaging operators are most often used,
\[
f(x,\alpha) = \int_{E_n}h_{\alpha}(u)f(x+u)\mathrm{d}u,
\]
where $h_{\alpha}$ is the weight function (the operator kernel), and $u\in E_n,$, $\alpha$ is a non-negative scalar parameter.

In the theory of generalized functions,  infinitely differentiable functions are
considered,
\[
f(x,\alpha)=\int_{-\alpha}^{\alpha} \omega_{\alpha}(u)f(x+u)du,\;\;\;
x, u\in E_1,
\]
where
\[
\omega_{\alpha}(u)=\begin{cases} C_{\alpha}\exp\left\{ - \frac{\alpha^2}{\alpha^2-|u|^2}\right\},&\quad |u|< \alpha,\\
0,&\quad|u\ge\alpha;
\end{cases}
\]
$C_{\alpha}$ is chosen from the condition
\[
\int_{-\alpha}^{\alpha}\omega_{\alpha}(u)\mathrm{d}u=1
\]

In the present book we study the smoothed functions of the following form:
\begin{equation}
f(x,\alpha)=\frac{1}{(2\alpha)^n}\int_{-\alpha}^{\alpha}...\int_{-\alpha}^{\alpha}f(x+y)\mathrm{d}y_1,...\mathrm{d}y_n.\quad \tag{2.1}
\end{equation}

The class of locally-lipschitz functions is extremely broad; it covers almost all classes of continuous functions studied in mathematical programming. At the same time, locally Lipschitz functions preserve some useful properties that allow us to introduce the notion of a "generalized gradient" and to construct numerical procedures for determining stationary points.
\begin{lemma}
\label{lem:2.1}
{Partial derivatives of a locally Lipschitz function}
\[
\partial f(x)/\partial x_i,\quad i=1\,,\,...\,,\,n.
\]
\textit{exist everywhere except for the set of measure zero.}
\end{lemma}

{\it P r o o f}. 
Since the function satisfies the Lipschitz condition, it follows that it is absolutely continuous in each variable $x_i (i=1,\,...\,,n)$ with the other variables fixed. Let us recall the definition of an absolutely continuous function.

Let a finite function $g(x)$ be given on the segment $[a,b]$. If for any $\varepsilon >0$ there exists such an $\delta >0$ that for any finite system of non-intersecting intervals $(a_1,b_1),\,...,\,,\,(a_m,\,b_m)$, where $\sum^{m}_{i=1}(b_i-a_i)<\delta$, the inequality 
$$
\sum_{i=1}^{m}\left|g\left(b_{i}\right)-g\left(a_{i}\right)\right|<\varepsilon
$$
holds, then the function $g(x)$ is called absolutely continuous on the segment $[a,b]$.

Absolutely continuous functions have the following property. 
Almost everywhere in the interval $[a, b]$ (i.e., except for a set of measure zero), there exists a finite derivative. 
Thus, for locally Lipschitz functions, there exists a partial derivative $\partial f(x) / \partial x_{i}$ almost everywhere for the variable $x_{i} (i=1, \ldots, n)$.
Hence, by the Fubini-Tonelli theorem on product measures [58], $f(x)$ has an ordinary partial derivative $\partial f(x) / \partial x_{i}(i=1, \ldots , n)$ almost everywhere in $E_n$.

It turns out that the result of Lemma 2.1 can be strengthened.
\begin{theorem}
\label{th:2.1}
{A locally Lipschitz function is differentiable almost everywhere, i.e. $\nabla f(x)$ exist everywhere except for a set of measure zero.}
\end{theorem}

This fact is known as Rademacher's theorem [173]; another proof of the theorem is contained in [55]. Based on this property, the concept of "generalized gradient $\partial f(x)$ of a locally Lipschitz function $f(x)$ at a point $x$" was originally introduced as a convex closure of a set of points $y$ of the form
$$
y=\lim _{k \rightarrow \infty} \nabla f\left(x^{k}\right),
$$
where $\left\{x^{k}\right\}$ is a sequence converging to $x$ and such that $f(x)$ is differentiable at every point of $x^{k}$.

In what follows, we will use another, in many cases more convenient, equivalent definition of the set $\partial f(x)$.

A locally Lipschitz function may not have the usual directional derivative
$$
f^{\prime}(x ; e)=\lim _{\lambda \downarrow 0}[f(x+\lambda e)-f(x)] / \lambda.
$$

\textit{The generalized derivative} of the function $f(x)$ at the point $x$ in the direction $e$, denoted by $f^{0}(x ; e)$, is defined by the formula
\[f^{0}(x ; e)=\lim _{\delta \downarrow 0} \sup _{\substack{\|v \mid\| \delta \le \delta \\ 0<\lambda \le \delta}}[f(x+v+\lambda e)-f(x+v)] / \lambda .\tag{2.2}\]
\begin{lemma}
\label{lem:2.2}
{The function $e \rightarrow f^{0}(x ; e)$ is convex.}
\end{lemma}

{\it P r o o f.} By Theorem 4.7 [109], a function $f(x)$ is convex if it is positively homogeneous and semi-additive, i.e.,
$$
\begin{gathered}
f(\lambda x)=\lambda f(x), \quad 0<\lambda<\infty \\
f(x+y) \le f(x)+f(y) \quad \forall x, y
\end{gathered}
$$

Positive homogeneity is obvious, and semi-additivity follows from the chain of relations
$$
\begin{aligned}
& f^0(x; e_1 + e_2) = \lim _{\delta \downarrow 0} \sup _{\substack{\|v\| \le \delta \\ 0<\lambda \le \delta}}[f\left(x+v+\lambda e_{1}+\lambda e_{2}\right)-f(x+v)]/\lambda \le \\
& \le \lim _{\delta \downarrow 0} \sup _{\substack{\|v\| \le \delta \\
0<\lambda \le \delta}}\left[f\left(x+v+\lambda e_{1}+\lambda e_{2}\right)-f\left(x+v+\lambda e_{2} \right)\right] / \lambda+ \\
& +\lim _{\delta \downarrow 0} \sup _{\substack{\|v\| \le \delta \\
0<\lambda \le \delta}}\left[f\left(x+v+\lambda e_{2}\right)-f(x+v)\right] / \lambda=f^{0}\left (x ; e_{1}\right)+f^{0}\left(x ; e_{2}\right) .
\end{aligned}
$$
\begin{definition}
\label{df:2.1}
{The generalized gradient} of the function $f(x)$ at the point $x$, denoted by $\partial f(x)$, is the subdifferential of the convex function $f^{0}(x ; \cdot)$ at zero,  i.e. $y \in \partial f(x)$ if
\[f^{0}(x ; e) \geq\left<y, e\right> \quad \forall e .\tag{2.3}\]
\end{definition}

Consequently, the set $\partial f(x)$ is convex and closed, and due to the inequality $\left|f^{0}(x ; e)\right| \le L\|e\|$ is bounded.

It follows from the inequality $(2.3)$ that if $f(x)$ is differentiable at a point $x$, then
$$
\nabla f(x) \in \partial f(x).
$$
\begin{theorem}
\label{th:2.2}
{The following formula is valid}
\[f^{0}(x, e)=\max_{y \in \partial f(x)}\left<y, e\right> .\tag{2.4}\]
\end{theorem}

{\it P r o o f.} According to Theorem 13.1 in [109], each convex set $M$ is defined by its support function $\sup _{y \in M}(y, e)$, where $e-$ is an arbitrary vector. This means that $h \in \bar{M}$ if and only if
\[\left<h, e\right> \le \sup _{y \in M}\left<y, e\right> \quad \forall e .\tag{2.5}\]
In addition, if on the sets $M_{1}$ and $M_{2}$
\[\sup _{y \in M_{1}}\left<y, e\right>=\sup _{y \in M_{2}}\left<y, e\right>,\tag{2.6}\]
then $\bar{M}_{1}=\bar{M}_{2}$. Therefore formula (2.4) follows from relations (2.3) and (2.5).
\begin{theorem}
\label{th:2.3}
{Let $A$ be a convex closed set. Then}
\[\partial f(x) \subset A,\]
if
\[f^{0}(x ; e) \le \max _{u \in A}\left<u, e\right> \quad \forall e \tag{2.7}.\]
\end{theorem}

{\it P r o o f.} If there exists some vector $y \notin A$, $y \in \partial f(x)$, then by the separation theorem there exists a vector $e$ such that
$$
\max _{u \in A}\left<e, u\right><\left<e, y\right>.
$$

But this inequality contradicts (2.7). Another equivalent definition of the set $\partial f(x)$ is as follows.
\begin{theorem}
\label{th:2.4}
{The set $\partial f(x)$ is the convex closure of the set of all points of the form}
$$
G=\{g=\lim _{k \rightarrow \infty} \nabla f\left(x^{k}\right)\},
$$
{where $\left\{x^{k}\right\}$ is a sequence converging to $x$ and such that $f(x)$ is differentiable at every point $x^{k}$.}
\end{theorem}

{\it P r o o f.} Denote $M=\operatorname{co}G$. It is easy to see that the set $M$ is convex and bounded. Its closedness follows from the closedness of the set $G$. Let $\left\{y^{k}\right\} \in G$ and $y^{k} \rightarrow y$. We need to show that $y \in G$. For each number $k$ there is a point $x^{k}$ such that
$$
\left\|y^{k}-\nabla f\left(x^{k}\right)\right\| \le 1 / k, \quad\left\|x^{k}-x\right\| \le 1 / k
$$
moreover, at the points $x=x^{k}$ the function $f(x)$ is differentiable. That's why
$$
\left\|y-\nabla f\left(x^{k}\right)\right\| \le\left\|y-y^{k}\right\|+\left\|y^{k}-\nabla f\left(x^{k}\right)\right\| \le\left\|y-y^{k}\right\|+1 / k,
$$
i.e.
$$
\left\|y-\nabla f\left(x^{k}\right)\right\| \rightarrow 0 \quad \text { as } \quad k \rightarrow \infty,
$$
whence it follows that $y \in G$.

Initially, in [145], the set of generalized gradients of a locally Lipschitz function was defined as $\partial f(x)=\operatorname{co}G$. With this definition of $\partial f(x)$, it was shown that
\[f^{0}(x, e)=\max _{y \in \operatorname{co}G}\left<y, e\right>. \tag{2.8}\]
Hence, comparing relations (2.4), (2.6), and (2.8), we obtain the assertion of the theorem.

An important property is that the point-to-set mapping $x \rightarrow \partial f(x)$ is closed.
\begin{theorem}
\label{th:2.5}
{Let}
$$
g\left(x^{k}\right) \in \partial f\left(x^{k}\right), \quad x^{k} \rightarrow x, \quad g\left(x^{k }\right)\rightarrow g
$$
{Then} $g \in \partial f(x)$.
\end{theorem}

{\it P r o o f.} By definition, we have
$$
f^{0}\left(x^{k}, e\right) \geqslant\left<g\left(x^{k}\right), e\right>.
$$
Therefore, there are $v^{k} \in E_{n}$ and $\lambda_{k}>0$ such that
$$
\frac{1}{\lambda_{k}}\left[f\left(x^{k}+v^{k}+\lambda_{k} e\right)-f\left(x^{k} +v^{k}\right)\right] \geqslant\left<g\left(x^{k}\right), e\right>-\frac{1}{l} ,
$$
where
$$
\left\|{v}^{k}\right\|+\lambda_{k}<1 / k
$$
Since $\left<g\left(x^{k}\right), e\right>$ converges to $\left<g, e\right>$, the last inequality implies the inequality
$$
f^{0}(x ; e) \geq\left<(g, e\right>.
$$
Therefore, by the definition of the generalized gradient, $g \in f(x)$, which was to be proved.
\begin{theorem}
\label{th:2.6}
{If $f(x)$ is a convex function, then the set $\partial f(x)$ coincides with the usual subdifferential for $f(x)$}.
\end{theorem}

{\it P r o o f.} It is known that a convex function $f(x)$ has a directional derivative $f^{\prime}(x ; e)$ for which
$$
\left.f^{\prime}(x; e\right)=\max _{y \in \hat{\partial} f(x)}\left<y, e\right>,
$$
where $\hat{\partial} f(x)$ is its subdifferential. Conditions (2.4), (2.6) imply that for the proof it suffices to show the equality
$$
f^{\prime}(x; e)=f^{0}(x ; e).
$$
The inequality
$$
f^{\prime}(x; e) \le f^{0}(x; e)
$$
is obvious. For a convex function, the mean value theorem holds,
\begin{eqnarray}
\frac{f(x+v+\lambda e) - f(x+v)}{\lambda}& =& \left<g(x+v+t \lambda e), e\right>,\;\; 0 \le t \le 1, \nonumber\\
g(x+v + t\lambda e) &\in& \hat{\partial} f(x+v+t \lambda e),\nonumber
\end{eqnarray}
whence, due to the closedness of the mapping $x \rightarrow \hat{\partial}f(x)$, 
the opposite inequality holds
$$
f^{0}(x ; e) \le \max _{y \in \hat{\partial}f(x)}\left<y, e\right>=f^{\prime}(x ; e).
$$

Let us now consider a question connected with the application of the usual operations of addition and multiplication of locally Lipschitz functions.
\begin{theorem}
\label{th:2.7}
{If the functions $f_{1}(x)$ and $f_{2}(x)$ satisfy the local Lipschitz condition, then}
$$
f_{3}(x)=f_{1}(x) \pm f_{2}(x),\;\;\; f_{4}(x)=f_{1}(x) f_{2}(x)
$$
{ also belongs to this class.}
\end{theorem}

{\it P r o o f.} Let $L_{1}, L_{2}$ be the Lipschitz constants of the functions $f_{1}, f_{2}$ with respect to an arbitrary bounded set $S$. We must show that the functions $f_{3}, f_{4}$ satisfy the Lipschitz condition in $S$.

For $f_{3}$ the expression is justified
$$
\begin{aligned}
f_{3}(x)-f_{3}(y) \mid= & \left|f_{1}(x) \pm f_{2}(x)-f_{1}(y) \mp f_{ 2}(y)\right| \le \\
& \le\left|f_{1}(x)-f_{1}(y)\right|+\left| \pm f_{2}(x) \mp f_{2}(y)\right| \le \\
& \le L_{1}\|x-y\|+L_{2}\|x-y\|=L\|x-y\|,
\end{aligned}
$$
where $L=L_{1}+L_{2}$.

Let $a_{1}$ and $a_{2}$ be the upper bounds of the functions $\left|f_{1}(x)\right|,\left|f_{2}(x)\right|$ on the set $ \bar{S}$, where $\bar{S}$ is the closure of $S ; L=\max \left\{L_{1}, L_{2}\right\}$. Then
$$
\begin{aligned}
& \left|f_{4}(x)-f_{4}(y)\right|=\left|f_{1}(x) f_{2}(x)-f_{1}(y) f_{ 2}(y)\right|= \\
& =\left|f_{1}(x) f_{2}(x)-f_{1}(x) f_{2}(y)+f_{1}(x) f_{2}(y)- f_{1}(y) f_{2}(y)\right| \le \\
& \le\left|f_{1}(x)\right| \times\left|f_{2}(x)-f_{2}(y)\right|+\left|f_{2}(y)\right| \times\left|f_{1}(x)-f_{1}(y)\right| \le \\
& \le a_{1} L_{2}\|x-y\|+a_{2} L_{1}\|x-y\| \le\left(a_{1}+a_{2}\right) L\|x-y\|,
\end{aligned}
$$
where $x, y \in S$.
\begin{theorem}
\label{th:2.8}
{The following inclusion formula is valid}
\[\partial\left(f_{1}(x)+f_{2}(x)\right) \subset \partial f_{1}(x)+\partial f_{2}(x).\tag{ 2.9}\]
\end{theorem}

{\it P r o o f.} We use Theorem 2.3. Let
$$
A=\partial f_{1}(x)+\partial f_{2}(x)
$$
Then
$$
\max _{u \in \partial_{1}(x)+\partial f_{2}(x)}\left<u, e\right>=f_{1}^{0}(x ; e)+f_{2} ^{0}(x ; e).
$$

It remains to note that
$$
\left(f_{1}+f_{2}\right)^{0}(x ; e) \le f_{1}^{0}(x ; e)+f_{2}^{0}(x; e) .
$$

It is easy to see that relation (2.9) can have a strict inclusion.

A very successful technique, which makes it possible to derive a number of other interesting properties of locally Lipschitz functions, consists in approximating the function $f(x)$ by a sequence of smoothed functions $f(x, \alpha)$ converging to $f(x)$ as $ \alpha \rightarrow 0$. Let us study the function $f(x, \alpha)$ of the form (2.1):
$$
f(x, \alpha)=\frac{1}{(2 \alpha)^{n}} \int_{-\alpha}^{\alpha} \ldots \int_{-\alpha}^{\alpha} f(x+y) d y_{1} \ldots d y_{n}
$$
where $\alpha>0$.

Below we will often use the following fact.

\noindent
{\bf Lemma {$2.2^*$}}
{\it The sequence of functions $f(x, \alpha)$ converges uniformly to $f(x)$ as $\alpha \rightarrow 0$ in any bounded region of space $E_{n}$.}

The proof of the lemma naturally follows from the relations
$$
f(x, \alpha)=\frac{1}{(2 \alpha)^{n}} \int_{-\alpha}^{\alpha} \cdots \int_{-\alpha}^{\alpha} [f(x+y)-f(x)] d y+f(x) \le f(x)+\sqrt{n} L \alpha,
$$
where $L$ is the Lipschitz constant of the function $f(x)$.

Let us formulate important differential properties of the function $f(x, \alpha)$.
\begin{lemma}
\label{lem:2.3}
{The gradient of the function $f(x, \alpha)$ at the point $x$ is determined by the formula}
$$
\nabla f(x, \alpha)=\frac{1}{(2 \alpha)^{n}} \int_{-\alpha}^{\alpha} \cdots \int_{-\alpha}^{\alpha} \nabla f(x+y) d y .
$$
\end{lemma}

The proof follows from Lebesgue's theorem on passing to the limit under the integral sign.
\begin{lemma}
\label{lem:2.4}
{The gradient of the function $f(x, \alpha)$ in any bounded region of space $E_{n}$ satisfies the Lipschitz condition with constant $L=C / \alpha$, where $C$ is some constant.}
\end{lemma}

{\it P r o o f.} The partial derivative $\partial f(x, \alpha) / \partial x_{1}$ is given by
\[
\begin{aligned}
& \frac{\partial f(x, \alpha)}{\partial x_{1}}=\frac{\partial}{\partial x_{1}} \frac{1}{(2 \alpha)^{n}} \int_{-\alpha}^{\alpha} \cdots \int_{-\alpha}^{\alpha} f(x+y) d y= \\
& =\frac{\partial}{\partial x_{1}} \frac{1}{(2 \alpha)^{n}} \int_{x_{1}-\alpha}^{x_{1}+\alpha} \cdots \int_{x_{n}-\alpha}^{x_{n}+\alpha} f(y) d y= \\
& =\frac{1}{(2 \alpha)^{n}} \int_{x_{2}-\alpha}^{x_{2}+\alpha} \cdots \int_{x_{n}-\alpha}^{x_{n}+\alpha}\left[f\left(x_{1}+\alpha, y_{2}, \ldots, y_{n}\right)-f\left(x_{1}-\alpha, y_{2}, \ldots, y_{n}\right)\right] d y_{2} \ldots d y_{n}= \\
& =\frac{1}{(2 \alpha)^{n}} \int_{-\alpha}^{\alpha} \ldots \int_{-\alpha}^{\alpha}\left[f\left(x_{1}+\alpha, x_{2}+z_{2}, \ldots, x_{n}+z_{n}\right)-\right. \\
& \left.-f\left(x_{1}-\alpha, x_{2}+z_{2}, \ldots, x_{n}+z_{n}\right)\right] d z_{2} \ldots d z_{n}, \quad 
\end{aligned} \tag{2.11}
\]
which follows from the formula for differentiating the integral with respect to the upper limit. Other partial derivatives are calculated similarly.

According to formula (2.11), the partial derivative $\partial f(x, \alpha) / \partial x_{1}$ satisfies the Lipschitz condition with constant $L=C / \alpha$, i.e.
$$
\left|\frac{\partial f(x, \alpha)}{\partial x_{1}}-\frac{\partial f\left(y, \alpha\right)}{\partial y_{1}} \right| \le \frac{C}{\alpha}\|x-y\| \text {. }
$$
This, in turn, implies that
$$
\|\nabla f(x, \alpha)-\nabla f(y, \alpha)\| \le \frac{C \sqrt{n}}{\alpha}\|x-y\| .
$$

Formulas (2.10), (2.11) show that there are two ways to calculate the gradient of a smoothed function: for one, $\nabla f(x, \alpha)$ is determined by the gradients of the function $f(x)$, for the other, through finite differences of $f (x)$. These formulas will prove useful in constructing numerical methods for minimizing Lipschitz functions. In accordance with them, two types of methods will be substantiated: finite-difference methods and methods of generalized gradient descent. In what follows, the following characteristic of the gradient of the function $f(x, \alpha)$ is often used.
\begin{lemma}
\label{lem:2.5}
{Let $\alpha_{k} \rightarrow 0$. For any $\varepsilon>0$ there are $\delta(x)>0$ and number $\bar{k}$ such that for all points $y$ such that $\|y-x\| \le \delta$, for $k \geqslant \bar{k}$ the inequality holds}
$$
\rho\left(\nabla f\left(y, \alpha_{k}\right), \partial f(x)\right) \equiv \min _{g \in \partial f(x)}\left\|\nabla f\left(y, \alpha_{k}\right)-g\right\| \le \varepsilon .
$$
\end{lemma}

{\it P r o o f.} Let us first show the following: for any $\varepsilon>0$ there exists $\bar{\delta}>0$ such that for all points $y$ at which the function $f(x)$ is differentiable and such that $\| y-x\| \le \bar{\delta}$, we have the inequality
$$
\rho(\nabla f(y), \partial f(x)) \le \varepsilon
$$

Assume the contrary: there exists $\varepsilon$ such that for any $\bar{\delta}$ there are points $\tilde{y},\|\tilde{y}-x\| \le \bar{\delta}$ for which $\rho(\nabla f(\tilde{y}),  \partial f(x))>\varepsilon$. Then the limit points of the sequence $\nabla f(\tilde{y})$ for $\bar{\delta} \rightarrow 0$ do not belong to $\partial f(x)$. Thus, we obtain a contradiction with the definition of the set $\partial f(x)$. We choose the value $\delta=$ $=\bar{\delta} / 2$ required in the condition of the lemma, and the number $\bar{k}$ such that $\alpha_{k} \le \bar{\delta} / (2 \sqrt{n})$ for $k \geqslant \bar{k}$. Applying Lemma 2.3, we obtain the assertion of the lemma.
\begin{corollary}
\label{cor:2.1}
$$
\lim _{k \rightarrow \infty} \rho\left(\nabla f\left(x, \alpha_{k}\right), \partial f(x)\right)=0 .
$$
\end{corollary}

Lemma 2.5 means that for $\alpha_{k} \rightarrow 0$ the gradients of the smoothed function $f\left(x, \alpha_{k}\right)$ calculated in a sufficiently small neighborhood of the point $x$, for sufficiently large $ k$ are close to the set of generalized gradients of the function $f(x)$ at the point $x$. This property plays an important role in the study of the convergence of numerical methods for minimizing Lipschitz functions.
\begin{lemma}
\label{lem:2.6}
{Let}
$$
\left\{x^{k}\right\} \rightarrow x, \quad \nabla f\left(x^{k}, \alpha_{k}\right) \rightarrow g, \quad \alpha_{k} \rightarrow 0, \quad k \rightarrow \infty .
$$
\textit{Then $g \in \partial f(x) $}.
\end{lemma}

The proof of this assertion follows directly from Lemma 2.5.

Let us continue studying the properties of locally Lipschitz functions. An analogue of Lagrange's mean value theorem holds for them.
\begin{theorem}
\label{th:2.9}
{Let the points $x, y$ belong to a convex set of $E_{n}$. Then there exist $g \in \partial f(x+t(y-x)), \quad t \in[0,1],$ such that}
\[f(y)-f(x)=\left<g, y-x\right>\tag{2.12}\]
(One can even assert [159] the existence of $t \in(0,1)$. 
\end{theorem}

The proof follows from the mean value formula for $f\left(x, \alpha_{k}\right)$:
$$f\left(y, \alpha_{k}\right)-f\left(x, \alpha_{k}\right)=\left<\nabla f\left(x+t_{k}(y-x), \alpha_{k}\right), y-x\right>, \quad 0 \le t_{k} \le 1.$$ 
Passing to the limit as $k \rightarrow \infty$ in this inequality and applying Lemma 2.6, we obtain relation (2.12).

The multivalued mapping $x \rightarrow \partial f(x)$ has an interesting feature.
\begin{theorem}
\label{th:2.10}
{The following statements are equivalent:}

\textit{a) The set $\partial f(x)$ consists of a single vector $g(x)$.}

\textit{b) At point $x$ there is a gradient $\nabla f(x)=g(x)$ which is continuous at $x$ with respect to the set of points where it is defined.}
\end{theorem}

{\it P r o o f.} Condition b) and Theorem 2.4 directly imply a). Let us show that b) follows from a). The existence of $\nabla f(x)$ at the point $x$ follows from relation (2.12) and the closedness of the mapping $x \rightarrow \partial f(x)$, while the continuity of $\nabla f(x)$ at the point $x$ follows directly from condition a).

The following property characterizes the set of generalized gradients of the maximum function.
\begin{theorem}
\label{th:2.11}
{Let $f_{i}(x) (i=1, \ldots, m)$ satisfy the local Lipschitz condition. Then the maximum function}
$$
\varphi(x)=\max _{1 \le i \le m} f_{i}(x)
$$
{belongs to this class and the floowing formula is valid:}
\[
\partial \varphi(x) \subset \operatorname{co}\left\{\partial f_{i}(x) \mid i \in I(x)=\left\{i \mid \varphi(x)= f_{i}(x)\right\}\right\} \tag{2.13}
\]
\end{theorem}

{\it P r o o f.} Let $x, y$ be arbitrary points of some bounded set for which the Lipschitz constant of the function $f_{i}(x)$ is $L_{i} (i=1, \ldots, m)$. Then
$$
\varphi(x)=f_{v}(x) \le f_{v}(y)+L_{v}\|x-y\| \le \varphi(y)+L\|x-y\|,
$$
where
$$
v \in I(x), \quad L=\max _{1 \le i \le m} L_{i}
$$
Since $x, y$ are arbitrary, the reverse inequality also holds. This implies that $\varphi(x)$ satisfies the local Lipschitz condition. 

From formula (2.12) it follows that
\[
\begin{aligned}
& \frac{f_{i}(x+\lambda e)-f_{i}(x)}{\lambda}=\left<g^{i}(x+t \lambda e), e\right>, \quad 0 \le t \le 1, \\
& g^{i}(x+t \lambda e) \in \partial f_{i}(x+t \lambda e), \quad i=1, \ldots, m .
\end{aligned}\tag{2.14}
\]

Since the set $\{1, \ldots, m\}$ is finite, there exists a subsequence $\left\{\lambda_{s}\right\} \rightarrow 0, s \rightarrow \infty,$ such that
\[
i_{\lambda_{s}} \in I\left(x+\lambda_{s} e\right), \quad i_{\lambda_{s}}=v .\tag{2.15}
\]

It is easy to see that $v \in I(x)$. If $v \notin I(x)$, then $f_{v}(x)<\varphi(x)$, and from the continuity of $f_{v}$, $\varphi$ for sufficiently small $\lambda_{s}$ it follows
$$
f_{v}\left(x+\lambda_{s} e\right)<\varphi\left(x+\lambda_{s} e\right),
$$
that contradicts (2.15).

Let $\nabla \varphi(x)$ exist at the point $x$. Since the mapping $x \rightarrow \partial f(x)$ is closed, passing to the limit as $\lambda_{s} \rightarrow 0$ in the equality
$$
\frac{\varphi\left(x+\lambda_{s}e\right)-\varphi(x)}{\lambda_{s}}=\frac{f_{v}\left(x+\lambda_{s}e \right)-f_{v}(x)}{\lambda_{s}},
$$
from relation (2.14), we obtain
\[
\nabla \varphi(x)=g^{v}(x) \in \partial f_{v}(x) .\tag{2.16}
\]

Let there be a sequence $x^{k} \rightarrow x$ of points where $\nabla \varphi(x)$ exists and $\lim\limits_{k \rightarrow \infty} \nabla \varphi\left(x^{ k}\right)=g(x)$. 
From relations (2.16) it follows that
$$
\nabla \varphi\left(x^{k}\right) \in \partial f_{v_{k}}\left(x^{k}\right), \quad v_{k} \in I\left( x^{k}\right) .
$$
That's why
$$
g(x) \in \partial f_{i}(x), \quad i \in I(x) .
$$
Hence,
$$
\partial \varphi(x) \subset \operatorname{co}\left\{\partial f_{i}(x) \mid i \in I(x)\right\},
$$
that was required to prove.

\smallskip
Let us show that in relation (2.13) a strict inclusion can hold. Let $f_{1}(x)=-|x|, f_{2}(x)=x$. Then
$$
f(x)=\max \left\{f_{1}(x), f_{2}(x)\right\}=x;
$$
$\partial f(0)=\{1\}$, while $\partial f_{1}(0)=[-1,1]$. 

The situation is similar for compound functions. Denote by $\left[v^{1} \ldots v^{m}\right]$ the $n \times m$-matrix, where $i$-th column is the vector $v^{i} \in E_{ n}$ for $i=1, \ldots, m$.
\begin{theorem}
\label{th:2.12}
{Let $f_{i}: E_{n} \rightarrow E_{1}, i=1, \ldots, m$ and let $F: E_{m} \rightarrow E_{1} are $ locally Lipschitz functions. For $x \in E_{n}$ we define}
$$
\begin{aligned}
& Y(x)=\left(f_{1}(x), \ldots, f_{m}(x)\right), \quad \Phi(x)=F(Y(x)), \\
& H(x)=c o\left\{g \in E_{n} \mid g=\left[g^{1} \ldots g^{m}\right] \omega\right\}, \\
& g^{i} \in \partial f_{l}(x), \quad i=1, \ldots, m, \quad \omega \in \partial F(Y(x)) .
\end{aligned}
$$
{Then the function $\Phi(x)$ satisfies the local Lipschitz condition and the formula holds}
\[
\partial \Phi(x) \subset H(x). \tag{2.17}
\]
\end{theorem}

{\it P r o o f.} By Theorem 2.9,
$$
\frac{F(Y(x+v+\lambda e))-F(Y(x+v))}{\lambda}=\left<\bar{\omega}, Y(x+v+\lambda e)-Y (x+v)\right>,
$$
where
$$
\bar{\omega} \in \partial F\left(Y(x+v)+t(Y(x+v+\lambda \varepsilon)-Y(x+v))\right). 
$$
Similarly,
$$
\begin{aligned}
f_{i}(x+v+\lambda e)-f_{i}(x+v) & =\lambda\left<\bar{g}^{i}, e\right>, \\
\bar{g}^{i} \in \partial f_{i}(x+v+t \lambda e), \quad i & =1, \ldots, m .
\end{aligned}
$$
Therefore, the definition of $f^{0}(x ; e)$ and Theorem 2.5 imply that
$$
f^{0}(x ; e) \le \max _{u \in H(x)}\left<u, e\right>.
$$
Next, applying Theorem 2.3, we obtain the required relation (2.17).

\smallskip
Formulas (2.9), (2.13), (2.17) show that it is impossible to construct gradient calculi for Lipschitz functions, as was the case for generalized differentiable functions. Thus, generalized gradient descent methods cannot be widely applied in Lipschitz programming. Therefore, finite-difference methods, which are explored in the subsequent chapters, are of great importance.

\section*{$\S$ 3. Necessary conditions for an extremum}
\label{Sec.3}
\setcounter{section}{3}
\setcounter{definition}{0}
\setcounter{equation}{0}
\setcounter{theorem}{0}
\setcounter{lemma}{0}
\setcounter{remark}{0}
\setcounter{corollary}{0}
\setcounter{example}{0}

Clarifying the necessary conditions for an extremum is important for several reasons. In many cases, these conditions help predict the structure of optimization methods. They play an important role in the study of the convergence of algorithms. It is expedient to derive such conditions that would be applicable to a wide class of problems, so as not to build their own relations in each individual case.

Necessary extremum conditions in minimization problems for generalized differentiable and Lipschitz functions, in contrast to convex functions, characterize only local properties. In other words, if the necessary extremum conditions are satisfied at the point $x=x^{*}$, then, without making any additional assumptions about the behavior of the function, we can hypothesize that at the point $x=x^{*}$ local extremum of $f(x)$ takes place. In this section, the necessary conditions for these classes are formulated in a form that is most convenient for studying the validity of the minimization methods considered below.

First, we briefly consider the necessary extremum conditions in mathematical programming problems with generalized differentiable functions. The presentation will be rather illustrative in order to show the fundamental possibility of deriving the necessary conditions for an extremum only from Definition 1.1. A detailed study is not necessary here, since it will be carried out below for a more general class of locally Lipschitz functions. The point is that generalized differentiable functions are locally Lipschitz. Moreover, as shown in Theorem 1.10, the generalized Clarke gradients of these functions are always pseudogradients of these functions as well; therefore, the necessary extremum conditions for problems with generalized differentiable functions in terms of generalized Clarke gradients coincide with the necessary extremum conditions in terms of pseudogradients.

Consider the following problem:
\[
f(x) \rightarrow \min , \quad h(x) \le 0, \tag{3.1}
\]
where $f(x)$, $h (x)$ are generalized differentiable functions. Let us define
a multivalued mapping $G: x \rightarrow G (x)$:
\begin{equation}\tag{3.2}
G(x)=\left\{
\begin{array}{ll}
G_f(x),& h(x)<0,\\
co\{G_f(x)\cup G_h(x)\},&h(x)=0,\\
G_h(x),& h(x)>0.
\end{array}
\right.
\end{equation}
\begin{theorem}
\label{th:3.1}
Let $x^*$ be the point of minimum of $f(x)$ under the constraint $h(x)\le 0.$ Then
\begin{equation}\tag{3.3}
0\in G(x)
\end{equation}
(John type condition). If the regularity condition is satisfied
\[
\inf\{\|g\|\,|\,g\in G_h(x),\,h(x)=0\}=\gamma>0,
\]
then there exist such
$\lambda^*\ge 0,$ $g_f^*\in G_f(x^*),$ $g_h^*\in G_h(x^*),$ that
\begin{equation}\tag{3.4}
g_f^*+\lambda^* g_h^*=0,\;\;\;\;\;\lambda^*h(x^*)=0
\end{equation}
(Kuhn-Tucker type condition).
\end{theorem}

{\it P r o o f.} Condition (3.4) follows obviously from (3.3). The regularity condition is introduced to ensure that the coefficient on the gradient of the function $f(x)$ satisfying 
(3.3) does not go to zero. The situation when it is equal to zero is usually not interesting, since the minimized function is not present in the necessary extremum conditions.

Let $h(x^*) = 0$ (the case $h(x^*) < 0$ is proved analogously). Suppose the contrary to the statement of the theorem. Then by convexity and closedness of 
$G(x^*)$,
\[
\rho(0,G(x^*))=\underset{g\in G(x^*)}{\min}\|g\|=\|d\|>0,
\]
and $d\in G(x^*).$ In the neighborhood of the point $x = x^*$ the following expantions are valid
\[
\begin{array}{lcl}
f(y)&=&f(x^*)+\left<g_f,y-x^*\right>)+o_f(x^*,y,g_f),\\
h(y)&=&\left<g_h,y-x^*\right>+o_h(x^*,y,g_h),
\end{array}
\]
where
\[
\begin{array}{lcl}
o_f(x^*,y,g_f)/\|y-x^*\|&\rightarrow &0,\\
o_h(x^*,y,g_h)/\|y-x^*\|&\rightarrow &0,
\end{array}
\]
uniformly in $y\rightarrow x^*$, and in $g_f\in G_f(y)$ and in $g_h\in G_h(y)$. It is obvious that the mapping $G$ is locally bounded and closed; hence it is semi-continuous from above. Consider $y = x^ * - td,$ $g\in G(y),$ $t\ge 0$. At small enough $t$, by virtue of the properties of the functions $o_f$ and $o_h$ and the semi-continuity from above of $G$, we have
\[
\begin{array}{c}
o_f(x^*,y,g_f)\le t\|d\|^2/3,\;\;\;\;\;\;o_h(x^*,y,g_h)\le t\|d\|^2/3,\\
\rho(0,G(x^*))=\underset{g^\prime\in G(x^*)}{\min}\|g-g^\prime\|\le\|d\|/3,\\
\left<g,y-x^*\right>\le -(2/3)\|d\|^2t.
\end{array}
\]
Then
\[
f(x^*-td)\le f(x^*)-t\|d\|^2/3,\;\;\;\;\;h(x^*-td)\le -t\|d\|^2/3,
\]
which contradicts the local optimality of $x^*$. The theorem is proved.

Let us now consider the necessary conditions for extrema of the Lipschitz  functions.
\begin{theorem}
\label{th:3.2}
For a locally Lipschitzian function $f(x),$ $x\in E_n,$ to have a local minimum at point $x= x^*,$ it is necessary that
\[ 0\in \partial f(x^*).  \]
\end{theorem}

{\it P r o o f.}
Suppose that $0\notin \partial f(x^*)$. By the separability theorem, there exists a vector $a$ and a number $\varepsilon >0$ for which 
$\left<a, g(x^*)\right> < -\varepsilon,$ where $g(x^*)$ is an arbitrary vector of the set $\partial f(x^*)$. Applying the Lagrange formula, we have
\[
f(y,\alpha_k)=f(x^*,\alpha_k)+\left<\nabla f(x^*+\theta(y-x^*),\alpha_k),y-x^*\right>,
\;\;\;0\le\theta\le 1.
\]

Consider the points $y = x^* + ta,$ $t >0.$ Then by virtue of Lemma 2.5, if $t$ is sufficiently small,
\[
\left<a,\nabla f(x^*+\theta(y-x^*),\alpha_k)\right>\le -\varepsilon/2,
\;\;\;k\ge\bar{k}.
\]
Hence,
\[
f(y,\alpha_k)\le f(x^*,\alpha_k)-\lambda\varepsilon/2. 
\]
Passing to the limit by $\alpha_k\rightarrow 0$ in this inequality, we obtain
\[
f(y)\le f(x^*)-\lambda\varepsilon/2.
\]
This relation contradicts the fact that $x^*$ is the point of local minimum of the function $f(x).$

Consider the problem
\begin{equation}\tag{3.5}
f_0(x)\rightarrow\min_x
\end{equation}
subject to constraints
\begin{equation}\tag{3.6}
f_i(x)\le 0,\;\;\;\;\; i=1,\ldots,m,
\end{equation}
where $f_i(x)$ $(i=0,\ldots,m)$ satisfy the local Lipschitz condition.
\begin{theorem}
\label{th:3.3}
If $x^*$ is a solution of the problem (3.5), (3.6), then there exist such numbers $\lambda_0$, $\lambda_i\ge 0$ $(i\in I(x^*))$ not all equal to zero that
\begin{equation}\tag{3.7}
\begin{array}{c}
\lambda_0 f_0^0(x^*;e)+\sum_{i\in I(x^*)}\lambda_i f_i^0(x^*;e)\ge 0
\;\;\;\;\;\forall e\in E_n ,\\
I(x^*)=\{i|\,f_i(x^*)=0\}.
\end{array}
\end{equation}
\end{theorem}

{\it P r o o f.} 
Suppose the opposite: there exists such a vector $v$ that for any numbers 
$\lambda_0$, $\lambda_i$, not all equal to zero, the inequality holds
\[
\lambda_0 f_0^0(x^*;v)+\sum_{i\in I(x^*)}\lambda_i f_i^0(x^*;v)< 0
\]

Hence, it follows that
$$
f_{0}^{0}\left(x^{*} ; v\right)<0,\;\;\; f_{i}^{0}\left(x^{*} ; v\right)<0, \;\;\;i \in I\left(x^{*}\right)
$$
Therefore, for sufficiently small numbers $t>0$
$$
f_{0}\left(x^{*}+t v\right)<f_{0}\left(x^{*}\right),\;\;\; f_{i}\left(x^{*}+t v\right)<f_{l}\left(x^{*}\right)=0, \;\;\;i \in I\left(x^{*}\right).
$$
Thus, for small $t>0$, the point $x^{*}+t v$ satisfies all the restrictions, since due to the continuity of $f_{i}(x) (i=0, \ldots, m)$, for small $t$ holds
$$
f_{i}\left(x^{*}+t v\right)<0, \;\;\;i \notin I\left(x^{*}\right)
$$
and
$$
f_{0}\left(x^{*}+t v\right)<f_{0}\left(x^{*}\right)
$$
This contradicts the fact that $x^{*}$ is a solution to problem $(3.5),(3.6)$. The theorem has been proven.

Relations (2.4) and (3.7) imply that there are $\lambda_{0}, \lambda_{i}\ge 0$, not all equal to zero, for which
\begin{equation}
\lambda_{0} \max _{g^{0} \in \partial f_{0}\left(x^{*}\right)}\left<g^{0}, e\right>+\sum_{i \in I\left(x^{*}\right)} \lambda_{i} \max _{g^{i} \in \partial f_{i}\left(x^{*}\right)}\left<g^{i}, e\right> \geq 0 \;\;\;\forall \;e \in E_{n}.\tag{3.8}
\end{equation}
\begin{lemma}
\label{lem:3.1}
If inequality (3.8) holds, then
\begin{equation}
0 \in \lambda_{0} \partial f_{0}\left(x^{*}\right)+\sum_{i \in I\left(x^{*}\right)} \lambda_{i} \partial f_{i}\left(x^{*}\right)\tag{3.9}
\end{equation}
\end{lemma}

{\it P r o o f.} Consider a convex set
$$
P=\left\{\lambda_{0} \partial f_{0}\left(x^{*}\right)+\sum_{i \in I\left(x^{*}\right)} \lambda_{i} \partial f_{i}\left(x^{*}\right)\right\}.
$$
Let $0 \notin P$. Then, by the separation theorem, there are vector $a$ and number $\varepsilon>0$ such that for any
$$
g^{0} \in \partial f_{0}\left(x^{*}\right),\;\;\; g^{i} \in \partial f_{i}\left(x^{*}\right), \;\;\;i \in I\left(x^{*}\right),
$$
the inequality holds
$$
\lambda_{0}\left<a, g^{0}\right>+\sum_{i \in I\left(x^{0}\right)} \lambda_{i}\left<a, g^{i}\right>\le -\varepsilon,
$$
that contradicts (3.9). 

The Kuhn-Tucker theorem takes place when an additional regularity condition is satisfied.
\begin{theorem}
\label{th:3.4}
Let the following condition be satisfied: there exists $d$ such that
\begin{equation}
f_{i}^{0}\left(x^{*} ; d\right)<0 \;\;\;\forall i \in I\left(x^{*}\right).
\tag{3.10}
\end{equation}
Then there are $\lambda_{i}\geq 0 \;(i=1, \ldots, m)$ such that
\begin{equation}
\begin{array}{r}
0 \in \partial f_{0}\left(x^{*}\right)+\sum_{i=1}^{m} \lambda_{i} \partial f_{i}\left(x^{*}\right), \\
\lambda_{i} f_{i}\left(x^{*}\right)=0,\; i=1, \ldots, m.
\end{array}
\tag{3.11}
\end{equation}
\end{theorem}

The proof obviously follows from relation (3.7).

\smallskip
Denote by $\Lambda$ the set of vectors $\lambda=\left(\lambda_{1}, \ldots, \lambda_{m}\right) \geq 0$ satisfying relation (3.11). From the definition of the sets $\partial f_{0}(x), \partial f_{i}(x) \;(i=1, \ldots, m)$ and the closedness of the mappings $x \rightarrow \partial f_{i}(x) \;(i=0,1, \ldots m)$ it follows that $\Lambda$ is convex and closed. As will be shown below, (3.10) is a necessary and sufficient condition for $\Lambda$ to be bounded, i.e., in this sense, it is the weakest possible regularity condition.
\begin{theorem}
\label{th:3.5}
Condition (3.10) is equivalent to the boundedness of the set $\Lambda$.
\end{theorem}

{\it P r o o f.} Let condition (3.10) be satisfied. Then it follows from the relation $(3.11)$ and the definition of $\partial f_{i}(x) \;(i=0,1, \ldots, m)$ that
$$
\begin{gathered}
f_{0}^{0}\left(x^{*} ; e\right)+\sum_{i=1}^{m} \lambda_{i} f_{i}^{0}\left(x ^{*};e\right) \geq 0 \;\;\;\forall e\, \in E_{n}, \\
\lambda_{i}\geq 0, \;\;\;\lambda_{i} f_{i}\left(x^{*}\right)=0,\;\; i=1, \ldots, m,
\end{gathered}
$$
which immediately implies that $\Lambda$ is bounded.

Now let $\Lambda$ be non-empty and bounded. We form the set
$$
M\left(x^{*}\right)=\operatorname{co} \bigcup_{i \in I\left(x^{*}\right)} \partial f_{i}\left(x^{* }\right)
$$
and show that $0 \notin M\left(x^{*}\right)$. Otherwise,
\begin{equation}\tag{3.12}
\begin{array}{l}
0=\sum_{i \in l\left(x^{*}\right)} \alpha_{i} h^{i}, \\
\alpha_{i }\ge 0,\;\;\; h^{i} \in \partial f_{i}\left(x^{*}\right), \;\;\;\sum_{i \in I\left(x^{*}\right)} \alpha_{i}=1 .
\end{array}
\end{equation}
Relation (3.11) implies that there exist such $\lambda_{i}\ge 0$,
$g^{0} \in \partial f_{0}\left(x^{*}\right)$, $g^{i} \in \partial f_{i}\left(x^{* }\right)\;\left(i \in I\left(x^{*}\right)\right)$, that
\begin{equation}\tag{3.13}
0=g^{0}+\sum_{i \in M\left(x^{*}\right)} \lambda_{i} g^{i}.
\end{equation}

Multiplying the first equality in (3.12) by $\sigma>0$ and adding it to $(3.13)$, we get
\begin{equation}\tag{3.14}
0=g^{0}+\sum_{i \in I\left(x^{*}\right)}\left(\lambda_{i} g^{i}+\sigma \alpha_{i} h^ {i}\right)
\end{equation}

Let's put
$$
\lambda_{i}(\sigma)=\lambda_{i}+\sigma \alpha_{i} .0
$$
for $i \in I\left(x^{*}\right)$ and $\lambda_{i}(\sigma)=0$ for $i \notin I\left(x^{*}\right).$ If $\lambda_{i}(\sigma)>0$, then
$$
P^{i}(\sigma)=\left[\lambda_{i}(\sigma)\right]^{-1}\left(\lambda_{i} g^{i}+\sigma \alpha_{i } h^{i}\right) \in \partial f_{i}\left(x^{*}\right).
$$

Therefore, from (3.14) we have
$$
0=g^{0}+\sum_{i=1}^{m} \lambda_{i}(\sigma) P^{i}(\sigma).
$$
Since there is $\alpha_{l}>0$ for $l \in I\left(x^{*}\right)$, then
$$
\lambda_{l}(\sigma)=\lambda_{l}+\sigma \alpha_{l} \rightarrow \infty \text { as } \sigma \rightarrow \infty,
$$
that contradicts the boundedness of $\Lambda$.

Since $0 \notin M\left(x^{*}\right)$, then by the separation theorem there are vector $d$ and number $\varepsilon>0$ such that for $i \in I\left(x^{ *}\right)$ we have $\left<d, g^{i}\right>\leq -\varepsilon$ for all $g^{i} \in \partial f_{i}\left(x^{*} \right).$ Therefore, $f_{i}^{0}\left(x^{*} ; d\right)<0\left(i \in I\left(x^{*}\right)\right).$ The theorem is proved.

\newpage
\begin{flushright}
CHAPTER 2
\label{Ch.2}

\textbf{MINIMIZATION OF LIPSCHITZ FUNCTIONS \\WITHOUT CALCULATION OF GRADIENTS
}
\underline{\hspace{15cm}}
\end{flushright}
\bigskip\bigskip\bigskip\bigskip\bigskip\bigskip

\section*{$\S$ 4. Finite-difference method for minimizing Lipschitz functions}
\label{Sec.4}
\setcounter{section}{4}
\setcounter{definition}{0}
\setcounter{equation}{0}
\setcounter{theorem}{0}
\setcounter{lemma}{0}
\setcounter{remark}{0}
\setcounter{corollary}{0}
\setcounter{example}{0}
\numberwithin{equation}{section}
In the case when it is required to minimize a non-smooth function, the movement along the generalized gradient or in the direction of the difference approximation of the gradient does not, in general, give a monotonic change either in the function itself or in the distance to the extremum point. In this regard, an interesting question arises of studying the convergence of methods similar to finite-difference or gradient methods for minimizing non-convex non-smooth functions.

From a practical point of view, the situation when information about gradients is not available is most typical in mathematical programming problems. When minimizing functions using gradient methods, the user needs to prepare ${n}$ subroutines for calculating the gradient components. If the dimension of the space is large, then this auxiliary procedure can take considerable time. In addition, in high-dimensional problems it is very difficult to formulate generalized gradient formulas. In smooth optimization, if the function is analytically complex, the components of the gradient are approximated by difference relations; for example,
\begin{equation}
\label{eqn:4.1}
\nabla f(x) \approx \sum_{i=1}^{n} \frac{f\left(x+\alpha e_{i}\right)-f(x)}{\alpha} e_{i}
\end{equation}

For non-smooth functions, the difference approximation (4.1) is inapplicable, since it leads to non-converging methods. This can be easily verified by considering the following example. Let each level line of a convex function be the boundary of a square. When calculating (4.1) in the lower left corner of the square, we get a zero value, although all subgradients at this point can be different from zero. At the same time, a slight change in scheme (4.1) leads to the fact that it can be used in the process of searching for the extremum of nonsmooth functions.

As for the gradient method,
\begin{equation}
\label{eqn:4.2}
x^{k+1}=x^{k}-\rho_{k} g\left(x^{k}\right),\;\;\; g\left(x^{k}\right) \in \partial f\left(x^{k}\right),
\end{equation}
its use for minimizing Lipschitz functions also does not guarantee the convergence of the sequence $\left\{x^{k}\right\}$ to stationary points of $f(x)$. For functions of this class, the direction $-g\left(x^{k}\right)$ does not lead to a decrease in the values of $f(x^{k})$ or the distance to the minimum point. Therefore, the question of choosing the direction of descent in methods like (4.2) is solved differently than for smooth and convex functions.

The substantiation of numerical methods for minimizing locally Lipschitz functions is carried out as follows. Consider the sequence of smoothed functions
$$
f(x, \alpha)=\frac{1}{(2 \alpha)^{n}} \int_{x_{1}-\alpha}^{x_{1}+\alpha} \ldots \int_{ x_{n}-\alpha}^{x_{n}+\alpha} f\left(y_{1}, \ldots, y_{n}\right) d y_{1} \ldots d y_{n}\nonumber
$$

As noted in $\S$ 2, the functions $f(x, \alpha)$ converge uniformly in any bounded domain of $E_{n}$ to $f(x)$ as $\alpha \rightarrow 0$. The question arises: is it possible to use the gradient of the function $f(x, \alpha)$ as a descent vector in the process of minimization, and how to change $\alpha$, tending it to zero, in order to obtain the solution of the original problem in the limit?

The transition to averaged functions $f(x, \alpha)$ is also important because it allows us to exclude the influence of small local singularities of the function $f(x)$ and is often useful for finding a global extremum in solving multiextremal problems. It may turn out that the function $f(x)$ has many local extrema, although there is also a pronounced global extremum. Minimization of the function $f(x)$ can lead to one of the local extrema, while minimization of the smoothed function can give a good approximation to the global minimum, if the parameter $\alpha$ is changed in the course of calculations, tending it to zero.

In $\S$ 2 it was shown that the partial derivative $\partial f(x, \alpha) / \partial x_{1}$ of the function $f(x, \alpha)$ is calculated by the formula
\begin{eqnarray}
 \frac{\partial f\left(x,\alpha\right)}{\partial x_{1}}
&=&
\frac{1}{(2 \alpha)^{n}} \int_{x_{ 2}-\alpha}^{x_{2}+\alpha} \ldots \int_{x_{n}-\alpha}^{x_{n}+\alpha}\left[f\left(x_{1} +\alpha, y_{2}, \ldots, y_{n}\right) \rightarrow\right. \nonumber\\
&& \left.-f\left(x_{1}-\alpha, y_{2}, \ldots, y_{n}\right)\right] d y_{2} \ldots d y_{n} \text { . }\label{eqn:4.3}
\end{eqnarray}
Derivatives
$\partial f(x, \alpha) / \partial x_{2}, \ldots, \partial f(x, \alpha) / \partial x_{n}$ are calculated similarly. 

The gradient of the function $f(x, \alpha)$ can be expressed differently:
\begin{eqnarray}
\label{eqn:4.4}
\begin{array}{c}
\nabla f(x, \alpha)=\frac{1}{(2 \alpha)^{h}} \int_{-\alpha}^{\alpha} \cdots \int_{-\alpha}^{\alpha} g(x+y) d y_{1} \ldots d y_{n}, \\
g(x+y) \in \partial f(x+y) .
\end{array}
\end{eqnarray}

The complexity of using $\nabla f(x, \alpha)$ in the process of minimizing $f(x)$ is related to the difficulty of calculating multidimensional integrals in formulas (4.3), (4.4). However, this problem can be solved by probabilistic methods.

Let us construct a random vector $H(x, \alpha)$ such that
\begin{equation}
\label{eqn:4.5}
\mathbb{E} H(x, \alpha)=\nabla f(x, \alpha)
\end{equation}
where $\mathbb{E}$ is the sign of mathematical expectation. It is easy to see that a random vector satisfying (4.5) is defined by the formula
\begin{equation}
\label{eqn:4.6}
H(x, \alpha)=\frac{1}{2 \alpha} \sum_{i=1}^{n}
\left[f\left(\tilde{x}_{1}, \ldots, x_{i}+\alpha, \ldots, \tilde{x}_{n}\right)
-f\left(\tilde{x}_{1}, \ldots, x_{i}-\alpha, \ldots, \tilde{x}_{n}\right)\right] e_{t}
\end{equation}
where $\tilde{x}_{i}(i=1, \ldots, n)$ are independent random variables uniformly distributed on the segments $\left[x_{i}-\alpha, x_{i}+\alpha\right](\alpha>0)$.

Taking a finite-difference approximation of the gradient of the smoothed function $f(x, \alpha)$, we obtain
\begin{equation}
\label{eqn:4.7}
H(x, \alpha)=\sum_{i=1}^{n} \frac{f\left(\tilde{x}+\Delta e_{i}\right)
-f(\tilde{x})}{\Delta} e_{i}
\end{equation}

Relation (4.7) differs from (4.1) in that it is calculated at a random point $x$ chosen (on the basis of a uniform distribution) in an $n$-dimensional cube centered at $x$ and side $2 \alpha$, and
$$
\mathbb{E} \sum_{i=1}^{n} \frac{f\left(\tilde{x}+\Delta e_{i}\right)
-f(\tilde{x})}{\Delta} e_{i}=\nabla f(x, \alpha)+b,
$$
where $\| b \|\le \sqrt{n}L \Delta /\alpha$.

If to apply formula (4.4), then the vector satisfying (4.5), is
\begin{equation}
\label{eqn:4.8}
H(x, \alpha)=g\left(\tilde{x}\right).
\end{equation}

Formulas (4.6)-(4.8) lead us to consider finite-difference and gradient methods for minimizing locally Lipschitz functions:
\begin{equation}
\label{eqn:4.9}
x^{k+1}=x^{k}-\rho_{k} H\left(x^{k}, \alpha_{k}\right),
\end{equation}
where the vector $H\left(x^{k}, \alpha_{k}\right)$ is determined by one of the formulas (4.6) - (4.8). In order not to introduce normalizing factors, we assume that the set of stationary points $X^{*}=\left\{x^{*} \mid 0 \in f^{*}(x)\right\}$ satisfying the necessary condition extreme, is bounded. Then, without loss of generality, we can say that the sequence $x^{k}$ defined by formulas (4.6)-(4.9) belongs to a bounded set.

It is easy to see that the sequence of points $x^{k}(k=0,1, \ldots)$ has a probabilistic nature. As is customary in probability theory, random vectors $x^{k}(\omega)(k=0,1, \ldots)$ are defined on the probability space $(\Omega, \Sigma, P)$ (which is the product of the spaces in which the points of $\tilde{x}$ are ejected, is usually called an elementary event of this space. In what follows, the dependence of $x^{k}$ on $\omega$ is usually omitted. 

Investigation of convergence of $\left\{x^ {k}\right\}$ to the set $X^{*} \subset E_{n}$ is carried out on the basis of the following sufficient conditions [91], slightly different from the more general conditions described in $\S$ 8. It turns out that the limit points of $\left\{x^{k}\right\}$ with probability 1 (almost surely) belong to some set $X^{*}$, if the following conditions are satisfied with probability 1:

I.
\begin{equation}\label{eqn:4.10}
\lim _{k \rightarrow \infty}\left\|x^{k+1}(\omega)-x^{k}(\omega)\right\|=0.
\end{equation}

II. There is a bounded set $X(\omega) \in E_{n}$ such that
\begin{equation}\label{eqn:4.11}
x^{k}(\omega) \in X(\omega), k=0,1, \ldots .
\end{equation}

III. For any subsequence $x^{s}(\omega)$ for which
$$
\lim _{s \rightarrow \infty} x^{s}(\omega)=x^{\prime}(\omega) \notin X^{*}, s \in S \subset\{1,2, \ldots\},\nonumber
$$
\noindent
there exists $\bar{\delta}(\omega)$ such that for all sufficiently large $S$ and all $\delta \in (0, \bar{\delta}(\omega)]$
\begin{equation}\label{eqn:4.12}
k(s)=\min \left\{r \mid x^{r}-x^{s}>\varepsilon, r>s\right\}<\infty .
\end{equation}

IV. There is a continuous function $W(x)$ for which
\begin{equation}\label{eqn:4.13}
\varlimsup_{s \rightarrow \infty} W\left(x^{k(s)}(\omega)\right)<\lim _{s \rightarrow \infty} W\left(x^{s}(\omega)\right)=W\left(x^{\prime}(\omega)\right) .
\end{equation}

V. The set
$$
W^{*}=\left\{W(x) \mid x \in X^{*}\right\}
$$
does not contain intervals.

The most important feature of  conditions I--V is that they do not require the entire sequence $\left\{W\left(x^{k}\right)\right\}$ to be monotone. Conditions III, IV mean that a certain monotonicity is observed only on some specially constructed subsequences. 

Let us now formulate conditions on the step size $\rho_{k}$ and the choice of parameters $\Delta_{k}, \alpha_{k}$ in formulas (4.6) - (4.8), under which the sequence $\left\{x^{ k}\right\}$ converges in the indicated sense to the set of points satisfying the necessary extremum condition.
\begin{theorem}\label{th:4.1}
Let the following conditions fulfill:
$$
\begin{gathered}
\sum_{k=0}^{\infty} \rho_{k}=\infty, \;\;\;\;\;
\sum_{k=0}^{\infty} \rho_{k}^{2}<\infty \\
\rho_{k} / \alpha_{k} \rightarrow 0,\;\;\; \Delta_{k} / \alpha_{k} \rightarrow 0,\;\;\;\left|\alpha_{k+1}-\alpha_{k}\right| / \rho_{k} \rightarrow 0,\;\;\; \alpha_{k} \rightarrow 0.
\end{gathered}
$$
Then all the limit points of the sequence $\left\{x^{k}\right\}$ defined by formulas (4.6)-(4.9) belong to the 
set $X^{*}=\left\{x^{*} \mid 0 \in \partial f\left(x^{*}\right)\right\}$
with probability 1.
\end{theorem}

{\it P r o o f.} Relation (4.10) follows directly from the conditions of the theorem and the boundedness of $H\left(x^{k}, \alpha_{k}\right)$. The main difficulty is the derivation of relations (4.11), (4.12). For simplicity of notation, we carry out the proof for the case (4.5).

For smoothed functions $f\left(x, \alpha_{k}\right)$ the Lagrange finite increment formula is valid:
$$
\begin{gathered}
f\left(x^{k+1}, \alpha_{k}\right)=f\left(x^{k}, \alpha_{k}\right)+\left<\nabla f\left(x^{k}+\tau\left(x^{k+1}-x^{k}\right), \alpha_{k}\right), x^{k+1}-x^{k}\right>, \\
0\leq  \tau \leq 1.
\end{gathered}
$$
Bounded values will be denoted by symbol $C$. Lemma 2.4 implies
\begin{eqnarray}
&&f\left(x^{k+1}, \alpha_{k}\right)=f\left(x^{k}, \alpha_{k}\right)\nonumber \\
&&\;\;\;+\left<\nabla f\left(x^{k}+\tau(x^{k+1}-x^{k}), \alpha_{k}\right)
-\nabla f(x^{k}, \alpha_{k})+\nabla f(x^{k}, \alpha_{k}), x^{k+1}-x^{k}\right> \nonumber \\
&&\leq f\left(x^{k}, \alpha_{k}\right)+C \rho_{k}^{2} / \alpha_{k}\nonumber\\
&&\;\;\;-\rho_{k}\left<\nabla f\left(x^{k}, \alpha_{k}\right), H\left(x^{k}, \alpha_{k}\right)-\nabla f\left(x^{k}, \alpha_{k}\right)+
\nabla f\left(x^{k}, \alpha_{k}\right)\right>\nonumber\\
&& \leq f\left(x^{k}, \alpha_{k}\right)-\rho_{k} \|\nabla f\left(x^{k}, \alpha_{k}\right)\|^{2}+C \rho_{k}^{2} / \alpha_{k}\nonumber\\
&&\;\;\;+\rho_{k}\left<\nabla f(x^{k}, \alpha_{k})-H(x^{k}, \alpha_{k}), \nabla f(x^{k}, \alpha_{k})\right>.\;\;\label{eqn:4.14}
\end{eqnarray}

Let us estimate the quantity $\left|f\left(x, \alpha_{k+1}\right)-f\left(x, \alpha_{k}\right)\right|$. Making a change of variables, we get
$$
\begin{gathered}
f\left(x, \alpha_{k}\right)=\frac{1}{2^{n}} \int_{-1}^{1} \ldots \int_{-1}^{1} f\left(x+\alpha_{k} z\right) d z_{1} \ldots d z_{n}, \\
f\left(x, \alpha_{k+1}\right)=\frac{1}{2^{n}} \int_{-1}^{1} \cdots \int_{-1}^{1} f\left(x+\alpha_{k+1} z\right) d z_{1} \ldots d z_{n}.
\end{gathered}
$$
Therefore, since $f(x)$ satisfies the Lipschitz condition, we have
$$
\left|f\left(x, \alpha_{k+1}\right)-f\left(x, \alpha_{k}\right)\right|\leq C\left|\alpha_{k}-\alpha_{k+1}\right|.
$$

Inequality (4.14) implies
\begin{eqnarray}
 f\left(x^{k+1}, \alpha_{k+1}\right) \quad 
&\le &
f\left(x^{k}, \alpha_{k}\right)-\rho_{k} \|\nabla f\left(x^{k}, \alpha_{k}\right)\|^{2}
\nonumber\\
&&+\rho_{k}\left<\nabla f\left(x^{k}, \alpha_{k}\right) 
-H\left(x^{k}, \alpha_{k}\right), \nabla f\left(x^{k}, \alpha_{k}\right)\right>\nonumber\\
&&+C\left|\alpha_{k}-\alpha_{k+1}\right|+C \rho_{k}^{2} / \alpha_{k}.
\label{eqn:4.15}
\end{eqnarray}

Let us show that condition (4.12) is satisfied. Let there be a subsequence $\left\{x^{s}(\omega)\right\} \rightarrow x^{\prime} \notin X^{*}$. Then one can find positive numbers 
$\bar{s}$ and $\bar{\delta}$ such that the $(2 \bar{\delta})$-neighborhoods of the points 
$x^{s} (s>\bar{s}, s\in S)$ do not intersect with $X^{*}$.

Construct a neighborhood of radius $\delta\leq \bar{\delta} / 2$ at the point $x^{s}\; (s\geq \bar{s})$. Let us assume that condition (4.12) is not satisfied. Then all points $x^{k}$ starting from some $k \geq \bar{s}$ are contained in the $\delta$-neighborhood of the point $x^{s}\; (s \in S)$. Therefore, Lemma 2.5 implies that for points $x^{k}$ located in a rather small $\delta$-neighborhood of the point $x^{s}\; (s\ge\bar{s})$, for sufficiently large numbers $k\ge\bar{s}$ the inequalities are fulfilled:
$$
\nabla f\left(x^{k}, \alpha_{k}\right)\geq\sigma> 0.
$$

According to the conditions of the theorem, starting from some number 
$k\ge\bar{s}$ we have
$$
\frac{C\left|\alpha_{k}-\alpha_{k+1}\right|}{\rho_{k}}+\frac{C \rho_{k}}{\alpha_{k}} \leq \frac{\sigma^{2}}{2} .
$$

Summing up the inequalities (4.15), we obtain
\begin{eqnarray}
&&f\left(x^{k+1}, \alpha_{k+1}\right)   
\leq f\left(x^{s}, \alpha_{s}\right)-\frac{\sigma^{2}}{2} \sum_{r=s}^{k} \rho_{r}
\nonumber\\
&&\;\;\;\;\;\;\;\;\;\;\;\;\;\;\;-\sum_{r=s}^{k} \rho_{r}\left[\frac{\sigma^{2}}{2}-\frac{C\left|\alpha_{r}-\alpha_{r+1}\right|}{\rho_{r}}-\frac{C \rho_{r}}{\alpha_{r}}\right] \nonumber\\
&&\;\;\;\;\;\;\;\;\;\;\;\;\;\;\;+\sum_{r=s}^{k} \rho_{r}\left<\nabla f\left(x^{r}, \alpha_{r}\right)-H\left(x^{r}, \alpha_{r}\right), \nabla f\left(x^{r}, \alpha_{r}\right)\right>\nonumber\\ 
&&\;\;\;\;\;\;\;\;\;\;\leq f\left(x^{s}, \alpha_{s}\right) 
-\frac{\sigma^{2}}{2} \sum_{r=s}^{k} \rho_{r}\nonumber\\ 
&&\;\;\;\;\;\;\;\;\;\;\;\;\;\;\;+\sum_{r=s}^{k} \rho_{r}\left<\nabla f\left(x^{r}, \alpha_{r}\right)-H\left(x^{r}, \alpha_{r}\right), \nabla f\left(x^{r}, \alpha_{r}\right)\right>.\label{4.16}
\end{eqnarray}

Let us use the following known fact [9].
\begin{lemma}\label{lem:4.1}
If $\left\{\varphi_{k}\right\}$ is a sequence of random variables such that
$$
\sum_{k=0}^{\infty} \mathbb{E} \varphi_{k}^{2}<\infty
$$
then the sequence of random variables
$$
\sum_{k=0}^{\infty} \left(\varphi_{k}-\mathbb{E}\left(\varphi_{k} \mid \varphi_{0}, \varphi_{1}, \ldots, \varphi_{k-1 }\right)\right)
$$
converges with probability 1.
\end{lemma}

In our case
$$
\begin{gathered}
\varphi_{k}=-\rho_{k}\left<H\left(x^{k}, \alpha_{k}\right), \nabla f\left(x^{k}, \alpha_{k }\right)\right>, \\
E\left(\varphi_{k} \mid \varphi_{0}, \varphi_{1}, \ldots, \varphi_{k-1}\right)=-\rho_{k} \|\nabla f\left(x ^{k}, \alpha_{k}\right)\|^{2} .
\end{gathered}
$$

Therefore, the following series converges with probability 1,
$$
\sum_{k=0}^{\infty} \rho_{k}\left<\nabla f\left(x^{k}, \alpha_{k}\right)-H\left(x^{k} , \alpha_{k}\right), \nabla f\left(x^{k}, \alpha_{k}\right)\right>.
$$

By the hypothesis of the theorem, the series $\sum_{k=0}^{\infty} \rho_{k}$ diverges, therefore, passing ${k}$ to the infinity in inequality (4.16), we obtain a contradiction with the boundedness of $\left\{f\left(x^{k}\right)\right\}$. This means that (4.12) holds:  there exists an index $k(s)<$ such that $k(s)=\min \left\{r \mid x^{r}-x^{s} >\delta, r>s\right\}$. 

It remains to derive condition (4.13). From the relations
$$
x^{k(s)}=x^{s}-\sum_{r=s}^{k(s)-1} \rho_{r} H\left(x^{r}, \alpha_{r }\right), 
\;\;\;\|x^{k(s)}-x^{s}\|>\delta,
$$
we get
$$
\sum_{r=s}^{k(s)-1} \rho_{r}>\frac{\delta}{C}
$$
Therefore, from inequality (4.16) for sufficiently large $s$, we have 
$$f\left(x^{k(s)}, \alpha_{k(s)}\right)\le f\left(x^{s}, \alpha_{s}\right)-\sigma^{2} \delta /(4 C),$$ 
whence, due to the uniform convergence of $f(x, \alpha) \rightarrow f(x)$, it follows
\begin{equation}\label{eqn:4.17}
\underset{s\rightarrow\infty}{\overline{\lim}}f\left(x^{k(s)}\right)
\le \underset{s\rightarrow\infty}{\lim}f\left(x^{s}\right)
-\sigma^2\delta/(4C),
\end{equation}
which proves (4.13). Thus, as function $W(x)$ we can take $f(x)$.

Based on the derived conditions (4.12), (4.13), we now show the convergence of the sequence $\left\{x^{k}\right\}$ to the set $X^{*}$. 

For any numbers $a$ and $b$ such that
$$
\overline{\lim } f\left(x^{k(s)}\right)<a<b<\lim _{s \rightarrow \infty} f\left(x^{s}\right)=f \left(x^{\prime}\right)
$$
the sequence $\left\{f\left(x^{k}\right)\right\}$ intersects the interval $(a, b)$ from left to right an infinite number of times. Then we can choose two sequences of points $\left\{x^{r}\right\}(r \in R),\left\{x^{p}\right\}(p \in P), R, P \subset\{1,2, \ldots, k, \ldots\}$ for which
$$
\begin{array}{r}
f\left(x^{r}\right)\le a, f\left(x^{r+1}\right)>a, \\
f\left(x^{p}\right) \ge b, f\left(x^{k}\right)>a, r<k<p .
\end{array}
$$
For simplicity, we assume that the sequence $\left\{x^{r}\right\}(r \in R)$ converges to some $x^{\prime \prime};$ otherwise via $\left\{x^{r}\right\}(r \in R)$ denote a convergent subsequence. Since $\rho_{k} \rightarrow 0$, then
$$
\lim _{r \rightarrow \infty} f\left(x^{r}\right)=a
$$
By virtue of the condition $\mathrm{V}$, the number $a$ can be chosen so that
$$
a \notin\left\{f(x) \mid x \in X^{*}\right\}
$$

For the sequence $\left\{x^{r}\right\}(r \in R)$, we perform the same procedures that have already been done for $x^{s}(S \in S)$. Let
$$
k(r)=\min \left\{t: \| x^{t}-x^r\|>\delta, t>r\right\}.
$$

Since the function $f(x)$ is continuous, we choose the number $\delta$ so small that for all sufficiently large $r \in R$, we get $r<k(r)<p$, i.e. .
$$
\begin{aligned}
& f\left(x^{r}\right)\le a<f\left(x^{k(r)}\right) \text {, } \\
& \varlimsup_{r \rightarrow \infty} f\left(x^{k(r)}\right) \ge \lim _{r \rightarrow \infty} f\left(x^{r}\right),
\end{aligned}
$$
which contradicts inequality (4.17) written for $\left\{x^{r}\right\}(r \in R)$. The theorem has been proven.
\begin{corollary}\label{cor:4.1}
With probability 1, there exists a subsequence $\left\{x^{s}\right\}(s \in S)$, for which $\nabla f\left(x^{s}, \alpha_s\right) \rightarrow 0$, $s\rightarrow\infty$.
\end{corollary}

Thus, if $\{ x^s \to x \}$, then by Lemma 2.6, $x \in X^*$

{\it P r o o f.} Suppose that there are $\varepsilon$ and $\overline{k}$, for which
\begin{equation*}
    \| \nabla f(x^k, \alpha_k) \| \geq \varepsilon, \quad k \geq \overline{k}.
\end{equation*}
Then for sufficiently large $k$, inequality (4.15) implies
\begin{equation*}
    f(x^{k+1}, \alpha_{k+1}) \leq f(x^k, \alpha_k) - \varepsilon^2 \rho_k / 2 + \rho_k \left<\nabla f(x^k, \alpha_k) - H(x^k, \alpha_k), \nabla f(x^k, \alpha_k)\right>.
\end{equation*}
Summing these inequalities over $k$ from $0$ to $N$ and tending $N$ to $\infty$, we obtain a contradiction with the boundedness of $\{ f(x^k) \}$, since with probability 1
\begin{equation*}
    \sum_{k=}^{\infty} \rho_k \left<\nabla f(x^k, \alpha_k) - H(x^k, \alpha_k), \nabla f(x^k, \alpha_k)\right> < \infty.
\end{equation*}

Let us formulate an important property of the set $\bar{X}^*$ of limit points of the sequence $\{ x^k \}$ defined by formulas (4.6)-(4.9).
\begin{corollary}
\label{cor:4.2}
{Let $x^1$ and $x^2$ be points from the set $\tilde{X}^*$. Then $f(x^1) = f(x^2)$, i.e. the sequence $\{ f(x^k) \}$ converges with probability 1.}
\end{corollary}
\begin{remark}
The sequence $\{ x^k \}$ can be made artificially bounded, if we consider the iterative process
\begin{equation*}
    x^{k+1} =
    \begin{cases}
        x^k - \rho_k H(x_k, \alpha_k), &f(x^k) \leqslant f(x^0) + C, \\
        x^0,                           &f(x^k)   >  f(x^0) + C, 
    \end{cases}
\end{equation*}
where $x^0$ is the initial approximation, $C > 0$ is some constant. Similar algorithms were studied in [91]. Applying conditions (4.12), (4.13), it is easy to show that with probability 1 the sequence $\{ x^k \}$ goes beyond the set $\{ x | f(x) \leqslant f(x^0) + C \}$ no more than finite number of times.
\end{remark}

    \section*{$\S$ 5. Construction of finite differences \\for the maximum function}
\label{Sec.5}
\setcounter{section}{5}
\setcounter{definition}{0}
\setcounter{equation}{0}
\setcounter{theorem}{0}
\setcounter{lemma}{0}
\setcounter{remark}{0}
\setcounter{corollary}{0}
\setcounter{example}{0}

From the point of view of practical applications, numerical methods that do not require gradient calculations are of great interest for minimizing the maximum function. On the one hand, this fact is explained by the fact that in problems with a large number of variables it is very difficult for the user to calculate generalized gradients, which are required in gradient-type methods. On the other hand, in the problem of minimizing the function
\begin{equation*} 
    f(x) = \max_{1 \le i \le r} f_i(x)
\end{equation*}
he needs to prepare $n \cdot r$ routines for computing the generalized gradients of $f_i(x) (i = 1, \cdots, r)$. Therefore, in problems of large dimensions, this procedure can take a significant amount of time.


The construction of finite difference methods based on the function $f(x)$ is inefficient, since they require the order of $n \cdot r$ to calculate the values of the functions $f_i(x) (i = 1, \ldots, r)$. The methods discussed below use $n + r + 1$ (or $2n + r$) calculations of $f_i(x)$; they are also applicable when $f_i(x)$ are non-differentiable.

Consider the problem:
\begin{equation*}\tag{5.1}
    f(x) = \max_{1 \le i \le r} f_i(x) \to \min_x
\end{equation*}
given that
\begin{equation*}\tag{5.2}
    x \in X
\end{equation*}
where $f_i(x)$ are convex functions, $X$ is a bounded closed convex set in $E_n$. The domain $X$ is such that the operation of projecting onto $X$ is simple.

Finite-difference methods for solving problems (5.1), (5.2) are defined by the relations
\begin{equation*}\tag{5.3}
    x^{k+1} = \pi_X(x^k - H^j (x^k, \alpha_k)),\;\;\; k=0,1,...,
\end{equation*}
where vector $H^j (x^k, \alpha_k)$ is chosen in one of two ways:
\begin{equation*}\tag{5.4}
    H^j (x^k, \alpha_k) = \frac{1}{2\alpha_k} \sum_{i=1}^{n} \left[ f_j(\tilde{x}_1^k, \ldots, x_j^k + \alpha_k, \ldots, \tilde{x}_n^k) - f_j(\tilde{x}_1^k, \ldots, x_i^k - \alpha_k, \ldots, \tilde{x}_n^k) \right] e_i,
\end{equation*}
\begin{equation*}\tag{5.5}
    H^j (x^k, \alpha_k) = \sum_{i=1}^{n} \frac{f_j(\tilde{x}^k + \Delta_k e_t) - f_j(\tilde{x}^k)}{\Delta_k} e_t;
\end{equation*}
$\pi_X(x)$ is the projection operator of point $x$ onto set $X$; $\tilde{x}_i^k (i = 1, \ldots, n)$ are independent random variables uniformly distributed on segments $\left[ x_i^k - \alpha_k, x_i^k + \alpha_k \right] $; index $j$ is such that
\begin{equation*}
     f_j(x^k) = f(x^k) = \max_{1 \leqslant i \leqslant r} f_i(x^k)
\end{equation*}

Thus, in method (5.3)-(5.5), as well as in algorithms of the subgradient type, the final difference is calculated for the function on which the maximum is achieved.
\begin{theorem}\label{th:5.1}
{Let the conditions be met}
\[
     \sum_{k=0}^{\infty} \rho_k = \infty, \;\;\; \sum_{k=0}^{\infty} \rho_k^2 < \infty, \;\;\; \frac{\Delta_k}{\alpha_k} \to 0, \;\;\;\quad \alpha_k \to 0.
\]
{Then, with probability 1, the limit points of the sequence $\{x^k\}$ belong to the set of solutions of problem (5.1), (5.2).}
\end{theorem}

{\it P r o o f.} In $\S$ 4 it is shown that the random vectors $H^j(x^k, \alpha_k)$, which determine the direction of descent in method (5.3), satisfy the condition
\begin{equation*}
    \mathbb{E}H^j(x^k, \alpha_k) | x^k = \nabla f_j(x^k, \alpha_k) + b_k,
\end{equation*}
where $f_j(x, \alpha$ are smoothed functions, $\| b_k \| \leq C\Delta_k / \alpha_k$. From the inequality
\begin{equation*}
     f_j(t, \alpha_k) - f_j(x, \alpha_k) \geq \left<\nabla f_j(x, \alpha_k), y - x\right>
\end{equation*}
it follows that $\nabla f_j(x, \alpha_k)$ converges uniformly to $\partial f_j(x)$ as $\alpha_k \to 0$, where $\partial f_j(x)$ is the subdifferential of the function $f_j( x)$. By virtue of the properties of the maximum function, for convex functions the following relation holds:
\begin{equation*}
    \partial f_j(x) \subset \partial f(x), \quad f_j(x) = f(x).
\end{equation*}

Let $x^* \in X^*, X^*$ be the set of solutions to problem (5.1), (5.2). Applying convergence conditions (4.10)-(4.13) with $W(x) = \displaystyle \min_{x^* \in X^*} \| x - x^* \|$ the assertion of the theorem is obtained directly from the inequalities
\begin{multline*}
     \| x^{k+1} - x^* \|^2 \leqslant \| x^k - \rho_k H^j(x^k, \alpha_k) - x^* \|^2 =\\
     = \| x^k - x^* \|^2 - 2\rho_k \left<\nabla f_j(x^k, \alpha_k) + b_k, x^k - x^*\right> +\\
     + 2\rho_k \left<\nabla f_j(x^k, \alpha_k) + b_k - H^j(x^k, \alpha_k), x^k - x^*\right> + \rho_k^2 \| H^j(x^k, \alpha_k) \|^2
\end{multline*}
due to the fact that with probability 1
\begin{equation*}
     \sum_{k=0}^{\infty} \rho_k \left<\nabla f_j(x^k, \alpha_k) + b_k - H^j(x^k, \alpha_k), x^k - x^*\right> < \infty.
\end{equation*}

\section*{$\S$ 6. Random finite difference directions}
\label{Sec.6}
\setcounter{section}{6}
\setcounter{definition}{0}
\setcounter{equation}{0}
\setcounter{theorem}{0}
\setcounter{lemma}{0}
\setcounter{remark}{0}
\setcounter{corollary}{0}
\setcounter{example}{0}
Determining the vectors $H(x^k, \alpha_k)$ in formulas (4.6), (4.7) requires $2n$ or $n+1$ calculations of the objective function values. If one calculation of the function $f(x)$ takes a significant amount of time, then the difference approximations (4.6), (4.7) may turn out to be unacceptable. In this case, it is advisable to use the ideas of random search algorithms.

Consider the difference approximation
\begin{equation*}
     z(x^k) = \sum_{i=1}^{p} \frac{f(\tilde{x}^k + \Delta_k\mu_j^k) - f(\tilde{x}^k)} {\Delta_k} \mu_j^k.\tag{6.1}
\end{equation*}
Here $\tilde{x}^k$ is a random point uniformly distributed in an $n$-dimensional cube centered at $x^*$ and with edge $2\alpha_k$; the components of the vector $\mu_i^k (i = 1, \ldots, p)$ are independent random variables uniformly distributed on $\left[ -1, 1 \right]$.

By the Lagrange formula, we have
\begin{multline*}
     (f(x^k + \Delta_k \mu_i^k, \alpha_k) - f(x^k, \alpha_k)) / \Delta_k = \left<\nabla f(x^k + \tau_k \Delta_k \mu_i^k, \alpha_k), \mu_i^k\right> =\\
     = \left<\nabla f(x^k, \alpha_k), \mu_i^k\right> + v_k, \quad 0 \leq \tau_k \leq 1,
\end{multline*}
where $f(x, \alpha)$ is a smoothed function, $v_k$ is some vector such that
\begin{equation*}
    \| v^k \| \leq C\Delta_k / \alpha_k.
\end{equation*}

\noindent
Since $\mathbb{E}(\mu_i^k)^2 = 2/3$, the conditional expectation of the vector $z(x^k)$ satisfies the relation
\begin{equation*}
     \mathbb{E}\{z(x^k) | x^k\} = \frac{2}{3}p \nabla f(x^k, \alpha_k) + b_k,
\end{equation*}
where
\begin{equation*}
     \| b_k \| \leqslant C\Delta_k / \alpha_k.
\end{equation*}
Therefore, to minimize locally Lipschitz functions, one can apply the random search method
\begin{equation*}
     x^{k+1} = x^k - \rho_k z(x^k),\tag{6.2}
\end{equation*}
in which the vector $z(x^k)$ is calculated by formula (6.1).
\begin{theorem}\label{th:6.1}
{Let the conditions be met}
\begin{equation*}
    \begin{gathered}
        \sum_{k=0}^{\infty} \rho_{k}=\infty, \quad \sum_{k=0}^{\infty} \rho_{k}^{2}<\infty, \quad \rho_{k} / \alpha_{k} \rightarrow 0, \\
        \Delta_{k} / \alpha_{k} \rightarrow 0, \quad\left|\alpha_{k}-\alpha_{k+1}\right| / \rho_{k} \rightarrow 0, \quad \alpha_{k} \rightarrow 0 .
    \end{gathered}
\end{equation*}
\emph{Then all limit points of the sequence $\left\{x^{k}\right\}$ defined by formulas (6.1), (6.2) belong to the set $X^* = \left\{x ^{*} \mid 0 \in \partial f\left(x^{*}\right)\right\};$ sequence $\left\{f\left(x^{k}\right)\right\}$ converges almost certainly.}
\end{theorem}

{\it P r o o f.} Let us derive an inequality which, together with conditions (4.12), (4.13), leads directly to the proof of Theorem 6.1.

For smoothed functions $f(x, \alpha)$ we have
\begin{multline*}
    f(x^{k+1}, \alpha_k) = f(x^k, \alpha_k) + \left<\nabla f (x^k + \tau_k (x^{k+1} - x^k), \alpha_k) \right. \\
    \left.- \nabla f\left(x^{k}, \alpha_{k}\right)+\nabla f\left(x^{k}, \alpha_{k}\right), x^{k+1}-x^{k}\right> \\
		\leq f(x^k, \alpha_k) + C \frac{\rho_k^2}{\alpha_k}
    -\rho_k \left<\nabla f(x^k, \alpha_k), z(x^k)\right>, \quad 0 \leq \tau_k \leq 1.
\end{multline*}
Then
\begin{multline*}
    f(x^{k+1}, \alpha_{k+1}) \leqslant \\
    \leqslant f(x^k, \alpha_k) - \frac{2}{3} p \rho_k \| \nabla f(x^k, \alpha_k) \|^2 + C\|\alpha_{k}-\alpha_{k+1} \| + C \frac{\rho_k^2}{\alpha_k} + \\
    + C{\rho_k} \frac{\Delta_k}{\alpha_k}+\rho_k \left<\frac{2}{3} p \nabla f(x^k, \alpha_k) + b_k - z(x^k), \nabla f(x^k, \alpha_k)\right>.
\end{multline*}

Taking into account the obtained inequality, the proof of Theorem 6.1 differs slightly from the proof of Theorem 4.1.

Method for minimizing a convex maximum function
\begin{equation*}
    f(x)=\max_{1 \le i \le r} f_i(x), \quad x \in X,
\end{equation*}
is determined by the sequence
\begin{equation*}
    x^{k+1}=\Pi_X \left(x^k - 
		\rho_k \sum_{i=1}^{p} \frac{f_{j} (\tilde{x}^k + \Delta_k \mu_i^k) - f_j \tilde{x}^k)}{\Delta_k} \mu_i^k\right),
\end{equation*}
where
\begin{equation*}
    f_j(x^k) = f(x^k) = \max_{1 \le i \le r} f_i(x^k).
\end{equation*}

\section*{$\S$ 7. Efficiency of the finite difference methods}
\label{Sec.7}
\setcounter{section}{7}
\setcounter{definition}{0}
\setcounter{equation}{0}
\setcounter{theorem}{0}
\setcounter{lemma}{0}
\setcounter{remark}{0}
\setcounter{corollary}{0}
\setcounter{example}{0}

Recently, the issues of the effectiveness of numerical optimization methods have become increasingly important. Due to the significant expansion of applications and the complication of optimization problems, it becomes necessary to have in the arsenal of computing tools, in a certain sense, optimal methods that would provide the solution of all problems of a certain class with the minimum possible laboriousness.

In [73], a theory of complexity of extremal problems of convex programming was constructed.

Consider the problem
\begin{equation}\label{eqn:7.1}
    \begin{gathered}
        f_0(x) \to \min, \\
        f_i(x) \le 0, \quad i=1, \ldots, m, \\
        x \in X,
    \end{gathered}
\end{equation}
where $f_j(x) (j = 0, 1, \ldots, m)$ are convex functions, $X$ is a bounded convex closed set of some finite-dimensional Banach space. Problem (7.1) is equipped with an oracle that calculates the values of the functions $f_j(x) (j = 0, 1, \ldots, m)$ (and the values of their gradients) at a point specified to it.

Numerical methods give an approximate solution of problem (7.1). This solution can be a point $x$ of the set $X$ or an indication that problem (7.1) is inconsistent (the latter fact is denoted by *). The relative error of the point $x \in X$ is taken as the error of the solution:
\begin{equation*}
    v(x) =
    \begin{cases}  
        1, \text { task } (7.1) \text { is compatible, } x = *, \\
         0, \text { task } (7.1) \text { is incompatible, } x = *, \\
        \max \{ \frac{f_0(x) - \bar{f}_0}{r_0}, \frac{f_1(x)}{r_1}, \ldots, \frac{f_m(x)}{r_m} \}, x \neq *,
    \end{cases}
\end{equation*}
where $\bar{f}_0$ is the optimal value of $f_0(x)$. Here $r_j > 0$ $(0 \le j \le m)$ are the normalizing factors, which are chosen, as a rule, such that the relative errors of any point $x \in X$ do not exceed 1.

The method for solving problem (7.1) is a set of rules for generating the next questions to the oracle, the moment of stopping and issuing the result. The complexity of solving a given problem by this method is defined as the number of steps of its work on it. In addition to the laboriousness, the method is characterized by its own error.


The complexity $N(v)$ of this class of problems as a function of $v$ is equal to the minimum possible complexity of the method that solves any problem of the class with an error not exceeding $v$.

It was shown in [73] that the complexity of convex problems admits the estimate
\begin{equation*}
    C_1 \leqslant \frac{N(v)}{1 + n \ln (1 / v)} \leqslant C_2,
\end{equation*}
where $C_1, C_2$ are positive constants. The right inequality is true for all $v<1$, while the left one is true asymptotically with respect to $v \to 0$.

In [56], the center of gravity method (CGM) was proposed, which has a laboriousness of $O(n \ln (1 / v))$ and is the optimal method for $v \to 0$. This method is based on drawing reference hyperplanes through the center of gravity of the region where the minimum is localized. The procedure for reducing the localization area by successive cuts is very effective, but this method cannot be used in practice due to the excessive complexity of the step.

In [141, 135], modified methods of the center of gravity were proposed (otherwise they are also called methods of ellipsoids), in which the initial localization region of the minimum.  a hemisphere, is immersed in an ellipsoid of minimum volume, the center of gravity of which is very easy to find. Then the ellipsoid is turned into a ball by changing the metric and the operation is repeated. The ellipsoid method can be interpreted in another way as a method of stretching space in the direction of a subgradient. This method is already practically implemented; its complexity is $O(n^2 \ln (1/v))$ with the number of $O(n^2)$ arithmetic operations at each step and the memory $O(n^2)$. Thus, the complexity of these methods is not related either to the presence (and number) of constraints on the problem, or to its smoothness and degree of conditionality.

It turns out that as the dimension $n$ of the problem grows, the complexity $N(v)$ depends on the geometric properties of $X$. If $X$ is a parallelepiped, then $N(v) \sim n \ln (1/v)$ for all $v < 1/4$ and all $n$, so that the complexity grows linearly with $n$. If $X$ is an ellipsoid, then as $n$ increases, the function $N(v)$ is bounded above by $O(v^{-2})$ regardless of the dimension, and for a given $v$ it is $O(v^ {-2})$ for all sufficiently large $n$.

Gradient methods were constructed in [73], which were called mirror descent methods for solving convex problems (7.1) on convex bodies $X$ of the type $L_p$-balls (in the space metric $L_p, 1 \leq p < \infty$) . The estimate of the complexity of the corresponding method has the form
\begin{equation*} 
    C(p)(1/v)^{\max (2, p)} \quad \text { for } \quad p > 1; 
		\quad  C(1) (1/v)^2 \ln n \mbox {   for } \quad p = 1 .
\end{equation*}

Consequently, convex problems on ellipsoids can be solved by methods whose complexity does not depend, generally speaking, on the dimension of the problem. In addition, in high-dimensional problems, the complexity of these methods cannot be improved.

If the initial space is Euclidean, then the mirror descent method turns into the usual subgradient method with a special step adjustment. If at the current point $x^{k} (k=0,1, \ldots)$  the constraints $f_i(x) \leq 0 \quad (i = 1, \ldots, m)$ are satisfied with accuracy
up to values of order $v$, then the descent at the point $x^k$ occurs in the direction of the vector $-g_0 (x^k ),\; g_0 (x^k) \in \partial f_0(x^k)$, otherwise, by the antisubgradient of the constraint function, which is not satisfied at the point $x^k$. Mirror descent methods are robust against errors in the original information.

In spaces other than Euclidean, gradient descent occurs in the dual space, after which the mapping to the original space is pierced. For example, if the original space is $L_p^{(n)}$ with the norm $\|x\|_p = (\sum_{i=1}^n |x_i|^p)^{1/p}$, then the main motion occurs in the space $L_q^{(n)}, 1/p + 1/q = 1,$ dual to $L_p^{(n)}.$

Mirror descent methods have been extended to stochastic programming problems, where the objective functions and constraints are not specified analytically and are probabilistic. The source of information about the problem to be solved (the oracle) reports at the right time the values of the functions and their subgradients at the point in question. The oracle may be mistaken, but so that the error does not exceed the specified accuracy of the solution, and its answers may be distorted by random noise.

The problem is subject to the following requirements $\mathrm{A}-\mathrm{B}$, which characterize random noise. Let $R > 0$ and $v_0$ describe the accuracy of the oracle.

{\bf A.} For all $x \in X, 0 \le j \le m$ the functionals
\begin{equation*}
     h_j(y) = \mathbb{E}_\theta \Psi_j(x, \theta) + \left<\mathbb{E}_\theta \xi_j(x, \theta), y - x\right> - v_0 R
\end{equation*}
satisfy the conditions
\begin{equation*} 
    f_j(y) \geq h_j(y), \quad y \in X, \quad f_j(x) \le h_j(x) + v_0 R.
\end{equation*}

Scalar functions $\Psi(x, \theta)$ are interpreted as observations of $f_j(x)$ at point $x$ under oracle noise $\theta,$ $\xi_j(x, \theta)$ are interpreted as observations of subgradients of functions $f_j (x)$.

{\bf B.} For all $j = 0, \ldots, m$ and $x \in X$
\begin{equation*}
     \mathbb{E}_\theta\| \xi_j(x, \theta) \|^r \le (R / (2\rho(X)))^r, \quad r > 1,
\end{equation*}
where $\rho(X)$ is the radius of the set $X$.

{\bf C.} For all $x \in X$
\begin{equation*}
     \mathbb{E}_\theta |\Psi_j(x, \theta) - f_j(x)|^r \le R^r
\end{equation*}

As normalizing factors, choose
\begin{equation*} 
    r_j =
    \begin{cases}
        R, & j = 0, \\
        R + \max [0, \min_{x \in X} f(x)], & j \ge 1.
    \end{cases}
\end{equation*}

Requirement {\bf A} means that a certain systematic error is allowed when calculating the functions and their subgradients. Requirements {\bf B} and {\bf C} limit the intensity of oracle disturbance.


In the Euclidean case, the stochastic mirror descent method for the problem 
$$\min f(x), x \in X$$
turns into the well-known method
\begin{equation*}\tag{7.2}
     x^{k+1} = \pi_X (x^k - \rho_k \xi (x^k, \theta^k))
\end{equation*}
with averaging
\begin{equation*}
     \bar{x}^k = \left( \sum_{s=1}^{k} \rho_s \right)^{-1} \sum_{s=1}^{k} \rho_s x^s .
\end{equation*}
Under the conditions $\sum_{k=0}^{\infty} \rho_k = \infty,$ $\rho_k \to 0$, the trajectory $\bar{x}^k$ converges to the optimum on average over the functional, and an estimate of the rate of convergence is given. If we choose the steps constant according to the formula $\rho_k = 4 (v - v_0) / R$, then the number of iterations that ensure the accuracy $v>v_{0}$ is equal to $N = \left] 1/(4(v - v_0 )^2) \right[ \left( \left] t \right[ \text { minimal integer not less than } t \right)$.

The complexity of these stochastic methods, when a given average error is reached, is the same as for mirror descent methods for deterministic problems. In other words, gradient methods of non-differentiable optimization almost do not react to noise in the supplied information.

Our goal is to show that the stochastic finite differences considered in $\S$ 4 -- 6 satisfy 
conditions A-B. Consider the smoothed function
\begin{equation*} 
    f(x, \alpha) = \frac{1}{(2 \alpha)^n} \int_{-\alpha}^{\alpha} \ldots \int_{-\alpha}^{\alpha} f(x+y) dy_1 \ldots d y_n.
\end{equation*}

Let $f(x)$ be a Lipschitz function with constant $L$. The study is carried out in the space $E_n$. Then
\begin{equation*} 
    | f(x, \alpha) - f(x) | \le \sqrt{n} L \alpha.\tag{7.3}
\end{equation*}

B $\S$ 4 it is shown that for the vector $H(x, \alpha)$ calculated by the formula
\begin{equation*}\tag{7.4}
     H(x, \alpha) = \frac{1}{2\alpha} \sum_{i=1}^{n} \left[ f(\tilde{x}_1, \ldots, x_i + \alpha, \ldots, \tilde{x}_n) - f(\tilde{x}_1, \ldots, x_i - \alpha, \ldots, \tilde{x}_n) \right] e_i,
\end{equation*}
the condition $\mathbb{E} H(x, \alpha) = \nabla f(x, \alpha)$ is satisfied, where $\tilde{x}_i \; (i = 1, \ldots, n)$ are independent random variables, uniformly distributed on segments $\left[ x_i - \alpha, x_i + \alpha \right]$.

From formula (7.4) it follows that
\begin{equation*}\tag{7.5}
     \mathbb{E} \| H(x, \alpha) \|^2 \le (\sqrt{nL})^2.
\end{equation*}

The smoothing parameter $\alpha$ due to (7.3) and the inequality for convex functions
\begin{equation*}
    f(y, \alpha) - f(x, \alpha) \ge \left<\nabla f(x, \alpha), y-x\right>
\end{equation*}

\noindent
can be chosen so that the requirements of A--B are satisfied. Thus, to the accuracy of random disturbance, satisfying the requirements of A--B, subgradient is modeled by the stochastic finite difference (7.4). Consequently, in the method (7.2) instead of $\xi\left(x^{k}, \theta^{k}\right)$ one can use a vector of finite differences $H(x, \alpha)$ (7.4):
\begin{equation} \tag{7.6}
    x^{k+1}=\pi_{X}\left(x^{k}-\rho_{k} H\left(x^{k}, \alpha\right)\right) .
\end{equation}
It follows from inequality (7.5) that the labor intensity of the method (7.6) tuned to the absolute error $v>v_{0}$ is $O\left(\sqrt{n L} /\left(v-v_{0}\right)\right)^{2}$; it is assumed that the computation difficulties of the function $f(x)$ and vector $H(x, \alpha)$ are of the same order. Taking a finite-difference approximation of the gradient, we obtain a random vector
\begin{equation} \tag{7.7}
    H(x, \alpha)=\sum_{i=1}^{n} \frac{f\left(\tilde{x}+\Delta e_{i}\right)-f(\tilde{x})}{\Delta} e_{i} .
\end{equation}
The conditional expectation of the vector (7.7) coincides with $\nabla f(x, \alpha)$ with the interference value $\sqrt{n} L \Delta / \alpha$, that is, again the parameters $\Delta$ and $\alpha$ can be chosen such that conditions A--B are satisfied.

If a function calculation takes a considerable amount of time, the finite-difference approximation discussed in B $\S$ 6 is used:
\begin{equation} \tag{7.8}
    H(x, \alpha)=\sum_{i=1}^{p} \frac{f\left(\tilde{x}+\Delta \mu_{i}\right)-f(\tilde{x})}{\Delta} \mu_{i} .
\end{equation}
Here $\tilde{x}_{i}$ $(i=1, \ldots, n)$ are uniformly distributed on segments $\left[x_{i}-\alpha\right.$, $\left.x_{i}+\alpha\right]$, the vector components $\mu_{i}$ are independent and uniformly distributed on $[-1,1]$ random variables.

From formula (7.8) it follows that
$$
\mathbb{E}\|H(x, \alpha)\|^{2} \leqslant(\sqrt{n p} L)^{2},
$$
where
$$
\mathbb{E} H(x, \alpha)=\frac{2}{3} p \nabla f(x, \alpha)+b,
$$
where $\|b\| \leq p L \Delta n^{3 / 2} / \alpha$. Therefore, the vector of finite differences (7.8) can be used as the descent direction in the method (7.6), but there are more stringent requirements with respect to the displacements of $\Delta$ here.

It is of great practical interest to construct finite-difference methods of minimization of the maximum function
\begin{equation} \tag{7.9}
    f(x)=\max _{1 \leqslant s \leqslant r} f_{s}(x), \quad x \in X,
\end{equation}
since the maximum operation is one of the main sources generating nonsmooth functions.

In $\S$ 5 we showed that the finite-difference method of minimizing function (7.9) is constructed similarly to the subgradient method:
$$
x^{k+1}=\pi_{X}\left(x^{k}-\rho_{k} H_{j}\left(x^{k}, \alpha\right)\right) \text {, }
$$
where the vectors $H_{j}\left(x^{k}, \alpha\right)$, determined by the formulas $(7.4),(7.7),(7.8)$, are computed by the function on which the maximum in formula (7.9) is reached. Thus, in many important cases the calculation of the subgradient can practically be replaced by the calculation of finite differences $(7.4),(7.7),(7.8)$.

One of the main results of the convex analysis is that the solution of (7.1) is equivalent to the minimization of a nonsmooth penalty function
\begin{equation} \tag{7.10}
    \Phi(x)=f_{0}(x)+\sum_{i=1}^{m} r_{i} \max \left(0, f_{i}(x)\right), \quad x \in X,
\end{equation}
where $r_{i}>0$ are penalty coefficients, F $(x)$ is a maximum function of a special kind; therefore, finite-difference methods are used to minimize it:
\begin{equation} \tag{7.11}
    \begin{gathered} 
    x^{k+1}=\pi_{X}\left(x^{k}-\rho_{k}\left(H_{0}\left(x^{k}, \alpha\right)+\sum_{i=1}^{m} r_{i}^{+} H_{i}\left(x^{k}, \alpha\right)\right)\right), \\
    r_{i}^{+}= \begin{cases}r_{i}, & f_{i}\left(x^{k}\right)>0, \\
    0, & f_{i}\left(x^{k}\right) \leqslant 0;\end{cases}
    \end{gathered}
\end{equation}
vectors $H_{i}\left(x^{k}, \alpha\right) (i=0,1, \ldots, m)$ are determined by one of  formulas (7.4), (7.7), (7.8).

Hence, the labor intensity of finite-difference methods for solving convex problems (7.1) is estimated by the value $O(m+1)\left(\sqrt{n} L /\left(v-v_{0}\right)\right)^{2}$. A study of the method (7.11) is given in $\S$ 18.

Thus, the finite-difference methods studied in this monograph can be widely used in solving non-smooth extremal problems.

Let us consider the convergence rate of finite-difference methods. The previous estimates of labor intensity were obtained for finite-difference methods with constant parameters $\rho_{k}$ and $\alpha_{k}$. Consider now the situation where $\rho_{k}, \alpha_{k} \rightarrow 0$ as $k \rightarrow \infty$, and the minimizing function $f(x)$ is strongly convex, i.e. satisfies the inequality
$$
f\left(\varepsilon x+(1-\varepsilon) y\right) \le \varepsilon f(x)+(1-\varepsilon) f(y)-\varepsilon(1-\varepsilon) \lambda\|x-y\|^{2}
$$
for any $x, y, \varepsilon \in[0,1] ; \lambda>0$ is a strong convexity parameter. We investigate the finite-difference minimization method for $f(x)$ on the convex set $X$ defined by the formula
\begin{eqnarray}
x^{k+1}&=&\pi_{X}
\left(
x^{k}-\frac{\rho_{k}}{2 \alpha_{k}} 
\sum_{i=1}^{n}
\left[
f\left(\tilde{x}_{1}^{k},\ldots,x_{i}^{k}+\alpha_{k},\ldots,\tilde{x}_{n}^{k}\right)
 \right.\right.\nonumber\\
&& \left.\left.\;\;\;\;\;\;\;\;\;\;\;\;\;\;\;\;\;\;\;\;\;\;\;\;\;\;\;\;\;\;
-f\left(\tilde{x}_{1}^{k}, \ldots, x_{i}^{k}-\alpha_{k}, \ldots, \tilde{x}_{n}^{k}\right)
\right] e_{i}
\right). \nonumber
\end{eqnarray}

The parameters $\rho_{k}, \alpha_{k}$ are chosen as follows:
$$
\begin{gathered}
\rho_{k}=b k^{-\beta}, \quad b>0, \quad 0<\beta \leqslant 1, \\
\alpha_{k}=m k^{-2 \mu}, \quad m>0, \quad \mu>0 .
\end{gathered}
$$
It is assumed that $x^{*} \in X$, where $x^{*}$ is the point of minimum of $f(x)$. From the results of [70] it follows that if $\beta=1, \mu>1$, $2b \lambda>1$, then the mean square rate and almost sure convergence rate are the highest, and $\mathbb{E}\left\|x^{k}-x^{*}\right\|^{2}$ has order $k^{-1}$.

\newpage
\begin{flushright}
CHAPTER 3
\label{Ch.3}

\textbf{GENERALIZED GRADIENT DESCENT METHODS}
\numberwithin{equation}{section} 
\underline{\hspace{15cm}}

\end{flushright}
\bigskip\bigskip\bigskip\bigskip\bigskip\bigskip

This chapter delves into the fundamentals of minimizing nonconvex nonsmooth functions, specifically focusing on the Generalized Gradient Method and the Local Steepest Descent Method. These methods employ the values and generalized gradients of the minimized functions to converge to the local extrema of the corresponding optimization problems. To establish convergence, generalized Lyapunov functions are utilized. Additionally, the challenges associated with finding the global extremum are discussed.

\section*{$\S$ 8. Convergence conditions for iterative algorithms of nonlinear programming}
\label{Sec.8}
\setcounter{section}{8}
\setcounter{definition}{0}
\setcounter{equation}{0}
\setcounter{theorem}{0}
\setcounter{lemma}{0}
\setcounter{remark}{0}
\setcounter{corollary}{0}

This section delves into the intricacies of proving the convergence of optimization algorithms, as discussed in the book. It introduces a general concept of "iterative algorithm," defines its convergence, and establishes necessary and sufficient conditions, as well as general sufficient conditions for algorithm convergence, all expressed in terms of the generalized Lyapunov function.

This study delves into the complexities of convergence analysis for numerical methods in nonconvex nonsmooth optimization, emphasizing the challenges compared to convex or smooth optimization methods due to the lack of a relaxation process. To address these challenges, a special technique based on the generalized Lyapunov function is employed. This technique, initially introduced in [87] for sufficient convergence conditions in nonlinear programming algorithms, has been refined by various authors. The text presents a modern version of these sufficient conditions, along with necessary and sufficient convergence conditions for comparison. It aims to consolidate and thoroughly examine the common foundation for convergence proofs of all methods discussed in the book. Consequently, specific methods in their respective sections will focus on the unique aspects of their convergence proofs, primarily verifying the general convergence conditions.

\textbf{1. Definitions.} 
\label{Sec.8.1}
Let $M$ be a complete metric space, $\rho(x, y)$ be the distance between points $x, y \in M, M^{n}=M \times \ldots \times M$ be the $n$-th Cartesian degree of $M, 2^{M}$ be the set of subsets of space $M$, $\bar{D}$ be the closure of the set $D \subset M$,
  $$
    \rho(x, D)=\inf \{p(x, y) \mid y \in D\}
  $$
 be the distance from the point $x$ to the set $D$.

\begin{definition}
\label{df:8.1}
{Iteration algorithm} $A=\left(A_{k}, D^{k}\right)_{k \ge l}$ is the set of such sets $D^{k} \subset M^{k+1}$ and mappings $A_{k}: D^{k} \rightarrow$ $ 2^{M}, k=l, l+1, \ldots$ such that for any $\left(x^{0}, \ldots, x^{k}\right) \in D^{k}$ it holds
$$
\left(x^{0}, \ldots, x^{k}\right) \times A_{k}\left(x^{0}, \ldots, x^{k}\right) \subset D^{k+1}, \quad k=l, l+1, \ldots
$$
\end{definition}
\begin{definition}
\label{df:8.2}
Any set of points $\left\{x^{k}\right\}_{k=0}^{s}$ from $M$ such that $\left(x^{0}, \ldots, x^{s}\right) \in D^{s}$, where $s \ge l$, will be called \emph{initial data for the algorithm}
$$
A=\left(A_{k}, D^{k}\right)_{k \ge l} .
$$

Any sequence of $\left\{x^{k}\right\}_{k=0}^{\infty}$ from $M$ such as $\left(x^{0}, \ldots, x^{k}\right) \in D^{k}$ for all $k \ge s$, where $s \ge l$, will be called \emph{sequence of initial data} for $A$.
\end{definition}
\begin{definition}
\label{df:8.3}
The sequence of points $\left\{x^{k}\right\}_{k=0}^{\infty}$ is called \emph{generated algorithm}
$$
A=\left(A_{k}, D^{k}\right)_{k \ge l}
$$
with initial data $\left(x^{0}, \ldots, x^{s}\right)$ if $\left(x^{0}, \ldots, x^{s}\right) \in D^{s}$ and
$$
x^{k+1} \in A_{k}\left(x^{0}, \ldots, x^{k}\right), \quad k=s, s+1, \ldots, \quad s \ge l .
$$
\end{definition}

Obviously, any sequence generated by the algorithm $\boldsymbol{A}$ is also a sequence of initial data for $A$.
\begin{definition}
\label{df:8.4}
{Algorithm definition area} $A=\left(A_{k}\right.$, $\left.D^{k}\right)_{k \ge l}$ is the set of
$$
D_{A}=\bigcup_{k=l}^{\infty} \bigcup_{i=0}^{k}\left\{x^{i} \mid\left(x^{0}, \ldots, x^{k}\right) \in D^{k}\right\}
$$
\end{definition}

The above definitions imply that an iterative algorithm produces infinite sequences of points, i.e., it works infinitely long. Finite algorithms can be formally considered as infinitely long if we assume that, starting from some point, they produce the same point.
\begin{definition}
\label{df:8.5}
A sequence of points $\left\{x^{k}\right\}_{k=0}^{\infty}$ is called  \emph{convergent to the set $D \subset M$ by distance} if
$$
\lim _{k \rightarrow \infty} \rho\left(x^{k}, D\right)=0 \text {. }
$$
\end{definition}
\begin{definition}
\label{df:8.6}
An algorithm $A$ is called \emph{convergent $\kappa$ to the set $X^{*} \subset M$ by distance} if any sequence generated by it converges to $X^{*}$ by distance.
\end{definition}

Often somewhat different definitions of the convergence of sequences and algorithms to a set are used.
\begin{definition}
\label{df:8.7}
The sequence $\left\{x^{k}\right\} \subset M$ is called compact if it belongs to some compact set. 
\end{definition}

An infinite sequence is compact if and only if any of its subsequences has limit points.

\begin{definition}
\label{df:8.8}
A compact sequence of points $\left\{x^{k}\right\}_{k=0}^{\infty}$ is called \emph{convergent to the set $X^{*} \subset M$ by limit points} if all its limit points belong to $X^{*}$.
\end{definition}
\begin{definition}
\label{df:8.9}
Algorithm $A$ is called \emph{compact} if for any initial data it generates compact sequences.
\end{definition}
\begin{definition}
\label{df:8.10}
A compact algorithm $A$ is called \emph{convergent to the set $X^{*} \subset M$ by limit points} if any sequence generated by it converges to $X^{*}$ by limit points.
\end{definition}

For a compact algorithm $A$ and a closed set $X^{*}$ the convergence of $A$ to $X^{*}$ by distance is equivalent to the convergence by limit points.
\begin{remark}
\label{rem:8.1}
Let the algorithm $A$, convergent to $X^{*}$, generate sequences $\left\{x^{k}\right\}$ in which, starting from some point, the same point is repeated. This assumption corresponds, as it were, to the finiteness of the algorithm $A$. Then it follows from the definitions of convergence that
$$
x=\lim _{k \rightarrow \infty} x^{k} \in \bar{X}^{*},
$$
in full agreement with the intuitive understanding of the convergence of finite algorithms.
\end{remark}

\textbf{2. Necessary and sufficient conditions for convergence of algorithms.} 
\label{Sec.8.2}
In nonconvex optimization, the convergence of iterative algorithms is often established by contrary reasoning. In a general form they can be given the following form.
\begin{theorem}
\label{th:8.1}
{For the algorithm $A=\left(A_{k}, D^{k}\right)_{k \ge l}$ to converge to the set $X^{*}$ (by distance), it is necessary and sufficient that there exists a continuous function $W(x)\left(x \in \bar{D}_{A}\right)$ such that for any sequence $\left\{x^{k}\right\}$ generated by $A$, the conditions are satisfied: }

\begin{enumerate}
  \item[1)] \emph{the numerical sequence $\left\{W\left(x^{k}\right)\right\}_{k=0}^{\infty}$ has a limit (or instead some condition equivalent to convergence of $\left.\left\{W\left(x^{k}\right)\right\}_{k=0}^{\infty}\right)$ is satisfied};

  \item[2)] \emph{with the formal assumption that the total $\left\{x^{k}\right\}$ does not converge to $X^{*}$ (by distance), it follows that 
  $\left\{W\left(x^{k}\right)\right\}$ has no limit (or some sufficient divergence condition  for $\left\{W\left(x^{k}\right)\right\}$ is satisfied).}
\end{enumerate}
\end{theorem}
{\it P r o o f.} Let us first prove sufficiency. For all sequences $\left\{x^{k}\right\}$ generated by the algorithm, the numerical sequence $ \left \{ W \left(x^{k} \right) \right \} $, according to condition 1), has a limit. Subdivide the whole set of sequences $\left\{x^{k}\right\}$ into subsets of convergent and nonconvergent to $X^{*}$. On sequences $\left\{x^{k}\right\}$ that do not converge to $X^{*}$, the numerical sequence $\left\{W\left(x^{k}\right)\right\}$ by condition 2) has no limit, and hence, the set of such sequences is empty. 

Let us prove the necessity. Let the algorithm $A$ converge to $X^{*}$. Take $ W(x)= \rho \left (x, X^{*} \right) $; the function $\rho \left (x, X^{*} \right )$ is continuous [46]. For all $ \left \{x^{k} \right \}$ generated by the algorithm $A$, by definition of convergence of $A$, we have
$$
\lim _{k \rightarrow \infty} \rho\left(x^{k}, X^{*}\right)=\lim _{k \rightarrow \infty} W\left(x^{k}\right)=0 \text {, }
$$
i.e. condition 1) of the theorem is satisfied. Condition 2) of the theorem refers to sequences $\left\{x^{k}\right\}$ that do not converge to the set $X^{*}$, whose set is empty. Therefore, the premise of condition 2) is false, and from the viewpoint of formal logic it is true (irrespective of the corollary). The theorem is proved.

Theorem 8.1 does not explicitly distinguish between minimization and maximization problems. It can be given a form oriented to minimization problems. Let us first prove one necessary and sufficient condition for convergence of numerical sequences.
\begin{lemma}
\label{lem:8.1}
{For a numerical sequence $\left\{W^{k}\right\}_{k=0}^{\infty}$ to have a finite limit, it is necessary and sufficient that:}
\begin{enumerate}
  \item[1)] \emph{$\left\{W^{k}\right\}_{k=0}^{\infty}$ was bounded;}
  \item[2)] \emph{if $\lim _{s \rightarrow \infty} W^{k_{s}}=W^{\prime}$, then for any $\varepsilon>0$ there is such a number $K(\varepsilon)$, that for all $k \ge K(\varepsilon)$ there be $W^{k} \le W^{\prime}+\varepsilon$.}
\end{enumerate}
\end{lemma}
{\it P r o o f.} Necessity follows directly from the definition of the limit of a numerical sequence. Let us prove sufficiency. Assume the contrary:
$$
\underline{W}=\underline{\lim} _{k \rightarrow \infty} W^{k}<\varlimsup_{k \rightarrow \infty} W^{k}=\overline{W},
$$
where $\underline{\lim}$ and $\varlimsup$ designate the lowest and the upperest limits of a sequence, respectively.
Choose the numbers $c$ and $d$ such that $\underline{W}<c<d<\overline{W}$. The sequence $\left\{W^{k}\right\}$ infinitely many times crosses the segment $[c, d]$ from $c$ to $d ;$ so there exist indices $k_{s}$ and $m_{s} (s=0,1, \ldots)$ such that $k_{s}<m_{s}$ and $W^{k_{s}} \le c<W^{k}<d \le W^{m_{s}}$, when $k \in\left(k_{s}, m_{s}\right)$ or $W^{k_{s}} \le c<d \le W^{k_{s}+1}$. Without losing generality, we can assume that
$$
\lim _{s \rightarrow \infty} W^{k_{s}}=W^{\prime} \le c.
$$
Take $\varepsilon<d-c$. According to condition 2) of the lemma, for all sufficiently large $k$ we have $W^{k} \le W^{\prime}+\varepsilon<d$, which contradicts the earlier construction. The lemma is proved.
\begin{theorem}
\label{th:8.2}
{For the algorithm $A=\left(A_{k}, D^{k}\right)_{k \ge l}$ to converge to a set $X^{*}$ (by distance), it is necessary and sufficient that a continuous function $W: \overline{D}_{A} \rightarrow E_{1}$ exists such that for any sequence $\left\{x^{k}\right\}$ generated by the algorithm, the conditions hold:}
\begin{enumerate}
  \item[0)] \emph{The numerical sequence $\left\{W\left(x^{k}\right)\right\}$ is bounded;}
  \item[1)] \emph{If some subsequence $\left\{x^{k_{s}}\right\}_{s=0}^{\infty}$ is not convergent to $X^{*}$ (by distance) and $\lim _{s \rightarrow \infty} W\left(x^{k_{s}}\right)=W^{\prime}$, then there is another subsequence $\left\{x^{l_{s}}\right\}_{s=0}^{\infty}$, such that}
\end{enumerate}
$$
W^{\prime}=\lim _{s \rightarrow \infty} W\left(x^{k_{s}}\right)>\varlimsup _{s \rightarrow \infty} W\left(x^{l_{s}}\right);
$$
\begin{enumerate}
  \item[2)] \emph{If $\lim _{s \rightarrow \infty} W\left(x^{k_{s}}\right)=W^{\prime}$, so for all $\varepsilon>0$ there is such a number $K(\varepsilon)$, that for all $k \ge K(\varepsilon)$ the  inequality is fulfilled:}
\end{enumerate}
$$
W\left(x^{k}\right) \le W^{\prime}+\varepsilon.
$$
\end{theorem}

{\it P r o o f.} Conditions 0), 2) of Theorem 8.1 are equivalent to the convergence of the sequence $\left\{W\left(x^{k}\right)\right\}$, and condition 1) is sufficient for the divergence of $\left\{W\left(x^{k}\right)\right\}$, so Theorem 8.2 follows from Theorem 8.1. The theorem is proved.

Now regroup the conditions of Theorem 8.2. Let us break condition 2) into two: the first one will refer to the subsequences $\left\{x^{k_{s}}\right\}$ that do not converge to $X^{*}$ (we then combine it with condition 1) of the theorem), and the second one will refer only to the subsequences $\left\{x^{k_{s}}\right\}$ that do converge to $X^{*}$. As a result, we obtain the following necessary and sufficient convergence conditions.
\begin{theorem}
\label{th:8.3}
{ For the algorithm $A=\left(A_{k}, D^{k}\right)_{k \ge l}$ to converge to $ X^{*}$ (by distance), it is necessary and sufficient that there exists a continuous function $W: \bar{D}_{A} \rightarrow E_{1}$ such that for any sequence $\left\{x^{k}\right\}$ generated by $A$, the following conditions are satisfied:}
\begin{enumerate}
  \item[0)] \emph{The numerical sequence $\left\{W\left(x^{k}\right)\right\}$ is bounded;}
  \item[1)] \emph{If a subsequence $\left\{x^{k_{s}}\right\}$ does not converge to $X^{*}$ (by distance) and \\$\lim _{s \rightarrow \infty} W(x^{k_{s}}) = W^{\prime}$, then for any $\varepsilon>0$ there be such indices $l_{s} \ge k_{s}$ that $\mathrm{W}\left(x^{k}\right) \le W^{\prime}+\varepsilon$ as $k \in\left[k_{s}, l_{s}\right), $ and}
\end{enumerate}
$$
W^{\prime}=\lim _{s \rightarrow \infty} W\left(x^{k_{s}}\right)>\varlimsup_{s \rightarrow \infty} W\left(x^{l_{s}}\right);
$$
\begin{enumerate}
  \item[2)] \emph{If a subsequence $\left\{x^{k_{s}}\right\}$ converges to $X^{*}$ (by distance) and $\lim _{s \rightarrow \infty}W\left(x^{k_{s}}\right)=W^{\prime}$, then for any $\varepsilon>0$ there is such a number $K(\varepsilon)$ that for all $k \ge K(\varepsilon)$ it takes place}
\end{enumerate}
$$
W\left(x^{k}\right) \le W^{\prime}+\varepsilon.
$$
\end{theorem}
{\it P r o o f.} Let us show that Theorems 8.2 and 8.3 are equivalent in the sense that if the function $W$ and the sequence $\left\{x^{k}\right\}$ generated by the algorithm satisfy conditions 0) - 2) of Theorem 8.2, then they also satisfy the conditions of Theorem 8.3, and vice versa. 

For the sequence $\left\{x^{k}\right\}$ and its subsequence $\left\{x^{k_{s}}\right\}$ such that
$$
\lim _{s \rightarrow \infty} W\left(x^{k_{s}}\right)=W^{\prime},
$$
let us introduce finite or infinite numbers
$$
m_{{s}}(\varepsilon)=\sup \left\{m:| W\left(x^{k}\right)-W^{\prime} \mid \le \varepsilon \text { for } k \in\left[k_{s}, m\right)\right\},
$$
if $\left|W\left(x^{k_{s}}\right)-W^{\prime}\right| \le \varepsilon$, and $m_{\mathrm{s}}(\varepsilon)=k_{\mathrm{s}}$ otherwise.

Let for $\left\{x^{k}\right\}$ and some function $W$  conditions 0) - 2) of Theorem 8.2 are satisfied. Obviously, then for the same $\left\{x^{k}\right\}$ and $W$ conditions 0), 2) of Theorem 8.3 are fulfilled. Let us show that condition 1) of Theorem 8.3 is also satisfied. Let there exist a subsequence $\left\{x^{k_{s}}\right\}$ that does not converge to $ X^{*}$, and
$$
\lim _{s \rightarrow \infty} W\left(x^{k_{s}}\right)=W^{\prime}.
$$
By condition 1) of Theorem 8.2, there exists $\varepsilon^{\prime}>0$ for which $m_{s}\left(\varepsilon^{\prime}\right)<\infty$. Let us show that for any $\varepsilon \in\left(0, \varepsilon^{\prime}\right]$ the inequality holds:
$$
W^{\prime}=\lim _{s \rightarrow \infty} W\left(x^{k_{s}}\right)>\varlimsup_{s \rightarrow \infty} W\left(x^{m_{s}(\varepsilon)}\right);
$$
from here the condition 1) of Theorem 8.3 will follow. Suppose the contrary; then for some $\varepsilon \le \varepsilon^{\prime}$ there are such arbitrarily large $s$ that
$$
W^{\prime}+\varepsilon \le W\left(x^{m_{s}(\varepsilon)}\right).
$$
This contradicts condition 2) of Theorem 8.2.

Now let $\left\{x^{k}\right\}$ be generated by the algorithm $A$ and some function $W(x)$ satisfy conditions 0) - 2) of Theorem 8.3. Obviously, then for the same $\left\{x^{k}\right\}$ and $W$ conditions 0), 1) of Theorem 8.2 are satisfied. Let us check the validity of condition 2) of Theorem 8.2. Let
$$
\lim _{s \rightarrow \infty} W\left(x^{k_{s}}\right)=W^{\prime}.
$$
This condition is obviously satisfied on the subsequences $ \left \{x^{k_{s}} \right \}$ that converge to $X^{*}$. Now assume that a subsequence $ \left \{x^{k_{s}} \right \}$ does not converge to $X^{*}$, and suppose that in this case condition 2) of Theorem 8.2 is not satisfied. Then for some $\varepsilon>0$ and for the infinite number of $s$ there exist such indices $m_{s} \ge k_{s}$ that
$$
W\left(x^{m_{s}}\right)>W^{\prime}+\varepsilon.
$$
Without loss of generality, we may assume that this inequality is satisfied for all $s$.

By virtue of condition 1) of Theorem 8.3, there exist indices $l_{s}$ for which $k_{s} \le l_{s} \le m_{s}$ and $W\left(x^{l_{s}}\right)<W^{\prime}$. Numerical sequences $\left\{W\left(x^{l_{s}}\right), \ldots\right.$ $\left.\ldots, W\left(x^{m_{s}}\right)\right\}$ $(s=0,1, \ldots)$ intersect the segment $\left[W^{\prime}, W^{\prime}+\varepsilon\right]$ from $W^{\prime}$ to $W^{\prime}+\varepsilon ;$
so there are indexes $p_s$ and $q_s$ $(l_s\le p_s<q_s\le m_s)$ such that
\begin{equation}
\label{eqn:8.1}
W(x^{p_s})\le W^\prime<W(x^k)<W^\prime+\varepsilon\le W(x^{q_s}),\;\;\;
k\in(p_s,q_s),\tag{8.1}
\end{equation}
or
\begin{equation}
W(x^{p_s})\le W^\prime<W^\prime+\varepsilon\le W(x^{q_s}),\;\;\;
q_s=p_s+1.\nonumber
\end{equation}
Without loss of generality, we can consider that
\[
\lim_{s\rightarrow\infty}W(x^{p_s})=W^{\prime\prime}\le W^\prime.
\]
Now, if $\{x^{p_s}\}$ converges to $X^*$, then inequalities (8.1) contradict condition 2) of Theorem 8.3, and if $\{x^{p_s}\}$ does not converge to $X^*$, then condition 1) of Theorem 8.3 applied to $\{x^{p_s}\}$ also contradicts the construction (8.1). The theorem is proved.

{\bf 3. Sufficient conditions for convergence.} 
\label{Sec.8.3}
First, we slightly weaken condition 2) of Theorem 8.3 and thus obtain the necessary conditions for convergence of algorithms. With the weakened (necessary) conditions there will be no convergence of algorithms in the sense of Definition 8.6, but it turns out that in this case there is a weakened convergence,  i.e., in some sense convergence of algorithms on a function (Theorem 8.4). Then we slightly strengthen the conditions of Theorem 8.3 and thus obtain sufficient conditions for convergence of the algorithms (Theorem 8.5).

Let us introduce the notation for the lowest and the upperest limits of a numerical sequence 
$\left\{W(x^k)\right\}$:
\begin{equation}
\label{eqn:8.2}
\underline{W}=\underset{k\rightarrow\infty}{\underline{\lim}} W(x^k),\;\;\;\;
\overline{W}=\underset{k\rightarrow\infty}{\overline{\lim}} W(x^k). \tag{8.2}
\end{equation}
By $W^*$, denote the set of limit points of the sequence $\{W(x^k)\}$, where $\{x^k\}$ is  an infinite sequences of points convergent to $X^*$. If $X^*$ is a compact in $M$, then
\[
W^*=\{W(x)|\,x\in X^*  \}.
\]
\begin{theorem}
\label{th:8.4}
For the algorithm $A = (A_k,D^k)_{k\geq l}$ to converge to a set $X^*$ (by distance), it is necessary that there exists such a continuous function 
$W: \bar{D}_A\rightarrow E_1$ that for any sequence $\{x^k\}$ generated by $A$, the conditions are fulfilled:

0) The numerical sequence $\{W(x^k)\}$ is bounded;

1) If the sequence $\{x^{k_s}\}$ does not converge to $X^*$ (in distance) and
$\underset{s\rightarrow\infty}{\lim}W(x^{k_s})=W^\prime$, then for any
$\varepsilon>0$ there exist such indexes $l_s\ge k_s$ that 
$W(x^k)\le W^\prime +\varepsilon$ for $k\in[k_s,l_s)$ and the inequality holds
true:
\[
W^\prime=\underset{s\rightarrow\infty}{\lim}W(x^{k_s})>
\underset{s\rightarrow\infty}{\overline{\lim}}W(x^{l_s});
\]

   2) If the subsequence $\left\{x^{k_{s}}\right\}$ converges to $X^{*}$ (in distance) and $ \lim_{s \rightarrow \infty} W\left( x^{k_{s}}\right)=W^{\prime}, then$
$$
\varlimsup_{s \rightarrow \infty}\left[W\left(x^{k_{s}+1}\right)-W\left(x^{k_{s}}\right)\right] \le 0 .
$$

If there exists a continuous function $W$ satisfying conditions 0)-2), then algorithm $A$ converges in function, i.e., for any sequence $\left\{x^{k}\right\}$ generated by $A$, the segment $[\underline{W}, \overline{W}]$ is embedded in the set $W^{* }$, and all subsequences $\left\{x^{k_{s}}\right\}$ such that $\lim _{s \rightarrow \infty} W\left(x^{k_{s}} \right)=\underline{W}$, converge to $X^{*}$.
\end{theorem}
{\it P r o o f.} Condition 2) of Theorem 8.4 follows from condition 2) of Theorem 8.3; therefore the conditions of Theorem 8.4 are necessary for the convergence of algorithm $A$. Let us prove their sufficiency for the convergence of the algorithm in terms of function values.

Let's show that if
$$
\lim _{s \rightarrow \infty} W\left(x^{k_{s}}\right)=\underline{W},
$$
then $\left\{x^{k_{s}}\right\}$ converges to $X^{*}$. Let's assume the opposite. Then condition 2) of the theorem is applicable to $\left\{x^{k_{s}}\right\}$, which states that there exists such a subsequence $\left\{x^{l_{s}}\right\} $ that
$$
\varlimsup_{s \rightarrow \infty} W\left(x^{l_{s}}\right) \le \underline{W} .
$$

This contradicts the definition of $\underline{W}$ as the lowest limit of the sequence $\left\{W\left(x^{k}\right)\right\}$. The required assertion is proved. At the same time, we showed that $\underline{W} \in W^{*}$.

Let us now prove that $[\underline{W}, \overline{W}) \subset W^{*}$. Let's assume the opposite. Then there are numbers $c \in[\underline{W}, \overline{W})$ and $c \notin W^{*}$. Since $\underline{W} \in W^{*}$, then $c>\underline{W}$. We choose a number $d$ so that $\underline{W}<c<d<\overline{W}$; then the segment $[c, d]$ is embedded in $(\underline{W}, \overline{W})$ and $c \notin W^{*}$. The sequence $\left\{W\left(x^{k}\right)\right\}$ intersects the segment $[c, d]$ from $c$ to $d$ an infinite number of times; so there are parts $\left\{x^{k_{s}}, \ldots, x^{m_{s}}\right\} (s=0,1, \ldots)$ of sequences $\left\{x^{k}\right\}$ such that either
$$
W\left(x^{k_{s}}\right) \le c<W\left(x^{k}\right)<d \le W\left(x^{m_{s}}\right) , \quad k \in\left(k_{s}, m_{s}\right)
$$
or
$$
W\left(x^{k_{s}}\right) \le c<d \le W\left(x^{k_{s}+1}\right) .
$$

Since $\left\{W\left(x^{k}\right)\right\}$ is bounded, we can assume without loss of generality that
$$
\lim _{s \rightarrow \infty} W\left(x^{k_{s}}\right)=W^{\prime} \le c
$$

Let us show that $\left\{x^{k_{s}}\right\}$ does not converge to $X^{*}$. Let's assume the opposite. Then, by virtue of (8.3) and condition 2) of the theorem, we have
$$
\lim _{s \rightarrow \infty} W\left(x^{k_{s}}\right)=\lim _{s \rightarrow \infty} W\left(x^{k_{s}+1} \right)=c,
$$
and since, by assumption, $x^{k_{s}} \rightarrow X^{*}$, then $c \in W^{*}$, which contradicts to the choice of $c$. Hence $\left\{x^{k_{s}}\right\}$ does not converge to $X^{*}$.

To the subsequence $\left\{x^{k_s}\right\}$ condition 1) of the theorem is applicable. Take $\varepsilon<d-c$. According to condition 1) of the theorem, there are indices $l_{s} \ge k_{s}$ such that $W\left(x^{k}\right) \le W^{\prime}+\varepsilon<d$ for $k \in\left[k_{s}, l_{s}\right)$, and we have
$$
W^{\prime}=\lim _{s \rightarrow \infty} W\left(x^{k_{s}}\right)>\varlimsup_{s \rightarrow \infty} W\left(x^{l_ {s}}\right)\text{. }
$$
But this means that $l_{s} \le m_{s}$, and the minimizing condition 1) of the theorem contradicts the above construction (8.3). So $[\underline{W}, \overline{W}) \subset W^{*}$.

It remains to note that the set $W^{*}$ is closed; so $[\underline{W}, \overline{W}) \subset W^{*}$ implies $[\underline{W}, \bar{W}] \subset W^{*}$. The theorem has been proven.

Let us now somewhat strengthen the necessary conditions of Theorem 8.4, namely, add one additional requirement to them. In this case, it turns out that the strengthened conditions will be sufficient for the convergence of the algorithms.

\begin{theorem}
\label{th:8.5}
For the algorithm $A=\left(A_{k}, D^{k}\right)_{k \geq l}$ to converge to a set $X^{*}$ (in distance), it suffices to have such a continuous function $ W: \overline{D}_{A} \rightarrow E_{1}$ that for any sequence $\left\{x^{k}\right\}$ generated by $A$ conditions 1), 2 ) of Theorem 8.3 are satisfied and the additional condition takes place:

3) the set $W^{*}$ does not contain segments (intervals).

In this case $\left\{x^{k}\right\}$ converges to $ X^{*}$ (in distance), and the numerical sequence $\left\{W\left(x^{k}\right)\right\}$ has a limit.
\end{theorem}
{\it P r o o f.} By Theorem 8.3 $[\underline{W}, \overline{W}] \subset W^{*}$. But according to the condition of the theorem, the linear set $W^{*}$ does not contain segments; so $\left\{W\left(x^{k}\right)\right\}$ has a limit
$$
\lim _{k \rightarrow \infty} W\left(x^{k}\right)=W^{\prime}.
$$

Assume that there exists a subsequence $\left\{x^{k_s}\right\}$ that does not converge to $X^{*}$. It's obvious that
$$
\lim _{s \rightarrow \infty} W\left(x^{k_{s}}\right)=W^{\prime}.
$$
To the subsequence $\left\{x^{k_{s}}\right\}$ now the condition 1) of Theorem 8.3 is applicable, which, in particular, means the divergence of $\left\{W\left(x^{k}\right)\right\}$. The resulting contradiction proves the theorem. Note that here we have actually verified the conditions of the general Theorem 8.1.

The necessary conditions of Theorem 8.3 and the sufficient conditions of Theorem 8.4 give a range containing necessary and sufficient conditions for the convergence of algorithms.

Let us specify Theorems $8.4,8.5$ for compact algorithms.
\begin{theorem}
\label{th:8.6}
For the convergence of a compact algorithm $A=\left(A_{k}, D^{k}\right)_{k \geq l} $ to a compact set $X^{*}$ (by limit points), there necessary the existence of a continuous function $ W: \overline{D}_{A} \rightarrow E_{1}$ such that for any sequence $\left\{x^{k}\right\}$ generated by $A$, the conditions are met:
\begin{enumerate}
   \item If $\left\{x^{k}\right\} \rightarrow x \notin X^{*}$, then for any $\varepsilon>0$ there are indices $l_{s} \ge k_{s }$ such that
\end{enumerate}
$$
\begin{gathered}
W\left(x^{k}\right) \le W(x)+\varepsilon \quad for \quad k \in\left[k_{s}, l_{s}\right), \\
W(x)=\lim _{s \rightarrow \infty} W\left(x^{k_{s}}\right)>\varlimsup_{s \rightarrow \infty} W\left(x^{l_{s }}\right) ;
\end{gathered}
$$
\begin{enumerate}
   \setcounter{enumi}{1}
   \item If $\left\{x^{k}\right\} \rightarrow x \in X^{*}$, then
\end{enumerate}
$$
\varlimsup_{s \rightarrow \infty}\left[W\left(x^{k_{s}+1}\right)-W\left(x^{k_{s}}\right)\right] \le 0 .
$$

If there exists a function $W$ satisfying conditions 1), 2), then the compact algorithm $A$ converges to $X^{*}$ with respect to the function values, i.e. the segment $[\underline{W}, \overline{W}]$ is embedded in the set $W^{*}=\left\{W(x) \mid x \in X^{*}\right\} $, and limit points $\left\{x^{k}\right\}$ minimal with respect to $W$, belong to $X^{*}$.
\end{theorem}
\begin{theorem}
\label{th:8.7}
For the convergence of a  compact algorithm $A=\left(A_{k}, D^{k}\right)_{k \geq l} $ to a compact set $X^{*}$ (by limit points), it suffices to have a continuous function $ W: \overline{D}_{A} \rightarrow E_{1}$ such that for any sequence $\left\{x^{k}\right\}$ generated by $A$, conditions 1), 2) Theorems 8.6 are fulfilled and the condition holds:
\begin{enumerate}
   \setcounter{enumi}{2}
   \item the set $W^{*}=\left\{W(x) \mid x \in X^{*}\right\}$ does not contain intervals.
\end{enumerate}

In this case, all limit points of the sequence $\left\{x^{k}\right\}$ belong to $X^{*}$ and the numerical sequence $\left\{W\left(x^{k}\right)\right\}$ has a limit.
\end{theorem}

Theorems 8.6, 8.7 give very general conditions for the convergence of iterative (compact) algorithms. Now let us narrow down the class of algorithms under consideration. Let us formulate stronger (but less general) sufficient conditions for the convergence of iterative algorithms. In conditions 1), 2) of Theorems 8.6, 8.7, the difference between the values of the function $W$ at two points acted as a measure of the proximity of these points. In the following theorem, we take the distance between points as such a measure.
\begin{theorem}
\label{th:8.8}
Let for a compact algorithm $A=\left(A_{k}, D^{k}\right)_{k\geq l}$ there is a continuous function $W: \overline{D}_{A} \rightarrow E_{1 }$ and a set $X^{*}$ such that for any sequence $\left\{x^{k}\right\}$ generated by $A$, the following condition is satisfied:
\begin{enumerate}
   \item If $\left\{x^{k_{s}}\right\} \rightarrow x \notin X^{*}$, then for any $\delta>0$ there are such indicesl $l_{\mathrm {s}} \ge k_{s}$ that $\rho\left(x, x^{k}\right) \le \delta \quad for \quad$ $k\in\left[k_{s}, l_{s}\right)$ and
\end{enumerate}
$$
W(x)=\lim _{s \rightarrow \infty} W\left(x^{k_{s}}\right)>\varlimsup_{s \rightarrow \infty} W\left(x^{l_s}\right);
$$
then $X^{*}$ is non-empty, and the limit points of the sequence $\left\{x^{k}\right\}$ minimal in the values of $W$, belong to $X^{*}$.

Let, in addition, the following condition be satisfied:

2) If $\left\{x^{k_{s}}\right\} \rightarrow x \in X^{*},$ then
$$
\lim _{s \rightarrow \infty} \rho\left(x^{k_{s}+1}, x^{k_{s}}\right)=0.
$$
Then the algorithm converges to $X^{*}$ with respect to the function values, i.e., the minimal (in the values of $W$) limit points of the sequence $\left\{x^{k}\right\}$  belong to $X^{*}$, and the semi-interval $[\varliminf_{k \rightarrow \infty} W\left(x^{k}\right), \varlimsup_{k \rightarrow \infty} W\left(x^{k}\right) ) $ is embedded in the set $W^{*}=\left\{W(x) \mid x \in X^{*}\right\}$.

If in addition to the previous conditions
\begin{enumerate}
   \setcounter{enumi}{2}
   \item The set ${W}^{*} \subset E_{1}$ does not contain intervals, then the algorithm converges to $ X^{*}$, i.e., all limit points $\left\{x^{k}\right\}$ belong to a connected subset of $X^{*}$, the numerical sequence $\{W(x^{k})\}$ has a limit, and
\end{enumerate}
$$
\lim _{k \rightarrow \infty} \rho\left(x^{k+1}, x^{k}\right)=0.
$$
\end{theorem}

Algorithms, satisfying the conditions of Theorem 8.8, can be called local, and the function $W(x)$ can naturally be called the generalized Lyapunov function for the algorithm $A$.

In the formulation of the theorem, condition 1), roughly speaking, means that if the algorithm starts from the point $x \notin X^{*}$, then it finds values of the function $W$ smaller than $W(x)$ in a small neighborhood of $x$. Condition 2)  is a relaxed requirement of
$$
\lim _{k \rightarrow \infty} \rho\left(x^{k+1}, x^{k}\right)=0,
$$
which is satisfied as a consequence of all conditions of Theorem 8.8.

The theorem asserts the convergence of the algorithm $A$ to connected subsets of the set $X^{*}$, which in practice correspond to local extrema (more precisely, to stationary sets) of the function $W(x)$.

Theorem 8.8 is the main tool for proving convergence of iterative nonconvex optimization algorithms in this book. If the conditions of Theorem 8.8 are satisfied, then, due to the continuity of $W(x)$, the conditions of Theorem 8.7 are also satisfied. However, the sufficient conditions of Theorem 8.8 imply a number of additional results on the convergence of algorithms.

{\it P r o o f of Theorem 8.8.} Let us prove the first assertion of the theorem. The set of limit points of the sequence $\left\{x^{k}\right\}$ generated by the algorithm is obviously closed. Since $W(x)$ is continuous, there exist limit points $\left\{x^{k}\right\}$ that are minimal with respect to values of $W$. Now suppose the contrary: there exists a minimal limit point $x=\lim _{s \rightarrow \infty} x^{k_{s}}$ not belonging to $X^{*}$. The minimizing condition 1) of the theorem is applicable to the subsequence $\left\{x^{k_{s}}\right\}$, which, in particular, means that there exists a lower limit point of $\left\{ x^{k}\right\}$, i.e., the point $x=\lim _{s \rightarrow \infty} x^{k_{s}}$ is not minimal. The obtained contradiction proves that $X^{*}$ is not empty, and all limit points of $\left\{x^{k}\right\}$ that are minimal with respect to $W$ belong to $X^{*}$. 

Let us prove the second assertion of the theorem. Let us show that the interval $[\underline{W}, \overline{W})$, where $\underline{W}, \overline{W}$ are defined by (8.2), is embedded in the set $W^{*}$. Assume the opposite, then there is $c$ such that $c \in[\underline{W}, \overline{W})$ and $c \notin W^{*}$. By the proved first assertion of the theorem, $\underline{W} \in W^{*}$; so $c>W$. Let's choose a number $d \in(c, \overline{W})$. Thus the interval $[c, d]$ is embedded in $(\underline{W}, \overline{W})$ and $c \notin W^{*}$. The numerical sequence $\left\{W\left(x^{k}\right)\right\}$ intersects the segment $[c, d]$ from $c$ to $d$ an infinite number of times, so there are such parts $ \left\{x^{k_{s}}, \ldots, x^{m_{s}}\right\} \quad(s=0,1, \ldots)$ of the sequence $\left\{x^{ k}\right\}$ that, either
\begin{equation}
W\left(x^{k_{s}}\right) \le c<W\left(x^{k}\right)<d \le W\left(x^{m_{s}}\right) , \quad k \in\left(k_{s}, m_{s}\right),\tag{8.4}
\end{equation}
or
$$
W\left(x^{k_{s}}\right) \le c<d \le W\left(x^{k_{s}+1}\right).
$$
The subsequence $\left\{x^{k_{s}}\right\}$ has limit points due to the compactness of the algorithm. Without loss of generality, we can assume that
$$
\lim _{s \rightarrow \infty} x^{k_{s}}=x^{\prime}.
$$

Let us show that $x^{\prime} \notin X^{*}$. Let's assume the opposite. Then by condition 2) of the theorem
$$
\lim _{s \rightarrow \infty} \rho\left(x^{k_{s}+1}, x^{k_{s}}\right)=0,
$$
and due to the continuity of $W(x)$
$$
\lim _{s \rightarrow \infty}\left[W\left(x^{k_{s}+1}\right)-W\left(x^{k_{s}}\right)\right]= 0,
$$
Together with (8.4), this gives
$$
W\left(x^{\prime}\right)=\lim _{s \rightarrow \infty} W\left(x^{k_{s}}\right)=\lim _{s \rightarrow \infty} W\left(x^{k_{s}+1}\right)=c \in W^{*},
$$
which contradicts the choice of $c$. Thus, $x^{\prime} \notin X^{*}$.

To the subsequence $\left\{x^{k_{s}}\right\}$ the minimizing condition 1) of the theorem is applicable. We choose $\delta$ such that
\begin{equation}
\max_y \left\{W(y) \mid \rho\left(y, x^{\prime}\right) \le \delta\right\} \le d.\tag{8.5}
\end{equation}
According to condition 1), there are indices $l_{s} \ge k_{s}$ such that for sufficiently large $s$ $\rho\left(x^{k}, x^{\prime}\right) \le \delta$ for $k \in\left[k_{s}, l_{s}\right)$, and
\begin{equation}
c \ge W\left(x^{\prime}\right)=\lim _{s \rightarrow \infty} W\left(x^{k_{s}}\right)>\varlimsup_{s \rightarrow \infty} W\left(x^{l_{s}}\right).\tag{8.6}
\end{equation}
By virtue of (8.5), for all $k \in\left[k_{s}, l_{s}\right)$ the inequality $W\left(x^{k}\right)<d$ holds, and then, according to the construction of $(8.4), l_{s} \le m_{s}$, whence $W\left(x^{l s}\right) \ge c$. Thus, we have obtained a contradiction between relation (8.6) and construction (8.4). The second assertion of the theorem is proved: $[\underline{W}, \overline{W}) \subset W^{*}$.

Let's prove the third assertion.
According to the second statement, $[\underline{W}, \overline{W}) \subset W^{*}$. But since $W^{*}$ does not contain segments, there is a limit
$$
\lim _{k \rightarrow \infty} W\left(x^{k}\right)=\varlimsup_{k \rightarrow \infty} W\left(x^{k}\right)= \varliminf_{s \rightarrow \infty} W\left(x^{k}\right).
$$

Now suppose there is a limit point
$$
x=\lim _{s \rightarrow \infty} x^{k_{s}},
$$
not contained in $X^{*}$. To the subsequence $\left\{x^{k_{s}}\right\}$ condition 1) of the theorem is applicable, which, in particular, means that $\left\{W\left(x^ {k}\right)\right\}$ has no limit. The obtained contradiction proves that all limit points of $\left\{x^{k}\right\}$ belong to $X^{*}$.

Note that here we have actually verified the conditions of the general Theorem 8.1.

Since now all limit points $\left\{x^{k}\right\}$ belong to $X^{*}$, it follows from condition 2) and the compactness of the algorithm that
$$
\lim _{k \rightarrow \infty} \rho\left(x^{k+1}, x^{k}\right)=0
$$

At the end of the section, in Lemma 8.2, we will show that then the set of limit points of  $\left\{x^{k}\right\}$ is connected and, therefore, belongs to a connected subset of $X^{*}$. The theorem has been proven.

The previous theorems, with varying degrees of detail, reflect the scheme of proving algorithm convergence by contradiction; therefore, some conditions of the theorems are formal in nature and psychologically inconvenient outside the context of proving convergence by contradiction. Let's reformulate Theorem 8.8 without the formal assumptions of contradiction.

Let $\left\{x^{k}\right\}$ be a sequence of initial data for the algorithm $A=\left(A_{k}, D^{k}\right)_{k \ge l}$. Denote by $\left\{x_{n}^{k}\right\}_{k=0}^{\infty}$ the sequence generated by $A$ with initial data $\left(x^{0}, \ldots, x^{n}\right)$, i.e.
$$
\left\{x_{n}^{k}\right\}_{k=0}^{\infty}=\left\{x_{n}^{k}=x^{k} \text { when } k \in[0, n] ; x_{n}^{k} \in A_{k-1}\left(x_{n}^{0}, \ldots, x_{n}^{k-1}\right) \text { with } k >n\right\} .
$$
\begin{theorem}
\label{th:8.9}
Let for a compact algorithm $A=\left(A_{k}\right.$, $\left.D^{k}\right)_{k \ge l}$ there is a continuous function $W: \overline{D} _{A} \rightarrow E_{1}$ and a set $X^{*}$ such that for any sequence of initial data $\left\{x^{k}\right\}$ the following condition holds:

1) If $\left\{x^{k_{s}}\right\} \rightarrow x \notin X^{*}$, then for any $\delta>0$ there are such indices $l_{s} \geq k_{s}$, that $\rho\left(x_{k_{s}}^{k}, x\right) \le \delta$ for $k \in\left[k_{s}, l_{ s}\right) $, and
$$
W(x)=\lim _{s \rightarrow \infty} W\left(x^{k_{s}}\right)>\varlimsup_{s \rightarrow \infty} W\left(x_{k_{s} }^{l_{s}}\right) \text {; }
$$
then $X^{*}$ is non-empty and the minimal (in values of ${W}$) limit points of any sequence $\left\{y^{k}\right\}$ generated by the algorithm, belong to $X^{* }$.

Let, in addition, the condition is met:
\begin{enumerate}
   \setcounter{enumi}{1}
   \item If $\left\{x^{k_{s}}\right\} \rightarrow x \in X^{*}$, then
\end{enumerate}
$$
\lim _{s \rightarrow \infty} \rho\left(x^{k_{s}}, x^{k_{s}+1}\right)=0.
$$
Then the algorithm $A$ converges to $ X^{*}$ with respect to the function values, i. e., for any sequence $\left\{y^{k}\right\}$ generated by $A$, the limit points of $\left\{y^{k}\right\}$ that are minimal in values of $W$, belong to $X ^{*}$, and the interval $[\lim _{k \rightarrow \infty} W\left(y^{k}\right), \varlimsup_{k \rightarrow \infty} W\left(y^{k }\right))$ is embedded in $W^{*}=\left\{W(x) \mid x \in X^{*}\right\}$.

Let in addition to the previous conditions:
\begin{enumerate}
   \setcounter{enumi}{2}
   \item The set $W^{*}$ does not contain segments. 
\end{enumerate}
Then the algorithm $A$ converges to $X^{*}$, i.e., for any sequence $\left\{y^{k}\right\}$ generated by $A$, all limit points of $\left\{y^{k}\right\}$ belong to a connected subset of $X^ {*}$, the number sequence $\left\{W\left(y^{k}\right)\right\}$ has a limit, and
$$
\lim _{k \rightarrow \infty} \rho\left(y^{k}, y^{k+1}\right)=0.
$$
\end{theorem}

{\it P r o o f.} Let $\left\{y^{k}\right\}$ be a sequence generated by algorithm $A$. Obviously, it is also a sequence of initial data for $A$; therefore, conditions 1), 2) of the theorem are applicable to it, where it is necessary to put $x^{k}=y^{k}$ and $x_{k_s}^{k}=y^{k }$. As applied to $\left\{y^{k}\right\}$, the conditions of Theorem 8.9 coincide with the analogous conditions of proved Theorem 8.8, that implies the validity of Theorem 8.9. The theorem has been proven.

Note that the conditions of Theorem 8.9 additionally imply that the local minima of $W$ are embedded in $X^{*}$; more precisely, if $x$ is a point where a local minimum of $W$ is reached, for which there exists a sequence of initial data $\left\{x^{k}\right\}$ that has $x$ as its a limit point, then $x $ belongs to $X^{*}$.

Indeed, otherwise, running the algorithm $A$ from the points $x^{k} \rightarrow x$, according to condition 1 ) of Theorem 8.9 , we will find values of $W$ smaller than $W(x)$ in an arbitrarily small
neighborhood of $x$, which contradicts the local optimality of $x$.

Thus, in Theorem 8.9 both the conditions for the convergence of optimization algorithms and the necessary conditions for the extremum of functions are tied.

To conclude this section, we prove a lemma on the connectedness of the sets of limit points of sequences, which was used in the proof of Theorems 8.8, 8.9.

First, we present some facts from [2] about the connectedness of sets in the metric space $M$.

\begin{definition}
\label{def:8.11}
A set $X \subset M$ is called connected if it cannot be represented as
$$
X=\left(\Phi_{1} \cap X\right) \cup\left(\Phi_{2} \cap X\right),
$$
where $\Phi_{1}$ and $\Phi_{2}$ are non-empty closed sets such that
$$
\left(\Phi_{1} \cap X\right) \cap\left(\Phi_{2} \cap X\right)=\emptyset.
$$
\end{definition}
\begin{definition}
\label{def:8.12}
A component of the set $X$ is called a connected subset that does not lie in another, wider connected subset of the set $X$.
\end{definition}

Each set $X \subset M$ splits into its connected components.

\begin{definition}
\label{def:8.13}
A finite sequence of points $x^{0}, \ldots$ $\ldots, x^{\text {s }}$ in a metric space $M$ is called an $\varepsilon$-chain (or an $\varepsilon$-chain connecting the point $x^{1}$ with point $x^{s})$ if $\rho(x^{i}, x^{i+1})<\varepsilon$ for $i \in[0, s)$. A set $X$ is said to be $\varepsilon$-linked if any two of its points can be connected by a $\varepsilon$-chain composed of the points of the set $X$. A set $X$ is called linked if it is $\varepsilon$-linked for any $\varepsilon>0$.
\end{definition}

Every connected set in a metric space is linked.

Every linked compact set is connected.

\begin{lemma}
\label{lem:8.2}
Let the sequence $\left\{x^{k}\right\}_{k=0}^{\infty}$ be embedded in some compact set and
$$
\lim _{k \rightarrow \infty} \rho\left(x^{k}, x^{k+1}\right)=0.
$$
Then the set $L$ of limit points $\left\{x^{k}\right\}$ is a connected compact set.
\end{lemma}

{\it P r o o f.} Let us show that $L$ is compact, i.e., any infinite sequence of points from $L$ has limit points. Let
$$
\left\{y^{n}\right\}_{n=0}^{\infty} \subset L.
$$
Since $\left\{y^{n}\right\}$ are limit points of $\left\{x^{k}\right\}$, there exist such points $x^{k_{n}}$ that
$$
\rho\left(x^{k_{n}}, y^{n}\right) \le 1 / n.
$$
By assumption, the subsequence $\left\{x^{k_{n}}\right\}$ has limit points; they will also be limiting for $\left\{y^{n}\right\}$. So, compactness is proved.

Let us show that $L$ is linked, hence the connectedness of $L$ will follow. Take an arbitrary $\varepsilon>0$. Let us show that
$$
\lim _{k \rightarrow \infty} \rho\left(x^{k}, L\right)=0 .
$$
Note that the function $\rho(x, L)$ is continuous. Let
$$
\varlimsup_{k \rightarrow \infty} \rho\left(x^{k}, L\right)=\lim _{s \rightarrow \infty} \rho\left(x^{k_{s}}, L \right), \quad \lim _{s \rightarrow \infty} x^{k_{s}}=x \in L .
$$
Then, due to the continuity of $\rho(x, L)$
$$
0 \le \varlimsup_{k \rightarrow \infty} \rho\left(x^{k}, L\right)=\lim _{s \rightarrow \infty} \rho\left(x^{k_{s} }, L\right)=\rho(x, L)=0,
$$

Let now for all $k \ge K$ be
$$
\rho\left(x^{k}, x^{k+1}\right)<\varepsilon / 3, \quad \rho\left(x^{k}, L\right)<\varepsilon / 3 .
$$
Let $a, b \in L$. Find an $\varepsilon$-chain in $L$ connecting $a$ and $b$. Since $a$ and $b-$ are limit points of $\left\{x^{k}\right\}$, then there are $k_{2}>k_{1}>K$ such that
$$
\rho\left(a, x^{k_{1}}\right)<\varepsilon / 3, \quad \rho\left(b, x^{k_{2}}\right)<\varepsilon / 3 .
$$
For points $\left\{x^{k}\right\}_{k=k_{1}}^{k_{2}}$, by construction, there are $y^{k} \in L$ such that
$$
\rho\left(x^{k}, y^{k}\right)<\varepsilon / 3, \quad k \in\left[k_{1}, k_{2}\right] .
$$
It's obvious that
$$
\rho\left(y^{k}, y^{k+1}\right) \le \rho\left(y^{k}, x^{k}\right)+\rho\left(x^ {k}, x^{k+1}\right)+\rho\left(x^{k+1}, y^{k+1}\right) \le \varepsilon
$$
for $k \in\left[k_{1}, k_{2}\right) ; \rho\left(a, y^{k_{1}}\right) \le 2 \varepsilon / 3, \rho\left(y^{k_{2}}, b\right) \le 2 \varepsilon / 3$.

Thus, $\left\{a, y^{k_{1}}, \ldots, y^{k_{2}}, b\right\}$ is the required $\varepsilon$-chain. Since $\varepsilon$ is arbitrary, $L$ is a linked set. The lemma is proven.

\section*{$\S$ 9. Generalized gradient decent method}
\label{Sec.9}
\setcounter{section}{9}
\setcounter{definition}{0}
\setcounter{equation}{0}
\setcounter{theorem}{0}
\setcounter{lemma}{0}
\setcounter{remark}{0}

In this section, the well-known generalized gradient method for minimizing convex functions is extended to problems of unconstrained minimization of nonconvex nonsmooth generalized differentiable functions.

Non-smooth optimization begins with the generalized gradient method. Initially, it was applied to minimize (with a constant step) convex piecewise linear functions [134], then the convergence conditions of the method were studied [40, 99], and numerous modifications and generalizations were constructed. We extend this important method to the case of minimizing a wide class of (nonconvex nonsmooth) generalized differentiable functions.

Recall that, in general, the generalized 
differentiable functions have no directional derivatives.

Consider the problem:
\begin{equation}
\label{eqn:9.1}
f(x) \rightarrow \min_x, \quad x \in E_{n},
\end{equation}
where $f(x)$ is a generalized differentiable function such that $f(x) \rightarrow \infty$ as $\|x\| \rightarrow \infty$.

Let us define the set of (pseudo-)stationary points of problem (9.1):
$$
X^{*}=\left\{x \in E_{n} \mid 0 \in G_{f}(x)\right\},
$$
where $G_{f}(x)$ is the set of pseudo-gradients of $f(x)$ at point $x$.

The set $X^{*}$ is closed because the pseudogradient map $G_{f}$ is closed. According to Theorem 3.1, all local minimum points of $f(x)$ are embedded in $X^{*}$.

We introduce the stationary value set
$$
F^{*}=\left\{f(x) \mid x \in X^{*}\right\}.
$$

The generalized gradient method has the form:
\begin{eqnarray}
x^{0} \in E_{n}, \quad x^{k+1}=x^{k}-\rho_{k} g^{k}, \quad g^{k} \in G_{ f}\left(x^{k}\right), \quad k=0,1, \ldots, \label{eqn:9.2}\\
0 \le \rho_{k} \le R, \quad \lim _{k \rightarrow \infty} \rho_{k}=0, \quad \sum_{k=0}^{\infty} \rho_{k }=+\infty .\label{eqn:9.3}
\end{eqnarray}

The following statements establish convergence of the generalized gradient method; they will be proved later on.
\begin{theorem}
\label{th:9.1}
Assume that the sequence $\left\{x^{k}\right\}_{k=0}^{\infty}$ generated by algorithm (9.2) - (9.3) is bounded. Then it converges to the solution of problem (9.1) in function, i.e., the limit points $\left\{x^{k}\right\}_{k=0}^{\infty}$ belong to $X^{*}$, and all limit points of the numerical sequence $\left\{f\left(x^{k}\right)\right\}_{k=0}^{\infty} $ constitute an interval lying in the set $F^{*}$. If the set $F^{*}$ does not contain intervals (for example, it is finite or countable), then the sequence $\left\{x^{k}\right\}_{k=0}^{\infty}$ converges to the solution of problem $(9.1)$, i.e., all limit points $\left\{x^{k}\right\}_{k=0}^{\infty}$ belong to a connected subset of $X^{*}$, and the numerical sequence $\left\{f\left(x^{k}\right)\right\}_{k=0}^{\infty}$ has a limit.
\end{theorem}

The theorem follows from Lemma 9.2, proved below, and the general Theorem 8.9. Since in our case $F^{*}$ is closed, then
$$
\left[ \underline\lim_{k \rightarrow \infty} f(x^{k}), \; \overline{\lim _{k \rightarrow \infty}} f(x^{k})\right] \subset F^{*} .
$$

The next lemma shows that if the steps $\rho_{k}$ in the generalized gradient method are small enough, then the generated sequences are bounded.
\begin{lemma}
\label{lem:9.1}
Let there exist $c>f\left(x^{0}\right)$ not belonging to $F^{*}$. Then for any $d>c$ there exists $\overline{R}$ such that for any sequence $\left\{\rho_{k}\right\}_{k=0}^{\infty}$ satisfying (9.3) with $ R \le \overline{R}$, sequences of points $\left\{x^{k}\right\}$ launched from $x^{0}$ according to $(9.2)$, do not leave the bounded region $\{x \mid f(x) \le d\}$.
\end{lemma}
\begin{remark}
\label{rem:9.1}
The number $c$ required in Lemma 9.1 certainly exists if the set $X^{*}$ (and hence $F^{*}$ ) is bounded or if $F^{*}$ does not contain intervals.
\end{remark}
\begin{remark}
\label{rem:9.2}
Practically, the bound 
$$
R=\sup _{k \ge 0} \rho_{k}
$$
is selected in the following way [88]. 
Suppose there exists a number $d > f(x^0)$ such that the interval $[f(x^0), d)$ does not lie entirely in $F^*$. Then, according to Lemma 9.1, there exists an $\bar{R}$ such that for $R \leq \bar{R}$, the sequences $\{x^k\}$ are uniformly bounded. Let us run algorithm (9.2) from the point $x^0$ with some $\{\rho_k\}_{k=0}^\infty$ and $R_0 = \sup_{k\geq 0}\rho_k$. If at some step it turns out that $f(x^k) > d$, then we interrupt the process, choose another sequence $\{\rho^\prime_k\}_{k=0}^\infty$ (with a smaller $R_1 = \sup_{k\geq 0}\rho^\prime_k$ and run algorithm (9.2) with new $\rho^\prime_k$ from the initial (or record) point $x^0$ and so on. If $R_s$ is set so that $R_s\rightarrow 0$ as $s\rightarrow\infty$, then sooner or later there will be $R_s < \bar{R}$, and the corresponding trajectory $\{x^k\}$ will not leave the bounded region $\{x | f(x) \leq d\}$.
\end{remark}

We split the proof of Lemma 9.1. into series of lemmas. First, we will prove some geometric lemma.
\begin{lemma}
\label{lem:9.2}
Let $ P $ be a convex subset in $ E_n $, where $ 0 < \gamma \leq || p || \leq \Gamma < + \infty $ for all $ p \in P$. Then for any set of vectors $ \{ p^r \in P | r = k, ..., m \} $ and for a set of numbers $ \{ \rho_r \in E_1 | \;r = k, ..., m, \rho_r \geq 0 \} $ such that
$$
\sum_{r=k}^{m} \rho_r \geq \sigma > 0,\;\;\; \sup_{k \leq r \leq m} \rho_r \leq R = \frac{\gamma^3 \sigma}{6 \Gamma^3},
$$
there exists an index $ l \in (k, m] $, for which
$$
\left<p^l, \frac{\sum_{r=k}^{l-1} \rho_r p^r}{\sum_{r=k}^{l-1} \rho_r}\right> \geq \frac{\gamma^2}{4}, \;\;\;\;\;\;\sum_{r=k}^{l-1} \rho_r \geq \frac{\gamma \sigma}{3 \Gamma}.
$$
\end{lemma}

{\it P r o o f.} Select indexes $t$ and $m^\prime\leq m$ such that:
$$
\sum_{r=k}^{t-1} \rho_r < \frac{\gamma \sigma}{3 \Gamma} \leq \sum_{r=k}^{t} \rho_r,
\;\;\;\;\;\; \sum_{r=k}^{m'-1} \rho_r < \sigma \leq \sum_{r=k}^{m'} \rho_r.
$$
Using proof by contradiction, we assume that for all $ l \in (t, m'] $
$$
\left<p^l, \frac{\sum_{r=k}^{l-1} \rho_r p^r}{\sum_{r=k}^{l-1} \rho_r}\right> < \frac{\gamma^2}{4}.
$$
Then we consider the following representation:
$$
\left\| \sum_{r=k}^{l} \rho_r p^r \right\|^2 = \left\| \sum_{r=k}^{l-1} \rho_r p^r \right\|^2 + 2 \rho_l \left<p^l, \sum_{r=k}^{l-1} \rho_r p^r\right> + \rho_l^2 \| p^l \|^2.
$$
Let us sum  these equations from $ t + 1 $ to $ m' $ along the index  $ l $, we get 
$$
\left\| \sum_{r=k}^{m^\prime} \rho_r p^r \right\|^2 = \left\| \sum_{r=k}^{t} \rho_r p^r \right\|^2 
+ 2 \sum_{k=t+1}^{m^\prime}\rho_l \left<p^l, \sum_{r=k}^{l-1} \rho_r p^r\right> 
+ \sum_{k=t+1}^{m^\prime}\rho_l^2 \| p^l \|^2.
$$
Under the made assumptions, we get the estimates: 
\begin{eqnarray}
    \gamma^2 \sigma^2 &=& \Gamma^2 \left(\sum_{r=k}^{t} \rho_r\right)^2 + \frac{\gamma^2}{2} \sum_{l=t+1}^{m'} \rho_l \sum_{r=k}^{l-1} \rho_r + \Gamma^2 \sup_{k \leq r \leq m'} \rho_r \sum_{l=t+1}^{m'} \rho_l \nonumber \\
    &\leq& \Gamma^2 \left(\frac{\gamma}{3 \Gamma} + \frac{\gamma^2}{6 \Gamma^2}\right)^2 \sigma^2 + \frac{\gamma^2}{2} \sigma ^ 2 + \Gamma^2 \frac{\gamma^2 \sigma}{6 \Gamma^2} \sigma \leq \frac{11}{12} \gamma^2 \sigma^2. \nonumber
\end{eqnarray}
The obtained contradiction proves the lemma.

The following lemma establishes the local minimizing property of the generalized gradient method: if the generalized gradient method is started from a non-stationary point $x$ with a sufficiently small step, then the method finds smaller values of the function $f$ in a small neighborhood of the starting point $x$ than $f(x)$. This property implies the fulfillment of condition 1) of Theorem 8.9 and basically ensures the convergence of the method.

\begin{lemma}
\label{lem:9.3}
Let a sequence of initial points $ \{ x^s \}_{s=0}^{\infty} $ converges to point $ x = \lim_{s \to \infty} x^s $. We will run the generalized gradient method (9.2) from the starting points $ x^s $ with the step values of $ \rho_s^k \;(k \geq k_s) $ until some moment $ n_s $. For each $ s, $ we get sequences $ \{ x_s^k \}_{k \geq k_s} $:
$$
x_s^{k_s} = x^s, \;\;\;x_s^{k+1} = x_s^{k} - \rho_s^k g_s^k, \;\;\;g_s^k \in G_j (x_s^{k}), \;\;\;k = k_s, k_s + 1, \ldots .
$$
Denote $ \rho_s = \sup_{k \geq k_s} \rho_s^k $ and $ \sigma_s = \sum_{k = k_s}^{n_s - 1} \rho_s^k $. 

If $ 0 \notin G_f (x) $, $ \sigma_s \geq \sigma > 0 $ and $ \lim_{s \to \infty} \rho_s = 0 $, then there exists such number $ \bar{\varepsilon} > 0 $ that for any $ \varepsilon \in (0, \bar{\varepsilon}] $ there are such indices $ l_s $, that for sufficiently large $ s $ holds
$\| x_s^k - x \| \leq \varepsilon $ for $ k \in [ k_s, l_s ]$, and:
\begin{equation*}\tag{9.4}
     f(x) = \lim_{s \to \infty} f(x^s) > \overline{\lim}_{s \to \infty} f(x^{l_s}).
\end{equation*}
\end{lemma}
{\it P r o o f.} Let
$$
\gamma_1 = \rho (0, G_f(x)) = \inf \{ \| g \| \; | \; g \in G_f(x) \}.
$$
By  upper semi-continuity of $ G_f $, there exists an $ \varepsilon_1 $-neighborhood of the point $ x $ such that for all $ g \in G_f(y) (\| y - x \| \leq \varepsilon_1) $ it holds:
$$
\rho (g, G_f(y)) \leq \frac{\gamma_1}{2}.
$$
Denote
$$
\Gamma_1 = \sup \{ \| g \| \; | \; g \in G_f(y), \| y - x \| \leq \varepsilon_1 \}.
$$
Due to the generalized differentiability of $ f(x) $ for $ c_f = \frac{\gamma_1^2}{16 \Gamma_1} $ there is $ \varepsilon_2 \leq \varepsilon_1 $ such that for $ \| y - x \| \leq \varepsilon_2 $ it holds
\begin{equation}\tag{9.5}
    f(y) = f(x) + \left<g, y - x\right> + o (x, y, g),
\end{equation}
where $ g \in G_f (y) $ and $ | o(x, y, g) | \leq c_f \| y - x \| $. We define $ \bar{\varepsilon} = \min \{ \varepsilon_2, \sigma \frac{\gamma_1}{2} \} $ and fix an arbitrary constant $ \varepsilon \leq \bar{\varepsilon} $. Let index $ s $ be  such that $ \| x^s - x \| \leq \frac{\varepsilon}{4} $ and $ \rho_s \leq \frac{\varepsilon}{4 \Gamma_1} $.

Let us show that the sequence $ \{ x_s^k \}_{k \geq k_s} $ leaves the $ \frac{\varepsilon}{2} $-neighbourhood of elements $ x^s $. Assume the opposite, that $ \| x_s^k - x^s \| \leq \frac{\varepsilon}{2} $ for all $ k \in [k_s, n_s] $. Then on the one hand, we get:
$$
\| x_s^{n_s} - x^s \| \leq \frac{\varepsilon}{2} \leq \sigma \frac{\gamma_1}{4}.
$$
But on the other hand:
$$
\| x_s^{n_s} - x^s \| = \left\| \sum_{k=k_s}^{n_s - 1} \rho_s^k g_s^k \right\| = 
\left\| \frac{\sum_{k=k_s}^{n_s - 1} \rho_s^k g_s^k}{\sum_{k=k_s}^{n_s - 1} \rho_s^k} \right\| \sum_{k=k_s}^{n_s - 1} \rho_s^k.
$$
For $ k \in [k_s, n_s] $ it holds
$$
\| x_s^k - x \| \leq \| x_s^k - x^s \| + \| x^s - x \| \leq \frac{3 \varepsilon}{4},
$$
thus
$$
\rho (g_s^k, G_f (x)) \leq \frac{\gamma_1}{2}.
$$
Since the set $ G_f(x) $ is convex and $ \rho (0, G_f(x)) = \gamma_1 $, then it follows that all $ g_s^k \;(k \in [k_s, n_s]) $ belong to a half-space separated from zero by quantity $ \frac{\gamma_1}{2} $. The same holds for convex combinations of vectors $ g_s^k (k \in [k_s, n_s]) $. As a result, we obtain $ \| x_s^{n_s} - x^s \| \geq \sigma \frac{\gamma_1}{2} $. The resulting contradiction proves that the sequences $ \{ x_s^k \} $ come out of the $ \frac{\varepsilon}{2}$-neighborhoods of points $ x^s $. 

We denote by $ m_s $ the time of the first exit of the sequence $ \{ x_s^k \}_{k \geq k_s} $ from the $ \frac{\varepsilon}{2} $-neighborhood of $ x^s $, i.e., for all $ k < m_s $:
$$
\| x_s^k - x^s \| \leq \frac{\varepsilon}{2}, \> \> \> \| x_s^{m_s} - x^s \| > \frac{\varepsilon}{2}
$$
Now the following estimate is true:
$$
\| x_s^{m_s} - x^s \| \leq \| x_s^{m_s - 1} - x^s \| + \rho_s^{m_s - 1} \| g_s^{m_s - 1} \| \leq \frac{3 \varepsilon}{4}
$$
and for all $ k \in [k_s, m_s] $:
$$
\| x_s^k - x \| \leq \| x_s^k - x^s \| + \| x^s - x \| \leq \varepsilon.
$$

For $ k \in [k_s, m_s] $ we substitute vectors of $ x_s^k $ and $ g_s^k $ in the equation (9.5):
\setcounter{equation}{5}
\begin{eqnarray} 
    f(x_s^k) &=& f(x) + \left<g_s^k, x_s^k - x\right> + o(x, x_s^k, g_s^k) \leq \nonumber \\
             &\leq& f(x) + \left<g_s^k, x_s^k - x^s\right>+ \left<g_s^k, x^s - x\right> + c_f \| x_s^k - x \| \leq  \nonumber \\
             &\leq& f(x) + \left<g_s^k, x_s^k - x^s\right> + c_f \| x_s^k - x \| + (\Gamma_1 + c_f ) \| x^s - x \|.\label{eqn:9.6}
\end{eqnarray}

Consider the scalar products
$$
\left<g_s^k, x_s^k - x^s\right> = - \left<g_s^k, \frac{\sum_{r=k_s}^{k-1} \rho_s^r g_s^r}{\sum_{r=k_s}^{k-1} \rho_s^r}\right> \sum_{r=k_s}^{k-1} \rho_s^r.
$$
To estimate them, we can apply Lemma 9.2, where we have to put
$$
P = co \{ G_j (y) \> | \> \| y - x \| \leq \varepsilon \}, \> p^r = g_s^r, \> k_s \leq r \leq m_s; \> \Gamma = \Gamma_1, \> \gamma = \frac{\gamma_1}{2}.
$$

The following conditions are satisfied:
$$
\sum_{r=k_s}^{m_s - 1} \rho_s^r \geq \frac{\| x_s^{m_s} - x^s \|}{\Gamma_1} \geq \frac{\varepsilon}{2 \Gamma_1} > 0, \> \> \> \lim_{s \to \infty} \rho_s = 0.
$$

Now, according to Lemma 9.2, for sufficiently large $ s $ there exist such indices $ l_s \leq m_s $ that
$$
\left<g_s^{l_s}, \frac{\sum_{k=k_s}^{l_s-1} \rho_s^k g_s^k}{\sum_{k=k_s}^{l_s-1} \rho_s^k}\right> \geq \frac{\gamma_1^2}{8}, \> \> \> \sum_{k=k_s}^{l_s-1} \rho_s^k \geq \frac{\epsilon \gamma_1}{12 \Gamma_1^2}.
$$

Let us rewrite inequality (9.6) for $ k = l_s $ taking into account the obtained estimates:
\begin{eqnarray}
    f(x_s^{l_s}) &= &- \frac{\gamma_1^2}{8} \sum_{k=k_s}^{l_s-1} \rho_s^k + c_f \Gamma_1 \sum_{k=k_s}^{l_s-1} \rho_s^k + (\Gamma_1 + c_f) \| x^s - x \| \nonumber \leq \\
                 &\leq& f(x) - \frac{\gamma_1^2}{192 \Gamma_1^2} \gamma_1 \epsilon + (\Gamma_1 + c_f) \| x^s - x \|.\label{eqn:9.7}
\end{eqnarray}

So, for all sufficiently small $ \epsilon $ for all sufficiently large $ s $ we have specified such indices $ l_s $ that $ \| x_s^k - x \| \leq \epsilon $ at $ k \in [k_s, l_s] $ and (9.7) is true. For other $ s $ the indices $ l_s $ are arbitrary. Then for the upper limit $ \overline{\lim}_{s \to \infty} f(x_s^{l_s}) $, by virtue of (9.7), the estimate (9.4) is valid. The lemma is proved.

{\it P r o o f of Lemma 9.1.} We need to show the boundedness under the conditions of Lemma 9.1 of the sequences $ \{ x^k \}_{k=0}^{\infty} $ generated by the generalized gradient method (9.2), (9.3) for sufficiently small $ R = \sup_{k \geq 0} \rho_k $.

Suppose the contrary: whatever $ R_s \to 0 (s \to \infty) $ is, there exists a sequence $ \{ x_s^k \}_{k=0}^{\infty} $ generated by (9.2), (9.3) with $ \sup_{k \geq 0}\; \rho_s^k \leq R_s $, going to infinity.

Choose a number $ d $ such that $ f(x^0) < c < d < \infty $. Denote 
$$
\Gamma = \sup \{ \> \| g \| \> | \> g \in G_f(y), f(y) \leq d \}.
$$
Since $ G_f $ is compact-valued and semi-continuous from above, it is bounded on the compact $ \{ \> x \> | \> f(x) \leq d \} $, i.e. $ \Gamma < \infty $. Every numerical sequence $ \{ f(x_s^k) \}_{k=0}^{\infty} \> (s = 0, 1, ...) $ intersects the segment $ [c, d] $ from $ c $ to $ d $; so there are parts $ \{ x_s^{k_s}, ..., x_s^{n_s} \} \> (s = 0, 1, ...) $ of sequences $ \{ x_s^k \}_{k \geq 0} \> (s = 0, 1, ...) $ such that either
\begin{equation}
\label{eqn:9.8}
    f(x_s^{k_s}) \leq c < f(x^k) < d  \leq f(x_s^{n_s}), \> \> \> k \in (k_s, n_s)
\end{equation}
or
$$
f(x_s^{k_s}) \leq c < d  \leq f(x_s^{k_s + 1})
$$
The sequence $ \{ f(x_s^{k_s}) \}_{s=0}^{\infty} $ belongs to the compact set $ \{ \> x \> | \> f(x) \leq c \} $; therefore it has limit points. Without changing the notation, we will assume that
$$
\lim_{s \to \infty} x_s^{k_s} = x'.
$$
By construction, $ f(x') \leq d $. Let us show that $ x' $ does not belong to $ X^{*} $. Indeed, suppose the contrary: $ x' \in X^{*} $. Then
$$
0 \leq \lim_{s \to \infty} \| x_s^{k_s + 1} - x_s^{k_s} \| \leq \lim_{s \to \infty} \Gamma R_s = 0.
$$
Hence, together with (9.8), it follows that
$$
\lim_{s \to \infty} f(x_s^{k_s}) = \lim_{s \to \infty} f(x_s^{k_s + 1}) = f(x') = c \in F^{*},
$$
which contradicts the choice of $ c $. So, $ x' \notin X^* $. For sufficiently small $ \varepsilon $ due to (9.8) and continuity of $ f(x) $ the sequence $ \{ x_s^{k_s} \}_{s=0}^{n_s} $ leaves the $ \varepsilon $-neighborhood of the point
$$
x' = \lim_{s \to \infty} x_s^{k_s}.
$$
Therefore, for sufficiently large $ s $
$$
\sum_{k=k_s}^{n_s-1} \rho_s^k \geq \frac{\| x_s^{n_s} - x_s^{k_s} \|}{\Gamma} \geq \frac{\varepsilon}{2 \Gamma}.
$$

We can assume that the generalized gradient method as if repeatedly starts from points $ x_s^{k_s} \to x' \notin X^* $. We can now see that lemma 9.3 can be applied to the initial data sequence $ \{ x_s^{k_s} \}_{s=0}^{\infty} $. Let us take such $ \varepsilon < \bar{\varepsilon} (x') $ that
$$
\max \{ f(x) \> | \> \| x' - x \| \leq \varepsilon \} \leq d,
$$
then the indices $ l_s $, from Lemma 9.3 are not larger than $ n_s $. But this means that the minimizing property (9.4) contradicts the construction (9.8). The resulting contradiction proves Lemma 9.1.

The generalized gradient method in the form (9.2), (9.3) has mainly theoretical value due to extremely slow convergence. There are many ways to improve it. One possibility is to normalize the pseudo-gradients in formulas (9.2).

Denote by $ \bar{g} $ the normalized vector $ g $:
\[ 
    \bar{g} = 
    \begin{cases}
        g / \| g \|, & g \neq 0, \nonumber \\
        0, & g= 0. \nonumber
    \end{cases}
\]

In method (9.2), we can substitute normalized pseudo-gradients $ \bar{g}^k $ instead of $ g^k $. For the method modified in this way, Theorem 9.1 and Lemmas 9.1, 9.3 remain valid.
\begin{lemma}
\label{lem:9.4}
For the normalized pseudo-gradient method the local minimizing property formulated in Lemma 9.3 (in which $ g^k $ should be replaced by $ \bar{g}^k $) remains valid.
\end{lemma}

{\it P r o o f.} Let us represent
$$
x_s^{k+1} = x_s^{k} - \rho_s^{k} \bar{g}_s^{k} = x_s^{k} - \bar{\rho}_s^{k} g_s^{k},
$$
where
\[ 
    \bar{\rho}_s^{k} = 
    \begin{cases}
        \rho_s^{k} / \| g_s^{k} \|, & g_s^{k} \neq 0, \nonumber \\
        \rho_s^{k}, & g_s^{k} = 0. \nonumber
    \end{cases}
\]
Let us show that lemma 9.3 applies to the sequence $ \{ x_s^{k} \}_{k \geq k_s} $ in a small neighborhood of point $ x $, from which lemma 9.4 will follow. Let us define 
$$
\gamma_1 = \rho (0, G_f(x)) > 0.
$$
Since $ G_f $ is upper semi-continuous, there exists an $ \varepsilon_1 $-neighborhood of point $ x $, such that for any $ g \in G_f(y) $ ($ \| y - x \| \leq \varepsilon_1 $) it take place the inequality $ \rho (g, G_f(x)) \leq \frac{\gamma_1}{2} $.

Let's denote as before
$$
\Gamma_1 = \sup \{ \> \| g \| \> | \> g \in G_f (y), \> \| y - x \| \leq \varepsilon_1 \}.
$$
Let's introduce the indices
$$
\bar{t}_s = \sup \{ \> t \> | \> \| x_s^k - x \| \leq \varepsilon_1, \> \forall k \in [k_s, t), t \in [k_s, n_s] \}.
$$
With $ k \in [k_s, \bar{t}_s) $ we have the following estimates
$$
\frac{\rho_s^k}{\Gamma_1} \leq \bar{\rho}_s^k \leq \frac{2 \rho_s^k}{\gamma_1}.
$$

We assign
$$
\bar{\rho}_s = \sup_{k_s \leq k < \bar{t}_s} \bar{\rho}_s^k, \> \> \> \bar{\sigma}_s = \sum_{k=k_s}^{\bar{t}_s - 1} \bar{\rho}_s^k.
$$
Obviously,
$$
\bar{\rho}_s \leq \frac{2 \rho_s}{\gamma_1} \to 0, \> \> \> s \to \infty.
$$

If $ \bar{t}_s = n_s $, then $ \bar{\sigma}_s \geq \frac{\sigma_s}{\Gamma_1} $. If $ \bar{t}_s < n_s $, then with significantly large value $ s $
$$
\bar{\sigma}_s \geq \frac{|| x_s^{\bar{t}_s} - x^s ||}{\Gamma_1} \geq \frac{\varepsilon_1}{2 \Gamma_1}.
$$

Thus,
$$
\bar{\sigma}_s \geq \min \{ \frac{\varepsilon_1}{2}, \sigma \} / \Gamma_1 = \bar{\sigma} > 0.
$$

It is now clear that lemma 9.3 applies to sequences $ \{ x_s^k \}_{k=s}^{\bar{t}_s} $ with sufficiently large s, from which lemma 9.4 follows.
\begin{theorem}
\label{th:9.2}
For the algorithm (9.2), (9.3) in which vectors $ g^k $ are replaced by normalized $ \bar{g}^k $, the results about convergence in the form of Theorem 9.1 and Lemma 9.1 remain valid.
\end{theorem}

Similarly to Theorem 9.1, Theorem 9.2 by virtue of Lemma 9.4 and
$$
\lim_{k \to \infty} \| x^{k+1} - x^k \| \leq \lim_{k \to \infty} \rho_k = 0,
$$
follows from the general theorem 8.9.
\begin{remark}
\label{rem:9.3}
An attentive reader may notice that the set of pseudo-stationary points $ X^* = \{ \> x \> | \> 0 \in G_f(x) \} $, to which the method (9.2) converges, may be much wider than the set of stationary points $ \{ \> x \> | \> 0 \in \partial f(x) \} $. This fact devalues the results on the convergence of the method (9.2). Therefore, consider the randomized generalized gradient method:
\begin{equation}
\label{eqn:9.9}
    x^0 \in E_n, \> \> \> x^{k+1} = x^k - \rho_k g^k, \> \> \> g^k \in G_f (\tilde{x}^k)
\end{equation}
\begin{equation}
\label{eqn:9.10}
    \| x^k - \tilde{x}^k \| \leq \varepsilon_k, \> \> \> \lim_{k \to \infty} \varepsilon_k = 0,
\end{equation}
where $ \tilde{x}^k $ is a random vector with a uniform distribution in the $ \varepsilon_k $-neighborhood of point $ x^k $. Under conditions (9.3), (9.10) the sequence $ \{ x^k \} $, according to Theorem 14.1, always converges to $ X^* $. But in this method at each iteration with probability one $ g^k \in \partial f(\tilde{x}^k) $, since $ G_f(x) $ by virtue of Theorem 1.12 coincides almost everywhere in $ E_n $ with $ \partial f(x) $. Thus in (9.9), (9.10) for almost all trajectories $ \{ x^k \} $ only the generalized Clarke gradients of the function $ f(x) $ are used. And since the mapping $ \partial f(x) $ is also pseudo-gradient for $ f(x) $, by virtue of Theorem 14.1 the method (9.9), (9.10) converges to the set $ \{ x \> | \> 0 \in \partial f(x) \} $ for almost all trajectories $ \{ x^k \} $. It remains to note that in practice there is no need to resort to randomization of the gradient choice at all, since it occurs at the expense of rounding errors anyway.
\end{remark}

\section*{$\S$ 10. The generalized gradient method in conditional optimization problems}
\label{Sec.10}
\setcounter{section}{10}
\setcounter{definition}{0}
\setcounter{equation}{0}
\setcounter{theorem}{0}
\setcounter{lemma}{0}
\setcounter{remark}{0}
\setcounter{corollary}{0}

In this section, the generalized gradient method is applied to local minimization of a generalized differentiable function under generalized differentiable inequality constraints and linear inequality constraints. 

First, the optimization problem with only one constraint-inequality is considered. The problem with several constraint-inequalities is easily reduced to it. To solve it, the generalized gradient method (of the type [99]) is used which uses pseudo-gradients of the target function inside the admissible domain and pseudo-gradients of the constraint outside it. 

Then we consider an optimization problem with several inequality constraints ordered by importance; it belongs, generally speaking, to multi-criteria optimization. A mathematical programming problem with constraint-inequalities falls into this class if we introduce an (arbitrary) order among the constraints. To solve the problem, the generalized gradient method is used, which consistently minimizes and achieves fulfillment of the constraints decreasing in importance. In this case the target function is considered to be the least important criterion and is minimized last.

To solve optimization problems containing linear constraint-inequalities we propose to use known methods of non-smooth optimization, including those considered in this section, in which generalized gradients are replaced by their projections onto a linear subspace corresponding to linear constraints. Consider the mathematical programming problem:
\begin{equation}
\label{eqn:10.1}
    f_0 (x) \to \min, \> \> \> x \in E_n
\end{equation}
subject to constraints
\begin{equation}
\label{eqn:10.2}
    f_i (x) \leq 0, \> \> \> i = 1, 2, ..., m
\end{equation}
where $ f_i (x) $ are generalized differentiable functions, $ G_i (x) $ are their pseudo-gradient mappings.

Let's reduce all constraints-inequalities of problem (10.1), (10.2) to a single one. This can be done in many ways. For example, the constraints (10.2) are equivalent to one constraint of the form
\begin{equation}
\label{eqn:10.3}
    h(x) \equiv \max_{1 \leq i \leq m} f_i (x) \leq 0
\end{equation}
or
\begin{equation}
\label{eqn:10.4}
    h(x) \equiv \sum_{i=1}^{m} \max(0, f_i (x)) \leq 0
\end{equation}

The first function is convenient because when using it in  optimization methods it is necessary to calculate only values of its constituent $ f_i (x) $ functions and one gradient of the function on which the maximum is reached.

For the second function, it is necessary to calculate values and gradients of all incoming functions. However, the second function is easier to minimize because it minimizes several functions $ f_i (x) $ at once. Therefore, if $ f_i (x) $ is simple and scaled, it is better to use the second function $ h (x) $. Note that in practice, in the case of the constraint (10.4), instead of zero in the right-hand side, one should put a small positive number in order to correctly determine the admissibility of points.

Naturally, the convolution of the original constraints into a single non-smooth constraint of the form (10.3) or (10.4) is designed for the worst case, i.e. it is assumed that $ f_i (x) (i = 1, 2, ..., m) $ are sufficiently general generalized differentiable functions. If the constraints $ f_i (x) \leq 0 $ have a special form, then to solve the problem we should  apply or develop special methods.

Recall that the calculation of the pseudo-gradients of complex functions (10.3), (10.4) should be done according to the rules of $\S$ 1.

\textbf{1. One inequality-constraint problem.} 
\label{Sec.10.1}
We will consider a problem of the form
\begin{equation}
\label{eqn:10.5}
    f(x) \rightarrow  \min_x
\end{equation}
subject to constraints
\begin{equation}
\label{eqn:10.6}
    h(x) \leq 0, \> \> \> x \in E_n
\end{equation}
where $ f(x) $ and $ h(x) $ are generalized differentiable functions.

  Without significant generality restriction, we can assume that
  $ h(x) \to \infty $ for $ \| x \| \to + \infty $. If this is not the case, we can always take instead of $ h(x) $ the function
  $$
    H(x) = \max (h(x), \|x\| - C),
  $$
  where $ C $ is a big enough constant.

  Let us define a multivalued mapping $ x \to G(x) $:
  \begin{equation*} 
    G(x) = 
    \begin{cases}
        G_f(x), & h(x) < 0, \nonumber \\
        co \{ G_f(x), G_h(x) \}, & h(x) = 0, \nonumber \\
        G_h(x), & h(x) > 0. \nonumber
    \end{cases}
  \end{equation*}

  According to Theorem 3.1, if $ x^* $ is a local minimum point of $ f(x) $ in the region $ \{ x \> | \> h(x) \leq 0 \} $, then $ 0 \in G(x^*) $ and $ h(x^*) \leq 0 $.

  The generalized gradient method for solving problem (10.5), (10.6) has the following form:
  \begin{equation}
	\label{eqn:10.7}
    x^0 \in E_n, \> \> \> x^{k+1} = x^k - \rho_k g^k, \> \> \> g^k\in G (x^k), \> \> \> k = 0, 1, ...,
  \end{equation}
  \begin{equation}
	\label{eqn:10.8}
    0 \leq \rho_k \leq R, \> \> \> \lim_{k \to \infty} \rho_k = 0, \> \> \> \sum_{k=0}^{\infty} \rho_k = \infty.
  \end{equation}

  Outside the admissible region, the method (10.7) minimizes the function $ h(x) $, and inside it the target function $ f(x) $.

  Instead of $ g^k $ in the method (10.7) we can use normalized values
  \begin{equation*} 
    \bar{g}^k =  
    \begin{cases}
        g^k / \| g^k \|, & g^k \neq 0, \nonumber \\
        0, & g^k = 0, \nonumber
    \end{cases}
  \end{equation*}
  in this case all statements about the convergence of the method remain valid, as in the problem without constraints.

  We define 
	
	$ D = \{ x \> | \> h(x) \leq 0 \} $ as the range of feasible vectors, 
	
	$ X^* = \{ x \> | \> x \in D, 0 \in G(x) \} $  the set of stationary points of the problem, 
	
	$ F^* = \{ f(x) \> | \> x \in X^* \} $ the set of values $ f(x) $ on the range $ X^* $, 
	
	$ X_h^* = \{ x \> | \> x \notin D, 0 \in G_h(x) \} $  the set of stationary points of the problem outside the range $ D $, 
	
	$ H^* = \{ h(x) \> | \> x \in X_h^* \} $ the set of values $ h(x) $ on the range $ X_h^* $.

  Let us formulate results about the convergence of the method. The following lemma establishes the boundedness of the sequence $ \{ x^k \}_{k=0}^{\infty} $.
\begin{lemma}
\label{lem:10.1}
	Let there exist a number $ c > \max \{ 0, h(x^0) \} $ that does not belong to $ H^* $. Then for any $ d > c $ there exists $ \bar{R} $ such that for any sequence $ \{ \rho_k \}_{k=0}^{\infty} $ satisfying (10.8) with $ R < \bar{R} $, the sequence $ \{ x^k \}_{k=0}^{\infty} $ launched from $ x^0 $ according to (10.7), do not leave the bounded region $ \{ x \> | \> h(x) \leq d \} $.
\end{lemma}

  Since outside the admissible region $ D $ the algorithm (10.7) minimizes the function $ h(x) $, Lemma 10.1 is similar to Lemma 9.1 and is proved exactly in the same way. Here also remarks 9.1, 9.2 concerning the selection of the value of $ R $ remain in force. In particular, the algorithm of the form
  \begin{equation*} 
    x^{k+1} =  
    \begin{cases}
        x^k - \rho_k g^k, \> \> g^k \in G(x^k), & \text{if} \> \> h(x^k) < d, \nonumber \\
        x^0, & \text{if} \> \> h(x^k) \geq d, \nonumber
    \end{cases}
  \end{equation*}
where  $ d > \max \{ 0, h(x^0) \} $, $ \{ \rho_k \}_{k=0}^{\infty} $ satisfy (10.8) with arbitrary $ R, $ and the interval $ [\max \{ 0, h(x^0) \}, d) $ does not lie entirely in $ H^* $, no more than a finite number of times goes out of the bounded region $ \{ x \> | \> h(x) \leq d \} $.
\begin{theorem}
\label{th:10.1}
	Suppose that the sequence $ \{ x^k \}_{k=0}^{\infty} $ is bounded. Then:
  \begin{enumerate}
    \item The minimal limit points $ \{ x^k \}_{k=0}^{\infty} $ (in the value of $ h(x)$) are proper to the set $ D \cup X_h^* $, and the interval $ [ \underline{\lim}_{k \to \infty} h(x^k), \overline{\lim}_{k \to \infty} h(x^k) ] $ is embedded in the set $ H^* \cup \{ h \in {R}^1 \> | \> h \leq 0 \} $;
    \item If there are limit points $ \{ x^k \}_{k=0}^{\infty} $ belonging to $ D $, then the minimal of them (by the value of $ f(x)$) belong to $ X^* $.
  \end{enumerate}
\end{theorem}
\begin{theorem}
\label{th:10.2}
	Let $ \{ x^k \}_{k=0}^{\infty} $ be bounded, the sets $ H^* $ and $ F^* $ contain no intervals. Then:
  \begin{enumerate}
    \item either all limit points of $ \{ x^k \}_{k=0}^{\infty} $ do not belong to the admissible region $ D $, but belong to $ X_h^* $, and there exists a limit $ \lim_{k \to \infty} h(x^k) > 0 $;
    \item or all limit points of $ \{ x^k \}_{k=0}^{\infty} $ belong to the admissible region $ D $, and then they all belong to the connected subset $ X^* $, and the sequence $ \{ f(x^k) \}_{k=0}^{\infty} $ has a limit.
  \end{enumerate}
\end{theorem}
\begin{corollary}
\label{cor:10.1}
	Let the sequence $ \{ x^k \}_{k=0}^{\infty} $ be bounded and the function $ h(x) $ has no stationary points outside the admissible region, i.e.
  \begin{equation*} 
    \{ x \> | \> 0\in G_h(x) \} \subset D.
  \end{equation*}
  Then all limit points of $ \{ x^k \}_{k=0}^{\infty} $ belong to the admissible region, and hence,
  \begin{enumerate}
    \item the minimal of them (by the value of $ f(x)$) belong to $ X^* $, and 
		$$ [ \underline{\lim}_{k \to \infty} f(x^k), \overline{\lim}_{k \to \infty} f(x^k) ] \subset F^*; $$
    \item if $ F^* $ contains no intervals, then all limit points $ \{ x^k \}_{k=0}^{\infty} $ belong to a connected subset of $ X^* $ and there exists a limit $ \lim_{k \to \infty} f(x^k) $.
  \end{enumerate}
	\end{corollary}

  Theorems 10.1, 10.2 and corollary 10.1 generalize Theorem 9.1 on convergence of the generalized gradient method to the minimization problem with constraints. Before proving the theorems, let us prove the minimizing property of the generalized gradient method under a constraint.
\begin{lemma}
\label{lem:10.2}
	Let the sequence of initial points $ \{ x^s \}_{s=0}^{\infty} $ converge to $ x = \lim_{s \to \infty} x^s $. Let us run the generalized gradient method (10.7) from starting points $ x^s $ with steps $ \rho_s^k $ ($ k \geq k_s $) up to some point $ n_s $. For each $ s $ we obtain sequences $ \{ x_s^k \}_{k=k_s}^{n_s} $:
  \[
    \begin{matrix}
      x_s^{k_s} = x^s, & x_s^{k+1} = x_s^k - \rho_s^k g_s^k \\
      g_s^k \in G(x_s^k), & k = k_s, k_s + 1, \ldots.
    \end{matrix}
  \]
  We denote 
  \[
    \begin{matrix}
      \rho_s = \sup_{k \geq k_s} \rho_s^k, & \sigma_s = \sum_{k=k_s}^{n_s - 1} \rho_s^k.
    \end{matrix}
  \]

  If $ 0 \notin G(x) $, $ \sigma_s \geq \sigma > 0 $ and $ \lim_{s \to \infty} \rho_s = 0 $, then there exists some constant $ \bar{\varepsilon} > 0 $ such that for any $ \varepsilon \in (0, \bar{\varepsilon}] $ there are indexes $ l_s \geq k_s$, such that $ || x_s^k - x || \leq \varepsilon $ for $ k \in [k_s, l_s]$, and
  \begin{enumerate}
    \item $ h(x) = \lim_{s \to \infty} h(x^s) > \overline{\lim}_{s \to \infty} h(x_s^{l_s}) $, $ h(x) > 0 $;
    \item $ f(x) = \lim_{s \to \infty} f(x^s) > \overline{\lim}_{s \to \infty} f(x_s^{l_s}) $, $ h(x) < 0 $;
    \item $ f(x) = \lim_{s \to \infty} f(x^s) > \overline{\lim}_{s \to \infty} f(x_s^{l_s}) $, $ 0 \geq \overline{\lim}_{s \to \infty} h(x_s^{l_s}) $, $ h(x) = 0 $.
  \end{enumerate}
	\end{lemma}

{\it  P r o o f.} Assertions 1), 2) of the lemma follow for  sufficiently small 
$\bar{\varepsilon}$ from Lemma 9.3.

  Consider the case $ h(x) = 0 $. Denote $ \gamma_1 = \rho (0, G(x)) $. It is easy to check that the mapping $ G $ is locally bounded and closed, and hence semi-continuous from above. Therefore there exists an $\varepsilon_1$-neighborhood of point $ x $ such that for all $ g \in G(y), (\| y - x \| \leq \varepsilon_1) $ the inequality holds
  \begin{equation*} 
    \rho (g, G(x)) \leq \frac{\gamma_1}{2}.
  \end{equation*}
  We denote 
  \begin{equation*} 
    \Gamma_1 = \sup \{ \| g \| \> | \> g \in G(y), \> \| y - x \| \leq \varepsilon_1 \} < \infty.
  \end{equation*}

  Due to the generalized differentiability of $ f(x) $ and $ h(x) $, for any constants $ c_f, c_h > 0 $ there is such an $ \varepsilon_2 \leq \varepsilon_1 $ that the relations are valid:
\begin{equation}
\label{eqn:10.9}
    \begin{matrix} 
      f(y) = f(x) + \left<g_f, y - x\right> + o_f (x, y, g_f), \\
      h(y) = h(x) + \left<g_h, y - x\right> + o_h (x, y, g_h),
    \end{matrix}
  \end{equation}
  where $ o_f (x, y, g_f) \leq c_f \| y - x \| $ and $ o_h (x, y, g_h) \leq c_h \| y - x \| $ if $ \| y - x \| \leq \varepsilon_2 $ and $ g_f \in G_f (y), g_h \in G_h (y) $. We will specify constants $ c_f $ and $ c_h $ later.

  Let us define
  \begin{equation*} 
    \bar{\varepsilon} = \min \left\{ \varepsilon_2, \sigma \frac{\gamma_1}{2} \right\}.
  \end{equation*}
  Let us fix an arbitrary $ \varepsilon \in (0, \bar{\varepsilon}]$. Let $ s $ satisfy the conditions $ \| x^s - x \| \leq \frac{\varepsilon}{4} $ and $ \rho_s \leq \frac{\varepsilon}{4 \Gamma_1} $.

  Let us show that the sequences $ \{ x_s^k \}_{k=k_s}^{n_s} $ exit from the $ \frac{\varepsilon}{2} $-neighborhoods of points $ x^s $. Suppose the contrary: $ || x_s^k - x^s || \leq \frac{\varepsilon}{2} $ for each $ k \in [k_s, n_s] $. Then, on the one hand, we have
  \begin{equation*} 
    \| x_s^{n_s} - x^s \| = \left\| \sum_{k=k_s}^{n_s - 1} \rho_s^k g_s^k \left/\right. \sum_{k=k_s}^{n_s - 1} \rho_s^k \right\| \sum_{k=k_s}^{n_s - 1} \rho_s^k.
  \end{equation*}
  Here for each $ k \in [k_s, n_s] $
  \begin{equation*} 
    \| x_s^k - x \| = \| x_s^k - x^s \| + \| x^s - x \|| \leq \frac{3 \varepsilon}{4};
  \end{equation*}
  thus
  \begin{equation*} 
    \rho (g_s^k, G(x)) \leq \frac{\gamma_1}{2}.
  \end{equation*}
  Since $ G(x) $ is convex and $ \rho (0, G(x)) = \gamma_1 $, then all $ g_s^k (k \in [k_s, n_s]) $ lie in a half-space departing from zero by $ \frac{\gamma_1}{2} $. The same holds for convex combinations of vectors $ g_s^k (k \in [k_s, n_s]) $. Thus, on the other hand, we have
  \begin{equation*} 
    \| x_s^{n_s} - x^s \| \geq \frac{\gamma_1}{2} \sum_{k=k_s}^{n_s - 1} \rho_s^k \geq \frac{\sigma \gamma_1}{2}.
  \end{equation*}

  The resulting contradiction proves that the sequences $ \{ x_s^k \}_{k=k_s}^{n_s} $ come out of the $ \frac{\varepsilon}{2} $-neighborhoods of points $ x_s $.

  Let us denote by $ m_s $ the moment of the first such exit, i.e.
  \begin{equation*} 
    m_s = \max \{ m \> | \> \| x_s^k - x^s \| \leq \frac{\varepsilon}{2} \> \> \> \forall k \in [k_s, m), m \in [k_s, n_s] \},
  \end{equation*}
  then
  \begin{equation*} 
    \| x_s^{m_s} - x \| \leq \rho_s^{m_s - 1} \| g_s^{m_s - 1} \| + \| x_s^{m_s - 1} - x^s \| + \| x^s - x \| \leq \varepsilon.
  \end{equation*}

  Let us define
  \begin{equation*} 
    g_s^k = \lambda_{fs}^k g_{fs}^k + \lambda_{hs}^k g_{hs}^k, \> \> \> g_{fs}^k \in G_f (x_s^k), \> \> \> g_{hs}^k \in G_h (x_s^k), \;\;\;\lambda_{fs}^k + \lambda_{hs}^k = 1.
  \end{equation*}

  For $ k \in [k_s, m_s] $ we substitute $ x_s^k $ and $ g_{fs}^k \in G_f (x_s^k) $ in the equation (10.9):
  \begin{eqnarray}
          f(x_s^k) &=& f(x) + \left<g_{fs}^k, x_s^k - x\right> + o_f (x, x_s^k, g_{fs}^k) 
					\nonumber\\
      &\leq& f(x) + \left<g_{fs}^k, x_s^k - x^s\right> + \left<g_{fs}^k, x^s - x\right> + c_f \| x_s^k - x \| \nonumber\\
      &\leq& f(x) + \left<g_{fs}^k, x_s^k - x^s\right> + c_f \| x_s^k - x^s \| + (\Gamma_1 + c_f) \| x^s - x \|.\nonumber\\ \label{eqn:10.10}
    \end{eqnarray}
  Similarly,
  \begin{equation}\label{eqn:10.11}
    h(x_s^k) \leq h(x) + \left<g_{hs}^k, x_s^k - x^s\right> + c_h \|x_s^k - x^s\| + (\Gamma_1 + c_h) \| x^s - x \|.
  \end{equation}

  Now consider scalar products:
  \begin{equation*} 
    \left<g_s^k, x_s^k - x^s\right> = - \left<g_s^k, \sum_{r=k_s}^{k-1} \rho_s^r g_s^r / \sum_{r=k_s}^{k-1} \rho_s^k\right> \sum_{r=k_s}^{k-1} \rho_s^k,
  \end{equation*}
  where $ g_s^r \in G(x_s^r) \;(r \in [k_s, m_s]) $. For their estimation we can apply Lemma 9.2, where we need to define
  \begin{equation*} 
    P = co \{ G(y) \> | \> \| y - x \| \leq \varepsilon_1 \};
  \end{equation*}
  \[
    \begin{matrix}
      p^r = g_s^r, & r \in [k_s, m_s], & \Gamma = \Gamma_1, & \gamma = \frac{\gamma_1}{2}.
    \end{matrix}
  \]

  The conditions are satisfied:
  \[
    \begin{matrix}
      \sum_{r=k_s}^{m_s - 1} \rho_s^r \geq \frac{\| x_s^{m_s} - x^s \|}{\Gamma_1} \geq \frac{\varepsilon}{2 \Gamma_1}, & \lim_{s \to \infty} \rho_s = 0.
    \end{matrix}
  \]

  Now, according to Lemma 9.2, for sufficiently large $ s $ one can find such indices $ l_s \leq m_s $ such that
  \begin{equation}\label{eqn:10.12}
    \left<g_s^{l_s}, \sum_{k=k_s}^{l_s - 1} \rho_s^k g_s^k / \sum_{k=k_s}^{l_s - 1} \rho_s^k\right> \geq \frac{\gamma_1^2}{8}, \;\;\;\;\;\; \sum_{k=k_s}^{l_s - 1} \rho_s^k \leq \frac{\varepsilon \gamma_1}{12 \Gamma_1^2}.
  \end{equation}

  Recall that we consider the case $ h(x) = 0 $. Let us show that if the constant $ c_h $ is chosen appropriately then for sufficiently large $c$ holds $ h(x_s^{l_s}) \leq 0 $, which proves the second part of statement 3) of the lemma. Indeed, suppose that $ h(x_s^{l_s}) > 0 $. Then $ g_s^{l_s} = g_{hs}^{l_s} $, and from (10.11) we obtain
  \begin{equation*}
    0 < h(x_s^{l_s}) \leq - \frac{\gamma_1^2}{8} \sum_{k=k_s}^{l_s - 1} \rho_s^k + c_h \Gamma_1 \sum_{k=k_s}^{l_s - 1} \rho_s^k + (\Gamma_1 + c_h) \| x^s - x \|.
  \end{equation*}

  If $ c_h \leq \gamma_1^2 / (16 \Gamma_1) $, then
  \begin{equation*}
    f(x_s^{l_s}) \leq f(x) - \frac{\gamma_1^2}{192 \Gamma_1^2} \gamma_1 \varepsilon + (\Gamma_1 + c_f) || x^s - x ||,
  \end{equation*}
that is impossible for sufficiently large $s$. This contradiction proves the second part of statement 3) of the lemma.

Now we have two possibilities for sufficiently large $s$:
$h(x_s^{l_s})<0$ and $h(x_s^{l_s})=0$. Under $h(x_s^{l_s})<0$ the equality
$g_s^{l_s}=g_{fs}^{l_s}$ holds true; so the inequality (10.10) can be rewritten at $k=l_s$ as
\[
f(x_s^{l_s})\le f(x)-\frac{\gamma_1^2}{8}\sum_{k=k_s}^{l_s-1}\rho_s^k
+c_f\Gamma_1\sum_{k=k_s}^{l_s-1}\rho_s^k+(\Gamma_1+c_f)\|x^s-x\|.
\]

If $c_f\le\gamma_1/(16\Gamma_1)$, then
\[
f(x_s^{l_s})\le f(x)-\frac{\gamma_1^2}{192\Gamma_1^2}+(\Gamma_1+c_f)\|x^s-x\|.
\]

  Hence, for $ s $ such that $ h(x_s^{l_s}) < 0 $, it follows that statement 3) of the lemma is true. 
	
	Consider the case $ h(x_s^{l_s}) = 0 $. From (10.11) at $ k = l_s $, we obtain
  \begin{equation*}
    \left<g_{hs}^{l_s}, \sum_{k=k_s}^{l_s - 1} \rho_s^k g_s^k\right> \leq c_h \left\| \sum_{k=s}^{l_s - 1} \rho_s^k g_s^k \right\| + (\Gamma_1 + c_h) \| x^s - x \|.
  \end{equation*}

  From the estimate (10.12), we have
  \begin{eqnarray}
    \lambda_{fs}^{l_s} \left<g_{fs}^{l_s}, \sum_{k=k_s}^{l_s - 1} \rho_s^k g_s^k\right> 
		&\geq& 
		\frac{\gamma_1^2}{8} \sum_{k=k_s}^{l_s - 1} \rho_s^k - \lambda_{hs}^{l_s} \left<g_{hs}^{l_s}, \sum_{k=k_s}^{l_s - 1} \rho_s^k g_s^k\right> \nonumber\\
       &\geq& \left(\frac{\gamma_1^2}{8} - c_h \Gamma_1\right) \sum_{k=k_s}^{l_s - 1} \rho_s^k \nonumber\\
			&&- (\Gamma_1 + c_h) \| x^s - x \|.\label{eqn:10.13}
 \end{eqnarray}

  Now find a suitable value for constant $ c_h $; all previous reasoning is true for $ c_h = \frac{\gamma_1^2}{16 \Gamma_!} $. If s is large enough, the right and left parts of (10.13) are positive, so given that $ 0 \leq \lambda_{fs}^{l_s} \leq 1 $, we obtain
  \begin{equation*}
    \left<g_{fs}^{l_s}, \sum_{k=k_s}^{l_s - 1} \rho_s^k g_s^k\right> \geq \frac{\gamma_1^2}{16} \sum_{k=k_s}^{l_s - 1} \rho_s^k - (\Gamma_1 + c_h) \| x^s - x \|.
  \end{equation*}

  Substituting this estimate into (10.10), we obtain for $ k = l_s $ and sufficiently large $ s $ it holds
  \begin{equation*}
    f(x_s^{l_s}) \leq f(x) - \frac{\gamma_1^2}{16} \sum_{k=k_s}^{l_s - 1} \rho_s^k + c_f \sum_{k=k_s}^{l_s - 1} \rho_s^k + (2 \Gamma_1 + c_f + c_h) \| x^s - x \|.
  \end{equation*}

  Now let's select the value of $ c_f $ ($ c_h $ was specified earlier). All the reasoning above is true for $ c_f = \frac{\gamma_1^2}{32} $. Then
  \begin{equation*}
    f(x_s^{l_s}) \leq f(x) - \frac{\gamma_1^2}{384 \Gamma_1^2} \gamma_1 \varepsilon + (2 \Gamma_1 + c_f + c_h) \| x^s - x \|.
  \end{equation*}

  Hence, for $ s $ such that $ h(x_s^{l_s})=0 $, the first part of the statement of the lemma also follows. For small $ s $ the indices of $ l_s $ are arbitrary. The lemma is proved.

  Having for the method (10.7) the minimizing property in the form of Lemma 10.2 and for bounded $ \{ x^k \}_{k=0}^{\infty} $ the property $ \lim_{k \to \infty} \| x^{k+1} - x^k \| = 0 $, we can deduce the validity of convergence theorems 10.1, 10.2 from general results on the convergence of algorithms in $\S$ 8.

{\it  P r o o f of Theorem 10.1.} Outside the admissible region $ D $, the method (10.7) minimizes the function $ h(x) $ as in the problem without constraints, with $ h(x) $ as the Lyapunov function and the set of solutions as the set $ D \cup X_h^* $. The function $ h(x) $ on $ D \cup X_h^* $ runs over values from the set $ \{ h \in E_1 \> | \> h \leq 0 \} \cup H^* $. Note that $ H^* $ is closed. Therefore the first part of Theorem 10.1 follows from Theorem 9.1 of convergence of the generalized gradient method in the problem without constraints, or, in other words, from (10.8), statement 1) of Lemma 10.2 and general Theorem 8.9.

  The second part of Theorem 10.1 follows from conditions 2), 3) of Lemma 10.2. Indeed, suppose the contrary. Let
  \begin{equation*}
    x = \lim_{s \to \infty} x^{k_s}
  \end{equation*}
  is one of the minimal (in value of $ f(x) $) limit points $ \{ x^k \}_{k=0}^{\infty} $ in a set $ D $ that does not belong to $ X^* $. We can assume that the algorithm (10.7) starts repeatedly from points $ x^{k_s} \> (s = 0, 1, \ldots) $; in this case we obtain sequences $ \{ x^k \}_{k \geq k_s} $, and, in this case, Lemma 10.2 is applicable.

  If $ x $ is an interior point of $ D $, it follows from assertion 2) of Lemma 10.2 that if $ \varepsilon $ is small enough, there exists an even smaller interior limit point $ \{ x^k \}_{k=0}^{\infty} $ of $ f $, which contradicts the minimality of $ x $.

  If $ x $ lies on the boundary $ D $, then by virtue of statement 3) of Lemma 10.2 again there exists a smaller admissible limit point $ \{ x^k \}_{k=0}^{\infty} $ of $ f $, which again contradicts the minimality of $ x $. Theorem 10.1 is proved.

{\it  P r o o f of Theorem 10.2.} Theorem 10.1 shows that the interval 
$$ \left[ \underset{k \to \infty}{\underline{\lim}} h(x^k), 
\underset{k \to \infty}{\overline{\lim}} h(x^k) \right] $$ 
is embedded in the set
  \begin{equation*}
    H^* \cup \{ h \in E_1 | h \leq 0 \}.
  \end{equation*}

  It implies that either the interval $ [ \underline{\lim}_{k \to \infty} h(x^k), \overline{\lim}_{k \to \infty} h(x^k) ] $ is fully embedded in the set $ \{ h \in E_1 | h \leq 0 \} $ and then the boundary elements $ \{ x^k \}_{k=0}^{\infty} $ belong to the $ D $, or there is a limit
  \begin{equation*}
    \lim_{k \to \infty} h(x^k) > 0,
  \end{equation*}
  and then the limit points $ \{ x^k \}_{k=0}^{\infty} $ do not belong to the admissible region $ D $. In this latter case all limit points $ \{ x^k \} $ belong to $ X_h^* $. Indeed, if this is not the case, then there exists a limit point
  \begin{equation*}
    \lim_{s \to \infty} x^{k_s} = x \notin X^*.
  \end{equation*}
  We can assume that  algorithm (10.7) as if repeatedly starts from points $ x^{k_s} \> (s = 0, 1, ...) $ and obtains sequences $ \{ x^k \}_{k \geq k_s} $. Lemma 10.2 applies to them. Minimizing property 1) of this lemma, in particular, denies the possibility of $ \lim_{k \to \infty} h(x^k) $. The resulting contradiction proves the first statement of Theorem 10.2.

  Let us prove the second statement of the theorem. As shown above, either all limit points of $ \{ x^k \}_{k=0}^{\infty} $ do not belong to $ D $, or all of them belong to $ D $. Suppose that they all belong to $ D $. Let us consider the values of $ f(x) $ at points $ x^k $, whether they are admissible or invalid. We may assume that method (10.7) minimizes the function $ f(x) $ on $ D $, and that the minimizing property (3) of Lemma 10.2 holds, and that $ \lim_{k \to \infty} \| x^{k+1} - x^k \| = 0 $ due to the boundedness of $ \{ x^k \} $. Then statement 2) of Theorem 10.2 follows from general Theorem 8.9.

\smallskip
\label{Sec.10.2}
\textbf{2. Multiple inequality constraints problem.} Considering  method (10.7) for solving problem (10.5), (10.6), we see that it first seems to minimize the function $ h(x) $ and then the function $ f(x) $ on the set $ \{ x \> | \> h(x) \leq 0 \} $. Nothing prevents us from increasing the number of sequentially minimized functions. Thus we come to the method of solving a mathematical programming problem with several constraints as a multi-criteria optimization problem.

  So, let's consider a mathematical programming problem:
\begin{equation}\label{eqn:10.14}
    \begin{matrix} 
      f_{m+1} (x) \to \min_x, & x \in E_n
    \end{matrix}
  \end{equation}
subject to the constraints
\begin{equation}\label{eqn:10.15}
    \begin{matrix} 
      f_j (x) \leq c_j, j = 1, 2, ..., m
    \end{matrix}
 \end{equation}
where we assume that constraints $ f_j (x) \leq c_j $ are ordered by their priority.

  A constraint with a smaller index will be considered more important than a constraint with a larger index, and hence we should achieve its fulfillment first, i.e. minimize the corresponding function $ f_j(x) $ to the level $ c_j $. Only when all the constraints are satisfied, we proceed to minimize the target function $ f_{m+1} (x) $, which can also be associated with some level $ c_{m+1} \geq - \infty $, i.e., essentially, do not distinguish between the constraint and the target function.

  At the same time, if there are many constraints, we may not achieve less important constraints, and the optimization process will stop at some constraint $ i^* \> (1 \leq i^* \leq m + 1) $. The corresponding point $ x^* $ will be the solution of the problem
\begin{equation}\label{eqn:10.16}
    f_{i^*} (x) \to \min_x
\end{equation}
subject to constraints
\begin{equation}\label{eqn:10.17}
    \begin{matrix} 
      f_j (x) \leq c_j, 1 \leq j < i^*
    \end{matrix}
 \end{equation}

  Thus, a solution of a possibly incompatible problem (10.14), (10.15) with ordered constraints is naturally understood as such an optimal solution $ x^* $ of a joint mathematical programming problem of the form (10.16), (10.17), that
  \begin{equation*}
    f_{i^*} (x^*) > c_{i^*}
  \end{equation*}
If the set $ x \> | \> f_1(x) \leq c_1 $ is empty, then we can assume that the problem (10.14), (10.15) has no solutions at all.

  Introducing order on a constraint set is one way to consider incompatible mathematical programming problems. 
	
	For problems with an ordered set of constraints we can also introduce a utility function. In this regard, it is convenient for us to refer to maximization problems.

  Let the maximization function be $ \phi_1 (x) $ on a set $ X \subset E_n $ to exceed the level $ c_1 $, and another maximization function $ \phi_2 (x) $ on a set $ \{ x \in X \> | \> \phi_1 (x) \geq c_1 \} $ to exceed the level $ c_2 $, and $ \phi_3 (x) \{ x \in X \> | \> \phi_1 (x) \geq c_1, \phi_2 (x) \geq c_2 \} $, etc.

  Note that in the series $ \phi_1 (x), ..., \phi_j (x), ..., \phi_m (x) $ the same function can be repeated several times under different indices $ j $, but with larger values of levels $ c_j $, which will reflect the decreasing importance of larger values of this function. We do not discuss here the problem of ordering target functions by importance.

  For the problem under consideration we construct a utility function $ \Phi (\phi) $, where $ \phi = (\phi_1, ..., \phi_m) \in E_m $. This function should be the larger the more constraints are fulfilled and the larger the value of the function in the first unfulfilled constraint. We will assume that $ \phi_j (x) \geq 0\; (j = 1, ..., m) $ when $ x \in X $. The following function has these properties:
  \[
    \Phi (\phi) = \begin{cases}
        \phi_1,                                & 0 \leq \phi_1 < c_1, \\
        \phi_2 + c_1,                          & c_1 \leq \phi_1, \> 0 \leq \phi_2 < c_2, \\
        \phi_3 + c_1 + c_2,                    & c_1 \leq \phi_1, c_2 \leq \phi_2, \> 0 \leq \phi_3 < c_3, \\
        & \cdots \\
        \phi_m + c_1 + c_2 + \cdots + c_{m-1}, & c_1 \leq \phi_1, \cdots, c_{m-1} \leq \phi_{m-1}, \> 0 \leq \phi_m < c_m, \\
        \> \> \> \> \> \> \> \> \> \> \> c_1 + c_2 + \cdots + c_{m}, & c_1 \leq \phi_1, \cdots, c_m \leq \phi_m.
    \end{cases}
  \]

  If there is a utility function, then the problem of multi-criteria optimization can be solved "directly" by substituting the criteria $ \phi_j (x) $ in the utility function $ \Phi (\phi) $ and maximizing the resulting function on $ x $. Everything depends on the properties of the function $ \Phi (\phi) $ and the information we have about it.

  In this particular case the function $ \Phi (\phi) $ is discontinuous, but it is also quasi-concave, i.e. sets $ \{ \phi \in E_m \> | \> \Phi (\phi) \geq c \} $ are convex (they have a representation  $ \{ \phi \in E_m \> | \> \phi_1 \geq c_1, ..., \phi_k \geq c_k, \phi_{k+1} \geq c \} $, where $ c_1 + c_2 + ... + c_k \leq c < c_1 + c_2 + ... + c_{k+1} $).

  If we substitute quasi-concave functions $ \phi(x) $ into $ \Phi(\phi) $, then we get a quasi-concave function $ \Phi(\phi (x)) $, due to the convexity of the set
  \begin{equation*}
    \{ x \> | \> \Phi(\phi (x)) \geq c \} = \{ x \> | \> \phi_1 (x) \geq c_1, ..., \phi_k (x) \geq c_k, \phi_k (x) \geq c_k, \phi_{k+1} (x) \geq c \},
  \end{equation*}
  where $ c_1 + c_2 + ... + c_k \leq c < c_1 + c_2 + ... + c_{k+1} $.

  The existing methods for maximizing quasi-convex functions [75] are justified under the assumption of their continuity. We consider a more general situation when functions $ \phi_j(x) $ are generalized differentiable, and construct a method for minimizing the function $ \Phi (\phi (x)) $.

  In order to maximize $ \Phi $ at the point $ \phi \in E_m $, it is necessary to maximize the first function $ \phi_j $ in the series $ \phi_1, \phi_2, ..., \phi_m $, which has not reached its level $ c_j $. With respect to the minimization problem (10.14), (10.15) this means that we have to minimize the function $ f_i(x) $ of the first of the unfulfilled constraints. Therefore consider the following method for solving problems (10.14), (10.15):
  \[
    \begin{matrix}
      x^0 \in E_n, \> \> \> x^{k+1} = x^k - \rho_k g^k, \> \> \> g^k \in G_{ik}(x^k), \;\;\;k = 0, 1, ... \\
      i_k = \min \{ j \> | \> 1 \leq j \leq m + 1, f_j(x^k) > c_j \;(c_{m+1} = - \infty) \}.
    \end{matrix}
  \]

  We may even consider a somewhat more general method:
  \begin{equation}\label{eqn:10.18}
      x^0 \in E_n, \;\;\; x^{k+1} = x^k - \rho_k g^k, \;\;\; g^k \in G(x^k),
 \end{equation}
\begin{equation}\label{eqn:10.19}
       0 \leq \rho_k \leq R, \;\;\; \lim_{k \to \infty} \rho_k = 0, \;\;\; \sum_{k=0}^{\infty} \rho_k = \infty,
\end{equation}
  where
\begin{equation}\label{eqn:10.20}
    G(x) = \text{co} \{ G_{i(x)} (x); \> G_j(x), \> j \in J(x) \}
  \end{equation}
\begin{equation}\label{eqn:10.21}
    i(x) \equiv \min \{ j \> | \> 1 \leq j \leq m + 1, \> f_j(x) \geq c_j, \> (c_{m+1} = - \infty) \}
  \end{equation}
\begin{equation}\label{eqn:10.22}
    J(x) \equiv \{ j \> | \> 1 \leq j < i(x), \> f_j(x) = c_j \},
  \end{equation}
  where $ G_j (x) $ is the pseudo-gradient mapping of the function $ f_j (x) $ of problems (10.14) and (10.15) ($ j = 1, 2, ..., m + 1 $).

  Let us formulate an analogue of Lemma 10.2 for this method. First, let us prove the semi-continuity from above of the mapping $ G(x) $.
\begin{lemma}
\label{lem:10.3}
	Function $ i(x) $ and its mapping $ x \to J(x), x \to i(x) \cup J(x) $ are semi-continuous from above.
	\end{lemma}

{\it   P r o o f.} According to the definition (10.21) of the index $ i(x) $, it takes place $ f_{i(x)} (x) > c_{i(x)} $. Due to continuity properties of $ f_{i(x)} $ the inequality $ f_{i(x)} (y) > c_{i(x)} $ holds for $ y $ from some neighborhood of $ x $. Here $ i(y) \leq i(x) $ and the function $ i(x) $ is semi-continuous from above.

  If $ f_j(x) \neq c_j $, then $ f_j(y) \neq c_j $ for $ y $ from some neighborhood of $ x $. It follows that $ J(y) \subset J(x) $, i.e. the mapping $ J(x) $ is semi-continuous from above (and closed). Let us prove that the mapping $ x \to i(x) \cup J(x) $ is closed; in this case this is equivalent to its semi-continuity from above. By virtue of the closeness of $ J(x) $ already proved, it suffices to show that if $ x^k \to x $, $ i(x^k) \in i $, then $ i \in i(x) \cup J(x) $. According to the above, $ i \leq i(x) $. Since $ f_{i(x^k)} (x^k) > c_{i(x^k)} $, then for big enough $ k $ values $ f_i (x^k) > c_i $ and due to the continuity of $ f_i (x) $ it takes place $ f_i (x) \geq c_i $. It follows that either $ i = i(x) $ or $ i \in J(x) $. The lemma is proved.
\begin{lemma}
\label{lem:10.4}
	The mapping $ x \to G(x) $, where $ G(x) $ is constructed according to (10.20)-(10.22), is locally bounded, closed, and therefore semi-continuous from above.
\end{lemma}
{\it  P r o o f.} Recall that the pseudo-gradient mappings $ G_j (x) $ of functions $ f_j (x) $ are locally bounded and closed ($ j = 1, 2, ..., m + 1 $). Therefore $ G(x) $ is also locally bounded. Let us prove the closeness of $ G $. Assume $ x^k \to x $, $ g^k \in G(x^k) $ and $ g^k \to g $, let us show that $ g \in G(x) $.

  Let us represent
  \begin{equation*}
    g^k = \lambda_{i(x^k)}^k g_{i(x^k)}^k + \sum_{j \in J(x^k)} \lambda_j^k g_j^k, \;\;\; \lambda_{i(x^k)}^k, \;\;\; \lambda_j^k \geq 0,
  \end{equation*}
  \begin{equation*}
    \lambda_{i(x^k)}^k + \sum_{j \in J(x^k)} \lambda_j^k = 1, \;\;\; g_{i(x^k)}^k \in G_{i(x^k)} (x^k), \;\;\; g_j^k \in G_j (x^k).
  \end{equation*}

  For sufficiently large $ k $, according to Lemma 10.3, $ J(x^k) \subset J(x) $; therefore it is possible to represent
  \begin{equation*}
    g^k = \lambda_{i(x^k)}^k g_{i(x^k)}^k + \sum_{j \in J(x^k)} \lambda_j^k g_j^k,
  \end{equation*}
  where
  \begin{equation*}
    \lambda_j^k = 0 \> \> \> \text{and} \> \> \> g_j^k \in G_j(x) \> \> \> \text{with} \> \> \> j \in J(x^k).
  \end{equation*}

  Due to the finiteness of the set of functions $ f_j(x) $, there exists a subset of indices $ \{ k_s \} $ such that $ i(x^{k_s}) = i $, $ J(x^{k_s}) \subset J(x) $, and
  \[
    \begin{matrix}
       \lambda_{i(x^{k_s})}^{k_s} = \lambda_i^{k_s} \to \lambda_i \geq 0, & \lambda_j^{k_s} \to \lambda_j, & \lambda_i + \sum_{j \in J(x)} \lambda_j = 1 \\
       g_{i(x^{k_s})}^{k_s} = g_i^{k_s} \to g_i, & g_j^{k_s} \to g_j, & j \in J(x).
    \end{matrix}
  \]

  Due to the closeness of mappings $ G_j $ ($ j = 1, 2, ..., m + 1 $) $ g_i \in G_i(x) $, $ g_j \in G_j(x) $ ($ j \in J(x) $). According to Lemma 10.3, either $ i = i(x) $ or $ i \in J(x) $; therefore
  \begin{equation*}
    g = \lim_{k \to \infty} g^k = \lambda_i g_i + \sum_{j \in J(x)} \lambda_j, \> \> \> g_j \in G(x).
  \end{equation*}
  Lemma is proven.
\begin{lemma}
\label{lem:10.5}
	Let the sequence of initial points $ \{ x^s \}_{s=0}^{\infty} $ converges to $ x = \lim_{s \to \infty} x^s $. We will run method (10.18) from initial points $x^s $ with steps $ \rho_s^k $ ($ k \geq k_s $) up to some point $ n_s $. For each $ s $, we obtain sequences of $ \{ x_s^k \}_{k=k_s}^{n_s} $:
  \[
    \begin{matrix}
       x_s^{k_s} = x^s, & x_s^{k+1} = x_s^k - \rho_s^k g_s^k, & g_s^k \in G(x_s^k), & k = k_s, k_s + 1, \ldots.
    \end{matrix}
  \]
  Let us define
  \[
    \begin{matrix}
       \rho_s = \sup_{k \geq s} \rho_s^k, & \sigma_s = \sum_{k=k_s}^{n_s - 1} \rho_s^k.
    \end{matrix}
  \]

  If $ 0 \notin G(x) $, $ \sigma_s \geq \sigma > 0 $ and $ \lim_{s \to \infty} \rho_s = 0 $, then there exists $ \bar{\varepsilon} (x) > 0 $ such that for any $ \varepsilon \in (0, \bar{\varepsilon}] $ there are such indices $ l_s $ that for large enough $ s $ $ || x_s^k - x || \leq \varepsilon $ for $ k \in [k_s, l_s] $, and
  \begin{enumerate}
    \item $ f_{i(x)} (x) = \lim_{s \to \infty} f_{i(x)} (x^s) > \overline{\lim}_{s \to \infty} f_{i(x)} (x_s^{l_s}) $, where $ i(x) = \min \{ j \> | \> 1 \leq j \leq m + 1, \> f_j(x) > c_j \> (c_{m+1} = - \infty) $\};
    \item $ c_j \geq \overline{\lim}_{s \to \infty} f_j (x_s^{l_s}) $, $ j \in J(x) $, where $ J(x) = \{ j \> | \> 1 \leq j \leq i(x), \> f_j(x) = c_j $\}.
  \end{enumerate}
	\end{lemma}

{\it  P r o o f.} The proof of this lemma basically repeats the proof of lemma 10.2. Denote by $ \gamma = \rho (0, G(x)) > 0 $. Because $ G(x) $ is semi-continuous from above (Lemma 10.4), there exists an $ \varepsilon_1 $-neighborhood of point $ x $ such that for all $ g \in G(y) $ ($ \| y - x \| \leq \varepsilon_1 $) the expression $ \rho (g, G(x)) \leq \frac{\gamma}{2} $ takes place. Let us denote by
  \begin{equation*}
    \Gamma = \sup \{ \| g \| \> | \> g \in G(y), \> \| y - x \| \leq \varepsilon_1 \} < \infty.
  \end{equation*}

  By Lemma 10.3, there exists an $ \varepsilon_2 $-neighborhood of point $ x $ such that for all $ y \> | \> \| y - x \| \leq \varepsilon_2 $ there is $ J(y) \subset J(x) $, $ i(y) \in i(x) \cup J(x) $.

  By virtue of the generalized differentiability of the functions $f_{j}(x)$, for $a=\gamma^{2} /(32 \Gamma)$ there is such an  $\varepsilon_{3} \le \varepsilon_{2}$ that
    \begin{equation} \label{eqn:10.23}
        f_{j}(y)=f_{j}(x)+\left<g_{j}, y-x\right>+o_{j}\left(x, y, g_{j}\right),
    \end{equation}
 where $\left|o_{j}\left(x, y, g_{j}\right)\right| \le a\|y-x\|$ at $\|y-x\| \le \varepsilon_{8}$ and $g_{j} \in G_{j}(y)$ $(j=1, \ldots, m+1)$. Let us set
    $$
    \bar{\varepsilon}=\min \left\{\varepsilon_{3}, \sigma \gamma / 2\right\} \text {. }
    $$
Fix an arbitrary $\varepsilon \in(0, \bar{\varepsilon}]$. Let $s$ be such that $\left\|x^{s}-x\right\| \le \varepsilon / 4$ and $\rho_{s} \le \varepsilon /(4 \Gamma)$.
    
    As in Lemma 10.2, the sequences $\left\{x_{s}^{k}\right\}_{k=k_{s}}^{n_{3}}$ go out of $(\varepsilon / 2)$-neighborhoods of the points $x_{s}$. The indices $m_{s}$ denote the moments of the first such exit. It holds $\left\|x_{\mathrm{s}}^{k}-x\right\| \le \varepsilon$ with $k \in\left[k_{\mathrm{s}}, m_{\mathrm{s}}\right]$.
    
    For each $s$ represent:
    $$
    \begin{gathered}
    g_{s}^{k}=\sum_{i \in i\left(x_{s}^{k}\right) \cup J\left(x_{s}^{k}\right)} \lambda_{s j}^{k} g_{s j}^{k}, \;\;\;g_{s j}^{k} \in G_{j}\left(x_{s}^{k}\right), \\
    \lambda_{s j}^{k} \ge 0, \;\;\;\sum_{i \in\left(x_{s}^{k}\right) \cup J\left(x_{s}^{k}\right)} \lambda_{s j}^{k}=1.
    \end{gathered}
    $$
    
    For the $k \in\left[k_{s}, m_{s}\right]$, let us substitute $x_{s}^{k}$ and $g_{s j}^{k}$, where
    $$
    j \in i\left(x_{s}^{k}\right) \cup J\left(x_{s}^{k}\right),
    $$
 in relations (10.23):
 \begin{equation} \label{10.24}
    f_{f}\left(x_{s}^{k}\right) \le f_{j}(x)+\left<g_{s j}^{k}, x_{s}^{k}-x^{s}\right>+a\left\|x_{s}^{k}-x^{s}\right\|+(\Gamma+a)\left\|x^{s}-x\right\|.
    \end{equation}
    
    Consider the scalar products:
     $$
    \left<g_{s}^{k}, x_{s}^{k}-x^{s}\right>=-\left<g_{s}^{k}, \sum_{r=k_{s}}^{k-1} \rho_{s}^{r} g_{s}^{r} \sum_{r=k_{s}}^{k-1} \rho_{s}^{r}\right> \sum_{r=k_{s}}^{k-1}\rho_{s}^{r}.
    $$
    
    As in the proof of Lemma 10.2, for sufficiently large $s$ there are indices such $l_{s} \le m_{s}$ that
     \begin{equation} \label{eqn:10.25}
    \left<g_{s}^{l_{s}}, \sum_{k=k_{s}}^{l_{s}-1} \rho_{s}^{k} g_{s}^{k} / \sum_{k=k_{s}}^{t_{s}-1} \rho_{s}^{k}\right> \ge \frac{\gamma^{2}}{8}, \quad \sum_{k=k_{s}}^{l_{s}-1} \rho_{s}^{k} \ge \frac{\varepsilon \gamma}{12 \Gamma^{2}} .
    \end{equation}
    
    Denote $i_{s}=i\left(x_{s}^{l_{s}}\right), J_{s}=J\left(x_{s}^{s_{s}}\right)$. Let us represent
    \begin{equation} \label{eqn:10.26}
    g_{s}^{l_{s}}=\lambda_{s i_{s}}^{l_{s}} g_{s i_{s}}^{l_{s}}+\sum_{j \in J_{s}} \lambda_{s j}^{l_{s}} g_{s j}^{l_{s}}.
    \end{equation}
    
    For $j \in J_{s} \subset J(x)$ and $k=l_{\mathrm{s}}$ from (10.24) we get
     $$
    c_{j}=f_{j}\left(x_{s}^{l_{s}}\right) \le c_{j}+\left<g_{s j^{\prime}}^{l_{s}} x_{s}^{l_{s}}-x^{s}\right>+a \Gamma \sum_{k=k_{s}}^{l_{s}-1} \rho_{s}^{k}+(\Gamma+a)\left\|x^{s}-x\right\|,
    $$
    from where
    \begin{equation} \label{eqn:10.27}
    -\left<g_{s j}^{l_{s}}, x_{s}^{l_{s}}-x^{s}\right> \le a \Gamma \sum_{k=k_{s}}^{l_{s}-1} \rho_{s}^{k}+(\Gamma+a)\left\|x^{s}-x\right\|.
    \end{equation}
    
    From (10.25) - (10.27), we get
     \begin{eqnarray}
    \lambda_{s_{s}}^{l s}\left<g_{s i_{s}}^{l_{s}},\sum_{k=k_{s}}^{l_{s}-1} \rho_{s}^{k} g_{s}^{k}\right> \ge \frac{\gamma^{2}}{8} \sum_{k=k_{s}}^{l_{s}-1} \rho_{s}^{k}-\sum_{j \in J_{s}} \lambda_{s j}^{l_{s}}\left<g_{sj}^{l_{s}}, \sum_{k=k_{s}}^{l_{s}-1} \rho_{s}^{k} g_{s}^{k}\right> \nonumber \\
     \ge \frac{\gamma^{2}}{8} \sum_{k=k_{s}}^{l_{s}-1} \rho_{s}^{k}-\sum_{j \in J_{s}} \lambda_{s j}^{l_{s}}\left(a \Gamma \sum_{k=k_{s}}^{l_{s}-1} \rho_{s}^{k}+(\Gamma+a)\left\|x^{s}-x\right\|\right) \nonumber \\
     \ge \frac{\gamma^{2}}{16} \sum_{k=k_{s}}^{l_{s}-1} \rho_{s}^{k}-(\Gamma+a)\left\|x^{s}-x\right\|.
 \label{eqn:10.28}
    \end{eqnarray}
    
    For sufficiently large $s$ the right and left parts of the (10.28) are positive, $0 \le \lambda_{s i_{s}}^{l_{s}} \le 1$; therefore
    $$
    \left<g_{s i_{s}}^{l_{s}}, \sum_{k=k_{s}}^{l_{s}-1} \rho_{s}^{k} g_{s}^{k}\right> \ge \frac{\gamma^{2}}{16} \sum_{k=k_{s}}^{l_{s}-1} \rho_{s}^{k}-(\Gamma+a)\left\|x^{s}-x\right\| .
    $$
    
    This evaluation is obviously carried out even when $J_{s}=\emptyset$. Substituting it into (10.24) for $k=l_{s}$ and $j=i_{s}$, we get
    \begin{eqnarray}
    f_{i_{s}}\left(x_{s}^{l_{s}}\right) \le f_{i_{s}}(x)-\frac{\gamma^{2}}{16} \sum_{k=k_{s}}^{l_{s}-1} \rho_{s}^{k}+a \Gamma \sum_{k=k_{s}}^{l_{s}-1} \rho_{s}^{k}+2(\Gamma+a)\left\|x^{s}-x\right\| \nonumber \\
    \le f_{i_{s}}(x)-\frac{\gamma^{2}}{32} \sum_{k=k_{s}}^{l_{s}-1} \rho_{s}^{k}+2(\Gamma+a)\left\|x^{s}-x\right\| \nonumber \\
    \le f_{i_{s}}(x)-\frac{\gamma^{2}}{384 \Gamma^{2}} \gamma \varepsilon+2(\Gamma+a)\left\|x^{s}-x\right\|.\label{eqn:10.29}
    \end{eqnarray}
    
    Let us now show that for sufficiently large values of  $s$, it holds
   $$
    f_{j}\left(x_{s}^{l_{s}}\right) \le c_{j}, \quad j \in J(x),
    $$
    from which the second assertion of the lemma will follow. Indeed, suppose the opposite: for arbitrarily large values of  $s$ there are such $j \in J(x)$ that $f_{j}\left(x_{3}^{l_{s}}\right)>c_{j}$. Let's define the indexes,
$$
    j_{s}=\min \left\{j \mid j \in J(x), f_{j}\left(x_{s}^{l_{s}}\right)>c_{j}\right\} .
    $$
It's obvious that
    $$
    i_{s}=i\left(x_{s}^{l_s}\right) \le j_{s}.
    $$
 According to the choice of the $\varepsilon \le \varepsilon_{2}$, it takes place
     $$
    i_{s} \in i(x) \cup J(x).
    $$
     Together with the previous inequality, this gives $i_{s} \in J(x)$. Hence, taking into account that
    $$
    f_{l_{s}}\left(x_{s}^{l_{s}}\right)>c_{i_{s}},
    $$
the opposite inequality follows: $i_{s} \ge j_{s}$. So we have shown that $i_{s}=j_{s}$. In this case, simultaneously $c_{i_{\mathrm{s}}}<f_{i_{s}}\left(x_{s}^{l_{s}}\right)$ and $f_{i_{s}}(x)=c_{i_{s}}$. Substituting these estimates into $(10.29)$, we obtain a contradiction for sufficiently large s:
    $$
    c_{i_{\mathrm{s}}} \le c_{i_{\mathrm{s}}}-\frac{\gamma^{2}}{384 \Gamma^{2}} \gamma e+(\Gamma+a)\left\|x^{s}-x\right\| .
    $$
    Hence, for sufficiently large values of the $s$
    $$
    f_{j}\left(x_{s}^{l_s} \right) \le c_{j} \quad \text { for } \quad j \in J(x) .
    $$
 The second assertion of the lemma is proved.
    
    Let us prove the first assertion of the lemma. We have
    $$
    i_{s}=i\left(x_{s}^{l_{s}}\right) \in i(x) \cup J(x) .
    $$
    It has already been shown that for sufficiently large $s$
    $$
  f_{j}\left(x_{s}^{l_{s}}\right) \le c_{j} \quad \text { for } \quad j \in J(x).
    $$
    
    Since
    $$
    f_{i_{s}}\left(x_{s}^{l_{s}}\right)>c_{i_{s}},
    $$
 then   for sufficiently large values of  $s$
    $$
    i_{s}=i\left(x_{s}^{l_{s}}\right)=i(x).
    $$
 Then the first assertion of the lemma follows from the inequality (10.29). Lemma is proven.
\begin{corollary}
\label{cor:10.2}    
		For point $x$ to be a local minimum point of (10.14), (10.15), the following conditions must be met:
    $$
    0 \in G(x), \quad f_{j}(x) \le c_{j}, \quad j=1,2, \ldots, m \text {. }
    $$
\end{corollary}
    
    Indeed, condition $f_{j}(x) \le c_{j}(j=1,2, \ldots, m)$ is an admissibility condition. Suppose that $x$ is a local minimum point, but $0 \notin G(x)$. We will repeatedly $(s=0,1, \ldots)$ run algorithm $(10.18)$ from the starting point $x_{s} \equiv x$ with steps $\rho_{s}^{k}\left(k \ge k_{s}\right)$ up to a certain moment $n_{\mathrm{s}}$ so that 
    $$
    \lim _{s \rightarrow \infty} \rho_{s}=\lim _{s \rightarrow \infty} \sup _{k_{s} \le k \le n_{s}} \rho_{s}^{k}=0, \quad \sum_{k=k_{s}}^{n_{s}-1} \rho_{s}^{k} \ge \sigma>0 .
    $$
    Then, according to Lemma 10.5, in an arbitrarily small neighborhood of  point $x$ the algorithm will find points with less than $f_{m+1}(x)$ values of the function $f_{m+1}$, which contradicts the local optimality of $x$. Thus, the minimizing property of Lemma 10.5 contains the necessary conditions for an extremum. This point has already been noted in the remark to Theorem 8.9.
    
    Since at the first stage algorithm (10.18) minimizes the function $h(x)=f_{1}(x)$, the sequences generated by the algorithm can be made bounded in the same way as it was done for method (10.7) in Lemma 10.1. Denote
    $$
    \begin{gathered}
    X_{i}^{*}=\left\{x \in E_{n} \mid 0 \in G(x) ; \quad f_{j}(x) \le c_{j}, \quad 1 \le j \le i ; \quad f_{i}(x)>c_{i}\right\}, \\
    F_{i}^{*}=\left\{f_{i}(x) \mid x \in X_{i}^{*}\right\} .
    \end{gathered}
    $$
\begin{theorem}
\label{th:10.3}    
		Assume that the sequence $\left\{x^{k}\right\}_{k=0}^{\infty}$ generated by the algorithm (10.18) - (10.22) is bounded, and the sets $F_{i}^{*}(i=1,2, \ldots, m+1)$ do not contain segments. Then the sequence $\left\{x^{k}\right\}_{k=0}^{\infty}$ converges to the solution of one of problems (10.16), (10.17), where $1 \le i^{*} \le m+1$, that is, all the limit points of $\left\{x^{k}\right\}_{k=0}^{\infty}$ are admissible for the problem (10.16), (10.17) and belong to the set $X_{i^{*}}^{*}$, and also there is a limit $\lim _{k \rightarrow \infty} f_{i^{*}}\left(x^{k}\right)>c_{i^{*}}$.
\end{theorem}
    
{\it    P r o o f.} At the first stage, the algorithm minimizes function $f_{1}(x)$. In this case, the role of the Lyapunov function is played by function $f_{1}(x)$ itself, the role of the set of solutions is played by the set
    $$
    X_{1}^{* *} \equiv X_{1}^{*} \cup\left\{x \in E_{n} \mid f_{1}(x) \le c_{1}\right\}
    $$
    
    By virtue of 
    $$
    \lim _{k \rightarrow \infty}\left\|x^{k+1}-x^{k}\right\|=0,
    $$
Lemma 10.5 (with $l(x)=1$), and the general Theorem 8.9, all limit points of $\left\{x^{k}\right\}_{k=0}^{\infty}$ belong to $X_{1}^{* *}$ and the half-interval $\left[\lim _{k \rightarrow \infty} f\left(x^{k}\right), \varlimsup_{k \rightarrow \infty} f\left(x^{k}\right)\right)$ is embedded in the set
    $$
    F_{1}^{*} \cup\left\{f \in E_{1} \mid f \le c_{1}\right\}
    $$
    
    Since the $F_{1}^{*}$ does not contain segments, then either there is a limit $\underset{k \rightarrow \infty}{\lim} f_{1}\left(x^{k}\right)>c_{1}$ and all limit points $\left\{x^{k}\right\}_{k=0}^{\infty}$ belong to $X_{1}^{*}$, or all limit points of $\left\{x^{k}\right\}_{k=0}^{\infty}$ belong to the set $\left\{x \in E_{n} \mid f_{1}(x) \le c_{1}\right\}$.
    
    If the latter case takes place, then we can assume that the algorithm minimizes function $f_{2}(x)$ under constraint $f_{1}(x) \le c_{1}$. Now the Lyapunov function is the $f_{2}(x)$, and the set of solutions is the set
    $$
    X_{2}^{* *}=X_{2}^{*} \cup\left\{x \in E_{n} \mid f_{1}(x) \le c_{1}, f_{2}(x) \le c_{2}\right\}
    $$

    By virtue of
    $$
    \lim _{k \rightarrow \infty}\left\|x^{k+1}-x^{k}\right\|=0
    $$
    Lemma 10.5 (with $l(x)=2$ ) and the general theorem 8.9, the sequence $\left\{x^{k}\right\}_{k=0}^{\infty}$ converges to the indicated set of solutions, and the semi-interval 
		$$\left[\left.\underset{k \rightarrow \infty}{\underline{\lim}} f_{2}\left(x^{k}\right), \varlimsup_{k \rightarrow \infty} f_{2}\left(x^{k}\right)\right.\right)$$ 
		is embedded in the set
    $$
    F_{2}^{*} \cup\left\{f \in E_{1} \mid f \le c_{2}\right\}
    $$
    
    Since $F_{2}^{*}$ does not contain intervals, then either there is a limit $\lim _{k \rightarrow \infty} f_{2}\left(x^{k}\right)>c_{2}$ and all limit points of $\left\{x^{k}\right\}_{k=0}^{\infty}$ belong to $X_{2}^{*}$, or $\varlimsup_{k \rightarrow  \infty} f_{2}\left(x^{k}\right) \le c_{2}$ and all limit points of $\left\{x^{k}\right\}_{k=0}^{\infty}$ belong to set $\left\{x \in E_{n} \mid f_{1}(x) \le c_{1}, f_{2}(x) \le c_{2}\right\}$ and so on. Theorem is proven.
\begin{corollary}
\label{cor:10.3}   
		If all stationary points of the mapping $G$ lie in the admissible region of the original problem (10.14), (10.15):
    $$
    \left\{x \in E_{n} \mid 0 \in G(x)\right\} \subset\left\{x \in E_{n} \mid f_{j}(x) \le c_{j}, 1 \le j \le m\right\},
    $$
    then under the conditions of Theorem 10.3 the sequence $\left\{x^{k}\right\}_{k=0}^{\infty}$ converges to the solution
    $$
    X^{*}=\left\{x \in E_{n} \mid 0 \in G(x), f_{j}(x) \le c_{j}, 1 \le j \le m\right\}
    $$
    of the original problem and there is a limit $\lim _{k \rightarrow \infty} f_{m+1}\left(x^{k}\right)$.
		\end{corollary}

\smallskip
    \textbf{3. Problems with linear constraints.} 
		\label{Sec.10.3}
		So far, we have considered problems of mathematical programming with inequality constraints. Let us now consider mathematical programming problems with general inequality constraints and linear equality constraints. Let it be necessary to solve the problem
    \begin{equation} \label{eqn:10.30}
        f(x) \rightarrow \min
    \end{equation}
subject to the constraints
    \begin{equation} \label{eqn:10.31}
    f_{j}(x) \le c_{j}, \quad j=1,2, \ldots, m,
    \end{equation}
    \begin{equation} \label{eqn:10.32}
    A x=b, \quad e \le x \le d,
    \end{equation}
    where $x=\left(x_{1}, \ldots, x_{n}\right) \in E_{n}, f(x)$ and $f_{j}(x)(j=1, \ldots, m)$ are generalized differentiable functions, $A$ is $m_{1} \times n$-matrix, $ b \in E_{m_{1}}$, $e$ and $d$ are $n$-dimensional vectors, possibly containing unbounded components.
    
    There are many possibilities for constructing algorithms for solving this problem; for example, linearization or reduced gradient methods can be generalized. However, in the context of the work done, the easiest way, apparently, is to use the method of projecting the generalized gradient onto a linear subspace given by constraints $A x=0$.
    
    Let a function $f(x)\left(x \in E_{n}\right)$ and a linear manifold
    $$
    L=\left\{x \in E_{n} \mid A x=b\right\}
    $$
    be given.

    Denote by the
    $$
    L_{0}=\left\{x \in E_{n} \mid A x=0\right\}
    $$
    the subspace associated with the $L$. We can consider the restriction of the function $f(x)$ on the $L$. If $f(x)$ is convex in the $E_{n}$, then it remains convex on a convex subset of $E_{n}$, in particular, on $L$. One can consider subgradients of a function $f(x)$ considered as a function on the $L$. We will show how to calculate these subgradients in terms of the usual subgradients $f(x)\left(x \in E_{n}\right)$. We have a subgradient inequality
    $$
    f(y) \ge f(x)+\left<g, y-x\right>, \quad x, y \in E_{n}, \quad g \in \partial f(x) .
    $$
    Let us decompose the vector $g$ in the subspace $L_{0}$ and the orthogonal complement to the $L_{0}$, that is, we represent $g=g_{0}+g^{\prime}$, where $g_{0}$ is the projection of $g$ onto $L_{0}$. When $x, y \in L$, it takes place $\left<g, y-x\right>=\left<g_{0}, y-x\right>$, and thus,
    
    $$
    f(y) \ge f(x)+\left<g_{0}, y-x\right> \text { for all } x, y \in L,
    $$
    that is $g_{0}$ is a subgradient of $f$ as a function onto $L$. The projection operator of space $E_{n}$ onto subspace $L_{0}$, as is known [120], has the form
    $$
    P=I-A\left(A^{T} A\right)^{-1} A^{T}
    $$
    where $I$ is the identity matrix of size $n \times n, A^{T}$ is the transposed matrix $A$.
    
    To minimize $f(x)$ on the $L$ it is now possible to use the known methods of convex optimization (the generalized gradient method, the ellipsoid method, methods with space dilation), in which the subgradients $f(x)$ should be used as functions on $L$. If $f(x)$ is to be minimized on $L$ under additional convex inequality constraints, then the functions in these constraints must also be considered as convex on $L$ with the corresponding subgradients. The same applies to constraints of the form $e \le x \le d$. They can also be replaced by a single constraint of the form
    $$
    h(x) \equiv \max _{1 \le i \le n} \max \left\{e_{i}-x_{i}, x_{i}-d_{i}\right\} \le 0
    $$
    or
    $$
    h(x) \equiv \sum_{i=1}^{n}\left(\max \left\{0, e_{i}-x_{i}\right\}+\max \left\{0, x_{i}-d_{i}\right\}\right) \le 0 .
    $$
    
    Projecting subgradients reduces the dimension of the working space, and thus the complexity of convex optimization methods decreases, since the estimates of complexity depend monotonically on the dimension [73].
    
    What was said above about convex functions carries over entirely to generalized differentiable functions. Let $f(x)$ be a generalized differentiable function, $G_{f}(x)$ be its pseudogradient mapping. Consider a multivalued mapping
    $$
    P G_{f}(x)=\left\{g_{0} \in L_{0} \mid g_{0}=P g, g \in G_{f}(x)\right\}
    $$
    where $P$ is the  projection operator of $E_{n}$ on $L_{0}$. Sets $P G_{f}(x)$, similar to $G_{f}(x)$, are non-empty, convex, closed, and bounded, and mapping $x \rightarrow P G_{f}(x)$ is upper semi-continuous. Let $x, y \in L$, $g \in G(y)$ and $g=P g+g^{\prime}$; then
    $$
    \begin{aligned}
    & o(x, y, P g)=f(y)-f(x)-\left<P g,y-x\right>= \\
    &=f(y)-f(x)-\left<g, y-x\right>=o(x, y, g)
    \end{aligned}
    $$
    satisfies the required smallness conditions (1.3). Thus, a function $f(x)$ that is generalized differentiable in $E_{n}$ is also generalized differentiable as a function in $L$ with pseudo-gradient mapping $P G_{f}(x)$.
    
    Now we can consider the restriction of the original problem (10.30)-(10.32) to the set $L$ and to solve it by applying the methods of the generalized gradient $(10.7)$ and the leading constraint method (10.18)-(10.22) considered earlier, in which instead of the pseudo-gradients of functions in the original space $E_{n}$ are substituted by their projections onto the subspace $L_{0}$. The starting point for methods must be taken from $L$. The constraints $e$ $\le x \le d$ should be taken into account on a par with other inequality constraints; in particular, they can be folded into one or ranked by importance.

\section*{$\S$ 11. Construction of relaxation methods for non-convex non-smooth optimization problems}
\label{Sec.11}
\setcounter{section}{11}
\setcounter{definition}{0}
\setcounter{equation}{0}
\setcounter{theorem}{0}
\setcounter{lemma}{0}
\setcounter{remark}{0}
\setcounter{corollary}{0}

In this section, we present a scheme for constructing relaxation algorithms for local minimization of generalized differentiable functions.

In relaxation methods for nonsmooth optimization, two interrelated problems must be solved: constructing a descent direction for the function and selecting a step size in that direction. We construct a relaxation method in which the key problem of finding a descent direction is reduced to a purely geometric task, which allows us to speak about its complexity and the search for an optimal solution algorithm. In addition, the method provides ample opportunities for optimizing the choice of step size. Finally, the proposed method stops after a finite number of iterations when the specified solution accuracy is reached. If the permissible deviation is made to tend to zero, then the obtained approximate solution will tend to the exact one. Some estimates of the computational complexity of the method on the class of convex problems are also given.

Let it be necessary to solve the problem
\begin{equation} \label{eqn:11.1}
 f(x) \rightarrow \min_x , \quad x \in E_{n};   
\end{equation}
$f(x)$ is a generalized differentiable function; $G_{f}(x)$ is its pseudo-gradient mapping.

According to Theorem 8.9, all local minimum points of $f(x)$ belong to the set (solutions)
$$
X^{*}=\left\{x \in E_{n} \mid 0 \in G_{f}(x)\right\}.
$$

The set $X^{*}$ will be considered as the exact solution of  problem (11.1). Along with it, we consider some approximate solutions $X_{\varepsilon\delta}^{*}$.

Let us introduce the notation
$$
\Delta(C, D)=\sup _{c \in C} \inf_{d\in D} \|c-d\|,
$$
the deviation of set $C$ from set $D$;
$$
\begin{aligned}
& \gamma_{\delta}(x)=\rho\left(0, \text { co } \underset{\{y|\|y-x \mid\| \le \delta\}}{\cup} G_{f}(y)\right) \text {; } \\
& X_{\varepsilon\delta}^{*}=\left\{x \in E_{n} \mid \gamma_{\delta}(x) \le \varepsilon\right\} .
\end{aligned}
$$
\begin{lemma}
\label{lem:11.1}
Function $\delta \rightarrow \gamma_{\delta}(x)$ is monotone decreasing along $\delta$, and function $(\delta, x) \rightarrow \gamma_{\delta}(x)$ is lower semicontinuous, that is, while $\left(\delta_{k}, x^{k}\right) \rightarrow (\delta, x)$ and $\gamma_{\delta_{k}}\left(x^{k}\right) \rightarrow \gamma$ the following inequality holds
$$
\gamma \ge \gamma_{\delta}(x).
$$
\end{lemma}

{\it P r o o f.} If $\delta$ increases, then the set with $\left\{G_{f}(y) \mid \| y-x \| \le \delta\right\}$ expands; hence $\gamma_{0}(x)$ is decreasing. Let us show that the function $(\delta, x) \rightarrow \gamma_{\delta}(x)$ is lower semicontinuous. Recall that
$$
\gamma _{\delta _ {\kappa }}(x) = \inf \left \{ \left \| \text{g} \right \| | \text{g} \in \text { co } \underset{\{y|\|y-x \mid\| \le \delta_k\}}{\cup}  G_{f}(y) \right \}.
$$
Therefore, for the numerical sequence $\varepsilon_{k} \rightarrow 0$ there are vectors $g_{i}^{k}$ and numbers $\lambda_{i}^{k}\;(i=1,2, \ldots, n+1 ;$ $ n$ is the space dimension) that
$$
\begin{gathered}
\left\|g^{k}\right\| \le \gamma_{\delta_{k}}\left(x^{k}\right)+\varepsilon_{k}, \quad g^{k}=\sum_{i=1}^{n+1} \lambda_{i}^{k} g_{i}^{k}, \\
g_{i}^{k} \in G_{f}\left(y_{i}^{k}\right), \quad\left\|y_{i}^{k}-x^{k}\right\| \le \delta_{k}, \quad \lambda_{i}^{k} \ge 0, \quad \sum_{i=1}^{n+1} \lambda_{i}^{k}=1 .
\end{gathered}
$$

Due to the local boundedness of $G_{f}$ there exists a subsequence of indices $\left\{k_{s}\right\}_{s=0}^{\infty}$ such that
$$
\lim _{s \rightarrow \infty} g_{i}^{k_{s}}=g_{i}, \quad \lim _{s \rightarrow \infty} y_{i}^{k_{s}}=y_{i} \quad \lim _{s \rightarrow \infty} \lambda_{i}^{k_{s}}=\lambda_{i} .
$$
It's obvious that
$$
\left\|y_{i}-y\right\| \le \delta, \quad \lambda_{i} \geq 0, \quad \sum_{i=1}^{n+1} \lambda_{i}=1,
$$
and due to the closedness of $G_{f}$ $g_{i} \in G_{f}\left(y_{i}\right)$. Let us introduce the notation
$$
\delta_{s}^{\prime}=\delta_{k_{s}}, \quad g=\sum_{i=1}^{n+1} \lambda_{i} g_{i} \in co \bigcup_{\{y: \|y-x\|\le \delta\}} G_{f}(y) .
$$
Passing in the inequality
$$
\left\|g^{k_{s}}\right\| \le \gamma_{\delta_{\mathrm{s}}^{\prime}}\left(x^{k_{s}}\right)+
\varepsilon_{k_{s}}
$$
to the limit in $s$, we obtain
$$
\|g\| \le \lim _{s \rightarrow \infty} \gamma_{\delta_{s}^{\prime}}\left(x^{k_{s}}\right)=\gamma.
$$
However $\gamma_{0}(x) \le\|g\|$. Thus,
$$
\gamma \ge\|g\| \ge \gamma_{\delta}(x),
$$
which was to be proved.
\begin{lemma}
\label{lem:11.2}
The multivalued mapping $(\varepsilon, \delta) \rightarrow X_{e \delta}^{*}$ is monotonically increasing, that is, $X_{\varepsilon_{1} \delta_{1}}^{*} \subset X_{\mathrm{\varepsilon}_{2} \delta_{2}}^{*}$ for $\varepsilon_{1} \le \varepsilon_{2}, \delta_{1} \le \delta_{2}$, and closed, that is, with $ \left(\varepsilon_{k}, \delta_{k}\right) \rightarrow(\varepsilon, \delta),$ $x^{k} \in X_{\varepsilon_{k} \delta_k}^{*}$ and $\left\{x^{k}\right\} \rightarrow x$ we have $x \in X_{\varepsilon \delta}^{*}$.
\end{lemma}

{\it P r o o f.} Let $\varepsilon_{1} \le \varepsilon_{2}, \delta_{1} \le \delta_{2}, \gamma_{\delta_{1}}(x) \le \varepsilon_{1}$. Due to the monotonic decrease of the $\gamma_{\delta}(x)$ on $\delta$, we have
$$
\gamma_{\delta_{2}}(x) \le \gamma_{\delta_{1}}(x) \le \varepsilon_{1} \le \varepsilon_{2};
$$
then $x \in X_{\varepsilon_2, \delta_2}^{*}$, which was to be proved. 
Let 
$$
\left(\varepsilon_{k}, \delta_{k}\right) \rightarrow(\varepsilon,\delta), \quad 
x^{k} \in X_{\varepsilon_{k}, \delta_{k}}^{*}, \quad\left\{x^{k}\right\} \rightarrow x.
$$
Then
$$
\gamma_{\delta_{k}}\left(x^{k}\right) \le \varepsilon_{k}.
$$
Passing to the limit with $k$, due to the lower semicontinuity of $\gamma_{\delta}(x)$ on $(\delta, x)$ we have
$$
\gamma_{\delta}(x) \le \lim _{k \rightarrow \infty} \gamma_{\delta_{k}}\left(x^{k}\right) \le \varepsilon,
$$
that is, $x \in X_{\varepsilon \delta}^{*}$, which was to be proved.
\begin{corollary}
\label{cor:11.1}
The set $X^{*}$, together with its $\delta$-neighbourhood, is embedded in $X_{\varepsilon\delta}^{*}$. Let $D$ be arbitrary closed bounded set in $E_{n}$. The multi-valued mapping $(\varepsilon, \delta) \rightarrow X_{\varepsilon \delta}^{*} \cap D$ is bounded and closed and hence upper semicontinuous. Therefore, with $X^{*} \cap D \neq \emptyset$ we have $X_{\varepsilon \delta}^{*} \cap D \neq \emptyset$ and
$$
\Delta\left(X_{\varepsilon \delta}^{*} \cap D, X^{*} \cap D\right) \rightarrow 0 \quad \mbox{ for } \quad(\epsilon, \delta) \rightarrow(0,0),
$$
that is, the set $X_{\varepsilon\delta}^{*} \cap D$ is an approximate solution to problem $(11.1)$.
\end{corollary}

The algorithm $\mathrm{Al}$ below  builds a sequence of approximations $\left\{x^{k}\right\}$ to the set $X_{\varepsilon\delta}^{*} \cap\left\{x \in E_{n} \mid f(x) \le f\left(x^{0}\right)\right\}$ and it converges to this set in a finite number of iterations. This algorithm allows a certain freedom of actions at some steps, that is, in essence, it contains a whole class of algorithms.

A l g o r i t h m A1.

S t e p  0. Choose the initial approximation $x^{0} \in E_{n}$, accuracy  $\varepsilon>0, \delta>0$ and parameters $\Lambda \ge \varepsilon,$ $\Delta \ge \delta,$ $\alpha,$ $\beta,$ $\theta$ $(0<\alpha, \beta, \theta<1)$, put $k=0$ (the beginning of the loop over $k$).

S t e p 1. Choose $\varepsilon_{0} \in[\varepsilon, \Lambda]$ and $\delta_{0} \in[\delta, \Delta]$, put $i=0$ (the beginning of the loop over $i$).

S t e p 2. Choose an initial trial direction $l^{0} \in E_{n},\left\|l^{0}\right\|>0$. Calculate trial point $y^{1}=x^{k}-\delta_{i} l^{0} /\left\|l^{0}\right\|$ and arbitrary generalized gradient $g^{1} \in G_{f}\left(y^{1}\right)$. Put $j=1$ (the beginning of the loop over $j$).

S t e p 3. Calculate new trial direction $l^{i}$ (the rule will be described below).

S t e p 4. If $\left\|l^{j}\right\| \le \varepsilon $ and $\delta_{i} \le \delta$, then stop.

S t e p 5. If $\left\|l^{j}\right\| \le \varepsilon_{i}$,  then  $\varepsilon_{i+1}=\alpha \varepsilon_{i}, \delta_{i+1}=\beta \delta_{i}$, replace $i$ with $i+1$, go to step 2.

S t e p 6. Calculate a new trial point $y^{j+1}=x^{k}-\delta_{i} l^{j} /\left\|l^{j}\right\|$ and arbitrary generalized gradient $g^{j+1} \in G_{f}\left(y^{j+1}\right)$.

S t e p 7. If $\left(g^{j+1}, l^{j}\right)<\theta \varepsilon_{i}^{2} / 2$, then replace $j$ with $j+1$, go to step 3.

S t e p 8. If $f\left(y^{j+1}\right)>f\left(x^{k}\right)-\theta \delta_{i} \varepsilon_{i}^{2} /\left(4\left\|l^{j}\right\|\right)$, then $\varepsilon_{l+1}=\varepsilon_{i}, \quad \delta_{i+1}=\beta \delta_{i}$, replace $i$ with $i+1$, go to step 2.

S t e p 9. Find the next approximation $x^{k+1}$ such that $f\left(x^{k+1}\right) \le f\left(y^{j+1}\right)$; in particular, we can take $x^{k+1}=y^{j+1}$. Replace $k$ with $k+1$, go to step 1.
\begin{remark}
\label{rem:11.1}
In the loop over $k$ the algorithm generates a sequence of approximations $x^{k}$, in the loop over $i$ the algorithm selects for a fixed $x^{k}$ the value of the trial step $\delta_{i}$, and in the loop over $j$ finds the direction of descent $l^{j}$ for fixed $x^{k}$ and $\delta_{i}$. At step 9, an arbitrary algorithm can be used to find a point $x^{k+1}$ that is better than $y^{i+1}$ (by the value of $f$); for instance, one-dimensional optimization of $f$ in the direction of $l^{i}$ can be performed.
\end{remark}

In algorithm A1, at step 3, rule L is used to construct test directions $l^{j}$ from the accumulated generalized gradients $g^{1}, g^{2}, \ldots$ $\ldots, g^{j}$. The rule L must satisfy the following two conditions:

Y1) $l^{j} \in \mathrm{co}\left\{g^{1}, \ldots, g^{j}\right\}, \quad j=1,2, \ldots$;

Y2) if the sequence $\left\{g^{i}\right\}_{j=1}^{\infty}$  is embedded in a convex set $P$ such that $0<\gamma \le\|g\| \le \Gamma<\infty$ for $g \in P$, then there is an index $m$ such that
$$
\left(g^{m+1}, l^{m}\right) \ge \theta \gamma^{2} / 2, \quad 0<\theta<1 .
$$

These conditions are satisfied, as will be shown in Lemmas $11.3-11.5$, for instance, by the following rules for constructing trial directions $l^{j}$.

L1. Let the numbers $\left\{\rho_{r}\right\}_{r=1}^{\infty}$ be such that $0 \le \rho_{r} \le \rho<\infty, \quad \sum_{r=1}^{\infty} \rho_{r}=\infty$.

Let
\begin{equation*}
    l^j = \sum_{r=1}^j \rho_r g^r / \sum_{r=1}^j \rho_r = \left(1 - \rho_j / \sum_{r=1}^j \rho_r\right) l^{j - 1} + \left(\rho_j / \sum_{r=1}^j \rho_r\right) g^j.
\end{equation*}
Particularly, with $ \rho_r \equiv \rho $ we have
\begin{equation*}
    l^j = \frac{1}{j} \sum_{r=1}^j \rho_r g^r = \left(1 - \frac{1}{j}\right) l^{j-1} + \frac{g^j}{j}.
\end{equation*}

L2. $ l^j \in \text{co} \{ l^{j-1}, g^j \} $ and $ \| l^j \| $ is minimal, $ l^1 = g^1 $.

L3. $ l^j \in \text{co} \{ g^1, ..., g^j \} $ and $ \| l^j \| $ is minimal.

L4. $ l^j \in \text{co} \{ l^{j-1}, g^i, i \in I_j \} $, where $ i \in I_j \subset \{ 1, ..., j \} $ and $ \| l^j \| $ is minimal.
\begin{theorem}
\label{th:11.1}
Let $ f(x) $ be a generalized differentiable function, the set $ D = \{ x \in E_n \> | \> f(x) \leq f(x^0) \} $ is bounded. Then  algorithm  A1 stops after a finite number of iterations at a point $ x^k $ belonging to the set $ X_{\varepsilon \delta}^{*} \cap D $.
\end{theorem}

{\it P r o o f.} By construction (see steps 8, 9) the algorithm under consideration is relaxational, so all approximations $ x^k $ belong to the set $ D $. If the algorithm stops after a finite number of iterations at some point $ x^k $ at step 5, then by the stopping condition the corresponding $ \| l^j \| \leq \varepsilon $ and $ \delta_i \leq \delta $, so $ x^k \in X_{\varepsilon \delta}^{*} \cap D $.

Let us prove the finiteness of the algorithm from the contrary. So let us assume that the algorithm works infinitely long. Let us show that then it produces an infinite sequence of points $ \{ x^k \}_{k=0}^{\infty} $, i.e. looping at steps 5, 7, 8 is impossible. At step 5 it is impossible because in this case sooner or later $ \varepsilon_i \leq \varepsilon $ and $ \delta_i \leq \delta $ will happen, and the algorithm stops. In step 7 it is impossible due to condition Y2 on rule $ L $. At step 8 it is impossible due to the properties of the functions in question. For example, for convex functions $ f(x) $, step 8 is always skipped, since at step 8 it is always takes place
\begin{eqnarray} 
    f(y^{j+1}) &\leq& f(x^k) + \left<g^{j+1}, y^{j+1} - x^k\right> \nonumber \\
		&=& f(x^k) - \delta_i \left<g^{j + 1}, \frac{l^j}{\| l^j \|}\right> \nonumber\\
		&\leq& f(x^k) - \theta \varepsilon_i^2 \delta_i / (2 \| l^j \|).\label{eqn:11.2}
\end{eqnarray}

The numerical sequence $ \{ f(x^k) \}_{k=0}^{\infty} $ is monotonically decreasing and bounded from below, so there is a finite limit $ \lim_{k \to \infty} f(x^k) $. But on the other hand, in the convex case it is always $ \varepsilon_i \geq \alpha^p \Lambda $ and $ \delta_i \geq \beta^p \Delta $, where the exponent of degree $ p $ is such that $ \alpha^p \Lambda \leq \varepsilon $ and $ \beta^p \Delta \leq \delta $. Moreover, $ || l^j || \leq \Gamma < + \infty $. Thus, $ f(x^k) $ decreases each time, according to (11.2), by at least some fixed value, and hence $ \{ f(x^k) \}_{k=0}^{\infty} $ cannot have a limit. In the convex case the theorem is proved.

Let us now show that if $ f(x) $ is generalized differentiable, then it is also impossible to loop at step 8. Indeed, if we loop at this step, $ \lim_{i \to \infty} \delta_i = 0 $ and $ \varepsilon_i \geq \gamma = \alpha^p \Lambda $, where $ p $ is such that $ \alpha^p \Lambda \leq \varepsilon $ and $ \beta^p \Delta \leq \delta $. Let us denote by
\begin{equation*}
    \Gamma = \sup \{ \| g \| \> | \> g \in G_j (y), f(y) \leq f(x^0) \}.
\end{equation*}
We have $ \Gamma < + \infty $ because of the boundedness on the compactness of the compact-valued upper semi-continuous mapping $ G_f $. For sufficiently small $ \delta_i $, due to the generalized differentiability of $ f(x) $, at step 8 we have
    \begin{eqnarray}{lcl}
        f(y^{j+1}) &\leq& f(x^k) + \left<g^{j+1}, y^{j+1} - x^k\right> + o(x^k, y^{j+1}, g^{j+1}) \nonumber \\
                   &\leq& f(x^k) - \delta_i \left<g^{j+1}, \frac{l^j}{\| l^j \|}\right> + \frac{\theta \gamma^2}{4 \Gamma} \| y^{j+1} - x^k \| \nonumber \\
                   &\leq& f(x^k) - \frac{\theta \varepsilon_i^2 \delta_i}{2 \| l^j \|} + \frac{\theta \gamma^2 \delta_i}{4 \Gamma}\nonumber\\
									&\leq& f(x^k) - \frac{\theta \varepsilon_i^2 \delta_i}{4 \| l^j \|},\nonumber
    \end{eqnarray}
i.e. it is not possible to loop at step 8.

So, if the algorithm does not stop, it produces an infinite sequence of approximations $ \{ x^k \} $. Since $ \{ f(x^k) \}_{k=0}^{\infty} $ is monotonically decreasing and bounded from below, there exists a limit $ \lim_{k \to \infty} f(x^k) $. Let $ \lim_{s \to \infty} x^{k_s} = x' $. By the condition of the transition from step 8 to step 9 we have
\begin{equation} \label{eqn:11.3}
    f(x^{k_s + 1}) \leq f(x^{k_s}) - \frac{\theta \varepsilon_i^2 \delta_i}{4 \| l^j \|}.
\end{equation}
Note that $ \| l^j \| \leq \Gamma < + \infty $ and $ \varepsilon_i \geq \gamma = \alpha^p \Lambda > 0 $, where $ p $ is such that $ \alpha^p \Lambda \leq \varepsilon $ and $ \beta^p \Delta \leq \delta $. Therefore, if we show that
\begin{equation} \label{eqn:11.4}
    \lim_{i \to \infty} \delta_i = \lim_{s \to \infty} \| y^{j + 1} - x^{k_s} \| > 0,
\end{equation}
then, passing in (11.3) to the limit on $s$, we obtain a contradiction with the convergence $ \{ f(x^k) \}_{k=0}^{\infty} $ and the theorem will be proved.

So let us show that (11.4) holds. Let us write down the expansion for $ f(x) $ at the point $ x' $:
\begin{equation*}
    f(y) = f(x') + \left<g_y, y - x'\right> + o(x', y, g_y).
\end{equation*}
According to the definition of a generalized differentiable function, for $ a = \frac{\theta \gamma^2}{12 \Gamma} $ there exists such  $ \delta' (a) $ that for all $ \{ y \> | \> \| y - x' \| \leq \delta' \} $ and all $ g_y \in G_f (y) $ it holds
\begin{equation*}
    o(x', y, g_y) \leq a \| y - x' \|.
\end{equation*}
For $ x^{k_s} $ and $ \{ y \> | \> \| y - x' \| \leq \delta' \} $ we can write
    \begin{eqnarray*}
        f(y) &\leq& f(x') + \left<g_y, y - x'\right> + a \| y - x' \| \\
             &\leq& f(x^{k_s}) + \left<g_y, y - x^{k_s}\right> + a \| y - x^{k_s} \| + f(x')  \\
             &&- f(x^{k_s}) + \left<g_y, x^{k_s} - x'\right> + a \| x^{k_s} - x' \|  \\
             &\leq& f(x^{k_s}) + \left<g_y, y - x^{k_s}\right> + a \| y - x^{k_s} \| + C \| x^{k_s} - x' \|,
    \end{eqnarray*}
where $ C = 2 \Gamma + a $. We define $ \delta'' = \frac{\delta'}{2} $ and show that for all sufficiently large $ k_s $
\begin{equation*}
    \| y^{j+1} - x^{k_s} \| \geq \beta \delta^{\prime\prime}.
\end{equation*}
We can find the value  $ K $ such that for all $ k_s \geq K $
\begin{equation*}
    \| x^{k_s} - x' \| \leq \min \{ \delta^{\prime\prime}, \alpha \beta \delta^{\prime\prime} / (2C) \}.
\end{equation*}
For $ k = k_s $ the algorithm at step 8 will produce variables $ y^{j+1} $, $ l^j $, $ \varepsilon_i $, $ \delta_i $, such that
\begin{equation*}
    \left<g^{j+1}, l^j\right> \geq \frac{\theta \varepsilon_i^2}{2}, \;\;\; \| y^{j+1} - x^{k_s} \| = \delta_i.
\end{equation*}

Let us show that if $ k_s \geq K $, then for all corresponding $ \delta_i $ the inequality $ \delta_i \geq \beta \delta'' $ is satisfied. Indeed, if this is not the case, then for some $ k_s \geq K $ and $ i $
\[
    \begin{matrix}
        \beta \delta^{\prime\prime} \leq \delta_i \leq \delta^{\prime\prime}, & \| y^{j+1} - x^{k_s} \| = \delta_i, & \left<g^{j+1}, l^j\right> \geq \frac{\theta \varepsilon_i^2}{2}.
    \end{matrix}
\]
Let us note that $ \| y^{j+1} - x' \| = \delta' $ and $ \varepsilon_i \geq \gamma > 0 $. Let us substitute these $ y^{j+1} $ and $ g^{j+1} $ into the above inequality:
    \begin{eqnarray}
        f(y^{j+1}) &\leq& f(x^{k_s}) + \left<y^{j+1} - x^{k_s}\right> + a \| y^{j+1} - x^{k_s} \| + \| x^{k_s} - x' \| \nonumber \\
        &\leq& f(x^{k_s}) - \delta_i \left<g^{j+1}, \frac{l^j}{\| l^j \|}\right> + a \delta_i + \frac{a \beta \delta^{\prime\prime}}{2} \nonumber  \\
        &\leq& f(x^{k_s}) - \frac{\theta \varepsilon_i^2 \delta_i}{2 \| l^j \|} + \frac{\theta \gamma^2 (2 \delta_i + \beta \delta^{\prime\prime})}{24 \Gamma} \nonumber\\
				&\leq& f(x^{k_s}) - \frac{\theta \delta_i \varepsilon_i^2}{4 \| l^j \|}.\nonumber
    \end{eqnarray}
But this means that at step 8 $ \delta_i $ cannot decrease further. Thus, for sufficiently large $ k_s $
\begin{equation*}
    \delta_i = \| y^{j+1} - x^{k_s} \| \geq \beta \delta^{\prime\prime},
\end{equation*}
which was required to show. The theorem is proved.
\begin{remark}
\label{rem:11.2}
At step 1 of the algorithm for $ k \geq 1 $, we can assume
\[
    \begin{matrix}
        \delta_0 = \max \{ \| x^k - x^{k-1} \|, \delta \} \\
        \varepsilon_0 = \max \{ \frac{f(x^{k-1}) - f(x^k)}{\| x^k - x^{k-1} \|}, \varepsilon \};
    \end{matrix}
\]
the theorem remains valid.
\end{remark}
\begin{remark}
\label{rem:11.3}
If at step 4, instead of "stop", you enter the division of $ \varepsilon $ and $ \delta $ by a number greater than one, then by virtue of Corollary 11.1 and Theorem 11.1 the iterative algorithm will converge to the set $ X^* $ by limit points.
\end{remark}
\begin{remark}
\label{rem:11.4}
The algorithm A1  can be extended to conditional optimization problems in the same way as it was done for the generalized gradient method in $\S$ 9.
\end{remark}

As can be seen from the proof of Theorem 11.1, the rule $ L $ for constructing test directions $ l^j $ solves the following problem.

{\it The direction finding problem.} Let there be a set $ P \subset E_n $ such that $ 0 < \gamma \leq || g || \leq \Gamma < + \infty $ for all $ g \in \text{co} P $. We obtain information about the set $ P $ as a sequence of vectors $ g^j \in P \;(j = 1, 2, ...) $. From  the  information $ {g^1, g^2,..., g^j} $ available at each moment $ j $, we have to construct vectors $ 
l^j \in \text{co} {g^1, g^2,..., g^j}\; (j = 1, 2, ...) $ so that to arrive to the inequality
\begin{equation*}
    \left<g^{j+1}, l^j\right> \geq \frac{\theta \gamma^2}{2}.
\end{equation*}

It is not immediately obvious that there are rules for constructing directions $ l^j $ that solve the direction finding problem. But we will show that rules L1-L4 solve the above problem, and give some estimates of the number of steps needed to solve it. Here the question of the optimal rule naturally arises, which we, however, leave open.
\begin{definition}
\label{def:11.1}
The labor complexity of the rule $ L $ on the direction finding problem, we will call the number
\begin{equation*}
    m_1 (P, L, \theta) = \sup_{\{ g^j \}_{j=1}^{\infty}\subset P} \min \{ j \> | \> \left<g^{j+1}, l^j\right> \geq \frac{\theta \gamma^2}{2}, \> 0 < \theta < 1 \}.
\end{equation*}
\end{definition}
\begin{lemma}
\label{lem:11.3}
Rule L1 solves the  direction finding problem. If in it $ \rho_r = \rho $, i.e.
\begin{equation*}
    l^j = \frac{1}{j} \sum_{r=1}^j g^r
\end{equation*}
then the following estimate for labor complexity holds:
\begin{equation*}
    m_1 \leq \frac{2}{1 - \theta} \frac{\Gamma^2}{\gamma^2}.
\end{equation*}
\end{lemma}

{\it P r o o f.} Let $ \{ g^j \}_{j=1}^{\infty} $ be an arbitrary sequence of vectors from the convex set $ P \subset E_n $ such that
\begin{equation*}
    0 < \gamma \leq \| g \| \leq \Gamma < + \infty, \> \> \> \text{for} \> \> \> g \in P.
\end{equation*}
Let for all $ j < m $ there be
\begin{equation*}
    \left<g^{j+1}, l^j\right> \geq \frac{\theta \gamma^2}{2}, \> \> \> 0 < \theta < 1.
\end{equation*}
Let's sum the equations
\begin{equation*}
    \left\| \sum_{r=1}^j \rho_r g^r \right\|^2 = \left\| \sum_{r=1}^{j-1} \rho_r g^r \right\|^2 + 2 \rho_j \left<g^j, \sum_{r=1}^{j-1} \rho_r g^r\right> + \rho_j^2 \| g^j \|^2
\end{equation*}
along $ j $  from $ 2 $ to $ m $:
\begin{equation*}
    \left\| \sum_{r=1}^m \rho_r g^r \right\|^2 = \rho_1^2 \| g^1 \|^2 + 2 \sum_{j=2}^m \rho_j \left<g^j, \sum_{r=1}^{j-1} \rho_r g^r\right> + \sum_{j=2}^m \rho_j^2 \| g^j \|^2.
\end{equation*}
The following estimates take place,
\begin{eqnarray}
      \gamma^2 \left(\sum_{r=1}^m \rho_r\right)^2 
				&\leq&
				\left\| \sum_{r=1}^m \rho_r g^r \right\|^2 \leq \rho_1 \Gamma^2 + \theta \gamma^2 \sum_{j=2}^m \rho_j \sum_{r=1}^{j-1} \rho_r + \rho \Gamma^2 \sum_{j=2}^m \rho_j \nonumber \\
        &\leq& \rho_1 \Gamma^2 + \theta \gamma^2 \left(\sum_{j=1}^m \rho_j\right)^2 + \rho \Gamma^2 \sum_{j=1}^m \rho_j,\nonumber
\end{eqnarray}
from where
\begin{equation*}
    (1  -\theta) \frac{\gamma^2}{\Gamma^2} \leq \frac{\rho_1}{(\sum_{j=1}^m \rho_j)^2} + \frac{\rho}{\sum_{j=1}^m \rho_j}.
\end{equation*}

Hence, it is clear that the number $ m $ cannot be arbitrary large, i.e., there exists $ m $ such that
\begin{equation*}
    \left<g^{m+1}, l^m\right> \geq \frac{\theta \gamma^2}{2}.
\end{equation*}
If $ \rho_r \equiv \rho $, then
\begin{equation*}
    (1 - \theta) \frac{\gamma^2}{\Gamma^2} \leq \frac{1}{m^2} + \frac{1}{m} \leq \frac{2}{m}
\end{equation*}
therefore $ m \leq \frac{2}{1 - \theta} \frac{\Gamma^2}{\gamma^2} $, so the lemma is proven.
\begin{lemma}\label{lem:11.4}
Rule L2 solves the direction finding problem and hence satisfies conditions Y1, Y2. Moreover, for the labor intensity of the form
\begin{equation*}
    m_2 (P, L2) = \sup_{\{ g^j \}_{j=1}^{\infty} \subset P} \min \{ j \> | \> \left<g^{j+1}, l^j\right> \geq \frac{1}{2} \> \| l^j \|^2 \}
\end{equation*}
the following estimate holds true:
\begin{equation*}
    m_2 \leq 1 + 8 \frac{\Gamma^2}{\gamma^2} \ln{\frac{\Gamma}{\gamma}}.
\end{equation*}
Since it follows from $ \left<g^{j+1}, l^j\right> \geq \frac{\| l^j \|^2}{2} $ that $ \left<g^{j+1}, l^j\right> \geq \frac{\theta \gamma^2}{2} $, $ 0 < \theta < 1 $, the same estimate is true for the labor complexity $ m_1 (P, L2, \theta) $ from Definition 11.1.
\end{lemma}
{\it P r o o f.} Let $ P $ be a convex set,
\begin{equation*}
    0 < \gamma \leq || g || \leq \Gamma < + \infty, \> g \in P; \> \> \> \{ g^j \}_{j=1}^{\infty} \subset P.
\end{equation*}
The vectors $ l^j $ are constructed recurrently:
\begin{equation*}
    l^1 = g^1; \> \> \> l^j \in \text{co} \{ l^{j-1}, g^j \},
\end{equation*}
where $ \| l^j \| $ is minimal. Let for all $ j \leq m, $
\begin{equation*}
    \left<g^j, l^{j-1}\right> \leq \frac{\| l^{j-1} \|^2}{2}
\end{equation*}
or
\begin{equation*}
    \left<g^j, \frac{l^{j-1}}{\| l^{j-1} \|}\right> \leq \frac{\| l^{j-1} \|}{2}.
\end{equation*}

Let us try to obtain estimates of the form
\begin{equation*}
    \| l^j \| \leq q \| l^{j-1} \|, \> \> \> 0 < q < 1,
\end{equation*}
provided that the projections of vectors $ g^j $ on the directions $ l^{j-1} $ do not exceed $ \frac{\| l^{j-1} \|}{2} $. It is necessary to consider several cases. If
\begin{equation*}
    \| g^j || \leq \frac{\| l^{j-1} \|}{2},
\end{equation*}
then
\begin{equation*}
    \| l^j \| \leq \| g^j \| \leq \frac{\| l^{j-1} \|}{2}.
\end{equation*}
If
\begin{equation*}
    \frac{\| l^{j-1} \|}{2} \leq \| g^j \| \leq \frac{\| l^{j-1} \|}{\sqrt{2}},
\end{equation*}
then from geometrical considerations we find
\begin{equation*}
    \| l^j \| \leq \| g^j \| \leq \frac{\| l^{j-1} \|}{\sqrt{2}}.
\end{equation*}
Now let us assume
\begin{equation*}
    \| g^j \| \geq \frac{\| l^{j-1} \|}{\sqrt{2}}.
\end{equation*}
We define $ x = \left<g^j, \frac{l^{j-1}}{\| l^{j-1} \|}\right> $. From geometrical considerations we can find
\begin{equation*}
    \| l^j \|^2 = \| l^{j-1} \|^2 \frac{\| g^j \|^2 - x^2}{\| g^j \|^2 - 2x \| l^{j-1} \| + \| l^{j-1} \|^2}.
\end{equation*}
The maximum of the right-hand side for $ x \leq \frac{\| l^{j-1} \|}{2} $ is reached at $ x = \frac{\| l^{j-1} \|}{2} $; so
\begin{equation*}
    \| l^j \|^2 = \| l^{j-1} \|^2 \left(1 - \frac{\| l^{j-1} \|^2}{4 \| g^j \|^2}\right) \leq \| l^{j-1} \|^2 \left(1 - \frac{\gamma^2}{4 \Gamma^2}\right).
\end{equation*}
Since
\begin{equation*}
    1 - \frac{\gamma^2}{2 \Gamma^2} > \frac{1}{\sqrt{2}},
\end{equation*}
then with
\begin{equation*}
    \left<g^j, \frac{l^{j-1}}{\| l^{j-1} \|}\right> \leq \frac{\| l^{j-1} \|}{2},
\end{equation*}
we always get
\begin{equation*}
    \| l^j \|^2 \leq \| l^{j-1} \|^2 \left(1 - \frac{\gamma^2}{4 \Gamma^2}\right).
\end{equation*}

For $ j=m $, we get
\begin{equation*}
    \gamma^2  \leq \| l^m \|^2 \leq \| l^1 \|^2 \left(1 - \frac{\gamma^2}{4 \Gamma^2}\right)^{m-1} \leq \Gamma^2 \left(1 - \frac{\gamma^2}{4 \Gamma^2}\right)^{m-1},
\end{equation*}
where
\begin{equation*}
    m \leq 1 + 8 \frac{\Gamma^2}{\gamma^2} \ln \frac{\Gamma}{\gamma}.
\end{equation*}
The lemma is proven.
\begin{lemma}
\label{lem:11.5}
Rules L3, L4 solve the direction finding problem and hence satisfy conditions Y1, Y2. For the labor complexities $ m_2 (P, L3) $ and $ m_2 (P, L4) $ the estimation from Lemma 11.4 for sure to be valid:
\begin{equation*}
    \max \{ m_2 (P, L3), m_2 (P, L4) \} \leq 1 + 8 \frac{\Gamma^2}{\gamma^2} \ln \frac{\Gamma}{\gamma}.
\end{equation*}
\end{lemma}

{\it P r o o f.} Let the direction finding problem be solved, and the vectors $ l^j $ are constructed according to the L3 or L4 rules. Let for all $ j \leq m $ the inequality be satisfied:
\begin{equation*}
    \left<g^j, l^{j-1}\right> \leq \frac{\| l^{j-1} \|^2}{2}.
\end{equation*}
Consider a vector $ \bar{l}^j $ constructed from $ l^{j-1} $ and $ g^j $ according to the L2 rule. Obviously, $ \| l^j \| \leq \| \bar{l}^j \| $. In the proof of Lemma 11.4, the inequality was obtained
\begin{equation*}
    \| \bar{l}^j \|^2 \leq \| l^{j-1} \|^2 \left(1 - \frac{\gamma^2}{4 \Gamma^2}\right),
\end{equation*}
from where, as in Lemma 11.4, it follows the estimate:
\begin{equation*}
    m \leq 1 + 8 \frac{\Gamma^2}{\gamma^2} \ln \frac{\Gamma}{\gamma}.
\end{equation*}
Lemma is proven.
\begin{lemma}
\label{lem:11.6}
If in the direction finding problem $ P = \{ a^1, ..., a^N \} $, then
\begin{equation*}
    m_2 (P, L3) \leq N + 1.
\end{equation*}
\end{lemma}

{\it P r o o f.} For the $ l^j $ calculated by the L3 rule, it is true
\begin{equation*}
    \left<g^j, l^j\right> \geq \| l^j \|^2.
\end{equation*}
Obviously, $ g^{N+1} \in \{ a^1, ..., a^N \} $; so $ l^{N+1} = l^N $ and $ \left<g^{N+1}, l^N\right> \leq \| l^N \|^2 $, i.e.
\begin{equation*}
    m_2 (P, L3) \leq N + 1.
\end{equation*}
The lemma is proven.

A l g o r i t h m A2. In algorithm A1, to calculate directions at step 3 and to check conditions at step 7, instead of vectors $ g^1, ..., g^{j+1} $, we can use normalized values $ \bar{g}^1, ..., \bar{g}^{j+1} $, where $ \bar{g} = \frac{g}{|| g ||} $ ($ \bar{g} = 0 $ if $ g = 0 $). Then step 8 of the algorithm should be replaced by the following step $8^\prime$.

S t e p $8^\prime$. If
\begin{equation*}
    f(y^{j+1}) > f(x^k) - \theta \delta_i \varepsilon_i^2 \| g^{j+1} \| / (4 \| l^j \|),
\end{equation*}
then $ \varepsilon_{i+1} = \varepsilon_i $, $ \delta_{i+1} = \beta \delta_i $, replace $ i $ by $ i + 1 $, go to step 2.

Proof of convergence of the algorithm thus modified to the set $ \bar{X}_{\varepsilon \delta} \cap \{ x \in E_n \> | \> f(x) \leq f(x^0) \} $, where
\[
    \begin{matrix}
        \bar{X}_{\varepsilon \delta} = \{ x \> | \> \bar{\gamma}_{\delta} (x) \leq \varepsilon \}, \\
        \bar{\gamma}_{\delta} (x) = \inf \{ \| \bar{g} \| \> | \> \bar{g} \in \text{co} \cup_{\{ y \> | \> \| y - x \| \leq \delta \}} \bar{G}_f (y) \}, \\
        \bar{G}_f (y) = \{ \bar{g} \in E_n \> | \> g \in G_f (y) \},
    \end{matrix}
\]
differs slightly from the convergence proof of Algorithm A1.

A l g o r i t h m A3. In lemmas 11.4, 11.5 we prove stronger properties of rules L2, L3 than those used in algorithm A1. Therefore, Algorithm A1 can be modified for the rules L2, L3 of construction of $ l^j $, namely, replace steps 7, 8 with the following:

S t e p $7^{\prime\prime}$. If $ \left<g^{j+1}, l^j\right> < \frac{|| l^j ||^2}{2} $, then replace $ j $ by $ j + 1 $, go to step 3. to step 3.

S t e p $8^{\prime\prime}$. If $ f(y^{j+1}) > f(x) - \delta_i \frac{|| l^j ||}{4} $, then $ \varepsilon_{i+1} = \varepsilon_i $, $ \delta_{i+1} = \beta \delta_i $, replace $ i $ by $ i + 1 $, go to step 2.

Theorem 11.1 remains valid for the A3 algorithm.

Let us estimate the labor complexity of the proposed algorithms in the case of convex function minimization. The labor intensity depends on the number of iterations (cycles) of the algorithm, which are calculated using the variables $ i $, $ j $, $ k $ (see Remark 11.1). Let us consider the simplest case. Let step 2 of the algorithms always take $ \varepsilon_0 = \delta $, $ \delta_0 = \delta $. Then in the convex case $i$-iterations are not performed and always $ \| l^j \| \geq \varepsilon $, $ \delta_i = \delta $.

We define
\[
    \begin{matrix}
        D = \{ x \> | \> f(x) \leq f(x^0) \}, \\
        \Gamma = \sup \{ \| g \| \> | \> g \in G_f(y), y \in D \};
    \end{matrix}
\]
it is obvious that $ \| l^j \| \leq \Gamma $. Thus, for algorithm A1, at $ k $-iterations the value of the convex function $ f(x) $ decreases by at least $ \frac{\theta \delta \varepsilon^2}{2 \Gamma} $. Hence, the total number of $k$-iterations does not exceed the value: 
\begin{equation*}
    \frac{2 \Gamma}{\theta \delta \varepsilon^2} [ f(x^0) - \min_{x \in D} f(x) ].
\end{equation*}
Let us denote by $ d $ the diameter of the region $ D $; then
\begin{equation*}
    [ f(x^0) - \min_{x \in D} f(x) ] \leq d \Gamma.
\end{equation*}
Therefore, the total number of $ k $-iterations does not exceed $ \frac{2 d \Gamma^2}{\theta \delta \varepsilon^2} $. Now everything depends on the number of $ j $-iterations at each $ k $-iteration. Let us denote
\[
    \begin{matrix}
        \Gamma_{\delta} (x) = \sup \{ \| g \| \> | \> g \in G_f(y), \| y - x \| \leq \delta \}, \\
        \gamma_{\delta} (x) = \inf \{ \| g \| \> | \> g \in \text{co} \{ G_f(y) \> | \> \| y - x \|\leq \delta \} \}.
    \end{matrix}
\]
The value $ Q_{\delta} (x) = \frac{\Gamma_{\delta} (x)}{\gamma_{\delta} (x)} $ characterizes the local gullibility of the function $ f(x) $. According to Lemma 11.3, at $ l^j = \frac{1}{j} \sum_{r=1}^j g^r $ the number of $ j $-iterations at each $ k $-iteration does not exceed the value $ 2 Q_{\delta}^2 (x) / (1 - \theta) $.

So, the total labor complexity of reaching the set $ X_{\varepsilon \delta}^{*} $ by algorithm A1, in which $ l^j = \frac{1}{j} \sum_{r=1}^j g^r $ and at step  9 $x^{k+1} = y^{j+1} $ is estimated as follows. The values of the function $ f(x) $ are not computed at all, and the required number of gradient computations is as follows
\begin{equation*}
    N_g^1 \leq \frac{4 d \Gamma^2}{\theta (1 - \theta) \delta \varepsilon^2} \sup_{{x} \notin X_{\varepsilon \delta}^{*}} Q_{\delta}^2 (x) \leq \frac{4 d \Gamma^4}{\theta (1 - \theta) \delta \varepsilon^4}.
\end{equation*}

Similar estimation for  algorithm A3, which uses the L2 procedure and at step  9 $x^{k+1} = y^{j+1} $, has the form
\begin{equation*}
    N_g^2 \leq 16 \frac{d \Gamma}{\delta \varepsilon} \left(\sup_{x\notin X_{\varepsilon \delta}^{*}} \{ Q_{\delta}^2 (x) \ln Q_{\delta} (x) \} + \frac{1}{8}\right) \leq 16 \frac{d \Gamma}{\delta \varepsilon} \left(\frac{\Gamma^2}{\gamma^2} \ln \frac{\Gamma}{\gamma} + \frac{1}{8}\right).
\end{equation*}

If the minimized function $ f(x) $ is piece-wise linear and consists of M pieces, then for the  algorithm A3, in which L3 is used and at step  9 $x^{k+1} = y^{j+1} $, the analogous estimate is of the form
\begin{equation*}
    N_g^3 \leq 16 \frac{d \Gamma}{\delta \varepsilon} M.
\end{equation*}

\section*{$\S$ 12. On the multiextremal optimization}
\label{Sec.12}
\setcounter{section}{12}
\setcounter{definition}{0}
\setcounter{equation}{0}
\setcounter{theorem}{0}
\setcounter{lemma}{0}
\setcounter{remark}{0}
\setcounter{corollary}{0}

Non-convex functions may have many local extrema that do not coincide with the absolute minimum of the function on the admissible set. The algorithms considered so far were intended only for finding local minima. The problem of global optimization of functions is complex and urgent. A review of the results obtained in this area is available in [114, 117, 165, 181]. In this paragraph, some additional results on three approaches in multivariate global optimization are described:
\begin{enumerate}
  \item method of combining local and non-local search;
  \item method of approximation of the initial multiextremal problem by multi-extremal problems of special kind;
  \item method of smoothing local extrema.
\end{enumerate}

  {\bf 1. Combined algorithms.} 
	\label{Sec.12.1}
	A natural way to solve real multi-extremal problems is to reasonably combine local and global optimization algorithms. For example, in combined algorithms [1, 12, 82] the stages of local and global algorithms alternate. In this case, each time the algorithms do not work to the end, and when some conditions are met, a switch from one to the other takes place. The method of enumeration of local extrema [28] can also be considered as a special case of this scheme. There are also other types of combining local and global algorithms [36].

  The convergence of the combined algorithms must be specifically proved, because without certain conditions they may not converge even to local extrema, despite the fact that they are based on convergent algorithms.

  We first give two combined algorithms in which local descent plays the main role in the sense of ensuring convergence to local extrema, and global search is used to find suitable starting points for the local algorithm. Then, conversely, we give an example of an algorithm in which global search plays the main role and local descent is used as an intermediate.

  Let us consider the problem:
  \begin{equation}\label{eqn:12.1}
    f(x) \to \min, \> \> \> x \in E_n,
  \end{equation}
where $ f(x) $ is a multiextremal generalized differentiable function.

  In the following combined algorithm, the generalized gradient method is used for local descent and any other algorithm can be used for non-local search.

{A l g o r i t h m K1.}

{\it S t e p 0.} Select an initial point $x^0\in E_n$, numbers $\sigma$, $\Delta>0$ and
sequences $\{\epsilon_k\}$, $\{\delta_k\}$, $\{\rho_k\}$ such that
$$
\epsilon_k\ge 0,\;\;\; \delta_k\ge 0, \;\;\;\rho_k\ge 0;
$$
\begin{equation}\label{eqn:12.2}
\lim_{k\rightarrow\infty}\epsilon_k=\lim_{k\rightarrow\infty}\delta_k=\lim_{k\rightarrow\infty}\rho_k=0,\;\;\;
\sum_{k=0}^\infty\rho_k =\infty.
\end{equation}
Set $k=0$, $\sigma_0=0$.

{\it S t e p 1.} Calculate the value $f(x^k)$ and an arbitrary vector $g^k\in G_f(x^k)$.

{\it S t e p 2.} If $f(x^k)\ge \min_{0\le r\le k}f(x^r)+\Delta$, then set 
$x^{k+1}=\mbox{argmin}\{f(x)| x=x^r, r\in[0,k]\}$, $\sigma_{k+1}=0$;
replace $k$ by $k+1$ and go to {\it Step 1}.

{\it S t e p 3.} If $\sigma_k\ge\sigma$ or $\|g^k\|\le\epsilon_k$, then find by some
algorithm $\tilde{A}$ a new starting point $\tilde{x}$ such that $f(\tilde{x})\le f(x^k)+
\delta_k$ and set $x^{k+1}=\tilde{x}$, $\sigma_{k+1}=0$; 
replace $k$ by $k+1$ and go to {\it Step 1}.

{\it S t e p 4.} Make a step of the gradient method,
$$
x^{k+1}=x^k-\rho_k g^k/\|g^k\|
$$
and set $\sigma_{k+1}=\sigma_k+\rho_k$. Replace $k$ by $k+1$ and go to {\it Step 1}.

\begin{remark}
\label{rem:12.1}
As the auxiliary algorithm $\tilde{A}$, it can be used, for example, any global search algorithm in the vicinity of point $x^k$. Algorithm $\tilde{A}$ should be efficient; 
in particular, for lack of a better point, we can take $\tilde{x} = x^k$ or
$\tilde{x} = \mbox{argmin } \{f(x) | x = x^r, r\in[0, k]\}$.
\end{remark}

\begin{theorem}\label{th:12.1}
Suppose that in problem (12.1) function f(x) is generalized differentiable and the set $D=\{x\in E_n|f(x) \le f(x^0) + \Delta\}$ is bounded, the sequence of points 
$\{x^k\}$ is constructed by  algorithm K1.
Then all limit points of $\{x^k\}$ belong to $D$, among the minimal ones in terms of the value of $f$ there are points belonging to $X^* =\{x\in E_n| 0\in G_f(x)\}$, and the interval 
$\left[\underline{\lim}_{k\rightarrow\infty}f(x^k),\overline{\lim}_{k\rightarrow\infty}f(x^k)\right]$
is embedded in the set $F^* = \{f(x) | x\in X^* \}.$ If $F^*$ does not contain intervals, then all limit points of $\{x^k\}$ belong to $X^*$, and the number sequence $\{f(x^k)\}$ has a limit.
\end{theorem}
{\it P r o o f.} Note that K1 actually combines three algorithms: gradient method $A_1$, 
return to the record point of the trajectory $A_2$,  and the search $A_3 = \tilde{A}$ for a new starting point $x$. For convenience, we introduce the symbolic function 
$A: k\rightarrow \{A_1, A_2, A_3\},$ 
which shows which of the algorithms, $A_{1}, A_{2}$ or $A_{3}$,  works at the $k$-th iteration, i.e. which algorithm generates the point $x^{k+1 }$. We single out the indices $k_{s}\;(s=0,1, \ldots)$ at which switching from one algorithm to another occurs, i.e., at moments $k_{s}$ either $\sigma_{k_{ s}}=0$, or $\sigma_{k_{s}} \ge \sigma$, or $\left\|g^{k_{s}}\right\| \le \varepsilon_{k_{s}}$, or $f\left(x^{k_{s}}\right) \ge \min _{0 \le k \le k_{s}} f\left( x^{k}\right)+\Delta$. By condition (12.2) and by the construction of the  algorithm $\mathrm{K} 1$, the sequence $\left\{k_{s}\right\}$ is infinite.

The  algorithm $\mathrm{K} 1$ can go beyond the bounded set $D$ by only one step of the gradient method, and since $\lim _{k \rightarrow \infty} \rho_{k}=0$, then all limit the points $\left\{x^{k}\right\}$ belong to D. Thus, the first statement of the theorem is proven.

Let us prove the second statement: among the minimal in value of $f$ limit points of $\left\{x^{k}\right\}$, there are points belonging to $X^{*}$. Let's assume the opposite. Let
\begin{equation}\label{eq:12.3}
\lim _{t \rightarrow \infty} x^{k_{t}}=x^{\prime} \notin X^{*}, \quad f\left(x^{\prime}\right)=\lim _{t \rightarrow \infty} f\left(x^{k_t}\right)=\lim _{k \rightarrow \infty} f\left(x^{k}\right).
\end{equation}
Let us find indices $k_{s(t)} \in\left\{k_{s}\right\}$ such that $k_{s(t)-1} \le k_{t}<k_{s(t )}\;(t=0,1, \ldots)$. Without loss of generality, we can assume that there exists $x^{\prime \prime}=\lim _{t \rightarrow \infty} x^{k_{s(t)}}$.
Let us show that
\begin{equation}\label{eq:12.4}
f\left(x^{\prime}\right)=\lim _{t \rightarrow \infty} f\left(x^{k_t}\right)=\varlimsup_{t \rightarrow \infty} \max _{k_{t} \le k \le k_{s(t)}} f\left(x^{k}\right).
\end{equation}

Let's consider several cases.

Let $A\left(k_{t}\right)=A_{1}$ for all $t$. Let us assume that (12.4) does not hold. Then, due to the continuity of $f(x)$,
$$
\varlimsup_{t \rightarrow \infty} \max _{k_{t} \le k \le k_{s(t)}}\left\|x^{k}-x\right\|>0
$$
and therefore,
\begin{equation}\label{eq:12.5}
\varlimsup_{t \rightarrow \infty} \sum_{k=k_{t}}^{k_{s}(t)-1} \rho_{k}>0 .
\end{equation}
It can be considered that the generalized gradient method, as it were, repeatedly starts from the points $x^{k_{t}}\;(t=0,1, \ldots)$. Under condition (12.5), Lemma 9.3 implies the existence of limit points of $\left\{x^{k}\right\}$ with $f$ values less than $f(x^{\prime})$, which contradicts ( 12.3).

If $A\left(k_{t}\right)=A_{2}$ for all $t$, then $k_{s(t)}=k_{t}+1$ and (12.4) is satisfied due to the fact that
$$
f\left(x^{k_{t}+1}\right) \le f\left(x^{k_{t}}\right)-\Delta .
$$
If $A\left(k_{t}\right)=A_{3}$ for all $t$, then $k_{s(t)}=k_{t}+1$ and (12.4) is valid due to conditions
$$
f\left(x^{k_{t}+1}\right) \le f\left(x^{k_{t}}\right)+\delta_{k_{t}}, \quad \lim _{t \rightarrow \infty} \delta_{k_{t}}=0.
$$
Now it is obvious that in the general case, when $A\left(k_{t}\right) \in\left\{A_{1}, A_{2}\right.$, $\left.A_{3} \right\}$, property (12.4) is satisfied.

Without loss of generality, we assume that $\lim_{t \rightarrow \infty} x^{k_{s(t)}}=x^{\prime \prime}$. By the contrary assumption, $x^{\prime \prime} \notin X^{*}$, and by (12.4), $f\left(x^{\prime}\right) \ge f\left(x^{\prime \prime}\right)$. We can assume that algorithms $A_{1}, A_{2}$, or $A_{3}$ start from the points $\left\{x^{k_{s(t)}}\right\}$. If $A\left(k_{s(t)}\right)=A_{1}$ for all $t$, then either $\sigma_{k_{s(t)+1}} \ge\sigma>0$ or 
$\|g^{k_{s(t)+1}}\|\le \epsilon_{k_{s(t)+1}}.$
In any case, 
$\underline{\lim}_{t\rightarrow\infty}\|x^{k_{s(t)+1}}-x^{\prime \prime} \|>0$, and hence,
$$
\sum_{k=k_{s(t)}}^{k_{s(t)+1}-1} \rho_{k} \ge c^{\prime}>0 .
$$

From this, Lemma 9.3 implies the existence of other limit points of the sequence $\left\{x^{k}\right\}$ with smaller values of $f$ than  $f(x^{\prime})$ and $f(x^{\prime \prime})$. We have arrived to a contradiction.

If $A\left(k_{s(t)}\right)=A_{2}$ for all $t$, then by virtue of the inequality $f\left(x^{k_{s(t)}+1} \right) \le$ $\le f\left(x^{k_{s}(t)}\right)-\Delta$ point $x^{\prime}$ is not the minimum limit point of $\left\{ x^{k}\right\}$ (by value of $f$ ). We got a contradiction.

If $A\left(k_{s(t)}\right)=A_{3}$ for all $t$, then $k_{s(t)+1}=k_{s(t)}+1$ . Without loss of generality, we can assume that $\lim x^{k_{s(t)}+1}=x^{\prime \prime \prime} \notin X^{*}$. Since $f\left(x^{k_{s(t)}+1}\right) \le f\left(x^{k_{s}(t)}\right)+\delta_{k_{s (t)}}$ and $\lim _{t \rightarrow \infty} \delta_{k_{s(t)}}=0$, then
\begin{equation}\label{eq:12.6}
f\left(x^{\prime}\right) \ge f\left(x^{\prime \prime}\right) \ge \lim _{t \rightarrow \infty} f\left(x^{k_{s}(t)+1}\right)=f\left(x^{\prime \prime\prime}\right).
\end{equation}
Now, according to the construction of the combined algorithm K1, at points $x^{k s(t)+1}$ only algorithms $A_{1}$ or $A_{2}$ can be used for sufficiently large $t$, which, as shown above , generate limit points $\left\{x^{k}\right\}$ with smaller values of $f$ than  $f(x^{\prime})$. Again we came to a contradiction.

It is easy to see that the general case, when $A\left(k_{t}\right) \in\left\{A_{1}, A_{2}, A_{3}\right\}$, reduces to the particular considered cases. Thus, the contradictions obtained above, prove the second statement of the theorem.

Let us prove the third assertion. Let's denote
$$
\underline{f}=\lim _{k \rightarrow \infty} f\left(x^{k}\right), \quad \bar{f}=\varlimsup_{k \rightarrow \infty} f\left(x^{k}\right) .
$$
Let us show that $[\underline{f}, \bar{f}) \subset F^{*}$; assume the opposite. Then there is a number $c$ such that $c \in[\underline{f}, \bar{f}]$ and $c \notin F^{*}$. According to the proven second statement of the theorem, $\underline{f}\in F^{*}$; hence $c>\underline{f}$. Let us choose a number $d$ so that $\underline{f}<c<d<\bar{f}$. The numerical sequence $\left\{f\left(x^{k}\right)\right\}$ crosses the interval $[c, d]$ from $c$ to $d$ an infinite number of times; so one can single out indices $\left\{k_{t}\right\}$ and $\left\{m_{t}\right\}$ such that $\lim _{t \rightarrow \infty} x^{ k_{t}}=x^{\prime}$ and either
\begin{equation}\label{eq:12.7}
f\left(x^{k_{t}}\right) \le c<f\left(x^{k}\right)<d \le f\left(x^{m_{t}}\right), \quad k \in\left(k_{t}, m_{t}\right)
\end{equation}
or
$$
f\left(x^{k_{t}}\right) \le c<d \le f\left(x^{k_{t}+1}\right) .
$$

For the subsequence $\left\{x^{k_t}\right\}$ we can prove the validity of relation (12.4) in the same way as was done above. From (12.4) and (12.7) it follows that for large $t$ $k_{s(t)}<m_{t}$ and
\begin{equation}\label{eq:12.8}
f\left(x^{\prime}\right)=\lim _{t \rightarrow \infty} f\left(x^{k_{t}}\right)=\lim _{t \rightarrow \infty} f\left(x^{k_{s(t)}}\right)=c.
\end{equation}

Without loss of generality, we will assume that $\lim _{t \rightarrow \infty} x^{k_{s}(t)}=x^{\prime \prime}$. From (12.8), due to the choice of $c \notin F^{*}$, it follows that $x^{\prime \prime} \notin X^{*}$. Now we can assume that the algorithm $A_{1}, A_{2}$ or $A_{3}$ starts from points $x^{k_{s(t)}}(t=0,1, \ldots)$. Let's consider a number of special cases.

If $A\left(k_{s(t)}\right)=A_{1}$ for all $t$, then by construction of the combined algorithm either $\sigma_{k_{s}(t)+1} \ge \sigma>0$, or $\left\|g^{k_{s(t)}+1}\right\| \le \varepsilon_{k_{s}(t)+1}$. From here, taking into account the fact that $x^{\prime \prime} \notin X^{*}$, it follows
$$
\lim _{t \rightarrow \infty}\left\|x^{k_{s(t)+1}}-x^{\prime \prime}\right\|>0, \quad \sum_{k=k_{s(t)}}^{{k_{s(t)+1}-1}} \rho_{k} \ge c^{\prime}>0 .
$$
We can assume that the gradient descent process starts multiple times from points of the sequence $\left\{x^{k_s(t)}\right\}$, which converges to $x^{\prime \prime} \notin X^ {*}.$ In this situation, the Lemma 9.3 applies to the sequence $\left\{x^{k}\right\}_{k \ge k_{s(t)}}$ $(t=0,1, \ldots)$. This lemma involves the parameter $\varepsilon$, which we choose so that 
$\max \left\{f(x) |\|x-x^{\prime \prime}\| \le \varepsilon\right\}<d$. Then, according to the lemma, there is a subsequence $\left\{x^{l_{t}}\right\}$ such that $\left\|x^{k}-x^{\prime \prime}\right\| \le \varepsilon$ for $k \in\left(k_{s(t)}, l_{t}\right)$ and
\begin{equation}\label{eq:12.9}
f\left(x^{\prime \prime}\right)=\lim _{t \rightarrow \infty} f\left(x^{k_{s(t)}}\right)>\overline{\lim _{t \rightarrow \infty}} f\left(x^{l_t}\right) .
\end{equation}
The resulting relation contradicts (12.7).

If $A\left(k_{s(t)}\right)=A_{2}$ for all $t$, then
\begin{equation}\label{eq:12.10}
f\left(x^{k_{s(t)}+1}\right) \le f\left(x^{k_{s}(t)}\right)-\Delta,
\end{equation}
which also contradicts (12.7).

Let $A\left(k_{s(t)}\right)=A_{3}$ for all $t$. Without loss of generality, we can assume that $\lim _{t \rightarrow \infty} x^{k_{s(t)+1}}=x^{\prime \prime \prime}$. For $\left\{x^{k_{s(t)+1}}\right\}$ relation (12.6) holds. From (12.6) and $(12.7)$ it follows that for sufficiently large $t$ $k_{s(t)+1}<m_{t}$ and
$$
f\left(x^{\prime}\right)=f\left(x^{\prime \prime}\right)=\lim _{t \rightarrow \infty} f\left(x^{k_{s(t)+1}}\right)=c
$$
whence, due to the choice of $c \notin F^{*}$, it follows that $x^{\prime\prime\prime} \notin X^{*}$. Now, by construction of the combined algorithm, for sufficiently large $t$ at points $x^{k_{s(t)+1}}$, only algorithms $A_{1}$ or $A_{2}$ can be used. These algorithms, as was already shown above for a similar situation (see (12.9), (12.10)), have minimizing properties that conflict with construction (12.7).

It is easy to see that the general case, when $A\left(k_{s(t)}\right) \in\left\{A_{1}, A_{2}, A_{3}\right\}$, reduces to the study of the considered special cases. Thus, the contradictions obtained above, prove that $[\underline{f}, \bar{f}] \subset F^{*}$. Now, due to the boundedness of the set $D$ and the closedness of the mapping $G_{f}$, the set $F^{*}$ is closed; therefore $[\underline{f}, \bar{f}] \subset F^{*}$. The third statement of the theorem has been proven.

Now let $F^{*}$ contains no intervals. Since $[\underline{f}, \bar{f}] \subset F^{*}$, then $\underline{f}=\bar{f}$, i.e. the number sequence $\left\{f\left (x^{k}\right)\right\}$ has a limit.

Let us show that all limit points of $\left\{x^{k}\right\}$ belong to $X^{*}$. Let's assume the opposite; then there is a subsequence $\left\{x^{k_t}\right\}$ such that $\lim _{t \rightarrow \infty} x^{k_t}=x^{\prime} \notin X^{* }$. Consider indices $\left\{k_{s(t)}\right\}$ such that $k_{s(t)-1}<k_{t} \le k_{s(t)}\;(t=0 ,1, \ldots)$. Without loss of generality, we can assume that $\lim _{t \rightarrow \infty} x^{k_{s(t)-1}}=x^{\prime \prime}$.

Note that $A\left(k_{s(t)-1}\right)=A_{2}$ holds only for a finite number of indices $t$, since otherwise we obtain a contradiction with the convergence of $\left\{ f\left(x^{k}\right)\right\}$.

If $A\left(k_{s(t)-1}\right)=A_{3}$ for an infinite set of numbers $t$, then for the same sufficiently large $t$ $A\left(k_{s( t)}\right)=A_{1}$. Recall that $x^{\prime} \notin X^{*}$. Since either $\sigma_{k_{s(t)+1}} \ge \sigma$ or $\left\|g^{k_{s(t)+1}}\right\| \le \varepsilon_{k_{s(t)+1}}$, then $\rho_{k_{s(t)}}+\ldots+\rho_{k_{s(t)+1}} \ge c^ {\prime}>0$. Lemma 9.3 is applicable to the sequence $\left\{x^{\left.k_{s(t)}\right\}}\right.$, the statement of which contradicts the already proven convergence of $\left\{f\left(x^ {k}\right)\right\}$.

Let $A\left(k_{s(t)-1}\right)=A_{1}$ for an infinite set of indices $t$; without loss of generality, we can assume that this holds for all $t$. It holds $x^{\prime \prime} \in X^{*},$since otherwise, applying  Lemma 9.3
to the sequence of initial points 
$\left\{ x^{k_{s(t)-1}}\right\}$, we obtain a contradiction with the already proven convergence of $\left\{f\left(x^{k}\right)\right\}$.

Since the set $X^{*}$ is closed, there exists $\delta>0$ such that $\rho\left(x^{\prime}, X^{*}\right) \ge \delta>0$.

Sequences $\left\{x^{k_{s}(t)-1}, \ldots, x^{k_t}\right\}\;(t=0,1, \ldots)$ move from $x^{\prime \prime}$ to $x^{\prime}$, and
$$
\begin{gathered}
\lim _{t \rightarrow \infty}\left\|x^{k_{s(t)-1}}-x^{\prime}\right\|=\left\|x^{\prime \prime}-x^{\prime}\right\|, \quad \lim _{t \rightarrow \infty}\left\|x^{k_{t}}-x^{\prime}\right\|=0, \\
\left\|x^{\prime \prime}-x^{\prime}\right\| \ge \rho\left(x^{\prime}, X^{*}\right) \ge \delta>0 .
\end{gathered}
$$
Number sequences 
$\left\{\left\|x^{k_{s(t)-1}}-x^{\prime}\right\|,
\left\|x^{k_{s(t)-1}+1}-x^{\prime}\right\|, \ldots,
\left\|x^{k_{t}}-x^{\prime}\right\|\right\}$ $(t=0,1, \ldots)$, 
intersect level $\delta$ an infinite number of times; therefore there is a subsequence $\left\{x^{k_{r(t)}}\right\}$ such that 
$k_{{s(t)}-1}<k_{r(t)}< k_{t}$ and for sufficiently large $t$ 
$\left\|x^{k_{r(t)}}-x^{\prime}\right\| \ge \delta$ and $\left\|x^{k_{r(t)}+1}-x^{\prime}\right\|<\delta$.

Without changing the notation, we will assume that $\lim _{t \rightarrow \infty} x^{k_{r(t)}}=x^{\prime \prime \prime}$. Since $\lim _{k \rightarrow \infty} \rho_{k}=0$, then
$$
\lim _{t \rightarrow \infty}\left\|x^{k_{r(t)}}-x^{\prime}\right\|=\left\|x^{\prime \prime}-x^{\prime}\right\|=\delta<\rho\left(x^{\prime}, X^{*}\right) .
$$
From here it follows that $x^{\prime \prime \prime} \notin X^{*}$.

Denote $\Gamma=\sup \left\{\|g\| \mid g \in G_{f}(y), y \in D\right\}$. For sufficiently large $t$, we have
$$
\frac{\delta}{2} \le\left\|x^{k_{t}}-x^{k_{r(t)}}\right\| \le \Gamma \sum_{k=k_{r(t)}}^{k_{t}-1} \rho_{k},
$$
hence,
$$
\sum_{k=k_{r(t)}}^{k_{t}-1} \rho_{k} \ge \frac{\delta}{2 \Gamma}.
$$
We can assume that the generalized gradient method starts repeatedly from the points $x^{k_{r(t)}} \rightarrow x^{\prime \prime \prime} \notin X^{*}$. In this situation, Lemma 9.3 is again applicable, from which, in particular, the divergence of the sequence $\left\{f\left(x^{k}\right)\right\}_{k=0}^{\infty}$ follows , which contradicts the previously established convergence of this sequence. The theorem is completely proven. 

A l o r i t h m $\mathrm{K} 2$.

Recall the relaxation algorithms A1--A3 from $\S$ 11. They generate a finite sequence of approximations to the set of approximate solutions $X_{\varepsilon \delta}^{*}$ of problems $(11.1),(11.2)$. At the current point, the direction $l$ of decrease of the minimized function and the point $y=x^{k}-\delta l /\|l\|$ are found in them. At step 9 in these algorithms, one can assume $x^{k+1}=y$, or one can, as noted, use some auxiliary algorithm to find a point $x^{k+1}$ such that $f\left( x^{k+1}\right) \le f(y)$. This auxiliary algorithm can be not only a linear search algorithm in the direction $l$, but also some (possibly heuristic) global search algorithm in the neighborhood of points $x^{k}$ or $y$.

Thus, in essence, the relaxation algorithms $\mathrm{Al}-\mathrm{A} 3$ from $\S 11$ are combined; in them, after each relaxation step, it is allowed to use a global search to find a point $x^{k+1}$ that is no worse than the already reached point $y$. In $\S$ 11 the local convergence of such a combined algorithm is proved.

Let us now consider a combined algorithm in which the main role is played by a limited global search in the vicinity of the current approximation.

Let it be necessary to solve the problem
$$
f(x) \rightarrow \min
$$
subject to constraints
$$
x \in D \subset E_{n}
$$
where $f(x)$ is a continuous function such that the sets $\left\{x \in E_{n} \mid x \in D\right.$, $f(x) \le c\}$ are bounded for any $c \in E_{1}$.

Let $K(x)$ be a multivalued mapping such that the mapping $x \rightarrow K(x) \cap D$ is lower semicontinuous, i.e., from $x^{k} \rightarrow x$ and $y \in K (x) \cap D$ implies the existence of $y^{k} \in K\left(x^{k}\right) \cap D$ such that $y^{k} \rightarrow y$.

For example, if $D=E_{n}$ and $K(x)$ are lower semicontinuous, then $K(x) \cap D=K(x)$ is lower semicontinuous. In particular, the mapping
$$
x \rightarrow K(x)=\left\{y \in E_{n} \mid y-x \in K\right\},
$$
where $K$ is an arbitrary set in $E_{n}$, semicontinuous both from above and below.
\begin{lemma}
\label{lem:12.1}
Let $D$ be aconvex set, $K_{r}(x)$ be a ball (in some norm) of radius $r$ centered at point $x$. Then the multivalued mapping $x \rightarrow K_{r}(x) \cap D$ is lower semicontinuous in $D$.
\end{lemma}

{\it P r o o f.} Let $x^{k} \rightarrow x, x^{k} \in D$ and $y \in K_{r}(x) \cap D$. Let us show that there exist $y^{k} \in K_{r}\left(x^{k}\right) \cap D$ such that $\left\{y^{k}\right\} \rightarrow y$.

Let's denote
$$
\begin{gathered}
t_{k}=\sup \left\{t \mid 0 \le t \le 1, \quad x^{k}+t\left(y-x^{k}\right) \in K_{r}\left(x^{k}\right)\right\}, \\
y^{k}=x^{k}+t_{k}\left(y-x^{k}\right)=\left(1-t_{k}\right) x^{k}+t_{k} y .
\end{gathered}
$$
Since $x^{k}, y \in D$, then $y^{k} \in D$. Thus, $y^{k} \in K_{r}\left(x^{k}\right) \cap D$.

Let us show that $\lim _{k \rightarrow \infty} t_{k}=1$; then $\left\|y^{k}-y\right\|=\left(1-i_{k}\right)\left\|y-x^{k}\right\| \rightarrow 0$ for $k \rightarrow \infty$, and the lemma would be proven. Assuming the contrary, then there is a subsequence $t_{k_{s}} \le 1-\varepsilon, \quad 0.5 \le \varepsilon<1$. Let's take $t^{\prime}=$ $=1-\varepsilon / 2, t_{k_{s}}<t^{\prime}$. Consider $y^{\prime}=x^{k_{s}}+t^{\prime}\left(y-x^{k_{s}}\right)$. The assessment holds true
$$
\left\|x^{k_{s}}-y^{\prime}\right\|=t^{\prime}\left\|y-x^{k_{s}}\right\| \le t^{\prime}\left(\|y-x\|+\left\|x^{k_{s}}-x\right\|\right)
$$
For $\left\|x^{k_{s}}-x\right\| \le \frac{1-t^{\prime}}{t^{\prime}}\|y-x\| $ it takes place $\left\|x^{k_{s}}-y\right\| \le\| y-x \| \le r$. Therefore, for sufficiently large $s$, there exists $t^{\prime}>t_{k_{s}}$ such that $y^{\prime}=x^{k_{s}}+t^{\prime }\left(y-x^{k_{s}}\right) \in K_{r}\left(x^{k_{s}}\right)$, that contradicts the choice of $$t_{k_{s}}=\sup \left\{t \mid 0 \le t \le 1, \quad x^{k_{s}}+t\left(y-x^{k_{s}}\right) \in K_{r}\left (x^{k_{s}}\right)\right\}.$$ The lemma is proven.

Consider the following idealized algorithm.

A l g o r i t h m $\mathrm{K} 3$.

Let an arbitrary point $x^{0} \in D$ and a number $r>0$ be given. Let us find the point $x^{1}$ of the global minimum $f(x)$ on the set $K_{r}\left(x^{0}\right) \cap D$. From the point $x^{1}$ we descend by any intermediate (in particular, the same) algorithm so that we get a point 
$x^{2}, f\left(x^{2}\right) \le f\left(x^{1}\right)$. 
Then again we find the point $x^{3}$ of the global minimum of $f(x)$ on the set $K_{r}\left(x^{2}\right) \cap$ $\cap D$. From it we will again carry out descent using some intermediate algorithm, we will obtain a point $x^{4}, f\left(x^{4}\right) \le f\left(x^{3}\right)$, etc.
\begin{theorem}
\label{th:12.2}
Let the sequence $\left\{x^{k}\right\}$ be generated by the algorithm $\mathrm{K} 3$, where the mapping $x \rightarrow K_{r}(x) \cap D$ is lower semicontinuous. Let $\left\{x^{k_{m}}\right\}$ denote the subsequence of those points $x^{k}$ that in $K_{r}\left(x^{k}\right)$-vicinity of which a global search was carried out. Then the sequence $\left\{x^{k}\right\}$ is bounded, and if $x^{\prime}$ is some limit point of $\left\{x^{k_{m}}\right \}$, then $x^{\prime}$ is the global minimum point of $f(x)$ on the set $K_{r}\left(x^{\prime}\right) \cap D$.
\end{theorem}

{\it P r o o f.} The K3 is a relaxation algorithm, i.e. the sequence $\left\{f\left(x^{k}\right)\right\}$ decreases monotonically. The sequence $\left\{x^{k}\right\}$ lies in the bounded set $\left\{x \mid x \in D, f(x) \le f\left(x^{0}\right )\right\}$; therefore $\left\{f\left(x^{k}\right)\right\}$ is bounded below. Therefore, there is a finite limit $\lim_{k\rightarrow\infty} f\left(x^{k}\right)$.

Now suppose the opposite: there is a limit point $x^{\prime}=\lim x^{k_{m_{s}}}$ of the sequence $\left\{x^{k_{m}}\right\}$, which is not the global minimum point of $f(x)$ on $K_{r}\left(x^{\prime}\right) \cap D$. Then there is a point $x^{\prime \prime} \in K_{r}\left(x^{\prime}\right) \cap D$ such that $f\left(x^{\prime \prime} \right) < f\left(x^{\prime}\right)$; it is obvious that $\left\|x^{\prime \prime}-x^{\prime}\right\| \le r$. Due to the lower semicontinuity of the mapping $x \rightarrow K_{r}(x) \cap D$, there is a sequence $y^{s} \rightarrow x^{\prime \prime}, y^{s} \in K\left( x^{k_{m_s}}\right) \cap D$. Let's choose $\varepsilon<f\left(x^{\prime}\right)-f\left(x^{\prime \prime}\right)$. Since $f\left(x^{k_{m_{s}}}\right) \rightarrow f\left(x^{\prime}\right), \quad f\left(y^{s}\right ) \rightarrow f\left(x^{\prime \prime}\right)$, then for sufficiently large $s$
$$
f\left(x^{k_{m_{s}}}\right)-f\left(y^{s}\right) \ge \varepsilon,
$$
whence $f\left(x^{k_{m_{s}}}\right) \ge f\left(y^{s}\right)+\varepsilon \ge f\left(x^{k_{m_s} +1}\right)+\varepsilon$, which contradicts the convergence of the sequence $\left\{f\left(x^{k}\right)\right\}$. The theorem has been proven.

The combined  algorithm K3 is idealized because it does not specify stopping conditions for the global search procedures. Let's fill this gap.

Let the global minimization algorithm for $f(y)$ on the set $K_{r}\left(x^{k}\right) \cap D$ generates a sequence of approximations $\left\{y^{s}\right\}$ and the estimate of the accuracy of the approximations is known:
$$
f\left(y^{s}\right)-f_{k}^{*} \le \varepsilon_{s} \rightarrow 0, \quad s \rightarrow \infty,
$$
where $f_{k}^{*}=\min \left\{f(y) \mid y \in K_{r}\left(x^{k}\right) \cap D\right\}$; for the algorithms considered in the next section of the paragraph, such estimate exists. Then we can restrict ourselves to finding $y^{s} \in K_{r}\left(x^{k}\right) \cap D$ such that
$$
\varepsilon_{s} \le \alpha\left(f\left(x^{k}\right)-f\left(y^{s}\right)\right), \quad \alpha \ge 0 .
$$
From the inequality
$$
f\left(y^{s}\right)-f_{k}^{*} \le \varepsilon_{s} \le \alpha\left(f\left(x^{k}\right)-f\left(y^{s}\right)\right)
$$
with accounting of $f\left(x^{k+1}\right) \le f\left(y^{s}\right)$, it follows
$$
f\left(x^{k+1}\right) \le f\left(y^{s}\right) \le \frac{\alpha}{\alpha+1} f\left(x^{k}\right)+\frac{1}{\alpha+1} f^*_k.
$$
This is sufficient for Theorem 12.2 to remain valid in the case under consideration.

{\bf 2. The approximation method.}
\label{Sec.12.2}
 In this section, we extend the method of nonconvex approximations of multi-extremal problems [29, 95, 96] to multidimensional multi-extremal problems with nonconvex constraints and outline some ways to solve the resulting auxiliary multi-extremal problems of a special kind.

Let $f(x)\left(x \in D \subset E_{n}\right)$ be a continuous function. Let us assume that there is a family $\varphi(y, x)(y \in D)$ of lower approximations of the function $f(x)$  such that
\begin{enumerate}
   \item $f(x) \ge \varphi(y, x)$ for all $y, x \in D$

   \item $f(y)=\varphi(y, y)$ for all $y \in D$

   \item $\varphi(y, x)$ are locally Lipschitz in $x$ uniformly in $y$, i.e. for any compact set $K \subset D$ there is a constant $L_{K}$ such that for any $ x, y, z \in K$ it holds
$$
|\varphi(y, x)-\varphi(y, z)| \le L_{K}\|x-z\| \text {. }
$$
\end{enumerate}
If $f(x) \quad(x \in D)$ is differentiable and its gradient $\nabla f(x)$ satisfies the Lipschitz condition in $D$ with constant $L$, then as a lower approximation to $f(\cdot)$ the tangent to graph of $f(\cdot)$ at point $y$ paraboloid can be taken:
\begin{equation}\tag{12.4'}
\varphi(y, x)=f(y)+\frac{1}{2 L}\|\nabla f(y)\|^{2}-\frac{L}{2}\left\|x-y-\frac{1}{L} \nabla f(y)\right\|^{2} .
\end{equation}
If $f(x)$ is Lipschitz in $D$ with Lipschitz constant $L$, then as lower approximations of $f(\cdot)$ we can take the functions
\begin{equation}\tag{12.5'}
\varphi(y, x)=f(y)-L\|y-x\| \text {. }
\end{equation}
Let us describe the method of non-convex approximations [96] for solving the following multiextremal problem:
\begin{equation}\tag{12.6'}
f(x) \rightarrow \min
\end{equation}
subject to
\begin{equation}\tag{12.7'}
x \in D \subset E_{n}
\end{equation}
where $f(x)$ is a continuous function on $D$ that has lower approximations $\varphi(y, x)$ satisfying conditions 1)-3), $D$ is a closed bounded set.

Let's take an arbitrary point $x^{0} \in D$. Let's find the point $x^{1} \in D$ with the lowest possible value of the function $\varphi_{0}(x)=\varphi\left(x^{0}, x\right)$. Then we find the point $x^{2} \in D$ with the lowest possible value of the function $\varphi_{1}(x)=$ $=\max \left(\varphi\left(x^{0}, x\right ), \varphi\left(x^{1}, x\right)\right)$, etc.

Let's denote
$$
\varphi_{k}(x)=\max _{0 \le i \le k} \varphi\left(x^{i}, x\right)
$$
In the described method, it is necessary to solve a sequence of special multiextremal problems of the form
\begin{equation}\tag{12.8'}
\varphi_{k}(x) \rightarrow \min
\end{equation}
subject to
\begin{equation}\tag{12.9'}
x \in D \subset E_{n}
\end{equation}
\begin{lemma}\label{le:12.2}. 
It takes place
$$
\lim _{k \rightarrow \infty}\left[f\left(x^{k}\right)-\varphi_{k-1}\left(x^{k}\right)\right]=0 .
$$
\end{lemma}
\begin{theorem}\label{th:12.3}
Let the points $x^{k}$ be the global minimum points of the functions $\varphi_{k-1}(x)\;(k=1,2, \ldots)$. Then all limit points of $\left\{x^{k}\right\}_{k=0}^{\infty}$ are global minimum points of $f(x)$ on $D$. 
\end{theorem}
Lemma 12.2 and Theorem 12.3 in a more general situation will be proved below.
\begin{remark}\label{re:12.2}
The quantity $\varepsilon_{k}=f\left(x^{k}\right)-\varphi_{k-1}\left(x^{k}\right)$ gives an estimate of the accuracy of the approximation $x^{k} $:
$$
f\left(x^{k}\right)-\min _{x \in D} f(x) \le f\left(x^{k}\right)-\varphi_{k-1}\left(x^{k}\right) .
$$
\end{remark}

Let us extend the described method to solve problems of the following type:
$$
f_{0}(x) \rightarrow \min
$$
subject to
$$
f_{i}(x) \le 0, \quad i=1,2, \ldots, m, \quad x \in D \subset E_{n},
$$
where $f(x)$ and $f_{i}(x)\;(i=1,2, \ldots, m)$ are continuous functions having lower approximations $\varphi_{i}(y, x)(i= 0,1, \ldots, m)$; the set $D$ is closed and bounded.

First, we will combine all inequality constraints into one (although this is not necessary), i.e. we will introduce the function
$$
h(x)=\max _{1 \le i \le m} f_{l}(x)
$$
and consider the equivalent problem:
\begin{equation}\tag{12.10'}
f(x) \rightarrow \min
\end{equation}
subject to
\begin{equation}\label{eq:12.11}
h(x) \le 0, \quad x \in D.
\end{equation}

Note that if
\begin{equation}\label{eq:12.12}
h(y)=\max _{1 \le i \le m} f_{i}(y)=f_{i^{*}(y)}(y)
\end{equation}
and $\varphi_{i^{*}(y)}(y, x)$ is a lower approximation of the function $f_{i^{*}(y)}(x)$, 
then
\begin{equation}\label{eq:12.13}
\psi(y, x)=\varphi_{i^{*}(y)}(y, x)
\end{equation}
is also a lower approximation of the function $h(\cdot)$. In what follows we will assume that the lower approximations $\psi(y, x)$ are known.

The method for solving problem (12.10), (12.11) constructs a sequence of points $\left\{x^{k}\right\} \subset D$ according to the following rule. The point $x^{0} \in D$ is arbitrary. Let the points $x^{0}, x^{1}, \ldots$ $\ldots, x^{k} \in D$ have already been constructed. Let's denote
$$
\begin{gathered}
\varphi_{k}(x)=\max _{0 \le i \le k} \varphi\left(x^{i}, x\right), \\
\psi_{k}(x)=\max _{0 \le i \le k}\left\{\psi\left(x^{i}, x\right) \mid h\left(x^{i}\right)>0\right\} .
\end{gathered}
$$
Then we find the point $x^{k+1}$ as a solution to the following multiextremal problem of a special form:
\begin{equation}\label{eq:12.14}
\varphi_{k}(x) \rightarrow \min
\end{equation}
subject to
\begin{equation}\label{eq:12.15}
\psi_{k}(x) \le 0, \quad x \in D.
\end{equation}
\begin{lemma}
\label{le:12.3}
If the original problem (12.10), (12.11) has no admissible solutions, then for some $k$ the approximating problem (12.14), (12.15) has no solutions. If the original problem has an admissible solution, then the sequence $\left\{x^{k}\right\}$ is infinite, and all limit points $\left\{x^{k}\right\}$ belong to the admissible region (12.11).
\end{lemma}

{\it P r o o f.} Let us show that if the sequence $\left\{x^{k}\right\}$ of admissible solutions to problems (12.14), (12.15) is infinite, then all limit points $\left\{x^{k}\right\} $ belong to the admissible domain of the problem (12.10), (12.11). In this case the region (12.11) is not empty. Thus, if domain (12.11) is empty, then the sequence $\left\{x^{k}\right\}$ cannot be infinite, i.e., for some $k$ the approximating problem (12.14), (12.15) has no feasible solutions.

So, let the sequence of points $\left\{x^{k}\right\}$ from the bounded set $D$ be infinite. Assume that there exists a subsequence $\left\{x^{k}\right\}$ such that
$$
\lim _{m \rightarrow \infty} x^{k_{m}}=x^{\prime}, \quad h\left(x^{\prime}\right)>0 .
$$
The sequences $\left\{\psi_{k}(x)\right\}(x \in D)$ are monotonically increasing and bounded from above, so there exists $\psi(x)=\lim _{k \rightarrow \infty} \psi_{k}(x)$. The functions $\psi\left(x^{k}, x\right)$, and subsequently $\psi_{k}(x)$, are Lipschitz in $D$ uniformly in $k$; therefore the function $\psi(x)$ is also Lipschitz in $D$.

Now, by construction. $\psi_{k_{m}-1}\left(x^{k_{m}}\right) \le 0$; passing here to the limit in $m$, we obtain $\psi\left(x^{\prime}\right) \le 0$. But on the other hand, $\psi\left(x^{k_{m}}\right)\ge \psi_{k_{m}}\left(x^{k_{m}}\right)= h\left(x^{k_{m}}\right)$ for sufficiently large $m$. Passing here to the limit in $m$, we have $\psi\left(x^{\prime}\right) \ge h\left(x^{\prime}\right)>0$. The resulting contradiction proves the lemma.
\begin{lemma}
\label{le:12.4}
Let the sequence $\left\{x^{k}\right\}$ of admissible solutions of problem (12.14), (12.15) be infinite. Then
$$
\lim _{k \rightarrow \infty}\left[f\left(x^{k}\right)-\varphi_{k-1}\left(x^{k}\right)\right]=0.
$$
\end{lemma}
{\it P r o o f}. Let's assume the opposite. Then there is \(\{x^{{k_m}}\}\) such that \(\mathop {\lim }\limits_{m \to \infty } {x^{{ k_m}}} = x'\) and
\[\mathop {\lim }\limits_{m \to \infty } \left[ f\left( {x^{k_m}} \right) - {\varphi _{{k_m} - 1 }}\left( {{x^{{k_m}}}} \right) \right] = \varepsilon > 0.\]
The sequences \(\left\{ {{\varphi _k}(x)} \right\}\) \((x \in D)\) are monotonically increasing and bounded above by the function \(f(x)\); so \(\varphi (x) = \mathop {\lim }\limits_{k \to \infty } {\varphi _k}(x)\) exists.
The functions \(\{{{\varphi _k}(x)}\}\) are obviously Lipschitz in D uniformly in k; so \(\varphi (x)\) is also Lipschitz in D. Then
\[\mathop {\lim }\limits_{m \to \infty } \left[ f\left( x^{{k_m}} \right) - {\varphi _{{k_m} - 1 }}\left( {{x^{{k_m}}}} \right) \right] = f(x') - \varphi (x') = \varepsilon > 0.\]
But, on the other hand, $\varphi \left( x^{{k_m}} \right) \ge \varphi _{{k_m}}\left( {{x^{{k_m} }}} \right) = f\left( {{x^{{k_m}}}} \right);$ passing here to the limit in $m$, we have \(\varphi (x') \ge f(x')\). The resulting contradiction proves the lemma.
The proof of Lemma 12.2 does not differ from the proof of Lemma 12.4.
\begin{theorem}\label{th:12.4}
Let the point \({x^{k + 1}}\) be the global minimum point of problem (12.14), (12.15); the sequence \(\left\{ {{x^k}} \right\}\) is infinite. Then all limit points of \(\left\{ {{x^k}} \right\}\) are points of the global solution of problem (12.10), (12.11).
\end{theorem}
{\it P r o o f.} Let's assume the opposite. Then there is a sequence 
\(\left\{ x^{{k_m}}\right\}\) such that\(\mathop {\lim }\limits_{m \to \infty } {x^{ {k_m}}} = x'\) and for some admissible point \(x''\) of problem (12.10), (12.11) we have \(f(x') > f(x'')\). By Lemma 12.4,
\[\mathop {\lim }\limits_{m \to \infty } {\varphi _{{k_m} - 1}}\left( {{x^{{k_m}}}} \right) = f(x ') > f(x'').\]
But, on the other hand, 
\({\varphi _{{k_m} - 1}}\left( x^{{k_m}} \right) \le f(x'')\) 
and \( \overline{\lim }_{m \to \infty } {\varphi _{{k_m} - 1}}\left( x^{{k_m}} \right) \le f(x'')\) The resulting contradiction proves the theorem.

Theorem 12.3 is proved similarly.

In the method under consideration, it is necessary to solve a sequence of special multiextremal problems (12.8), (12.9) or (12.4), (12.5). Let's outline some ways to solve them. To do this, we will have to transform the extremal problems.

Two problems will be called equivalent if there exists an algorithm that allows to reconstruct the solution of one problem from the solution of the other. Equivalent problems either have solutions simultaneously or do not have solutions simultaneously.

Two minimization problems with feasible solutions are known to be equivalent if there exists an algorithm that, for any feasible point of one problem, can identify a feasible point of the other with a no greater objective function value.

Indeed, in this case, the optimal values of the objective functions of the problems coincide, and the existing correspondences between the points of the feasible regions allow us to reconstruct the solutions of one problem based on the solution of the other.

Consider auxiliary problems (12.8), (12.9).

Let \(f(x)\) in (12.6) admit lower approximations of the form (12.4) with some constant L. In particular, if \(f(x)\) is a maximum function, then lower approximations can be constructed using formulas of the type ( 12.12), (12.13).

Now problem (12.8), (12.9) has the form:
\[\mathop {\max }\limits_{1 \le i \le k} \left({a_i} - \frac{L}{2}{\left\| {x - \left. {{b^i}} \right\|} \right.^2}\right) \to \mathop {\min }\limits_x \]
subject to
\[x \in D \subset {E_n},\]
where
\[{a_i} = f({x^i}) + \frac{1}{{2L}}{\left\| {\nabla f({x^i})} \right\|^2},\] \[{b^i} = {x^i} + \frac{1}{L}\nabla f({x ^i}).\]

This task is equivalent to the following one:
\[\Phi (x,{z_1}) \equiv {z_1} \to \mathop {\min }\limits_{x,{z_1}} \]
under restrictions
\[L\left<{b^i},x\right> + {a_i} - \frac{L}{2}{\left\| {{b^i}} \right\|^2} - \frac{L}{2}{\left\| x \right\|^2} \le{z_1},\]
\[i = 1,2,...,k;\;\;\; x \in D.\]
Making the change of variables \(z = {z_1} + \frac{L}{2}{\left\| x \right\|^2}\), we get the following equivalent problem:
\begin{equation}\label{eqn:12.16}
        z - \frac{L}{2}{\left\| x \right\|^2} \to \mathop {\min }\limits_{x,z}\tag{12.16}
\end{equation}
subject to the constraints
\begin{equation}\label{eqn:12.17}
        L\left<{b^i},x\right> - z \le \frac{L}{2}{\left\| {{b^i}} \right\|^2} - {a_i}, 
				\;\;\;i = 1,2,...,k;\tag{12.17}
\end{equation}
\begin{equation}\label{eqn:12.18}
        x \in D.\tag{12.18}
\end{equation}

If the set $D$ is given by linear constraints, then we have a multi-extremal problem of minimizing a concave quadratic function on a set defined by linear constraints. In this case, all local minima of the problem (12.16)–(12.18) are located at the vertices of the polyhedron (12.17), (12.18). Some algorithms for minimizing a concave function on a linear polyhedron based on directed enumeration of vertices and cutting off unpromising areas of the polyhedron are considered in [8, 126, 188].

Note that finding local minima of the problem (12.16)-(12.18) is relatively easy.
Let us construct the sequence \(\left\{ {{y^s}} \right\}\), ending at some local minimum of the problem. The point \({y^0} \in D\) is arbitrary. Let the point \({y^s}\) be already constructed. We linearize the objective function (12.16) at the point \({y^s}\) and find \({y^{s + 1}}\) as a solution to the following linear programming problem:
\begin{equation}\label{eqn:12.19}
    z - L\left<{y^s},y\right> \to \mathop {\min }\limits_{y,z}\tag{12.19} 
\end{equation}
under restrictions
\begin{equation}\label{eqn:12.20}
    L\left<{b^i},y\right> - z \le \frac{L}{2}{\left\| {{b^i}} \right\|^2} - {a_i}, 
		\;\;\;i = 1,2,...,k;\tag{12.20}
\end{equation}
\begin{equation}\label{eqn:12.21}		
    y \in D\tag{12.21}
\end{equation}
Since only the objective function changes in the sequence of problems (12.19) - (12.21), the optimal solution of one problem is admissible for the next one. Here the choice of point \({y^0} \in D\) is important.

Having descended to the local minimum of the problem (12.16)—(12.18), it is necessary to find a point with smaller than the one already achieved by the value of the objective function (12.6) and then descend again to a new local extremum, etc.

In essence, we have described the idea of a combined algorithm for solving problem (12.16)—(12.18).

Unfortunately, in other cases it is not possible to reduce auxiliary problems (12.8), (12.9) or (12.14), (12.15) to a known form by simple transformations. However, we will make some transformations of these tasks.

Let $f(x)$ and $h(x)$ in (12.10), (12.11) admit lower approximations of the form (12.4), (12.13) with constants \({L_f}\) and \({L_h}\). Then the auxiliary problem (12.14), (12.15) has the form.

\[\mathop {\max }\limits_{1 \le i \le k} \left[ {{a_i} - \frac{1}{2}{L_f}{{\left\| {x - {b^i}} \right\|}^2}} \right] \to \mathop {\min }\limits_x \]
under restrictions
\[\mathop {\max }\limits_{\left\{ {i|1 \le i \le k,h({x^i}) > 0} \right\}} \left[ {{c_i} - \frac{1}{2}{L_h}{{\left\| {x - {d^i}} \right\|}^2}} \right] \le 0,\;\;\;
x \in D,\]      
where
\[{a_i} = f({x^i}) + \frac{1}{{2{L_f}}}{\left\| {\nabla f({x^i})} \right\|^2},\;\;\;{b^i} = {x^i} + \frac{1}{{{L_f}}}\nabla f({x^i}),\]
\[{c_i} = h({x^i}) + \frac{1}{{2{L_h}}}{\left\| {\nabla {f_{j*({x^i})}}({x^i})} \right\|^2},\;\;\;{d^i} = {x^i} + \frac{1}{{{L_h}}}\nabla {f_{j^*({x^i})}}({x^i}),\]
\[h({x^i}) = \mathop {\max }\limits_{1 \le j \le m} {f_j}({x^i}) = {f_{j^*({x^i})}}({x^i}).\]

This task is equivalent to the following:
\begin{equation}\label{eqn:12:22}
    W(x,t,z) \equiv z \to \mathop {\min }\limits_{x,t,z}\tag{12.22}
\end{equation}
under restrictions
\begin{equation}\label{eqn:12:23}
{L_f}\left<{b^i},x\right> - \frac{1}{2}{L_f}t - z \le \frac{1}{2}{L_f}{\left\| {{b^i}} \right\|^2} - {a_i},\;\;\;         i = 1,2,...,k;\tag{12.23}
\end{equation}
\begin{equation}\label{eqn:12:24}                          
{L_h}\left<{d^j},x\right> - \frac{1}{2}{L_h}t \le \frac{1}{2}{L_h}{\left\| {{d^j}} \right\|^2} - {c_j},\;\;\;
 j \in \left\{ {i|1 \le i \le k,h({x^i}) > 0} \right\};\tag{12.24}
\end{equation}
\begin{equation}\label{eqn:12:25}           
x \in D,\;\;\;{\left\| x \right\|^2} \ge t > 0.\tag{12.25}
\end{equation}                                      
Here the variables are $x, t, z$ and there is only one non-linear concave constraint ${\left\| x \right\|^2} \ge t > 0.$

Let now the minimized multiextremal function $f(x)$ in problem (12.6), (12.7) satisfy the Lipschitz condition with constant $L$ in $D$ and the lower approximations \(\varphi ({x^i},x)\) have form (12.5). In this case, auxiliary problems (12.8), (12.9) have the form
\[\mathop {\max }\limits_{1 \le i \le k} ({f_i} - L\left\| {x - {x^i}} \right\|) \to \mathop {\min }\limits_x \]
subject to
\[x\in D,\]
where ${f_i} = f({x^i}).$
This task is equivalent to the following:
\[V(x,z) \equiv z \to \mathop {\min }\limits_{x,z} \]
under restrictions
\[{f_i} - L\left\| {x - {x^i}} \right\| \le z,\;\;\; i = 1,2,...,k;\]
\[z \le \mathop {\min }\limits_{1 \le i \le k} {f_i},\;\;\;x \in D.\]

In turn, this task is equivalent to the following:
\[V(x,z) \equiv z \to \mathop {\min }\limits_{x,z} \]
under restrictions
\[2{L^2}\left<{x^i},x\right> - 2{f_i}z + {z^2} - {L^2}{\left\| x \right\|^2} \le c_i^2,\;\;\; i = 1,2,...,k;\]
\[z \le \mathop {\min }\limits_{1 \le i \le k} {f_i},\;\;\;x \in D,\]
where ${c_i} = {L^2}{\left\| {{x^i}} \right\|^2} - f_i^2.$

Introducing an additional variable $t$, we arrive at a new equivalent problem:
\begin{equation}\label{eqn:12.26}
		U(x,z,t) \equiv z \to \mathop {\min }\limits_{x,t,z}\tag{12.26}
\end{equation}
under restrictions
\begin{equation}\label{eqn:12.27}
   2{L^2}\left<{x^i},x\right> - 2{f_i}z + {z^2} - Lt \le {c_i},\;\;\;i = 1,2,... ,k;\tag{12.27}
\end{equation}
\begin{equation}\label{eqn:12.28}
    z \le \mathop {\min }\limits_{1 \le i \le k} {f_i},\;\;\;x \in D.\tag{12.28}
\end{equation}
\begin{equation}\label{eqn:12.29}
    {\left\| x \right\|^2} \ge t \ge 0.\tag{12.29}
\end{equation}
Here the variables are $x, z, t$ and there is only one non-linear concave constraint \({\left\| x \right\|^2} \ge t \ge 0;\)
Consider the transformed auxiliary problems (12.22)-(12.25), (12.26)-(12.29). Let us denote the optimal value of the objective function in them by \({z^*}\). Linearizing at $y$ the non-convex constraint \(t - {\left\| x \right\|^2} \le 0\), we get \(t - 2\left<y,x\right> + {\left\| y \right\|^2} \le 0\). Substituting the latter for the constraint \(t - {\left\| x \right\|^2} \le 0\) into problems (12.22)-(12.25), (12.26)-(12.29), we obtain problems of linear and convex programming.
Denote the optimal solution of such problems by \({z^*}(y)\). Since \(t - {\left\| x \right\|^2} \le t - 2\left<y,x\right> + {\left\| y \right\|^2}\), then \({z^*} \le {z^*}(y).\)
It occurs:
\[\left\{ {(t,x) \in {E_{n + 1}}|x \in D,\;t - {{\left\| x \right\|}^2} \le 0} \right\} = \]
\[ = \bigcup\limits_{y \in D} {\left\{ {(t,x)|x \in D,\;t - 2\left<y,x\right> + {{\left\| y \right\|}^2} \le 0} \right\}} ,\]
so \({z^*} = \mathop {\inf }\limits_{y \in D} {z^*}(y).\) Thus, taking a grid of points \(\left\{ {{y^s}} \right\} \subset D\) and solving a set of transformed auxiliary problems with the non-convex constraint \(t - {\left\| x \right\|^2} \le 0\) linearized at the points \({y^s}\), we obtain an estimate for the optimal value \({z^*}\). For different \({y^s}\), problems (12.22)-(12.25) with a linearized constraint \(t - {\left\| x \right\|^2} \le 0\) differ only by one linear constraint, so it is reasonable to solve them by the dual simplex method.

\textbf{3. The smoothing method.} 
\label{Sec.12.3}
Smoothing or averaging of functions has long been used in mathematics. In particular, in this book it is widely used to construct local minimization methods for Lipschitz functions. For the purposes of global optimization, it has been used since [138]. The point is that smoothing (averaging) of a function eliminates small local extrema, while preserving the overall picture of the relief of the function. By minimizing the smoothed function, we hope to skip the local extrema of the original function and get into the region of the global minimum.

Let it be necessary to solve the problem:
\[f(x) \to \mathop {\min }\limits_{x \in {E_n}}, \]
where $f(x)$ is a continuous function such that \(f(x) \to + \infty \) as \(\left\| x\right\| \to + \infty. \)
Denote a cube with edge $h$ and center $x$ as
\[{K_h}(x) = \left\{ {y \in {E_n}|\mathop {\max }\limits_{1 \le i \le n} |{y_i} - {x_i}| \le h/2} \right\},\]
\({S_h}(x)\) is its surface, $N (y)$ is the outward normal to \({S_h}(x)\) at point
\(y \in {S_h}(x).\)

Consider the smoothed (averaged) functions:
\[{f_h}(x) = \frac{1}{{{h^n}}}\int\limits_{{K_h}(x)} {f(y)d{y_1}d{y_2}.. .d{y_n}} \]

The gradient \({f_h}(x)\) when f(x) is continuous is calculated by the surface integral
\[\nabla {f_h}(x) = \frac{1}{{{h^n}}}\int\limits_{{S_h}(x)} {f(y)N(y)} dS.\]

The global optimization method based on smoothing is as follows.

Starting from an arbitrary point \({x^0} \in {E_n}\), we minimize the smoothed function \({f_{{h_1}}}(x)\) with a sufficiently large smoothing parameter \({h_1}\) till the accuracy \[\left\| {\nabla {f_{{h_1}}}({x^1})} \right\| \le {\varepsilon _i}.\] Then, starting from \({x^1}\), we minimize the function \({f_{{h_2}}}(x)\) with a smaller smoothing parameter \({h_2} < {h_1}\) till the accurateness \( \left\| {\nabla {f_{{h_2}}}({x^2})} \right\| \le {\varepsilon _2}\), etc.

In this method
\[\mathop {\lim }\limits_{k \to \infty } {h_k} = \mathop {\lim }\limits_{k \to \infty } {\varepsilon _k} = 0\]
When the function $f(x)$ is continuously differentiable, the mapping \((x,h) \to \nabla {f_h}(x)\) is continuous. Therefore the sequence \(\{ {x^k}\} \) cannot converge by its limit points to non-stationary points of $f(x)$. Using Lemma 2.5, a similar property can also be proved for locally Lipschitz functions $f(x)$.

However, what are the global properties of a method? In this connection, the following observation is interesting [81].

Let $f(x)$ depends on a scalar variable $x$. Consider smoothed functions (\(h \ge 0\)):
\[{f_h}(x) = \frac{1}{h}\int\limits_{x - h/2}^{x + h/2} {f(y)dy,} \;\;\;{f_0 }(x) = f(x).\]
Their gradients in this case have the form:
\[\nabla {f_h}(x) = \left[ {f\left( {x + \frac{h}{2}} \right) - f\left( {x - \frac{h}{2} } \right)} \right]{h^{ - 1}},\;\;\;h > 0.\]

Let us construct a picture of stationary points of smoothed functions \({f_h}(x)\) in the plane \({E_x} \times {E_h} = \left\{ {(x,h)|x \in {E_x},h \in {E_h}} \right\}\), i.e. we construct the set \(T = \left\{ {\left( {x,h} \right) \in {E_x} \times {E_h}|\nabla {f_h}(x) = 0} \right\}\).

Note that if \(f({x^1}) = f({x^2})\) \(({x^1} \ne {x^2})\), then the point \((( {x^1} + {x^2})/2,\) \(|{x^2} - {x^1}|) \in T\). Thus, to construct $T$, it suffices to trace all segments parallel to the \({E_x}\)-axis, the ends of which lie on the graph of $f(x)$. By continuously moving these segments up or down and keeping track of their length, it is easy to build the set $T$.

The picture of the set $T$ is most transparent when all extrema of  $f(x)$ are different. Then $T$ consists of continuous non-intersecting lines or bands associated with extrema on the \({E_x}\)-axis. Moreover, any line starting at a local extremum \(({x^*},0)\) (on the axis \({E_x}\)) is bounded and closes to another local extremum \(({x^{**} },0)\)(on the \({E_x}\) axis). Lines do not come to stationary non-extremal points. And only the line starting at the global minimum \(\left( {{x_{\min }},0} \right)\) goes to infinity in the plane \({E_x} \times {E_h}\).

If not all extrema of $f(x)$ are different, then the picture of $T$ is somewhat more complicated. However, as before, $T$ splits into connected limited components, some connecting local extrema, and an unbounded connected component connecting all global minima of $f(x)$ and infinity, as well as main maxima between adjacent global minima.

Minimizing \({f_h}(x)\) with a sufficiently large smoothing parameter $h$, we obviously get to a stationary line from $T$ leading to a global minimum. However, it is not easy to follow this line to the global extremum, since it can behave quite bizarrely, in particular, it can have kinks. The points on this line can be local minima, local maxima, or simply stationary in $x$ for the corresponding smoothed functions \({f_h}(x)\).

\newpage
\begin{flushright}
CHAPTER 4
\label{Ch.4}

\textbf{NONMONOTONE METHODS WITH AVERAGING OF\\
DESCENT DIRECTIONS}

\underline{\hspace{15cm}}

\end{flushright}
\bigskip\bigskip\bigskip\bigskip\bigskip\bigskip

Recently, much attention has been attracted by relaxation methods of non-smooth optimization, which, on the one hand, follow from the conjugate gradient methods of smooth optimization, and, on the other hand, generalize the known methods for minimizing the maximum function.
Their main feature is that the direction of descent is determined by the averaged subgradients calculated in some neighborhood of the current point. Such algorithms are also called \(\varepsilon \)-subgradient methods. In these methods, the main and very complex issue is the construction of the direction of a significant decrease in the function.

In this chapter, we study nonmonotonic methods for minimizing nondifferentiable functions by averaging finite differences or generalized gradients. Nonmonotonic methods are generally no more complicated than the iterative procedures discussed in previous chapters. Their construction does not require additional memory and complex parameter adjustments. It should be noted that the optimal fast convergent ravine step method [74] was constructed in the class of nonmonotone methods. In Section 14 it will be shown that in its structure it is close to the methods of averaged gradients. In what follows, averaging procedures of the form
\[{z^{k + 1}} = (1 - {a_k}){z^k} + {a_k}H({x^k},{\alpha _k})\]
will be widely used in Chapter 5 when constructing numerical methods for minimizing Lipschitz functions with constraints, and also in Chapter 7 for solving stochastic extremal problems.

\section*{$\S$ 13. Methods with averaging descent directions }
\label{Sec.13}
\setcounter{section}{13}
\setcounter{definition}{0}
\setcounter{equation}{0}
\setcounter{theorem}{0}
\setcounter{lemma}{0}
\setcounter{remark}{0}
\setcounter{corollary}{0}
\numberwithin{equation}{section}
In Chapter 2, finite-difference methods for minimizing Lipschitz functions were studied. They have the disadvantages of conventional one-step procedures. To make the process of searching for an extremum more regular, it is advisable to average the descent direction vectors that are calculated at previous iterations. Direction averaging is advantageous in problems of minimizing ravine-type functions, since the average direction indicates movement along the ravine.

Let it be required to minimize a locally Lipschitz function $f(x)$. Define a sequence of points
\begin{equation}
\label{eqn:13.1}
{x^{k + 1}} = {x^k} - {\rho _k}\sum\limits_{t = {r_k}}^k {H({x^t},{\alpha _i})},
\;\;\;x^0\in E_n, \;\;\; k=0,1,\ldots. 
\end{equation}
Here \(\left\{ {{r_k}}\leq k \right\}\)  is a sequence of moments (at which the summation resumes) that determines the number of averaged vectors. The vectors \(H({x^k},{\alpha _k})\)  are defined by one of the following formulas:
\[
H({x^k},{\alpha _k}) = \frac{1}{{2{\alpha _k}}}\sum\limits_{i = 1}^n {\left[ {f(\tilde x_1^k,...,x_i^k + {\alpha _k},...,\tilde x_n^k) - } \right.} \]
\begin{equation}
\label{eqn:13.2}
\left. { - f(\tilde x_1^k,...,x_i^k - {\alpha _k},...,\tilde x_n^k)} \right]{e_i}, 
\end{equation} 
\begin{equation} \label{eqn:13.3}
    H({x^k},{\alpha _k}) = \sum\limits_{i = 1}^n {\frac{{f({{\tilde x}^k} + {\Delta _k}{e_i}) - f({{\tilde x}^k})}}{{{\Delta _k}}}} {e_i},
\end{equation} 
\begin{equation} \label{eqn:13.4}                             
    H({x^k},{\alpha _k}) = \sum\limits_{i = 1}^p {\frac{{f({{\tilde x}^k} + {\Delta _k}\mu _i^k) - f({{\tilde x}^k})}}{{{\Delta _k}}}} \mu _i^k,
\end{equation} 
similar to formulas (4.6), (4.7), (6.1), respectively.

In the generalized gradient descent method
\begin{equation} \label{eqn:13.5}
            H({x^k},{\alpha _k}) = g({\tilde x^k}).
\end{equation}
            
In practical calculations, averaging the vectors \(H({x^k},{\alpha _k})\) over  points contained in the \({\varepsilon _k}\)-neighborhood of the point \({x^{{r_k }}}\), i.e. the next restore occurs when
\[\left\| {{x^k} - {x^{{r_k}}}} \right\| > {\varepsilon _k}.\]
In order not to take into account too old information, ${\varepsilon _k}$ tends to zero. In what follows, we will assume that the number of averaged vectors is limited by a certain value M, i.e., $k - {r_k} \le M.$
\begin{theorem}\label{th:13.1}
Let the the following conditions hold true:
\[\sum\limits_{k = 0}^\infty  {{\rho _k}}  = \infty,\;\;\; 
\sum\limits_{k = 0}^\infty  {\rho _k^2}  < \infty,\;\;\; 
\frac{{{\rho _k}}}{{{\alpha _k}}} \to 0,\;\;\; 
\frac{{{\Delta _k}}}{{{\alpha _k}}} \to 0\]
\[\frac{{\left| {{\alpha _k} - {\alpha _{k + 1}}} \right|}}{{{\rho _k}}} \to 0,
\;\;\;{\alpha _k} \to 0,\;\;\;
\mathop {\lim }\limits_{k \to \infty } \frac{{{\alpha _{k - 1}}}}{{{\alpha _k}}} = 1,\]
\[
\underset{k \to \infty }{\overline {\lim }} \frac{{{\rho _k}}}{{{\rho _{k + 1}}}} \le C < \infty.
\]
Then the limit points of the sequence \(\left\{ {{x_k}} \right\}\) with probability 1 belong to the set 
\[{X^*} = \left\{ {{x^*}|0 \in \partial f({x^*})} \right\},\]
and if $\{f(x)|x\in X^*\}$ does not contain intervals, the sequence \(\left\{ {f({x^k})} \right\}\) converges to some limit with probability 1.
\end{theorem}

{\it P r o o f}. The proof largely coincides with the proof of Theorem 4.1. From the mean value theorem, we have
    \[f(x^{k + 1}),{\alpha _k}) = f({x^k},{\alpha _k}) + \left<\nabla f({x^k} + \tau ( {x^{k + 1}} - {x^k}),{\alpha _k}),{x^{k + 1}} - {x^k}\right> - \]
    \[ - \frac{1}{{k - {r_k} + 1}}\left<\sum\limits_{t = {r_k}}^k \nabla f({x^t},{\alpha _t}), {x^{k + 1}} - {x^k}\right> + \]
    \[ + \frac{1}{{k - {r_k} + 1}}\left<\sum\limits_{t = {r_k}}^k \nabla f({x^t},{\alpha _t}), {x^{k + 1}} - {x^k}\right>,\;\;\;0 \le \tau \le 1.\]

Let us estimate the quantities
\[\left\| {\nabla f({x^k} + \tau ({x^{k + 1}} - {x^k}),{\alpha _k}) - \nabla f({x^i},{\alpha _i})} \right\|,\;\;\;
i = {r_k},...,k.\]
Notice, that
\[\left\| {\nabla f(y,{\alpha _k}) - \nabla f(y,{\alpha _{{r_k}}})} \right\| \le C\left[ {\frac{{\left| {{\alpha _k} - {\alpha _{{r_k}}}} \right|}}{{{\alpha _{{r_k}}}}} + \frac{{\left| {{\alpha _k} - {\alpha _{{r_k}}}} \right|}}{{{\alpha _k}}}} \right],\]
\[\left\| {\nabla f(y,{\alpha _k}) - \nabla f(y,{\alpha _{k - 1}})} \right\| \le C\left[ {\frac{{\left| {{\alpha _k} - {\alpha _{k - 1}}} \right|}}{{{\alpha _{k - 1}}}} + \frac{{\left| {{\alpha _k} - {\alpha _{k - 1}}} \right|}}{{{\alpha _k}}}} \right],\]
\[\left| {f(y,{\alpha _{k + 1}}) - f(y,{\alpha _k})} \right| \le C\left| {{\alpha _{k + 1}} - {\alpha _k}} \right|.\]
That's why
\[\left\| {\nabla f({x^k} + \tau ({x^{k + 1}} - {x^k}),{\alpha _k}) - \nabla f({x^i},{\alpha _i})} \right\| = \]
\[ = \left\| {\nabla f({x^k} + \tau ({x^{k + 1}} - {x^k}),{\alpha _k}) - \nabla f({x^i},{\alpha _k}) + \nabla f({x^i},{\alpha _k}) - } \right.\]
\[\left. { - \nabla f({x^i},{\alpha _i})} \right\| \le \frac{C}{{{\alpha _k}}}\sum\limits_{t = i}^k {{\rho _t} + C\left[ {\frac{{\left| {{\alpha _k} - {\alpha _i}} \right|}}{{{\alpha _i}}} + \frac{{\left| {{\alpha _k} - {\alpha _i}} \right|}}{{{\alpha _k}}}} \right]} ,\]
\[\left\| {\nabla f({x^k} + \tau ({x^{k + 1}} - {x^k}),{\alpha _k}) - \nabla f({x^k},{\alpha _k})} \right\| \le \frac{{C{\rho _k}}}{{{\alpha _k}}}.\]
Hence, it follows that
\[f({x^{k + 1}},{\alpha _{k + 1}}) \le f({x^k},{\alpha _k}) + C\left| {{\alpha _k} - {\alpha _{k + 1}}} \right| + \]
\[ + C\frac{{{\rho _k}}}{{{\alpha _k}}}\sum\limits_{t = {r_k}}^k {{\rho _t} + C} {\rho _k}\sum\limits_{t = {r_k}}^{k - 1} {\left[ {\frac{{\left| {{\alpha _t} - {\alpha _k}} \right|}}{{{\alpha _t}}} + \frac{{\left| {{\alpha _t} - {\alpha _k}} \right|}}{{{\alpha _k}}}} \right]}  - \]
\[ - \frac{{{\rho _k}}}{{k - {r_k} + 1}}\left< {\sum\limits_{t = {r_k}}^k {\nabla f({x^t},{\alpha _t})} ,\sum\limits_{t = {r_k}}^k {(H({x^t},{\alpha _t})}  - \nabla f({x^t},{\alpha _t}) + } \right.\]
\[\left. { + \nabla f({x^t},{\alpha _t}))} \right> \le f({x^k},{\alpha _k}) + C\left| {{\alpha _k} - {\alpha _{k + 1}}} \right| + C\frac{{{\rho _k}}}{{{\alpha _k}}}\sum\limits_{t = {r_k}}^k {{\rho _t}}  + \]
\begin{eqnarray}
&&
+C\rho_k\sum_{t=r_k}^{k-1}
\left[\frac{|\alpha_t-\alpha_k}{\alpha_t}+\frac{|\alpha_t-\alpha_k}{\alpha_k}  \right]  
-\frac{\rho_k}{k-r_k+1}\left\|\sum_{t=r_k}^{k}\nabla f(x^t,\alpha_t)   \right\|^2  \nonumber\\
&&
+\frac{\rho_k}{k-r_k+1}
\left<   
\sum_{t=r_k}^{k}\nabla f(x^t,\alpha_t),
\sum_{t=r_k}^{k}\left(\nabla f(x^t,\alpha_t)-H(x^t,\alpha_t) \right)
\right>.\label{13.6}
\end{eqnarray}
To prove convergence, it is enough to check the fulfillment of
of conditions (4.12), (4.13).

Let $x^s\rightarrow x^\prime \notin X^*$ ($s\in S$) and the $2\bar{\delta}$-neighborhoods of the points $x_s$ ($s\ge \bar{s})$ do not intersect with $X^*$. Suppose that condition (4.12) is not satisfied, i.e.
all points $x^k$ ($k\ge s\ge \bar{s}$) are contained in the $\delta$-neighborhood of  point $x^s$ , $\delta \le\bar{\delta}/2$. 
Since $\rho_k\rightarrow 0$, for sufficiently large $s$ points 
$x^{s-M}, x^{s-M+1}, \ldots, x^{s-1}$ also belong to the $\delta$-neighborhood of $x_s$. Therefore, it follows from Lemma 2.5 that for sufficiently large $s$ and sufficiently small $\delta$
the inequality holds true:
\[
\left\|\sum_{t=r_k}^{k}\nabla f(x^t,\alpha_t) \right\|\ge\sigma>0.
\]
From the conditions of the theorem, it follows
\[
\frac{\rho_{r_k}}{\alpha_k}=\frac{\rho_{r_k}}{\rho_{r_k+1}}\times
\ldots\times\frac{\rho_{k}}{\alpha_k}\rightarrow 0;
\]
similarly,
\[
\frac{\rho_{r_k+1}}{\alpha_k}\rightarrow 0,\ldots,
\frac{\rho_{k-1}}{\alpha_k}\rightarrow 0.
\]
and
\[
\frac{|\alpha_{i}-\alpha_k|}{\alpha_i}\rightarrow 0,\ldots,
\frac{|\alpha_{i}-\alpha_k|}{\alpha_k}\rightarrow 0,\;\;\;
i=r_k,\ldots,k.
\]
Hence, for sufficiently large $k$
\[
C|\alpha_k-\alpha_{k+1}|+\frac{C}{\alpha_k}\sum_{t=r_k}^{k}\rho_t+
C\sum_{t=r_k}^{k-1}
\left[
\frac{|\alpha_{t}-\alpha_k|}{\alpha_t}+
\frac{|\alpha_{t}-\alpha_k|}{\alpha_k}
\right]
\le \frac{\sigma}{2M}.
\]
Then by virtue of convergence with probability 1 of the series
\[
\sum_{k=0}^{\infty}
\left<   
\sum_{t=r_k}^{k}\nabla f(x^t,\alpha_t),
\sum_{t=r_k}^{k}\left(\nabla f(x^t,\alpha_t)-H(x^t,\alpha_t) \right)
\right>
\]
and from inequality (13.6) at large $s$ it follows
\begin{equation}
\label{eqn:13.7}
f(x^k,\alpha_k)\le f(x^s,\alpha_k)-\frac{\sigma}{2M}\sum_{t=s}^{k-1}\rho_t,
\end{equation}
which at $k\rightarrow\infty$ contradicts the boundedness of 
$\{f(x_k, a_k)\}.$ Hence,
\[
k(s)=\min\{r|\,\|x^r-x^s\|>\delta,\;r>s\}<\infty.
\]
By definition,
\[
x^k(s)\notin U_\delta(x^s)=\{x|\,\|x-x^s\|\leq\delta\},
\]
but since $\rho_k\rightarrow 0$, for sufficiently large $s$ holds
\[
x^{k(s)}\in U_{2\delta}(x^s),
\]
i.e., inequality (13.7) remains valid for $k=k(s)$. From the relation
\[
x^{k(s)}=x^s-\sum_{t=s}^{k(s)-1}\rho_t\left(\sum_{p=r_t}^{t}H(x^p,\alpha_p)  \right)
\]
it follows that
\[
\sum_{t=s}^{k(s)-1}\rho_t>\frac{\delta}{CM}.
\]
Therefore, from inequality (13.7) we obtain
\[
f(x^{k(s)},\alpha_{k(s)})\le f(x^s,\alpha_s)-\frac{\sigma\delta}{2CM^2},
\]
from where (4.13) directly follows. The theorem is proved.

As the direction of motion, we can take the vector $d_k$, which
belongs to the convex hull of vectors 
$H (x^{r_k},\alpha_{r_k}),\ldots,H({x}^k,\alpha_k)$, i.e.
\begin{equation}
\label{eqn:13.8}
x^{k+1}=x^k-\rho_k d^k,
\end{equation}
\[
d^k=\sum_{t=r_k}^k\lambda_t^kH(x^t,\alpha_t),\;\;\;
\sum_{t=r_k}^k\lambda_t^k=1,\;\;\; \lambda_t^k\ge 0.
\]
\begin{theorem}
\label{th:13.2}
Let the conditions of Theorem 13.1 be satisfied. Then
the limit points of the sequence $\{x^k\}$ (13.8) with probability 1 belong to the set 
\[
X^*=\{x^*|\,0\in\partial f(x^*)\},
\]
the sequence $\{f(x^k)\}$ converges with probability 1.
\end{theorem}
{\it P r o o f}. The proof is practically similar to the previous one, since
vector belonging to the convex hull of vectors $H (x^{r_k},\alpha_{r_k}),\ldots,H({x}^k,\alpha_k)$ can be represented it in the form (13.3). Therefore,
\begin{eqnarray}
f\left(x^{k+1},\alpha_{k+1}\right)
&\le&
f\left(x^{k},\alpha_{k}\right)
+\sum_{t=r_k}^k\lambda_t^k
\left<\nabla f\left(x^k+\tau(x^{k+1}-x^k),\alpha_k\right)\right.\nonumber\\
&&\left.-\nabla f(x^t,\alpha_t),x^{k+1}-x^k\right>
+\sum_{t=r_k}^k\lambda_t^k\left<\nabla f(x^t,\alpha_t),x^{k+1}-x^k  \right>.\nonumber
\end{eqnarray}

From here it follows
\begin{align}
  & f({{x}^{k+1}},{{a}_{k+1}})\le f({{x}^{k}},{{a}_{k}})+C\left| {{a}_{k}}-{{a}_{k+1}} \right|+ \nonumber\\ 
 & +C\frac{{{\rho }_{k}}}{{{\alpha }_{k}}}\sum\limits_{t={{r}_{k}}}^{k}{{{\rho }_{t}}+C{{\rho }_{k}}\sum\limits_{t={{r}_{k}}}^{k}{[\frac{\left| {{a}_{t}}-{{a}_{k}} \right|}{{{a}_{t}}}}}+\frac{\left| {{a}_{t}}-{{a}_{k}} \right|}{{{a}_{k}}}-\nonumber \\ 
 & -{{\rho }_{k}}{{\left\| \sum\limits_{t={{r}_{k}}}^{k}{\lambda _{t}^{k}\nabla f({{x}^{t}},{{a}_{t}}} )\right\|}^{2}}+{{\rho }_{k}}\left<\sum\limits_{t={{r}_{k}}}^{k}{\lambda _{t}^{k}}\nabla f({{x}^{t}},{{a}_{t}}),\right.\nonumber
\end{align}
\begin{equation}\label{eqn:13.9}
\left.\sum\limits_{t={r}_{k}}^{k}\lambda _{t}^{k}(\nabla f({{x}^{t}},{{a}_{t}})-H({{x}^{t}},{{a}_{t}}))\right>. 
\end{equation}

Thus, we obtained the inequality similar to (13.6).

In relaxation methods of nonsmooth optimization of the type (13.8), as $d^k$ is usually taken the nearest to zero vector from the convex hull of accumulated vectors, i.e., a special quadratic programming problem is solved additionally. If monotonicity of the process is not required, then, as practical calculations show, methods of the type (13.8) have no advantage in convergence rate compared to methods (13.1) with the simplest averaging and without memorization of information.

Note that the conditions of Theorem 13.1 are also satisfied for the step multipliers
\[\rho _{k}^{\prime}={{\rho }_{k}}{{n}_{k}},\;\;\;\rho _{k}^{\prime\prime}={{\rho }_{k}}/{{n}_{k}},\;\;\;1\le {{n}_{k}}\le N<\propto .\]
Therefore, a rough one-dimensional minimization can be performed in the chosen direction of motion. Usually in numerical calculations the vectors $H (x^k, a^k)$, as well as the descent vector itself in (13.1), are chosen to be normalized.

Averaging of descent directions is also advantageous in conditional minimization methods. Let
\[X={{X}_{1}}\cap {{X}_{2}},\;\;\;{{X}_{1}}=\left\{ x\left| {{f}_{i}}(x)\le 0, \right. \right.\;\;\;i=1,...,\left. m \right\},\]
$X_2$ is a convex bounded and closed set, and $X_2$ has a non-empty intersection with the set of interior points of $X_1$. The method of minimizing $f_0 (x)$ at $x \in X$ with subgradient averaging is defined by the sequence of points:
\begin{equation}\label{eqn:13.10}
{{x}^{k+1}}=\pi_{{X}_{2}}({{x}^{k}}-{{\rho }_{k}}{\bar{x}^{k}}), 
\end{equation}
\[{\bar{x}^{k}}=\sum\limits_{t={{r}_{k}}}^{k}{{{g}^{{{\nu}_{t}}}}({{x}^{t}});}\]
\[
{{\nu}_{t}}=\left\{ 
\begin{array}{lc}
0,&\underset{1\le i\le m}{\mathop{\max }}\,{{f}_{i}}({{x}^{t}})\le 0;\\
j,&\underset{1\le i\le m}{\mathop{\max }}\,{{f}_{i}}({{x}^{t}})={{f}_{j}}({{x}^{t}})>0.
\end{array} 
\right.
\]
Here ${{f}_{i}}(x)$ are convex functions, ${{g}^{i}}(x)\in \partial {{f}_{i}}(x)$
$(i=0,1,.... ,m),$ ${\pi}_{{X}_{2}}(x)$ is the projection operator of a point $x$ on the set $X_2$, at points ${x}^{{r}_{k}}$ the summation resumes.
\begin{theorem}
\label{th:13.3}
Let
\[{{\rho }_{k}}\to 0,\;\;\;\sum\limits_{k=0}^{\propto }{{{\rho }_{k}}=\infty .}\]
Then the limit points of the sequence $\left\{ {{x}^{k}} \right\}$ belong to the set
\[{{X}^{*}}=\left\{ {{x}^{*}}\in\operatorname{arg}\underset{x\in X}{\mathop{\min }}\,{{f}_{0}}(x) \right\}.\]
\end{theorem}
{\it P r o o f.} Let us denote 
\[{{Q}_{\tau }}=\left\{ x\left| x\in {{X}_{1}},{{f}_{0}}(x)\le \tau \right. \right\}.\]
Let $X_{\delta }^{*}\subset {{Q}_{\tau }}$ be the set of interior points of the set ${{Q}_{\tau }}$ that are at least $\delta$-away from the boundary of ${{Q}_{\tau }}$. When $\delta $ is small enough, the set $X_{\delta }^{*}$ is nonempty. We first show that the limit points of the sequence $\left\{ {{x}^{k}} \right\}$ belong to the set ${{Q}_{\tau }}$; hence, by taking $\tau $ to ${{f}_{0}}({{x}^{*}})$, we obtain the statement of the theorem.

Let there exist a sequence
 \[{{x}^{s}}\to {{x}^{\prime}}\not\subset {{Q}_{\tau }},\;\;\;  s\in S.\]
There are such $\bar{\mathop{s}}\,$ and $\bar{\varepsilon}\,$ that, starting from some number $s\ge \bar{\mathop{s}},\, $ $2\bar{\varepsilon}$-neighborhoods of the points ${{x}^{s}}$ do not intersect with ${{Q}_{\tau }}$. Let us show that for any
$\varepsilon\le \bar{\varepsilon }/2$ and large enough $s$
\[k(s)=\min \left\{ r\left| \left\| {{x}^{r}}-{{x}^{s}} \right\|>\varepsilon  \right. \right\}<\infty .\]
Let's assume the opposite:
\[{{x}^{k}}\in {{U}_{\varepsilon }}({{x}^{s}}),\;\;\;s\ge \bar{s}.\,\]
Since ${{{\rho }_{k}}}\to 0$, then for sufficiently large $s$ the points ${{x}^{s-M}},...,{{x}^{s-1}}$ also belong to ${{U}_{\varepsilon}}({{x}^{s}})$. Let 
${{x}^{k}}\in {{X}_{1}};$ then
\[0>-\sigma \ge \underset{y\in X_{\delta }^{*}}{\mathop{\max }}\,({{f}_{0}}(y)-{{f}_{0}}({{x}^{k}})).\]
If \[{{x}^{k}}\notin {{X}_{1}},\] then
\[0>-\sigma \ge \underset{y\in X_{\delta }^{*}}{\mathop{\max }}\,({{f}_{j}}(y)-{{f}_{j}}({{x}^{k}})),\;\;\;j=1,...,m.\]
Therefore, it follows from these inequalities for sufficiently small $\varepsilon$ that
\[\underset{y\in X_{\delta }^{*}}{\mathop{\max }}\,\sum\limits_{t={{r}_{k}}}^{k}{({{g}^{{{\nu}_{t}}}}({{x}^{t}}),y-{{x}^{k}})\le -\frac{\sigma }{2}.}\]
Let us define
\[W(x)=\underset{y\in X_{\delta }^{*}}{\mathop{\min }}\,{{\left\| {{x}^{k}}-y \right\|}^{2}}={{\left\| {{x}^{k}}-y({{x}^{k}}) \right\|}^{2}},\;\;\;y({{x}^{k}})\in X_{\delta }^{*}.\]
Then
\begin{align}
  & W({{x}^{k+1}})\le {{\left\| {{x}^{k+1}}-y({{x}^{k}}) \right\|}^{2}}={{\left\| \pi_{{X}_{2}}({{x}^{k}}-{{\rho }_{k}}{\bar{x}^{k}})-y({{x}^{k}}) \right\|}^{2}}\le  
	\nonumber\\ 
 & \le {{\left\| {{x}^{k}}-{{\rho }_{k}}{\bar{x}^{k}}-y({{x}^{k}}) \right\|}^{2}}=W({{x}^{k}})+2{{\rho }_{k}}\left<{\bar{x}^k},y({{x}^{k}})-{{x}^{k}}\right>+ 
\nonumber\\ 
 & +\rho _{k}^{2}{{\left\| {\bar{x}^{k}} \right\|}^{2}}\le W({{x}^{k}})+C\rho _{k}^{2}-{{\rho }_{k}}\sigma . \nonumber 
\end{align}
For large $s$, we have $C{{\rho }_{k}}\le \sigma /2$, therefore
\[0\le W({{x}^{k}})\le W({{x}^{s}})-\frac{\sigma }{2}\sum\limits_{t=s}^{k-1}{{{\rho }_{t}}}.\]
Passing to the limit on $k\to \infty $, we obtain a contradiction with the condition of boundedness from below of $W(x)$. It is easy to show that
\[\underset{s\to \infty }{\overset{\_\_\_\_\_}{\mathop{\lim }}}\,W({{x}^{k(s)}})<\underset{s\to \infty }{\mathop{\lim }}\,W({{x}^{s}}).\]
Thus, condition (4.13) is satisfied. Note that $W(x)$ on ${{Q}_{\tau }}$ takes a continuum of values, but this is not an obstacle to the convergence proof because $W(x)$ on ${{Q}_{\tau }}$, takes values from $[0, \delta ]$. Hence, the limit points of the sequence $\left\{ {{x}^{k}} \right\}$ belong to ${{Q}_{\tau }}$, from where by pointing $\tau$ to ${{f}_{0}}({{x}^{*}})$ , we obtain the required statement of the theorem.

The finite-difference method for minimizing ${{f}_{0}}(x)$ at $x\in X$ is as follows:
\begin{equation}\label{eqn:13.11}
{{x}^{k+1}}=\pi_{{X}_{2}}({{x}^{k}}-{{\rho }_{k}}{\bar{x}^{k}}),
\end{equation} 
where
\[
{\bar{x}^{k}}=\sum\limits_{t={{r}_{k}}}^{k}{{{H}^{{{\nu}_{t}}}}({{x}^{t}},{{a}_{t}});}
\]
\[
{{\nu}_{t}}=\left\{ 
\begin{array}{lcl}
0,&\underset{1\le i\le m}{\mathop{\max }}\,{{f}_{i}}({{x}^{t}})\le 0,\\
j,&\underset{1\le i\le m}{\mathop{\max }}\,{{f}_{i}}({{x}^{t}})={{f}_{j}}({{x}^{t}})>0;
\end{array} \right.
\]
vectors ${H}^{{{{\nu}_{t}}}}({{x}^{t}},{{a}_{t}})$ are defined by formulas (13.2)-(13.4).
\begin{theorem}
\label{th:13.4}
Let the following conditions be satisfied:
\[\sum\limits_{k=0}^{\infty }{{{\rho }_{k}}=\infty ,}\;\;\;                        \sum\limits_{k=0}^{\infty }{\rho _{k}^{2}<\infty ,}\;\;\;                     
\frac{{{\Delta }_{k}}}{{{\alpha }_{k}}}\to 0,\;\;\;{{\alpha }_{k}}\to 0.\]
Then with probability 1 the limit points of the sequence $\{x^k\}$ belong to the set $X^*$.
\end{theorem}

The proof of the theorem is similar to the previous one.

\section*{$\S$ 14. Methods of averaged generalized gradients}
\label{Sec.14}
\setcounter{section}{14}
\setcounter{definition}{0}
\setcounter{equation}{0}
\setcounter{theorem}{0}
\setcounter{lemma}{0}
\setcounter{remark}{0}
\setcounter{corollary}{0}
\numberwithin{equation}{section}
In this paragraph we study a class of non-relaxation methods for minimizing generalized differentiable functions, in which the direction of motion is a convex combination of some number of generalized gradients obtained at previous iterations, and the step is adaptively adjusted. The averaging of previous gradients can be viewed as a simple technique that improves the performance of the gradient method. It is shown that special cases of the averaged gradient method are the well-known heavy ball [100] and gully step [14, 74] methods generalized to minimize (non-convex non-smooth) generalized differentiable functions under constraints. It is well known that in the smooth convex case, these methods are far superior to the gradient method. Moreover, these methods can be used for global minimization of functions under constraints.

Let us construct the averaged gradient method to solve the following problem:
\begin{equation}
\label{eqn:14.1}
f(x)\rightarrow\min_x
\end{equation}
subject to
\begin{equation}
\label{eqn:14.2}
h(x)\le 0
\end{equation}
where $f(x)$ and $h(x)$ are generalized differentiable functions, $G_f (x)$ and $G_h (x)$ are their pseudogradient mappings, $h (x)\rightarrow+\infty$ at $\|x\|\rightarrow\infty$. The extension of the method to problems with more general constraints can in principle be done in the same way as it was done for the generalized gradient method in Section 10. Let us denote by
\begin{equation}\label{eqn:14.3}
G(x)=\left\{
\begin{array}{ll}
G_f(x),&h(x)<0,\\
co\{G_f(x),G_h(x)\},&h(x)=0,\\
G_h(x),&h(x)>0;
\end{array}
\right.
\end{equation}
\[{{G}_{\delta }}(x)=co\left\{ G(y)\left| \left\| y-x \right\|\le \delta  \right. \right\}.\]
We will treat the convergence of the method of averaged gradients as a consequence of a certain stability of the generalized gradient method, to the study of which we will now proceed.

{\bf 1. Stability of the generalized gradient method.} 
\label{Sec.14.1}
Consider the following method:
\begin{equation}
\label{eqn:14.4}
{{x}^{0}}\in {{E}_{n}},\;\;\;{{x}^{k+1}}={{x}^{k}}-{{\rho }_{k}}{{P}^{k}},\;\;\;
\;\;\;{{P}^{k}}={{g}^{k}}+{{\vartriangle }^{k}},\;\;\;k=0,1,...,
\end{equation}
\begin{equation}
\label{eqn:14.5}
{{g}^{k}}\in {{G}_{{{\delta }_{k}}}}({{x}^{k}}),
\end{equation}
\begin{equation}
\label{eqn:14.6}
0\le {{\rho }_{k}}\le \rho ,\;\;\;0\le {{\delta }_{k}}\le \delta ,\;\;\;            \left\| {{\Delta }^{k}} \right\|\le \Delta ,
\end{equation}
\begin{equation}
\label{eqn:14.7}
\underset{k\to \infty }{\mathop{\lim }}\,{{\rho }_{k}}=\underset{k\to \infty }{\mathop{\lim {{\delta }_{k}}}}\,=\underset{k\to \infty }{\mathop{\lim }}\,\left\| {{\Delta }^{k}} \right\|=0,\;\;\;\sum\limits_{k=0}^{\infty }{{{\rho }_{k}}=+\infty .}
\end{equation}
This method differs from the generalized gradient method (9.2), (9.3) or (10.7), (10.8) in that, without losing convergence, we can take the current gradient $g^k$ not at $x^k$ but at some nearby point $y: \| y - x^k \| \le {{{\delta }_{k}}}$. Moreover, an arbitrary convex combination of the gradients computed with error at points neighboring $x^k$ can be taken as the direction of motion $P_k$. This property, together with the general conditions (14.7), expresses a certain stability of the generalized gradient method (9.2), (9.3) and (10.7), (10.8).

Note that this stability differs from the real computational stability of the generalized gradient method of the form:
\begin{equation}
\label{eqn:14.8}
{{x}^{k+1}}={{x}^{k}}-{{\rho }_{k}}{{g}^{k}},\;\;\;{{g}^{k}}\in G({{x}^{k}}),\;\;\;k=0,1,\ldots. 
\end{equation}

In practical calculations by any method, errors of the values included in the method entry are inevitable. A method is stable if at small errors it arrives at a small neighborhood of the solution. For example, a computing machine always performs actions with rational numbers of limited digit capacity. If the computation occurs within the bounded digit grid, it is absolutely accurate. However, if any operation on numbers within the bounded digit grid results in the result of the operation going outside the bounded digit grid, then it is rounded or truncated. This is where rounding errors creep into the calculation. This leads to the fact that, for example, in the generalized gradient method the next point $x^{k+1}$ is calculated not by (14.8) but by the formula:
\begin{equation}
\label{eqn:14.9}
{{x}^{k+1}}=\overline{\bar{x}^{k}-\overline{\bar{\rho }_{k}\cdot \bar{g}_{k}}},\;\;\;{g}^{k}\in G(\bar{x}^{k}), 
\end{equation}
where the dash denotes rounding (i.e., in this formula, values without a dash denote numbers whose machine representations may extend beyond the bit grid, and values with a dash denote rational numbers within the bit grid). If in (14.9) the rounding errors at each iteration are small compared to the step size , then the study of the convergence of the algorithm (14.9) can be reduced to the study of the convergence of procedures of the form (14.4), (14.7) [85].

The following lemma establishes the local stability of the generalized gradient method.
\begin{lemma}
\label{lem:14.1}
Let the sequence of initial points $\left\{{ {{x}^{s}}} \right\}$ converges to $x=\underset{s\to \infty }{\mathop{\lim }}\,{{x}^{s}}$. For each $s$, consider sequences $\left\{ x_{s}^{k} \right\}_{k={{k}_{s}}}^{{{n}_{s}}}$ such that
\begin{equation}
\label{eqn:14.10}
x_{s}^{{{k}_{s}}}={{x}^{s}},\;\;\;x_{s}^{k+1}=x_{s}^{k}-\rho _{s}^{k}P_{s}^{k},
\end{equation}
\begin{equation}
\label{eqn:14.11}
P_{s}^{k}=Q_{s}^{k}+\Delta _{s}^{k},\;\;\;Q_{s}^{k}\in {{G}_{{{\delta }_{{{k}_{s}}}}}}(x_{s}^{k}),\;\;\;{{k}_{s}}\le k\le {{n}_{s}}.
\end{equation}
Let us denote
\begin{equation}
{{\rho }_{s}}=\underset{{{k}_{s}}\le k\le {{n}_{s}}}{\mathop{\sup }}\,\rho _{s}^{k},\;\;\;\;\;{{\delta }_{s}}=\underset{{{k}_{s}}\le k\le {{n}_{s}}}{\mathop{\sup }}\,\delta _{s}^{k},\nonumber
\end{equation}
\begin{equation}
{{\Delta }_{S}}=\underset{{{k}_{s}}\le k\le {{n}_{s}}}{\mathop{\sup }}\,\left\| \Delta _{s}^{k} \right\|,\;\;\;\;\;{{\sigma }_{s}}=\sum\limits_{k={{k}_{s}}}^{{{n}_{s}}-1}{{{\rho }_{k}}.}\nonumber
\end{equation}
If $0\notin G(x)$ and ${{{\sigma }_{s}}}\ge \sigma >0$, then there exists $\bar{\varepsilon }\,=\bar{\varepsilon }(x,\sigma )$ such that for any 
$\varepsilon\in(0,\bar{\varepsilon })$ there are $\bar{\rho }=\bar{\rho }(x,\varepsilon ),$ $\bar{\delta }=\bar{\delta }(x,\varepsilon )$ and 
$\bar{\delta } =\bar{\Delta }(x,\varepsilon )$ such that for 
$\left\{ x_{s}^{k} \right\}_{k}^{{n}_{s}}={k}_{s}$ defined by formulas (14.10), (14.11), where $\delta _{s}^{k}\le \bar{\delta }$, $\rho _{s}^{k}\le \bar{\rho }$, and $\Delta _{s}^{k}\le \bar{\Delta }$ there are  indexes $l_s$ such that $\left\| {{x}^{k}}-x \right\|\le \varepsilon$ for $k\in \left[ {{k}_{s}},{{l}_{s}} \right]$ and
\begin{equation}
1) f(x)=\underset{s\to \infty }{\mathop{\lim }}\,f({{x}^{s}})>\underset{s\to \infty }{\overset{\_\_\_\_\_}{\mathop{\lim }}}\,f(x_{s}^{_{_{{{l}_{s}}}}}),\;\;\;h(x)<0;
\nonumber
\end{equation}
\begin{equation}
2) h(x)=\underset{s\to \infty }{\mathop{\lim }}\,h({{x}^{s}})>\underset{s\to \infty }{\overset{\_\_\_\_\_}{\mathop{\lim }}}\,h(x_{s}^{{{l}_{s}}}),\;\;\;h(x)>0;
\nonumber
\end{equation}
\begin{equation}
3) f(x)=\underset{s\to \infty }{\mathop{\lim }}\,f({{x}^{s}})>\underset{s\to \infty }{\overset{\_\_\_\_\_}{\mathop{\lim }}}\,f(x_{s}^{{{l}_{s}}}),\nonumber
\end{equation}
\begin{equation}
C(x)\overset{\_}{\mathop{\delta }}\,\ge \underset{s\to \infty }{\overset{\_\_\_\_\_}{\mathop{\lim }}}\,h(x_{s}^{{{l}_{s}}}),\;\;\;  C(x)>0,\;\;\;h(x)=0.\nonumber
\end{equation}
\end{lemma}
{\it P r o o f}. Recall that ${{G}_{\delta }}(x)=co\left\{ G(y)\left| \left\| y-x \right\|\le \delta \right. \right\}.$ Due to semi-continuity from above of $G(x)$, there exists an ${{\varepsilon }_{1}}$-neighborhood of point $x$ such that ${{\gamma }_{1}}=\rho (0,{{G}_{{\varepsilon}_{1}}}(x))>0.$ Let us denote 
\[\Gamma_f=\sup \left\{ \left\| g \right\|\left| g\in {{G}_{f}}(y),\left\| y-x \right\| \right.\le {{\varepsilon }_{1}} \right\},\]
\[\Gamma_h=\sup \left\{ \left\| g \right\|\left| g\in {{G}_{h}}(y),\left\| y-x \right\|\le {{\varepsilon }_{1}} \right. \right\},\]
\[\Gamma=\max ({{\Gamma}_{f}},{{\Gamma}_{h}})+{{\gamma }_{1}}<+\infty .\]
Let us take ${{\varepsilon }_{2}}=\min ({{\varepsilon }_{1}},\sigma {{\gamma }_{1}})$. By virtue of generalized differentiability of $f (x)$ and $h (x)$, for any constants $C_f$ and $C_h$, there is ${{{\varepsilon }_{3}}\le {{\varepsilon }_{2}}}$ such that when $\left\| y-x \right\|\\\le {{\varepsilon }_{3}}$, it holds
\[f(y)=f(x)+\left<{{g}_{f}},y-x\right>+{{O}_{f}}(x,y,{{g}_{1}}),\]
\[h(y)=h(x)+\left<{{g}_{h}},y-x\right>+{{O}_{h}}(x,y,{{g}_{h}}),\]
where
$\left| {{O}_{f}}(x,y,{{g}_{f}}) \right|\le {{C}_{f}}\left\| y-x \right\|$, 
$\left| {{O}_{h}}(x,y,{{g}_{h}}) \right|\le {{C}_{h}}\left\| y-x \right\|$, 
${{g}_{f}}\in {{G}_{f}}(y)$, ${{g}_{h}}\in {{G}_{h}}(y)$. 
The constants $C_f$ and $C_h$, will be specified later. If $h(x)\ne 0$, then there is ${{\varepsilon }_{4}}\le {{\varepsilon }_{3}}$ such that for all $y$, 
$\left\| y-x \right\|\le {{\varepsilon }_{4}}$, it holds $\left| h(y)-h(x) \right|\le \left| h(x) \right|/2$. If $h(x)=0$, then we assume 
${{{\varepsilon }_{4}}={{\varepsilon }_{3}}}$.  

Let's put 
$\bar{\varepsilon }={{\varepsilon }_{4}}$. 
Let's fix an arbitrary 
$\varepsilon \in (0,\bar{\varepsilon} ]$. 
Let's take 
$\bar{\delta }=\varepsilon /4$, $\bar{\rho }_{1}=\varepsilon /(4\Gamma)$, 
$\bar{\Delta }_{1}={{\gamma }_{1}}/2\,$. 
Let for $s\ge S$ holds
$\left\| {{x}^{s}}-x \right\|\le \varepsilon /4$, 
${\delta}_{s}\le \bar{\delta }_{1}$, 
${{\Delta }_{s}}\le \bar{\Delta }_{1}$, 
${{\rho }_{s}}\le \bar{\rho }_{1}$. 

We denote by ${{m}_{s}}={{m}_{s}}(\varepsilon )$ the maximal index $k\in [{{k}_{s}}, {{n}_{s}}]$ such that for all $r\in [{{k}_{s}},k)$, $\left\| x_{s}^{r}-x \right\|\le \varepsilon /2$ holds.

Let us show that $\varepsilon /2\le \left\| x_{s}^{{{{m}_{s}}}}-x \right\|\le 3\varepsilon /4$. Let us prove the left inequality first. Indeed, if 
$\left\| x_{s}^{{{m}_{s}}}-x \right\|\le \varepsilon /2$, then ${{m}_{s}}={{{n}_{s}}}$, 
and we get a contradiction: 
${\varepsilon }_{2}>3\varepsilon /4\ge \left\| x_{s}^{{n}_{s}}-{{x}^{s}} \right\|\ge \sigma {{\gamma }_{1}}$. Next,
\[\left\| x_{s}^{{{m}_{s}}}-x \right\|\le \left\| x_{s}^{{{m}_{s}}-1}-x \right\|+\rho _{s}^{{{m}_{s}}-1}\left\| P_{s}^{{{m}_{s}}-1} \right\|\le 3\varepsilon /4.\]
From here it follows
\[\frac{3\varepsilon }{4}\ge \left\| \sum\limits_{k={{k}_{s}}}^{{{m}_{s}}-1}{{{\rho }_{k}}{{P}^{k}}} \right\|\ge \frac{{{\gamma }_{1}}}{2}\sum\limits_{k={{k}_{s}}}^{{{m}_{s}}-1}{{{\rho }_{k}},}\]
\[\sum\limits_{k={{k}_{s}}}^{{{m}_{s}}-1}{{{\rho }_{k}}\le a=\frac{3\varepsilon }{2{{\gamma }_{1}}}.}\]
Let us consider $g_{s}^{k}\in {{G}_{\delta _{s}^{k}}}(x_{s}^{k})$. The following  representation holds true:
\begin{eqnarray}
g_{s}^{k}&=&\sum\limits_{i=1}^{n+1}\lambda_{s}^{{k}_{i}}g_{s}^{k_i}
	=\sum\limits_{i=1}^{n+1}\lambda_{s}^{{k}_{i}}a_{s}^{{k}_{i}}g_{{f}{s}}^{{k}_{i}}+
	\sum\limits_{i=1}^{n+1}\lambda_{s}^{{k}_{i}}\beta_{s}^{{k}_{i}}g_{{h}{s}}^{{k}_{i}} \nonumber\\ 
 & =&\left(\sum\limits_{i=1}^{n+1}\mu_{s}^{{k}_{i}}g_{{f}{s}}^{{k}_{i}}\right)\sum\limits_{i=1}^{n+1}\lambda_{s}^{{k}_{i}}a_{s}^{{k}_{i}}
+\left(\sum\limits_{i=1}^{n+1}\nu_{s}^{{k}_{i}}g_{{h}{s}}^{{k}_{i}}\right)
\sum\limits_{i=1}^{n+1}\lambda_{s}^{{k}_{i}}\beta_{s}^{{k}_{i}} \nonumber\\ 
 & =&\bar{a}_{s}^{k}\bar{g}_{{f}{s}}^{k}+
\bar{\beta}_{s}^{k}\bar{g}_{{h}{s}}^{k}, \label{14.12} 
\end{eqnarray}
where $\lambda _{s}^{{{k}_{i}}}$, $a_{s}^{{{k}_{i}}}$, $\beta _{s}^{{{k}_{i}}}$, $\mu _{s}^{{{k}_{i}}}$, $\nu _{s}^{{{k}_{i}}}$, $\bar{a}_{s}^{k}$, $\bar{\beta}_{s}^{k}\ge 0,$
\[\sum\limits_{i=1}^{n+1}{\lambda _{s}^{{{k}_{i}}}}=\sum\limits_{i=1}^{n+1}{\mu _{s}^{{{k}_{i}}}}=\sum\limits_{i=1}^{n+1}{\nu _{s}^{{{k}_{i}}}}=a_{s}^{{{k}_{i}}}+\beta _{s}^{{{k}_{i}}}=\bar{a}_{s}^{k}+\bar{\beta}_{s}^{k}\,=1,\]
\[g_{s}^{ki}\in G(x_{s}^{{{k}_{i}}}),\;\;\;g_{fs}^{ki}\in Gf(x_{s}^{{{k}_{i}}}),\;\;\;            g_{hs}^{ki}\in {{G}_{h}}(x_{s}^{{{k}_{i}}}),\]
\[\bar{g}_{fs}^{k}\in co\left\{ {G}_{f}(y)\left| \left\| y-x_{s}^{k} \right\|\le \delta_{s}^{k} \right. \right\}\,,\]
\[\bar{g}_{hs}^{k}\in co\left\{ {G}_{h}(y)\left| \left\| y-x_{s}^{k} \right\|\le \delta _{s}^{k} \right. \right\}\,.\]
If $\left\| {{x}^{s}}-x \right\|\le \varepsilon /4$, ${{\delta }_{s}}\le \varepsilon /4$, ${{k}_{s}}\le k\le {{m}_{s}}$, $1\le i\le n+1$, then
$\left\| x_{s}^{{{k}_{i}}}-x \right\|\le \left\| x_{s}^{{{k}_{i}}}-x_{s}^{k} \right\|+\left\| x_{s}^{k}-x \right\|\le \varepsilon \le {{\varepsilon }_{3}}.$ So the following
expansion holds true
\begin{eqnarray}
   f(x_{s}^{{{k}{i}}})&=&f(x)+\left<g_{fs}^{ki},x_{s}^{ki}-x\right>+{{O}_{f}}(x,x_{s}^{ki},g_{fs}^{ki})  \nonumber\\ 
 & \le& f(x)+\left<g_{f_s}^{ki},x_{s}^{ki}-{{x}^{s}}\right>+{{C}_{f}}\left\| x_{s}^{ki}-{{x}^{s}} \right\|+({{\Gamma}_{f}}+{{C}_{f}})\left\| {{x}^{3}}-x \right\|. \nonumber 
\end{eqnarray}
In this expansion, $x_{s}^{ki}\;(1\le i\le n+1)$ can be approximated by $x_{s}^{k}$; the inequality holds
\begin{eqnarray}
  f(x_{s}^{k})&\le& f(x)+\left<g_{fs}^{ki},x_{s}^{k}-{{x}^{s}}\right>+{{C}_{f}}\left\| x_{s}^{k}-{{x}^{s}} \right\| \nonumber\\ 
 && +(2{{\Gamma}_{f}}+{{C}_{f}}){{\delta }_{s}}+({{\Gamma}_{f}}+{{C}_{f}})\left\| {{x}^{s}}-x \right\|. \nonumber 
\end{eqnarray}
Let us multiply these inequalities by $\mu _{s}^{ki}$ and sum over $i$:
\begin{eqnarray}
   f(x_{s}^{k})&\le& f(x)+(\bar{g}_{fs}^{k},x_{s}^{k}-{{x}^{s}})+{{C}_{f}}\left\| x_{s}^{k}-{{x}^{s}} \right\|\, \nonumber\\ 
 && +(2{{\Gamma}_{f}}+{{C}_{f}}){{\delta }_{s}}+({{\Gamma}_{f}}+{{C}_{f}})\left\| {{x}^{s}}-x \right\|. \label{14.13} 
\end{eqnarray}
Let us represent
\[x_{s}^{k}=\bar{x}_{s}^{k}\,+z_{s}^{k},\]
\[\bar{x}_{s}^{k}={{x}^{s}}-\sum\limits_{r={{k}_{s}}}^{k-1}{\rho _{s}^{r}g_{s}^{r},}\,\;\;\;\;\;
z_{s}^{k}=-\sum\limits_{r=s}^{k-1}{\rho _{s}^{r}\Delta _{s}^{r}.}\]
The estemate holds true for $k\in [{{k}_{s}},{{m}_{s}}]:$
\[\left\| z_{s}^{k} \right\|\le \left(\underset{s\le r\le {{n}_{s}}}{\mathop{\sup }}\,\Delta _{s}^{r}\right)\sum\limits_{r={{k}_{s}}}^{{{m}_{s}}-1}{\rho _{s}^{r}\le a{{\Delta }_{s}}.}\]
We replace $x_{s}^{k}$ in the right-hand side of (14.13) by $\bar{x}_{s}^{k}\,$:
\begin{eqnarray}
   f(x_{s}^{k})&\le& f(x)+\left<\bar{g}_{fs}^{k}\,,\bar{x}_{s}^{k}-{{x}^{s}}\right>+{{C}_{f}}       \left\| \bar{x}_{s}^{k}-{{x}^{s}} \right\|\, \nonumber\\ 
 && +(2{{\Gamma}_{f}}+{{C}_{f}}){{\delta}_{s}}+({{\Gamma}_{f}}a+{{C}_{f}}){{\Delta}_{s}}\nonumber\\
&&+({{\Gamma}_{f}}+{{C}_{f}})\left\| {{x}^{s}}-x \right\|. \label{14.14} 
\end{eqnarray}
Similarly,
\begin{eqnarray}
   h(x_{s}^{k})&\le& h(x)+\left<\bar{g}_{hs}^{k},\bar{x}_{s}^{k}\,\,-{{x}^{s}}\right>+{{C}_{h}}\left\| \bar{x}_{s}^{k}\,-{{x}^{s}} \right\|+ \nonumber\\ 
 && +(2{{\Gamma}_{h}}+{{C}_{h}}){{\delta }_{s}}+({{\Gamma}_{h}}a+{{C}_{h}}){{\Delta }_{s}}\nonumber\\
&&+({{\Gamma}_{h}}+{{C}_{h}})\left\| {{x}^{s}}-x \right\|. \label{14.15} 
\end{eqnarray}
Now consider the scalar product
\[\left<g_{s}^{k},\bar{x}_{s}^{k}\,-{{x}_{s}}\right>=-\left<g_{s}^{k},\sum\limits_{r={{k}_{s}}}^{k-1}{\rho _{s}^{r}g_{s}^{r}/\sum\limits_{r={{k}_{s}}}^{k-1}{\rho _{s}^{r}}}\right>\sum\limits_{r={{k}_{s}}}^{k-1}{\rho _{s}^{r}.}\]
For $s\ge S$ and $k\in [{{k}_{s}},{{m}_{s}}]$ it takes place $g_{s}^{k}\in {{G}_{{{\varepsilon}_{1}}}}(x)$, so ${{\gamma }_{1}}/2\le \left\| P_{s}^{k} \right\|\le \Gamma$, and
\[\sum\limits_{k={{k}_{s}}}^{{{m}_{s}}-1}{\rho _{s}^{k}\ge \left\| x_{s}^{{{m}_{s}}}-{{x}^{s}} \right\|/\Gamma\ge \varepsilon /(4\Gamma).}\]
Obviously, to evaluate the scalar products $\left<g_{s}^{k},\bar{x}_{s}^{k}\,-{{{x}^{s}}}\right>$ at $s\ge S$ we apply Lemma 9.2, where we need to put
\[P={{G}_{{{\varepsilon}_{1}}}}(x),\;\;\;            {{p}^{k}}=g_{s}^{k},\;\;\;            \sigma =\varepsilon /(4\Gamma),\;\;\;              \rho =\varepsilon \gamma _{1}^{2}/(24{{\Gamma}^{3}}).\]
Let's put ${\bar{\rho }}\,=\min (\bar{\rho }_{1}\,,\rho )$. If $s\ge S$ and ${{\rho }_{s}}\le \bar{\rho }\,$, then, according to Lemma 9.2, there are indices ${{l}_{s}}\le {{m}_{s}}$ for which $\left\| x_{s}^{k}-x \right\|\le \varepsilon $ at $k\in [{{k}_{s}},{{l}_{s}}]$, and
\[\left<g_{s}^{{l}_{s}},\sum\limits_{k={{k}_{s}}}^{{{l}_{s}}-1}\rho _{s}^{k}g_{s}^{k}/\sum\limits_{k={{k}_{s}}}^{{l}_{s}-1}\rho _{s}^{k}\right>\ge \frac{\gamma _{1}^{2}}{4},\]
\begin{equation}
\label{eqn:14.16}
\sum\limits_{t={{k}_{s}}}^{{{l}_{s}}-1}{\rho _{s}^{k}\ge \frac{\varepsilon {{\gamma }_{1}}}{12{{\Gamma}^{2}}}}.
\end{equation}

Three cases need to be considered: $h(x)<0$, $h(x)>0,$ and $h(x)=0$.

Let $h(x)<0$. For $s$ sufficiently large and $k\in [{{k}_{s}},{{m}_{s}}]$ it holds $g_{s}^{k}=\bar{g}_{fs}^{k}\,$, so  inequality (14.14) can be rewritten for $k={{{l}_{s}}}$  as
\begin{eqnarray}
   f(x_{s}^{{{l}_{s}}})&\le& f(x)-\left(\frac{\gamma _{1}^{2}}{4}-{{\Gamma}C_{f}}\right)\sum\limits_{k={{k}_{s}}}^{{{l}_{s}}-1}{\rho _{s}^{k}}\nonumber \\ 
 && +(2{{\Gamma}_{f}}+{{C}_{f}}){{\delta }_{s}}+({{\Gamma}_{f}}a+{{C}_{f}}){{\Delta }_{f}}\nonumber\\
&&+({{\Gamma}_{f}}+{{C}_{f}})\left\| {{x}^{s}}-x \right\|. \label{14.17} 
\end{eqnarray}
Now let us instantiate $C_f$ (and hence $\overset{-}{\mathop{\varepsilon }}\,$): take ${{C}_{f}}=\gamma _{1}^{2}/(8\Gamma)$. We have $\epsilon\in(0,\bar{\varepsilon }\,]$ and $\bar{\rho }\,=\min (\bar{\rho }_{1}\,,\rho )$. Let's take
 \[\bar{\delta }_{f}\,=\frac{\varepsilon \gamma _{1}^{3}}{300{{\Gamma}^{2}}(2{{\Gamma}_{f}}+{{C}_{f}})},\;\;\;\;\;           
{{\bar{\Delta }}\,}_{f}=\frac{\varepsilon \gamma _{1}^{3}}{300{{\Gamma}^{2}}({{\Gamma}_{f}}a+{{C}_{f}})}.\]
All the reasoning carried out earlier will turn from conditionally true to unconditionally true. As a result, from (14.16), (14.17) when ${\rho }_{s}\le \bar{\rho }\,$, ${{\delta }_{s}}\le \bar{\delta }_{f}\,$ and ${{\Delta }_{s}}\le \bar{\Delta }_{s}\,$, we obtain
\begin{equation}
\label{14.18}
f(x_{s}^{{{l}_{s}}})\le f(x)-\frac{\varepsilon \gamma _{1}^{3}}{300{{\Gamma}^{2}}}+({{\Gamma}_{f}}+{{C}_{f}})\left\| {{x}^{s}}-x \right\|,
\end{equation}
whence follows the statement 1) of the lemma.

Assertion 2) is proved analogously.

Now let $h(x)=0$. Let us prove statement 3) of the lemma.

We show that, given an appropriate choice of constants $C_h$, $\bar{\delta }\,$ and $\bar{\Delta }\,$ for sufficiently large $s$, there will be 
$\min \left\{ h(y)|\, \left\| y-x_{s}^{{l}_{s}} \right\|\le {\delta }_{s} \right\}\le 0$, which entails $h(x_{s}^{{l}_{s}})\le 2{\Gamma}_{h}{\delta }_{s}=C(x){\delta }_{s}$ and proves the second part of statement 3) of the lemma. Indeed, suppose that $\min \left\{ h(y)| \left\| y-x_{s}^{{l}_{s}} \right\|\le {{\delta }_{s}} \right\}>0$. Then 
$g_{s}^{{l}_{s}}=\bar{g}_{h_s}^{{l}_{s}}\,$, and from (14.15), we obtain
\begin{eqnarray}
   0&\le& h(x_{s}^{{{l}_{s}}})\le -\left(\frac{\gamma _{1}^{2}}{4}-{{C}_{h}}\right)\sum\limits_{k={{k}_{s}}}^{{{l}_{s}}-1}{\rho _{s}^{k}+(2{{\Gamma}_{h}}+{{C}_{h}}){{\delta }_{s}}}\nonumber \\ 
 && +({{\Gamma}_{h}}a+{{C}_{h}}){{\Delta }_{s}}+({{\Gamma}_{h}}+{{C}_{h}})\left\| {{x}^{s}}-x \right\|. \nonumber 
\end{eqnarray}
For ${{C}_{h}}\le \gamma _{1}^{2}/(8\Gamma)$ and
 \[\bar{\delta }_{f}\,=\frac{\varepsilon \gamma _{1}^{3}}{300{{\Gamma}^{2}}(2{{\Gamma}_{h}}+{{C}_{h}})},\;\;\;\;\;           
{{\bar{\Delta }}\,}_{f}=\frac{\varepsilon \gamma _{1}^{3}}{300{{\Gamma}^{2}}({{\Gamma}_{h}}a+{{C}_{h}})},\]
for sufficiently large $s$ we obtain a contradiction:
\[0\le -\frac{\varepsilon \gamma _{1}^{3}}{300{{\Gamma}^{2}}}+({{\Gamma}_{h}}+{{C}_{h}})\left\| {{x}^{s}}-x \right\|.\]
The second part of statement 3) of the lemma is proved.

Now, depending on the location of $x_{s}^{{{{l}_{s}}}}$ we consider two cases:

a) $\max \left\{ h(y)\left| \left\| y-x_{s}^{{{l}_{s}}} \right\|\le {{\delta }_{s}} \right. \right\}<0$; then $g_{s}^{{l}_{s}}=\bar{g}_{fs}^{{l}_{s}}\,$, the inequalities (14.14), (14.17) and (14.18) hold true, whence at ${{C}_{f}}\le \gamma _{1}^{2}/(8\Gamma)$ (and corresponding ${\varepsilon \le \bar{\varepsilon }}\,$), and when ${\Delta }_{s}\le \bar{\delta }_{f}\,$ and ${\Delta }_{s}\le \bar{\Delta }_{s}\,$, it follows the first part of statement 3) of the lemma;

b) $\max \left\{ h(y)\left| \left\| y-x_{s}^{{l}_{s}} \right\|\le {{\delta }_{s}} \right. \right\}\ge 0$; then $h(x_{s}^{{l}_{s}})\ge -2{\Gamma}_{h}{\delta }_{s}$ and from (14.15) with $k={{l}_{s}}$, we obtain
\begin{eqnarray}
   \left<\bar{g}_{hs}^{{l}_{s}},\,\sum\limits_{k={{k}_{s}}}^{{{l}_{s}}-1}{\rho _{s}^{k}g_{s}^{k}}\right>&\le& {{C}_{h}}\left\| \sum\limits_{k={{k}_{s}}}^{{{l}_{s}}-1}{\rho _{s}^{k}g_{s}^{k}} \right\|+(4{{\Gamma}_{h}}+{{C}_{h}}){{\delta }_{s}}
	\nonumber\\ 
 && +({{\Gamma}_{h}}a+{{C}_{h}}){{\Delta }_{s}}+({{\Gamma}_{h}}+{{C}_{h}})\left\| {{x}^{s}}-x \right\|. \label{14.19} 
\end{eqnarray}
From (14.12), (14.16), (14.19), it follows
\begin{eqnarray}
   \bar{\alpha}_{s}^{{l}_{s}}\,\left<\bar{g}_{fs}^{{l}_{s}}\,,\sum\limits_{k={{k}_{s}}}^{{{l}_{s}-1}}{\rho_{s}^{k}g_{s}^{k}}\right>
	&\ge& \frac{\gamma _{1}^{2}}{4}\sum\limits_{k={{k}_{s}}}^{{{l}_{s}}-1}\rho _{s}^{k}-\bar{\beta}_{s}^{{{l}_{s}}}\left(\,\sum\limits_{k={{k}_{s}}}^{{{l}_{s}}-1}{\rho _{s}^{k}g_{s}^{k}}\right)  \nonumber\\ 
 & \ge& \left(\frac{\gamma _{1}^{2}}{4}-{{C}_{h}}\right)\sum\limits_{k={{k}_{s}}}^{{{l}_{s}}-1}{\rho _{s}^{k}-(4{{}_{h}}+{{C}_{h}}){{\delta }_{s}}} \nonumber\\ 
 && -({{\Gamma}_{h}}a+{{C}_{h}}){{\Delta }_{s}}
-({{\Gamma}_{h}}+{{C}_{h}})\| {{x}^{s}}-x\|. \label{14.20} 
\end{eqnarray}
Now let us choose the value of $C_h$; all the reasoning we have done is true when ${{\Gamma}_{h}}=\gamma _{1}^{2}/(8\Gamma)$. If
 \[\bar{\delta }_{s}\,\le\bar{\delta}\le
\frac{\varepsilon \gamma _{1}^{3}}{200{{\Gamma}^{2}}(4{{\Gamma}_{h}}+{{C}_{h}})},\;\;\;\;\;
{{\bar{\Delta }}\,}_{s}\le\bar{\Delta}\le
\frac{\varepsilon \gamma _{1}^{3}}{200{{\Gamma}^{2}}({{\Gamma}_{h}}a+{{C}_{h}})},\]
then for sufficiently large $s$ from (14.20) we obtain:
\begin{eqnarray}
\left<\bar{g}_{fs}^s,\sum_{k=k_s}^{l_s-1}\rho_s^kg_s^k\right>
&\ge&
\frac{\gamma_1^2}{8}\sum_{k=k_s}^{l_s-1}\rho_s^k
-\left(4\Gamma_h+C_h\right)\delta_s \nonumber\\
&&
-\left(\Gamma_ha+C_h\right)\Delta_s-\left(\Gamma_h+C_h\right)
\|x^s-x\|.\label{14.21}
\end{eqnarray}
From (14.14), (14.21) for $k=l_s$ and for sufficiently large $s$, we obtain
\begin{eqnarray}
f\left(x_s^{l_s}\right)&\le&f(x)-\left(\frac{\gamma_1^2}{8}-\Gamma C_f\right)\sum_{k=k_s}^{l_s-1}\rho_s^k+
\left(\Gamma_fa+\Gamma_ha+C_f+C_h\right)\Delta_s  \nonumber\\
&&
+\left(2\Gamma_f+4\Gamma_h+C_f+C_h\right)\delta_s
+\left(\Gamma_f+\Gamma_h+C_f+C_h\right)\Delta_s\|x^s-x\|.\nonumber
\end{eqnarray}
Let us now concretize $C_f$ (since $C_h$ has been already stated earlier, thus we concretize also $\bar{\varepsilon}$): let $C_f = \gamma_1^2/(16\gamma)$. We have
$0\le\varepsilon\le\bar{\varepsilon}$ and $\bar{\rho}=\min\left(\bar{\rho_1},\rho\right)$. 
Now let us define $\bar{\delta}$ and $\bar{\Delta}$. Let us put
 \[\bar{\delta}=
\frac{\varepsilon \gamma _{1}^{3}}{200{{\Gamma}^{2}}
(2{\Gamma}_{f}+4{{\Gamma}_{h}}+C_f+{{C}_{h}})},\;\;\;\;\;           
\bar{\Delta}=
\frac{\varepsilon \gamma _{1}^{3}}{200{{\Gamma}^{2}}
({\Gamma}_{f}a+{{\Gamma}_{h}}a+C_f+{{C}_{h}})}.\]

All the reasoning we did earlier will turn from conditionally true to unconditionally true.  As a result, at $\rho_s\le\bar{\rho}$, $\delta_s\le\bar{\delta}$, and $\Delta_s\le\bar{\Delta}$, we obtain
\[
f\left(x_s^{l_s}\right)\le f(x)-\frac{\varepsilon\gamma^3}{600\Gamma^2}
+(\Gamma_f+\Gamma_h+C_f+C_h)\|x^s-x\|,
\]
whence also follows the first part of statement 3) of the lemma. Thus, we specified $\bar{\varepsilon}(x,\sigma)$ and $\bar{\rho}(x,\varepsilon)$, $\bar{\delta}(x,\varepsilon)$, $\bar{\Delta}(x,\varepsilon)$ ($\varepsilon\in(0,\bar{\varepsilon}]$) under which statements 1) - 3) of the lemma are satisfied.
\begin{theorem}
\label{th:14.1}
Let in problem (14.1)-(14.2) $f(x)$ and $h(x)$ 
are generalized differentiable functions, and $h(x)\rightarrow+\infty$ for $\|x\|\rightarrow+\infty$. Let the set $H^*=\{h(x)|\,0\in G_h(x), h(x)>0\}$ contains no intervals. Then for sufficiently small $\rho$, $\delta$, and $\Delta$ from (14.6), the sequence (14.4) - (14.7) is bounded, and

either all limit points of $\{x^k\}$ do not belong to the admissible region $D=\{x\in\mathbb{E}_n|\,h(x)\le 0\}$, and then they all belong to 
$X_h^*=\{x|\,0\in G_h(x),\,h(x)>0\}$ and there exists a limit 
$\lim_{k\rightarrow\infty}h(x^k) > 0$;

or all limit points of $\{x^k\}$ belong to $D$, and then the minimal ones (by value of $f$) belong to $X^* = \{x\in\mathbb{E}_n|\,0\in G(x), h(x)\le 0\}$, and the segment $[\underline{\lim}_{k\rightarrow\infty}f(x^k),\,\overline{\lim}_{k\rightarrow\infty}f(x^k)]$ is nested in the set 
$F^*=\{f(x)|\,x\in X^*\}$. If the set $F^*$ contains no intervals, then all the limit points $\{x^k\}$ belong to a connected subset of $X^*$ and
there exists a limit $\lim_{k\rightarrow\infty}f(x^k)$.
\end{theorem}

This convergence theorem for method (14.4) - (14.7) is proved in the same way as Lemmas 9.1, 10.1 and Theorems 10.1, 10.2 on the convergence of the generalized gradient method (10.7), (10.8)  due to the validity of Lemma 14.1, similar to Lemma 10.2.

\smallskip
{\bf 2. The method of averaged gradients.} 
\label{Sec.14.2}
The method is designed to solve  problem (14.1), (14.2). It produces an infinite sequence of points $\{x^k\}$ according to the following rules:
\begin{equation}
\label{eqn:14:22}
x^0\in {E}_n \mbox{ is  arbitrary},
\end{equation}
\begin{equation}
\label{eqn:14:23}
x^{k+1}=x^k-\rho_k\nu_kP^k,\;\;\;P^k=Q^k+\Delta^k,\;\;\;\rho_k\ge 0,\;\;\;
\nu_k>0,
\end{equation}
\begin{equation}
\label{eqn:14:24}
Q^k=\sum_{r=r_k}^k\lambda_{kr}g^r,\;\;\;g^r\in G(x^r),\;\;\;\lambda_{kr}\ge 0,
\end{equation}
where $x^k$ is the current approximation, $P^k$ is the direction of motion, $Q^k$ is an average gradient, $\Delta^k$ is some small additive, 
$\rho_k = \rho_k(x^0,g^0,\ldots,x^k,g^k)$ is an adjustable step multiplier, $\nu_k = \nu_k(x^0,g^0,\ldots,x^k,g^k)$ is an adjustable normalizing multiplier, $\lambda_k = \lambda_k(x^0,g^0,\ldots,x^k,g^k)$ is an adjustable averaging multiplier.

The following conditions are met.

For any compact $K\subset{E}_n$ there are numbers 
$\nu_K>0$, $\lambda_K>0$, $N_K>0$, $\Lambda_K>0$, such that if
$\{x^k\}\subset K$ and $g^k\in G(x^k)$, then for all $k$ holds
\begin{equation}
\label{eqn:14:25}
0<\nu_K\le\nu_k(x^0,g^0,\ldots,x^k,g^k)\le N_K<+\infty,
\end{equation}
\begin{equation}
\label{eqn:14:26}
0<\lambda_K\le\sum_{r=r_k}^k\lambda_{kr}(x^0,g^0,\ldots,x^k,g^k)\le 
\Lambda_K<+\infty.
\end{equation}

For any bounded sequence $\{x^k\}$ and sequence $\{g^k\}$ 
$(g^k\in G(x^k))$ the following holds
\begin{equation}
\label{eqn:14:27}
\lim_{k\rightarrow\infty}\Delta_k(x^0,g^0,\ldots,x^k,g^k)=0,
\end{equation}
\begin{equation}
\label{eqn:14:28}
\max_{r_k\le r\le k}\|x^k-x^r\|\le\varepsilon_k(x^0,g^0,\ldots,x^k,g^k),
\end{equation}
\begin{equation}
\label{eqn:14:29}
\lim_{k\rightarrow\infty}\varepsilon_k(x^0,g^0,\ldots,x^k,g^k)=0,
\end{equation}
(function $\varepsilon_k$ majorizes $\max_{r_k\le r\le k}\|x^k-x^r\|$),
\begin{equation}
\label{eqn:14:30}
\lim_{k\rightarrow\infty}\rho_k(x^0,g^0,\ldots,x^k,g^k)=0.
\end{equation}

Functions $\rho_k(x^0,g^0,\ldots,x^k,g^k)$ have the following property: if
$\{x^k\rightarrow x\}$, $g^k\in G(x^k)$ and $0\notin G(x),$ then
\begin{equation}
\label{eqn:14:31}
\sum_{k=0}^\infty\rho_k(x^0,g^0,\ldots,x^k,g^k)=+\infty.
\end{equation}
\begin{remark}
\label{rem:14.1}
The vector $Q^k$ in (14.23) is some average of the previous generalized gradients, which explains the name of the method.
 It is easy to see that in the case of unconditional minimization of the function $f(x)$, the vector $Q^k$ is some generalized gradient of the of the averaged function
\[
\phi_k(x)=\sum_{r=r_k}^k\lambda_{kr}f(x-(x^k-x^r)).
\]
Therefore, the method of averaged gradients can also be considered as a method for solving a specific nonstationary (limit) extremal problem, in which a sequence of functions $\phi_k(x)$ is minimized.
\end{remark}
\begin{remark}
\label{rem:14.2}
In the method of averaged gradients (14.22) - (14.24) we can use normalized generalized gradients
\[
\bar{g}^r=g^r/(\|g^r\|+\varepsilon),\;\;\;\;\;\epsilon>0.
\]
In this case, the averaged gradients $Q^k$ have the form
\[
Q^k=\sum_{r=r_k}^k\lambda_{kr}g^r=
\sum_{r=r_k}^k\frac{\mu_{kr}}{\|g^r\|+\varepsilon}g^r,
\]
where
\[
\mu_{kr}\ge 0,\;\;\;\;\;\sum_{r=r_k}^k\mu_{kr}=1.
\]
We can also normalize the directions of motion $P^k$; for this purpose we should take normalizing multipliers of the form
\[
\nu_k(x^0,g^0,\ldots,x^k,g^k)=\frac{1}{\|P^k\|+\epsilon},
\;\;\; \epsilon>0.
\]
Due to the boundedness of the mapping G on the compacta, the method of averaged gradients with normalization of generalized gradients and vectors $P^k$ satisfies conditions (14.25)-(14.31).
\end{remark}

The following lemma establishes a local minimizing property of the method of averaged gradients, which basically ensures its convergence.

\begin{lemma}
\label{lem:14.2}
Suppose that the sequence $\{x^k\}$ generated by the algorithm (14.22) - (14.31) is bounded and there is a limit point 
$x^\prime = \lim_{s\rightarrow\infty} x^{k_s}$ such that $0\notin G(x).$ Then it can be found $\bar{\epsilon}>0$ such that for any 
$\epsilon\in (0,\bar{\epsilon})$ there exist indices $l_s\ge k_s$, for which (for sufficiently large $s$) $\|x^k -x^\prime\|\le\epsilon$ at 
$k\in [k_s,l_s]$, and it holds true:

1) $f(x^\prime)=\lim_{s\rightarrow\infty}f(x^{k_s})>
\overline{\lim}_{s\rightarrow\infty}f(x^{l_s}), \;\;\;h(x)<0;$

2) $h(x^\prime)=\lim_{s\rightarrow\infty}h(x^{k_s})>
\overline{\lim}_{s\rightarrow\infty}h(x^{l_s}), \;\;\;h(x)>0;$

1) $f(x^\prime)=\lim_{s\rightarrow\infty}f(x^{k_s})>
\overline{\lim}_{s\rightarrow\infty}f(x^{l_s}), \;\;\;
0\ge \overline{\lim}_{s\rightarrow\infty}h(x^{l_s}),
\;\;\;h(x)=0.$
\end{lemma}
{\it P r o o f}. Let us consider that $\{x^k\}\subset K$, where $K$ is a compact.
Denote $\gamma=\rho(0,G(x^\prime))>0$,
\[
\Gamma=\sup_{k\ge 0}\|g^k\|<+\infty,\;\;\;
G_\epsilon(x)=co\{G(y)|\,\|y-x\|\le\epsilon\}.
\]
By virtue of semi-continuity from above of the convex-valued mapping $G$, 
there exists $\epsilon^\prime>0$ such that 
$\rho(0,G_{\epsilon^\prime}(x^\prime))\ge\gamma/2$. Let us represent
\[
x^{k+1}=x^k-\rho_k\nu_k P^k=x^k-\bar{\rho}_k\bar{P}^k,
\]
where
\[
\bar{P}^k=\bar{Q}^k+\bar{\Delta}^k=P^k/\sum_{r=r_k}^k\lambda_{kr},\;\;\;
\bar{Q}^k=Q^k/\sum_{r=r_k}^k\lambda_{kr},
\]
\[
\bar{\Delta}^k=\Delta/\sum_{r=r_k}^k\lambda_{kr},\;\;\;
\bar{\rho}^k=\rho_k\nu_k\sum_{r=r_k}^k\lambda_{kr}.
\]
Let us show that the whole sequence $\{x^k\}$ cannot converge to the point $x.$ Suppose the contrary. Then for sufficiently large $k$ by virtue of (14.29) 
$\bar{Q}^k\in G_{\epsilon_k}(x^k)\subset G_{\epsilon^\prime}(x^\prime)$, and by virtue of (14.26), (14.27), $\|\bar{\Delta}^k\|\ge \gamma/4.$  Hence for sufficiently large $t$ and $s$, we obtain the relation
\[
2\epsilon^\prime\ge\|x^t-x^k\|=\left\|\sum_{k=s}^{t-1}\bar{\rho}_k\bar{P}^k \left/\sum_{k=s}^{t-1}\bar{\rho}_k \right.\right\|
\sum_{k=s}^{t-1}\bar{\rho}_k
\ge \frac{1}{4}\gamma\nu_K\lambda_K\sum_{k=s}^{t-1}{\rho}_k,
\]
where the right hand part, by virtue of the assumption made and condition (14.31), tends to infinity as $t\rightarrow +\infty$. The obtained contradiction proves that the whole sequence cannot converge to $x^\prime$. Therefore, there exists $\epsilon^{\prime\prime}<\epsilon^\prime$ such that indices
\[
m_s(\epsilon^{\prime\prime})=\sup\{m\ge k_s|\,\|x^k-x^\prime\|
\le\epsilon^{\prime\prime}, \mbox{ for } k\in[k_s,m)\}<+\infty.
\]
Consider the sequence $\{x^k_s=x^k, k\in[k_s,m_s(\epsilon^{\prime\prime})\}$. Let us show that Lemma 14.1 applies to them for sufficiently large $s$. Indeed, $x_s^{k_s}=x^{k_s}\rightarrow x^\prime,$ $0\notin G(x^\prime);$
\[
x_s^{k+1}=x_s^k-\bar{\rho}_k\bar{P}^k,\;\;\;
k\in[k_s,m_s(\epsilon^{\prime\prime}));
\]
\[
\bar{P}^k=\bar{Q}^k+\bar{\Delta}^k;\;\;\;
\|\bar{\Delta}^k\|\le\|\Delta^k\|/\lambda_K\rightarrow 0,\;\;\;
k\rightarrow\infty;
\]
\[
\bar{Q}^k\in G_{\epsilon_k}(x^k)\subset G_{\delta_s}(x^k);\;\;\;
\delta_s=\sup_{k\ge k_s}\epsilon_k\rightarrow 0,\;\;\;
s\rightarrow\infty;
\]
\[
\bar{\rho}_k=\rho_k\nu_k\sum_{r=r_k}^k\lambda_{kr}
\le N_K\Lambda_K\rho_k\rightarrow 0,\;\;\;
k\rightarrow\infty.
\]
For $s$ and $k\ge k_s$ such that $\|x^{k_s}-x^\prime\|\le \epsilon^{\prime\prime}/2$ and
$\|\bar{\Delta}^k\|\le \|{\Delta}^k\|/\lambda_K\le \Gamma$, it takes place
\[
\sigma_s=\sum_{k=k_s}^{m_s(\epsilon^{\prime\prime})-1}\bar{\rho}_k
\ge \|x^{m_s(\epsilon^{\prime\prime})}-x^{k_s}\|\left/
\max_{k_s\le k\le m_s(\epsilon^{\prime\prime})}\|\bar{P}^k\|\right.
=\frac{\epsilon^{\prime\prime}}{8\Gamma}=\sigma^\prime>0.
\]
Thus, Lemma 14.1 applies to the sequences 
$\{x^k_s=x^k, k\in[k_s,m_s(\epsilon^{\prime\prime})\}$
$(s=0,1,\ldots)$,
whence, given that $\lim_{s\rightarrow\infty}\delta_s=0$ and 
$\lim_{s\rightarrow\infty} \sup_{k\ge k_s} \rho_k = 0,$ 
the statements of this lemma follow. The lemma is proved.

The following lemma establishes boundedness conditions for sequences 
$\{x^k\}$ constructed by the method of averaged gradients.
\begin{lemma}
\label{lem:14.3}
Let $h(x)\rightarrow+\infty$ as $\|x\|\rightarrow+\infty$, and there exists $c>\max\{0,h(x^0)\}$ not beloning to
$H^*=\{h(x)|\,0\in G_h(x),\,h(x)>0\}.$
Let the sequence $\{x^k\}$ be built according to (14.22) -- (14.24)
under conditions (14.27) -- (14.31). Then for any $d>c$ there are
numbers $R$, $V$, and $\Delta$ such that if 
${\sup}_{k\geq 0}\,\rho_k\le R,\,$
${\sup}_{k\geq 0}\,\epsilon_k\le V,\,$
${\sup}_{k\geq 0}\,\|\Delta^k\|\le \Delta,\,$
then the sequence $\{x^k\}$ does not go out of the bounded set
$\{x|\,|h(x)\le d\}.$
\end{lemma}
{\it P r o o f.} Suppose the opposite. Then for some $d>c$ we find sequences 
$\{x^k\}$ starting for each $s$ from  point $x^0$ according to (14.22) - (14.31), where 
${\sup}_{k\geq 0}\,\rho_s^k\le R_s\rightarrow 0,\,$
${\sup}_{k\geq 0}\,\epsilon_s^k\le V_s\rightarrow 0,\,$
and ${\sup}_{k\geq 0}\,\|\Delta_s^k\|\le \Delta_s
\rightarrow 0,\,$ $(s\rightarrow\infty)$
, and $\{x^k\}$ going to infinity; the index $s$ denotes the number of the running sequence. Let us denote
\[
K=\{x|\,h(x)\le d\},\;\;\;\Gamma=\sup\{\|g\||\, g\in G(y), y\in K\}<+\infty.
\]
Let us select indices $k_s$ and $t_s$ such that
\begin{equation}
\label{14.32}
h(x_s^k)<d \mbox{ for } k\in[0,k_s], 
\end{equation}
\begin{equation}
\label{14.33}
h(x_s^{k_s})\le c<h(x^k)<d\le h(x_{t_s}), \;\;\; k\in(k_s,t_s]. 
\end{equation}
The sequence $\{x_s^{k_s}\}$ is bounded. Without changing the notation, we will assume that ${\lim}_{s\rightarrow\infty}x_s^{k_s}=x^\prime.$ Let us represent
\[
x_s^{k+1}=x_s^k-\rho_s^k\nu_s^k P_s^k=x_s^k-\bar{\rho}_s^k\bar{P}_s^k,
\;\;\;k\ge k_s,
\]
where
\[
\bar{\rho}_s^k=\rho_s^k\nu_s^k\sum_{r=r_k}^k\lambda_{kr},
\]
\[
\bar{P}_s^k=\bar{Q}_s^k+\bar{\Delta}_s^k=
P^k\left/\sum_{r=r_k}^k\lambda_{kr}(x_s^0,g_s^0,\ldots,x_s^k,g_s^k),\right.
\]
\[
\bar{Q}_s^k={Q}_s^k\left/\sum_{r=r_k}^k\lambda_{kr}\right.,\;\;\;
\bar{\Delta}_s^k={\Delta}_s^k\left/\sum_{r=r_k}^k\lambda_{kr}\right..
\]
Let us show that to the sequence $\{x_s^k, k_s\le k\le t_s\}$ 
$(s = 0,1,...)$ lemma 14.1 applies.

By virtue of (14.32), (14.33) and (14.27), (14.28), when $k < t_s$ the following holds
\begin{eqnarray}
\bar{\rho}_s^k\le \nu_K\Lambda_K\rho_s^k\le N_K\Lambda_KR_s\rightarrow 0,
\;\;\;s\rightarrow+\infty;\label{14.34}\\
\|\bar{\Delta}_s^k\|\le\|{\Delta}_s^k\|/\lambda_K\le\Delta_s/\lambda_K
\rightarrow 0,
\;\;\;s\rightarrow+\infty;\label{14.35}\\
\|\bar{P}_s^k\|\le\|\bar{Q}_s^k\|+\|\bar{\Delta}_s^k\|
\le\Gamma+\frac{1}{\lambda_K}\sup_{s\ge 0}\Delta_s
<+\infty;\label{14.36}\\
\underset{r_k\le r\le k}{\mathop{\max}}\|x_s^k-x_s^r\|\le\epsilon_s^k
\le V_s\rightarrow 0, \;\;\;\, 
s\rightarrow+\infty;\label{14.37}
\end{eqnarray}
\[
\bar{Q}_s^k\in co\{G(y)|\,\|y-x_s^k\|\le V_s\}.
\]
By virtue of (14.34), (14.36) 
$\underset{s\rightarrow\infty}{\mathop{\lim}}\|x^{k_s+1}-x^{k_s}\|=0$, whence, given (14.33), it follows
\[
h(x^\prime)=\underset{s\rightarrow\infty}{\mathop{\lim}}h(x_s^{k_s})
=\underset{s\rightarrow\infty}{\mathop{\lim}}h(x_s^{k_s+1})
=c\notin H^*.
\]
Thus $0\notin G_h(x^\prime)$, $h(x^\prime)>0$.

Let us find $\delta>0$ such that 
$\max \{h(y)|\,\|y - x\|\le\delta\}<d.$  By virtue of
(14.33) and the continuity of $h(y)$, the sequences 
$\{x_s^k, k_s\le k\le l_s\}$ $(s = 0,1,...)$ emerge from the $\delta$-neighborhood of $x^\prime$. Let us denote by $m_s$ the moments of the first such exit: $m_s\le l_s$. For all $k\in[k_s,m_s]$ it is satisfied  
(14.36); hence for sufficiently large $s$ it takes place
$\sum_{k=k_s}^{m_s-1}\bar{\rho}_s^k\ge\sigma>0$.

Thus, Lemma 14.1 applies to the sequences $\{x_s^k, k_s\le k\le m_s\}$. However, the minimizing property 2) of this lemma contradicts the construction (14.33). The lemma is proved.

The following technique is used for actual selection of the method parameters.

Let the estimates be known:
\[
R_0\ge {\sup}_{k\geq 0}\;\rho_k^0,\;\;\;
V_0\ge {\sup}_{k\geq 0}\;\epsilon_k^0,\;\;\;
\Delta_0\ge {\sup}_{k\geq 0}\;\|\Delta_k^0\|,
\]
which hold  for all sequences generated by the method.
Let us run the method from  point $x^0$ as long as $h(x^k) < d$. If at some moment $\tau_1$ it turns out that $h(x^{\tau_1})>d$, we return to the initial point, choose $\{\rho_k^1\}$, $\{\epsilon_k^1\}$,
$\{\Delta_1^k\}$ with smaller estimates:
\[
R_1\ge {\sup}_{k\geq 0}\;\rho_k^1,\;\;\;
V_1\ge {\sup}_{k\geq 0}\;\epsilon_k^1,\;\;\;
\Delta_1\ge {\sup}_{k\geq 0}\;\|\Delta_k^1\|,
\]
and run the method again from point $x^0$ with new parameters. If at some moment $\tau_2$ $h(x^{\tau_2})>d,$ then we again return to the initial point, choose $\{\rho_k^2\}$, $\{\epsilon_k^2\}$,
$\{\Delta_2^k\}$
with even smaller estimates:
\[
R_2\ge {\sup}_{k\geq 0}\;\rho_k^2,\;\;\;
V_2\ge {\sup}_{k\geq 0}\;\epsilon_k^2,\;\;\;
\Delta_2\ge {\sup}_{k\geq 0}\;\|\Delta_k^2\|,
\]
and run the method again from point $x^0$ with new parameters, etc.

If 
$\underset{s\rightarrow\infty}{\mathop{\sup}}R_s
=\underset{s\rightarrow\infty}{\mathop{\sup}}V_s
=\underset{s\rightarrow\infty}{\mathop{\sup}}\Delta_s=0$,
then, by virtue of Lemma 14.3, the method 
after a finite number of restarts will not leave the bounded
region $\{x | h(x)\le d\}$.
\begin{theorem}
\label{th:14.2}
Let in the problem (14.1), (14.2) the functions $f(x)$ and
$h(x)$ are generalized differentiable, with $h(x)\rightarrow +\infty$ at $\|x\|\rightarrow\infty$, and the set $H^* = \{h (x) | 0\in G_h(x), h(x) >0\}$ contains no intervals. Let the sequence $\{x^k\}$ constructed according to (14.22) - (14.31) is bounded. In this case 

either all limit points of $\{x^k\}$ do not belong to the admissible region $D = \{x | h (x)\le 0\}$; they belong to a connected subset of the set
\[
X_h^*=\{x|\,0\in G_h(x), h(x)>0\}
\]
and there exists a limit 
$\underset{s\rightarrow\infty}{\mathop{\lim}}h(x^k)>0;$

or all limit points $\{x^k\}$ belong to the admissible region $D$, then the minimal ones belong to the set $X^*=\{x|\,0\in G(x), h(x)\le 0\}$ and the segment 
$[\underset{s\rightarrow\infty}{\mathop{\underline{\lim}}}f(x^k),
\underset{s\rightarrow\infty}{\mathop{\overline{\lim}}}f(x^k)]$ 
is nested in the set $F^*=\{f(x) |\, x\in X^* \}.$ If thus $F^*$ contains no intervals, then all limit points $\{x^k\}$ belong to a connected subset of the set $X^*$ and there exists a limit 
$\underset{s\rightarrow\infty}{\mathop{\lim}}f(x^k).$
\end{theorem}

{\it P r o o f.} By the assumption of the theorem, the sequence $\{x^k\}$ belongs to some compact $K$. The minimal in $h$ limit points $\{x^k\}$ belong to $X_h^*\cup D$, for otherwise, due to property 2) of Lemma 14.2, they cannot be minimal. 

Denote 
$\underline{h}=\underset{s\rightarrow\infty}{\mathop{\underline{\lim}}}h(x^k)$,
$\overline{h}=\underset{s\rightarrow\infty}{\mathop{\overline{\lim}}}h(x^k)$,
$H_-^*=H^*\cup \{h\in \mathbb{E}_1|\,h\le 0\}$.
Let us show that $[\underline{h},\overline{h})\subset H_-^*$. Suppose the contrary. It was previously proved that $\underline{h}\in H_-^*$. Therefore there exists $c$ such that $\underline{h}<c<\overline{h}$ and 
$c\notin H_-^*$. Let us choose $d$ such that 
$\underline{h} < c < d < \overline{h}.$ 
The sequence $\{h(x^k)\}$ intersects the segment $[c, d]$ from $c$ to $d$ an infinite number of times; therefore there exist indices $k_s$, and $t_s$ such that 
$\underset{s\rightarrow\infty}{\mathop{{\lim}}}x^k=x^\prime$ 
and
\begin{equation}
\label{eqn:14.38}
h(x^{k_s})\le c<h(x^k)<d\le h(x^{t_s}),\;\;\; k\in(k_s.t_s).
\end{equation}
In view of (14.30), it takes place 
$h(x^\prime)=\underset{s\rightarrow\infty}{\mathop{{\lim}}}h(x^{k_s})
=\underset{s\rightarrow\infty}{\mathop{{\lim}}}h(x^{k_s+1})=c;$
hence $x^\prime\notin X_h^*\cup D$. Let us take $\epsilon^\prime$ such that
\[
\max\{h(y)|\,\|y-x^\prime\|\le\epsilon^\prime\}<d.
\]
Then for $\epsilon\le\epsilon^\prime$ the minimizing property 2) of Lemma 14.2 contradicts the construction (14.38). The obtained contradiction proves that 
$[\underline{h},\overline{h})\subset H_-^*.$ 
Since $H_-^*$ is closed, then $[\underline{h},\overline{h}]\subset H_-^*.$

If $H^*$ does not contain intervals, then it follows from the proof that either there exists a limit 
$\underset{s\rightarrow\infty}{\mathop{{\lim}}}h(x^{k})>0$ 
and then all limit points of $\{x^k\}$ do not belong to $D$, or 
$\overline{h}=\overline{\underset{s\rightarrow\infty}{\mathop{{\lim}}}}h(x^k)\le 0$ 
and then all limit points of $\{x^k\}$ belong to $D$. In the first case, all limit points of $\{x^k\}$ belong to $X_h^*$, for otherwise the divergence of $\{h(x^k)\}$ follows from property 2) of Lemma 14.2.

Let all limit points $\{x^k\}$ belong to $D$. Now we can assume that the function $f(x)$ on $D$ is minimized. The minimal in $f(x)$ limit points of $\{x^k\}$ belong to $X^*$, because otherwise the contradiction follows from property 3) of Lemma 14.2, namely, the existence of even smaller limit points of $\{x^k\}$. Just as it was done for $\{h(x^k)\}$, we can prove that the semi-interval 
$[\underset{s\rightarrow\infty}{\mathop{\underline{\lim}}}f(x^k),
\underset{s\rightarrow\infty}{\mathop{\overline{\lim}}}f(x^k))$ 
is embedded in $F^*$. By virtue of the boundedness of $D$ and the closedness of $G$, the set $F^*$ is closed; therefore 
$[\underset{s\rightarrow\infty}{\mathop{\underline{\lim}}}f(x^k),
\underset{s\rightarrow\infty}{\mathop{\overline{\lim}}}f(x^k)]
\subset F^*$. 

If $F^*$ does not contain intervals, then it follows from the proof that there exists 
$\underset{s\rightarrow\infty}{\mathop{{\lim}}}f(x^k)$. 
Then all limit points of $\{x^k\}$ belong to $X^*$, for otherwise we obtain a contradiction since property 3) of Lemma 14.2 implies the divergence of $\{f(x^k)\}.$

Due to the boundedness of the sequence $\{x^k\}$ and conditions (14.25)-(14.30), it is true that 
$\underset{s\rightarrow\infty}{\mathop{{\lim}}}\|x^{k+1}-x^k\|=0$.
Then, by Lemma 8.3, the set of limit points of the sequence $\{x^k\}$ is connected, i.e., $\{x^k\}$ converges to a connected subset of the set 
$X^*$. The proof of the theorem is completed.

\smallskip
{\bf 3. Adjustment of parameters in the method of averaged gradients.} 
\label{Sec.14.3}
Relations (14.25) - (14.31) represent general conditions of convergence of the method of averaged gradients. To be concretized are the parameters $\rho_k, \epsilon_k, r_k, \lambda_{kr}$. 
The values $\lambda_{kr}$ determine the method of averaging the previous gradients, $r_k$ and $\epsilon_k$ determine the number of considered (memorized) previous gradients, $\rho_k$ determine the step size along the chosen average direction. All these parameters must satisfy the general conditions (14.25)-(14.31).

Let us consider some adjustments of the step $\rho_k$.

Note that if $\rho_k$ satisfies conditions (14.30), (14.31), then 
$\rho_k^\prime=M_k\rho_k$ $(1/N\le M_k\le N<\infty)$ also satisfies (14.30), (14.31), i.e., one-dimensional minimization can be done in the direction of $P^k$, if, for example, the minimized function is decreasing.

{\it Step adjustment R1.} The simplest is the program step adjustment, when the sequence $\{\rho_k\}$ is set in advance such that
\[
\rho_k\ge 0,\;\;\;
\underset{s\rightarrow\infty}{\mathop{{\lim}}}\rho_k=0,\;\;\;
\sum_{k=0}^\infty\rho_k=\infty.
\]

{\it Step adjustment R2.}
Consider an adaptive (stepwise) step adjustment in which the indicator for the decreasing value of $\rho_k$ is the smallness of the distance 
$\rho(0,\{g^k,g^{k-1},\ldots,g^{s_k}\}$  
from  zero point to the convex hull of the accumulated generalized gradients. Let the numerical sequences $\{M_k\}$, $\{u_k\}$, $\{\theta_k\}$ are such that $u_k\ge 0$, $M_k>0$, $\theta_k>0$;
\[
\underset{k\rightarrow\infty}{\mathop{{\lim}}}M_k
=\underset{k\rightarrow\infty}{\mathop{{\lim}}}\theta_k=0,\;\;\;
\underset{k\ge 0}{\mathop{{\sup}}}\;u_k<\infty,\;\;\;
\sum_{k=0}^\infty u_k=\infty.
\]
Let in addition to conditions (14.25) - (14.31), 
$\underset{s\rightarrow\infty}{\mathop{{\lim}}}r_k=\infty.$
In this step adjustment, $\rho_k=M_{t_k}u_k$, and indices $t_k$ and $s_k$ change as follows:

$s_0=t_0=0$;

if $\rho(0,co\{g^k,g^{k-1},\ldots,g^{s_k}\})\ge \theta_{t_k}$,
then $s_{k+1}=s_k$, $t_{k+1}=t_k$;

if $\rho(0,co\{g^k,g^{k-1},\ldots,g^{s_k}\})< \theta_{t_k}$,
then $s_{k+1}=k$, $t_{k+1}=t_k+1$.

Let us prove convergence of the averaged gradient method with adjustable step rule R2. For this purpose it is enough to check conditions (14.30), (14.31). Conditions (14.25)-(14.29) are assumed to be satisfied.

Suppose that $x^k\rightarrow x$, $0\notin G(x)$. The sequences $\{t_k\}$ and, hence, $\{s_k\}$ cannot tend to infinity, for in such a case, when $k$ is sufficiently large, a contradiction is obtained:
\[
\frac{1}{2}\rho(0,G(x))\le
\rho(0,co\{g^{s_{k+1}},g^{s_{k+1}-1},\ldots,g^{s_k}\})< \theta_{t_k}
\rightarrow 0.
\]
Thus, for sufficiently large $k\ge s$, it holds $t_k = t$ and $s_k = s;$ therefore
\[
\sum_{k=0}^\infty\rho_k
=\sum_{k=0}^{s-1}\rho_k+M_t\sum_{k=s}^\infty u_k=\infty.
\]
Condition (14.31) is verified. In fact, it means the impossibility of $x^k \rightarrow x, 0 \notin G(x)$.

Let's check now the condition  (14.30). For simplicity, we assume $\Delta^k \equiv 0$. Assume that $\left\{x^k\right\}$ belongs to some compact $K$. Let us show that $t_k \rightarrow \infty$, from which it follows that $\lim _{k \rightarrow \infty} \rho_k=0$. Assume the opposite: $\lim _{k \rightarrow \infty} t_k=t<\infty$ and, therefore, $\lim _{k \rightarrow \infty} \mathrm{s}_k=s<\infty$. Then for $k \geq s$, $t_k=t$ and $s_k=s$; so
$$
\rho\left(0, \operatorname{co}\left\{g^k, g^{k-1}, \ldots, g^s\right\}\right) \geq \theta_t \text { and } \left\{x^k\right\} \subset K .
$$

Since $r_k \rightarrow \infty$, then for sufficiently large $k \geq k^{\prime}$ $r_k \geq s$. Then we get
$$
\left\|x^k-x^{k^{\prime}}\right\|=\left\|\sum_{r=k^{\prime}}^{k-1} \rho_r P^r\right \| \geq \lambda_K \theta_t M_s \sum_{r=k^{\prime}}^{k-1} u_r \rightarrow \infty,
$$
which contradicts the boundedness of $\left\{x^k\right\}$. The condition (14.30) has been checked. The convergence of the method of averaged gradients with the adjustment of the step R2 is proved.

If $P^k=\sum_{r=r_k}^k \mu_{k r} \bar{g}^r$, where
$$
\sum_{r=r_k}^k \mu_{k r}=1, \quad \mu_{k r} \geq 0, \quad \bar{g}^r=g^r /\left(\left\|g ^r\right\|+\varepsilon\right), \quad \varepsilon>0,
$$
then $\overline{g}^r$ can be used instead of $g^r$ in the regulation of $\mathrm{R} 2$.

To control the step $\rho_k$, it is possible to apply a regulation close to R2 from [90], the peculiarity of which is that the distance $\rho\left(0, co\left\{g^k, g^{k-1}, \ldots, g^s\right\}\right)$, firstly, is not calculated exactly, and, secondly, its approximate value is only corrected from iteration to iteration.

{\it Step adjustment R3.} This regulation is based on the following reasoning. Let first  $\rho_k=$ const, then $\sum_{k=0}^{\infty} \rho_k=+\infty$.
If the trajectory of the method falls into the neighborhood of the solution, then the distance alters around some value; therefore $\left\|x^k-x^s\right\| / \sum_{r=s}^{k-1} \rho_r \rightarrow 0$. This circumstance is an indication of reaching the neighborhood of the solution and signals the need to reduce the step.

Let $\left\{M_k\right\},\;\left\{\theta_k\right\},\;\left\{u_k\right\}$ be numerical sequences such that $M_k>0,\; \theta_k>0,\; u_k \geqslant 0$
$$
\lim _{k \rightarrow \infty} M_k=\lim _{k \rightarrow \infty} \theta_k=0, \quad \sup _{k \geqslant 0} u_k<\infty, \quad \sum_{k= 0}^{\infty} u_k=\infty .
$$

In this regulation, the step $\rho_k=M_{t_k} u_k$, and the indices $t_k$ and $s_k$ change as $k$ increases as follows:
$$
\begin{aligned}
& t_0=s_0=0 ; \\
& \text { if }\left\|x^k-x^s \right\| / \sum_{r=s_k}^{k-1} \rho_k \ge \theta_{t k} \text {, then } s_{k+1}=s_k,\; t_{k+1}=t_k \text {; } \\
& \text { if }\left\|x^k-x^s \right\| \sum_{r=s_k}^{k-1} \rho_r<\theta_{t k}, \text { then } s_{k+1}=k,\; t_{k+1}=t_k+1.
\end{aligned}
$$

{\it Step adjustment R4.} Numerous experiments show that the regulation of the step $\mathrm{R}3$ can be significantly improved if instead of the relation $\|x^k-x^{s_k}\|/\sum_{r=s_k}^{k-1} \rho_r$ to use the value
$$
R_k=\min _{s_k \le t<k}\left(\left\|x^k-x^t\right\| / \sum_{r=t}^{k-1} \rho_r\right) \text {. }
$$
Here it is necessary to remember the points of the trajectory $x^{s_k}, x^{s_k+1}, \ldots, x^k$. 

{\it Step adjustment R5.} In order not to remember the trajectory in the regulation $\mathrm{R} 4$, it is possible, without losing efficiency, to use  instead of $R_k$ the value
$$
T_k=\min \left(\left\|x^k-x^{s_k}\right\| / \sum_{r=s_k}^{k-1} \rho_r,\left\|x^k-x^{k_{ \min }}\right\| /\sum_{r=k_{\min }}^{k-1} \rho_r\right),
$$
where the indices $k_{\operatorname{mln}}$ are such that
$$
\begin{aligned}
& h\left(x^{k_{\mathrm{min}}}\right)=\min _{s_k \leqslant r<k} h\left(x^r\right), 
& \text { if } \min _{s_k \le r<k} h\left(x^r\right)>0 \text {, and } \\
& f\left(x^{k_{\mathrm{min}} }\right)=\min _{\left\{r \mid s_k \le r<k, h\left(x^ r\right) \le 0\right\}} f\left(x^r\right), 
& \text { if } \min _{s_k \le r<k} h\left(x^r\right) \le 0 \text {. } 
\end{aligned}
$$

To prove the convergence of the method of averaged gradients for step adjustments R3 - R5, it is enough to check conditions (14.30), (14.31). Let's check them for adjustment R4.

Suppose $\lim _{k \rightarrow \infty} x^k=x, 0 \notin G(x)$. Obviously, in this case
$\left\{x^k\right\}$ belongs to some compact $K$. For simplicity, we consider $\Delta^k\equiv 0$. Let's denote the sequence
$$
\left\{k_l\right\}_{l \ge 0}=\left\{k \mid R_k<\theta_{t_k}\right\}=\left\{k \mid s_{k+1}= k\right\} .
$$

Indices $t_k$ cannot tend to $+\infty$, because in this case $s_k \rightarrow \infty$ results in a contradiction: on the one hand,
$$
R_{k_l}=\min _{k_{l-1} \leqslant k<k_l}\left(\left\|x^{k_l}-x^k\right\| / \sum_{r=k}^ {k_l-1} \rho_k\right)<\theta_{t_l} \rightarrow 0, \quad l \rightarrow \infty,
$$
and on the other hand, for sufficiently large $k_l$
$$
R_{k_l} \geq \frac{1}{2} \lambda_K \rho(0, G(x))>0 .
$$

Therefore, for $k \geqslant k^{\prime}$, all $t_k=t<\infty$, so
$$
\sum_{k=0}^{\infty} \rho_k=\sum_{k=0}^{k^{\prime}-1} \rho_k+M_t \sum_{k=k^{\prime}}^ {\infty} u_k=+\infty .
$$

Condition (14.31) is verified. In fact, it means the impossibility of $x^k \rightarrow x, 0 \notin G(x)$.

Let's check condition (14.30). Let $\left\{x^k\right\}$ be bounded. We will show that $\lim _{k \rightarrow \infty} t_k=\infty$; from here it follows that $\lim _{k \rightarrow \infty} \rho_k=0$. Assume the opposite: $\lim _{k \rightarrow \infty} t_k=t<\infty$. It also follows that $\lim _{k \rightarrow \infty} s^k=s<\infty$. Then for $k \geq s$ holds
$$
\left\|x^k-x^s\right\| \geqslant \theta_t \sum_{r=s}^{k-1} \rho_r=\theta_t M_t \sum_{r=s}^{k-1} u_r \rightarrow \infty,
$$
which contradicts the boundedness of $\left\{x^k\right\}$. The condition (14.30) has been checked.

Let us consider some methods of averaging gradients, i.e. consider the adjustment of the averaging coefficients $\lambda_{kr}$ in the method (14.22)--(14.24).

Averaging the previous gradients is used to give the algorithm anti-aliasing properties. If the algorithm begins to dangle across the ravine of the minimized function, then the appropriate average of the previous gradients can indicate an acceptable direction along the ravine.

{\it Averaging procedure P1.} As the direction $P^k$ in  method (14.22) - (14.24), one can take the vector closest to the origin of the coordinates of the convex hull of the previous gradients $\left\{g^{r_k}, \ldots, g^k\right\} $ . Such $P^k$ is the direction of the fastest possible descent at point $x^k$ for the approximating function
$$
\bar{f}(x)=f\left(x^k\right)+\max _{r_k \leq r \leq k}\left<g^r, x-x^k\right>.
$$
\newpage
In this case, the construction of $P^k$ directions is reduced to solving a sequence of related problems of quadratic programming of a special kind. Algorithms from [33,170,186] can be used to solve these problems.

{\it Averaging procedure P2.} The direction $P^k$ in method (14.22) - (14.24) can be calculated recursively:
$$
P^0=g^0, \quad P^k=\left(1-\alpha_k\right) P^{k-1}+\alpha_k g^k, \quad 0 \leq \alpha_k \leq 1,
$$
moreover, periodically it is necessary to carry out the so-called recovering, that is, to set $ \alpha_k = 1 $. We denote by $r_k$ the nearest moment of recovery, i.e. $r_k=\max \left\{r \mid r \le k, \alpha_r=1\right\}$. Then $P^k=$ $=\sum_{r=r_k}^{k-1} \lambda_{k r} g^r$, where
$$
\begin{gathered}
\lambda_{k r_k}=\prod_{s=r_k}^k\left(1-\alpha_s\right), \quad \lambda_{k k}=\alpha_k, \\
\lambda_{k r}=\alpha_r \prod_{s=r+1}^k\left(1-\alpha_s\right) ; \quad \sum_{r=r_k}^k \lambda_{k r}=1 .
\end{gathered}
$$
Here $\alpha_k$ can be chosen, for example, from the condition
$$
\alpha_k=\mbox{argmin}_{0 \leq \alpha \leq 1}\left\|(1-\alpha) P^{k-1}+\alpha g^k\right\|,
$$
here $ \alpha_k $ is easily calculated analytically.

{\it Averaging procedure P3.} One can use the arithmetic average of the previous gradients:
$$
P^k=\frac{1}{k-r_k+1} \sum_{r=r_k}^k g^{r} .
$$
If the function $r_k$ increases monotonically, then
$$
P^k=\left(1-\frac{1}{k-r_k+1}\right) P^{k-1}+\frac{1}{k-r_k+1} g^k.
$$

In procedures P1 - P3, one can use normalized generalized gradients $\bar{g}^r$, it is especially useful in procedure P3.

The key problem in the method of averaged gradients is the regulation of the values $r_k$ and $\varepsilon_k$, which determine the number of previous gradients that are taken into account (averaged, remembered). It follows from numerical experiments that the number of $k-r_k+1$ gradients being averaged must be increased as $k$ grows; at the same time, they must be taken from smaller and smaller neighborhoods of points $x^k$.

For example, for step adjustments R2 - R4 and averaging procedure P3, the number $r_k$ can be adjusted as follows:
$$
\begin{array}{ll}
r_{k+1}=r_k, \quad \text { if } & \sum_{r=r_k}^k \rho_r \leq M_{t_k} ; \\
r_{k+1}=k, \quad \text { if } & \sum_{r=r_k}^k \rho_r>M_{t_k} .
\end{array}
$$
At the same time, we can consider, for example, $u_k=\left(k-r_k+1\right)^{1/2}$.

\smallskip
\textbf{4. The heavy ball (Polyak) method.} 
\label{Sec.14.4}
A special case of the method of averaged gradients is the well-known heavy ball method [100]. If a smooth function $f(x)$ is minimized, then the latter has the form:
\begin{equation}
\label{eqn:14.39}
\begin{gathered}
x^0 \in E_n, \quad x^1=x^0-\alpha_0 \nabla f\left(x^0\right), \\
x^{k+1}=x^k-\alpha_k \nabla f\left(x^k\right)+\beta_k\left(x^k-x^{k-1}\right), \quad k \ge 1,
\end{gathered}
\end{equation}
where $\nabla f(x)$ is the gradient of the function $f$ at the point $x$; $\alpha_k \geq 0, \beta_k \geq 0$ are step multipliers. The introduction of the inertial term $\beta_k\left(x^k-x^{k-1}\right)$ into the iterative process leads to acceleration of convergence, the motion trajectory becomes smooth and passes along the bottom of the ravine of the minimized function. The process (14.39) has a noticeable advantage over the gradient method  on smooth furrow functions [100].

The heavy ball method can be used for global minimization of functions, because with the help of an inertial term it can skip local minima.

Let's extend this method to the task of minimizing (nonconvex nonsmooth) generalized differentiable functions under constraints.

Let the problem (14.1), (14.2) be solved. For problems with more general constraints, the method is extended in the same way as it was done for the generalized gradient method of $\S$ 10.
Let
\begin{eqnarray}
x^0 \in E_n, & x^{k+1}=x^k-\rho_k P^k, \quad \rho_k \geqslant 0, \label{eqn:14.40}\\
P^0=g^0, & P^k=\left(1-\gamma_k\right) P^{k-1}+\gamma_k g^k,\label{eqn:14.41} \\
g^k \in G\left(x^k\right), & 0 \leq \gamma_k \leq 1, \quad k=0,1, \ldots,\label{eqn:14.42}
\end{eqnarray}
where $x^k$ is the current approximation, $P^k$ is the direction of movement, $\rho_k$ is the step multiplier, $\gamma_k$ is the averaging parameter; the mapping $G(x)$ is defined in (14.3). The following conditions are assumed:
\begin{equation}
\lim _{k \rightarrow \infty} \rho_k=0, \quad \sum_{k=0}^{\infty} \rho_k=+\infty, \quad \lim _{k \rightarrow \infty} \frac {\rho_k}{\gamma_k}=0.\label{eqn:14.43}
\end{equation}

Method (14.40) – (14.42) can be reduced to form (14.39):
$$
\begin{aligned}
x^{k+1}&=&x^k-\rho_k \gamma_k g^k-\rho_k\left(1-\gamma_k\right) P^{k-1}\;\;\;\;\;\;\;\;\;\; \;\;\;\;\; \;\; \\
& =&x^k-\rho_k \gamma_k g^k+\frac{\rho_k}{\rho_{k-1}}\left(1-\gamma_k\right)\left(x^k-x^{k-1} \right).
\end{aligned}
$$
It follows from relation (14.41)
$$
P^k=\sum_{r=0}^k \mu_r^k g^r ; \quad \mu_r^k \geq 0, \quad \sum_{r=0}^k \mu_r^k=1,
$$
where
$$
\mu_0^k=\prod_{i=1}^k\left(1-\gamma_i\right), \quad \mu_k^k=\gamma_k ; \quad \mu_r^k=\gamma_r \prod_{i=r+1}^k\left(1-\gamma_i\right).
$$
Consider $\left(0 \leq r_k \leq k\right)$
$$
P^k=\sum_{r=0}^k \mu_r^k g^r=\sum_{r=r_k}^k \mu_r^k g^r+\sum_{r=0}^{r_k-1} \mu_r ^k g^r.
$$
Let's define
\begin{eqnarray}
\bar{Q}^k&=&\sum_{r=r_k}^k \mu_r^k g^r / \sum_{r=r_k}^k \mu_r^k, \quad \bar{\rho}_k=\rho_k \sum_{r=r_k}^k \mu_r^k, \label{eqn:14.44}\\
\bar{\Delta}^k&=&\sum_{r=0}^{r_k-1} \mu_r^k g^r / \sum_{r=r_k}^k \mu_r^k .
\label{eqn:14.45}
\end{eqnarray}

Let's rewrite method (14.40) in the following form:
\begin{equation}
\label{eqn:14.46}
x^{k+1}=x^k-\rho_k P^k=x^k-\bar{\rho}_k\left(\bar{Q}^k+\bar{\Delta}^k\right) .
\end{equation}

The process (14.46) is represented in the form of methods of averaged gradients. The boundedness of the sequences (14.40), (14.46) can be ensured in the same way as it was done for the method of the generalized gradient $\S$ 9,10 and for the method of averaged gradients, that is, either by choosing sufficiently small steps $\rho_k$, or by using the mechanism of returning the starting point. Let us choose the indices $r_k$ (14.44), (14.45) so that the conditions (14.25) - (14.31) of the convergence of the method of averaged gradients are fulfilled.
\begin{lemma}
\label{lem:14.4}
Under conditions (14.43), there exist indices $r_k \leqslant k$ such that $\lim _{k \rightarrow \infty} r_k=+\infty$ and, at the same time,
\begin{equation}
\label{eqn:14.47}
\lim _{k \rightarrow \infty} \sum_{r=r_k}^k \rho_r=0, \quad \lim _{k \rightarrow \infty} \sum_{r=r_k}^k \gamma_k=+\infty.
\end{equation}
\end{lemma}
{\it P r o o f.} First, we construct a sequence of indices $\left\{r_k^{\prime}\right\}$ such that
$$
r_k^{\prime} \leq k, \quad \lim _{k \rightarrow \infty} r_k^{\prime}=+\infty, \quad \lim _{k \rightarrow \infty} \sum_{r= r_k^{\prime}}^k \rho_k=+\infty .
$$
This can be done, for example, as follows.

Let's put $\rho=\sup _{r \geq 0} \rho_r,$ $c_0=0,$ $r_0^{\prime}=0,$ $k=0$.

Let $c_k$ and $r_k^{\prime}$ already be calculated; now let us calculate $c_{k+1}$ and $r_{k+1}^{\prime}$.

If $\sum_{r=r_k^\prime}^k \rho_r<c_k$, then we set $c_{k+1}=c_k, r_{k+1}^{\prime}=r_k^{\prime}$, otherwise $c_{k+1}=c_k+2 \rho,$ $r_{k+1}^{\prime}=r_k^{\prime}+1$.

It is clear that $\lim _{k \rightarrow \infty} r_k=+\infty,$ $\lim _{k \rightarrow \infty} c_k=+\infty$, and
$$
\sum_{r = r_k^\prime}^k \rho_r \geq c_k-\rho \rightarrow+\infty \quad \text { at } \quad k \rightarrow+\infty.
$$

We denote $\delta_k=\sup _{r \geq r_k^\prime} (\rho_r / \gamma_r) \rightarrow 0,$ $ k \rightarrow+\infty$. Now we will find the indices $r_k$ from the condition
$$
\sum_{r=r_k}^k \gamma_r \le \min \left(\frac{1}{\sqrt{\delta_k}}, \;\sum_{r=r_k^{\prime}}^k \gamma_r\right)<\sum_{r=r_k-1}^k \gamma_r .
$$
Here  
$$\sum_{r=r_k}^k \gamma_r \leq \sum_{r=r_k^{\prime}}^k \gamma_r ;\text { \;\;\;  so\;\; } r_k \ge r_k^\prime.
$$
For sufficiently large  $k$
$$\sum_{r=r_k^{\prime}}^k \gamma_r \ge \sum_{r=r_k^{\prime}}^k \rho_r \rightarrow+\infty, \quad k \rightarrow+\infty, 
$$
So 
$$\sum_{r=r_k}^k \gamma_r \ge \min \left(\frac{1}{\sqrt{\delta_k}}, \sum_{r=r_k^{\prime}}^k \gamma_r\right)-\gamma_{r_k-1} \rightarrow+\infty, \quad k \rightarrow \infty . 
$$
At the same time, we have
$$
\sum_{r=r_k}^k \rho_r \leq \sup _{r_k^{\prime} \leq r \leq k} \frac{\rho_r}{\gamma_r} \sum_{r=r_k}^k \gamma_r \leq \sqrt{\delta_k} \rightarrow 0, \quad k \rightarrow+\infty .
$$
The lemma is proved.

Now, for process (14.44)-(14.46), we will check the conditions (14.25) – (14.31) for the convergence of the method of averaged gradients.
Note that
$$
\sum_{r=r_k}^k \mu_r^k=\sum_{r=r_k}^{k-1} \alpha_r \prod_{i=r+1}^k\left(1-\gamma_i\right)+\gamma_k=1-\prod_{i=r_k}^k\left(1-\gamma_i\right), 
$$
$$
\sum_{r=0}^{r_k-1} \mu_r^k=\prod_{l=r_k}^k\left(1-\gamma_i\right), 
$$
$$
\prod_{i=r_k}^k\left(1-\gamma_i\right)=\exp \left\{\sum_{i=r_k}^k \ln \left(1-\gamma_i\right)\right\} \leq \exp \left\{-\sum_{i=r_k}^k \gamma_i\right\} .
$$
Under conditions (14.43)-(14.45), (14.47) and due to boundedness of $\left\{x^k\right\}$, it holds
$$
\lim _{k \rightarrow \infty} \bar{\rho}_k \leq \lim _{k \rightarrow \infty} \rho_k=0; 
$$
$$
\sum_{k=0}^{\infty} \bar{\rho}_k \geqslant \sum_{k=0}^{k^{\prime}} \bar{\rho}_k+\frac{1}{ 2} \sum_{k=k^{\prime}+1}^{\infty} \rho_k=+\infty \;\;\;
\left( \sum_{r=r_k}^{k} \mu_r^k \ge \frac{1}{2} \mbox{  for  } k\geq k^\prime\right); 
$$
\begin{eqnarray}
\max _{r k \leq s \leq k}\left\|x^k-x^s\right\|&=&\sup _{r_k \leq s \leq k}\left\|\sum_{r=s}^ {k-1} \rho_r \sum_{t=0}^r \mu_t^r g^t\right\| \nonumber \\
&\leq&\left(\sup _{0 \leq t \leq k}\left\|g^t\right\|\right) \sum_{r=r_k}^k \rho_r \rightarrow 0, \quad k \rightarrow \infty; \nonumber
\end{eqnarray}
\begin{eqnarray}
\left\|\bar{\Delta}^k\right\| &\leq&\left(\sup _{0 \leq r<r_k}\left\|g^r\right\|\right) \sum_{r=0}^{r_k-1} \mu_r^k\left( 1-\sum_{r=0}^{r_k-1} \mu_r^k\right) \nonumber \\
&=&\left(\sup _{0 \leq r \leq r_k}\left\|g^r\right\|\right) \prod_{i=r_k}^k\left(1-\gamma_i\right) / \left(1-\prod_{i=r_k}^k\left(1-\gamma_i\right)\right) \rightarrow 0, \quad k \rightarrow \infty .\nonumber
\end{eqnarray}

Thus, under the assumption of boundedness of the sequence $\left\{x^k\right\}$ under conditions (14.43), all conditions (14.27) - (14.31) for the convergence of method (14.46) and, consequently, method (14.40) - (14.42) are satisfied. Therefore, the assertions of Theorem 14.2 are valid for the heavy ball method. The convergence of process (14.46) also follows from the convergence of method (14.4) - (14.7).

As with the averaged gradient method, normalized generalized gradients can be used in procedure (14:40) - (14:42) (see Remark 14.2).

We will describe the mechanism of returning to the starting point, which guarantees the boundedness of the heavy ball method.

Let the sequences $\left\{\rho_k\right\}_{k \geq 0},\left\{\gamma_k\right\}_{k\geq 0}$ satisfy conditions (14.42), (14.43). Let's run the heavy ball method from point $x^0$ with parameters $\left\{\rho_k\right\}_{k \geq 0},\left\{\gamma_k\right\}_{k \geq 0}$. If at some moment $\tau_1$ $h\left(x^{\tau_1}\right)>d$ happens, then we return to the starting point and redefine $x^{\tau_1}=x^0$. Let's run the method again from the point $x^{\tau_1}=x^0$ from the moment $\tau_1$ with the parameters $\left\{\rho_k\right\}_{k \geq \tau_1},\left\{\gamma_k\right\}_{k>\tau_1}$. If at some moment $\tau_2$ it turns out that $h\left(x^{\tau_2}\right)>d$, then we return to the point $x^0$ and redefine $x^{\tau_2}=x^0$. Let's run the method again from the point $x^{\tau_2}=x^0$ from the moment $\tau_2$ with parameters $\left\{\rho_k\right\}_{k \geqslant \tau_1},\left\{\gamma_k\right\}_{k>\tau_2}$. And so on.

Denote $\left\{\gamma_k^{\prime}\right\}=\left\{1, \gamma_1, \ldots, \gamma_{\tau_1-1}, 1, \gamma_{\tau_1+1}, \ldots, \gamma_{\tau_2-1}, 1, \gamma_{\tau_2+1}, \ldots\right\}$.

Obviously, the sequences $\left\{\rho_k\right\},\left\{\gamma_k^{\prime}\right\}$ satisfy condition (14.43); therefore, there are indices $r_k \rightarrow+\infty$ that satisfy $(14.47)$ (where $\gamma_k^{\prime}$ stands instead of $\gamma_k$).

Formally, it can be assumed that the points $\left\{x^{\tau_s}, x^{\tau_s+1}, \ldots, x^{\tau_{s+1}-1}\right\}$ are generated by the method of averaged gradients of the form $(14.44)-(14.46)$, in which only gradients $g^r$ with $r \geqslant \tau_s$ are actually taken into account, because $\gamma_{\tau_s}^{\prime}=1$. It is obvious that $R_s=\sup _{k \geqslant \tau_s} \bar{\rho}_k \rightarrow 0,$ $\Delta_s=\sup _{k \geqslant \tau_s}\left\|\bar{\Delta} ^k\right\| \rightarrow 0$ if $\tau_s \rightarrow \infty$.

Let us denote $\Gamma=\sup \{\|g\| \mid g \in G(y), h(y) \leqslant d\}$. Let's define indices
$\left\{r_k^{\prime}\right\}$ such that $r_k^{\prime}=r_k$, if $\max \left\{\tau_s \mid \tau_s \leqslant k\right\} \leqslant r_k$, and $r_k^{\prime}=\max \left\{\tau_s \mid r_k \leqslant\right.$
$\left.\leqslant \tau_s \leqslant k\right\}$, otherwise. Let $\varepsilon_k^{\prime}=\max _{r_k^{\prime} \leqslant r \leqslant k}\left\|x^k-x^r\right\|$.

For averaged gradients $\bar{Q}^k$ in this case $\bar{Q}^k \in \operatorname{co}\left\{G(y)\|\| y-x^k \| \leqslant \varepsilon_k^{\prime}\right\}$. 
The estimate is true,
$$
\varepsilon_k^{\prime}=\max _{r_k^{\prime} \leqslant r \leqslant k}\left\|x^k-x^r\right\| \leqslant \Gamma \sum_{r=r_k^{\prime}}^k \rho_r \leqslant \Gamma \sum_{r=r_k}^k \rho_r \rightarrow 0, \quad k \rightarrow \infty .
$$
Therefore, $V_s=\sup _{k \geq \tau_s} \varepsilon_k^{\prime} \rightarrow 0$ if $\tau_s \rightarrow \infty$.

Therefore, by virtue of Lemma 14.3, the heavy ball method leaves the bounded region $\{x \mid h(x) \leq d\}$ only a finite number of times.

\smallskip
\textbf{5. The gully step (Nesterov) method.} 
\label{Sec.14.5}
As a special case of the method of averaged gradients considered above, is the well-known method of the gully step. It was proposed in [14], then significantly generalized and investigated in $[74,76]$. In the smooth convex case, the method has a high rate of convergence (of the order of $O\left(k^{-2}\right)$, where $k$ is the number of iterations) and is optimal for tasks of this class.

The gully step method can be used globally for global minimization of functions. It passes over the ravines of the minimized function and has a certain inertia (see (14.51)), which allows skipping local minima. Therefore, with its help, by adjusting the step multipliers, you can roughly review the valleys of the function in order to detect areas of local minima. Then local descents can be made from the found promising points.

We will extend this method to the problem of minimizing (non-convex non-smooth) generalized differentiable functions under constraints. Let's write down the method of the rigging step for solving the problem (14.1), (14.2). On problems with more general constraints, it is extended in the same way as it is done for the generalized gradient method $\S$ 10 . The method produces two sequences $\left\{x^k\right\}$ and $\left\{y^k\right\}$ according to the following rules:
\begin{eqnarray}
&&x^0=y^0=\bar{x} \in E_n, \label{eqn:14.48}\\
&&y^{k+1}=x^k-\rho_k g^k, \quad g^k \in G\left(x^k\right), \label{eqn:14.49}\\
&&x^{k+1}=y^{k+1}+\lambda_k\left(y^{k+1}-y^k\right), \quad k=0,1, \ldots
\label{eqn:14.50}
\end{eqnarray}
Here, from the point $x^k$, a gradient descent step is made with a step multiplier $\rho_k\geq 0$, and along the direction $(y^{k+1}-y^k)$, a ramping step is made with a step multiplier $\lambda_k \geq 0$. The sets $G(x)$ are defined in (14.3). Normalized generalized gradients can be used.

Note that in the method $(14.48)-(14.50)$ the gradients of the task functions are not calculated at the points $y^k$; therefore, we will present the method in a form that does not contain auxiliary points $y^k$. Excluding $y^k$ and $y^{k+1}$ from (14.50) using (14.49), we obtain
\begin{equation}\label{eqn:14.51}
\begin{aligned}
& x^0 \in E_n, \quad x^1=x^0-\left(1+\lambda_0\right) \rho_0 g^0 \\
& x^{k+1}=x^k-\left(1+\lambda_k\right) \rho_k g^k+\lambda_k \rho_{k-1} g^{k-1}+\lambda_k\left( x^k-x^{k-1}\right) \\
& g^k \in G\left(x^k\right), \quad g^{k-1} \in G\left(x^{k-1}\right), \quad k=0,1 , \ldots
\end{aligned}
\end{equation}

In the right-hand side of formula (14.51) there is an inertial term $\lambda_k\left(x^k-x^{k-1}\right)$, which allows the sequence $\left\{x^k\right\}$ generated by the method skip local extremum.

From these relations it follows,
$$
\begin{aligned}
& x^{k+1}=x^k-\rho_k g^k-\sum_{r=0}^k \lambda_r \lambda_{r+1} \ldots \lambda_k \rho_r g^r, \\
& x^{k+1}-x^s=-\left(1+\lambda_s+\lambda_s \lambda_{s+1}+\ldots+\lambda_s \lambda_{s+1} \ldots \lambda_k\right) \rho_s g^s- \\
& -\left(1+\lambda_{s+1}+\lambda_{s+1} \lambda_{s+2}+\ldots+\lambda_{s+1} \ldots \lambda_k\right) \rho_{s+1} g^{s+1}-\ldots \\
& \ldots-\left(1+\lambda_k\right) \rho_k g^k+\left(\lambda_s+\lambda_s \lambda_{s+1}+\ldots+\lambda_s \lambda_{s+1} \ldots \lambda_k\right)\left(x^s-x^{s-1}\right)+ \\
& +\left(\lambda_s+\lambda_s \lambda_{s+1}+\ldots \lambda_s \lambda_{s-1} \ldots \lambda_k\right) \rho_{s-1} g^{s+1}.
\end{aligned}
$$
Let's introduce the notation:
$$
\begin{gathered}
\lambda_{ji}=\lambda_i \lambda_{i+1} \ldots \lambda_j \\
\sigma_{q p}=\lambda_p+\lambda_p \lambda_{p+1}+\ldots+\lambda_p \lambda_{p+1} \ldots \lambda_q .
\end{gathered}
$$
In these notations, we have
$$
\begin{gathered}
x^{k+1}=x^{k}-\rho_{k} g^{k}-\sum_{r=0}^{k} \lambda_{k r} \rho_{r} g^{r}, \\
x^{k+1}-x^{s}=-\sum_{r=s}^{k}\left(1+\sigma_{k r}\right) \rho_{r} g^{r}+\sigma_{k s}\left(x^{s}-x^{s-1}\right)+\sigma_{k s} \rho_{s-1} g^{s-1}.
\end{gathered}
$$

First, let's consider the gully step method with the so-called recovery. We assume that at certain points (recovery points) $k_{s}$, the gully step $\lambda_{k_{s}}$ is equal to zero. For convenience, we set $k_{0}=0$. Then, we have:
\begin{eqnarray}
x^{k+1}=x^{k}-\boldsymbol{\rho}_{k} g^{k}-\sum_{r=k_{s}}^{k} \lambda_{k r} \rho_{r} g^{r}=x^{k}-\bar{\rho}_{k} \bar{P}_{k}, \quad k_{s} \le k<k_{s+1},\label{eqn:14.52}\\
\bar{\rho}_{k}=\rho_{k}+\sum_{r=k_{s}}^{k} \lambda_{k r} \rho_{r}, \quad \bar{P}^{k}=\left(\rho_{k} g^{k}+\sum_{r=k_{s}}^{k} \lambda_{k r} \rho_{r} g^{r}\right) / \bar{\rho}_{k} .\nonumber
\end{eqnarray}

At $k=k_{s}$, we have $x^{k_{s}+1}=x^{k_{s}}-\rho_{k_{s}} g^{k_{s}},\quad
\bar{\rho}_{k_s}=\rho_{k_{s}}, \quad \bar{P}^{k_{s}}=g^{k_{s}}$. Now it is clear that this method belongs to the class of averaged gradient methods (with an appropriate choice of parameters). The boundedness of the method can be ensured in the same way as it is done for the averaged gradient method (14.23): return the process to the initial point if it goes too high on the function $h(x)$, and perform a restoration at this moment.

Let $\rho_k$ satisfy conditions (14.30) and (14.31). Let $r_k=\max_{s\geq 0}\{k_s|k_s\leq k\}$. Assume that the following conditions are satisfied:
\begin{equation}
\max_{r \in\left[r_{k}, k\right]} \sigma_{k r} \le C<+\infty, \quad \lim_{k \rightarrow \infty} \sum_{r=r_{k}}^{k} \rho_{r}=0, \label{eqn:14.53}
\end{equation}
which is the case, for example, when $\lambda_{k} \leq \lambda<+\infty$ and $k-r_{k} \leq const$.

Since $\lambda_{k r} \leq \sigma_{k r} \leq C$ for $r \in\left[r_{k}, k\right]$, then
$$
\sum_{r=r_{k}}^{k} \lambda_{k r} \rho_{r} \leq \sum_{r=r_{k}}^{k} \sigma_{k r} \rho_{r} \leq C \sum_{r=r_{k}}^{k} \rho_{r} \rightarrow 0, \quad k \rightarrow \infty .
$$
From here, taking into account (14.53), we obtain
$$
\lim _{k \rightarrow \infty}\left(\rho_{k}+\sum_{r=r_{k}}^{k} \lambda_{k r} \rho_{r}\right)=\lim _{k \rightarrow \infty} \sum_{r=r_{k}}^{k}\left(1+\sigma_{k r}\right) \rho_{r}=0.
$$
And thus, all convergence conditions (14.25) -- (14.31) of the method of averaged gradients (14.52) are satisfied.

Let's consider the gully step method (14.48) -- (14.50) without recovering. In the transformed form, it has the following form:
\begin{eqnarray}
x^{k+1}&=&x^{k}-\left(1+\lambda_{k}\right) \rho_{k} g^{k}+\lambda_{k} \rho_{k-1} g^{k-1}+\lambda_{k}\left(x^{k}-x^{k-1}\right)\nonumber \\
& =&x^{k}-\rho_{k} g^{k}-\sum_{r=0}^{k} \lambda_{k r} \rho_{r} g^{r}.\nonumber
\end{eqnarray}

For simplicity, let's consider programmatic parameter adjustments. Suppose
\begin{equation}
0 \leq \lambda_{k} \leq \lambda<1, \quad \rho_{k} \geq 0 . \quad \lim_{k \rightarrow \infty} \rho_{k}=0, \quad \sum_{k=0}^{\infty} \rho_{k}=+\infty .\label{eqn:14.54}
\end{equation}
Note that for the convex smooth case, optimal multipliers $\lambda_{k}$ have been found in [74]; they have the properties $\lambda_{k} \leq 1$ and $\lim_{k \rightarrow \infty} \lambda_{k}=1$.

Let us represent
$$
\sum_{r=0}^{k} \lambda_{k r} \rho_{r} g^{r}=\sum_{r=r_{k}}^{k} \lambda_{k r} \rho_{r} g^{r}+\lambda_{k r_{k}} \sum_{r=0}^{r_{k}-1} \lambda_{\left(r_{k}-1\right) r} \rho_{r} g^{r}.
$$

Let us assume that there exist indices $r_{k} \leq k$ such that $\lim_ {k \rightarrow \infty} r_{k}=+\infty$,
\begin{equation}
\lim_ {k \rightarrow \infty} \sum_{r=r_{k}}^{k} \rho_{r}=0, \quad \lim_ {k \rightarrow \infty} \sum_{r=r_{k}}^{k} \lambda^{r_{k}-r} \rho_{r}=+\infty .
\label{eqn:14.55}
\end{equation}

We will show that this assumption holds, for example, if $\rho_{k}=c / k^{\alpha}$, where $0<\alpha<1$. We have
$$
\sum_{r=r_{k}}^{k} \rho_{r} \leq \frac{c}{1-\alpha}\left[(k+1)^{1-\alpha}-r_{k}^{1-\alpha}\right] \leq \frac{c}{1-\alpha} \frac{k+1-r_{k}}{k^{\alpha}} .
$$
For monotonically decreasing $\left\{\rho_{k}\right\}$,  using summation by parts, we obtain  $(q=1 / \lambda>1)$
$$
\sum_{r=r_{k}}^{k} \lambda^{r_{k}-r} \rho_{r}=\sum_{r=r_{k}}^{k} q^{r-r_{k}} \rho_{r} \geq \frac{1}{q-1}\left(q^{k+1-r_{k}} \rho_{k}-\rho_{r_k}\right) .
$$
Let us choose $r_k$ in such a way that $r_k \rightarrow +\infty$ and $\sum_{r=r_{k}}^{k} \rho_{r}$ decreases sufficiently slowly as $k \rightarrow \infty$, for example, we take $r_{k} \sim k+1-k^{\beta}, 0<\beta<\alpha<1$. Then the conditions (14.55) are guaranteed to hold.

Let us denote
$$
\begin{gathered}
\bar{\rho}_{k}=\rho_{k}+\sum_{r=r_{k}}^{k} \lambda_{k r} \rho_{r}, \quad \bar{Q}^{k}=\left(\rho_{k} g^{k}+\sum_{r=r_{k}}^{k} \lambda_{k r} \rho_{r} g^{r}\right) / \bar{\rho}_{k}, \\
\bar{\Delta}^{k}=\frac{\lambda_{k r_{k}}}{\bar{\rho}_{k}} \sum_{r=0}^{r_{k}-1} \lambda_{\left(r_{k}-1\right) r} \rho_{r} g^{r} .
\end{gathered}
$$
Let us represent
\begin{equation}
x^{k+1}=x^{k}-\rho_{k}\left(\bar{Q}^{k}+\bar{\Delta}^{k}\right) \label{eqn:14.56}
\end{equation}

It is clear that this method falls into the framework of the method of averaged gradients and also the method (14.4), (14.5). Let us check the convergence conditions (14.28) - (14.33). Due to (14.54), (14.55), the following relationships hold for bounded sequences:
\begin{eqnarray}
& \lim _{k \rightarrow \infty} \bar{\rho}_{k} \le \lim _{k \rightarrow \infty}\left(\rho_{k}+\sum_{r=r_{k}}^{k} \rho_{r}\right)=0, \label{eqn:14.57} \\
& \sum_{k=0}^{\infty} \bar{\rho}_{k} \ge \sum_{k=0}^{\infty} \rho_{k}=+\infty, \label{eqn:14.58}\\
& \left\|\sum_{r=0}^{r_{k}-1} \lambda_{\left(r_{k}-1\right) r} \rho_{r} g^{r}\right\| \le \frac{1}{1-\lambda} \sup _{r \ge 0} \rho_{r} \sup _{r \ge 0}\left\|g^{r}\right\|, \label{eqn:14.59} \\
& \lim _{k \rightarrow \infty} \bar{\Delta}^{k}=0, \;\;\;\lim _{r \rightarrow \infty}\left\|x^{r}-x^{r-1}\right\|=0,\label{eqn:14.60}
\end{eqnarray}
\begin{eqnarray}
&\max _{r_{k} \le r \le k}\left\|x^{k}-x^{r}\right\|\leq \nonumber\\
& \le \sum_{r=r_{k}}^{k}(1\left.+\sigma_{k r}\right) \rho_{r}\left\|g^{r}\right\|+\sigma_{k_{r k}}\left\|x^{r_{k}}-x^{r_{k}-1}\right\|+\sigma_{k r_{k}} \rho_{r_{k}-1}\left\|g^{r_{k}-1}\right\| \nonumber \\
& \le \frac{1}{1-\lambda} \sup _{r \ge 0}\left\|g^{r}\right\| \sum_{r=r_{k}}^{k} \rho_{r}+\frac{1}{1-\lambda}\left\|x^{r_{k}}-x^{r_{k}-1}\right\|\nonumber \\
&+\frac{\lambda}{1-\lambda} \rho_{r_{k}-1}\left\|g^{r_{k}-1}\right\| \rightarrow 0, \;\;\;k \rightarrow \infty. \label{eqn:14.61}
\end{eqnarray}

Thus, for the method (14.56) under the assumption of boundedness of $\left\{x^{k}\right\}$ all conditions (14.7) and (14.25) -- (14.31) are satisfied; therefore, the statements of Theorems 14.1 and 14.2 are valid for the method of gully steps. The same holds if the normalized generalized gradients are used in the method of ridge steps (14.48) -- (14.50).

Let us describe in more detail the mechanism of returning to the starting point, which guarantees the boundedness of the method of gully steps.

Let the sequences $\left\{\rho_{k}\right\}_{k \geq 0},\left\{\lambda_{k}\right\}_{k>0}$ satisfy the conditions (14.54). We run the gully step method from the point $x^{0}$ until $h\left(x^{k}\right) \leq d$. If at some moment 
$h\left(x^{\tau_{1}}\right)>d$, we return to the initial point and redefine $x^{\tau_{1}}=x^{0}$. Then we run the method again from the point $x^{\tau_{1}}=x^{0}$ from the moment $\tau_{1}$ with the parameters $\left\{\rho_{k}\right\}_{k \geq \tau_{1}},\left\{\lambda_{k}\right\}_{k>\tau_{1}}$. If at some moment $\tau_{2}$ we have $h\left(x^{\tau_{2}}\right)>d$, we return to the initial point and redefine $x^{\tau_{2}}=x^{0}$. Then we run the method again from the point $x^{\tau_{2}}=x_{0}$ from the moment $\tau_{2}$ with the parameters $\left\{\rho_{k}\right\}_{k} \geq \tau_{2},\left\{\lambda_{k}\right\}_{k>\tau_{2}}$. And so on.

Let us denote $\left\{\lambda_{k}^{\prime}=0, \lambda_{1}, \ldots, \lambda_{\tau_{1}-1}, 0, \lambda_{\tau_{1}+1}, \ldots, \lambda_{\tau_{2}-1}, 0, \ldots\right\}$.

Obviously, the sequences $\{\rho_k\}_{k\geq 0}, \{\lambda_k'\}_{k\geq 0}$ satisfy the conditions (14.54), (14.55), where $\lambda_k$ is replaced by $\lambda_k^\prime$ and the indices $r_k \rightarrow +\infty$. Formally, one can assume that the points ${x^{\tau_s}, x^{\tau_s+1}, \ldots, x^{\tau_{s+1}-1}}$ are generated by the method of averaged gradients of the form (14.56), in which only the gradients $g'$ with $r\geq \tau_s$ are actually taken into account, since $\lambda_{\tau_s}'=0$. Obviously, $R_s=\sup_{k\geq \tau_s} \bar{\rho}_k \rightarrow 0$, if $\tau_s\rightarrow \infty$. By (14.60), $\Delta_s=\sup_{k\geq \tau_s} |\Delta^k| \rightarrow 0$, if $\tau_s\rightarrow \infty$.

Let $\Gamma=\sup \{\|g\| | g \in G(y),\; h(y) \leq d\}$. Define the indices $\left\{r_{k}^{\prime}\right\}$ such that 
$r_{k}^{\prime}=r_{k}$ if 
$\max \left\{\tau_{s} \mid \tau_{s} \leq k\right\} \leq r_{k}$, and 
$r_{k}^{\prime}=\max \left\{\tau_{s} \mid r_{k} \leq \tau_{s} \leq k\right\}$, 
otherwise. Let $\varepsilon_{k}^{\prime}=\underset{r_{k}^{\prime} \leq r \leq k}\max \left|x^{k}-x^{r}\right|$. For the averaged gradients $\bar{Q}^{k}$, it holds that 
$\bar{Q}^{k} \in \operatorname{co} \{G(y) \mid | \|y -x^{k} \| \leqslant \varepsilon_{k}^{\prime}\}$. 
It is obvious that relations (14.57) -- (14.60) are satisfied, as well as the following estimates, similar to (14.61):

if $r_{k}^{\prime}=\tau_{s}$, then
$$
\varepsilon_{k}^{\prime} \leqslant \sum_{r=\tau_{s}}^{k}\left(1+\sigma_{k r}\right) \rho_{r}\left\|g^{r}\right\| \leqslant \frac{\Gamma}{1-\lambda} \sum_{r=\tau_{s}}^{k} \rho_{r} \leqslant \frac{\Gamma}{1-\lambda} \sum_{r=r_{k}}^{k} \rho_{r} \rightarrow 0 ;
$$

if $r_{k}^{\prime}>\tau_{s}$, then
$$
\varepsilon_{k}^{\prime} \leqslant \frac{\Gamma}{1-\lambda} \sum_{r=r_{k}}^{k} \rho_{r}+\sigma_{k r_{k}}\left\|x^{r_{k}}-x^{r_{k}-1}\right\|+\sigma_{k r_{k}} \rho_{r_{k}-1}\left\|g^{r_{k}-1}\right\| \rightarrow 0 .
$$
Furthermore, for $k>\tau_{s}$ we have:
$$
\left\|x^{k}-x^{k-1}\right\| \leqslant \bar{\rho}_{k}\left(\left\|\bar{Q}^{k}\right\|+\left\|\bar{\Delta}^{k}\right\|\right) \leqslant\left(\Gamma+\sup _{k \geq 0}\left\|\bar{\Delta}^{k}\right\|\right) \sup _{k \geq \tau_{s}} \bar{\rho}_{k} .
$$

From here it follows that $V_{s}=\underset{k \geq \tau_{s}}\sup\ \varepsilon_{k}^{\prime} \rightarrow 0$, if $\tau_{s} \rightarrow+\infty$. Therefore, due to Lemma 14.3, the gully step method only exits the bounded region $\{x \mid h(x) \leq d\}$ a finite number of times.

In conclusion of this section, we note that in constructing methods of averaged gradients, we relied on the stability (Lemma 14.1 and Theorem 14.1) of the original generalized gradient method. Similarly, having any stable base method of the form $x^{k+1} \in x^{k}-\rho_{k} G\left(x^{k}\right)$, where $G(x)$ is no longer a gradient mapping, one can construct methods of averaged directions and analogues of heavy ball and gully step methods.

\newpage
\begin{flushright}
\textbf{CHAPTER 5}
\label{Ch.5}

\textbf{SOLUTION OF EXTREMAL PROBLEMS WITH \\LIPSCHITZ FUNCTIONS UNDER CONSTRAINTS}
\end{flushright}
\hrule
\bigskip\bigskip\bigskip\bigskip\bigskip\bigskip

This chapter delves into more intricate minimization problems involving Lipschitz functions with constraints, compared to those discussed in Chapter 2. For problems with linear constraints, it explores analogues of the well-known methods of conditional and reduced gradients in nonlinear programming. For problems with nonlinear constraints, it examines a generalized version of the method of feasible directions and penalty function methods. It is demonstrated that exact penalties for Lipschitz functions, as in the case of convex functions, possess a global property.

\section*{$\S$ 15. The conditional gradient method}
\label{Sec.15}
\setcounter{section}{15}
\setcounter{definition}{0}
\setcounter{equation}{0}
\setcounter{theorem}{0}
\setcounter{lemma}{0}
\setcounter{remark}{0}
\setcounter{corollary}{0}
The method of conditional gradient was originally proposed for finding extremums of smooth functions in the presence of constraints. The most direct approach to solving the problem
\begin{align*}
f(x) \rightarrow \min_x \tag*{(15.1)}
\end{align*}
subject to
\begin{align*}
x \in D \tag*{(15.2)}
\end{align*}
where $D$, is a convex, bounded, and closed set, consists in linearization of the function $f(x)$: for a given feasible point $x^k \in D$, a solution $\bar{x}^k$ is found for the problem, which has the same constraints (15.2), and as the objective function, a linear approximation of the function $f(x)$ at the point $x^k$ is used. Then, it is assumed that the direction $d^k = \bar{x}^k - x^k$ is a suitable direction for finding the minimum. Therefore, to determine the direction of movement $d^k$, the problem is solved:
\begin{align*}
(\nabla f(x), y) \rightarrow \min \tag*{(15.3)}
\end {align*}
subject to
\begin{align*}
y \in D \tag*{(15.4)}
\end {align*}

In the case when the set $D$ is formed by linear constraints, the problem (15.3), (15.4) can be solved using known methods of linear programming.

It turns out that the general scheme of the conditional gradient method can also be applied to solving nonsmooth extremal problems. The main difference lies in constructing an auxiliary sequence of vectors $z^{k}$ such that:
\begin{align*}
z^{k+1}=z^{k}+a_{k}\left(H\left(x^{k}, \alpha_{k}\right)-z^{k}\right) \tag *{(15.5)}
\end{align*}
where $a_{k}$ are some non-negative multipliers and $H\left(x^{k}, \alpha_{k}\right)$ is determined by one of the formulas:
\begin{align*}
 &\frac{1}{2 \alpha_{k}} \sum_{i=1}^{n}\left[f\left(\tilde{x}_{i}^{k}, \ldots, x_{l}^{k}+\alpha_{k}, \ldots, \tilde{x}_{n}^{k}\right)\right.  \left.- \ f\left(\tilde{x}_{1}^{k}, \ldots, x_{i}^{k}-\alpha_{k}, \ldots, \tilde{x}_{n}^{k}\right)\right] e_{i} , \tag*{(15.6)}\\
&\sum_{i=1}^{n} \frac{\left.f\left(\tilde{x}^{k}+\Delta_{k} e_{i}\right)-f \tilde{x}^{k}\right)}{\Delta_{k}} e_{i}, \tag*{(15.7)}
\end{align*}
where $\tilde{x}^k$ are independent random variables uniformly distributed in an $n$-dimensional cube with center $x^k$ and edge length $2\alpha_k$.

The direction $z^{k+1}$ is a convex combination of $z^{k}$ and the vector $H\left(x^{k}, \alpha_{k}\right)$, which makes it similar to the choice of descent direction in conjugate gradient methods. More information about the properties of the averaging operation (15.5) is provided in Chapter 7.

Finite difference methods for solving problem (15.1), (15.2) with a Lipschitz function $f(x)$ are defined by the equations
\begin{align*}
x^{k+1}=x^{k}+\rho_{k}\left(\bar{x}^{k}-x^{k}\right), \quad 0 \leq \rho_{k} \leq 1 \tag*{(15.8)}
\end{align*}
where $\bar{x}^{k}$ is the solution of the problem:
\begin{align*}
\underset{x}{\min }\left<z^{k}, x\right>, \tag*{(15.9)}
\end{align*}
\begin{align*}
x \in D \tag*{(15.10)}
\end{align*}
If there exists a way to compute the generalized gradient of the function $f(x)$, then in formula $(15.5)$ one can use
\begin{align*}
H\left(x^{k}, \alpha_{k}\right)=g\left(\tilde{x}^{k}\right) \tag*{(15.11)}
\end{align*}
where $g\left(\tilde{x}^{k}\right) \in \partial f\left(\tilde{x}^{k}\right)$.

Let the domain $D$ be defined by inequalities.
$$
f_{i}(x) \leq 0, \quad i=1, \ldots, m,
$$
and $f_{i}(x)$ be convex functions. Assume that the regularity condition is satisfied. Then, the set $X^{*}$ of stationary points of the problem $(15.1)$, $(15.2)$ can be expressed as follows: $x^{*} \in X^{*}$ if there exist vectors $g\left(x^{k}\right) \in \partial f\left(x^{k}\right), g^{i}\left(x^{k}\right) \in \partial f_{i}\left(x^{k}\right)$ such that
$$
\begin{aligned}
& g\left(x^{*}\right)=\sum_{i=1}^{m} \lambda_{i} g^{i}\left(x^{*}\right), \quad \lambda_{i} \leq 0, \\
& \lambda_{i} f_{i}\left(x^{*}\right)=0, \quad i=1, \ldots, m .
\end{aligned}
$$
\begin{theorem}
\label{th:15.1}
{Suppose the following conditions are satisfied:}
$$
\begin{aligned}
& \sum_{k=0}^{\infty} \rho_{k}=\infty, \quad \sum_{k=0}^{\infty} a_{k}^{2}<\infty, \quad \frac{\rho_{k}}{\alpha_{k} a_{k}} \rightarrow 0, \\
& \sum_{k=0}^{\infty} \frac{\rho_{k}^{2}}{\alpha_{k}^{2}}<\infty, \quad \frac{\left|\alpha_{k}-\alpha_{k+1}\right|}{\rho_{k}} \rightarrow 0, \quad \alpha_{k} \rightarrow 0 .
\end{aligned}
$$
{Then the sequence $f\left(x^{k}\right)$ converges with probability 1 and limit points of $x^{k} $ belong to $X^{*}$ with probability 1.}
\end{theorem}

{\it P r o o f.} Following the proof technique described in $\S$ 4, it suffices to show that conditions (4.12) and (4.13) are satisfied. As will become clear from the following, it is useful to consider the auxiliary problem
\begin{align*}
\left<g(y), x\right> \rightarrow \min \tag{15.12}
\end{align*}
subject to the constraints:
\begin{align*}
f_{i}(x) \leq 0, \quad i=1, \ldots, m, \tag{15.13}
\end{align*}
where $g(y) \in \partial f(y)$ and the point $y$ satisfies the constraints (15.18). Let $\bar{y}$ be a solution to this problem. For any $x$ with the property (15.13)
$$
\left<g(y), \bar{y}-x\right> \leq 0.
$$

If $y \notin X^{*}$, then
\begin{align*}
\left<g(y), \bar{y}-y\right> \leq \sigma(y)<0 . \tag*{(15.14)}
\end{align*}
Indeed, assume that $\left<g(y), \bar{y}-y\right>=0$, or in other words, $y$ is a solution to problem (15.12), (15.13). According to the necessary conditions for an extremum of problem (15.12), (15.13),
$$
\begin{aligned}
& g(y)=\sum_{i=1}^{m} \lambda_{i} g^{i}(y), \quad \lambda_{i} \leq 0, \\
& \lambda_{i} f_{i}(y)=0, \quad i=1, \ldots, m .
\end{aligned}
$$
But these relationships mean that $y \in X^{*}$. This leads to a contradiction. In fact, it takes place a more general fact:
\begin{align*}
\max _{g(y) \in \partial f(x)}\left<g(y), \bar{y}-y\right> \leq \hat{\sigma}(y)<0. \tag*{(15.15)}
\end{align*}
Let us prove this. Suppose there exist sequences $\left\{\bar{y}^{s}\right\}$ and $\left\{g^{s}(y)\right\} \in \partial f(y)$, such that
$$
\left<g^{s}(y), \bar{y}^{s}-y\right> \rightarrow 0,
$$
where $\bar{y}^{s}$ is a solution of the problem $\min_{x\in D} \left<g^{s}(y),x\right>$. For simplicity, let us denote $g^{s}(y) \rightarrow \hat{g}(y) \in \partial f(y)$. From inequality (15.14), it follows that
$$
\left<\hat {g}(y), \hat{y}-y\right> \leq \sigma(y)<0,
$$
where
$$
\left<\hat{g}(y), \hat{y}\right>=\min _{x \in D} \left<\hat{g}(y), x\right> .
$$
Choose an index $s$ large enough so that the following inequalities hold:
$$
\begin{gathered}
\left<g^{s}(y), \bar{y}^{s}-y\right> \geq \sigma(y) / 2, \\
\left\|g^{s}(y)-\hat{g}(y)\right\| \leq|\sigma(y)| /(8 C),
\end{gathered}
$$
where $\|x\| \leq C$ for all $x \in D$. Then from the relation
$$
\left<\hat {g}(y)-g^{s}(y), \hat{y}-y\right>+\left<g^{s}(y), \hat{y}-y\right> \leq \sigma(y),
$$
it follows that
$$
\left<g^{s}(y), \hat{y}-y\right> \leq \frac{3}{4} \sigma(y).
$$
Thus,
$$
\left<g^{s}(y), \hat{y}-y\right> \leq \frac{3}{4} \sigma(y), \quad\left<g^{s}(y), \overline{y}^{s}-y\right> \geq \frac{\sigma(y)}{2}.
$$
We have arrived at a contradiction, since
$$
\left<g^{s}(y), y^{s}\right> \leq\left<g^{s}(y), \hat{y}\right>.
$$

Thus, the inequality (15.15) holds. Hence, for points $y$ belonging to a sufficiently small neighborhood of a point $y^{\prime} \notin X^{*}$, we have that
\begin{align*}
\max \{\left<g(y), \bar{y}-y\right> \mid g(y) \in \partial f(y)\} \leq \hat{\sigma}<0 . \tag*{(15.16)}
\end{align*}
For smooth functions we have
$$
\begin{aligned}
f\left(x^{k+1}, \alpha_{k}\right) & =f\left(x^{k}, \alpha_{k}\right)+\left<\nabla f\left(x^{k}+\tau\left(x^{k+1}-x^{k}\right), \alpha_{k}\right), \quad x^{k+1}-x^{k}\right>= \\
=f\left(x^{k}, \alpha_{k}\right) & +\rho_{k}\left<\nabla f\left(x^{k}+\tau\left(x^{k+1}-x^{k}\right), \alpha_{k}\right)-\nabla f\left(x^{k}, \alpha_{k}\right), \bar{x}^{k}-x^{k}\right>+ \\
& +\rho_{k}\left<\nabla f\left(x^{k}, \alpha_{k}\right), \bar{x}^{k}-x^{k}\right> \leqslant f\left(x^{k}, \alpha_{k}\right)+C \rho_{k}^{2} / \alpha_{k}+ \\
& +\rho_{k}\left<\nabla f\left(x^{k}, \alpha_{k}\right), \bar{x}^{k}-x^{k}\right>,
\end{aligned}
$$
where
\begin{align*}
f\left(x^{k+1}, \alpha_{k+1}\right) \leqslant f\left(x^{k}, \alpha_{k}\right)+C \rho_{k}^{2} / \alpha_{k} & +C\left|\alpha_{k}-\alpha_{k+1}\right|+ \\
& +\rho_{k}\left<\nabla f\left(x^{k}, \alpha_{k}\right), \bar{x}^{k}-x^{k}\right>.\tag*{(15.17)}
\end{align*}

Assume that there exists a subsequence $x^{s} \rightarrow x^{\prime} \notin X^{}$ $(s \in S)$. Let condition (4.12) does not hold, i.e., all points $x^{k}(k \geq s)$ belong to a sufficiently small $\delta$-neighborhood of the point $x^{s}$, which does not intersect with $X^{*}$. From inequality (15.16) it follows that
$$
\left<g\left(x^{s}\right), \hat{x}^{s}-x^{s}\right> \leq \hat{\sigma}<0,
$$
where $\hat{x}^s$ is the solution of problem (15.12), (15.13). For sufficiently small $\delta$, we have
$$
\left<g\left(x^{s}\right), \hat{x}^{s}-x^{k}\right> \leq \hat{\sigma} / 2.
$$

Let us now note an interesting property of the operation $z^{k}$: if the conditions of the theorem are satisfied, then with probability 1,
$$
\left(z^{k}-\nabla f\left(x^{k}, \alpha_{k}\right)\right) \rightarrow 0, \quad k \rightarrow \infty .
$$
The proof of this fact is given in Chapter 7.

According to Lemma 2.5, starting from some index $k$, the distance from the vector $\nabla f\left(x^{k}, \alpha_{k}\right)$ to the set $\partial f\left(x^{s}\right)$ can be made arbitrarily small if the points $x^{k}$ belong to a sufficiently small neighborhood of $x^{s}$. Then, due to the fact that $g\left(x^{s}\right)$ is an arbitrary vector from $\partial f\left(x^{s}\right)$, for sufficiently large $k$ we have:
$$
\left<z^{k}, \hat{x}^{s}-x^{k}\right> \leq \hat{\sigma} / 4.
$$
Therefore,
$$
\left<z^{k}, \bar{x}^{k}-x^{k}\right> \leq \hat{\sigma} / 4,
$$
where $\bar{x}^{k}$ is the solution of problems (15.9), (15.10). Therefore, starting from some number $k$, we have
$$
\left<\nabla f\left(x^{k}, \alpha_{k}\right), \bar{x}^{k}-x^{k}\right> \leq \hat{\sigma} / 8, \quad k \geq \bar{s} .
$$

Thus, in inequality (15.17), we estimated the quantity $\left<\nabla f\left(x^{k}, \alpha_{k}\right), \bar{x}^{k}-x^{k}\right>$. It follows from this that
\begin{align*}
f\left(x^{k+1}, \alpha_{k+1}\right) \leqslant f\left(x^{s}, \alpha_{s}\right)+\frac{\hat{\sigma}}{16} \sum_{r=s}^{k} \rho_{r}, \tag*{(15.18)}
\end{align*}
since for sufficiently large $k$
$$
C\left[\frac{\rho_{k}}{\alpha_{k}}+\frac{\left|\alpha_{k}-\alpha_{k+1}\right|}{\rho_{k}}\right] \leqslant \frac{|\hat{\sigma}|}{16} .
$$

Taking the limit as $k \rightarrow \infty$ in inequality (15.18), we arrive at a contradiction with the assumption that $f(x)$ is bounded. Hence, condition (4.12) holds:
$$
k(s)=\min \left\{r \mid\left\|x^{r}-x^{s}\right\|>\delta, r>s\right\}<\infty .
$$
From the relation
$$
x^{k(s)}=x^{s}+\sum_{r=s}^{k(s)-1} \rho_{r}\left(\bar{x}^{r}-x^{r}\right),
$$
it follows that
$$
\sum_{r=s}^{k(s)-1} \rho_{r} \geqslant \frac{\delta}{C}.
$$
Therefore, from inequality (15.18) and the uniform convergence of $f(x,\alpha) \rightarrow f(x)$, condition (4.13) follows.
$$
\overline{\lim }_{s \rightarrow \infty} f\left(x^{k(s)}\right)<\lim _{s \rightarrow \infty} f\left(x^{s}\right).
$$
The theorem is proved.

It should be noted that if we substitute 
$H\left(x^{k}, \alpha_{k}\right)$ or $g\left(\tilde{x}^{k}\right)$ directly instead of $z^{k}$ in equations (15.9), the sequence 
$\left\{x^{k}\right\}$ 
may not converge to $X^{*}$.

The method (15.8) -- (15.10) has broader capabilities compared to the projection method onto the feasible set. The projection operation is equivalent to solving an optimization problem with a quadratic objective function. In the conditional gradient method, the projection operation is replaced by solving an optimization problem with a linear objective function. Therefore, this method is useful when it is difficult to perform the projection operation. It is also advantageous to use this method in problems with constraints that have a block structure, where the auxiliary linear programming problem can be solved relatively easily.

\section*{$\S$ 16. The reduced gradient method}
\label{Sec.16}
\setcounter{section}{16}
\setcounter{definition}{0}
\setcounter{equation}{0}
\setcounter{theorem}{0}
\setcounter{lemma}{0}
\setcounter{remark}{0}
\setcounter{corollary}{0}
The method of reduced gradient was proposed by P. Wolfe for solving nonlinear programming problems with linear constraints.
\begin{align*}
\min f(x)\tag*{(16.1)}
\end{align*}
subject to constraints
\begin{align*}
A x=b, \quad x \geqslant 0
\tag*{(16.2)}
\end{align*}
where $A$ is an $(m \times n)$-matrix of rank $m$, and $b$ is a vector in $E_m$.

The method is based on reducing the number of basic variables using non-basic independent variables, resulting in a new, reduced problem containing only non-basic variables. In the reduced problem with simple constraints, one of the unconditional minimization methods can be applied, and the next descent direction is chosen so that the new approximation satisfies the constraints (16.2). The method of conjugate gradients or quasi-Newton method is most often used for this method. The reduced gradient scheme generalizes the simplex procedure for linear programming and widely uses techniques for working with sparse matrices [71].

Let us show how the method of reduced gradient can be extended to the problem with a Lipschitz continuous function $f(x)$. Suppose that any $m$ columns of the matrix $A$ are linearly independent, and if in the equation $A x=b$ the non-zero components of the vector $x$ correspond to linearly independent columns of the matrix $A$, then there are exactly $m$ such components. Then each point in the set $S=\{x \mid A x=b, x \geq 0\}$ has at least $m$ positive components.

Let a point $x$ satisfy (16.2). The matrix $A$ can be represented as $[B, N]$ and the vector $x$ as $[x^{T}_{B}, x^{T}_{N}]^ T$, where $B$ is a nonsingular matrix of size $m\times m, x_{B}>0$. The vector $x_{B}$ is called \textit{basic}. The components of the nonbasic vector $x_{N}$ can be either positive or zero. Let
\[ g(x) = [g^{T}_{B}(x), \quad g^{T}_{N}(x)]^T, \]
where $g_{B}(x)$ is the generalized gradient of the function $f(x)$ with respect to the variables of the basis vector $x_{B}$, and $g_{N}(x)$ is the generalized gradient with respect to the nonbasic vectors variables $x_{N}$.

Let's represent the direction vector as $[d^{T}_{B}, d^{T}_{N}]^T$. Note that the equality $0 = Ad = Bd_{B} + Nd_{N}$ holds if for any $d_{N}$ we put $d_{B} = - B^{-1}Nd_{N}$. Let's call
\[ r^{T} = [r^{T}_{B}, r^{T}_{N}] = g^{T}(x) - g^{T}_{B}(x )B^{-1}A = g^{T}_{N}(x) - g^{T}_{B}(x)B^{-1}N \]
\textit{reduced gradient}. Let us explore $\left<g(x), d\right>$:
\[ \left<g(x), d\right> = \left<g_{B}(x), d_{B}\right> + \left<g_{N}(x), d_{N}\right> = [g^{T}_{N }(x) - g^{T}_{B}(x)B^{-1}N]d_{N} = \left<r_{N}, d_{N}\right>. \]

Choose $d_{N}$ so that $\left<r_{N}, d_{N}\right> < 0$ and $d_{j} \geq0 $ if $x_{j} =0 $. For each out-of-basic component $j$, we set $d_{j} = - r_{j}$ if $r_{j} \leq 0$, and $d_{j} = - x_{j}r_{j} $, if $r_{j} > 0$. This ensures that the inequality $d_{j} \geq0 $ holds if $x_{j} = 0 $. Moreover, $\left<g(x), d\right> \leq 0$, and the strict inequality holds if $d_{N}\neq0$.

\begin{lemma}
\label{lem:16.2}
{Let $d$ be chosen as described above, $g(x)$ is an arbitrary vector from the set $\partial f(x)$. Then, if $d=0$, then the Kuhn--Tucker conditions are satisfied at the point $x$.}
\end{lemma}

{\it P r o o f}. First of all, we note that the vector $d$ is possible direction, since
\[ Ad = Bd_{B} + Nd_{N} = B(- B^{-1}Nd_{N}) + Nd_{N} = 0. \]
If the variable $x_{j}$ is basic, then by assumption $x_{j}>0$. If $x_{j}$ is off-basic, then the component $d_{j}$ can be negative only if $x_{j}>0$. Thus $d_{j}\geq0$ if $x_{j}=0$, and hence, $d$ is a possible direction. Moreover, \[ \left<g(x), d\right>=\sum_{j \in I} r_{j} d_{j}, \] where $I$ is the set of indices of nonbasic variables. By the definition of $d_{j}$, it is obvious that either $d = 0$ or $\left<g(x), d\right> < 0$.

Let the Kuhn--Tucker conditions be satisfied at the point $x$. Then there are vectors $g(x)\in \partial f(x), u = (u^{T}_{B}, u^{T}_{n})^{T} \geq0$, and $ v$ such that
\begin{equation}
     [g_{B}^{T}(x), g_{N}^{T}(x)]+v^{T}[B, N]-[u_{B}^{T}, u_{N }^{T}]=0, \tag{16.3}\label{eq:16.3}
\end{equation}
\[ (u_{B}, x_{B})=0, \quad(u_{N}, x_{N})=0 . \]
Since $x_{B}>0$ and $u_{B}\geq0$, it follows from the equality $(u_{B}, x_{B}) = 0$ that $u_{B} = 0$. It follows from the first equality in $(16.3)$ that $\upsilon^{T} = - g^{T}_{B} B^{-1}$. Substituting this expression into the second equality in $(16.3)$, we get
\[ u^{T}_{N} = g^{T}_{N} - g^{T}_{B}B^{-1}.\]
In other words, $u_{N}=r_{N}$. Thus, the Kuhn--Tucker conditions are reduced to the relations $r_{N} \geq0$ and $(r_{N}, x_{N}) = 0$. However, by definition, the equality $d=0$ is true if and only if $r_{N}\geq0$ and $(r_{N}, x_{N}) = 0$. Therefore, if $d=0$, then the Kuhn--Tucker conditions are satisfied at the point $x$.

The reduced gradient method for solving problem $(16.1), \: (16.2)$, when $f(x)$ is a Lipschitz function, is defined as follows.

{\it Initial stage}. Choose an admissible point $x^{1}$ that satisfies the conditions $Ax=b, x\geq0$. Put $k=1$ and go to step $1$.

{\it Step 1}. Choose $I_{k}$, the set of $m$ indices  with the largest components of the vector $x^{k}$:
\[ x_{B}=\left\{x_{i}, i \in I_{k}\right\}, \quad B=\left\{a^{j}, j \in I_{k}\right\}, \quad N=\left\{a^{j}, j \notin I_{k}\right\}. \]
Calculate the reduced gradient
\[ r_{N}^{T}=\left(z_{N}^{k}\right)^T-\left(z_{B}^{k}\right)^T B^{-1} N .\]
Determine the direction of descent $d_{k}$:
\begin{equation}\tag{16.4}
\begin{aligned}
& d_{j}^{k}=\left\{\begin{array}{lll}
-r_{j}, & j \notin I_{k}, & r_{j} \leq 0, \\
-x_{j} r_{j}, & j \in I_{k}, & r_{j}>0,
\end{array}\right. \\
& d_{B}^{k}=-B^{-1} N d_{N} .
\end{aligned}
\end{equation}

{\it Step 2}. The sequences of points $x^{k}, z^{k}$ are recalculated according to the formulas
\begin{equation}\tag{16.5}
\begin{aligned}
x^{k+1}= & x^{k}+\sigma_k d^{k} \\
z^{k+1}= & z^{k}+a_{k}\left(H\left(x^{k}, \alpha_{k}\right)-z^{k}\right) \\
& \sigma_{k}=\min \left\{\rho_{k}, \lambda,\right\}
\end{aligned}
\end{equation}
\begin{equation}\tag{16.6}
\lambda=\left\{\begin{array}{ll}
\min\limits_{1 \leq j \leq n}\left\{-\frac{x_{j}^{k}}{d_{j}^{k}} \mid d_{j}^{k}<0\right\},&\exists \;\;d_{j}^{k} < 0, \\
\infty, & \mbox{all}\;\; d^{k}_j \geq 0.  
\end{array}\right.
\end{equation}
Set $k = k + 1$ and go to step 1.

\begin{theorem}
\label{th:16.1}
{Let the conditions of Theorem $15.1$ be satisfied. Then  the sequence $f(x^{k})$ converges with probability $1$, and the limit points of $x^{k}$ satisfy the Kuhn--Tucker conditions with probability $1$.}
\end{theorem}

{\it P r o o f}. To prove  convergence of the reduced gradient method $(16.5)$, the following two results should be taken into account. The first one is that
\begin{equation}\tag{16.7}
\left<z^{k}, d^{k}\right> < \gamma < 0
\end{equation}
for all points $x^{k}$ (with sufficiently large numbers $k$) that are contained in a neighborhood of $y$, which does not intersect $X^{*}$, the set of Kuhn--Tucker points. The inequality $(16.7)$ follows from Lemmas 2.5, 16.1, and properties of $z^{k}$ and $\partial f(x)$.

The next important point of the proof is that at the above points the step factor $\sigma_{k}$ is equal to $\rho_{k}$. Denote by $I_{s}$, the set of indices  of $m$ the largest components of the vector $x_{s}$ used to calculate $d_{s}$. Suppose there exists a subsequence $x^{s} \to \hat{x} \: (s \in S \subset \{1, ... , k, ...\})$ such that
\[ \inf \{\delta_{s}\mid s \in S\} = 0,\]
where
\[ \delta_{s} =\min\left[\min\limits_{i}\left\{-\frac{x_{j}^{k}}{d_{j}^{k}} \mid d_ {j}^{k}<0\right\}, \infty\right]. \]
Then there exists a set of indices $S^{\prime} \subset S$ such that $\delta_{s} = - x_{p}^{s}/d_{p}^{s}$ converge to zero as $ s \in S^{\prime}$, and also $I^{s} = \hat{I}$, where $x_{p}^{s} > 0, d_{p}^{s} > 0$, $p$ is an element of the set ${1, 2,..., n}$, $\hat{I}$ is the set of indices of  $m$ the largest components of the vector $\hat{x} $. Note that, according to $(16.4)$, the sequence $\{d_{p}^{s}\} \: (s \in S^{\prime})$ is bounded, and since $\delta_{s} \to 0 \: (s \in S^{\prime})$, then $x_{p}^{s} \to 0$. Thus, $\hat{x_{p}} = 0$, i.e. $p \notin \hat{I}$. But $I_{s} = \hat{I}$ for $s \in S^{\prime}$; hence $p \notin I_{s}$. Since $d_{p}^{s} < 0$, (16.4) implies $d_{p}^{s} = - x_{p}^{s}r_{p} ^{s}$. Hence follows
\[ \delta_{s} = - x_{p}^{s}/d_{p}^{s} = 1/r_{p}^{s}. \]
This shows that $r_{p}^{s} \to \infty$, which is impossible because $r^{s}$ is bounded. Thus, there exists $\delta > 0$ such that the point $x^{k} + \lambda d^{k}$ is admissible for all $\lambda \in [0, \delta]$ and for all numbers $k $, starting from some.

In what follows, the proof is similar to that given for Theorem 15.1.

\section*{$\S$ 17. The method of feasible directions}
\label{Sec.17}
\setcounter{section}{17}
\setcounter{definition}{0}
\setcounter{equation}{0}
\setcounter{theorem}{0}
\setcounter{lemma}{0}
\setcounter{remark}{0}
\setcounter{corollary}{0}
Let us now turn to the study of a numerical method for solving the general problem of nonlinear programming: to find
\begin{equation}\tag{17.1}
\min f_{0}(x)
\end{equation}
under conditions
\begin{equation}\tag{17.2}
f_{i} (x) \leq 0, \quad i = 1, ..., m,
\end{equation}
where $f_{0}(x)$ is a Lipschitz function, $f_{i}(x)$ are continuously differentiable function-ins. Here the ideas of the Topkis-Beinott possible directions method (see $[6]$) and the conditional gradient method discussed in $\S$ 15 are used.

Recall that a non-zero vector $d$ is called a possible direction at a point $x \in D$ ($D$ is a non-empty set in $E_{n}$) if there exists $\delta > 0$ such that $x + \lambda d \in D$ for all $\lambda \in [0, \delta]$. The method for solving the original problem is determined by the relation
\begin{equation}\tag{17.3}
x^{k+1} = x^{k} + \gamma_{k}d^{k}.
\end{equation}
The vector $x^k$ is a valid point. Here, both active and inactive constraints are taken into account when determining the direction of movement. In this it differs from the G. Seutendijk method of possible directions, in which the search for the direction of descent occurs on the basis of almost active constraints. A possible direction $d^{k}$ is found from the solution of the linear program problem
\begin{equation}\tag{17.4}
\min_{d,y} y
\end{equation}
under conditions
\begin{gather*} \tag{17.5}
\left<z^{k}, d\right> - y \leq 0, \\
f_{i}(x^{k}) + \left<\nabla f_i(x^{k}), d\right> - y \leq 0, \quad i = 1, ..., m, \\
-1 \leq d_{j} \leq 1, \quad j = 1, ..., n.
\end{gather*}

The sequence $\{z^{k}\}$ is constructed using the averaging operation according to the formula
\[ z^{k+1} = z^{k} + a_{k} (H(x^{k}, \alpha_{k}) - z^{k}), \]
and the vector $H(x^{k}, \alpha_{k})$ is determined by one of the formulas: $(15.6), \: (15.7), \: (15.11).$

It remains to determine the step $\gamma_{k}$. Let the domain $D$ defined by the constraints $(17.2)$ be a bounded set satisfying the regularity condition. Let's put
\[ \rho_{k}^{\prime} = \max \{ \lambda \mid x^{k} + \lambda d^{k} \in D \}. \]
Then in the formula $(17.3)$ $\gamma_{k} = \min \{\rho_{k}^{\prime}, \rho_{k}\}$, where $\rho_{k}$ is chosen from the relations
\begin{equation}\tag{17.6}
\sum_{k=0}^{\infty} \rho_{k} = \infty, \quad \rho_{k} \to 0.
\end{equation}

Let us define the set $X^{*}$ of stationary points of problem $(17.1), (17.2)$: $x^{*} \in X^{*}$ if there exist vectors $g^{\nu} (x ^{*}) \in \partial f_{\nu}(x^{*}) \: (\nu = 0, 1, ..., m)$ and numbers $\lambda_{i} \leq 0 \ : (i = 1, ..., m)$, which
\[ g^{0}(x^{*}) = \sum_{i=1}^{m}\lambda_{i}g^{i}(x^{*}), \]
\[ \lambda_{i} f_{i}(x^{*}) = 0, \quad (i= 1, ..., m).\]

\begin{theorem}
\label{th:17.1}
{Let the conditions of Theorem 15.1 be satisfied, as well as the relation (17.6). Then the limit points of the sequence $\{x^{k}\}$ belong to the set $X^{*}$ with probability $1$ and the sequence $f_{0}(x^{k})$ converges almost surely.}
\end{theorem}

{\it P r o o f}. As in the study of convergence of the conditional gradient method, we consider the auxiliary problem
\begin{equation}\tag{17.7}
\min_{d,y} y
\end{equation}
under conditions
\begin{gather*} \tag{17.8}
\left<g^{0}(x), d\right> - y \leq 0, \\
f_{i}(x) + \left<\nabla f_{i}(x), d\right> - y \leq 0, \quad i = 1, ..., m, \\
-1 \leq d_{j} \leq 1, \quad j = 1, ..., n.
\end{gather*}
where $g^{0}(x) \in \partial f_{0}(x).$

Let $x \notin X^{*}$. Then it is easy to show that $\bar{y}(x) < 0$, where $\bar{y}(x)$ is the solution to problem $(17.7), (17.8)$. This implies that for points $u$ that belong to a sufficiently small neighborhood of the point $x \notin X^{*}$, the inequality holds:
\begin{equation}\tag{17.9}
\max_{g^{0}(u) \in \partial f_{0}(u)} \bar{y}(u) < \varepsilon < 0.
\end{equation}

The next important point of the proof is that for all points
\[ x^{k} \in U_{2\hat{\delta}}(x^{s}) \equiv \{x \mid \lVert x - x^{s} \rVert \leq 2\bar{ \delta}\} \]
for sufficiently large $s$ and small $\bar{\delta}$ such that $x^{s} \to x^{\prime} \notin X^{*} \quad (s \in S)$, the choice step factor in method (17.3)  leads to $\gamma_{k} = \rho_{k}$, i.e.
satisfies the relations (17.6).

By virtue of the relation $(17.9)$, there is a non-zero vector $\bar{d}$ that satisfies the inequalities
\begin{gather*}
\left<g^{0}(x^{k}), \bar{d}\right> \leq \varepsilon, \\
f_{i}(x^{k}) + \left<\nabla f_{i}(x^{k}), \bar{d}\right> \leq \varepsilon, \quad i = 1, ..., m, \\
-1 \leq \bar{d_{j}} \leq 1, \quad j = 1, ..., n.
\end{gather*}

Using the result of Lemma $2.5$, we can choose the vector $g^{0}(x^{s}) \in \partial f_{0}(x^{s})$ and the numbers $s$ and $\bar{\delta }$ such that the inequalities
\begin{gather*}
\left<\nabla f_{0}(x^{k}, \alpha_{k}), \bar{d}\right> \leq \varepsilon/2, \\
f_{i}(x^{k}) + \left<\nabla f_{i}(x^{k}), \bar{d}\right> \leq \varepsilon/2, \quad i = 1, ..., m,\\
-1 \leq \bar{d_{j}} \leq 1, \quad j = 1, ..., n.
\end{gather*}

With probability $1$ it holds:
\[ (z^{k} + \nabla f_{0}(x^{k}, \alpha_{k})) \to 0;\]
therefore, for the optimal solution $d^{k}$ of problem $(17.4), \: (17.5)$
\begin{gather*}\tag{17.10}
\left<z^{k}, \: d^{k}\right> \leq \varepsilon/4, \\
f_{i}(x^{k}) + \left<\nabla f_{i}(x^{k}), \: d^{k}\right> \leq \varepsilon/4, \quad i = 1, .. ., m, \\
-1 \leq d_{j}^{k} \leq 1, \quad j = 1, ..., n.
\end{gather*}

Since the functions $f_{i}(x) \: (i = 1, ..., m)$ are continuously differentiable, (17.10) implies the existence of $\tau$ such that for all $\lambda \in [0, \tau]$
\begin{equation}\tag{17.11}
f_{i}(x^{k}) + \left<\nabla f_{i}(x^{k} + \lambda d^{k}), \: d^{k}\right> \leq \varepsilon/8.
\end{equation}
Since $f_{i}(x^{k}) < 0$ for all $i$, by the mean value theorem we get
\begin{equation}\tag{17.12}
\begin{aligned}
f_{i}(x^{k} + \lambda d^{k}) = f_{i}(x^{k}) + \lambda\left<\nabla f_{i}(x^{k} + \mu_ {ik}\lambda d^{k}), \: d^{k}\right> = \\
(1 - \lambda)f_{i}(x^{k}) + \lambda[f_{i}(x^{k}) + \left<\nabla f_{i}(x^{k} + \mu_{ ik}\lambda d^{k}), \: d^{k}\right>] \leq \\
\leq \lambda[f_{i}(x^{k}) + \left<\nabla f_{i}(x^{k} + \mu_{ik}\lambda d^{k}), \: d^{ k}\right>],
\end{aligned}
\end{equation}
where $0\leq \mu_{ik} \leq 1$. Since $\lambda\mu_{ik} \in [0, \tau]$, from $(17.10)$ and $(17.11)$
it follows
\[ f_{i}(x^{k} + \lambda d^{k}) \leq \lambda\varepsilon/8.\]
This means that the point $x^{k} + \lambda d^{k}$ is valid for any $\lambda \in [0, \tau]$. Thus, for sufficiently large $k$, the step size is $\gamma_{k} = \rho_{k}$. In what follows, the proof is based on the derivation of conditions (4.12), (4.13) and is carried out to some extent already in the standard way.

\section*{$\S$ 18. Penalty functions Methods}
\label{Sec.18}
\setcounter{section}{18}
\setcounter{definition}{0}
\setcounter{equation}{0}
\setcounter{theorem}{0}
\setcounter{lemma}{0}
\setcounter{remark}{0}
\setcounter{corollary}{0}
\textbf{1. Finite-difference methods for non-smooth penalties.} 
\label{Sec.18.1}
One of the efficient ways to solve extremal problems with constraints is to reduce these problems with the help of penalty functions to problems of unconstrained minimization. The most interesting are the local penalty functions, which reduce the original problem
\begin{equation}\tag{18.1}
f_{0}(x) \to \min
\end{equation}
under conditions
\begin{equation}\tag{18.2}
f_{i}(x) \leq 0, \quad i = 1, ..., m.
\end{equation}
\begin{equation}\tag{18.3}
x \in X
\end{equation}
to an unconstrained minimization of a nonsmooth function.

In convex programming $[37]$, the penalty function has the form
\begin{equation}\tag{18.4}
\Phi(x) = f_{0}(x) + \sum_{i=1}^{m} r_{i}\max \{0, f_{i}(x)\}, \quad r_{i } > 0.
\end{equation}

Let $(x^{*}, u^{*})$ be a saddle point of the Lagrange function
\[ L(x, u) = f_{0}(x) + \sum_{i = 1}^{m} u_{i}f_{i}(x),\]
\[x \in X, \quad u \geq 0, \]
i.e., $(x^{*}, u^{*})$ satisfies the relations:
\begin{equation}\tag{18.5}
\begin{aligned}
f_{i}(x^{*}) + \sum_{i = 1}^{m} u_{i}f_{i}(x^{*}) \leq f_{0}(x^{*} ) + \sum_{i = 1}^{m} u_{i}^{*}f_{i}(x^{*}) \leq f_{0}(x) + \sum_{i = 1}^{m} u_{i}^{*}f_{i}(x), \\
x \in X, \quad u\geq 0.
\end{aligned}
\end{equation}

Conditions $(18.5)$ can also be written like this:
\[ L(x^{*}, u^{*}) = \min_{x \in X}\max_{u \geq 0}L(x, u) = \max_{u \geq 0}\min_ {x \in X}L(x, u).\]

It is known that if $(x^{*}, u^{*})$ is a saddle point of the function $L(x, u), \: x \in X, \: u \geq 0$, then $x^ {*}$ is the optimal solution of problem $(18.1) - (18.3)$. The converse statement is true for the problem of convex programming under the assumption that the constraints $(18.2), (18.3)$ satisfy the Slater regularity condition, i.e., there exists a vector $\bar{x} \in X$, for which $f_{i}(\bar{x}) < 0 \: (i = 1, ..., m)$. The importance of the Slater condition is related to the fact that it implies that the set of components of the vector $u$ is bounded. If $f_{i}(x)  \; (i = 0, 1, ..., m) $ are convex functions, $X$ is a convex set, then for the function $\Phi (x)\;(18.4)$ with $r_{i} > u_{i}^{*}$ the minimization problem is equivalent to (18.1) - (18.3) in the sense that their optimal values and optimal sets coincide.

It should be noted that even for smooth $f_{i}(x)$ the function $\Phi (x)$ is not  smooth; therefore, to minimize it, it is necessary to apply methods of nondifferentiable optimization. Subgradient methods for minimizing $\Phi (x)$ are well known; therefore, we will consider finite-difference methods.

In what follows, $X$ is a bounded closed set onto which one can easily project; $\Phi (x) $ is a special type of maximum function, so the following finite difference method can be used to minimize it:
\begin{equation}\tag{18.6}
\begin{gathered}
x^{k+1}=\pi_X\left(x^k-\rho_h\left[H^0\left(x^k, \alpha_k\right)+\sum_{i=1}^m r_i^{ +} H^i\left(x^k, \alpha_k\right)\right]\right), \\
r_i^{+}= \begin{cases}0, & f_i\left(x^k\right) \leq 0, \\
r_i, & f_i\left(x^k\right)>0,\end{cases}
\end{gathered}
\end{equation}
where the vectors $H^i(x^k, \alpha_k)$ are determined by one of the formulas $(15.6), (15.7)$.
\begin{theorem}
\label{th:18.1}
{Let the conditions be met}
\[ \sum_{k = 0}^{\infty}\rho_{k} = \infty, \quad \sum_{k = 0}^{\infty} \rho_{k}^{2} < \infty, \quad \frac{\Delta_{k}}{\rho_{k}} \to \infty, \quad \alpha_{k} \to 0. \]
{Then with probability $1$ the limit points of the sequence $\{x^{k}\}$ belong to the set $X^{*}$ of minima of the problem} $(18.1)-(18.3)$.
\end{theorem}
{\it P r o o f}. It is required to show that conditions  (4.12), (4.13) are satisfied. Let $x^{s} \to x^{\prime} \notin X^{*} \: (s \in S)$. One can specify positive numbers $\Bar{s}$ and $\Bar{\delta}$ such that for all points
\[ x \in U_{2\delta(x^{s})} \equiv \{ x \mid \lVert x - x^{s} \rVert \leq 2\bar{\delta}\}, \quad s \leq \bar{s_{i}}, \]
the inequalities hold:
\[ \Phi (x) - \Phi (x^{*}) \geq \sigma > 0, \quad x^{*}\in X^{*}. \]
Then, due to the uniform convergence of the functions $f_{i} (x, \alpha_{k})$ to $f_{i}(x) \; (i = 0, 1, ..., m)$ at these points for sufficiently large $k$, we have
\[f_0\left(x^k, \alpha_k\right)+\sum_{i=1}^m r_i^{+} f_i\left(x^k, \alpha_k\right)-f_0\left(x^*, \alpha_k\right)-\sum_{i=1}^m r_i^{+} f_i\left(x^*, \alpha_k\right) \geq \frac{\sigma}{2}>0.\]
That's why
\[\left<\nabla f_0 (x^{k}, \alpha_{k}) + \sum_{i=1}^m r_i^{+} f_i(x^k, \alpha_k), \: x^{* } - x^{k}\right> \leq - \frac{\sigma}{2}.\]

Let us show that $k (s) < \infty$, where
\[ k (s) = \min \{ r \mid \lVert x^{r} - x^{s} \rVert > \delta, r > s \}, \quad \delta \leq \bar{\delta }/2 .\]

Assume the contrary: $x^{k} \in U_{\delta}(x^{s})$ for all $k > s$. For simplicity, we assume that the vectors $H^{i}(x^{k}, \alpha_{k})$ are defined by the formulas $(15.6)$. Then the following chain of inequalities takes place:
\begin{gather}
\left\|x^{k+1}-x^*\right\|^2 \leq\left\|x^k-\rho_k\left[H^0\left(x^k, \alpha_k\right)+\sum_{l=1}^m r_i^{+} H^i\left(x^k, \alpha_k\right)\right]-x^*\right\|^2= \nonumber\\
=\left\|x^k-x^*\right\|^2-2 \rho_k\left<H^0\left(x^k, \alpha_k\right)+\sum_{i=1}^m r_l^{+} H^l\left(x^k, \alpha_k\right), x^k-x^*\right>+ \nonumber\\
\quad+2 \rho_k^2\left\|H^0\left(x^k, \alpha_k\right)+\sum_{i=1}^m r_{i}^{+}H^i\left(x^k, \alpha_k\right)\right\|^2 \leq\left\|x^k-x^*\right\|^2- \nonumber\\
\quad-2 \rho_k\left<\nabla f_0\left(x^k, \alpha_k\right)+\sum_{i=1}^m r_i^{+} \nabla f_i\left(x^k, \alpha_k\right), x^k-x^*\right>+ \nonumber\\
+2 \rho_k\left<\nabla f_0\left(x^k, \alpha_k\right)-H^0\left(x^k, \alpha_k\right)+\right.\sum_{i=1}^m r_i^{+}\left[\nabla f_i\left(x^k, \alpha_k)-\right.\right. \nonumber\\
\left.\left.-H^i\left(x^k, \alpha_k\right)\right], x^k-x^*\right>+C \rho_k^2 \leq\nonumber\\
\leq\left\|x^k-x^*\right\|^2-\sigma \rho_k
+2 \rho_k\left<\nabla f_0\left(x^k, \alpha_k\right)-H^0\left(x^k, \alpha_k\right)+\right.\nonumber\\
+\sum_{i=1}^m r_i^{+}\left[\nabla f_i\left(x^k, \alpha_k\right) 
-H^i\left(x^k, \alpha_k\right)\right], x^k-x^*\left.\right>+C \rho_k^2. \tag{18.7}
\end{gather}
Let's put
\[ W(x^{k}) = \min_{x^{*} \in X^{*}} \|x^k-x^*\|^2 = \|x^k-x^{*}(k )\|^2. \]
Then from $(18.7)$ it follows that
\begin{equation}\tag{18.8}
\begin{aligned}
& W\left(x^{k+1}\right) \leq W\left(x^k\right)-\sigma \rho_k+2 \rho_k\left<\nabla f_0\left(x^k, \alpha_k\right)-H^0\left(x^k, \alpha_k\right)+\right. \\
& \quad+\sum_{i=1}^m r_i^{+}\left[\nabla f_i\left(x^k, \alpha_k\right)-H^i\left(x^k, \alpha_k\right)\right], x^k-x^{*}\left.\right>+C \rho_k^2 .
\end{aligned}
\end{equation}
Summing $(18.8)$ over $k$ and letting $k$ tend to infinity, we get a contradiction, because with probability $1$
\begin{equation}\tag{18.9}
\begin{aligned}
& \sum_{k=0}^{\infty} \rho_k\left<\right.\nabla f_0\left(x^k, \alpha_k\right)-H^0\left(x^k, \alpha_k\right)+\sum_{i=1}^m r_i^{+}\left[\nabla f_i\left(x^k, \alpha_k\right)-\right. \\
& \left.\left.-H^i\left(x^k, \alpha_k\right)\right], x^k-x^*\right>+\sum_{k=0}^{\infty} C \rho_k^2<\infty,
\end{aligned}
\end{equation}
i.e. condition $(4.12)$ is proved.

From the properties of the projection operation, we have
\[\delta<\left\|x^{k(s)}-x^s\right\| \leq \sum_{r=s}^{k(s)-1}\left\|x^r-x^s\right\| \leq C \sum_{r=s}^{k(s)-1} \rho_r .\]
It follows from relation $(18.9)$ that for sufficiently large $s$
\begin{equation}\tag{18.10}
\begin{aligned}
&\sum_{r=s}^{k(s)}\left\{2 \rho _ { r } \left<\nabla f_0\left(x^r, \alpha_r\right)\right.\right.  -H^i\left(x^r, \alpha_r\right)+\sum_{i=1}^m r_l^{+}\left[\nabla f_i\left(x^r, \alpha_r\right)- \right.\\
& \left.\left.\left.-H^i\left(x^r, \alpha_r\right)\right], x^r-x^*\right>+C \rho_r^2\right\} \leq \frac{\sigma}{2} \sum_{r=s}^{k(s)-1} \rho_r .
\end{aligned}
\end{equation}
For sufficiently large $s$ $\;x^{k(s)} \in U_{2\delta}(x^{s})$; so from  inequality $(18.10)$, we have
\[ W x^{k(s)} \leq W x^{s} - \sigma\delta/(2C),\]
whence the condition $(4.13)$ follows. The theorem has been proven.

{\it Note}. 
Non-smooth functions $\Phi\;(x)$ $(18.4)$ have an advantage over penalty functions
\[ F(x)=f_0(x)+\sum_{i=1}^m r_i\left(\max \left\{0, f_l(x)\right\}\right)^2,\]
since the minimization of $F (x)$ under the conditions $x \in X$ gives an approximate solution to problem $(18.1) - (18.3)$, which tends to the exact solution only as $r_{i} \to \infty$.

The non-smooth penalty function $\Phi(x)$ $(18.4)$ is, as a rule, a ravine-type function, so to minimize it it is expedient to use procedures for averaging the vectors $H\;(x^{k}, \alpha_{k}) \; (k = 0, 1, ..., m)$, which were considered in Ch. $4$.

\textbf{2. Global features of non-smooth penalties.} 
\label{Sec.18.2}
Consider now the problem
\begin{equation}\tag{18.11}
f_{0}(x) \to \min
\end{equation}
under conditions
\begin{equation}\tag{18.12}
f_{i}(x) \leq 0, \quad i = 1, ..., m,
\end{equation}
where $f_{i}(x)$ are Lipschitz functions. The question arises: will the function $(18.4)$ be an exact penalty function in this case? We investigate the following penalty function:
\begin{equation}\tag{18.13}
P(x, \lambda)=f_0(x)+\lambda \sum_{i=1}^m \max \left(0, f_i(x)\right), \quad \lambda>0 .
\end{equation}
Let $f_{0}(x) \to \infty$ for $\|x \|\to \infty, - \infty < \inf f_{0}(x) < \infty$, and let in minimum points $x^{* }$ of problem $(18.11), (18.12)$ the regularity conditions are satisfied: there exists $d$ such that
\[f_i^0\left(x^* ; d\right)<0,\quad \forall i \in I\left(x^*\right),\]
\[I\left(x^*\right) \equiv\left\{i \mid f_i(x^*)=0\right\} .\]

Denote by $Y$ the set of minimum points of problem $(18.11), (18.12)$.
\begin{theorem}
\label{18.2}
{There is $\bar{\lambda}>0$ such that for} $\lambda \geq \bar{\lambda}$
\[ \arg \min P(x, \lambda) \in Y.\]
\end{theorem}

{\it P r o o f (by contradiction).} We choose an arbitrarily large $\lambda>0$. Let $x_\lambda$ be the minimum point of $P(x, \lambda)$, where $x_\lambda \notin Y$. It is natural to assume that $x_\lambda$ does not satisfy the constraints $(18.12)$, because otherwise $x_\lambda \in Y$. Therefore, there exists an index $j$ such that $f_{j}(x_\lambda) > 0$. By assumption, there is a sequence
\[\{x_{\lambda_{s}}\}, \quad s \to \infty, \quad \max_{j\leq i \leq m} f_{i}(x_{\lambda_{s}}) > 0.\]
Assume without loss of generality that $x_{\lambda_{s}} \to x^{\prime}$ as $s \to \infty$. If, in turn, the point $x^{\prime}$ also does not satisfy the constraints $(18.12)$, then
\[ P (x_{\lambda_{s}}, \lambda_{s}) \to \infty, \quad s \to \infty. \]
That's why
\[ f_{i}(x^{\prime}) \leq 0, \quad i = 1, ..., m.\]
Since
\[ P (x_{\lambda_{s}}, \lambda_{s}) < f_{0}(x^{*}), \quad x^{*} \in Y, \]
\setcounter{equation}{13}
then $x_{\lambda_{s}} \rightarrow x^{\prime} \in Y$. Since $x_{\lambda_{s}}-$ minimum point $P\left(x, \lambda_{s}\right)$, it follows from the necessary extremum condition that
\begin{align}
&0 \in \partial P\left(x_{\lambda_{s}}, \lambda_{s}\right)=& \nonumber\\
&=\partial\left(f_{0}\left(x_{\lambda_{s}}\right)+\lambda_{s} \sum_{i\in I^{+}\left(x_{\lambda_{s}}\right)} f_{i}\left(x_{\lambda_{s}}\right)
+\lambda_{s} \sum_{i\in I\left(x_{\lambda_{s}}\right)} \max \left(0, f_{i}\left(x_{\lambda_{s}}\right)\right)\right),\tag{18.14}
\end{align}
\begin{eqnarray}
I^{+}\left(x_{\lambda_{s}}\right) \equiv\left\{i \mid f_{i}\left(x_{\lambda_{s}}\right)>0\right\}, \;\;\;
I\left(x_{\lambda_{s}}\right) \equiv\left\{i \mid f_{i}\left(x_{\lambda_{s}}\right)=0\right\} .\nonumber
\end{eqnarray}

Due to the properties of the generalized gradients of the sum and maximum of Lipschitz functions, formula (18.14) implies that
\begin{align}
&0 \in \partial f_{0}\left(x_{\lambda_{s}}\right)+\lambda_{s} \sum_{i \in I^+\left(x_{\left.\lambda_{s}\right)}\right)} \partial f_{i}\left(x_{\lambda_{s}}\right)+\sum_{i \in I\left(x_{\left.\lambda_{s}\right)}\right.} \bar{\lambda}_{i} \partial f_{i}\left(x_{\lambda_{s}}\right), \quad \tag{18.15}\\
&\bar{\lambda}_{i}={\lambda}_{s} \mu_{i}, \quad 0 \leqslant \mu_{i} \leqslant 1 .\nonumber
\end{align}

Note that for sufficiently large $s$ the inclusion holds:
$$
I\left(x_{\lambda_{s}}\right) \cup I^{+}\left(x_{\lambda_{s}}\right) \subset I\left(x^{\prime}\right)=\left\{i \mid f_{i}\left(x^{\prime}\right)=0\right\} \text {. }
$$

The relation $(18.15)$ means that there are vectors $g^{i}\left(x_{\lambda_{s}}\right) \in \partial f_{i}\left(x_{\lambda_{s}}\right)$ $(i=0,1, \ldots, m)$, for which
\begin{align}
g^{0}\left(x_{\lambda_{s}}\right)+\lambda_{s} \sum_{i \in I^+\left(x_{\lambda_{s}}\right)} g^{i}\left(x_{\lambda_{s}}\right)+\sum_{i \in I\left(x_{\lambda_{s}}\right)} \bar{\lambda}_{i} g^{i}\left(x_{\lambda_{s}}\right)=0.  \tag{18.16 }
\end{align}
Since the set $\{1, \ldots, m\}$ is finite, then we can choose a subsequence $\lambda_{s}$ for which $f_{i}\left(x_{\lambda_{s}}\right) >0$. It follows from the closedness property of mappings $x \rightarrow \partial f_{i}(x)$ $(i=0,1, \ldots, m)$ that
$$
\begin{aligned}
g^{0}\left(x_{\lambda_{s}}\right) & \rightarrow \bar{g}^{0}\left(x^{\prime}\right) \in \partial f_{0}\left(x^{\prime}\right), \\
g^{j}\left(x_{\lambda_{s}}\right) & \rightarrow \bar{g}^{j}\left(x^{\prime}\right) \in \partial f_{j}\left(x^{\prime}\right) .
\end{aligned}
$$

Therefore, passing to the limit in $s \rightarrow \infty$ in relation (18.16), we obtain a contradiction with the boundedness of the Lagrange multipliers at the points of the set $Y$. The theorem has been proven.

\section*{$\S$ 19. Finite-difference methods for constrained minimization of  Lipschitz functions}
\label{Sec.19}
\setcounter{section}{19}
\setcounter{definition}{0}
\setcounter{equation}{0}
\setcounter{theorem}{0}
\setcounter{lemma}{0}
\setcounter{remark}{0}
\setcounter{corollary}{0}
Note that the choice of penalty coefficients $r_{i}$ of the exact penalty function $\Phi(x)$ (18.4) is not an easy task, because the components $u_{i}^{*}\; (i=1, \ldots, m )$ of a saddle point of the Lagrange function are known only in exceptional cases. The proposed method, on the one hand, is close to the penalty function method, and on the other hand, it uses the idea of the well-known linearization method [105]. Its peculiarity lies in the fact that the penalty coefficient is calculated analytically.

Consider the problem:
\begin{equation}
\label{eqn:19.1}
f(x) \rightarrow \min_x \quad 
\end{equation}
given that
\begin{equation}
\label{eqn:19.2}
h(x) \le 0, \quad 
\end{equation}
where $f(x)$ and $h(x)$ are locally Lipschitz functions. Condition (19.2) does not narrow the problem, since the system $f_{i}(x) \leq 0$ $(i=1, \ldots, m)$ with the help of the function $h(x)=\max _{1 \leq i \leq m}f_i(x)$ reduces to $(19.2)$.

The method for solving problem (19.1), (19.2) is determined by the relation
\begin{equation}
\label{eqn:19.3}
x^{k+1}=x^{k}+\rho_{k} d^{k}, \quad 
\end{equation}
where the vector $d^{k}$ is the solution of the quadratic programming problem
\begin{equation}
\label{eqn:19.4}
\left(z_{f}^{k}, d\right)+\frac{1}{2}\|d\|^{2} \rightarrow \min_d \quad 
\end{equation}
with an additional restriction
\begin{equation}
\label{eqn:19.5}
\left(z_{h}^{k}, d\right)+h\left(x^{k}\right) \leqslant 0,\quad 
\end{equation}
which is present when the condition $h\left(x^{k}\right) \geq 0$ is satisfied. Here $x^{0}, z_{f}^{0}, z_{h}^{0}$ are arbitrary initial approximations,
\begin{equation}
\label{eqn:19.6}
\begin{array}{l}
z_{f}^{k+1}=z_{f}^{k}+a_{k}\left(H_{f}\left(x^{k}, \alpha_{k}\right)-z_{f}^{k}\right),\quad \\
\\ 
z_{h}^{k+1}=z_{h}^{k}+a_{k}\left(H_{h}\left(x^{k}, \alpha_{k}\right)-z_{h}^{k}\right), 
\end{array}   
\end{equation}
the vectors $H_{f}\left(x^{k}, \alpha_{k}\right), H_{h}\left(x^{k}, \alpha_{k}\right)$ are defined by the formulas ( 15.6), (15.7). Further, it will be shown that the direction vector $d^{k}$ is calculated as
$$
\begin{aligned}
d^{k}=-z_{f}^{k}-u_{k} z_{h}^{k},
\end{aligned}
$$
where the coefficient $u_{k}$ is determined analytically.

To prove the convergence of method (19.3) - (19.6), it is more convenient to assume that $d^{k}$ is a solution to problem (19.4), (19.5). Let us define the set $X^{*}$ of solutions to problem $(19.1),(19.2): x^{*} \in X^{*}$ if there are such $g_{f}\left(x^{*} \right) \in$ $\in \partial f\left(x^{*}\right), g_{h}\left(x^{*}\right) \in \partial h\left(x^{ *}\right), \lambda \geq 0$, that
\begin{equation}
\label{eqn:19.7}
\begin{array}{l}
g_{f}\left(x^{*}\right)+\lambda g_{h}\left(x^{*}\right)=0, \quad \\
\lambda h\left(x^{*}\right)=0, \quad h\left(x^{*}\right) \leq 0 .
\end{array}
\end{equation}

Relations (19.7) satisfy the necessary conditions of extremum  if the regularity condition is satisfied: there exists $d$ such that
\begin{equation}
\label{eqn:19.8}
h^{0}\left(x^{*} ; d\right)<0 . \quad h\left(x^{*}\right)=0 \quad 
\end{equation}

From (19.8), in particular, it follows that $0 \notin \partial h\left(x^{*}\right)$ for $h\left(x^{*}\right)=0$. We assume that the set $X^{*}$ is bounded and all points of $x^{k}$ belong to a bounded domain (see the remark on Theorem 4.1).

Consider the auxiliary problem
\begin{equation}
\label{eqn:19.9}
\left<g_{f}(y), d\right>+\frac{1}{2}\|d\|^{2} \rightarrow \min_d
\end{equation}
with an additional restriction
$$
\left<g_{h}(y), d\right>+h(y) \le 0 ,\quad 
$$
if
$$
h(y) \ge 0, \quad g_{f}(y) \in \partial f(y), \quad g_{h}(y) \in \partial h(y) .
$$
We assume that problem (19.9) is solvable for any $y$ belonging to an arbitrary bounded set from $E_{n}$.

Let $d(y)$ be a solution to problem (19.9). Let us show that if $y \notin X^{*},$ then $\operatorname{d}(y) \neq$ $\neq 0$. Suppose the contrary, $d(y)=0$. It is easy to see that in this case $h(y) \leq 0$. The vector $d(y)$ is a solution to problem (19.9) if and only if there exists $\lambda(y) \geq 0$ such that
$$
\begin{gathered}
g_{f}(y)+d(y)+\lambda(y) g_{h}(y)=0, \\
\lambda(y)\left<g_{h}(y), d(y)\right>+h(y)=0, \quad h(y) \ge 0 .
\end{gathered}
$$
Therefore, these relations for $d(y)=0$ mean that the necessary extremum conditions (19.7) are satisfied at the point $y$. Therefore, if $y \notin X^{*}$, then $d(y) \neq 0$. Since the sets $\partial f(x), \partial h(x)$ are closed and the mappings $x \rightarrow \partial f(x), x \rightarrow \partial h(x)$ are closed, this implies that for points $x$ located in a sufficiently small neighborhood of the point $x^{\prime} \notin X^{*}$, the following inequalities hold:
$$
\inf \|d(x)\| \geqslant c\left(x^{\prime}\right)>0,
$$
where the infimum is taken over all possible $g_{f}(x) \in \partial f(x)$, $g_{h}(x) \in \partial h(x)$.
\begin{theorem}
\label{th:19.1}
{Let the conditions be met:}
$$
\begin{gathered}
\sum_{k=0}^{\infty} \rho_{k}=\infty, \quad \sum_{k=0}^{\infty}\left(\frac{\rho_{k}}{\alpha_{k}}\right)^{2}<\infty, \quad \frac{\rho_{k}}{\alpha_{k} a_{k}} \rightarrow 0, \\
\sum_{k=0}^{\infty} a_{k}^{2}<\infty, \quad \frac{\left|\alpha_{k}-\alpha_{k+1}\right|}{\rho_{k}} \rightarrow 0, \quad \alpha_{k} \rightarrow 0, \quad \frac{\Delta_{k}}{\alpha_{k}} \rightarrow 0 .
\end{gathered}
$$
{Then with probability 1  limit points of the sequence } $\left\{x^{k}\right\}$ t{from} (19.3) - (19.6) { belong to the set} $X^{ *}$, the sequence $\{f\left(x^{k}\right)\}$ {converges almost surely.}
\end{theorem}
{\it P r o o f.} For the smoothed function $h\left(x, \alpha_{k}\right)$, we have
$$
\begin{gathered}
h\left(x^{k+1}, \alpha_{k}\right)=h\left(x^{k}, \alpha_{k}\right)+\\
+\left<\nabla h\left(x^{k}+\beta\left(x^{k+1}-x^{k}\right), \alpha_{k}\right), x^{k+1}-x^{k}\right>= \\
=h\left(x^{k}, \alpha_{k}\right)+\left<\nabla h\left(x^{k}, \alpha_{k}\right), x^{k+1}-x^{k}\right>+ \\
+\left<\nabla h\left(x^{k}+\beta\left(x^{k+1}-x^{k}\right), \alpha_{k}\right)-\nabla h\left(x^{k}, \alpha_{k}\right), x^{k+1}-x^{k}\right>, \\
\quad 0 \leq \beta \leq 1 .
\end{gathered}
$$

It follows that if $h\left(x^{k}\right) \ge 0$, then
$$
\begin{gathered}
h\left(x^{k+1}, \alpha_{k+1}\right) \le\\
\le h\left(x^{k}, \alpha_{k}\right)+C \frac{\rho_{k}^{2}}{\alpha_{k}}+C\left|\alpha_{k}-\alpha_{k+1}\right|+\rho_{k}\left<\nabla h\left(x^{k}, \alpha_{k}\right)-z_{h}^{k}, d^{k}\right>+ \\
+\rho_{k}\left<z_{h}^{k}, d^{k}\right> \le h\left(x^{k}, \alpha_{k}\right)-\rho_{k} h\left(x^{k}\right)+C\left|\alpha_{k}-\alpha_{k+1}\right|+C \frac{\rho_{k}^{2}}{\alpha_{k}}+C \rho_{k} \gamma_{k},
\end{gathered}
$$
where
$$
\gamma_{k}=\left|\left(\nabla h\left(x^{k}, \alpha_{k}\right)-z_{h}^{k}, d^{k}\right)\right|.
$$

In Ch. 7 it will be shown that under the conditions of this theorem with probability 1
$$
\gamma_{k} \rightarrow 0, \quad\left(z_{f}^{k}-\nabla f\left(x^{k}, \alpha_{k}\right)\right) \rightarrow 0 .
$$

It is easy to see from conditions (4.12), (4.13) that all limit points of the sequence $\left\{x^{k}\right\}$ satisfy the constraint $h(x) \le 0$.

For the function $f\left(x, \alpha_{k}\right)$, we have
$$
\begin{gathered}
f\left(x^{k+1}, \alpha_{k}\right) =\\
=f\left(x^{k}, \alpha_{k}\right)+\left<\nabla f\left(x^{k}+\beta\left(x^{k+1}-x^{k}\right), \alpha_{k}\right),  x^{k+1}-x^{k}\right>= \\
=f\left(x^{k}, \alpha_{k}\right) +\\
+\left<\nabla f\left(x^{k}+\beta\left(x^{k+1}-x^{k}\right), \alpha_{k}\right)-\nabla f\left(x^{k}, \alpha_{k}\right), x^{k+1}-x^{k}\right>+ \\
+\left<\nabla f\left(x^{k}, \alpha_{k}\right), x^{k+1}-x^{k}\right> \leq f\left(x^{k}, \alpha_{k}\right)+C \frac{\rho_{k}^{2}}{\alpha_{k}}+ \\
+\rho_{k}\left<\nabla f\left(x^{k}, \alpha_{k}\right)-z_{f}^{k}, d^{k}\right>+\rho_{k}\left<z_{f}^{k}, d^{k}\right>, \quad 0 \le \beta \le 1 .
\end{gathered}
$$
Hence, it follows that
\begin{eqnarray}
f\left(x^{k+1}, \alpha_{k+1}\right) 
&\le& 
f\left(x^{k}, \alpha_{k}\right)  +C \frac{\rho_{k}^{2}}{\alpha_{k}}+C\left|\alpha_{k}-\alpha_{k+1}\right|+ \nonumber\\
&& +\rho_{k}\left<z_{f}^{k}, d^{k}\right>+\rho_{k}\left<\nabla f\left(x^{k}, \alpha_{k}\right)-z_{f}^{k}, d^{k}\right>.\label{eqn:19.10}
\end{eqnarray}

The proof of convergence of the sequence $\left\{x^{k}\right\}$ to $X^{*}$ is carried out from conditions (4.12), (4.13). Let condition (4.12)  be not satisfied, i.e., all points $x^{k}\;(k>\bar{s}, \bar{s} \in S)$ are contained in a sufficiently small $\sigma$-neighborhood of the point $x^{s}$, which does not intersect $X^{*}, x^{s} \rightarrow x^{\prime} \notin X^{*}\;(s \in S)$. If $h\left(x^{k}\right)<0$, then $z_{f}^{k}=-d^{k}$. Therefore, for sufficiently large $k$ $d^{k} \neq 0$, that follows from Lemma 2.5 and the conditions that with probability 1
$$
\left(z^{k}-\nabla f\left(x^{k}, \alpha_{k}\right)\right) \rightarrow 0, \quad 0 \notin \partial f\left(x^{\prime}\right) .
$$

Let $h\left(x^{k}\right) \ge 0$. Since $d^{k}$ is a solution to problem (19.4), (19.5), there exists $\lambda_{k} \ge 0$ such that
\begin{equation}
\label{eqn:19.11}
\begin{array}{l}
z_{f}^{k}+d^{k}+\lambda_{k} z_{h}^{k}=0,  \\
\lambda_{k}\left(\left<z_{h}^{k}, d^{k}\right>+h\left(x^{k}\right)\right)=0 .
\end{array}
\end{equation}
From here it follows
\begin{equation}
\label{eqn:19.12}
\left<z_{f}^{k}, d^{k}\right>=\lambda_{k} h\left(x^{k}\right)-\left\|d^{k}\right\|^{2}. \quad 
\end{equation}

Consider the problem
\begin{equation}
\label{eqn:19.13}
\begin{array}{l}
\left<g_{f}\left(x^{\prime}\right), d\right>+\frac{1}{2}\|d\|^{2} \rightarrow \min_d, \quad \\
\left<g_{h}\left(x^{\prime}\right), d\right>+h\left(x^{\prime}\right) \leqslant 0, \quad h\left(x^{\prime}\right) \geqslant 0, \quad 
\end{array}
\end{equation}
where
$$
g_{f}\left(x^{\prime}\right) \in \partial f\left(x^{\prime}\right), \quad g_{h}\left(x^{\prime}\right) \in \partial h\left(x^{\prime}\right).
$$
It suffices to study the case $h\left(x^{\prime}\right)=0$. As noted above, $d\left(x^{\prime}\right) \neq 0$, where $d\left(x^{\prime}\right)$ is the solution to problem (19.13). Besides,
$$
\left<g_{f}\left(x^{\prime}\right), d\left(x^{\prime}\right)\right>+\frac{1}{2}\left\|d\left(x^{\prime}\right)\right\|^{2} \leq \varphi\left(x^{\prime}\right)<0,
$$
since $d\left(x^{\prime}\right)=0$ is not a solution to problem (19.13).

Let us show that there exists an admissible vector of problem (19.4), (19.5) slightly different from $d\left(x^{\prime}\right)$. Let
$$
\begin{gathered}
\tilde{g}_{f}\left(x^{\prime}\right)=\mbox{argmin} _{y \in \partial f (x^{\prime})}\left\|z_{f}^{k}-y\right\|, \\
\tilde{g}_{h}\left(x^{\prime}\right)=\mbox{argmin}_{y \in \partial h\left(x^{\prime}\right)}\left\|z_{h}^{k}-y\right\| .
\end{gathered}
$$

We will look for a solution $d^{k}$ of problem $(19.4),(19.5)$ in the form
$$
d=d\left(x^{\prime}\right)+\Delta d, \quad \Delta d \equiv-\mu \frac{\tilde{g}_{h}\left(x^{\prime}\right)}{\left\|\tilde{g}_{h}\left(x^{\prime}\right)\right\|},
$$
since under the regularity condition it holds $0 \notin \partial h\left(x^{\prime}\right)$. Then 
$$\left<z_{h}^{k}, d\right>+h\left(x^{k}\right)=$$
$$
=\left<z_{h}^{k}-\tilde{g}_{h}\left(x^{\prime}\right), d\left(x^{\prime}\right)\right>+\left<\tilde{g}_{h}\left(x^{\prime}\right), d\left(x^{\prime}\right)\right>-\mu \frac{\left<z_{h}^{k}, g_{h}\left(x^{\prime}\right)\right>}{\left\|\tilde{g}_{h}\left(x^{\prime}\right)\right\|}+h\left(x^{k}\right) .
$$
Since with probability 1,
$$
\left(z_{h}^{k}-\nabla h\left(x^{k}, \alpha_{k}\right)\right) \rightarrow 0,
$$
for sufficiently large $k$ and small $\delta$ the quantity $\left<z_{h}^{k}-\tilde{g}_{h}\left(x^{\prime}\right), d\left(x^{\prime}\right)\right>$ is also small. Since $\left<\tilde{g}_{h}\left(x^{\prime}\right), d\left(x^{\prime}\right)\right> \le 0$, then one can choose such a small $\mu>0$ that the inequality
holds,
$$
\left<z_{h}^{k}, d\right>+h\left(x^{k}\right) \le 0.
$$

Let us now estimate the value of the objective function (19.4) at the point
$d=d\left(x^{\prime}\right)+\Delta d:$
$$
\begin{aligned}
\left<z_{f}^{k}, d\right>+ & \frac{1}{2}\|d\|^{2}= \\
= & \left<z_{f}^{k}-\tilde{g}_{f}\left(x^{\prime}\right), d\right>+\left<\tilde{g}_{f}\left(x^{\prime}\right), d\left(x^{\prime}\right)+\Delta d\right>+\frac{1}{2}\left\|d\left(x^{\prime}\right)\right\|^{2}= \\
= & \left<\tilde{g}_{f}\left(x^{\prime}\right), d\left(x^{\prime}\right)\right>+\left<\tilde{g}_{f}\left(x^{\prime}\right), \Delta d\right>+\frac{1}{2}\left\|d\left(x^{\prime}\right)\right\|^{2}+\left<d\left(x^{\prime}\right), \Delta d\right>+ \\
& +\frac{1}{2}\|\Delta d\|^{2}+\left<z_{f}^{k}-\tilde{g}_{f}\left(x^{\prime}\right), d\right> .
\end{aligned}
$$

For the same reasons, for sufficiently large numbers $k$
$$
\left<z_{f}^{k}, d\right>+\frac{1}{2}\| d\|^{2}<\frac{\varphi\left(x^{\prime}\right)}{2}<0 .
$$

This implies that, starting from some number $k$, we have
$$
\left\|d^{k}\right\|>\sigma>0.
$$
From relations (19.10), (19.12), it follows that
$$f\left(x^{k+1}, \alpha_{k+1}\right) \le f\left(x^{k}, \alpha_{k}\right)+\rho_{k} \lambda_{k} h\left(x^{k}\right)-$$
$$
-\rho_{k}\left\|d^{k}\right\|^{2}+C \frac{\rho_{k}^{2}}{\alpha_{k}}+C\left|\alpha_{k}-\alpha_{k+1}\right|+\rho_{k}\left<\nabla f\left(x^{k}, \alpha_{k}\right)-z_{f}^{k}, d^{k}\right> .
$$
It follows from the conditions of the theorem that starting from some number $k$, for sufficiently small $\delta$
$$
\lambda_{k} h\left(x^{k}\right)+C\left[\frac{\rho_{k}}{\alpha_{h}}+\frac{\left|\alpha_{k}-\alpha_{k+1}\right|}{\rho_{k}}\right]+\left<\nabla f\left(x^{k}, \alpha_{k}\right)-z_{f}^{k}, d^{k}\right> \le \frac{\sigma^{2}}{2} .
$$
That's why
$$
f\left(x^{k+1}, \alpha_{k+1}\right) \le f\left(x^{s}, \alpha_{s}\right)-\frac{\sigma^{2}}{2} \sum_{r=s}^{k} \rho_{r}.
$$

Passing to the limit in $k \rightarrow \infty$, we obtain a contradiction. Therefore, condition (4.12) is satisfied:
$$
k(s) \equiv \min \left\{r \mid\left\|x^{r}-x^{s}\right\|>\delta, r>s\right\}<\infty .
$$
From the relation
$$
x^{k(s)}=x^{s}+\sum_{r=s}^{k(s)-1} \rho_{r} d^{r}
$$
it follows that
$$
\sum_{r=s}^{k(s)-1} \rho_{r}>\frac{\delta}{C},
$$
that's why
$$
\varlimsup_{s \rightarrow \infty} f\left(x^{k(s)}\right) \le \lim _{s \rightarrow \infty} f\left(x^{s}\right)-\sigma^{2} \delta /(2 C)
$$
i.e., condition (4.13) is satisfied.

It is easy to see that problem (19.4), (19.5) is solved analytically. If $h\left(x^{k}\right)<0$, then $d^{k}=-z_{f}^{k}$. If $h\left(x^{k}\right) \geq 0$, then relations (19.11) imply that
$$
d^{k}=-z_{f}^{k}-u_{k} z_{h}^{k},
$$
where
$$
u_{k}=\max \left\{0, \frac{h\left(x^{k}\right)-\left<z_{f}^{k}, z_{h}^{k}\right>}{\left\|z_{h}^{k}\right\|^{2}}\right\}.
$$

\section*{$\S$ 20. Arrow-Hurwicz finite-difference method with averaging}
\label{Sec.20}
\setcounter{section}{20}
\setcounter{definition}{0}
\setcounter{equation}{0}
\setcounter{theorem}{0}
\setcounter{lemma}{0}
\setcounter{remark}{0}
\setcounter{corollary}{0}
Gradient methods for minimizing convex nonsmooth functions based on the well-known Arrow-Hurwitz procedure [137] have been studied in great detail by now. In this section, we study the finite difference method for solving the problem
\begin{equation}
\label{eqn:20.1}
f_{0}(x) \rightarrow \min_x \quad  \tag{20.1}
\end{equation}
under restrictions
\begin{equation}
\label{eqn:20.2}
f_{i}(x) \le 0, \quad i=1, \ldots, m, \quad \tag{20.2}
\end{equation}
\begin{equation}
\label{eqn:20.3}
x \in X, \quad \tag{20.3}
\end{equation}
where $f_{i}(x)\;(i=1, \ldots, m)$ are convex functions, $X$ is a convex bounded closed set from $E_{n}$.
We introduce the Lagrange function
$$
L(x, u)=f_{0}(x)+\sum_{i=1}^{m} u_{i} f_{i}(x) .
$$
Consider the following method for solving problem (20.1) - (20.3):
\begin{equation}
\label{eqn:20.4}
x^{k+1}=\pi_{X}\left(x^{k}-\rho_{k}\left[H^{0}\left(x^{k}, \alpha_{k}\right)+\sum_{i=1}^{m} u_{i}^{k} H^{i}\left(x^{k}, \alpha_{k}\right)\right]\right), \quad \tag{20.4}
\end{equation}
\begin{equation}
\label{eqn:20.5}
u^{k+1}=\pi_{U}\left(u^{k}+\rho_{k} L_{u}\left(x^{k}, u^{k}\right)\right). \quad \tag{20.5}
\end{equation}

Here $L_{u}(x, u)$ is the gradient of the function $L(x, u)$ with respect to the variables $u$. The boundedness of the optimal dual variables $u$ follows from the regularity conditions, in particular, the Slater condition; therefore, the set $U$ onto which the projection is carried out in formula (20.5) is also assumed to be convex, closed, and bounded.

The finite difference vectors $H^{i}\left(x^{k}, \alpha_{k}\right)(i=0,1, \ldots, m)$ are determined by one of the formulas similar to (15.6), (15.7). The convergence of the trajectories $\left\{x^{k}\right\},\left\{u^{k}\right\}$ to the set of saddle points of the Lagrangian $L(x, u)$ will be studied in the Cesaro sense. Namely, let
\begin{equation}
\label{eqn:20.6}
\hat{x}^{k}=\sum_{s=0}^{k} \rho_{s} x^{s} / \sum_{s=0}^{k} \rho_{s}, \quad \tag{20.6} 
\end{equation}
\begin{equation}
\label{eqn:20.7}
 \hat{u}^{k}=\sum_{s=0}^{k} \rho_{s} u^{s} / \sum_{s=0}^{k} \rho_{s^{\circ}} \quad \tag{20.6}
\end{equation}

It turns out that, under natural assumptions, the limit points of the sequences $\left\{\hat{x}^{k}\right\},\left\{\hat{u}^{k}\right\}$ belong to the set of saddle points of $L(x, u)$.

It is important to note that averaging (20.6), (20.7) makes it possible to avoid imposing an additional condition on the strict convexity of the objective function $f_{0}(x)$ or a different choice of step factors in the coordinates $x$ and $u$.
\begin{theorem}
\label{th:20.1}
{Let the conditions be met:}
\begin{equation}
\label{eqn:20.8}
\sum_{k=0}^{\infty} \rho_{k}=\infty, \quad \sum_{k=0}^{\infty} \rho_{k}^{2}<\infty, \quad \frac{\Delta_{k}}{\alpha_{k}} \rightarrow 0, \\
\quad \alpha_{k} \rightarrow 0, \quad \sum_{k=0}^{\infty} \rho_{k} \alpha_{k}<\infty .\quad \tag{20.8}
\end{equation}
{Then with probability 1, the limit points of the sequences} $\left\{\hat{x}^{k}\right\},\left\{\hat{u}^{k}\right\}$ \textit {belong to the set of saddle points} of $L(x, u)$.
\end{theorem}

{\it P r o o f}. The proof will be carried out for the case when the vectors $H^{i}\left(x^{k} \alpha_{k}\right)$ $(i=0,1, \ldots, m)$ are defined by the formula (15.6) .

Let us introduce notation:
$$
\begin{gathered}
v^{k}=H^{0}\left(x^{k}, \alpha_{k}\right)+\sum_{i=1}^{m} u_{i}^{k} H^{i}\left(x^{k}, \alpha_{k}\right), \\
\gamma_{k}=H^{0}\left(x^{k}, \alpha_{k}\right)-\nabla f_{0}\left(x^{k}, \alpha_{k}\right)+\\
+\sum_{i=1}^{m} u_{i}^{k}\left(H^{i}\left(x^{k}, \alpha_{k}\right)-\nabla f_{i}\left(x^{k}, \alpha_{k}\right)\right),
\end{gathered}
$$
$\left(X^{*}, u^{*}\right)$ is the set of saddle points of the function $L(x, u)$. The inequalities hold true
\begin{eqnarray}
\left\|x^{k+i}-x^{*}\right\|^{2} 
\le\left\|x^{k}-\rho_{k} v^{k}-x^{*}\right\|^{2}=\left\|x^{k}-x^{*}\right\|^{2}+2 \rho_{k}\left(v^{k}, x^{*}-x^{k}\right)+\rho_{k}^{2}\left\|v^{k}\right\|^{2} \nonumber \\
\le\left\|x^{k}-x^{*}\right\|^{2}+2 \rho_{k}\left<\nabla f_{0}\left(x^{k}, \alpha_{k}\right)+\sum_{i=1}^{m} u_{i}^{k} \nabla f_{i}\left(x^{k}, \alpha_{k}\right), x^{*}-x^{k}\right>\nonumber \\
+2 \rho_{k}\left<\gamma_{k}, x^{*}-x^{k}\right>+C \rho_{k}^{2} \le\left\|x^{k}-x^{*}\right\|^{2}
+2 \rho_{k}\left(L\left(x^{*}, u^{k}\right)-L\left(x^{k}, u^{k}\right)\right)
\nonumber \\
\quad+C \rho_{k} \alpha_{k}+2 \rho_{k}\left<\gamma_{k}, x^{*}-x^{k}\right>+C \rho_{k}^{2} \le\left\|x^{k}-x^{*}\right\|^{2}
+2 \rho_{k}\left(\bar{L}\left(x^{*}\right)-L\left(x^{k}, u^{k}\right)\right)
\nonumber\\
+C \rho_{k} \alpha_{k}+C \rho_{k}^{2}+2 \rho_{k}\left<\gamma_{k}, x^{*}-x^{k}\right>, \quad \text{( 20.9) } \nonumber\\
\bar{L}(x)=\max _{u \in U} L(x, u) .\;\;\;\;\;\;\;\;\;\;\;\;\;\;\;\;\;\;\;\;\nonumber
\end{eqnarray}

For an arbitrary vector of dual variables $u$, the following inequalities hold:
\begin{eqnarray}
\left\|u^{k+1}-u\right\|^{2} &\le&\left\|u^{k}+\rho_{k} L_{u}\left(x^{k}, u^{k}\right)-u\right\|^{2} \nonumber\\
&\le&\left\|u^{k}-u\right\|^{2}+2 \rho_{k}\left(L\left(x^{k}, u^{k}\right)-L\left(x^{k}, u\right)\right)+C \rho_{k}^{2}. \;\;\;\text{( 20.10) }\nonumber
\end{eqnarray}

Adding the inequalities (20.9),(20.10), we get
$$
\begin{aligned}
\left\|x^{k+1}-x^{*}\right\|^{2} & +\left\|u^{k+1}-u\right\|^{2} \leqslant\left\|x^{k}-x^{*}\right\|^{2}+\left\|u^{k}-u\right\|^{2}+ \\
&+2 \rho_{k}\left(\bar{L}\left(x^{*}\right)-L\left(x^{k}, u\right)\right)+C \rho_{k} \alpha_{k}+ \\
&+C \rho_{k}^{2}+2 \rho_{k}\left<\gamma_{k}, x^{k}-x^{k}\right> .
\end{aligned}
$$

After summing over $k=0,1, \ldots, s$, we have
$$
\begin{aligned}
& \left(\sum_{k=0}^{s} \rho_{k}\right)^{-1} \sum_{k=0}^{s} \rho_{k}\left(L\left(x^{k}, u\right)-\bar{L}\left(x^{*}\right)\right) \leqslant \\
& \leqslant\left[\left\|x^{1}-x^{*}\right\|^{2}+\left\|u^{1}-u\right\|^{2}-\left\|x^{s+1}-x^{*}\right\|^{2}-\left\|u^{s+1}-u\right\|^{2}+\right. \\
& \left.+\sum_{k=0}^{s} C\left(\rho_{k}^{2}+\rho_{k} \alpha_{k}\right)+2 \sum_{k=0}^{s} \rho_{k}\left<\gamma_{k}, x^{*}-x^{k}\right>\right]\left(2 \sum_{k=0}^{s} \rho_{k}\right)^{-1} \text {. }
\end{aligned}
$$

From the definition of $\hat{x}^{k}, \;\hat{u}^{k}$ and the convexity of $L$ in $x$, we get
$$
\begin{aligned}
&L\left(\hat{x}^{k}, u\right)-\bar{L}\left(x^{*}\right)\leq \\
&\le\left[C+\sum_{k=0}^{s} C\left(\rho_{k}^{2}+\rho_{k} \alpha_{k}\right)+2 \sum_{k=0}^{s} \rho_{k}\left<\gamma_{k}, x^{*}-x^{k}\right>\right]\left(2 \sum_{k=0}^{s} \rho_{k}\right)^{-1}.
\end{aligned}
$$

Thus, $\bar{L}\left(\hat{x}^{k}\right) \rightarrow \bar{L}\left(x^{*}\right)$ with probability 1 because with probability 1 the series $\sum_{k=0}^{\infty} \rho_{k}\left<\gamma_{k}, x^{*}-x^{k}\right>$ converges. In other words, the limit points of the sequence $\left\{\hat{x}^{k}\right\}$ belong to the set $X^{*}$ with probability 1. Similarly, it is shown that with probability 1,
$$
\underline{L}\left(\hat{u}_{k}\right) \rightarrow \max _{u} L(u), \quad \underline{L}(u) \equiv \min _{x \in X} L(x, u) \text {. }
$$
The theorem has been proven.

\newpage

\begin{flushright}
CHAPTER 6
\label{Ch.6}

\textbf{RANDOM LIPSCHITZ AND RANDOM \\GENERALIZED DIFFERENTIABLE FUNCTIONS}
 
\underline{\hspace{15cm}}

\end{flushright}
\bigskip\bigskip\bigskip\bigskip\bigskip\bigskip

In this chapter, we study properties of random Lipschitz and random generalized differentiable functions that are necessary to justify methods of stochastic programming. The measurability of some generalized gradient mappings with respect to a set of deterministic and random variables is proved, and some results on the permutation of the signs of mathematical expectation and generalized differentiation are given.

\section*{$\S$ 21. Measurable multivalued mappings}
\label{Sec.21}
\setcounter{section}{21}
\setcounter{definition}{0}
\setcounter{equation}{0}
\setcounter{theorem}{0}
\setcounter{lemma}{0}
\setcounter{remark}{0}
\setcounter{corollary}{0}

We need some information from the theory of (measurable) multivalued mappings.

Let $E_n$ be $n$-dimensional Euclidean space; $2^{E_n}$ be the set of subsets of $ E_n$; $\mathcal{B} \text { (or } \mathcal{B}_{E_n})$ be $\sigma$-algebra of Borel subsets of 
$E_{n}$ ; $\Theta$ be a set; $\Sigma$ be $\sigma$-algebra of measurable subsets of $\Theta$ ; $\mathcal{B} \times \Sigma$ be the minimal algebra in $E_{n} \times \Theta$, containing all sets of the form $B \times T$, where $B \in \mathcal{B}, T \in \Theta$.

\textbf{1. Multivalued mappings.} 
\label{Sec.21.1}
Let us give some general information about multivalued mappings (see [60, 77] for details).
\begin{definition}
\label{def:21.1}
{A multi-valued mapping} from the space $X$ to the space $Y$ is any mapping
$$
G: X \rightarrow 2^{Y}, \quad \operatorname{dom} G=\{x \in X \mid G(x) \neq \varnothing\}
$$

A mapping $G$ is called closed-valued (convex-valued, compact-valued) if its values $G(x)$ are closed (convex, compact) sets in $Y$.
\end{definition}
\begin{definition}
\label{def:21.2}
A multivalued mapping $G: E_{n} \rightarrow 2^{E_{m}}$ is called \textit{closed} if its graph 
$$
\Gamma=\left\{(x, y) \in E_{n+m} \mid x \in \operatorname{dom}G, y\in G(x)\right\}
$$
is a closed set in $E_{n+m}$, i.e. from, $x^{k} \rightarrow x$ and $g^{k} \rightarrow g\left(g^{k} \in G\left(x^{k}\right)\right)$, it follows $g \in G(x)$.
\end{definition}
\begin{definition}
\label{def:21.3}
A multivalued mapping $G: E_{n} \rightarrow 2^{E_{m}}$ is called \textit{upper semicontinuous} if for any point $x \in \operatorname{dom} G$ and any neighborhood $U$ of the zero  in $E_{m}$ there is a neighborhood  $V$ of the zero in $E_{n}$ such that $G(x+V) \subset G(x)+U$, i.e., if for any point $x \in \operatorname{dom} G$ and for any $\varepsilon>0$ there is $\delta>0$ such that $\rho(g$, $G(x))<\varepsilon$ for all $g \in G(y),\|y-x\| \leq\delta$.
\end{definition}
\setcounter{page}{211}
\begin{lemma}
\label{lem:21.1}
A multivalued mapping $G: E_n\rightarrow 2^{E_n}$ is upper semi-continuous if and only if for any $x\in dom G \subset E_n$ and for any $x^k \rightarrow x$, $x^k\in \mbox{ dom }G$ the deviation
\[
\Delta\left(G(x^k),G(x)\right)=
\underset{g^k\in G(x^k)}{\mathop{\sup}} \,\underset{g\in G(x)}{\mathop{\inf}}
\|g-g^k\|\rightarrow 0,\;\;\;
k\rightarrow\infty. 
\]
\end{lemma}
\begin{theorem}
\label{th:21.1}
Let $G: E_n\rightarrow 2^{E_n}$ be a multi-valued mapping. Then:

the closed mapping $G$ is closed-valued;

a closed-valued upper semi-continuous mapping $G$ is closed;

a closed locally bounded (bounded in each compact) mapping $G$ is upper semi-continuous;

a finite (bounded at each point) upper semi-continuous  mapping $G$ is locally bounded (bounded in each compact);

in particular, a compact-valued upper semi-continuous mapping is locally bounded.
\end{theorem}

{\bf 2. Measurable multivalued mappings.} 
\label{Sec.21.2}
Consider multivalued mappings from a measurable space $\left(\Theta,\Sigma \right)$ in $E_n$. Let
\[
G:\Theta\rightarrow 2^{E_n},\;\;\;\mbox{dom }G=\{\theta\in\Theta|\,G(\theta)\neq\emptyset \}.
\]
\begin{definition}
\label{def:21.4}
A single-valued mapping $g:\Theta\rightarrow 2^{E_n}$ is called a section (selector) of a multi-valued mapping $G:\Theta\rightarrow 2^{E_n}$, if 
$g(\theta)\in G(\theta)$ for all $\theta\in\Theta$.
\end{definition}
\begin{definition}
\label{def:21.5}
A multivalued mapping $G:\Theta\rightarrow 2^{E_n}$ is called measurable, if for any closed set $C\subset E_n$ it takes place
\[
G^{-1}(C)=\{\theta\in\Theta|\,G(\theta)\cap C\neq \emptyset\}\in\Sigma.
\]
\end{definition}
\begin{definition}
\label{def:21.6}
A multivalued mapping $G:\Theta\rightarrow 2^{E_n}$ 
is normal if it is measurable and the sets $G(\theta)$ are closed for all 
$\theta\in\mbox{dom }G.$
\end{definition}
\begin{theorem}
\label{th:21.2}
[151, 144, 175]. The following statements concerning the closed-valued mapping $G:\Theta\rightarrow 2^{E_n}$ are equivalent:

a) $G$ is measurable;

b) for any compact $K\subset E_n$, the set
\[
G^{-1}(K)=\{\theta\in\Theta|\,G(\theta)\cap K\neq \emptyset\}\in\Sigma.
\]
c) for any open set $B\subset E_n$
\[
G^{-1}(B)=\{\theta\in\Theta|\,G(\theta)\cap B\neq \emptyset\}\in\Sigma.
\]
d) there exists a Castaign representation [144] for $G$, i.e. $\mbox{dom }G$
is measurable, and there exists a countable family of measurable sections  
$g_l: \Theta\rightarrow E_n$ $(l=\,1,2,...)$ of the mapping $G$ such that the set 
$\{g_l(\theta),\; l=1,2,\ldots\}$ is dense in $G(\theta)$ for $\theta\in\Theta$.
\end{theorem}
\begin{lemma}
\label{lem:21.2}
A closed mapping $G:\Theta\rightarrow 2^{E_n}$ is
$\mathcal{B}_{E_n}$-measurable (Borelian) and hence has sections which are Borelian functions.
\end{lemma}

Indeed, for any compact $K\subset E_n$, the set 
$\{x\in E_n|\, G(x)\cap K\neq\emptyset\}$
is closed and hence Borelian; then by Theorem 21.2, the mapping $G$  is $\mathcal{B}_{E_n}$-measurable and has Borel sections.

\begin{lemma}
\label{lem:21.3}
[49]. Let a multivalued mapping $G: x\rightarrow G(x)\subset E_n$ is measurable; then the mapping $\mbox{co } G: x\rightarrow\mbox{co }G(x)\subset E_n$
is also measurable.
\end{lemma}
\begin{definition}
\label{def:21.7}
Any function $\varphi:E_n\times\Theta\rightarrow E_1\cup\{+\infty\}$  will be called an {\it integrant}.

Let us denote by $\varphi_\theta$ a function from $E_n$ to $E_1$ such that 
$\varphi_\theta(x)=\varphi(x,\theta)$; its epi-graph is by definition the set
\[
\mbox{epi }\varphi_\theta=\{(x,\alpha)\in E_{n+1}|\,\alpha\ge \varphi_\theta(x)\}.
\]

The integrant $\varphi$  is called normal if the multivalued mapping $\theta\rightarrow \mbox{epi }\varphi_\theta$ is normal.

The integrant $\varphi$  is called convex if the functions $\varphi_\theta(x)$ are convex for all $\theta\in\Theta$.
\end{definition}
\begin{definition}
\label{def:21.8}
The function $f: E_n\times \Theta\rightarrow E_1$ is called Carath\'{e}odory one if it is measurable in $\theta\in\Theta$ for every $x\in E_n$ and 
is continuous in $x\in E_n$ for all $\theta\in\Theta$.
\end{definition}
\begin{lemma}
\label{lem:21.4}
[110]. The Carath\'{e}odory function $f$ is a normal integrant.
\end{lemma}
\begin{lemma}
\label{lem:21.5}
[175]. If $f$ is a normal convex integrant,
then the multivalued mapping
\[
\theta\rightarrow K(\theta)=\{x\in E_n|\,\varphi(x,\theta)\le 0\}
\]
is normal.
\end{lemma}
\begin{definition}
\label{def:21.9}
Mapping $g:E_n\times\Theta\rightarrow E_m$ is called
$\mathcal{B}\times\Sigma$-measurable (or $\mathcal{B}_{E_n}\times\Sigma$-measurable) if
$\{(x,\theta)\in E_n\times\Theta|\, g(x.\theta)\le c\}\in \mathcal{B}_{E_n}\times\Sigma$ 
for any $c\in E_m$.
\end{definition}
\begin{definition}
\label{def:21.10}
A mapping $g$ is called {\it Carath\'eodorian}, if it is measurable in 
$\theta\in\Theta$  for every $x\in E_n$ and continuous in $x\in E_n$ for all 
$\theta \in\Theta$.
\end{definition}
\begin{lemma}
\label{lem:21.6}
[151]. The Carathéodory function $f$ and, as a consequence, Carathéodory mapping $g$ is $\mathcal{B}\times \Sigma$-measurable.
\end{lemma}
\begin{remark}
\label{rem:21.1}
Further, we will also have to deal with mappings of the form 
$G: E_n\times\Theta\rightarrow 2^{E_m}$ and integrants of the form 
$\varphi: E_m\times(E_n\times\Theta)\rightarrow E_1\cup\{+\infty\}$. 
The above definitions and statements obviously apply to them.

For example, a multivalued mapping $G: E_n\times\Theta\rightarrow 2^{E_m}$ is called $\mathcal{B}\times \Sigma$-measurable, if for any closed set 
$C\subset E_n$
\[
G^{-1}(C)=\{(x,\theta)|\,G(x,\theta)\cap C\neq \emptyset\}\in\mathcal{B}\times\Sigma.
\]

The integrand $\varphi(g, x, \theta)$: $ E _n \times ( E _n \times \Theta) \to E_1\cup \{+\infty\}$ is called normal, if the multivalued mapping
\begin{equation} \notag
	( x , \theta) \to epi\,\varphi_{x,\theta} = \{(g, \alpha) \in E _{m+1} | \alpha \geq \varphi (g, x, \theta)\}
\end{equation}
is
$\mathcal{B} \times \Sigma$-measurable and the sets $epi\,\varphi_{x,\theta}$ are non-empty and closed for all $(x, \theta)$.
\end{remark}

\begin{remark}
\label{rem: 21.2} 
$\mathcal{B} \times \Sigma$ -measurable mapping $G( x , \theta)$: $E_n \times \Theta \to 2^{E _m}$ is measurable in $\theta$ for fixed x due to the similar property measurable of functions [10, p. 126] and the Castaing representation of measurable set-valued mappings (Theorem 21.2).
\end{remark}
\begin{definition}
\label{def:21.11}
{The (Lebesgue) integral of a set-valued mapping} $G : \Theta \to 2^{E _n}$ is the set of integrals of all possible integrable sections of the mapping $G$:
\begin{equation} \notag
	\underset{\Theta}{\int} G (\theta) P (d\theta) = \left\{ \underset{\Theta}{\int} g (\theta) P (d\theta) \;|\; g (\theta) \in G (\theta) \right\}.
\end{equation}
\end{definition}
\begin{theorem}
\label{th:21.3}
[144]. {Let $(\Theta, \Sigma, P)$ be a measurable space with finite positive measure $P$, a multivalued mapping $G: \Theta \to 2^{E _n}$ is normal;  $sup \{ \| g\| \:|\: g \in G (\theta) \} \leq L (\theta)$, where $L (\theta)$ is an integrable function. Then the set $\int_\Theta G (\theta) \text{P} (\theta)$ is compact in $E _n$, and the set $\int_\Theta\mbox{co } G (\theta) \text{P} (\theta)$ is convex and compact in $E _n$, and if the measure $P$ is continuous, then}
\begin{equation} \notag
	\underset{\Theta}{\int} G (\theta) P (d\theta) = \underset{\Theta}{\int} co\,G (\theta) P (d\theta).
\end{equation}
\end{theorem}
\begin{lemma}
\label{lem:21.7}
{Let $(\Theta, \Sigma, P)$ be a measurable space with finite positive measure $P$, a multivalued mapping $G (\text{x}, \theta)$: $E _n \times \Theta \to 2^ {E _m}$ is upper semicontinuous in x in K $\subset \text{E} _n$ for all $\theta$ and is measurable in $\theta$ for every x $\in K$, the sets $G(x, \theta)$ are non-empty compact sets for all $(x, \theta$), $sup \{ \|g\| \:|\: g \in G ( x , \theta), x \in K \} \leq L _K (\theta)$, where $L _K (\theta)$ is an integrable function. Then the sets}
\begin{equation}  \notag
	G (x) = \underset{\Theta}{\int} G (x, \theta) P (d\theta)
\end{equation}
{are non-empty compacts (i.e. bounded closed sets)
and the multivalued mapping $x \to G(x)$ is upper semicontinuous in $K$.}
\end{lemma}

{\it P r o o f}. The compactness of $G(x)$ follows from Theorem 21.3. Let us prove the upper semicontinuity of the mapping $G$. Denote
\begin{equation} \notag
	\Delta (B, A) = \underset{b \in B}{\sup} \: \underset{a \in A}{\inf} \: \|\: b - a \:\|,
	\;\;A, B \subset E _n .
\end{equation}


It is necessary to show that $\underset{y \to x}{\lim} \Delta (G(y), G(x)) = 0$. By Definition 21.11, for any $g _z \in G (z)$ there exists a measurable section $g _z (\theta)$ of $G ( z , \cdot): \Theta \to 2^{E _m}$ such that $g _z = \underset{\Theta}{\int} g _z (\theta) P (d\theta)$. That's why
\begin{equation} \tag{21.1}
\begin{array}{c}
	\Delta (G ( y ), G (x)) = \underset{g _y \in G (y)}{\sup}\;\; \underset{g _x \in G (x)}{\inf} \,\|\, g _y - g _x \,\|\, \leq \\
	\leq \underset{g _y (\theta) \in G (y, \theta)}{\sup}\;\; \underset{g _x (\theta) \in G (x , \theta)}{\inf} \: \underset{\Theta}{\int} \,\|\, g _y (\theta) - g _x (\theta) \,\|\, P (d\theta),
\end{array}
\end{equation}
where $g _x (\theta), g _y (\theta)$ are measurable (and integrable) sections of the maps $G (x, \cdot), G (y, \cdot)$, respectively. Let us show that in (21.1) inf we can introduced under the integral sign, i.e.,
\begin{equation} \tag{21.2}
	\underset{g _x (\theta) \in G (x , \theta)}{\inf} \: \underset{\Theta}{\int} \,\|\, g _y (\theta) - g _x (\theta) \,\|\, P (d\theta)
	= \underset{\Theta}{\int} \underset{g _x (\theta) \in G (x , \theta)}{\inf} \,\|\, g _y (\theta) - g _x (\theta) \,\|\, P (d\theta).
\end{equation}

Let us show that in (21.2) the left-hand side is not less than the right-hand side. For an integrable function $g _x (\theta) \in G (x , \theta)$,
\begin{equation} \tag{21.3}
\|\, g _y (\theta) - g _x (\theta) \,\|\, \geq \underset{g _x \in G (x , \theta)}{\inf} \|\, g _y (\theta) - g _x  \,\|,
\end{equation}
moreover, the right hand side of (21.3) is measurable and integrable. Indeed, let $\{ g^{m}_x (\theta), \; m = 1,2,\dots \}$ be a countable dense family of sections of $G(x, \cdot)$. Due to the continuity of $\|\, g _y (\theta) - g _x (\theta) \,\|$ in g$_x$,
\begin{equation} \tag{21.4}
\begin{array}{lcl}
\underset{g _x (\theta) \in G (x , \theta)}{\inf} \,\|\, g _y (\theta) - g _x (\theta) \,\|
&=& \underset{m}{\inf} \,\|\, g_y (\theta) - g^{m}_x (\theta)\,\|  \\
&=&\underset{M \to \infty}{\lim} \, \underset{0 \leq m \leq M}{\min} \,\|\, g _y (\theta) - g^{m}_x (\theta )\,\|.
\end{array}
\end{equation}

The function $\min_{\{0 \leq m \leq M\}} \|\, g _y (\theta) - g^{m}_x (\theta) \,\|$ is measurable, since the Lebesgue sets (with $c\in E_1$)
\begin{equation} \notag
\{ \theta \,|\, \underset{0 \leq m \leq M}{\min} \,\|\, g _y (\theta) - g^{m}_x (\theta) \,\| \, \leq c\}
= \underset{0 \leq m \leq M}{\cup} \:
\{ \theta \,| \: \|\, g_y (\theta) - g^{m}_x (\theta) \,\|\, \leq c \}
\end{equation}
are measurable.

Then the left hand side of (21.4) is measurable as the limit of measurable functions. Besides,
\begin{equation} \notag
\underset{g _x (\theta) \in G (x , \theta)}{\inf} \,||\, g _y (\theta) - g^{m}_x (\theta) \,||\, \leq 2 L _K (\theta).
\end{equation}

Thus, the right hand side in (21.3) is integrable. We integrate (21.3) and take inf over g$_x (\theta) \in G (x, \theta)$ on the left-hand side; then we get
\begin{equation} \tag{21.5}
	\underset{g _x (\theta) \in G (x , \theta)}{\inf} \underset{\Theta}{\int} \,\|\, g _y (\theta) - g^{m}_x (\theta) \,\|\, P (d\theta)
	 \geq \underset{\Theta}{\int} \underset{g _x (\theta) \in G (x , \theta)}{\inf} \|\, g _y (\theta) - g^{m}_x (\theta) \,\|\, P (d\theta).
\end{equation}


We now prove the reverse inequality. Let us show that inf on the right-hand side of (21.5) is actually reached not pointwise, but on measurable sections of $G(x,\cdot)$. Let $\{g^{m}_x (\theta), m = 1,2,\dots\}$ be a countable dense family of measurable sections of $G(x,\cdot)$. Take an arbitrary $\varepsilon > 0$. Consider the set
\begin{equation} \notag
\Theta^{m}_\varepsilon = \{ \theta \,| \: \|\, g _y (\theta) - g^{m}_n (\theta) \,\|\, \leq \underset{g _x \in G (x, \theta)}{\inf} \,\|\, g _y (\theta) - g _x (\theta) \,\|\, + \varepsilon\}.
\end{equation}

We define the function $h^{m}_{x, \varepsilon}: \Theta \to E _n$ as follows:
\begin{equation} \notag
h^{m}_{x, \varepsilon} =
\begin{Bmatrix}
g^{1}_{x}, & \theta \in \Theta^{1}_{\varepsilon}, \\
g^{2}_{x}, & \theta \in \Theta^{2}_{\varepsilon} \setminus \Theta^{1}_{\varepsilon}, \\
\dots & \dots \\
g^{m}_{x}, & \theta \in \Theta^{m}_{\varepsilon} \setminus \cup_{i=1}^{m-1} \Theta^{i}_{\varepsilon} \\
g^{1}_{x}, & \theta \in \Theta \setminus \cup_{i=1}^{m} \Theta^{i}_{\varepsilon}.
\end{Bmatrix}
\end{equation}
Obviously, the function $h^{m}_{x, \varepsilon}$ is a measurable section of $G(x,\cdot)$, and on the set $T^{m}_{\varepsilon} = \overset{m} {\underset{i=1}\cup} \Theta^{i}_{\varepsilon}$ holds
\begin{equation} \notag
||\, g _y (\theta) - h^{m}_{x, \varepsilon} (\theta) \,||\, \leq \underset{g _x \in G (x, \theta)} {\inf} \,||\, g _y (\theta) - g _x (\theta) \,||\, + \varepsilon.
\end{equation}
Let us show that $P^{m}_{\varepsilon} = P (\theta \setminus T^{m}_{\varepsilon}) \to 0 \;(m \to \infty)$. Indeed, $\{T^{m}_{\varepsilon}\}_{m \geq 1}$ is a monotonically increasing sequence of measurable sets, and $P (T^{m}_{\varepsilon}) \leq P (\Theta) < +\infty$. Therefore, the limit set $T^{\infty}_{\varepsilon} = \overset{\infty}{\underset{i=1}\cup} \Theta^{i}_{\varepsilon}$ is measurable and $P^ {\infty}_{\varepsilon} = P(T^{\infty}_{\varepsilon}) \leq P(\Theta)$. If $P^{\infty}_{\varepsilon} < P(\Theta)$, then on the set $\Theta \setminus T^{\infty}_{\varepsilon}: P(\Theta \setminus T^{ \infty}_{\varepsilon}) > 0$ for all $m$ holds
\begin{equation} \notag
\|\, g _y (\theta) - g^{m}_{x} (\theta) \,\|\, > \underset{g _x \in G (x, \theta)}{\inf} \,\|\, g _y (\theta) - g _x (\theta) \,\|\, + \varepsilon.
\end{equation}

This means that the family $\{g^{m}_{x} (\theta)\}^{\infty}_{m=1}$ on the set $\Theta \setminus T^{\infty}_{ \varepsilon}$ is not dense in $G(x,\theta)$. The resulting contradiction proves that $P^{m}_{\varepsilon} = P(\Theta \setminus T^{\infty}_{\varepsilon}) \to 0 \;(m \to \infty)$. Then
\begin{equation} 
\begin{array}{c}
\underset{g _x \in G (x, \theta)}{\inf} \underset{\Theta}{\int} ||\, g _y (\theta) - g^{m}_{x} (\theta) \,||\, P (d\theta) \leq \underset{\Theta}{\int} ||\, g_y (\theta) - h^{m}_{x, \varepsilon} ( \theta) \,||\, P (d\theta) = \\
= (\underset{T^{m}_{\varepsilon}}{\int} + \underset{\Theta \setminus T^{m}_{\varepsilon}}{\int}) \,||\, g _y (\theta) - h^{m}_{x, \varepsilon} (\theta) \,||\, P (d\theta) \leq \\
\leq (\underset{\Theta}{\int} - \underset{\Theta \setminus T^{m}_{\varepsilon}}{\int}) \underset{g _x \in G (x, \theta)}{\inf} ||\, g _y (\theta) - g_{x} (\theta) \,||\, P (d\theta)   \\
+ \varepsilon+ \underset{\Theta \setminus T^{m}_{\varepsilon}}{\int} \,||\, g _y (\theta) - h^{m}_{x, \varepsilon} (\theta) \,||\, P (d\theta);\tag{21.6}
\end{array}
\end{equation}
moreover, the residual integrals on the right-hand side of (21.6) over the sets $\Theta \setminus T^{\infty}_{\varepsilon})$ tend to zero $m \to +\infty$ due to the absolute continuity of the Lebesgue integral. Passing in (21.6) to the limit first in $m \to$ $+\infty$ and then in $\varepsilon \to 0$, we obtain
\begin{equation} \tag{21.7}
	\underset{g _x \in G (x, \theta)}{\inf} \underset{\Theta}{\int} ||\, g _y (\theta) - g^{m}_{x} (\theta) \,||\, P (d\theta) 
	\leq
	\underset{\Theta}{\int} \underset{g _x \in G (x, \theta)}{\inf} ||\, g _y (\theta) - g_{x} ||\, P (d\theta).
\end{equation}
From (21.5) and (21.7) the required equality (21.2) follows. From (21.1) and (21.2) we get
\begin{equation*} \tag{21.8}
\Delta(G(y), G(x)) \leq \sup _{g_{y}(\theta) \in G(x, \theta)} \int_{\Theta}\;\; \inf _{g_x \in G(x, \theta)}\left\|g_{y}(\theta)-g_{x}\right\| P(d\theta).
\end{equation*}
It can proved exactly in the same way that sup can be introduced under the integral sign in (21.8). It is only necessary to show that the function $\inf \{\|g_{y}-g_{x}\|\,|\, g_{x} \in G(x, \theta)\}$ is continuous in $g_{y}$. Indeed,
\begin{equation*} \notag
\begin{array}{c}
\inf _{g_{x} \in G(x, \theta)} \,||\,g_{y}^{2}-g_{x}\,|| - ||\,g_{y}^{1}-g_{y}^{2}\,||\, \leq \\
\leq \inf _{g_{x} \in G(x, \theta)}\|g_{y}^{1}-g_{x}\| \leq \inf _{g_{x} \in G(x, \theta)}\|g_{y}^{2}-g_{x}\|+\|g_{y}^{1}-g_{y}^{2}\| .
\end{array}
\end{equation*}
So,
\begin{equation*} \notag
\Delta(G(y), G(x)) \leq \int_{\theta} \Delta(G(y, \theta), G(x, \theta)) P(d \theta) .
\end{equation*}
Since $\Delta(G(y, \theta), G(x, \theta)) \leq 2 L_{K}(\theta), L_{K}(\theta)$ is integrable and $\Delta(G (y, \theta), G(x, \theta)) \to 0$ for $y \to x(\theta \in \Theta)$, then by Lebesgue's theorem on passage to the limit under the integral sign, $\Delta( G(y), G(x)) \to 0$ for $y \to x$. The lemma is proven.

\section*{$\S$ 22. Random Lipschitz functions}
\label{Sec.22}
\setcounter{section}{22}
\setcounter{lemma}{0}
\setcounter{theorem}{0}
\setcounter{corollary}{0}
\setcounter{remark}{0}
\setcounter{equation}{0}
In this section, we study the Clarke subdifferential of random Lipschitz functions. It is shown that it is a multi-valued mapping measurable in the set of deterministic and random variables. This result is used to prove the measurability of the trajectories of some stochastic optimization methods. A formula related to subdifferential calculus is obtained, namely, that the subdifferential of the mathematical expectation of a random Lipschitz function is embedded in the mathematical expectation of the subdifferential of this function (which is understood as an integral of a multivalued mapping).


Note that for convex functions, the Clarke subdifferential coincides with the subdifferential of a convex function (the set of subgradients). Thus, we have the inclusion of the subdifferential of the mathematical expectation of a random convex function in the mathematical expectation of its subdifferential. The reverse inclusion (for the proven measurability of the subdifferential) follows trivially from the subgradient inequality for convex functions. Thus, in the convex case, we have a well-known result: the equality of the subdifferential of the mathematical expectation and the mathematical expectation of the subdifferential. The proof of this fact for functions in abstract spaces, information of a historical and bibliographic nature can be found in [48].
\begin{lemma}
\label{lem:22.1}
{Let $(\Theta, \Sigma)$ be a measurable space, the function f: $E_{n} \times \Theta \to E_{1}$ be measurable in $\theta \in \Theta$ for every $x \in E_{n}$ and is Lipschitzian in $x$ in every compact set $K \subset E_{n} $ for all $\theta \in \Theta,$ $\partial_{x} f(x, \theta)$ be Clarke's subdifferential of the function $f(\cdot, \theta)$ at point x. Then the multivalued mapping $\partial_{x} f: E_{n} \times \Theta \to 2^{E_{n}}$ is $\mathcal{B} \times \Sigma$-measurable, and therefore has $\mathcal{B} \times \Sigma$-measurable sections. The mapping $\partial_{x} f(x, \cdot): \theta \to 2^{E_{n}}$ is measurable in $\theta$ for fixed x.}
\end{lemma}
{\it P r o o f}. Recall that
\begin{equation} \notag
\begin{array}{c}
\partial f(x, \theta) = \{g \in E_{n} \mid\left<g, d\right> \leq f^{0}(x, \theta ; d)\; \forall\; d \in E_{n} \}, \\
f^{0}(x, \theta; d) =\underset{\delta\rightarrow 0}{\lim} \;\;\underset{\|y\|\leq\delta, 0<\lambda\leq\delta}{\sup}
[f(x + y + \lambda d, \theta) - f(x + y, \theta)] / \lambda,
\end{array}
\end{equation}
where $f^{0}(x, \theta ; d)$ is the Clarke derivative of the function $f(\cdot, \theta)$ in the direction $d \in E_{n}$. Let us show that $f^{0}(\cdot, \cdot ; d)$  is $\mathcal{B}\times \Sigma$-measurable. To do this, it suffices to show that  the function
\begin{equation} 
\label{eqn:22.1}
	f_{\delta}^{0}(x, \theta ; d)=\sup _{\substack{\|y\| \leq \delta \\ 0<\lambda \leq \delta}}[f(x+y+\lambda d, \theta)-f(x+y, \theta)] / \lambda 
\end{equation}
is  $\mathcal{B}\times \Sigma$-measurable.

Let $\{y_{l}\}$ and $\{\lambda_{m}\}$ be countable everywhere dense sets in
\begin{equation} \notag
\{y \in E_{n}|\;\| y\| \leq \delta\} \text{  and  }\{\lambda \in E_{1} \mid 0<\lambda \leq \delta\}.
\end{equation}
Then
\begin{equation} \notag
\begin{aligned}
	f_{\delta}(x, \theta ; d) & =\sup _{\substack{\| y \| \leq \delta \\
	0<\lambda \leq \delta}}[f(x+y+\lambda d, \theta)-f(x+y, \theta)] / \lambda= \\
	& =\lim _{\varepsilon \to+0} \sup _{\substack{\|y\| \leq \delta \\
	\varepsilon \leq \lambda \leq \delta}}[f(x+y+\lambda d, \theta)-f(x+y, \theta)] / \lambda= \\
	= & \lim _{\varepsilon \to+0} \sup _{\substack{\|y_{l}\| \leq \delta \\
	\varepsilon \leq \lambda_{m} \leq \delta}}[f(x+y_{l}+\lambda_{m} d, \theta)-f(x+y_{l}, \theta)] / \lambda_{m}= \\
	= & \sup _{l, m}[f(x+y_{l}+\lambda_{m} d, \theta)-f(x+y_{l}, \theta)] / \lambda_{m}= \\
	& =\lim _{N \to \infty} \underset{l,m \leq N}{\max} [f(x+y_{l}+\lambda_{m} d, \theta)-f(x+y_{l}, \theta)] / \lambda_{m} .
\end{aligned}
\end{equation}


Note that the functions $f(x+y_{l}+\lambda_{m} d, \theta)$ and $f(x+y_{l}, \theta)$ are Carathéodorian on $E_{n} \times \Theta$ and hence  are $\mathcal{B}\times \Sigma$-measurable (Lemma 21.6). Function
\begin{equation} \notag
\underset{l, m \leq N}{\max} [f(x+y_{l}+\lambda_{m} d, \theta) - f(x+y_{l}, \theta)] / \lambda_{m},
\end{equation}
obviously,  is $\mathcal{B}\times \Sigma$-measurable. Then $f_{\delta}^{0}(\cdot, \cdot; d)$ is measurable as the limit of measurable functions. Consequently, $f^{0}(\cdot, \cdot; d)$ is $\mathcal{B} \times \Sigma$-measurable as well. Note that $f^{0}(x, \theta; \cdot)$ is continuous. Let us represent
\begin{equation} \notag
\partial_{x} f(x, \theta)=\{g \in E_{n} \mid \psi(g, x, \theta) \equiv \max _{\|d\| \leq 1}[\left<g, d\right>-f^{0}(x, \theta ; d)] \leq 0\} .
\end{equation}
The function $\psi(g, \cdot, \cdot)$ is $\mathcal{B} \times \Sigma$-measurable for each $g$ (similarly to $f_{\delta}^{0}(\cdot, \cdot ; d))$, is convex and continuous in $g$; by Lemma 21.4, $\psi(g, x, \theta)$ is a normal convex integrand on $E_{n} \times(E_{n} \times \Theta)$. Now, by Lemma 21.5, the multivalued mapping
\begin{equation} \notag
(x, \theta) \to \partial_{x} f(x, \theta)=\{g \in E_{n} \mid \psi(g, x, \theta) \leq 0\}
\end{equation}
is $\mathcal{B} \times \Sigma$-measurable and, by Theorem 21.2, has $\mathcal{B} \times \Sigma$-measurable sections. By Remark 21.2, the measurability of $\partial_{x} f(x, \theta)$ in two variables $(x, \theta)$ implies the measurability in one variable $\theta$. The lemma is proven.
\begin{corollary}
\label{cor:22.1}
{If  function $f(x, \theta)$ is convex in $x$ for each $\theta$ and measurable in $\theta$ for each $x$, then its subgradient mapping (it coincides with the Clarke subdifferential) $ \partial_{x} f(x, \theta)$  is $\mathcal{B}\times \Sigma$-measurable, and the mapping $\partial_{x} f(x, \cdot)$ is measurable.}
\end{corollary}
\begin{theorem}
\label{th:22.1}
{Let $(\Theta, \Sigma, P)$ be a space with finite positive measure (for example, a probability space), the function $f: E_{n} \times \Theta \to E_{1}$ be measurable in $ \theta \in \Theta$ for every $x \in E_{n}$ and Lipschitz in $x$ in every compact set $K \subset E_{n}$ with a finite Lipschitz constant $L_{K} (\theta)$ for all $\theta \in \Theta$, $\partial_{x} f(x, \theta)$ be the Clarke subdifferential of the function $f(\cdot, \theta)$ at the point $x$. Then $\partial_{x} f(x, \cdot): \Theta \to 2^{E_{n}}$ is a convex-valued normal map. If, in addition, $f(x, \theta)$ and $L_{K}(\theta)$ are integrable, then the function $F(x)=\int_{\theta} f(x, \theta) P( d \theta)$ is Lipshitzian in $K$ with constant $L_{K}=\int_{\Theta} L_{K}(\theta) P(d \theta)$ and the inclusion
holds}
\begin{equation} \tag{22.2}
\partial F(x) \subset \int_{\Theta} \partial_{x} f(x, \theta) P(d \theta).
\end{equation}
\end{theorem}
{\it P r o o f}. Recall that the sets $\partial_{x} f(x, \theta)$ are convex and compact $(\S 2)$. The measurability of the mapping $\partial_{x} f(x, \cdot): \Theta \to 2^{E_{n}}$ was proved in Lemma 22.1. Hence the mapping $\partial_{x} f(x, \cdot)$ is normal.

The fact that $F(x)$ is Lipschitzian is obvious: for $x, y \in K$
\begin{equation} \notag
|F(y)-F(x)|\leq \int_{\Theta}| f(y, \theta)-f(x, \theta) \mid P(d \theta) \leq\|y-x\| \int_{\Theta} L_{K}(\theta) P(d \theta).
\end{equation}


Recall that
\begin{equation} \tag{22.3}
\begin{gathered}
\partial F(x)=\{g \in E_{n} \mid\left<g, d\right> \leq F^{0}(x ; d)\; \forall\; d \in E_{n}\}, \\
F^{0}(x ; d)=\lim _{\delta \to+0} \sup _{\substack{\|<\|y\|\leq \delta \\
0<\lambda \leq \delta}}[f(x+y+\lambda d)-f(x+y)] / \lambda .
\end{gathered}
\end{equation}

Let's prove the inequality
\begin{equation} \tag{22.4}
F^{0}(x ; d) \leq \int_{\Theta} f^{0}(x ; \theta) P(d \theta).
\end{equation}
Let $f_{\delta}^{0}(x, \theta; d)$ be defined by formula (22.1). For $y, \lambda$ such that $\|y\| \leq \delta$ and $0<\lambda \leq \delta$,
\begin{equation} \notag
f_{\delta}^{0}(x, \theta ; d) \geq [f(x+y+\lambda d, \theta)-f(x+y, \theta)] / \lambda.
\end{equation}
The functions $f_{\delta}^{0}$ and $f^{0}$  are $\mathcal{B}\times \Sigma$-measurable and hence measurable (see Remark 21.2). In addition, $|f_{\delta}^{0}(x, \theta ; d)| \leq L_{K}(\theta)$ and $|f^{0}(x, \theta ; d)| \leq L_{K}(\theta)$; so they are integrable. Then
\begin{equation} \notag
\begin{gathered}
\int_{\Theta} f_{\delta}^{0}(x, \theta ; d) P(d \theta) \geq [f(x+y+\lambda d)-f(x+y)] / \lambda, \\
\int_{\Theta} f_{\delta}^{0}(x, \theta ; d) P(d \theta) \geq \sup _{\substack{\| y\| \leq \delta \\
0<\lambda \leq \delta}}[f(x+y+\lambda d)-f(x+y)] / \lambda .
\end{gathered}
\end{equation}
In the last inequality, we pass to the limit with respect to $\delta \to 0$. By virtue of Lebesgue's theorem, the limit can be introduced under the integral sign; as a result, we obtain inequality (22.4). Denote $\mathbb{E} \partial_{x} f(x, \theta)=\int_{\Theta} \partial_{x} f(x, \theta) P(d \theta)$ and $\mathbb{E} f^{0}( x, \theta ; d)=\int_{\Theta} f^{0}(x, \theta ; d) P(d \theta)$.

Let's prove the equality
\begin{equation} \tag{22.5}
	\mathbb{E}\partial_{x} f(x, \theta)=\{g \in E_{n} \mid\left<g, d\right> \leq \mathbb{E}f^{0}(x, \theta ; d) \quad \forall d \in E_{n}\}.
\end{equation}

First we prove that the left-hand side of (22.5) is embedded in the right-hand side. Let $g \in \mathbb{E}\partial_{x} f(x, \theta)$. According to Definition 21.10, there exists an integrable function $g(\theta)$ such that $g(\theta) \in \partial_{x} f(x, \theta)$ and $g=\int_{\theta} g( \theta) P(d \theta)$. Then by the definition of $\partial_{x} f(x, \theta)$, $\left<g(\theta), d\right> \leq f^{0}(x, \theta ; d)$ for any $d \in E_{n}$. Integrating this inequality, we get
\begin{equation} \notag
\left<g, d\right> \leq \mathbb{E} f^{0}(x, \theta ; d) \quad \forall d \in E_{n} .
\end{equation}

Consequently, $g$ belongs to the right-hand side of (22.5), which was to be proved.

Now suppose that equality (22.5) does not hold, i.e., there exists $g^{\prime} \notin \mathbb{E} \partial_{x} f(x, \theta)$ and
\begin{equation} \notag
g^{\prime} \in\{g \in E_{n} \mid \left<g, d\right> \leq \mathbb{E} f^{0}(x, \theta ; d) \quad \forall d \in E_{n }\}.
\end{equation}
By Theorem 21.3, $\mathbb{E} \partial_{x} f(x, \theta)$ is a convex compact set; so there exists $g^{\prime \prime} \in \mathbb{E} \partial_{x} f(x, \theta)$ such that
\begin{equation} \notag
\|g^{\prime}-g^{\prime \prime}\|=\rho(g^{\prime}, \mathbb{E} \partial_{x} f(x, \theta)).
\end{equation}
Let $d=g^{\prime}-g^{\prime \prime}$. By virtue of the above assumption, for any $g \in \mathbb{E} \partial_{x} f(x, \theta)$
\begin{equation} \tag{22.6}
\left<g, d\right><\left<g^{\prime}, d\right> \leq \mathbb{E} f^{0}(x, \theta ; d).
\end{equation}
Let us show that this cannot be true. Let's define
\begin{equation} \notag
G_{d}(x, \theta)=\{g \in E_{n} \mid g \in \partial_{x} f(x, \theta),\;\left<g, d\right>=f^{0}( x, \theta ;d)\}.
\end{equation}
Let us show that $\theta \to G_{d}(x, \theta)$ is a measurable mapping. Indeed,
\begin{eqnarray*} 
G_{d}(x, \theta) & =&\left\{\right. g \in E_{n} \mid \psi(g, x, \theta)=\sup _{\| l \|\leq 1}[\left<g, l\right>-f^{0}(x, \theta ; l)] \leq 0, \\
&&\;\;\;\;\;\;\;\;\;\;\;\;\;\;\;\;\;\;\;\;\;\;\;\;\psi_{d}(g, x, \theta) \equiv  \left<g, d\right>-f^{0}(x,\theta;d)=0\}= \\
& =&\{g \in E_{n} \mid \varphi(g, x, \theta)=\max (\psi(g, x, \theta),-\psi_{d}(g, x, \theta)) \leq 0\} .
\end{eqnarray*}

The functions $\psi(g, x, \theta),-\psi_{d}(g, x, \theta), \varphi(g, x, \theta)$ are convex and continuous in $g$ and measurable in $ \theta$, since the functions $\left<g,l\right>-f^{0}(x, \theta ; l)$ have this property, and the operations sup and max do not violate these properties. Hence $\varphi(\cdot, x, \cdot)$ is a normal convex integrand on $E_{n} \times \Theta$; by Lemma 21.5, the mapping $G_{d}(x, \cdot)$ is measurable. By Theorem 21.2, there exists a measurable (and integrable) section $g(\theta) \in G_{d}(x, \theta)$, i.e., $g(\theta) \in \partial_{x} f (x, \theta)$ and $\left<g(\theta), d\right>=f^{0}(x, \theta ; d)$. Integrating these expressions, we get
\begin{equation} \notag
\begin{gathered}
	g=\int_{\Theta} g(\theta) P(d \theta) \in  \partial_{x} f(x, \theta), \\
	\left<g, d\right>=\left<\int_{\Theta} g(\theta) P(d \theta), d\right>=\mathbb{E} f^{0}(x, \theta ; d),
\end{gathered}
\end{equation}
which contradicts (22.6). Thus, equality (22.5) is proved. Now from (22.3) - (22.5) the required inclusion (22.2) follows.
\begin{remark}
\label{rem:22.1}
Inequality (22.4) was proved in [180], and the inclusion (22.2) was established under the following assumptions:
\begin{enumerate}
   \item There are directional derivatives $d \in E_{n}$ of functions $f(\cdot, \theta):$
\end{enumerate}
$$
f^{\prime}(x, \theta ; d)=\lim _{\lambda \to+0}[f(x+\lambda d, \theta)-f(x, \theta)] / \lambda
$$
\begin{enumerate}
   \setcounter{enumi}{1}
   \item It holds $f^{\prime}(x, \theta ; d)=f^{0}(x, \theta ; d)$.
\end{enumerate}
Since the Clarke derivative $f^{0}(x, \theta ; d)$ is convex in $d$, it is actually assumed in [180] that the directional derivative $f^{\prime}(x, \theta ;)$ is convex in $d$, i.e., the local convexity of $f(\cdot, \theta)$ [49].
\end{remark}

\section*{$\S$ 23. Random generalized differentiable functions 
and the calculus of stochastic generalized gradients}
\label{Sec.23}
\setcounter{section}{23}
\setcounter{definition}{0}
\setcounter{equation}{0}
\setcounter{theorem}{0}
\setcounter{lemma}{0}
\setcounter{remark}{0}
\setcounter{corollary}{0}

In this section, we study random generalized differentiable functions and their random pseudogradient mappings. It is shown that the class of generalized differentiable functions is closed under the operation of taking the expectation. In this case, the expectation of a pseudo-gradient mapping (as an integral of a multi-valued mapping) of a random function is itself some pseudo-gradient representation of the mathematical expectation of this function. In this respect, generalized differentiable functions are similar to convex functions.

Due to these properties, and also due to the stability of the generalized gradient method, it turns out that it is possible to extend the method of stochastic generalized gradients (stochastic quasi-gradients [41]) to stochastic programming problems with generalized differentiable functions. This will be done in the next paragraph. An important role is played by the stochastic generalized gradients of the problem functions, which are measurable in the set of deterministic and random variables sections of random pseudogradient mappings. In this section, we study the existence and construction of such sections.

\textbf{1. Random generalized differentiable functions}.
\label{Sec.23.1}
\begin{theorem}
\label{th:23.1}
{Let $(\Theta, \Sigma, P)$ be a measurable space with a finite positive measure (for example, a probability space). The function $f: E_{n} \times \Theta \to E_{1}$ is generalized differentiable with respect to $x \in E_{n}$ for $\theta \in \Theta$ and integrable with respect to $\theta$ for each $ x$, $G_{f}: E_{n} \times \Theta \to 2^{E_{n}}$ is its $\theta$-measurable for each $x$ pseudo-gradient (in $x$) for each $\theta$ mapping ( for example, $G_{f}(x, \theta)=\partial_{x} f(x, \theta)$ ). Let for each compact set $K \subset E_{n}$ there exists an integrable function $L_{K}(\theta)$ such that}
\begin{equation} \notag
	\sup \{\|g\||\, g \in G_{f}(x, \theta), x \in K\} \leq L_{K}(\theta) .
\end{equation}
{Then function $F(x)=\int_{\Theta} f(x, \theta) P(d \theta)$ is generalized differentiable and $x \to G_{F}(x)=\int_{\Theta } G_{f}(x, \theta) P(d \theta)$ is some pseudogradient mapping for F(x).}
\end{theorem}

{\it P r o o f}. It needs to be shown

firstly, that the sets $G_{F}(x)$ are non-empty, convex and compact, this follows from Theorems 21.2, 21.3;

second, that the mapping $G_{F}: E_{n} \to 2^{E_{n}}$ is upper semicontinuous, as follows from Lemma 21.7;

thirdly, that the expansion holds true:
\begin{equation} \notag
	F(y)=F(x)+\left<g, y-x\right>+o(x, y, g),
\end{equation}
where $o(x, y^{k}, g^{k}) /\|y^{k}-x\| \to 0$ for any sequences $y^{k} \to x$ and $g^{k} \in G_{F}(y^{k})$.

By Definition 21.11, for each $g \in G_{F}(y)$ there exists an integrable section $g_{y}(\theta)$ of $G_{f}(y, \cdot)$ such that $g=\mathbb{E} g_{y}(\theta)=\int_{\theta} g_{y}(\theta) P(d \theta).$ Due to the generalized differentiability of $f(\cdot, \theta)$ the expansion holds true:
\begin{equation} \tag{23.1}
f(y, \theta)=f(x, \theta)+\left<g, y-x\right>+o(x, y, g ; \theta) \text {, }
\end{equation}
where $o(x, y, g ; \theta) /\|y-x\| \to 0$ for arbitrary $y \to x$ and $g \in G_{f}(y, \theta)$. Substitute the section $g_{y}(\theta)$ into (23.1) :
\begin{equation} \tag{23.2}
	f(y, \theta)=f(x, \theta)+\left<g_{y}(\theta), y-x\right>+o(x, y, g_{y}(\theta) ; \theta)
\end{equation}
Obviously, $o(x, y, g_{y}(\theta) ; \theta)$ is integrable under integrable $g_{y}(\theta)$. Integrating (23.2), we obtain
\begin{equation} \notag
o(x, y, g)=F(y)-F(x)-\left<g, y-x\right>=\int_{\Theta} o(x, y, g_{y}(\theta) ; \theta) P (d\theta),
\end{equation}
where $g=\mathbb{E} g_{y}(\theta)$.

Let now arbitrary sequences $y^{k} \to x$ and $g^{k} \in G_{F}(y^{k})$ be given. Let us show that $o(x, y^{k}, g^{k}) /\|y^{k}-x\| \to 0.$ Indeed, by the definition of the integral, for each $g^{k} \in G_{F}(y^{k})$ there exists an integrable section $g^{k}(\theta)$ of the mapping $G_{f}(y ^{k}, \cdot)$ such that $g^{k}=\int_{\Theta} g^{k}(\theta) P(d \theta)$. Let $\|y^{k}-x\| \leq \varepsilon$. Since $|o(x, y^{k}, g^{k}(\theta) ; \theta)| /\|y^{k}-x\| \leq L_{\varepsilon}(\theta)$, where $L_{\varepsilon}(\theta)$ is some integrable function, and $o(x,y^{k}, g^{k} (\theta) ; \theta) /\|y^{k}-x\| \to 0$ for $k \to \infty$ for all $\theta$, then by Lebesgue's theorem on passage to the limit under the integral sign,
\begin{equation} \notag
	\lim _{k \to \infty} \frac{o(x, y^{k}, g^{k})}{\|y^{k}-x\|}=\int_{\Theta} \lim _{k \to \infty} \frac{o(x, y^{k}, g^{k}(\theta) ; \theta)}{\|y^{k}-x\|} P(d \theta)=0 .
\end{equation}
The theorem has been proven.

\textbf{2. Calculus of stochastic generalized gradients.} 
\label{Sec.23.2}
Consider the problem of stochastic programming
\begin{equation} \tag{23.3}
	F(x) \equiv \int_{\Theta} f(x, \theta) P(d \theta) \to \min_x,
\end{equation}
where $f(x, \theta)$ satisfies the conditions of Theorem 23.1. Let $G_{f}(x, \theta)$ be the pseudogradient (in $x$) mapping required under the conditions of Theorem 23.1 for the function $f(\cdot, \theta)$. For example, according to Theorem 1.10 and Lemma 22.1, it can be the Clarke subdifferential $\partial_{x} f(x, \theta)$; further examples of suitable $G_{f}(x, \theta)$ will be given. Let $g(x, \theta)$ be a  $\mathcal{B}\times \Sigma$-measurable section of $G_{f}(x, \theta)$. Under our assumptions, it exists because, by virtue of Lemma 22.1, Theorem 21.2, and Theorem 1.10, there exists a  $\mathcal{B}\times \Sigma$-measurable section $g(x, \theta) \in \partial_{x} f(x, \theta ) \subset G_{f} (x, \theta)$.

The simplest stochastic algorithm for solving problem (23.3) has the form
\begin{align*}
x^0 \in E_n, \quad x^{k+1}=x^k-\rho_k g(x^k, \theta^k), \quad k=0,1, \ldots, \tag*{(23.4)}
\end{align*}
where $x^0$ is the starting point, $x^k$ and $x^{k+1}$ are successive approximations to the solution of the problem, $\{\rho_k\}$ are step multipliers, $\{\theta^k\}$ are independent realization of the random parameter $\theta \in \Theta,$ $g\left(x^k, \theta^k\right)$ is the gradient of the integrand function $f\left(x, \theta^k\right)$ at the point $x^k$, which is a stochastic estimate of the gradient of the minimized function $F(x)$.

Let $\left(\Omega, \Sigma_\omega\right)$ denote a measurable space which is a direct countable product of measurable spaces $(\Theta, \Sigma)$. We denote by $\omega$ the sequence $\left\{\theta^k\right\}$; obviously, $\omega \in \Omega$. Let us define the mapping $\pi_k: \Omega \rightarrow \Theta$, which is the projector of the space $\Omega=\Theta \times \Theta \times \ldots$ onto the $k$-th multiplier $\Theta$, i.e., $\pi_k(\omega)=\theta^k$. We will consider the trajectory of points $\left\{x^k\right\}$ as functions of $\omega$. To analyze the trajectories $\left\{x^k(\omega)\right\} \quad(\omega \in \Omega)$ methods of random process theory can be applied, one must first be sure that $\left\{x^k(\omega)\right\}$ are measurable as functions on $\left(\Omega, \Sigma_\omega\right)$.

The measurability of the trajectory $\left\{x^k(\omega)\right\}$ of the algorithm (23.4) is established recurrently. The function $x^0(\omega)=x^0$ is measurable as a constant. If $x^k(\omega)$ is measurable, then $g^k(\omega)=g\left(x^k(\omega), \pi_k(\omega)\right)$ is measurable as a superposition of the measurable mapping $z^k=\left(x^k, \pi_k\right): \Omega \rightarrow E_n \times \Theta$ and $\mathcal{B} \times \Sigma$-measurable mapping $g: E_n \times \Theta \rightarrow E_n$. Hence, the vector-function $x^{k+1}(\omega)$ is also measurable.

Thus, when constructing stochastic optimization methods, the cross sections of generalized gradient mappings of functions under the sign of mathematical expectation, measurable jointly in deterministic and random variables, play an essential role. Let us consider the issues of existence and construction of such cross sections.

Let $(\Theta, \Sigma, P)$ be a probability space.

\begin{definition}
\label{def:23.1}
Function $f: E_n \times \times \Theta \rightarrow E_1$ is called random generalized differentiable (continuously differentiable, convex, Lipschitz), if it is measurable on $\theta \in E_n$ for every $x \in E_n$ and generalized differentiable (continuously differentiable, convex, Lipschitz) on $x \in E_n$ for all $\theta \in \Theta$.
\end{definition}
\begin{definition}
\label{def:23.2}
A stochastic pseudogradient (generalized gradient) of a random generalized differentiable function $f(x, \theta)$ is any $\mathcal{B} \times \Sigma$-measurable section $g(x, \theta)$ of its any pseudogradient mapping $G_f(x, \theta)$.
\end{definition}
\begin{lemma}
\label{lem:23.1}
Every random generalized differentiable function has stochastic pseudogradients (generalized gradients).
\end{lemma}

{\it P r o o f.} By Theorem 1.10, $\partial_x f(x, \theta) \subset G_f(x, \theta)$, and, by Lemma 22.1 , the mapping $\partial_x f(x, \theta)$ is $\mathcal{B} \times \Sigma$-measurable; therefore (Theorem 21.2) it has $\mathcal{B} \times \Sigma$-measurable sections.

\begin{lemma}
\label{lem:23.2}
{Let mapping $G_f(x, \theta)$ be $\mathcal{B} \times \Sigma$-measurable, then its
selection
$$
g(x, \theta)=\arg \min \left\{\|g\| \mid g \in G_f(x, \theta)\right\}
$$
is also $\mathcal{B} \times \Sigma$-measurable.}
\end{lemma}

The assertion of the lemma follows from the results of [175, Corollary 4.3].

\begin{lemma}
\label{lem:23.3}
{The gradient $g(x, \theta)$ of a random continuously differentiable function $f(x, \theta)$ is $ \mathcal{B} \times \Sigma$-measurable. If $f(x, \theta)$ is integrable over $\theta \in \theta$ for every $x$ and for any compact $K \subset E_n$ holds $\sup \{\| g(x$, $\theta) \| \mid y \in K\} \le L_K(\theta)$, where $L_K(\theta)$ is an integrable function, then the function $F(x) \equiv \int_{\Theta} f(x, \theta) P(d\theta)$ is continuously differentiable, $\quad g(x) $ $\equiv \int_{\Theta} g(x, \theta) P(d\theta)$ is its gradient. }
\end{lemma}

{\it P r o o f.} The gradient $g(x, \theta)$ is continuous on $x$ and measurable on $\theta$; by Lemma 21.6 it is $\mathcal{B} \times \Sigma$-measurable. Let $K=\{y\,|\,\| y-x \| \le \varepsilon\}$; by the lemma assumption $\sup \{\|g(x, \theta)\| \,|\, \| y-x \| \le \varepsilon\} \le L_K(\theta),$ and the function $L_K(\theta)$ is integrable. By the mean value theorem, the following inequalities hold for a continuously differentiable function $f(\cdot, \theta)$:
$$
\begin{gathered}
\frac{|f(y, \theta)-f(x, \theta)-\left<g(x, \theta), y-x\right>|}{\|y-x\|} \le 2 L_K(\theta), \\
\|g(y, \theta)-g(x, \theta)\| \le 2 L_K(\theta) .
\end{gathered}
$$
Now by virtue of Lebesgue's theorem on the limit transition under the sign of the integral, we have
$$
\begin{gathered}
\lim _{y \rightarrow x} \frac{F(y)-F(x)-\left<g(x), y-x\right>}{\|y-x\|}= \\
=\int_\Theta \lim _{y \rightarrow x} \frac{f(y, \theta)-f(x, \theta)-\left<g(x, \theta), y-x\right>}{\|y-x\|} P(d \theta)=0 , \\
\lim _{y \rightarrow x}\|g(y)-g(x)\| \le \int_{\Theta} \lim _{y \rightarrow x}\|g(y, \theta)-g(x, \theta)\| P(d \theta)=0,
\end{gathered}
$$
that means that $F(x)$ is continuously differentiable. The lemma is proved.

\begin{lemma}
\label{lem:23.4}
{Suppose that function $f_\alpha(x, \theta) \;\left(\alpha \in K \subset E_m, \quad x \in E_n\right.$, $\theta \in \Theta)$ is measurable in $\theta$ and continuous in $(\alpha, x)$ together with its gradient in $x$ $g_\alpha(x,\theta);$ $K$ is a compact in $E_m$. Let}
$$
\begin{aligned}
f(x, \theta) & =\max \left\{f_\alpha(x, \theta) \mid \alpha \in K\right\}, \\
A(x, \theta) & =\left\{\alpha \in K \mid f_\alpha(x, \theta)=f(x, \theta)\right\}.
\end{aligned}
$$
{Then the function $f(x, \theta)$ is generalized differentiable in $x$, the set of pseudogradients at a point $x$ is given by the formula},
$$
G_f(x, \theta)=\operatorname{co}\left\{g \in E_n \mid g=g_\alpha(x, \theta), \alpha \in A(x, \theta)\right\},
$$
{multivalued mappings $(x, \theta) \rightarrow G_f(x, \theta)$ and $(x, \theta) \rightarrow A(x, \theta)$ are $\mathcal{B} \times \Sigma$-measurable, and if $\alpha(x, \theta)$ is  a $\mathcal{B} \times \Sigma$-measurable section of $A(x, \theta)$, then $g_{\alpha(x, \theta)}(x, \theta)$ is $\mathcal{B} \times \Sigma$-measurable cross section of $G_f(x, \theta)$. If, in addition, $f(x, \theta)$ is integrable over $\theta$ and for any compact $B \subset E_n$ holds}
$$
\sup \left\{\left\|g_\alpha(x, \theta)\right\| \mid \alpha \in K, x \in B\right\} \le L_{K, B}(\theta),
$$
{where $L_{K, B}(\theta)$ is integrable, then the function $F(x)=\int_\Theta f(x, \theta) P(d \theta)$ is generalized differentiable, $G_F(x)=\int_\Theta G_f(x, \theta) P(d \theta)$ is some pseudo-gradient mapping for $F(x)$.}
\end{lemma}
{\it P r o o f.} The function $f(x, \theta)$ is weakly convex in $x$ [91]; therefore it is generalized differentiable (Theorem 1.3) and $G_f(x$, $\theta)$ is its pseudogradient mapping.

It is not difficult to show that for any compact $D \subset E_m$ the function 
$$\max \left\{F_\alpha(x, \theta) \mid \alpha \in K \cap D\right\}$$ 
is continuous in $x$ and measurable in $\theta$, i.e., it is Carath\'eodory one and hence (Lemma 21.6), is $\mathcal{B} \times \Sigma$-measurable. Then the mapping $A(x, \theta)$ is $\mathcal{B} \times \Sigma$-measurable by virtue of $\mathcal{B} \times \Sigma$-measurability of the sets
$$
\{(x, \theta) \mid A(x, \theta) \cap D \neq \emptyset\}=\left\{(x, \theta) \mid \max _{\alpha \in K \cap D} f_\alpha(x, \theta)=f(x, \theta)\right\} .
$$

The mapping $g_\alpha(x, \theta)$ is continuous in $\alpha$, continuous in ${x}$, and measurable in $\theta$ (as the limit of measurable finite differences); by Lemma 21.6, it is $\mathcal{B}_{E_m} \times\left(\mathcal{B}_{E_n} \times \Sigma\right)$-measurable. If $\alpha(x, \theta)$ is $\mathcal{B} \times \Sigma$-measurable section of $A(x, \theta)$, then the section $g_{\alpha(x, \theta)}(x, \theta)$ of the mapping $G_f(x, \theta)$ is $\mathcal{B} \times \Sigma$-measurable as a superposition of measurable mappings.

Let $\left\{\alpha_l(x, \theta)\right\}_{l=1}^{\infty}$ be a countable dense family of measurable sections of $A(x, \theta)$; then $\left\{g_{\alpha_l(x, \theta)}(x, \theta)\right\}_{l=1}^{\infty}$ is a countable dense family of measurable sections of the mapping
$$
G_f^{\prime}(x, \theta)=\left\{g_\alpha(x, \theta) \mid \alpha \in A(x, \theta)\right\}.
$$
By Theorem 21.2, $ G_f^{\prime}(x, \theta)$ is $\mathcal{B} \times \Sigma$-measurable, and by Lemma 21.2 $G_f(x, \theta)=\operatorname{co} G_f^{\prime}(x, \theta)$ is also $\mathcal{B} \times \Sigma$-is measurable.

The generalized differentiability of $F(x)$ and the pseudogradientness of $G_F(x)$ follow from Theorem 23.1. The lemma is proved.

\begin{remark}
\label{rem:23.1}
The mapping $A(x, \theta)$ and the cross section $\alpha(x, \theta) \in A(x, \theta)$ are exact solutions to the problem $\max \left\{f\alpha(x, \theta) \mid \alpha \in K\right\}$. In practice, this problem is solved by some algorithm (either exactly or approximately), i.e., for given $(x, \theta)$ some solution $\tilde{\alpha}(x, \theta)$ is constructed and then substituted into $g_\alpha(x, \theta)$. To guarantee $\mathcal{B} \times \Sigma$-measurability of the function $g_{\tilde{\alpha}(x, \theta)}(x, \theta)$, we must specifically check whether the result $\tilde{\alpha}(x, \theta)$ of the algorithm is $\mathcal{B} \times \Sigma$-measurable.
\end{remark}

Let $f_0(y, \theta)\left(y \in E_m\right)$ and $f_i(x, \theta)\left(x \in E_n, i=1, \ldots, m\right)$ be random generalized differentiable functions; $G_0(y, \theta)$ and $G_i(x, \theta) \quad(i=1,2, \ldots, m)$ are their pseudogradient mappings; $L_0\left(\theta, K_0\right)$ and $L_i(\theta, K)$ are the corresponding Lipschitz constants on sets $K_0 \subset E_m, K \subset E_n$, and
$$
\begin{gathered}
\sup \left\{\|g\| \mid g \in G_0(y, \theta), y \in K_0\right\} \le L_0\left(\theta, K_0\right), \\
\sup \left\{\|g\| \mid g \in G_i(x, \theta), x \in K\right\} \le L_i(x, K) .
\end{gathered}
$$
Let us introduce the vector function $z(x, \theta)=\left(f_1(x, \theta), \ldots, f_m(x, \theta)\right)$ and $n \times m$-matrix $\left[g^1 \ldots g^m\right]$. Consider the function
\begin{align*}
f(x, \theta)=f_0(z(x, \theta), \theta)=f_0\left(f_1(x, \theta), \ldots, f_m(x, \theta), \theta\right).  \tag*{(23.5)}
\end{align*}

Its properties are studied further in lemmas $23.5-23.8$.
\begin{lemma}
\label{lem:23.5}
{The function $f(x, \theta)$ defined by formula (23.5) is a random generalized differentiable function, the mapping $x \rightarrow G_f(x, \theta),$\ where}
$$
\begin{gathered}
G_f(x, \theta)=\operatorname{co}\left\{g \in E_n \mid g=\left[g^1 \ldots g^m\right] g^0,\right. \\
\end{gathered}
$$
\begin{align*}
\left.g^0 \in G_0(z(x, \theta), \theta), g^i \in G_i(x, \theta), i=1, \ldots, m\right\},  \tag*{(23.6)}
\end{align*}
is the pseudogradient one in $x$ for $f(x, \theta)$.
\end{lemma}

Indeed, the function $f_0(y, \theta)$ is Carath\'{e}odory one and hence $\mathcal{B}_{E_m} \times \Sigma$-measurable, and the vector-function $z(x, \cdot)$ is measurable by convention (for a fixed $x$ ). Then $f(x, \cdot)$ is measurable in $\theta$ at fixed $x$ as a superposition of the above measurable functions. The generalized differentiability of $f(\cdot, \theta)$ in $x$ and the pseudogradientness of $G_f(x, \theta)$ at each $\theta$ are proved in Theorem 1.6 .
\begin{lemma}
\label{lem:23.6}
{Let the functions $f_0(y, \theta)$ and $f_i(x, \theta)\;(i=1,2, \ldots, m)$ are bounded on $K_0 \times \Theta$ and $K \times \Theta$, respectively, and the corresponding Lipschitz constants $L_0\left(\theta, K_0\right)$ and $L_i(\theta, K)$ are bounded on $\Theta$. Then the compound function $f(x, \theta)$ satisfies the Lipschitz condition on $\mathrm{\text {K}}$ and is integrable with its Lipschitz constant, and the function $F(x)=\int_{\Theta} f(x, \theta) P(d \theta)$ is generalized differentiable and $G_F(x)=\int_\Theta G_f(x, \theta) P(d \theta)$ is its some pseudogradient mapping. } 
\end{lemma}

Indeed, the Lipschitz property of the compound function $f(x, \theta)$ is easily verified, the integrability of $f(x, \theta)$ and $G_f(x, \theta)$ follows from their measurability (Lemmas $23.5,23.7$) and boundedness, and the generalized differentiability of $F(x)$ and the pseudogradientness of $G_F(x)$ follow from Theorem 23.1. 
\begin{lemma}
\label{lem:23.7}
{Suppose $g^0(y, \theta)$ is a $\mathcal{B}_{E_m} \times \Sigma$-measurable cross section of $G_0(y, \theta)$, and $g^i(x, \theta)$ is a $\mathcal{B}_{E_n} \times \Sigma$-measurable cross section of $G_i(x, \theta)$ $\{i=1,2, \ldots, m)$. Then the function}
$$
g(x, \theta)=\left[g^1(x, \theta) \ldots g^m(x, \theta)\right] g^0(z(x, \theta), \theta)
$$
{is a $\mathcal{B} \times \Sigma$-measurable cross section of $G_f(x,\theta)$ (see (23.6)). If $G_0(y,\theta)$ and $G_i(x,\theta)$ are jointly measurable over their variables, then also $G_f(x, \theta)$ is $ \mathcal{B} \times \Sigma$-measurable}.
\end{lemma}

{\it P r o o f.} Note that the measurability of a vector-function is equivalent to its coordinate measurability; therefore, to prove the measurability of $g(x, \theta)$ it suffices to show the measurability of $g^0(z(x, \theta), \theta)$. The vector-function $z(x, \theta)$ is Carath\'{e}odory one, and by virtue of Lemma 21.6 it is $\mathcal{B} \times \Sigma$-measurable; then $g^0(z(x, \theta), \theta)$ is $\mathcal{B} \times \Sigma$-measurable as a superposition of measurable functions.

Let us prove $\mathcal{B} \times \Sigma$-measurability of the mapping $(x, \theta) \rightarrow G_f(x, \theta)$. Let us denote
$$
\begin{gathered}
G_f^{\prime}(x, \theta)=\left\{g \in E_n \mid g=\left[g^1 g^2 \ldots g^m\right] g^0,\right. \\
\left.g^0 \in G_0(z(x, \theta), \theta), g^i \in G_i(x, \theta), i=1, \ldots, m\right\},
\end{gathered}
$$
i.e., $\operatorname{co} \ G_f^{\prime}(x, \theta)=G_f(x, \theta)$. By virtue of Lemma 21.3, it suffices to show the measurability of $(x, \theta) \rightarrow G_f^{\prime}(x, \theta)$. Let us show the existence of a countable dense family of measurable sections of the mapping $G^{\prime}$ over the set of variables; by Theorem 21.2, the measurability of $G^{\prime}$ follows from here. By Theorem 21.2, there exist countable dense  $\mathcal{B} \times \Sigma$-measurable families of sections $\left\{g_{\mu_0}^0, \mu_0=1,2 \ldots\right\}$ and $\left\{g_{\mu_i}^i, \mu_i=1,2, \ldots\right\}$ of mappings $G_0$ and $G_i$, respectively. Obviously, the function
$$
g_\mu(x, \theta)=\left[g_{\mu_1}^1(x, \theta) \ldots g_{\mu_m}^m(x, \theta)\right] g_{\mu_0}^0(z(x, \theta), \theta)
$$
where $\mu=\left(\mu_0, \mu_1, \ldots, \mu_m\right)$, is a $\mathcal{B} \times \Sigma$-measurable cross section of $G^{\prime}$. The family $\left\{g_\mu\right\}$ is countable. We show that the set of values of $\left\{g_\mu(x, \theta)\right\}$ is dense in $G_f^{\prime}(x, \theta)$ for all $(x, \theta)$. Let us take an arbitrary $g^{\prime} \in G_f^{\prime}(x, \theta)$. By construction, there exist $g_0^{\prime} \in G_0(z(x, \theta), \theta)$ and $g_i^{\prime} \in G_l(x, \theta)$ such that $g^{\prime}=\left[g_1^{\prime} \ldots g_m^{\prime}\right] g_0^{\prime}$. By virtue of the density of families of sections $g_{\mu_i}^i$, there exist $g_{\mu_0^k}^0(z(x, \theta), \theta) \rightarrow g_0^{\prime}$ and $g_{\mu_i^k}^i(x, \theta) \rightarrow g_i^{\prime}$. Then
$$
g_{\mu^k}(x, \theta)=\left[g_{\mu_1^k}^1(x, \theta) \ldots g_{\mu_m^k}^m(x, \theta)\right] g_{\mu_0^k}^0(z(x, \theta), \theta) \rightarrow g^{\prime} .
$$
Lemma 23.7 is proved.

\begin{lemma}
\label{lem:23.8}
{Consider the generalized differentiable (see Theorem 1.5) max-function}
$$
f(x, \theta)=\max _{1 \le i \le m} f_i(x, \theta)
$$
{As shown in $\S$ 1, its pseudogradient mapping in $x$ is given by the formulas}
\begin{align*}
G_f(x, \theta) & =\operatorname{co}\left\{g \in E_n \mid g_i \in G_i(x, \theta), \quad i \in I(x, \theta)\right\}, \tag*{(23.7)}
\end{align*}
\begin{align*}
I(x, \theta) & =\left\{i \mid 1 \le i \le m, f_i(x, \theta)=f(x, \theta)\right\}. \tag*{(23.8)}
\end{align*}

{If $g_i(x, \theta)\;(i=1, \ldots, m)$ are $\mathcal{B} \times \Sigma$-measurable sections of $G_{\imath}(x, \theta)$, and $i(x, \theta)$ is $\mathcal{B} \times \Sigma$-measurable section of the mapping}
\noindent$(x, \theta) \rightarrow I(x, \theta),$ then
$$
g(x, \theta)=g_{i(x, \theta)}(x, \theta)
$$
{is $\mathcal{B} \times \Sigma$-the measurable cross section of the mapping $G_f(x, \theta)$.

If the mappings $G_i(i=1, \ldots m)$ are $\mathcal{B} \times \Sigma$-measurable, then $G_f$ is $\mathcal{B} \times \Sigma$-measurable.}
\end{lemma}
\begin{remark}
\label{rem:23.2}
To construct a measurable section of the mapping $G_f$ in Lemma 23.8, we need $\mathcal{B} \times \Sigma$-measurable sections of the mapping $(x, \theta) \rightarrow I(x, \theta)$.

Let us show that the mapping $(x, \theta) \rightarrow I(x, \theta)$ is $ \mathcal{B} \times \Sigma$-measurable. The finite set $I(x, \theta)$ is closed. For any closed $B $ $\subset E_1$, the set
$$
\{(x, \theta) \mid I(x, \theta) \cap B \neq \emptyset\}=\left\{(x, \theta) \mid \max _{\substack{i \in B \\ 1 \le i \le m}} f_i(x, \theta)=f(x, \theta)\right\}
$$
is $\mathcal{B} \times \Sigma$-measurable, since $f(x, \theta)$ and $f_i(x, \theta)$ are Carath\'{e}odory ones. Hence, the mapping $(x, \theta) \rightarrow I(x, \theta)$ is $ \mathcal{B} \times \Sigma$-measurable, and by Theorem 21.2 there exist $\mathcal{B} \times \Sigma$-measurable sections of the mapping $(x, \theta) \rightarrow I(x, \theta)$. Consider, for example,
$$
i(x, \theta)=\min \{i \mid i \in I(x, \theta)\}
$$

The function $i(x, \theta)$ is $ \mathcal{B} \times \Sigma$-measurable since the Lebesgue sets $\left(c \in E_1\right)$
$$
\{(x, \theta) \mid i(x, \theta) \le c\}=\{(x, \theta) \mid I(x, \theta) \cap\{r \mid r \le c\} \neq \emptyset\} 
$$
are $\mathcal{B} \times \Sigma$-measurable.
\end{remark}
\begin{remark}
\label{rem:23.3}
If in the conditions of Lemma 23.8 functions $f_i(x$, $\theta)\;(i=1,2, \ldots, m)$ satisfy the Lipschitz condition in $x$ in compact $K \subset E_n$ with constants $L_i(\theta)$, then the maximum function $f(x, \theta)$ satisfies the Lipschitz condition in $K$ with constant $L(\theta)=\max _{1 \le i \le m} L_i(\theta)$. If, in addition, $f_i(x, \theta)$ are integrable with its Lipschitz constants, then the maximum function is integrable together with its Lipschitz constant, and by virtue of Theorem 23.1, the function $F(x)=\int_\Theta f(x$,$\theta) \ P (d\theta)$ is generalized differentiable.
\end{remark}

{\it P r o o f of Lemma 23.8.} The cross section $g(x, \theta)$ is $\mathcal{B} \times \Sigma$-measurable, since  the Lebesgue sets $\left(c \in E_n\right)$ are $\mathcal{B} \times \Sigma$-measurable,
$$
\{(x, \theta) \mid g(x, \theta) \le c\}=\bigcup_{k=1}^m\left(\{(x, \theta) \mid i(x, \theta)=k\} \cap\left\{(x, \theta) \mid g_k(x, \theta) \le c\right\}\right) .
$$

Let us prove $\mathcal{B} \times \Sigma$-measurability of $G_f$ if $G_i(i=1, \ldots, m)$ are $\mathcal{B} \times \Sigma$-measurable. By virtue of lemma 21.3, it suffices to show the measurability of the mapping
$$
(x, \theta) \rightarrow G_f^{\prime}(x, \theta)=\left\{g \in E_n \mid g \in G_i(x, \theta), i \in I(x, \theta)\right\} .
$$
Let $\left\{g_i^{\nu_i}(x, \theta)\right\}_{\nu_i=1}^{\infty} \;(i=1, \ldots, m)$ are countable dense families of $\mathcal{B} \times \Sigma$-measurable sections of $G_i$. It suffices to consider the case $m=2$.

Let us define functions:
$$
\begin{aligned}
h_1^{\nu_1 \nu_2}(x, \theta) & = \begin{cases}g_1^{\nu_1}(x, \theta), & f_1(x, \theta) \ge f_2(x, \theta), \\
g_2^{\nu_2}(x, \theta), & f_1(x, \theta)<f_2(x, \theta);\end{cases} \\
h_2^{\nu_1 \nu_2}(x, \theta) & = \begin{cases}g_1^{\nu_1}(x, \theta), & f_1(x, \theta)>f_2(x, \theta), \\
g_2^{\nu_2}(x, \theta), & f_1(x, \theta) \le f_2(x, \theta) .\end{cases}
\end{aligned}
$$

Recall that $f_1(x, \theta)$ and $f_2(x, \theta)$ are Carath\'{e}odory ones; therefore $h_1^{\nu_1 \nu_2}(x, \theta)$ and $h_2^{\nu_1 \nu_2}(x, \theta)$ are $ \mathcal{B} \times \Sigma$-measurable. Obviously, the set $\left\{h_1^{\nu_1 \nu_2}(x, \theta)\right.$, $\left.h_2^{\nu_1 \nu_2}(x, \theta)\right\}_{\nu_1, \nu_2=1}^{\infty}$ is a countable dense family of measurable sections of the mapping $G^{\prime}$. The lemma is proved.

\newpage

\begin{flushright}
\text{CHAPTER 7}
\label{Ch.7}

\textbf{SOLVING STOCHASTIC EXTREMAL PROBLEMS}
\end{flushright}
\hrule
\vspace{50pt}

This chapter focuses on numerical methods for solving stochastic problems in which there are no continuous derivatives of the objective and constraint functions, and the function values are not computed exactly and are probabilistic in nature.

The following extremal problem is considered:
$$
F_0(x) \equiv \mathbb{E} f_0(x, \theta) \rightarrow \min _x
$$
subject to
\begin{align*}
F_l(x)=\mathbb{E} f_i(x, \theta) \le 0, \quad i=1,2, \ldots, m, \tag*{(C)} \\
x \in X \subset E_n,\;\;\;\;\;\;\;\;\;\;\;\;\;\;\;\;\;\;\;\;
\end{align*}
where $\theta$ is an elementary event of some probability space $(\Theta, \Sigma, P),$ $\mathbb{E}$ is the sign of mathematical expectation.

The main difficulty in solving the problem $(\mathrm{C})$ is that for a fixed $x$ the computation of $F_i(x)\;(i=0,1, \ldots, m)$ requires considerable effort, since it involves the computation of multivariate integrals. Instead of values  $F_i(x)\;(i=0,1, \ldots, m)$, it is usually possible to observe at a given $x$ the values of random variables $f_i(x, \theta)$. Therefore, nonlinear programming methods cannot be widely applied to problems (C), and stochastic methods are developed to solve them.

The numerical methods studied in this chapter allow solving problem (C), when the distribution laws of $\theta$ are unknown, but there is a way to compute the random variables $f_i(x, \theta)\;(i=0,1, \ldots$ $\ldots, m)$. The same methods are also applicable, when the probabilistic properties of $\theta$ are given, but the computation of functions $F_i(x)\;(i=0,1, \ldots, m)$ is either impossible (the dependence of functions $f_i(x, \theta)$ on $\theta$ is unknown) or too complicated.

There are direct and indirect methods of stochastic programming. If the probabilistic properties of $\theta$ are known and the analytical dependence $f_i(x, \theta)$ on $\theta$ is given, then the mathematical expectations in (C), at least in principle, can be calculated. Then problem (C) formally ceases to be stochastic. In indirect methods, one tries to obtain the dependencies $F_i(x)(i=0,1, \ldots, m)$ and then apply nonlinear programming methods. In other words, instead of a stochastic problem, its deterministic version is considered. Therefore, the success of indirect methods application depends to a large extent on
 the probabilistic properties $\theta$ and properties
of the functions $f_{i}(x, \theta)$.

Methods that are based on the information about the values $f_{i}(x, \theta)$
or analogues of their gradients are called direct methods of stochastic programming.
When solving practical problems, functional dependencies $f_{i}(x, \theta)$,
as a rule, are nonlinear, and the vectors $x$ and $\theta$ have a large dimension.
Therefore, finding the analytical form of the functions $F_{i}(x, \theta)$
most often is not possible, since this is associated with the definition of mathematical 
expectations, which are multidimensional Lebesgue integrals.

If the probability distribution $\theta$ is given, then the calculation of $f_{i}(x, \theta)$ 
is reduced to modeling random variables with a given distribution law.
This process is usually carried out on a computer with the help of transformations of one 
or several independent values of a random variable uniformly distributed in the interval $[0,1]$.

An essential feature of direct methods is that they are applicable to solving 
problems in which the distribution laws of $\theta$ are unknown.
In the process of optimizing complex objects, situations may arise when 
the formulation of the probabilistic properties of $\theta$ presents a significant difficulty 
and it is necessary to solve the problem without an analytical study of the probabilistic model.
Direct methods require only observations $\theta_{1}, \theta_{2}, \ldots, \theta_{k}, \ldots$ 
of parameter $\theta$. Therefore, these methods can be widely used to optimize systems 
based on simulation. Direct methods are the main subject of study in this chapter.

Let's do more general remarks about the technology of solving stochastic programming 
problems.
Firstly, the stochasticity of the problem may not be clearly manifested in the entire 
admissible region, but only, say, in some neighborhood of the solution, i.e. far 
from the solution, the variance of the objective function gradients can be small 
compared to the gradients themselves, and the problem is almost deterministic. 
Therefore, one should first apply non-linear programming methods that are 
resistant to errors in the gradients used. If the analytical form of the random 
functions of the problem is known, then, having observed one or more implementations 
of the random parameters of the problem, one can first replace the integral, the mathematical
 expectation, with an integral sum and use non-linear programming methods 
(non-smooth optimization) to minimize this sum. Along the way, one needs to periodically 
monitor the variance of the gradients of the objective function.
If it becomes comparable to the gradients themselves, then stochastic analogues 
of stable non-linear programming methods should be used. If the variance becomes 
larger than the gradients themselves, then one should either reduce this variance 
by additional observations,
or switch to a simple stochastic gradient method, since complex methods 
in such a situation become no better than a simple one.

\section*{$\S$ 24. Methods of averaged stochastic gradients}
\label{Sec.24}
\setcounter{section}{24}
\setcounter{definition}{0}
\setcounter{equation}{0}
\setcounter{theorem}{0}
\setcounter{lemma}{0}
\setcounter{remark}{0}
\setcounter{corollary}{0}
In [41], to solve the problem of convex stochastic programming, a method 
of the stochastic approximation type was proposed, the method of stochastic quasi-gradients, 
which uses the gradients of a random function under the sign of mathematical expectation and 
the operation of projecting onto a convex set. 
In [73], this method was improved and extended to the infinite-dimensional case. 
In [91], the method of sochastic quasigradients was extended to stochastic 
programming problems with weakly convex functions (see $\S$ 1). 
In this section, we consider the method of stochastic generalized gradients for solving 
a stochastic programming problem with generalized differentiable functions under constraints. 
Methods with averaging of stochastic gradients are also considered, in particular, 
stochastic analogues of the heavy ball and ravine step methods (see $\S$ 14).

Consider the problem of non-convex stochastic programming:
\begin{equation}
F(x) \equiv \mathbb{E} f(x, \theta) \rightarrow \min _{x} \tag{24.1}
\end{equation}
under constraints
\begin{equation}
h(x) \le 0, \quad x \in E_{n}, \tag{24.2}
\end{equation}
where $\theta \in \Theta ;(\Theta, \Sigma, P)$ is a probability space; $f(x, \theta): E_{n} \times \Theta \rightarrow$ $\rightarrow E_{1}$ is a random (measurable in $\theta$ for every $x$ ) generalized differentiable (in $x$ for all $\theta$ ) integrable together with its own Lipschitz constant $L_{K}(\theta)$ for each compact set $K \subset E_{n}$; the function $h: E_{n} \rightarrow E_{1}$ is generalized differentiable, and $h(x) \rightarrow+\infty$ as $\|x\| \rightarrow+\infty$.

By Theorem 23.1, the function $F$ is generalized differentiable. Suppose that some measurable (in $\theta$ for each $x$ ) pseudogradient (in $x$ for all $\theta$ ) mapping $G_{f}(x, \theta)$ of the function $f(x, \theta)$, as well as the $\mathcal{B} \times \Sigma$-measurable section $g_{f}(x, \theta)$ of the mapping $(x, \theta) \rightarrow G_{f}(x, \theta)$ is known. As shown in $\S$ 23, the required $G_{f}$ and $g_{f}$ exist. By Theorem 23.1, $G_{F}(x)=\int_{\Theta} G_{f}(x, \theta) P(d \theta)$ is some pseudogradient mapping for $F(x)$ and $g_{ f}(x)=\mbox{E} g_{f}(x, \theta) \in G_{F}(x)$. Assume that the pseudogradient mapping $G_{h}(x)$ and its Borel section $g_{h}(x) \in G_{h}(x)$ are also known. Since $G_{h}:E_{n} \rightarrow 2^{E_{n}}$ is closed, the required section exists (Lemma 21.2).

Denote
\begin{equation}
\tag{24.3}
\begin{array}{r}
G(x)= \begin{cases}G_{F}(x), & h(x)<0, \\
\operatorname{co}\left\{G_{F}(x), G_{h}(x)\right\}, & h(x)=0, \\
G_{h}(x), & h(x)>0\end{cases} \\ 
g(x, \theta)= \begin{cases}g_{f}(x, \theta), & h(x)<0, \\
g_{h}(x, \theta), & h(x) \ge 0 .\end{cases}\;\;\;\;\;\;\;\;\;\;\;\;\;\;\;\;\;
\end{array}
\end{equation}

The function $g(x, \theta)$ is $ \mathcal{B} \times \Sigma$-measurable, since the sets $\left(c \in E_{n}\right)$

$\{(x, \theta) \mid g(x, \theta) \le c\}=\left(\left\{(x, \theta) \mid g_{f}(x, \theta) \le c\right\} \cap\right.$

$\;\;\;\;\;\;\;\;\;\;\;\;\;\;\;\;\;\;\;\;\cap\{(x, \theta) \mid h(x)<0\}) \cup\left(\left\{(x, \theta) \mid g_{h}(x, \theta) \le c\right\} \cap\{(x, \theta) \mid h(x) \ge 0\}\right)$

\noindent
are $\mathcal{B} \times \Sigma$-measurable.

Denote the set of pseudostationary points of problem (24.1), (24.2) 
$$X^{*}=\left\{x \in E_{n} \mid 0 \in G(x)\right\}.$$

Let us introduce the notation:

$\Theta^{i}=\Theta \times \ldots \times \Theta$ be the $i$-th Cartesian power of $\Theta$;

$\theta^{i}=\left(\theta_{0}, \ldots, \theta_{i-1}\right)$ be an element of the space $\Theta^{i}$;

$\Sigma^{i}=\Sigma \times \ldots \times \Sigma=\sigma\left\{A_{0} \times \ldots \times A_{i-1} \mid A_{k} \in \Sigma, 0 \le k<i\right\}$ be the $\sigma$-algebra, which is a product of $i$ $\sigma$-algebras  $\Sigma$;

$P^{i}=P \times \ldots \times P$ be the probability measure on $\left(\Theta^{i}, \Sigma^{i}\right)$, which is the product of $i$ measures $P $; $\left(\Theta^{i} ; \Sigma^{i}, P^{i}\right)$ be the corresponding probability space;

$\Omega=\Theta^{\infty}=\Theta \times \Theta \times \ldots$ be the countable direct product of $\Theta$ spaces;

$\Sigma_{\omega}=\Sigma^{\infty}=\Sigma \times \Sigma \times \ldots=\sigma\left\{A_{0} \times A_{1} \times \ldots \mid A_ {k} \in \Sigma, k \ge 0\right\}$ be the countable product of $\sigma$-algebras $\Sigma$.

Measures $P^{i}$ on $\left(\Theta^{i}, \Sigma^{i}\right)$, according to Ionescu Tulcea's theorem (see [133]), can be uniquely extended to some probabilistic measure $P_{\omega}$ on $\left(\Omega, \Sigma_{\omega}\right)$ so that
$$
P^{i}\left\{A \mid A \in \Sigma^{i}\right\}=P_{\omega}\left\{A \times \Theta \times \Theta \times \ldots \mid A \in \Sigma^{i}\right\} .
$$

Then $\left(\Omega, \Sigma_{\omega}, P_{\omega}\right)$ is a probability space.

\textbf{1. The method of averaged stochastic gradients.} 
\label{Sec.24.1}
The method is designed to solve problem (24.1), (24.2); it generates a random sequence $\left\{x^{k}(\omega)\right\}$ of approximations to the solution of the problem according to the following relations:
\begin{align}
x^{0}(\omega)&=x^{0} \in E_{n}, \tag{24.4}\\
x^{k+1}(\omega)&=x^{k}(\omega)-\rho_{k} P^{k}(\omega), \tag{24.5}\\
P^{k}(\omega)&=\sum^{k} \lambda_{k r} g^{r}(\omega) / H_{r}(\omega), \tag{24.6}\\
g^{r}(\omega)&=g\left(x^{r}(\omega), \pi_{r}(\omega)\right), \quad \pi_{r}(\omega)=\theta_{r}, \tag{24.7}\\
H_{r}(\omega)&=H_{r}\left(x^{0}(\omega), \ldots, x^{r}(\omega), g^{0}(\omega), \ldots, g^{r-1}(\omega)\right) \ge \nu>0, \tag{24.8}
\end{align}
where $\omega=\left(\theta_{0}, \theta_{1}, \ldots\right) \in \Omega;$ $\theta_{r}(r=0,1, \ldots)$ are independent realizations of the random parameter $\theta \in \Theta$; the mapping $g(x, \theta)$ is defined in (24.3) and is bounded in $x$ in each compact set $K \subset E_{n}$ uniformly in $\theta$; $H_{r}\left(x^{0}, \ldots, x^{r}, g^{0}, \ldots, g^{r-1}\right)$ is a Borel function of its arguments, bounded in every compact in $\left(E_{n}\right)^{2 r+1} ;$ $v$ is some positive constant; $\rho_{k}, r_{k}, \lambda_{k r}$ are deterministic functions of $k$ and $r$ satisfying the following conditions:
\begin{equation}
\tag{24.9}
\lambda_{k r} \ge 0, \quad \sum_{r=r_{k}}^{k} \lambda_{k r}=1
\end{equation}
(for convenience, we assume that $\lambda_{k r}=0$ for $r \in\left[0, r_{k}\right)$ );
\begin{align}
& 0 \le \rho_{k} \le \rho, \quad \lim _{k \rightarrow \infty} \sum_{r=r_{k}}^{k} \rho_{k}=0, \quad \sum_{k=0}^{\infty} \rho_{k}=+\infty; \tag{24.10} \\
& \sum_{r=0}^{\infty}\left(\sum_{k=r}^{\infty} \lambda_{k r} \rho_{r}\right)^{2}=\sum_{k , s=0}^{\infty} \rho_{k} \rho_{s} \sum_{r=\max \left(r_{k}, r_{s}\right)}^{\min (k, s)} \lambda_{k r} \lambda_{k s}<+\infty. \tag{24.11}
\end{align}
\begin{remark}
\label{rem:24.1}
It is assumed that relations (24.4)--(24.8) hold for all $\omega \in \Omega$. The random process $\left\{x^{k}(\omega)\right\}$ is considered in the probability space $\left(\Omega, \Sigma_{\omega}, P_{\omega}\right)$. The vectors $g^{r}(\omega)$ are called stochastic (in $\Omega$ ) generalized gradients of the problem functions at  points $x^{r}(\omega)$, the quantities $H_{r}(\omega)$ are normalizing coefficients, and the directions $P^{k}(\omega)$ can be called averaged stochastic gradients.
\end{remark}
\begin{remark}
\label{rem:24.2}
The measurability of the trajectory $\left\{x^{k}(\omega)\right\}_{k=0}^{\infty}$ is established recursively. The function $x^{0}(\omega)=x^{0}$ is obviously measurable. 
If $x^{0}, \ldots, x^{k}(\omega), g^{0}(\omega), \ldots, g^{k-1}(\omega)$ are measurable, then $ g^{k}(\omega)$ is 
measurable as a superposition of the measurable mapping $z^{k}(\omega)=\left(x^{k}(\omega), \pi_{k}(\omega)\right ): \Omega \rightarrow E_{n} \times \Theta$ and $\mathcal{B} \times \Sigma$-measurable mapping $g: E_{n} \times \theta \rightarrow E_{n}$ 
(at this point, the $\mathcal{B} \times \Sigma$-measurability of the  sections $g_{j}(x, \theta)$ and the Borel property of $g_{h}(x)$ are used). 
The quantity $H_{r}(\omega)$ is measurable as a superposition of a measurable mapping $y^{r}(\omega)=\left(x^{0}(\omega), \ldots, x^{r}(\omega), 
g^{0}(\omega), \ldots, g^{r-1}(\omega)\right): \Omega \rightarrow\left(E_{n}\right)^{2 r+ 1}$ and the Borel function $H_{r}:\left(E_{n}\right)^{2 r+1} \rightarrow E_{1}$. This implies that $x^{k+1}(\omega)$ is also measurable.
\end{remark}
\begin{remark}
\label{rem:24.3}
Outside the admissible region $\{x \mid h(x) \le 0\}$ of problem $(24.1),(24.2)$, the method of averaged stochastic gradients (24.4) - (24.11) actually works as a deterministic method of averaged gradients (14.22) - (14.33). Therefore, its boundness (for any $\omega$ ) can be ensured in the same way as in the purely deterministic case, i.e., either by choosing sufficiently small steps $\rho_{k} \le \rho$, or by using the return to starting point (see $\$$ 14).
\end{remark}
\begin{remark}
\label{rem:24.4}
If there exists an estimate $\left\|g_{f}(x, \theta)\right\| \le L_{f}(x)$ for all $\theta$, where $L_{f}(x)$ is a continuous function, then certainly $g_{f}(x, \theta)$ is bounded in $x$ in each compact subset of $E_{n}$ uniformly in $\theta$.
\end{remark}
\begin{remark}
\label{rem:24.5}
Let the estimate $\left\|\mbox{E} g_{f}(x, \theta)\right\| \le H(x)$ is known, where $H(x)$ is a continuous function; then the normalizing coefficients can be taken in the following form $(v>0)$ :
$$
H_{r}(\omega)=\left\{\begin{aligned}
H\left(x^{r}(\omega)\right)+\nu, & \;\;\;h\left(x^{r}(\omega)\right)<0, \\
\left\|g_{h}\left(x^{r}(\omega)\right)\right\|+\nu, & \;\;\;h\left(x^{r}(\omega)\right) \ge 0 .
\end{aligned}\right.
$$
Alternatively, one can assume
$$
H_{0}(\omega)=\text { const, } \quad H_{r}(\omega)=\left\|P^{r-1}(\omega)\right\|+\nu, \quad \nu>0, \quad r=1,2, \ldots
$$
\end{remark}
\begin{remark}
\label{rem:24.6}
Condition (24.11) is satisfied, if
$$
\sum_{k=0}^{\infty}\left(\rho_{k} \sum_{r=r_{k}}^{k} \rho_{r}\right)<+\infty,
$$
which takes place, for example, when $k-r_{k} \le$ const $(k)$ or when $\rho_{k}=c / k^{\alpha}$, $k-\left(c_{ 1}+c_{2} k^{\beta}\right) \le r_{k} \le k,$ $0.5<\alpha \le 1,$ $0<\beta<2 \alpha-1 ;$ $ c_{1}, c_{1}, c_{2}>0$. For the indicated $\alpha, \beta$ relation (24.10) also holds.
\end{remark}
\begin{remark}
\label{rem:24.7}
It is possible, following [73, 143], to speed up the convergence of the method by averaging the minimizing sequence $\left\{x^{k}(\omega)\right\}_{k=0}^{\infty}$, for example, in the following way:
\begin{align}
\bar{x}^{0}(\omega) & =x^{0}(\omega), \notag{}\\
\bar{x}^{k}(\omega) & =\left(1-\tau_{k}\right) \bar{x}^{k-1}(\omega)+\tau_{k} x ^{k}(\omega), \quad k \ge 1, \tag{24.12}\\
\tau_{0} & =1, \quad 0 \le \tau_{k} \le 1, \quad \sum_{k=0}^{\infty} \tau_{k}=+\infty. \tag{24.13}
\end{align}

If $\left\{\rho_{k}\right\}$ satisfy (24.10), then by the Abel-Dini lemma
\footnote{Fikhtengol'ts G. M. A course of differential and integral calculus. Vol. 2. - Moscow: Nauka, 1966. 290 p.}
conditions $(24.13)$ are satisfied by the sequence
\begin{equation}
\tau_{k}=\rho_{k} / \sum_{r=0}^{k} \rho_{r}, \quad k \ge 0 .\tag{24.14}
\end{equation}
\end{remark}
\begin{lemma}
\label{lem:24.1}
Let $X^{*}=\bigcup_{j \in J} X_{j}^*$, where $X_{j}^{*}\;(j \in J)$ are connected sets. Then $\left\{x^{k}(\omega)\right\}$ converges under conditions (24.13) to the set $\bar{X}^{*}=\bigcup_{j \in J}\mbox{ co }X_{j}^{*}.$
\end{lemma}
{\it P r o o f.} Denote $\Phi_{r}^{k}=\prod^{k}_{t=r}\left(1-\tau_{t}\right)$ and, for convenience, set $\Phi_{k+1}^{k}=1 $. Note that $\Phi_{r}^{k} \le \exp \left(-\sum_{t=r}^{t} \tau_{t}\right)$ and, under condition (24.13), $\lim _{k \rightarrow \infty} \Phi_{r}^{k}=0$ for any fixed $r$. For any $s \in(0, k]$ the representation is true:
\begin{align}
&\bar{x}^{k}(\omega)=\sum_{r=0}^{k} \tau_{r} \Phi_{r+1}^{k} x^{r}(\omega )=\sum_{r=s}^{k} \mu_{r}^{k} x^{r}(\omega) \notag{} \\
&\;\;\;\;\;\;\;\;\;\;+\Phi_{s}^{k}\left(\sum_{r=0}^{s-1} v_{r}^{k} x^{r}(\omega)-\sum_{r=s }^{k} \mu_{r}^{k} x^{r}(\omega)\right); \tag{24.15}
\end{align}
$$
\begin{array}{cc}
\mu_{r}^{k}=\tau_{r} \Phi_{r+1}^{k} /\left(1-\Phi_{s}^{k}\right), & r \in[ s, k] ; \\
\nu_{r}^{k}=\tau_{r} \Phi_{r+1}^{k} / \Phi_{s}^{k}, \quad & r \in[0, s) ; \\
\sum_{r=0}^{s-1} \nu_{r}^{k}=\sum_{r=s}^{k} \mu_{r}^{k}=1 .
\end{array}
$$

As will be shown (Theorem 24.1), under conditions (24.9)-(24.11) the sequence $\left\{x^{k}(\omega)\right\}_{k=0}^{\infty}$ 
converges to some connected subset $X_{j(\omega)}^{*}$ of the solution set $X^{*}$. 
In (24.15) the first term $\sum_{r=s}^{k} \mu_{r}^{k} x^{r}(\omega)$ for sufficiently large $s$ and all $k \ge s $ will be in an 
arbitrarily small neighborhood of the convex hull of $X_{i(\omega)}^{*}$, and the second term tends to zero as $k \rightarrow+\infty$. 
Therefore $\left\{\bar{x}^{k}\right\}_{k=0}^{\infty}$ converges to an arbitrarily small neighborhood of $\mbox{co }X_{j(\omega)}^{*}$. The lemma is proven.
\begin{remark}
\label{rem:24.8}
If problem (24.1), (24.2) also has linear constraints $N x=b$, where $N$ is an $m \times n$ matrix and $b$ is an $m$-dimensional vector, 
then we can use the method of projection of the (stochastic) gradient onto the linear subspace $L_{0}=\{x \mid N x=0\}$, 
as was done in $\S$ 10 for the generalized gradient method. In $\S$ 10 it was shown that the projection of the pseudo-gradient of the function $F(x)$ onto 
the subspace $L_{0}$ is the pseudogradient of this function, considered only on the linear manifold $L_{b}=\{x \mid N x=b\}$. 
The designed pseudo-gradients can be used in any  methods for optimization of $F(x)$ and operating in the $L_{b}$ manifold. In view of the computability 
of the projection operation on $L_{0}$ and the operation of taking the mathematical expectation, the projected stochastic gradients 
of the function $F(x)=\mbox{E} f(x, \theta)$ will be the stochastic gradients of the function $F(x)$, already considered as a function on the manifold $L_{b}$. 
Therefore, using the designed gradients, one can run the considered stochastic programming methods in the linear manifold $L_{b}$ and solve problems with linear constraints.
\end{remark}
\begin{theorem}
\label{th:24.1}
Let the problem of non-convex stochastic programming (24.1), (24.2) be solved by the method of averaged stochastic gradients $(24.4)-(24.11)$. Let the set $H^{*} = \left\{h(x) \mid 0 \in G_{h}(x), h(x)>0\right\}$ contains no intervals. Assume that the sequence $\left\{x^{k}(\omega)\right\}_{k=0}^{\infty}$ is bounded for any $\omega;$ this is the case, for example, for sufficiently small $\rho=\sup _{k \ge 0} \rho_{k}$ and $R = \sup_{k \ge 0}\sum_{r=r_{k}}^{k} \rho_{ r}$. 
Then: 

either all limit points of $\left\{x^{k}(\omega)\right\}$ do not belong to the admissible region $D=\{x \mid h(x) \le 0\}$, but belong to the set $ X_{h}^{*}=\{x \mid 0 \in G_{h}(x), h(x)>0\}$ and there exists a limit $\lim _{k \rightarrow \infty} h\left(x^{k}(\omega)\right)>0$;

or all limit points of $\left\{x^{k}(\omega)\right\}$ belong to the admissible domain D; in this case, the smallest of them (with respect to the value of $F$) $P_{\omega}$-almost surely ($P_{\omega}$-a.s.) belong to the set $X^{*} =\left\{x \in E_{n} \mid 0 \in G(x), \quad h(x) \le 0\right\}$ and the interval 
$$\left[\varliminf_{k \rightarrow \infty} F(x^{k}(\omega)), \varlimsup_{k \rightarrow \infty} F(x^{ k}(\omega))\right]$$ 
($P_{\omega}$-a.s.) is included in the set $F^{*}= 
\{f(x) \mid x \in X^{*}\}$. If at the same time $F^{*}$ does not contain intervals, then all limit points of $\{x^{k}(\omega)\}$  $P_{\omega}$-a.s. belong to the connected component of $X^{*}$ and $P_{\omega}$-a.s. there is a limit $\lim_{k \rightarrow \infty} F(x^{k}(\omega))$.
\end{theorem}

\begin{lemma}
\label{lem:24.2}
Let there exist $c >\max \left(0, h\left(x^{0}\right)\right)$ not belonging to $H^{*}$. 
Then for any $d>c$ there are sufficiently small $\rho^{\prime}$ and $R^{\prime}$ such that for any sequence $\left\{\rho_{k}\right\}$, 
satisfying (24.10) with $\rho \le \rho^{\prime}$ and $R \le R^{\prime}$, for all $\omega$ the sequence $\left\{x^{h}( \omega)\right\}$ launched from $x^{0}$ 
according to (24.4) -- (24.9) do not leave the bounded region $\{x \mid h(x) \le d\}$.
\end{lemma}

{\it P r o o f.} The lemma is valid because outside the domain $D_{\varepsilon}=\{x \mid h(x) \ge \varepsilon>0\}$ for sufficiently large $k$ the method of averaged 
stochastic gradients minimizes the deterministic function $h( x)$ and therefore works like a deterministic average gradient method. Indeed, suppose the contrary; 
then for some $d>c$ there are $\left\{\rho_{\mathrm{s}}^{k}\right\}(s=0,1, \ldots)$ with $\sup _{k \ge 0} \rho_{\mathrm{s}}^{k}=\rho_{\mathrm{s}}^{\prime} \rightarrow 0$ 
and $\sup _{k \ge 0} \sum_{r=r_{k}}^{k} \rho_{\mathrm{s}}^{k}=R_{\mathrm{s}}^{\prime} \rightarrow 0$ such that the corresponding sequences $ \left\{x_{s}^{k}\left(\omega^{s}\right)\right\}$ ($\left.\omega^{s} \in \Omega\right)$, 
launched from the point $x^{0}$ for $s=0,1, \ldots$ $\left(\omega^{s} \in \Omega\right)$, go to infinity. 
Let us single out indices $k_{s}\left(\omega^{s}\right)$ and $t_{s}\left(\omega^{s}\right)$ such that for $k \in\left( k_{s}, t_{s}\right)$
\begin{equation}
h\left(x_{s}^{k_{s}}\left(\omega^{s}\right)\right) \le c<h\left(x_{s}^{k}\left(\omega^{s}\right)\right)<d \le h\left(x_{s}^{t_s}(\omega^{s})\right). \tag{24.16}
\end{equation}

The sequence $\left\{x_{s}^{k_{s}}\left(\omega^{s}\right)\right\}_{s=0}^{\infty}$ is bounded; without loss of generality, we assume that $\lim _{s \rightarrow \infty} x_{s}^{k_{s}}\left(\omega^{s}\right)=x^{\prime}$ . Under assumptions $(24.10)$ and due to (24.16), we have
$$
\lim _{s \rightarrow \infty} h\left(x_{s}^{k_{s}}\left(\omega^{s}\right)\right)=\lim _{s \rightarrow \infty } h\left(x_{s}^{k_{s}+1}\left(\omega^{s}\right)\right)=c \notin H^{*}
$$
so $x^{\prime} \notin X_{h}^{*}$. Let us represent
$$
x_{s}^{k+1}\left(\omega^{s}\right)=x_{s}^{k}\left(\omega^{s}\right)-\rho_{s}^ {k} P_{s}^{k}\left(\omega^{s}\right)=x_{s}^{k}\left(\omega^{s}\right)-\bar{\rho}_{s}^{ k} \bar{P}_{s}^{ k}\left(\omega^{s}\right),
$$
where
$$
\begin{aligned}
\bar{\rho}_{k} &= \rho_{k} \sum_{r=r_{k}}^{k} \frac{\lambda_{k r}}{H_{r}(\omega^{\alpha})}, \\
\bar{P}^{k}(\omega^{s}) &= \sum_{r=r_{k}}^{k} \frac{\lambda_{k r}}{H_{r}(\omega^{s})} \frac{g^{r}(\omega^{s})}{\sum_{r=r_{k}}^{k} \frac{\lambda_{k r}}{H_{r}(\omega^{s})}}.
\end{aligned}
$$

For sufficiently large $s$, for $k \in\left[k_{s}, t_{s}\right]$, and $r \ge r_{k}$ it takes place $g^{r}\left(\omega^{s}\right) = g_{h}\left(x_{s}^{r}\left(\omega^{s}\right)\right) \in G_{h}\left( x_{s}^{r}\left(\omega^{s}\right)\right)$, 
and the averaged stochastic gradient method works like a deterministic averaged gradient method. In this case,
$$
\begin{gathered}
\left\|g^{r}(\omega)\right\| \le \Gamma=\sup \{\|g(x, \theta)\|\,|\, h(x) \le d\}<+\infty, \\
\bar{P}_{s}^{k}\left(\omega^{s}\right) \in \mbox{co} \left\{G_{h}(y)\|\,|\,\| y-x^{k}\left(\omega^{s}\right) \| \le \delta_{s}^{k}\right\}, \\
\delta_{s}^{k}=\frac{\Gamma}{\nu} \rho_{s}^{r} \le \varepsilon_{s}=\frac{\Gamma}{\nu} \sup_{k \ge k_{s}} \sum_{r=r_{k}}^{k} \rho_{s}^{r} \rightarrow 0, \quad s \rightarrow \infty.
\end{gathered}
$$

Let us take $\delta>0$ such that $\max \{h(y) \mid\|y-x\| \le \delta\}<d$. Due to (24.16) and continuity of $h(y)$, the sequence $\left\{x_{s}^{k}\left(\omega^{s}\right)\right\}_{k=k_{s }}^{t_{s}}(s=0,1, \ldots)$ 
comes out of the $\delta$-neighborhood of the point $x^{\prime}$. Denote by $m_{s}$ the time of the first such exit. 
For sufficiently large $s$ and $k \in\left[k_{s}, m_{s}\right)$ it holds $ \left\|\bar{P}_{s}^{k}\right\| \le \Gamma<+\infty$ and $\sum_{k=k_{s}}^{m_{s}-1} \bar{\rho}_{s}^{k} \ge \sigma>0$. To sequences $\left\{x_{s}^{k}\left(\omega^{s}\right), k \in\left[k_{s}, m_{s}\right]\right\}$ 
Lemma 14.1 can be applied but its minimizing property 2) contradicts the construction (24.16). The lemma is proven.
\begin{lemma}
\label{lem:24.3}
Assume that the sequences $\left\{x^{k}(\omega)\right\}$ constructed according to (24.4) - (24.11) are bounded uniformly in $\omega$. Let $\bar{g}^{r}(\omega)=\mathbb{E}_{\omega} \{g^{r}(\omega) \mid x^{0}(\omega), \ldots, x^{r} (\omega)\}$ be the conditional expectation in $\omega$ of $g^{r}(\omega)$ for fixed $x^{0}(\omega), \ldots, x^{r}(\omega )$; define the quantities
\begin{equation}
\xi_{0}^{k}(\omega)=\sum_{t=0}^{k-1} \rho_{t} \sum_{r=r_{t}}^{t} \frac{\lambda_{t r}}{H_{r}(\omega)}\left(g^{r}(\omega)-\bar{g}^{r}(\omega)\right). \tag{24.17}
\end{equation}

Then the random sequence $\{\xi_{0}^{k}(\omega)\}_{k=0}^{\infty}$ $P_{\omega}$- a.s. has a limit.
\end{lemma}

{\it P r o o f.} As a rule, statements of this kind follow easily from martingale theory [133]. The complication here is that the sequence $\left\{\xi_{0}^{k}(\omega)\right\}_{k=0}^{\infty}$ is not a martingale with respect to the flow of $\sigma $-algebras generated by random variables $\left\{x^{k}(\omega)\right\}_{k=0}^{\infty}$. We artificially construct a sequence $\left\{\bar{\xi}_{0}^{k}(\omega)\right\}_{k=0}^{\infty}$ close to $\left\{\xi_{0}^{k}(\omega)\right\}_{k=0}^{\infty}$ and which is a martingale with respect to the flow of $\sigma$-algebras generated by $\left\{x^ {k}(\omega)\right\}_{k=0}^{\infty}$. Let us make the transformations:
$$
\begin{aligned}
& \xi_{0}^{k}=\sum_{t=0}^{k-1} \rho_{t} \sum_{r=r_{t}}^{t} \lambda_{t r}\left(g^{r}-\bar{g}^{-r}\right) / H_{r}=\sum_{r=0}^{k-1}\left(\sum_{t=r}^{k- 1} \lambda_{t r} \rho_{t}\right)\left(g^{r}-\bar{g}^{r}\right) / H_{r}= \\
&=\sum_{r=0}^{k-1}\left(\sum_{t=r}^{\infty} \lambda_{t r} \rho_{t}\right)\left(g^{r }-\bar{g}^{r}\right) / H_{r}-\sum_{r=0}^{k-1}\left(\sum_{t=k}^{\infty} \lambda_ {t r} \rho_{t}\right)\left(g^{r}-\bar{g}^{r}\right) / H_{r} .
\end{aligned}
$$
The random sequence
$$
\xi_{0}^{k}(\omega)=\sum_{r=0}^{k-1}\left(\sum_{t=r}^{\infty} \lambda_{t r} \rho_{ t}\right)\left(g^{r}-\bar{g}^{r}\right) / H_{r}
$$
forms a martingale with respect to $\left\{x^{k}(\omega)\right\}_{k=0}^{\infty}$. Let
$$
\Gamma=\sup_{\omega \in \Omega} \sup_{k \ge 0}\left\|g^{k}(\omega)\right\|<\infty .
$$
For random variables
$$
\delta^{k}(\omega)=\sum_{r=0}^{k-1}\left(\sum_{t=k}^{\infty} \lambda_{t r} \rho_{t}\right)\left(g^{r}-\bar{g}^{r}\right) / H_{\mathbf{r}}
$$
the estimates hold:
$$
\begin{aligned}
&\left\|\delta^{k}\right\| \le \frac{2 \Gamma}{\nu} \sum_{r=0}^{k}\left(\sum_{t=k}^{\infty} \lambda_{t r} \rho_{t}\right ) 
 =\frac{2 \Gamma}{\nu} \sum_{t=k}^{\infty} \rho_{t} \sum_{r=0}^{k} \lambda_{t r} \\
& \;\;\;\;\;\;\;\;=\frac{2 \Gamma}{\nu} \sum_{t=k}^{\infty} \rho_{t} \sum_{r=r_{t}}^{k} \lambda_{t r} \le \frac{2 \Gamma}{\nu} \sum_{t=k}^{q(k)} \rho_{t} \le \frac{2 \Gamma}{\nu} \sum_{t=r_{q(k)}}^{q(k)} \rho_{t},
\end{aligned}
$$
where the quantities $r_{t}$ and $r_{q(k)}$ are obtained from $r_{k}$ by replacing $k$ with $t$ and $q(k)$ respectively; 
$q(k)=\sup \left\{t \mid t \ge k, r_{t} \le k\right\}$. Since $q(k) \rightarrow+\infty$ as $k \rightarrow \infty$, then by condition (24.10) $\left\|\delta^{k}\right\| \rightarrow 0$ for $k \rightarrow \infty$.

It takes place
$$
\mathbb{E}\left\|\bar{\xi}_{0}^{k}(\omega)\right\|^{2} \le \frac{4 \Gamma^{2}}{\nu^{2}} \sum_{r=0}^{\infty}\left(\sum_{t=r}^{\infty} \lambda_{t r} \rho_{t}\right)^{2}<+\infty ;
$$
hence $\quad \mbox{E}\left\|\bar{\xi}_{0}^{k}(\omega)\right\| \le 1+\mbox{E}\left\|\bar{\xi}_{0}^{k}\right\|^{2}<+\infty, \quad$ and since $\quad\left\{\bar{\xi}_{0}^{k}(\omega)\right\}_{k=0}^{\infty}$ is a martingale, then the sequence $\left\{\bar{\xi}_ {0}^{k}(\omega)\right\}$, and hence the sequence $\left\{\xi_{0}^{k}(\omega)\right\} P_{\omega}$- a.s. has a limit  (see [133, p. 496]). The lemma is proven.
\begin{lemma}
\label{lem:24.4}
Let $\omega$ be such that $\left\{\xi_{0}^{k}(\omega)\right\}_{k=0}^{\infty}$ has a limit. Suppose that $\lim _{s \rightarrow \infty} x^{k_{s}}(\omega)=x(\omega) \notin X^{*}$. Denote
$$
m_{s}(\varepsilon, \omega)=\sup \left\{m \mid\left\|x^{k}(\omega)-x\right\| \le \varepsilon \quad \mbox{ for } \quad k \in\left[k_{s}, m\right)\right\}
$$
Then there exists $\bar{\varepsilon}(\omega)$ such that for any $\varepsilon \in(0, \bar{\varepsilon}]$ there are indices $l_{s}(\omega) \in\left [k_{s}(\omega), m_{s}(\varepsilon, \omega)\right]$ such, that
\begin{enumerate}
   \item $F(x(\omega))=\lim _{s \rightarrow \infty} F\left(x^{k_{s}}(\omega)\right)>
	\overline{\lim} F\left(x^ {l_{s}}(\omega)\right), \quad h(x(\omega))<0$;
   \item $h(x(\omega))=\lim _{s \rightarrow \infty} h\left(x^{k_{s}}(\omega)\right)>
	\overline{\lim} h\left(x^ {l_{s}}(\omega)\right), \quad h(x(\omega))>0$;
   \item $F(x(\omega))=\lim _{s \rightarrow \infty} F\left(x^{k_{s}}(\omega)\right)>\varlimsup_{s \rightarrow \infty} F\left(x^{l_{s}}(\omega)\right), \quad 0 \ge \varlimsup_{s \rightarrow \infty} h\left(x^{l_{s}}(\omega)\right)$, $h(x(\omega))=0$.
\end{enumerate}
\end{lemma}

{\it P r o o f.} The following representation holds:
$$
\begin{aligned}
x^{k+1}(\omega)&=x^{k_{s}}(\omega)  -\sum_{t=k_{s}}^{k} \rho_{t} \sum_{r= r_{t}}^{t} \lambda_{t r} g^{r} / H_{r}= \\
&=  x^{k_{s}}(\omega)-\sum_{t=k_{s}}^{k} \rho_{t} \sum_{r=r_{t}}^{t} \lambda_ {t r} \bar{g}^{r} / H_{r}-\xi_{k_{s}}^{k+1}(\omega)= \\
&=  x^{k_{s}}(\omega)-\sum_{t=k_{s}}^{k} \rho_{t} \bar{P}^{t}(\omega)-\xi_ {k_{s}}^{k+1}(\omega)=\bar{x}_{k_{s}}^{k+1}(\omega)-\xi_{k_{s}}^{ k+1}(\omega),
\end{aligned}
$$
where
\begin{align}
\bar{x}_{k_{s}}^{k+1}(\omega)&=\sum_{t=k_{s}}^{k} \rho_{t} \bar{P}^{t }(\omega)=\bar{x}_{k_{s}}^{k}(\omega)-\rho_{k} \bar{P}^{k}(\omega), \tag{24.18}\\
\bar{P}^{t}(\omega)&=\sum_{r=r_{t}}^{t} \lambda_{t r} \bar{g}^{r}(\omega) / H_{r }(\omega), \tag{24.19}\\
\bar{g}^{r}(\omega)&=\mathbb{E}_{\omega}\{ g^{r}(\omega) \mid x^{0}(\omega), \ldots, x^{r}(\omega)\}, \tag{24.20}\\
\xi_{n}^{m}(\omega)&=\sum_{t=n}^{m-1} \rho_{t} \sum_{r=r_{t}}^{t} \lambda_{t r }\left(g^{r}(\omega)-\bar{g}^{r}(\omega)\right) / H_{r}(\omega). \notag{}
\end{align}

Instead of $\left\{x^{k}(\omega)\right\}_{k=0}^{\infty}$, we will study the behavior of close sequences $\left\{x_{k_{s}}^{ k}(\omega)\right\}_{k \ge k_{s}}(s=0,1, \ldots)$ generated at a fixed $\omega$ by the deterministic method of averaged gradients (24.18) - (24.20), which uses the pseudo-gradients of the functions $F$ and $h$ taken not at the points $\bar{x}^{k}(\omega)$, but at close points $x^{k}(\omega)$. Then
$$
\left\|\bar{x}_{k_{s}}^{k}(\omega)-x^{k}(\omega)\right\| \le\left\|\xi_{k_{s}}^{k}(\omega)\right\| \le \sup _{k \ge k_{s}}\left\|\xi_{k_{s}}^{k}(\omega)\right\|=\delta_{k_{s}}(\omega )
$$
and $\lim _{s \rightarrow \infty} \delta_{k_{s}}(\omega)=0$, since the series $\xi_{0}^{\infty}(\omega)$ converges by the assumption of the lemma. Notice, that
$$
\begin{aligned}
& \left|F\left(x^{k}(\omega)\right)-F\left(\bar{x}_{k_{s}}^{k}(\omega)\right)\right | \le L_{F}\left\|\xi_{k_{s}}^{k}(\omega)\right\|, \\
& \left|h\left(x^{k}(\omega)\right)-h\left(\bar{x}_{k_{s}}^{k}(\omega)\right)\right | \le L_{h}\left\|\xi_{k_{s}}^{k}(\omega)\right\|,
\end{aligned}
$$
where $L_{F}$ and $L_{h}$ are the Lipschitz constants of functions $F$ and $h$ for a compact set containing $\left\{x^{k}(\omega)\right\}_{k =0}^{\infty}$. This implies that the differences $\mid F\left(x^{k}(\omega)\right)-F\left(\bar{x}_{k_{s}}^{k}(\omega)\right) \mid$ and $\left|h\left(x^{k}(\omega)\right)-h\left(\bar{x}_{k_{s}}^{k} (\omega)\right)\right|$ can be made arbitrarily small uniformly in $k$ for sufficiently large $s$.

\numberwithin{equation}{section}
\setcounter{section}{24}
\setcounter{equation}{20}
Let $\lim _{s \rightarrow \infty} \bar{x}_{k_{s}}^{k_{s}}(\omega)=\lim _{s \rightarrow \infty} x^{k_{s}}(\omega)=x(\omega), \quad 0 \notin G(x)$. Let us show that the sequences $\left\{\bar{x}_{k_{s}}^{s}(\omega)\right\}_{k \ge k_{s}}$ have the following local minimizing property.

There exists $\bar{\varepsilon}(\omega)$ such that for any $\varepsilon \in(0, \bar{\varepsilon}(\omega)]$ there exist indices $l_{s}(\omega) \ge k_{s}$ such that for sufficiently large $s$, it holds $\left\|\bar{x}_{k_{s}}^{k}(\omega)-x(\omega)\right\| \le \varepsilon$ at $k \in\left[k_{s}, l_{s}\right)$ and
\begin{enumerate}
  \item $F(x(\omega))=\lim _{s \rightarrow \infty} F\left(\bar{x}_{k_s}^{k_{s}}(\omega)\right)>\varlimsup_{s \rightarrow \infty} F\left(\bar{x}_{k_{s}}^{l_{s}}(\omega)\right), 
    \quad h(x(\omega))<0$;

  \item $h(x(\omega))=\lim _{s \rightarrow \infty} h\left(\bar{x}_{k_{s}}^{k_{s}}(\omega)\right)>\overline{\lim} _{s \rightarrow \infty} h\left(\bar{x}_{k_{s}}^{l_{s}}(\omega)\right), 
  \;\;\;
  \quad h(x(\omega))>0$;

  \item $F(x(\omega))=\lim _{s \rightarrow \infty} F\left(\bar{x}_{k s}^{k^{s}}(\omega)\right)>\varlimsup_{s \rightarrow \infty} F\left(\bar{x}_{k_{s}}^{l_{s}}(\omega)\right), 
  \\
  0\ge\lim _{s \rightarrow \infty} h\left(\bar{x}_{k_{s}}^{l_{s}}(\omega)\right), \quad h(x(\omega))=0$;
\end{enumerate}

The specified minimizing property of sequences
$\left\{\bar{x}_{k_{s}}^{k}(\omega)\right\}_{k \ge k_{s}}(s=0,1, \ldots)$ is analogous to the property of Lemma 14.2 and means simply the stability of the averaged gradient method $(24.18)-(24.20)$ with respect to the errors of the quantities included in the record of the method, namely with respect to the errors at the points where the pseudo-gradients of functions are taken. From this property of stability of the deterministic averaged gradient method, which in turn follows from the stability of the generalized gradient method (Lemma 14.1), due to the closeness of the sequences $\left\{x^{k}(\omega)\right\}_{k \ge k_{\mathrm{s}}}$ and $\left\{\bar{x}_{k_{\mathrm{s}}}^{k}(\omega)\right\}_{k \ge k_{\mathrm{s}}}$ the statement of the lemma follows.

So let us show that Lemma 14.1 is applicable to the sequences $\left\{\bar{x}_{k_{\mathrm{s}}}^{k}(\omega)\right\}_{k \geq k_{s}}(s=0,1, \ldots)$  that will imply the required minimizing property of these sequences. Let us denote 
$$
\begin{gathered}
\Gamma=\sup _{k \ge 0}\left\|g^{k}(\omega)\right\|, \quad \rho_{k_{s}}^{k}=\rho_{k}, \quad k \ge k_{s} ; \quad \bar{\rho}_{s}=\sup _{k \ge k_{\mathrm{s}}} \rho_{k} ; \\
\delta_{k_{\mathrm{s}}}^{k}=\left\|\xi_{k_{s}}^{k}\right\|+\frac{\Gamma}{\nu} \sum_{r=r_{k}}^{k} \rho_{r}, \quad \delta_{s}=\sup _{k \ge k_{s}}\left\|\bar{\xi}_{k_{s}}^{k}\right\|+\frac{\Gamma}{\nu} \sup _{k \ge k_{s}} \sum_{r=r_{k}}^{k} \rho_{r}
\end{gathered}
$$
Note that
$$
\begin{gathered}
\lim _{s \rightarrow \infty} \bar{\rho}_{s}=\lim _{s \rightarrow \infty} \delta_{s}=0 \text { and } \sum_{k \ge k_{s}} \rho_{k_{\mathrm{s}}}^{k}=+\infty, \\
\bar{P}^{k}(\omega) \in \operatorname{co}\left\{G(y) \mid\left\|y-\bar{x}_{k_{s}}^{k}(\omega)\right\| \le \delta_{k_{s}}^{k}\right\}, \quad k \ge k_{s} .
\end{gathered}
$$
So all the conditions of Lemma 14.1 are satisfied, hence the minimizing property of the sequences $\left\{\bar{x}_{k_{s}}^{k}(\omega)\right\}_{k \ge k_{s}}$ and $\left\{x^{k}(\omega)\right\}_{k \ge k_{s}}$ holds true. The lemma is proved.

{\it P r o o f of Theorem 24.1.} The boundedness of all sequences $\left\{x^{k}(\omega)\right\}_{k=0}^{\infty}(\omega \in \Omega)$ for small enough $\rho=\sup _{k \ge 0} \rho_{k}$ and $R=\sup _{k \ge 0} \sum_{r=r_{k}}^{k} \rho_{r}$  follows from Lemma 24.2. 
It follows from (24.5) - (24.10) that $\lim _{k \rightarrow \infty}\left\|x^{k+1}(\omega)-x^{k}(\omega)\right\|=0$, for all $\omega \in \Omega$. Denote by $\Omega^{\prime}$ the set of those $\omega$, for which the sequence $\left\{\xi_{0}^{k}(\omega)\right\}_{k=0}^{\infty}$ (cf. (24.17))  has a limit. By virtue of Lemma 24.3, $ \quad P_{\omega}\left(\Omega^{\prime}\right)=1$. Now for each $\omega \in \Omega^{\prime}$, using Lemma 24.4, we can prove the statement of Theorem 24.1 exactly as we did this for the deterministic method in Theorem 14.2. The theorem is proved.

Let us briefly consider stochastic analogues of other methods of nonconvex non-smooth optimization.

{\bf 2. The stochastic generalized gradient method.} 
\label{Sec.24.2}
The method is given by the following relations:
$$
\begin{gathered}
x^{0}(\omega)=x^{0} \in E_{n}, \\
x^{k+1}(\omega)=x^{k}(\omega)-\rho_{k} g^{k}(\omega) / H_{k}(\omega), \\
g^{k}(\omega)=g\left(x^{k}(\omega), \pi_{k}(\omega)\right), \quad \pi_{k}(\omega)=\theta_{k}, \quad k=0,1, \ldots, \\
0 \le \rho_{k} \le \rho, \quad \sum_{k=0}^{\infty} \rho_{k}=+\infty, \quad\sum_{k=0}^{\infty} \rho_{k}^{2}<+\infty .
\end{gathered}
$$
Here $\omega, g^{k}(\omega), H_{k}(\omega)$ are the same as in (24.4) - (24.8). Obviously, this method is a special case of the averaged stochastic gradient method; therefore, Theorem 24.1 also holds for it. To speed up the convergence of the method, the sequence $\left\{x^{k}(\omega)\right\}$ can be averaged according to (24.12), (24.14).

{\bf 3. The stochastic heavy ball method}. 
\label{Sec.24.3}
The method has the form:
\begin{eqnarray}
&x^{0}(\omega)=x^{0} \in E_{n}, \quad P^{0}(\omega)=g^{0}(\omega), \label{eqn:24.1}\\
&x^{k+1}(\omega)=x^{k}(\omega)-\rho_{k} P^{k}(\omega), \label{eqn:24.2}\\
&P^{k}(\omega)=\left(1-\gamma_{k}\right) P^{k-1}(\omega)+\gamma_{k} g^{k}(\omega) / H_{k}(\omega), \label{eqn:24.3}\\
&g^{k}(\omega)=g\left(x^{k}(\omega), \pi_{k}(\omega)\right), \pi_{k}(\omega)=\theta_{k}, \label{eqn:24.4}\\
&0<\gamma \le \gamma_{k} \le 1, \quad 0 \le \rho_{k+1} \le \rho_{k}, \quad \sum_{k=0}^{\infty} \rho_{k}=+\infty, \quad\sum_{k=0}^{\infty} \rho_{k}^{2}<+\infty.
\label{eqn:24.5}
\end{eqnarray}
Here $\omega, g^{k}(\omega), H_{h}(\omega), \pi_{k}(\omega)$ have the same meaning and satisfy the same conditions as in the method (24.4) -- (24.8) of averaged stochastic gradients. Note that $\lim _{k \rightarrow \infty} \rho_{k} / \gamma_{k}=0$.

Recall that the equation of motion (24.22) can be transformed to the following form:
$$
\begin{aligned}
x^{k+1}(\omega) & =x^{k}(\omega)-\rho_{k} \gamma_{k} g^{k}(\omega) / H_{k}(\omega)-\rho_{k}\left(1-\gamma_{k}\right) P^{k-1}(\omega)= \\
& =x^{k}(\omega)-\rho_{k} \gamma_{k} g^{k}(\omega) / H_{k}(\omega)+\frac{\rho_{k}}{\rho_{k-1}}\left(1-\gamma_{k}\right)\left(x^{k}(\omega)-x^{k-1}(\omega)\right) .
\end{aligned}
$$

For this method, the convergence assertions in the form of Theorem 24.1 are fulfilled. Let us comment this assertion.

First, the stochastic heavy-ball method (24.21) -- (24.25) can be represented as (14.44) -- (14.46). Then outside the region $D_{\varepsilon}=\{x \mid h(x) \ge \varepsilon>0\}$ for large enough $k$ it will work as a deterministic averaged gradient method, and its boundedness uniformly over $\omega$ can be provided by the return-to-start point mechanism described in $\S$ 14.

Second, the method can be treated as some modified degerministic averaged gradient method.

The following representations hold true:
$$
\begin{gathered}
P^{k}(\omega)=\sum_{r=0}^{k} \lambda_{k r} g^{r} / H_{r}, \quad \lambda_{k r} \ge 0, \quad \sum_{r=0}^k \lambda_{k r}=1, \\
\lambda_{k 0}=\prod_{i=1}^{k}\left(1-\gamma_{i}\right), \quad 
\lambda_{k k}=\gamma_{k} \quad \lambda_{k r}
=\gamma_{r} \prod_{i=r+1}^k\left(1-\gamma_{i}\right), \\
0<r<k .
\end{gathered}
$$

Now let's consider the trajectory:
$$
\begin{aligned}
& x^{k+1}(\omega)=x^{k}(\omega)- \sum_{t=k_{\mathrm{s}}}^{k} \rho_t \sum_{r=0}^{t} \lambda_{t r} g^{r} / H_{r} \\
&=x^{k_{s}}(\omega)-\sum_{t=k_{s}}^{k} \rho_{t} \sum_{r=r_{t}}^t \lambda_{t r} \bar{g}^{r} / H_{r}-\xi_{k_s}^{k+1}(\omega) \\
&=x^{k_{s}}(\omega)-\sum_{t=k_{s}}^{k} \rho_{t} \bar{P}^{t}(\omega)-\xi_{k_{\mathrm{s}}}^{k+1}(\omega)=x_{k_{\mathrm{s}}}^{k+1}(\omega)-\xi_{k_{\mathrm{s}}}^{k+1}(\omega),
\end{aligned}
$$
where
$$
\begin{aligned}
\bar{x}_{k}^{k+1}(\omega) & =\sum_{t=k_{s}}^k \rho_{t} \bar{P}^{t}(\omega)=\bar{x}_{k_{\mathrm{s}}}^{k}(\omega)-\rho_{k} \bar{P}^{k}(\omega), \\
\bar{P}^{t}(\omega) & =\sum_{r=0}^t \lambda_{t r} \bar{g}^{r}(\omega) / H_{r}(\omega), \\
\bar{g}^{r}(\omega) & =\mathbb{E}_{\omega}\{ g^{r}(\omega) \mid x^{0}(\omega), \ldots, x^{r}(\omega)\}, \\
\xi_{k_s}^{k+1}(\omega) & =\sum_{r=k_{s}}^k \rho_{t} \sum_{r=0}^t \lambda_{t r}\left(g^{r}(\omega)-\bar{g}^{r}(\omega)\right) / H_{r}(\omega).
\end{aligned}
$$

Let's represent
$$
\bar{P}^{k}=\sum_{r=0}^{k} \lambda_{k r} \bar{g}^{r} / H_{r}
=\sum_{r=r_{k}}^k \lambda_{k r} \bar{g}^{r} / H_{r}+\sum_{r=0}^{r_k-1} \lambda_
{kr} \bar{g}^{r} / H_{r}.
$$

Denote
$$
\begin{gathered}
Q^{k}=\sum_{r=r_{k}}^{k} \lambda_{k r} \frac{\bar{g}^{r}}{H_{r}} / \sum_{r=r_{k}}^{k} \lambda_{k r}, \quad \rho_{k}^{\prime}=\rho_{k} \sum_{r=r_{k}}^{k} \lambda_{k r} \\
\Delta^{k}=\sum_{r=0}^{r_k-1} \lambda_{k r} \frac{\bar{g}^{r}}{H_{r}} / \sum_{r=r_{k}}^{k} \lambda_{k r}
\end{gathered}
$$

Since $\lim _{k \rightarrow \infty} \rho_{k} / \gamma_{k} \le \lim _{k \rightarrow \infty} \rho_{k} / \gamma=0$, then there exist $r_{k}$ satisfying the conditions of Lemma 14.4 .

Suppose that the sequences $\left\{x^{k}(\omega)\right\}$ are bounded uniformly over $\omega$. Denote $\Gamma=\sup _{\omega \in \Omega} \sup _{k \ge 0}\left\|g^{k}(\omega)\right\|<\infty$.

Let us choose $r_{k}$ according to Lemma 14.4, then $\lim _{k \rightarrow \infty} \Delta^{k}=0$. With these notations in mind, we rewrite
\begin{equation}
\label{24.26}
\bar{x}_{k_{s}}^{k+1}(\omega)=\bar{x}_{k_{s}}^{k}(\omega)-\rho_{k} \bar{P}^{k}(\omega)=\bar{x}_{k_{s}}^{k}(\omega)-\rho_{k}^{\prime}\left(Q^{k}(\omega)+\Delta^{k}(\omega)\right), 
\end{equation}
\[
x^{k+1}(\omega)=\bar{x}_{k_{s}}^{k+1}(\omega)-\xi_{k_{s}}^{k+1}(\omega) .
\]

For each $\omega$, relations (24.26) can be seen as a deterministic averaged gradient method, which uses pseudogradients of problem functions taken not at points $\bar{x}_{k_{\mathrm{s}}}^{k}(\omega)$, but at close points $x^{k}(\omega)$, and averaging the gradients with some error $\Delta^{k}(\omega)$.

If we show that $\sup _{k>k_{s}}\left\||\xi_{k_s}^{k}(\omega)\right\| \rightarrow 0 \quad P_{\omega}$-a.s. as $s \rightarrow \infty$, then using the stability of the generalized gradient method (Theorem 14.1), one can prove the minimizing property and convergence $P_{\omega}$-a.s. of the stochastic heavy ball method exactly as it is done for the averaged stochastic gradient method in Lemma 24.4 and Theorem 24.1.

Thus, it has to be shown that the sequence $\left\{\xi_{0}^{k}(\omega)\right\}$
$P_{\omega}$-a.s. converges, where
$$
\xi_{0}^{k}(\omega)=\sum_{t=0}^{k-1} \rho_{t} \sum_{r=0}^{t} \lambda_{t r}\left(g^{r}(\omega)-\overline{g}^{r}(\omega)\right) / H_{r}(\omega), \quad k \ge 0 .
$$

Just as in Lemma 24.3, let us represent
$$
\begin{aligned}
& \xi_{0}^{k}=\sum_{t=0}^{k-1} \rho_{t} \sum_{r=0}^{t} \lambda_{t r}\left(g^{r}-\bar{g}^{r}\right) / H_{r}
=\sum_{r=0}^{k-1}\left(\sum_{t=r}^{k-1} \rho_{t} \lambda_{t r}\right)\left(g^{r}-\bar{g}^{r}\right) / H_{r}= \\
&=\sum_{t=0}^{k-1}\left(\sum_{t=r}^{\infty} \lambda_{t r} \rho_{t}\right)\left(g^{r}-\bar{g}^{r}\right) / H_{r}- \sum_{t=0}^{k-1}\left(\sum_{t=k}^{\infty} \lambda_{t r} \rho_{t}\right)\left(g^{r}-\bar{g}^{r}\right) / H_{r},
\end{aligned}
$$
where
$$
\begin{aligned}
& \sum_{t=r}^{\infty} \lambda_{t r} \rho_{t} \le \rho_{r} \sum_{t=r}^{\infty}(1-\gamma)^{t-r}=\rho_{r} / \gamma<+\infty; \\
& \sum_{t=k}^{\infty} \lambda_{t r} \rho_{t} \le \rho_{k} \sum_{t=r}^{\infty}(1-\gamma)^{t-r}=\rho_{k} / \gamma<+\infty \text {  at } r \le k .
\end{aligned}
$$

The random sequence
$$
\bar{\xi}_{0}^{k}(\omega)=\sum_{r=0}^{k-1}\left(\sum_{t=r}^{\infty} \lambda_{t r} \rho_{t}\right)\left(g^{r}-\bar{g}^{r}\right) / H_{r}
$$
forms a martingale with respect to the flow of $\sigma$-algebras generated by random variables $\left\{x^{k}(\omega)\right\}_{k=0}^{\infty}$. It takes place
$$
\begin{aligned}
\mathbb{E}_{\omega}\left\|\bar{\xi}_{0}^{k}(\omega)\right\|^{2} \le \frac{4 \Gamma^{2}}{\nu^{2}} & \sum_{r=0}^{\infty}\left(\sum_{t=r}^{\infty} \lambda_{t r} \rho_{t}\right)^{2} \le \frac{4 \Gamma^{2}}{\nu^{2}} \sum_{r=0}^{\infty} \rho_{r}^{2}\left(\sum_{t=r}^{\infty} \lambda_{t r}\right)^{2} \\
& \le \frac{4 \Gamma^{2}}{\nu^{2}} \sum_{r=0}^{\infty} \rho_{r}^{2}\left(\sum_{t=r}^{\infty}(1-\gamma)^{t-r}\right)^{2} \le \frac{4 \Gamma^{2}}{\gamma^{2} \nu^{2}} \sum_{r=0}^{\infty} \rho_{r}^{2}<+\infty .
\end{aligned}
$$

Therefore, the martingale $\left\{\bar{\xi}_{0}^{k}(\omega)\right\}_{k=0}^{\infty}$ has a limit $P_{\omega}$-a.s. For the random vector
$$
\delta^{k}(\omega)=\sum_{r=0}^{k-1}\left(\sum_{t=k}^{\infty} \lambda_{t r} \rho_{t}\right)\left(g^{r}-\bar{g}^{r}\right) / H_{r}
$$
the estimates are valid: 
$$
\begin{aligned}
&\left\|\delta^{k}(\omega)\right\| \le \frac{2 \Gamma}{\nu} \rho_{k} \sum_{r=0}^{k-1} \sum_{t=k}^{\infty} \lambda_{t r} \le \\
& \le \frac{2 \Gamma}{\nu} \rho_{k} \sum_{r=0}^{k-1} \sum_{t=k}^{\infty}(1-\gamma)^{t -r} \le \frac{2 \Gamma}{\gamma \nu} \rho_{k} \rightarrow 0, \;\;\;k \rightarrow \infty .
\end{aligned}
$$

Thus, the sequence $\left\{\xi_{0}^{k}(\omega)\right\}_{k=0}^{\infty}$ has a limit $P_{\omega}$-a.s. The convergence of the stochastic heavy ball method is validated.

\textbf {4. Stochastic gully step method.} 
\label{Sec.24.4}
The deterministic gully step method is discussed in detail in $\S$ 14. Its stochastic counterpart for solving $(24.1),(24.2)$ is:
\begin{equation}
\label{eqn:24.27}
x^{0}(\omega)=y^{0}(\omega)=x^{0} \in E_{n}, 
\end{equation}
\begin{equation}
\label{eqn:24.28}
y^{k+1}(\omega)=x^{k}(\omega)-\rho_{k} g^{k}(\omega), 
\end{equation}\begin{equation}
\label{eqn:24.29}
g^{k}(\omega)=g\left(x^{k}(\omega), \pi_{k}(\omega)\right), \quad \pi_{k}(\omega)=\theta_{k},
\end{equation}
\begin{equation}
\label{eqn:24.30}
x^{k+1}(\omega)=y^{k+1}(\omega)+\lambda_{k}\left(y^{k+1}(\omega)-y^{k}(\omega)\right), \quad k=0,1, \ldots,
\end{equation}
where the multipliers of the gradient step $\rho_{k}$ and the gully step $\lambda_{k}$ satisfy the conditions:
\begin{equation}
\label{eqn:24.31}
0 \le \rho_{k+1} \le \rho_{k} \le \rho, \;\;\;\sum_{k=0}^{\infty} \rho_{k}=+\infty, \;\;\;\sum_{k=0}^{\infty} \rho_{k}^{2}<+\infty, 
\end{equation}
\begin{equation}
\label{eqn:24.32}
0 \le \lambda_{k} \le \lambda<1 .
\end{equation}

Here $\theta_{k}, \omega, g^{k}(\omega), \pi_{k}(\omega)$ have the same meaning as for the averaged stochastic gradient method $(24. 4)-(24.8)$, and $x^{k}(\omega), y^{k}(\omega)$, $\rho_{k}, \lambda_{k}$ are similar to the corresponding values in the gully step method $(14.48)-(14.50)$.

Suppose that there exist indices $r_{k} \le k$ such that $\lim _{k \rightarrow \infty} r_{k}=+\infty$ and
\begin{equation}
\label{eqn:24.33}
\lim _{k \rightarrow \infty} \sum_{r=r_{k}}^{k} \rho_{r}=0, \quad \lim _{k \rightarrow \infty} \sum_{r=r_{k}}^{k} \lambda^{r_{k}-r} \rho_{r}=+\infty .
\end{equation}

In $\S$ 14 it was shown that such $r_{k}$ exists for $\rho_{k}=\mathrm{const} / k^{\alpha}$, $0<\alpha<1$.

Note that if conditions (24.33) are satisfied for some $r_{k}$, then it takes place
\begin{equation}
\label{eqn:24.34}
\lim _{k \rightarrow \infty} \sum_{r=0}^{k} \lambda^{k-r} \rho_{r}=0 \text {. }
\end{equation}
Indeed, let's represent:
$$
\sum_{r=0}^{k} \lambda^{k-r} \rho_{r}=\sum_{r=r_{k}}^{k} \lambda^{k-r} \rho_{k}+\lambda^{k-r_{k}} \sum_{r=0}^{r_{k}-1} \lambda^{r_{k}-r} \rho_{r} .
$$
The estimates hold true:
$$
\begin{aligned}
& \sum_{r=r_{k}}^{k} \lambda^{k-r} \rho_{r} \le \sum_{r=r_{k}}^{k} \rho_{k}, \quad \sum_{r=0}^{r_{k}-1} \lambda^{r_{k}-r} \rho_{r} \le \frac{\lambda \rho}{1-\lambda}<+\infty, \\
& \sum_{r=0}^{k} \lambda^{k-r} \rho_{r} \le \sum_{r=r_{k}}^{k} \rho_{r}\left(1+\frac{\lambda^{k-r_{k}}}{\sum_{r=r_{k}}^{k} \lambda^{k-r_{r}} \rho_{r}} \frac{\lambda \rho}{1-\lambda}\right) \rightarrow 0, \;\;\;
 k \rightarrow \infty \text {. }
\end{aligned}
$$
As in the deterministic case, the method (24.27)--(24.32) can be reduced to the form
$$
\begin{aligned}
& x^{0}(\omega)=x^{0} \in E_{n}, \quad x^{1}(\omega)=x^{0}-\left(1+\lambda_{0}\right) \rho_{0} g^{0}(\omega), \\
& x^{k+1}(\omega)=x^{k}(\omega)-\left(1+\lambda_{k}\right) \rho_{k} g^{k}(\omega)+ \\
& +\lambda_{k} \rho_{k-1} g^{k-1}(\omega)+\lambda_{k}\left(x^{k}(\omega)-x^{k-1}(\omega)\right) \text {. }
\end{aligned}
$$

It follows from these relations that
\begin{equation}
\label{24.35}
x^{k+1}(\omega)=x^{k}(\omega)-\rho_{k} g^{k}(\omega)-\sum_{r=0}^{k} \lambda_{r} \lambda_{r+1} \ldots \lambda_{k} \rho_{r} g^{r}(\omega) \text {. }
\end{equation}

Thus, the method $(24.27)-(24.30)$ essentially works like a stochastic averaged gradient method. Since $\lambda_{r} \le \lambda<1$, stochastic gradients $g^{r}(\omega)$ with $r$ close to $k$, i.e., stochastic gradients taken at points close to $x^{k}(\omega)$, play the main role in the right-hand side of (24.35).

For the method (24.27)-(24.30), Theorem 24.1 is also valid. The proof of this fact is done in the same way as the proof of convergence of the stochastic heavy-ball method. In particular, the boundedness of the sequence $\left\{x^{k}(\omega)\right\}$ uniformly over $\omega$ can be ensured using the return-to-initial-point mechanism described in $\S$ 14 .

Let us represent relation (24.35) as
\begin{equation}
\label{eqn:24.36}
x^{k+1}(\omega)=x^{k}(\omega)-\rho_{k} P^{k}(\omega), 
\end{equation}
$$
\begin{gathered}
P^{k}(\omega)=\sum_{r=0}^{k} \lambda_{k r} g^{r}(\omega), \\
\lambda_{k k}=1+\lambda_{k}, \quad \lambda_{k r}=\frac{\rho_{r}}{\rho_{k}} \lambda_{r} \lambda_{r+1} \ldots \lambda_{k} \le \frac{\rho_{r}}{\rho_{k}} \lambda^{k-r+1}, \quad 0 \le r<k .
\end{gathered}
$$

Starting at some point $k_{s}$, along with (24.36), we will consider the close process
$$
\bar{x}_{k_{s}}^{k_{s}}(\omega)=x^{k_{s}}(\omega)
$$
\begin{eqnarray}
 \bar{x}_{k_{s}}^{k+1}(\omega)&=&\bar{x}_{k_{s}}^{k}(\omega)-\rho_{k} \bar{P}^{k}(\omega)=x^{k_{s}}(\omega)-\sum_{t=k_{s}}^{k} \rho_{t} \bar{P}^{t}(\omega)
\nonumber\\
&=&x^{k_{s}}(\omega)-\rho_{k} \bar{g}^{k}(\omega)-\sum_{r=0}^{k} \lambda_{r} \lambda_{r+1} \ldots \lambda_{k} \rho_{r} \bar{g}^{r}(\omega), \label{24.37}
\end{eqnarray}
$$
\begin{aligned}
& \bar{P}^{k}(\omega)=\sum_{r=0}^{k} \lambda_{k r} \bar{g}^{r}(\omega), \\
& \bar{g}^{r}(\omega)=\mathbb{E}_{\omega} \{g^{r}(\omega) \mid x^{0}(\omega), \ldots, x^{k}(\omega)\}.
\end{aligned}
$$

There is the following relationship between $x^{k+1}(\omega)$ and $\bar{x}_{k_{s}}^{k+1}(\omega)$:
\begin{equation}
\label{24.38}
x^{k+1}(\omega)=\bar{x}_{k_{s}}^{k+1}(\omega)-\xi_{k_{s}}^{k+1}(\omega), 
\end{equation}
\[
\xi_{k_{s}}^{k+1}(\omega)=\sum_{t=k_{s}}^{k} \rho_{t} \sum_{r=0}^{t} \lambda_{t r}\left(g^{r}-\overline{g}^{r}\right) .
\]

Thus, we can consider the process (24.37) for each $\omega$ as if generated by the deterministic gully step method with the difference that the gradients for this method are not taken at points $\bar{x}_{k_{s}}^{k}(\omega)$ but at close points $x^{k}(\omega)$, with 
$\left\|\bar{x}_{k_{s}}^{k}(\omega)-x^{k}(\omega)\right\|=\left\|\xi_{s}^{k}(\omega)\right\|$. 
If the series $\xi_{0}^{\infty}(\omega)$ converges, then 
$\sup _{k \ge k_{s}}\left\|\xi_{k_{s}}^{k}(\omega)\right\| \rightarrow 0$ at $s \rightarrow \infty$.

The relation (24.37) can be reduced to the form (14.56) and then, using the stability of the generalized gradient method (Theorem 14.1), one can prove the minimizing property and $P_{\omega}$-a.s. convergence of the stochastic gully step method in exactly the same way as for the averaged stochastic gradient method in Lemma 24.4 and Theorem 24.1.

So, it remains to show that the series $\xi_{0}^{\infty}(\omega)$ has a limit $P_{\omega}$-a.s. Like in Lemma 24.3, let us present
$$
\xi_{0}^{k}(\omega)=\sum_{r=0}^{k-1}\left(\sum_{t=r}^{k-1} \lambda_{t r} \rho_{t}\right)\left(g^{r}-\bar{g}^{r}\right)=\bar{\xi}_{0}^{k}(\omega)-\delta^{k}(\omega),
$$
where
$$
\begin{aligned}
\bar{\xi}_{0}^{k}(\omega) & =\sum_{r=0}^{k-1}\left(\sum_{t=r}^{\infty} \lambda_{t r} \rho_{t}\right)\left(g^{r}-\bar{g}^{r}\right), \
&\delta^{k}(\omega)  =\sum_{r=0}^{k-1}\left(\sum_{t=k}^{\infty} \lambda_{t r} \rho_{t}\right)\left(g^{r}-\bar{g}^{r}\right) .
\end{aligned}
$$

Suppose that the sequence $\left\{x^{k}(\omega)\right\}_{k=0}^{\infty}$ is bounded uniformly over $\omega$. Denote $\Gamma=\sup _{\omega \in \Omega} \sup _{r \ge 0}\left\|g^{r}(\omega)\right\|<\infty$. The random sequence $\left\{\bar{\xi}_{0}^{k}(\omega)\right\}$ forms a martingale with respect to the flow of $\sigma$-algebras generated by random vectors $\left\{x^{k}(\omega)\right\}$. 
The estimates are valid: 
\begin{eqnarray}
\mathbb{E}_{\omega}\left\|\bar{\xi}_{0}^{k}(\omega)\right\|^{2} 
&\le& 4 \Gamma^{2} \sum_{r=0}^{\infty}\left(\sum_{t=r}^{\infty} \lambda_{t r} \rho_{t}\right)^{2} \nonumber\\
&\le &4 \Gamma^{2} \sum_{r=0}^{\infty} \rho_{r}^{2}\left(1+\sum_{t=r}^{\infty} \lambda^{t-r}\right)^{2} \nonumber\\
&\le& 4 \Gamma^{2}\left(1+\frac{1}{1-\lambda}\right)^{2} \sum_{r=0}^{\infty} \rho_{r}^{2}<\infty.\label{24.39}
\end{eqnarray}
Therefore, the martingale $\left\{\xi_{0}^{k}(\omega)\right\}$ has a limit $P_{\omega}$-a.s.

For the random variables $\delta^{k}(\omega)$ it is true:
\begin{eqnarray}
\left\|\delta^{k}(\omega)\right\| 
&\le& 
2 \Gamma \sum_{r=0}^{k-1} \sum_{t=k}^{\infty} \lambda_{t r} \rho_{t}  \le 2 \Gamma \sum_{r=0}^{k-1} \rho_{r} \sum_{t=k}^{\infty} \lambda^{t-r+1} \nonumber\\
& \le &2 \Gamma \sum_{r=0}^{k} \lambda^{k-r} \rho_{r} \rightarrow 0, \quad k \rightarrow \infty \quad \text { (see } \quad \text { (24.34)). }\nonumber
\end{eqnarray}

Thus, the sequence $\left\{\xi_{0}^{k}(\omega)\right\}$ has a limit $P_{\omega}$-a.s. The convergence of the stochastic gully step method is validated.

\section*{$\S$ 25. Finite difference methods in stochastic programming}
\label{Sec.25}
\setcounter{section}{25}
\setcounter{definition}{0}
\setcounter{equation}{0}
\setcounter{theorem}{0}
\setcounter{lemma}{0}
\setcounter{remark}{0}
\setcounter{corollary}{0}

\textbf{1. The stochastic approximation method.} 
\label{Sec.25.1}
This is one of the first direct methods for solving a stochastic problem on an unconditional extremum. It was first proposed in [174] to find the root of the function $F(x)=\mathbb{E} f(x, \theta)$ in the space $E_{n}$. Then in [155] similar procedures were used to minimize $F(x)$.

The Kiefer-Wolfowitz stochastic approximation method proposed for minimizing $F(x)$ is defined by the relation
\begin{equation}
\label{25.1}
x^{k+1}=x^{k}-\rho_{k} \sum_{i=1}^{n} \frac{f\left(x^{k}+\Delta_{k} e_{i}, \theta^{k}\right)-f\left(x^{k}-\Delta_{k} e_{i}, \theta^{k}\right)}{\Delta_{k}} e_{i}
\end{equation}
where $\theta^{k}$ are independent observations of $\theta$, $e_i$ are unit coordinate vectors. The convergence of the sequence (25.1) to the minimum of $F(x)$ is usually investigated under the assumption that $F(x)$ is a continuously differentiable function.

A slight modification of stochastic approximation procedures leads to stochastic methods for minimizing nonsmooth functions of probabilistic nature. The construction of these methods is related to function smoothing by a technique that was studied in Ch. 1.

Let for a fixed $\theta$ the function $f(x, \theta)$ satisfy the local Lipschitz condition, that is, for some bounded set $X$ in $E_{n}$ there exists a constant $L(\theta)\in [0, \infty)$ such that
$$
|f(x, \theta)-f(y, \theta)| \le L(\theta) \| x-y \|
$$
for any points $x, y \in X$. Then assume that $\mathbb{E} L(\theta)<\infty$. 
Let us consider the quantities
$$
F(x, \alpha)=\frac{1}{(2 \alpha)^{n}} \int_{x_{1}-\alpha}^{x_{1}+\alpha} \ldots \int_{x_{n}-\alpha}^{x_{n}+\alpha} F\left(y_{1}, \ldots, y_{n}\right) d y_{1} \ldots d y_{n},
$$
where $F(x)=\mathbb{E} f(x, \theta)$.

From the results of Ch. 2, it is easy to see that the random vector $\xi(x, \alpha)$ defined by the formula
\begin{eqnarray}
\xi(x, \alpha)&=&\frac{1}{2 \alpha} \sum_{i=1}^{n}\left[f (\tilde{x}_{1}, \ldots, x_{i}+\alpha, \ldots, \tilde{x}_{n}, \theta^{k})\right.\nonumber\\
&& \left.\;\;\;\;\;\;\;\;\;\;\;\;-f(\tilde{x}_{1}, \ldots x_{i}-\alpha, \ldots, \tilde{x}_{n}, \theta^{k})\right] e_{i}\label{25.2}
\end{eqnarray}
satisfies the condition
$$
\mathbb{E}(\xi(x, \alpha) \mid x)=\nabla{F}(x, \alpha),
$$
where $\tilde{x}_{i}\;(i=1, \ldots, n)$ are independent random variables uniformly distributed on segments $\left[x_{i}-\alpha, x_{i}+\alpha\right],$ $ \theta^{k}$ are independent observations of $\theta$. 

Analogously to formulas (4.7), (6.1), let us introduce finite differences of the following form:
\begin{equation}
\label{25.3}
\xi(x, \alpha)=\sum_{i=1}^{n} \frac{f\left(\tilde{x}+\Delta e_{i}, \theta^{k}\right)-f\left(\tilde{x}, \theta^{k}\right)}{\Delta} e_{i}, 
\end{equation}
\begin{equation}
\label{25.4}
\xi(x, \alpha)=\sum_{i=1}^{p} \frac{f\left(\tilde{x}+\Delta \mu_{i}^{k}, \theta^{k}\right)-f\left(\tilde{x}, \theta^{k}\right)}{\Delta} \mu_{l}^{k},
\end{equation}
for which
$$
\mathbb{E}\{\xi(x, \alpha) \mid x\}=\nabla F(x, \alpha)+b, \quad\|b\| \le C \Delta / \alpha.
$$

Hence, the stochastic finite-difference minimization methods for $F(x)$, generalizing the stochastic Kiefer-Wolfowitz approximation procedure (25.1), are defined by the relations
\begin{equation}
\label{25.5}
x^{k+1}=x^{k}-\rho_{k} \xi(x^{k}, \alpha_{k}),
\end{equation}
where the vector $\xi\left(x^{k}, \alpha_{k}\right)$ is calculated by using one of the formulas (25.2) -- (25.4). 
\begin{theorem}
\label{th:25.1}
Let the conditions hold:
$$
\begin{gathered}
\sum_{k=0}^{\infty} \rho_{k}=\infty, \quad \sum_{k=0}^{\infty} \rho_{k}^{2}<\infty, \quad \frac{\rho_{k}}{\alpha_{k}} \rightarrow 0, \quad \frac{\Delta_{k}}{\alpha_{k}} \rightarrow 0, \\
\alpha_{k} \rightarrow 0, \quad \frac{\left|\alpha_{k}-\alpha_{k+1}\right|}{\rho_{k}} \rightarrow 0, \quad \mathbb{E}_\theta L^{2}(\theta)<\infty .
\end{gathered}
$$
Then all limit points of the sequence $\left\{x^{k}\right\}$ with probability 1 belong to the set $X^{*} \equiv\left\{x^{*} | 0 \in \partial F\left(x^{*}\right)\right\}$; the sequence $\left\{F\left(x^{k}\right)\right\}$ converges with probability 1 .
\end{theorem}

The proof is similar to the one given for Theorem 4.1. The relation (4.10) is not quite obvious. From the condition $\mathbb{E} L^{2}(\theta)<\infty$, it follows that with probability 1,
$$
\sum_{k=0}^{\infty} \rho_{k}\left(\xi(x^{k}, \alpha_{k})-\nabla F(x^{k}, \alpha_{k})\right)<\infty,
$$
whence, with probability 1,
$$
\rho_{k}\left(\xi(x^{k}, \alpha_{k})-\nabla F(x^{k}, \alpha_{k})\right) \rightarrow 0 .
$$

A similar fact holds for formulas (25.3) through (25.4). Therefore, from the relation
$$
x^{k+1}=x^{k}-\rho_{k} \nabla F(x^{k}, \alpha_{k})+\rho_{k}\left(\nabla F(x^{k}, \alpha_{k})-\xi(x^{k}, \alpha_{k})\right),
$$
it follows (4.10). The derivation of conditions (4.12), (4.13) is done in the same way as in Theorem 4.1 .

\numberwithin{equation}{section}
\setcounter{section}{25}
\setcounter{equation}{5}
	\setcounter{page}{241}
	\setcounter{theorem}{1}
	\setcounter{equation}{5}
	\textbf{2. A stochastic game problem.} 
	\label{Sec.25.2}
	The stochastic approximation method (25.5) can be used to solve more complex problems when the classical method (25.1) is inapplicable. Stochastic game-type problems arise when choosing optimal actions under conditions where several factors influence the outcome of the decision, which can be divided into three groups: one group is controlled by the first player, the other by the second player, whose interests are opposite to those of the first player, and the third consists of random variables $\theta$.

	Let the first player makes a decision $x$ from the admissible set $X$, the second player makes a decision $y$ from the admissible set $Y$, and the state of nature $\theta$ is an elementary event of some probability space. Assume that the second player makes their choice after the first player, and he/she knows the choice  $x$ made by the first player and the state of nature $\theta$. If $f(x, y, \theta)$ represents the payment that the second player receives from the first one when nature is in state $\theta$, then the expected loss for the first player is
	\begin{equation}
	\label{eqn:25.6}
		F(x) = \mathbb{E} \max _{y \in Y} f(x, y, \theta), \quad x \in X \text {.} 
	\end{equation}
	Then the first player should choose $x$ that minimizes the function $F(x)$, subject to the condition $x \in X$.
	
	The difficulty of solving this problem lies in the fact that only in rare cases an explicit form of $F(x)$ can be found. Moreover, the function $F(x)$ is not, in general, continuously differentiable even with a sufficiently smooth $f(x, y, \theta)$.
	
	Assume that $X, Y$ are convex, bounded, and closed sets, and let there exist such a $y(x, \theta)$ for every $x$ and $\theta$ that
	$$
	f(x, y(x, \theta), \theta)=\max _{u \in Y} f(x, y, \theta),
	$$
	where $f(x, y, \theta)$ is a convex function with respect to $x$. Then, to minimize $F(x)$ (25.6), the following method is applied: 
	\begin{equation}
	\label{eqn:25.7}
	x^{k+1}=\pi_X\left(x^k-\rho_k \eta(x^k, \alpha_k)\right),  
	\end{equation} 
	in which the vector $\eta\left(x^k, \alpha_k\right)$ is determined by one of the following formulas:
	\begin{eqnarray}
		& \eta\left(x^k, \alpha_k\right)=\frac{1}{2 \alpha_k} \sum_{i=1}^n\left[f\left(\tilde{x_1^k}, \ldots, x_i^k+\alpha_k, \ldots, \tilde{x}_n^k, y\left(x^k, \theta^k\right), \theta^k\right)-\right. \nonumber\\
		& \left.-f\left(\bar{x}_1^k, \ldots, x_i^k-\alpha_k, \ldots, \bar{x}_n^k, y\left(x^k, \theta^k\right), \theta^k\right)\right] e_i,  \label{25.8}\\ 
		& \eta\left(x^k, \alpha_k\right)=\sum_{i=1}^n \frac{f\left(\bar{x}^k+\Delta_k e_i, y\left(x^k, \theta^k\right), \theta^k\right)-f\left(\tilde{x}^k, y\left(x^k, \theta^k\right), \theta^k\right)}{\Delta_k} e_i \text {, } \label{25.9}\\
		& \eta\left(x^k, \alpha_k\right)=\sum_{i=1}^p \frac{f\left(\bar{x}^k+\Delta_k \mu_i^k, y\left(x^k, \theta^k\right), \theta^k\right)-f\left(\tilde{x}^k, y\left(x^k, \theta^k\right), \theta^k\right)}{\Delta_k} e_i . \label{25.10}
	\end{eqnarray}
	
\begin{theorem}
\label{th:25.2}
	{Let the following conditions be satisfied:}
	$$
	\sum_{k=0}^{\infty} \rho_k=\infty, \quad \sum_{k=0}^{\infty} \rho_k^2<\infty, \quad \frac{\Delta_k}{\alpha_k} \rightarrow 0, \quad \alpha_k \rightarrow 0, \quad \mathbb{E}_\theta L^2(\theta)<\infty .
	$$
	{Then, with probability 1, the limit points of the sequence $\left\{x^k\right\}$ belong to the set of minima of function} (25.6). 
	\end{theorem}
	The proof is almost identical to the proof of Theorem 25.1. 
	It should be noted that:
	$$
	\mathbb{E}\left\{\eta(x^k, \alpha_k) \mid x^k, y(x^k, \theta^k), \theta^k\right\}=\nabla f\left(x^k, y(x^k, \theta^k), \theta^k, \alpha_k\right)+b_k,
	$$
	where $f(x, y(x, \theta), \theta, \alpha)$ is a smoothed function with respect to $x$ of the form $(2.1)$, $\|b_k\| \le C L\left(\theta^k\right) \Delta_k / \alpha_k$. The convergence analysis of the method is primarily based on the inequality:
	\begin{equation}
		F(x)-F\left(x^k\right) \ge \mathbb{E}\left<g\left(x^k, y(x^k, \theta^k), \theta^k \mid x^k\right), x-x^k\right>, \label{eqn:25.11}
	\end{equation}
	where $g\left(x^k, y(x^k, \theta^k), \theta^k\right) \in \partial f\left(x^k, y(x^k, \theta^k), \theta^k\right)$.
	Obviously, that
	$$
	\begin{aligned}
		& f\left(x, y(x, \theta^k), \theta^k\right)-f\left(x^k, y(x^k, \theta^k), \theta^k\right) \ge \\
		& \ge f\left(x, y(x^k, \theta^k), \theta^k\right)-f\left(x^k, y(x^k, \theta^k), \theta^k\right) \ge \\
		& \ge\left<g\left(x^k, y(x^k, \theta^k), \theta^k\right), x-x^k\right> .
	\end{aligned}
	$$
	Taking the conditional expectation of both sides of this inequality, we obtain (25.11).
	
	\textbf{3. Convergence rate.} 
	\label{Sec.25.3}
	Let the assumptions formulated in $\S$ 7 hold, and $F(x)$ is a strongly convex function,
	$$
	\mathbb{E}\left(F(x)-f(x, \theta)\right)^2 \le \sigma^2, \quad x \in X,
	$$
	$X$ is a convex compact set, and $L$ is the Lipschitz constant of the function $F(x)$ on the set $D \supset X$. We examine the behavior of the following method:
	\begin{equation}
		x^{k+1}=\pi_X\left(x^k-\rho_k \xi(x^k, \alpha_k)\right), \label{eqn:25.12}
	\end{equation}
	where $\xi\left(x^k, \alpha_k\right)$ is determined by formula (25.2).

\begin{theorem}
\label{th:25.3}	
	{Let the conditions (7.13) be satisfied, and, in addition, $\beta - 4\mu > 0$. Then, for the sequence (25.12), the following estimate holds:}
	$$
	\begin{gathered}
		\varphi(k) \;\mathbb{E}\left\|x^k-x^*\right\|^2 \le C, \\
		\varphi(k)=\left\{\begin{array}{lll}
			k^\psi, & \beta<1 & \text { or } \beta=1, \quad 2 \lambda b>\psi, \\
			k^\psi / \ln k, & \beta=1, & 2 b \lambda=\psi, \\
			k^{2 b \lambda}, & \beta=1, & 2 b \lambda<\psi, \\
		\end{array}\right.
	\end{gathered}
	$$
	$$
	\begin{aligned}
		& \psi=\min \{\mu, \beta-4 \mu\}, \\
	\end{aligned}
	$$
	
	\begin{flushleft}
		\textit{moreover, when} $2 b \lambda>\psi$, $\beta=1$
	\end{flushleft}
	$$
	\begin{aligned}
		& \lim _{p \rightarrow \infty} \sup _{k \ge p} k^\psi\, \mathbb{E}\left\|x^k-x^*\right\|^2 \le d(2 b \lambda-\psi)^{-1}, \\
	\end{aligned}
	$$
	\begin{flushleft}
		\textit{where}
	\end{flushleft}
	$$
	\begin{aligned}
		& d= \begin{cases}d_1=2 L b \sqrt{L \lambda^{-1} n \sqrt{n}}, & \text { if } \mu<1 / 5, \\
			d_2=n(\sigma b / m)^2, & \text { if } \mu>1 / 5, \\
			d_3=d_1+d_2, & \text { if } \mu=1 / 5 .\end{cases}
	\end{aligned}
	$$
\end{theorem}

	The convergence rate of the method is characterized almost surely by the next theorem.
\begin{theorem}
\label{th:25.4}	
	{Let the conditions of Theorem 7.4 be fulfilled, as well as $\psi+\beta+1>0$. Then,}
	$$
	\mathrm{P}\left\{\sup _{t \ge k} t^\nu\left\|x^t-x^*\right\| \ge a\right\} \le C\left(k^\nu / a\right)^2\left[\varphi_1(k)\right]^{-1}, \quad a>0
	$$
	\begin{flushleft}
		\textit{where}
	\end{flushleft}
	$$
	\begin{gathered}
		0 \le \nu<\min \left\{\frac{\psi+\beta-1}{2}, b \lambda\right\}, \\
		\varphi_1(k)=\left\{\begin{array}{lll}
			\varphi(k), & \text{if} & \beta=1, \\
			k^{\psi+\beta-1}, & \text {if} & \beta<1 .
		\end{array}\right.
	\end{gathered}
	$$
\end{theorem}	
	From theorems $25.3$ and $25.4$, it follows that if $\beta=1$, $10b\lambda L>3$, and $\mu=\frac{1}{5}$, then both the root-mean-square convergence rate and the almost sure convergence rate are the highest. Then
	$$
	\lim _{t \rightarrow \infty} \sup _{k \ge t} k^{1 / 5} \mathbb{E}\left\|x^k-x^*\right\|^2 \le d_3\left(\frac{2 b \lambda L}{3}-\frac{1}{5}\right)^{-1} .
	$$
	
	In the case when the descent direction is chosen according to formula (25.4), the mean-square deviation of the $k$-th approximation from the minimum point with optimally chosen parameters has an order of $k^{-1/6}$.
	
	Theorems 7.1--7.4, 25.3, 25.4 show that random noise in the objective function significantly affects the convergence rate of finite-difference methods.
	
	\section*{$\S$ 26. The averaging operation}
	\label{Sec.26}
	\setcounter{section}{26}
\setcounter{definition}{0}
\setcounter{equation}{0}
\setcounter{theorem}{0}
\setcounter{lemma}{0}
\setcounter{remark}{0}
\setcounter{corollary}{0}
	
	The averaging operation of random functions $f(x, \theta)$, gradients, and finite-difference approximations is very useful in solving a wide range of stochastic programming problems. As a rule, the structure of the methods studied in this book is as follows: a certain sequence of points is constructed $x^{k+1}=x^k+\rho_k d\left(x^k\right)$, satisfying the condition $x^k \in X$, where ${X}$ is a bounded set in $E_n$, and the vector $d\left(x^k\right)$ (in general, random) defines the direction of motion from the point $x^k$.
	
	Let's consider the sequence
	\begin{equation}
		z^{k+1}=z^k+a_k\left(f(x^k, \theta^k)-z^k\right), \label{eqn:26.1}
	\end{equation}
	\noindent where $z^0$ is the initial approximation, $a_k$ are positive multipliers, and $\theta^k$ are independent observations of the parameter $\theta$. Note that if $a_k=1/(k+1), z^0=0$, then
	$$
	z^{k+1}=\frac{f\left(x^0, \theta^0\right)+\ldots+f\left(x^k, \theta^k\right)}{k+1}.
	$$
	Therefore, it can be said that $z^k$ is obtained by the averaging operation. The difference between the averaging operation (26.1) and the usual averaging performed within the framework of the law of large numbers is that the probability distribution of the random variable $f(x, \theta)$ depends on the multidimensional parameter $x$, which changes during the process of searching for the minimum of $F(x)$.
	
	First, we will show that if the random variable $f(x, \theta)$ is bounded for all $x, \theta$, i.e., $|f(x, \theta)| \le C$ ($C$ is some constant), then the averaged estimates $z^k\;(k=0,1, \ldots)$ are also uniformly bounded. From a practical point of view, this requirement is not essential since it is automatically satisfied when implementing numerical methods on computers. Later on, the boundedness condition for $f(x, \theta)$ will be dropped.
\begin{lemma}
\label{lem:26.1}	
	{Let us assume that}
	$$
	|f(x, \theta)| \le C, \quad \sum_{k=0}^{\infty} a_k=\infty, \;\;\;a_k \rightarrow 0.
	$$
	\textit{Then, in the iterative procedure (26.1), the values $z^k$ are uniformly bounded.}
\end{lemma}	
	{\it P r o o f.} Without loss of generality, we assume that $a_k < 1$. Therefore
	$$
	\left|z^{k+1}\right| \le\left|z^k\right|\left(1-a_k\right)+C a_k
	$$
	From here, assuming convention that
	$$
	\prod_{l=k+1}^k\left(1-a_l\right)=1,
	$$
	we obtain the inequality
	\begin{equation}
		\left|z^{k+1}\right| \le\left|z^0\right| \prod_{s=0}^k\left(1-a_s\right)+C \sum_{s=0}^k a_s \prod_{l=s+1}^k\left(1-a_l\right) . \label{26.2}
	\end{equation}
	Let us show that the first term tends to zero as $k\to\infty$. From the expansion:
	$$
	e^{-a_k}=1-\frac{a_k}{1 !}+\frac{a_k^2}{2 !}-\frac{a_k^3}{3 !}+\cdots
	$$
it follows
	$$
	e^{-a_k} \ge 1-a_k.
	$$
	
	\noindent Therefore
	$$
	\varphi \equiv \exp \left\{-\sum_{s=0}^k a_s\right\} \ge \prod_{s=0}^k\left(1-a_s\right).
	$$
	The quantity $\varphi$ tends to zero since $\sum_{k=0}^\infty a_k = \infty$, implying that $\prod_{s=0}^k (1-a_s) \rightarrow 0$.
	
	Let us represent the second term of the sum (26.2) as
	$$
	\sum_{s=0}^k a_s \prod_{l=s+1}^k\left(1-a_l\right)=\sum_{s=0}^k\left(\prod_{l=s+1}^k\left(1-a_l\right)-\prod_{l=s}^k\left(1-a_l\right)\right) .
	$$
	Then we have
	$$
	\sum_{s=0}^k a_s \prod_{l=s+1}^k\left(1-a_l\right)=1-\prod_{l=0}^k\left(1-a_l\right)<1 .
	$$
	The lemma is proven.
	
	If we do not assume the boundedness of $f(x, \theta)$, then the averaged estimates $z^k$ can take arbitrarily large values. In this case, $z^k$ is defined by the formula
	$$
	z^{k+1}=\pi_Z\left(z^k+a_k\left(f(x^k, \theta^k)-z^k\right)\right),
	$$
	where $Z$ is a bounded set. In the future, this situation is not specifically discussed, and for the sake of simplicity, we assume that the estimates $z^{k+1}$ are determined by formula (26.1). We formulate the main statement of this section. Let us denote by $B_k$ the $\sigma$-algebra induced by the random vectors $x^0, z^0, \ldots, x^k, z^k$.
\begin{theorem}
\label{th:26.1}	
	{Let the following conditions be satisfied}
	$$
	\begin{gathered}
		\sum_{k=0}^{\infty} a_k=\infty, \quad \sum_{k=0}^{\infty} a_k^2<\infty, \quad \frac{\rho_k}{a_k} \rightarrow 0, \\
		\mathbb{E}\left\|d(x^k)\right\|^2<\infty, \quad \mathbb{E}_\theta f^2(x, \theta)<\infty,
	\end{gathered}
	$$
	and let the function $F(x)=\mathbb{E}_\theta f(x, \theta)$ satisfy the Lipschitz condition in the region $X$. Then the limit of the sequence $\left\{z^k-F\left(x^k\right)\right\}$ is zero with probability 1.
\end{theorem}
	
	{\it P r o o f.} The proof is carried out by contradiction, using the scheme described in $\S$ 4. First, we show that with probability 1,
	$$
	\left|z^{k+1}-z^k\right| \rightarrow 0.
	$$
	Note that
	$$
	\mathbb{E}\left\{f(x^k, \theta^k \mid B_k\right\}=F(x^k);
	$$
	therefore, with probability 1 from the conditions of the theorem it follows:
	$$
	\sum_{k=0}^{\infty} a_k\left(F(x^k)-f(x^k, \theta^k)\right)<\infty.
	$$
	
	\noindent Therefore, from the necessary condition for the convergence of a series, it follows that
	$$
	a_k\left(f(x^k, \theta^k)-F(x^k)\right) \rightarrow 0 .
	$$
	Therefore, from the equations
	$$
	\begin{aligned}
	&	\boldsymbol{z}^{k+1}=z^k+a_k\left(f(x^k, \theta^k)-z^k\right)= \\
		& =z^k+a_k\left(F(x^k)-z^k\right)+a_k\left(f(x^k, \theta^k)-F(x^k)\right)
	\end{aligned}
	$$
	it follows that with probability 1,
	$$
	\left|z^{k+1}-z^k\right| \rightarrow 0.
	$$
	
	Let us assume that there exists a convergent subsequence $\{z^s(\omega)\}$:
		$$
	\left(z^s(\omega)-F\left(x^s(\omega)\right) \rightarrow d \neq 0 .\right.
	$$
	One can specify such positive numbers $\bar{s}$ and $\bar{\delta}$ that the $2\bar{\delta}$-neighborhoods of the points $\left(z^s-F\left(x^s\right)\right), s \ge\bar{s}$, do not contain zero. Let us show that there exist  indexes $k(s)<\infty$, where
	$$
	\begin{gathered}
		k(s)=\min \left\{r \|\left(z^r-F\left(x^r\right)\right)-\left(z^s-F\left(x^s\right)\right) \mid>\delta, r>s\right\}, \\
		\delta \le \bar{\delta} / 2 .
	\end{gathered}
	$$
	
	Let us derive the inequality from which this result follows. We denote all bounded quantities by the symbol $C$. Omitting some calculations, we obtain
	\begin{eqnarray}
		\left(z^{k+1}-F(x^{k+1})\right)^2  
		&\le&\left(z^k-F(x^k)+F(x^k)-F(x^{k+1})+a_k\left(f(x^k, \theta^k)-z^k\right)\right)^2  \nonumber\\
		&\le&\left(z^k-F(x^k)\right)^2+2\left(F(x^k)-F(x^{k+1})\right)^2+2 a_k^2\left(f(x^k, \theta^k)-z^k\right)^2 \nonumber\\
		&&+2\left(z^k-F(x^k)\right)\left(F(x^k)-F(x^{k+1})\right) \nonumber\\
		&&+2 a_k\left(z^k-F(x^k)\right)\left(f(x^k, \theta^k)-z^k\right) \nonumber\\
		&\le&\left(z^k-F(x^k)\right)^2 
		+C \rho_k^2+C a_k^2+C \rho_k+T_k\nonumber\\
		&&-2 a_k\left(z^k-F(x^k)\right)^2, \label{26.3}
	\end{eqnarray}
where
	$$
	\begin{aligned}
		& T_k=C \rho_k^2\left(\|d(x^k)\|^2-\mathbb{E}\{\|d(x^k)\|^2 \mid B_k\}\right)+4 a_k^2\left(f^2(x^k, \theta^k)-\right. \\
		&\left.\quad-\mathbb{E} \{f^2(x^k, \theta^k) \mid B_k\}\right)+C \rho_k\left(\|d(x^k)\|-\mathbb{E}\{\|d(x^k)\| \mid B_k\}\right)+ \\
		&+2 a_k\left(z^k-F(x^k)\right)\left(f(x^k, \theta^k)-F(x^k)\right) .
	\end{aligned}
	$$
	Note that the inequality
	$$
	\begin{aligned}
		& \mathbb{E}\left\{\|d(x^k)\|^2 \mid B_k\right\}<\infty \\
	\end{aligned}
	$$
	implies
	$$
	\begin{aligned}
		& \mathbb{E}\left\{\|d(x^k)\| \mid B_k\right\}<\infty,
	\end{aligned}
	$$
	\noindent since
	\[
	\|d(x^k)\| \le 1+\frac{1}{4}\left\|d(x^k)\right\|^2.
	\]
	From the theorem's conditions, we obtain $\sum_{k=0}^{\infty} T_k<\infty$ with probability 1. According to the assumption, for the points $(z^k-F(x^k))$, belonging to the $\delta$-neighborhood of the point ${z}^s-F(x^s),$ $k \ge \bar{s}$, the following relation holds:
	\[
	\left(z^k-F(x^k)\right)^2>\bar{\delta}^2.
	\]
	By the theorem's condition, the sum $C\left(\rho_k^2+a_k^2+\rho_k\right)$, starting from a certain number $\bar{s}$ $(k \ge \bar{s})$, does not exceed $\bar{\delta}^2 a_k / 2$. Therefore, summing up inequality (26.3), we get
	\begin{equation}
		\left(z^k-F(x^k)\right)^2 \le \left(z^s-F(x^s)\right)^2+\sum_{r=s}^{k-1} T_r-\sum_{r=s}^{k-1} a_r .\label{26.4}
	\end{equation}
	Letting $k$ to tend to infinity, we obtain a contradiction with the conditions of boundedness of $F(x)$ in the region $X$ and boundedness of $z^k$. Consequently, there exists a finite number
	\[
	k(s)=\min \left\{r|\;|\left(z^r-F(x^r)\right)-\left(z^s-F(x^s)\right) \mid>\delta, \;r>s\right\} .
	\]
	
	From the relations
	\[
	\begin{gathered}
		\left|z^{k(s)}-z^s\right| \le  \left|\sum_{k=s}^{k(s)-1}\left[a_k\left(F(x^k)-z^k\right)+a_k\left(f( x^k, \theta^k)-F(x^k)\right)\right]\right|, \\
		\left|F(x^{k(s)})-F(x^s)\right| \le C \sum_{k=s}^{k(s)} \rho_k\left\|d(x^k)\right\|,
	\end{gathered}
	\]
	as well as from the convergence of the series
	$$
	\begin{aligned}
		& \sum_{k=0}^{\infty} \rho_k\left(\|d(x^k)\|-\mathbb{E}\{\|d(x^k)\| \mid B_k\}\right), \\
		& \sum_{k=0}^{\infty} a_k\left(f(x^k, \theta^k)-F(x^k)\right),
	\end{aligned}
	$$
	it follows that for sufficiently large numbers $s$
		$$
	C \sum_{k=s}^{k(s)-1} a_k>\frac{\delta}{2}.
	$$
	By construction,
	$$
	\left|z^{k(s)}-F(x^{k(s)})-z^s+F(x^s)\right|>\delta,
	$$
	\noindent but for sufficiently large $s$,
	$$
	\left|z^{k(s)}-F(x^{k(s)})-z^s+F(x^s)\right|<2 \delta,
	$$
	and also
	$$
	\sum_{k=s}^{k(s)-1} T_k<\frac{\bar{\delta} \delta}{4 C};
	$$
	therefore, all reasoning conducted while deriving inequality (26.4) remains valid for $k=k(s)$. Thus,
	\begin{equation}
		\left(z^{k(s)}-F(x^{k(s)})\right)^2 \le\left(z^s-F(x^s)\right)^2-\bar{\delta}^2 \delta /(4 C).\label{26.5}
	\end{equation}
	
	Taking the limit in (26.5), we obtain
		$$
	\varlimsup_{s \rightarrow \infty}\left(z^{k(s)}-F(x^{k(s)})\right)^2 \le d^2-\bar{\delta}^2 \delta /(4 C) .
	$$
	For positive numbers $a, b$ such that
	$$
	\varlimsup_{s \rightarrow \infty}\left(z^{k(s)}-F(x^{k(s)})\right)^2<a<b<d^2,
	$$
	the sequence $\left(z^k-F(x^k)\right)^2$ intersects the interval  $(a, b)$ from left to right infinitely many times. Therefore, there exists a sequence of pairs $\left\{x^r\right\}(r \in R),\left\{x^p\right\}(p \in P)$, for which
	$$
	\begin{gathered}
		\left(z^r-F\left(x^r\right)\right)^2 \le a, \quad\left(z^p-F\left(x^p\right)\right)^2 \ge b, \\
		a<\left(z^k-F\left(x^k\right)\right)^2<b, \quad r<k<p .
	\end{gathered}
	$$
	With the subsequence $\left\{\left(z^r-F\left(x^r\right)\right)^2\right\}$  we proceed in the same way as with $\left\{\left(z^s-F\left(x^s\right)\right)^2\right\}$ $(s \in S)$.  Let
	$$
	k(r)=\min \left\{t \,|\,\left(z^t-F\left(x^t\right)\right)-\left(z^r-F\left(x^r\right)\right) \mid>\delta, \;t>r\right\} .
	$$
	
	For sufficiently large $r$ and sufficiently small $\delta$  the relation $r<k(r)<p$ holds; therefore
	$$
	\left(z^{k(r)}-F(x^{k(r)})\right)^2>\left(z^r-F(x^r)\right)^2,
	$$
	which contradicts to inequality (26.5) written for the subsequence $\left\{\left(z^r-F(x^r)\right\}\right.$. The obtained contradiction proves the theorem.
	
	A result, analogous to Theorem 26.1, is valid for the averaging of $\nabla f(x, \theta)$. Let
	$$
	\begin{aligned}
		& y^{k+1}=y^k+a_k\left(\xi(k)-y^k\right), \\
		& \mathbb{E}\{ \xi(k) \mid B_k\}=\nabla F(x^k)+b_k,
	\end{aligned}
	$$
	where $b_k$ is some random noise that appears, for example, when considering the finite-difference approximations
	$$
	\xi(k)=\sum_{i=1}^n \frac{f\left(x^k+\alpha_k e_i, \theta^k\right)-f\left(x^k\right)}{\alpha_k} e_i.
	$$
	
	\noindent If the functions $f(x, \theta)$ are continuously differentiable, then the vector $\xi(k)$ can be set equal to $\nabla f\left(x^k, \theta^k\right),$ $\{\theta^k\}$ are independent observations of the parameter $\theta$.
\begin{theorem}
\label{th:26.2}	
	{Let the conditions of Theorem 26.1 be satisfied},
	$$
	\begin{gathered}
		\mathbb{E}\left\|b_k\right\| \le r_k, \quad r_k \rightarrow 0, \\
		\mathbb{E}\left\|b_k\right\|^2<\infty, \quad \mathbb{E}\|\xi(k)\|^2<\infty,
	\end{gathered}
	$$
	\textit{the gradient of the function $F(x)$ satisfies the Lipschitz condition in $X$. Then with probability 1,}
	$$
	\left(y^k-\nabla F(x^k)\right) \rightarrow 0, \quad k \rightarrow \infty .
	$$
\end{theorem}
	
	The proof of this theorem with minor modifications is carried out in the same way as for Theorem 26.1.
	
	In minimization problems with the objective function $F(x)=\mathbb{E} f(x, \theta)$, where $f(x, \theta)$ satisfies the local Lipschitz condition with constant $L(\theta)$, the averaging operation
	$$
	v^{k+1}=v^k+a_k\left(\xi(x^k, \alpha_k)-v^k\right)
	$$
	is characterized by the fact that with probability 1,
	$$
	\left(v^k-\nabla F(x^k, \alpha_k)\right) \rightarrow 0, \quad k \rightarrow \infty .
	$$
	Here $F(x, \alpha)$ is a smoothed function, and the vector $\xi\left(x^k, \alpha_k\right)$ is determined by formulas (25.2) -- (25.4). If $f(x, \theta)$ is a convex function with respect to $x$, then
	\begin{equation}
		\xi(x^k, \alpha_k)=g(\tilde{x}^k, \theta^k), \quad g(\tilde{x}^k, \theta^k) \in \partial f(\tilde{x}^k, \theta^k) . \quad\label{26.6}
	\end{equation}
\begin{theorem}
\label{th:26.3}	
	{Let the conditions of Theorem 26.1 be satisfied},
	$$
	\begin{gathered}
		\sum_{k=0}^{\infty}\left(\frac{\rho_k}{\alpha_k}\right)^2<\infty, \quad \frac{\rho_k}{\alpha_k a_k} \rightarrow 0, \quad \alpha_k \rightarrow 0, \\
		\frac{\Delta_k}{\alpha_k} \rightarrow 0, \frac{\left|\alpha_k-\alpha_{k+1}\right|}{\alpha_k a_k} \rightarrow 0, \quad \mathbb{E} L^2(\theta)<\infty .
	\end{gathered}
	$$
	\textit{Then with probability 1,}
	$$
	\left(v^k-\nabla F(x^k, \alpha_k)\right) \rightarrow 0, \quad k \rightarrow \infty .
	$$
\end{theorem}	
{\it 	P r o o f.} For simplicity of notation, we assume that
	$$
	{\mathbb{E}\{\xi }(x^k, \alpha_k) \mid B_k\}=\nabla F(x^k, \alpha_k) .
	$$
	The following chain of inequalities is fulfilled:
	$$
	\begin{aligned}
		&\left\|v^{k+1}-\nabla F(x^{k+1}, \alpha_{k+1})\right\|^2=\left\| v^k-\nabla F(x^k, \alpha_k)+a_k\left(\xi(x^k, \alpha_k)-v^k\right)+\right. \\
		&\left.+\left(\nabla F(x^k, \alpha_k)-\nabla F(x^{k+1}, \alpha_k)\right)+\left(\nabla F(x^{k+1}, \alpha_k)-\nabla F(x^{k+1}, \alpha_{k+1})\right) \right\|^2 \le
	\end{aligned}
	$$
	$$
	\begin{gathered}
		\le\left\|v^k-\nabla F(x^k, \alpha_k)\right\|^2+2 a_k\left<v^k-\nabla F(x^k, \alpha_k), \xi(x^k, \alpha_k) \pm\right. \\
		\left. \pm \nabla F(x^k, \alpha_k)-v^k\right>+2\left<v^k-\nabla F(x^k, \alpha_k), \nabla F(x^k, \alpha_k)-\nabla F(x^{k+1}, \alpha_k)\right>+ \\
		+2\left<v^k-\nabla F(x^k, \alpha_k), \nabla F(x^{k+1}, \alpha_k)-\nabla F(x^{k+1}, \alpha_{k+1})\right>+ \\
		+3 a_k^2\left\|\xi(x^k, \alpha_k)-v^k\right\|^2+3\left\|\nabla F(x^k, \alpha_k)-\nabla F(x^{k+1}, \alpha_k)\right\|^2+ \\
		+3\left\|\nabla F(x^{k+1}, \alpha_k)-\nabla F(x^{k+1}, \alpha_{k+1})\right\|^2 \le\left\|v^k-\nabla F(x^k, \alpha_k)\right\|^2- \\
		\quad-2 a_k\left\|v^k-\nabla F(x^k, \alpha_k)\right\|^2+2 a_k\left<v^k-\nabla F(x^k, \alpha_k),\right. \\
		\left.\xi(x^k, \alpha_k)-\nabla F(x^k, \alpha_k)\right>+C \frac{\rho_k}{\alpha_k}\left\|d(x^k)\right\|+\frac{C\left|\alpha_k-\alpha_{k+1}\right|}{\alpha_k}+ \\
		+C a_k^2\left\|\xi(x^k, \alpha_k)\right\|^2+C a_k^2+C\left(\frac{\rho_k}{\alpha_k}\right)^2\left\|d(x^k)\right\|^2+C \frac{\left|\alpha_k-\alpha_{k+1}\right|^2}{\alpha_k^2} \le \\
		\le\left\|v^k-\nabla F(x^k, \alpha_k)\right\|^2-2 a_k\left\|v^k-\nabla F(x^k, \alpha_k)\right\|^2
		\end{gathered}
		$$
		\begin{equation}
		\label{eqn:26.7}
		+C \frac{\rho_k}{\alpha_k}+C \frac{\left|\alpha_k-\alpha_{k+1}\right|}{\alpha_k}+C a_k^2+C \frac{\rho_k^2}{\alpha_k^2}+C \frac{\left|\alpha_k-\alpha_{k+1}\right|^2}{\alpha_k^2}+T_k,
	\end{equation}
where
	$$
	\begin{aligned}
		T_k= & 2 a_k\left<v^k-\nabla F(x^k, \alpha_k), \xi(x^k, \alpha_k)-\nabla F(x^k, \alpha_k)\right>+ \\
		+C & \frac{\rho_k}{\alpha_k}\left(\|d(x^k)\|-\mathbb{E}\{\|d(x^k)\| \mid B_k\}\right)+C a_k^2\left(\|\xi(x^k, \alpha_k)\|^2-\right. \\
		& \left.-\mathbb{E}\{\|\xi(x^k, \alpha_k)\|^2 \mid B_k\}\right)+C \frac{\rho_k^2}{\alpha_k^2}\left(\|d(x^k)\|^2-\mathbb{E}\{\|d(x^k)\|^2 \mid B_k\}\right) .
	\end{aligned}
	$$
	
	From the conditions of the theorem, it follows that $\sum_{k=0}^{\infty} T_k<\infty$. In deriving (26.7) the following inequality was used:
	$$
	\left(c_1+c_2+\ldots+c_n\right)^2 \le n\left(c_1^2+\ldots+c_n^2\right) .
	$$
	The further proof is conducted similarly that of Theorem 26.1
	
	\section*{$\S$ 27. Stochastic optimization based on the averaging operation}
	\label{Sec.27}
\setcounter{section}{27}
\setcounter{definition}{0}
\setcounter{equation}{0}
\setcounter{theorem}{0}
\setcounter{lemma}{0}
\setcounter{remark}{0}
\setcounter{corollary}{0}
	In this section, numerical methods for minimizing an expectation function are investigated:
	\begin{equation}
		F(x)=\mathbb{E} f(x, \theta)\label{27.1}
	\end{equation}
subject to the condition
	\begin{equation}
		x \in X,\label{27.2}
	\end{equation}
	where $X$ is a bounded set from $E_n$.

\textbf{1. Methods with averaging descent directions.}
\label{Sec.27.1} 
\setcounter{section}{27}
\numberwithin{equation}{section}
Let us first consider the problem of minimizing a smooth convex function $F (x)$ (27.1) with conditions (27.2), where $X$ is a convex set.

The method is defined by recurrent relations:
\begin{equation}
x^{k+1} = \pi_X(x^k - \rho_k y^k), \quad\label{27.3}
\end{equation}
\begin{equation}
y^{k+1} = y^k + a_k(\xi(k) - y^k), \quad\label{27.4}
\end{equation}
where
\[
\mathbb{E}\{\xi(k)|\,B_k\} = \nabla F(x^k) + b_k,
\]
$\pi_X (x)$ is the projection operator of point $x$ onto the set $X$. If
the function $f(x,\theta)$ for each $\theta$ is continuously differentiable with respect to $x$, then $\xi(k)$ can be taken as $\xi(k)=\nabla f(x^k,\theta ^k)$, where $\{\theta^k\}$ are independent observations of the parameter $\theta$.

The expediency of using the averaging operation is due to the fact that it is necessary to reduce the influence of noise on the search process and make it more regular. The peculiarity of the method (27.3), (27.4) is that at sufficiently large number $k$ the vector indicating the direction of motion is close to the antigradient of the target function. However, this advantage affects the initial minimization stage for a sufficiently distant initial approximation $x^0$. The averaging operation $y^k$ incorporates all the old information about the random gradients computed in previous iterations; in other words, the $y^k$ procedure has an inertia that degrades the properties of the algorithm in the neighborhood of the solution. In [100], it is shown that near the solution, methods with averaging descent directions are inferior in convergence rate to single-step procedures, e.g., stochastic quasi-gradient method.
\begin{theorem}
\label{th:27.1}
{Let the following conditions be fulfilled:}
\[
\sum^{\infty}_{k=0}\rho_k=\infty,\;\medskip\sum^{\infty}_{k=0}a^2_k<\infty,\;\frac{\rho_k}{a_k}\xrightarrow{}0,
\]
\[
\mathbb{E}||\xi(k)||^2 < \infty,\;\mathbb{E}||b_k||^2<\infty,\;\mathbb{E}||b_k||\leq r_k,
\]
\[
r_k\xrightarrow{}0,
\]
\textit{$\nabla F(x)$ satisfies the Lipshitz condition in $X$. Then with probability 1,
limit points of the sequence $\{x^k\}$ belong to the set
$X^*$ of solutions of problem (27.1), (27.2) and with probability 1},
\[
F(x^k)\xrightarrow{}F(x^*), \; x^*\in X^*.
\]
\end{theorem}
{\it P r o o f}. It is enough to show that conditions (4.12), (4.13) are fulfilled. Assume that there is a subsequence $x^s\xrightarrow{}x'\notin X^*\;(s\in S)$. It is possible to specify such positive numbers $\bar{s}$ and $\bar{\delta}$ that $2\bar{\delta}$-neighborhoods of points $x^s(s\geq\bar{s})$ are not intersect with $X^*$.
Let us construct neighborhoods of radius $\delta \leq\bar{\delta}/2$ at each point $x^s(s\geq\bar{s})$.
Further, for convenience, all limited values are denoted by the symbol $C$.\\
The inequality holds:
\begin{eqnarray}
\|x^{k+1} - x^*\|^2&\leq&\|x_k - \rho_k y^k - x^*\|^2=\|x^k-x^*\|^2 
\nonumber
\\
&&- 2\rho_k\left<y^k,x^k-x^*\right>+\rho^2_k\|y^k\|^2\leq\|x^k-x^*\|^2 
\nonumber
\\
&&- 2\rho_k\left<y^k-\nabla F(x^k), x^k-x^*\right>+2\rho_k\left<\nabla F(x^k),x^* -x^k\right>,\quad
\label{27.5}
\end{eqnarray}
where $x^* \in X^*$, the point $x^k$ belongs to the $\delta$-neighborhood of $x^s$ $(s\geq\bar{s})$. Then there is $\epsilon > 0$ such that
\[
F(x^k)-F(x^*)\geq\epsilon,
\]
\[
\left<\nabla F(x^k), x^*-x^k\right>\leq F(x^*)-F(x^k)\leq -\epsilon.
\]
It follows from Theorem 26.2 that the value
\[
-2\left<y^k - \nabla F(x^k), x^*-x^k\right>+C\rho^2_k,
\]
starting with some number $k$, does not exceed $\epsilon$. Therefore from (27.5)
it follows that
\[
\|x^{k+1}-x^*\|^2\leq\|x^k-x^*\|^2-\epsilon\rho_k.
\]

Let's put
\[
W(x^k)=\min_{x^*\in X^*}||x^k-x^*||^2=||x^k-x^*(k)||^2.
\]
Then from the previous inequality, we have
\begin{equation}
W(x^{k+1})=\|x^{k+1}-x^*(k)\|^2\leq W(x^k)-\epsilon\rho_k.\quad\label{27.6}
\end{equation}
Summing the inequality (27.6), we get
\begin{equation}
W(x^k)\leq W(x^s) - \epsilon\sum^{k-1}_{r=s}\rho_r, \quad\label{27.7}
\end{equation}
that for $k\rightarrow\infty$ contradicts the boundedness of $W(x)$ on the closed
bounded set $U_\delta(x^s)$. Therefore, the condition (4.12) is fulfilled,
 i.e. there is a finite number $k(s)$ such that
\[
k(s)=\min_{r>s} \{r\,|\,\|x^r-x^s\|>\delta\}.
\]
Inequality (27.7) is also valid for $k = k(s)$. From the relation
\[
\delta\leq\|x^{k(s)}-x^s\|\leq\sum^{k(s)-1}_{r=s}\rho_r\|y^r\|
\]
it flows
\[
\sum^{k(s)-1}_{r=s}>\frac{\delta}{C},
\]
that upon substitution in (27.7) gives
\[
W(x^{k(s)}\leq W(x^s)-\epsilon\delta/C,
\]
and this is no other than (4.13). The theorem is proved.

If $f(x, \theta)$ is a convex function for each $\theta$, then the method of solving problem (27.1)--(27.2) is determined by the relations:
\begin{equation}
x^{k+1}=\pi_X\left(x^k-\rho_k v^k\right), \quad\label{27.8}
\end{equation}
\begin{equation}
v^{k+1}=v^k + a_k(\xi(x^k,\alpha_k) - v^k), \quad
\label{27.9}
\end{equation}
vector $\xi(x^k,\alpha_k)$ is calculated by formulas (25.2) -- (25.4), (26.6).
\begin{theorem}
\label{th:27.2}
{Let the conditions of Theorem 26.3 be fulfilled, and  }
\[
\sum^{\infty}_{k=0}\rho_k = \infty.
\]
Then with probability 1, the limit points of the sequence $\{x^k\}$
belong to the set $X^*$ of solutions to problem (27.1), (27.2) and with probability 1,
\[
F(x^k)\xrightarrow{}F(x^*).
\]
\end{theorem}
{\it P r o o f}. The following relations hold true:
\[
\|x^{k+1} - x^*\|^2\leq\|x^k-\rho_k v^k-x^*\|^2=
\]
\[
=\|x^k-x^*\|^2-2\rho_k\left<v^k,x^k-x^*\right>+\rho^2_k\|v\|^2\leq
\]
\[
\leq\|x^k-x^*\|^2-2\rho_k\left<v^k-\nabla F(x^k,\alpha_k),x^k-x^*\right>-
\]
\[
-2\left<\nabla F(x^k,\alpha_k),x^k-x^*\right>+C\rho^2_k,
\]
where $F(x,\alpha)$ is a smoothed function. The following proof is carried out analogously to the proof of Theorem 27.1 due to the fact that with probability 1,
\[
(v^k-\nabla F(x^k,\alpha_k))\rightarrow{}0.
\]

If we apply the averaging operation for the unconditional minimization of the function $F(x)=\mathbb{E}f(x,\theta)$, where $f(x,\theta)$ is Lipschitz function, we obtain the method
\begin{equation}
x^{k+1}=x^k-\rho_k v^k, \quad\label{27.10}
\end{equation}
\begin{equation}
v^{k+1} = v^k+\alpha_k(\xi(x^k,\alpha_k)-v^k).\quad
\label{27.11}
\end{equation}

In relations (27.11), the vectors $\xi(x^k,\alpha_k)$ are determined by formulas (26.2) -- (26.4).
\begin{theorem}
\label{th:27.3}
{Let the conditions of Theorem 26.3 be fulfilled, and let $\sum^{\infty}_{k=0}\rho_k=\infty$. Then, with probability 1, the limit points of $\{x^k\}$ belong to the set}
\[
X^*=\{x^*|0\in\partial F(x^*)\},
\]
$F(x^k)$ converges with probability 1.
\end{theorem}
The proof of convergence of method (27.10), (27.11) is based on
the fact that with probability 1,
\[
(v^k-\nabla F(x^k,\alpha_k))\rightarrow{}0.
\]
The statement follows from inequalities
\[
F(x^{k+1},\alpha_k)=F(x^k,\alpha_k)+\left<\nabla F(x^k+\tau(x^{k+1}-x^k),\alpha ^k),x^{k+1}-x^k\right>=
\]
\[
=F(x^k,\alpha_k)+\left<\nabla F(x^k+\tau(x^{k+1}-x^k),\alpha^k)\pm\right.
\]
\[
\left.\pm \nabla F(x^k,\alpha_k),x^{k+1}-x^k\right>\leq F(x^k, \alpha_k)+
\]
\[
+C\rho^2_k/\alpha_k-\rho_k\left<\nabla F(x^k,\alpha_k),v^k\right>, \quad(0\leq\tau\leq 1),
\]
and also from the proof of Theorem 4.1.

\textbf{2. Stochastic conditional gradient methods.} 
\label{Sec.27.2}
The operation of projecting a point $y$ onto a set $X$ is equivalent to solving an extremal problem (with a quadratic objective function) of minimizing $\|y-x\|^2$ over $x\in X$. The conditional gradient method differs from the previous algorithms in that instead of the projection operation, a linearization operation is considered, which is equivalent to solving the problem with
linear objective function. Therefore, this method is effectively applicable when it is difficult to carry out the projection operation.

Let's consider a method for minimization of a continuously differentiable 
function $F(x)=\mathbb{E}f(x,\theta)$ over $x\in X$ assuming that there is
a possibility to observe a random vector $\xi(k)$ for which
\[
\mathbb{E}\{\xi(k)|B_k\}=\nabla F(x^k)+b_k.
\]

Let's define sequences of points
\begin{equation}
x^{k+1}=x^k+\rho_k(\bar{x}^k-x^k), \quad 0\leq\rho_k\leq 0 \quad\label{27.12}
\end{equation}
\begin{equation}
y^{k+1}=y^k+\alpha_k(\xi(k)-y^k), \quad\label{27.13}
\end{equation}
where $\bar{x}^k$ is the solution to the problem
\begin{equation}
\min_{x\in X}\left<y^k,x\right>.\quad\label{27.14}
\end{equation}

Simple examples show that the method in which the vector $\bar{x}^k$
is a solution to the problem
\[
\min_{x\in X}\left<\xi(k),x\right>,
\]
may not converge to the set of solutions of problem (27.1), (27.2).
\begin{theorem}
\label{th:27.4}
{Let the conditions be fulfilled:}
\[
\sum^{\infty}_{k=0}\rho_k=0,\quad \sum^{\infty}_{k=0}a^2_k<\infty,\quad \frac{\rho_k}{a_k }\xrightarrow{}0,
\]
\[
\mathbb{E}\|\xi(k)\|^2<\infty,\quad \mathbb{E}\|b_k\|^2<\infty,\quad \mathbb{E}\|b_k\|\leq r_k,\quad r_k\xrightarrow{}0 ,
\]
\textit{the gradient of the function $F(x)$ satisfies the Lipshitz condition in $X$. Then
with probability 1, the limit points of the sequence $\{x^k\}$ (27.12) --
(27.14) belong to the set $X^*$ satisfying the necessary
conditions of extremum for problem (27.1), (27.2), the sequence $\{F(x^k)\}$
converges with probability 1.}
\end{theorem}
{\it P r o o f}. The proof largely coincides with the proof of the conditional gradient method
 considered in $\S$ 15. From the mean value theorem, we have
\begin{eqnarray}
F(x^{k+1})&=&F(x^k)+\left<\nabla F(x^k+\tau(x^{k+1}-x^k),x^{k+1}-x^k)\right>
\nonumber\\
&=&F(x^k)+\rho_k\left<\nabla F(x^k+\tau(x^{k+1}-x^k))-\nabla F(x^k),\bar{x}^k-x^k\right>
\nonumber\\
&&+\rho_k\left<\nabla F(x^k),\bar{x}^k-x^k\right>\nonumber\\
&\leq& F(x^k)+C\rho^2_k+\rho_k\left<\nabla F(x^k),\bar{x}^k-x^k\right>,\quad\label{27.15}
\end{eqnarray}
\[
0\leq\tau\leq 1
\]
In inequality (27.15), we need to estimate the value $\rho_k\left<\nabla F(x^k),\bar{x}^k-x^k\right>$.
To do this, consider the auxiliary problem of minimization over $x$ the
functions $\left<\nabla F(z), x\right>$ under the conditions $x\in X$, $z\in X$. Let $\bar{z}$ be the solution of the given problem. For points $z$ belonging to a sufficiently small neighborhood of the point $z^\prime\notin X^*$, the inequality holds:
\begin{equation}
\left<\nabla F(z),\bar{z}-z\right>\leq\sigma<0.\quad\label{27.16}
\end{equation}

Assume that $x^s\xrightarrow{}x^\prime\notin X^*$ $(s\in S)$. Let the condition (4.12) is not fulfilled, i.e., all points $x^k\;(k\geq s)$ are contained in a sufficiently small $\delta$-neighborhood of the point $x^s$, which does not intersect with $X^*$ . From the inequality (27.16), we have
\[
\left<\nabla F(x^k),\tilde{x}^k-x^k\right>\leq\sigma<0,
\]
where $\tilde{x}^k$ is the solution of the previous auxiliary problem. Since with probability 1,
\[
(y^k-\nabla F(x^k))\xrightarrow{}0,
\]
then for sufficiently large $k\geq s$,
\[
\left<y^k,\tilde{x}^k-x^k\right>\leq\sigma/2.
\]
Therefore,
\[
\left<y^k,\bar{x}^k-x^k\right>\leq\sigma/2,
\]
where $\bar{x}^k$ is the solution of problem (27.14). Hence, by virtue of the previous remark, we get
\[
\left<\nabla F(x^k),\bar{x}^k-x^k\right>\leq\sigma/4.
\]
Summing inequality (27.16) and taking into account that $\rho_k\xrightarrow{}0$, for sufficiently large $s$, we have
\begin{equation}
F(x^k)\leq F(x^s)+\frac{\sigma}{8}\sum^{k-1}_{r=s}\rho_r. \label{27.17}
\end{equation}
Therefore, passing to the limit on $k\xrightarrow{}\infty$, we obtain the contradiction to the boundedness of $F(x)$. Thus, condition (4.12) is fulfilled,
i.e. there exists
\[
k(s)=\min_{r>s}\{r\,|\,\|x^r-x^s\|>\delta\}.
\]
It follows directly from the inequality (27.17) that
\[
\overline{\lim}_{x\rightarrow\infty}F(x^{k(s)})<\lim_{x\rightarrow{}\infty}F(x^s),
\]
and this means that condition (4.13) is valid. The theorem is proved.

Let us now discuss the expectation minimization problem:
\[
F(x)=\mathbb{E}f(x,\theta)\rightarrow\min_{x\in X},
\]
where $f(x,\theta)$ satisfies in $D\supset X$ the Lipshitz condition with a constant
$L(\theta)$. The considered method is a stochastic variant of the conditional gradient method described in $\S$ 15. The following sequences of points are defined:
\begin{equation}
x^{k+1}=x^k+\rho_k(\bar{x}^k-x^k),\quad\label{27.18}
\end{equation}
\begin{equation}
v^{k+1}=v^k+\alpha_k(\xi(x^k, \alpha_k)-v^k),\quad
\label{27.19}
\end{equation}
\begin{equation}
\left<v^k,\bar{x}^k\right>=\min_{x\in X}\left<v^k,x\right>.\quad\label{27.20}
\end{equation}
Here $\xi(x^k, \alpha_k)$ are random vectors calculated by formulas
(25.2) -- (25.4). If $f(x,\theta)$ is a convex function in $x$, then
\[
\xi(x^k, \alpha_k)=g(\tilde{x}^k, \theta_k), \quad g(\tilde{x}^k, \theta_k)\in\partial f(\tilde{x }^k, \theta_k).
\]
\begin{theorem}
\label{th:27.5}
{Let the conditions be fulfilled:}
\[
\sum^{\infty}_{k=0}\rho_k=\infty,\quad
\sum^{\infty}_{k=0}a^2_k<\infty,\quad
\sum^{\infty}_{k=0}\left(\frac{\rho_k}{\alpha_k}\right)^2<\infty,\quad
\alpha_k\xrightarrow{}0,
\]
\[
\frac{\rho_k}{\alpha_k a_k}\xrightarrow{}0,\quad
\frac{|\alpha_k-\alpha_{k+1}|}{\alpha_k a_k}\xrightarrow{}0,\quad
\mathbb{E}L^2(\theta)<\infty.
\]
\textit{Then the number sequence $\{F(x^k)\}$ converges with probability 1 and the limit points of the sequence $\{x^k\}$ (27.18) -- (27.20) with probability 1 belong to the set $X^*$ that satisfies the necessary conditions of
extremum for problem (27.1), (27.2).}
\end{theorem}

The proof completely coincides with the proof of Theorem 15.1.

\smallskip
\textbf{3. The rate of convergence of the stochastic conditional gradient method.}
\label{Sec.27.3}
We give an estimate of the rate of convergence of the method (27.12) -- (27.14) without proof. Let $C_L, C_\xi, C_b, C_\rho, C_a, r, t, \beta$ be positive constants,
\[
\alpha=\min(r-t,\frac{t}{2},\beta), \quad
0<t<r\leq 1,
\]
and the gradient of function $F(x)$ satisfies the inequality
\[
\|\nabla F(x)-\nabla F(y)\|\leq C_L\|x-y\|, \quad x,y\in X.
\]

Let the following conditions be fulfilled for $k\geq N$:
\[
\rho_k=C_\rho k^{-r},\quad a^k=C_a k^{-t},
\]
\[
\alpha N^{-(1-r)}<C_\rho\leq N^r, \quad
2\alpha N^{t-1}<2C_a\leq N^t,
\]
\[
\|\xi(k)\|\leq C_\xi, \quad \|b_k\|\leq C_b k^{-\beta}.
\]
The number $N$ is introduced in order to select the constants $C_\rho$ and $C_a$ large enough, since the conditions $\rho_h\leq 1$, $a_k\leq 1$ are used in the proof of the theorem.
\begin{theorem}
\label{th:27.6}
{If $k\geq N$, then the estimate holds true}
\begin{equation}
\mathbb{E}(F(x^k)-F(x^*))\leq C_1 k^{-\alpha},\quad x^*\in X^*, \quad \label{27.21}
\end{equation}
where
\[
C_1 =\max\left\{(\mathbb{E}F(x^N)-F(x^*))N^\alpha,\;
(2C_x C_a C_\rho + \frac{1}{2}C_L C^2_x C^2_\rho)(C_\rho-\alpha N^{-(1-r)})^{-1}\right\},
\]
\[
C_x=\max_{x,y\in X}\|x-y\|.
\]
\end{theorem}

We choose the constants $r$ and $t$ so that the exponent $\alpha$ in formula (27.21) is maximal, i.e., we solve the problem
\[
\alpha\equiv\min\left(r-t,\frac{t}{2},\beta\right)\xrightarrow{}\max, \quad 0<t<r\leq 1.
\]
It is easy to notice that the optimal solution is:
\[
r=1,\quad t=\alpha/3,\quad \alpha=\min(1/3,\beta).
\]

\textbf{4. Examples}.
\label{Sec.27.4}
We illustrate the application of the conditional stochastic gradient  method for solving some stochastic programming problems.

\textit{Multi-product task.} Let's assume that we need to do
an order $x=(x_1,...,x_n)$ for the supply of $n$ heterogeneous products, provided that the demand is characterized by a random vector $\theta=(\theta_1,...,\theta_n)$
and the coefficient of interchangeability $\lambda_{js}$ of the $s$-th demanded product by  the $j$-th supply product (the order vector may consist of products not included in the demand vector). We present the value of $x_j$ as $x_j=\sum^r_{s=1}x_{js},$ and
consider the loss function $f(x, \theta)$. If the loss function takes into account
only costs associated with the underuse of some products (in the case of $\sum^n_{j=1}\lambda_{js}x_{js}>\theta_s$) and losses from the deficit of others (in the case of $\sum^n_{ j=1}\lambda_{js}x_{js}\leq\theta_s$), then it is natural to assume that
\[
f(x,\theta)=\sum^r_{s=1}\max\Big\{\alpha_s\Big(\sum^n_{j=1}\lambda_{js}x_{js}-\theta_s\Big),\;
\beta_s\Big(\theta_s-\sum^n_{j=1}\lambda_{js}x_{js}\Big)\Big\},
\]
where $\alpha_s$ is the cost of storing a unit of  $s$-th product, $\beta_s$ is the cost of its deficit. Then the mathematical expectation of losses is equal to
\begin{equation}
F(x)=\sum^r_{s=1}\mathbb{E}\max\Big\{\alpha_s\Big(\sum^n_{j=1}\lambda_{js}x_{js}-\theta_s\Big) ,\;
\beta_s\Big(\theta_s-\sum^n_{j=1}\lambda_{js}x_{js}\Big)\Big\}.\quad\quad\quad
\label{27.22}
\end{equation}
We need to find  vector $x$ that minimizes (27.22), under the conditions:
\begin{equation}
\sum^n_{j=1}a_jx_j\leq a, \quad x_j=\sum^r_{s=1}x_{js},\quad x_{js}\geq 0,\quad\quad\label{27.23}
\end{equation}
where $a_j$ is the volume of one unit of $j$-th product, and $a$ is the warehouse capacity.

Note that the linear programming problem (27.20) has the
trivial solution: all $x_{js}$ are equal to zero, except for $x_{js}=a/a_j$, for which the minimum is achieved:
\[
\min_{j,s}\Big\{v_{js}\frac{a}{a_j}\Big\}<0,\quad j=1,...,n, \;\;\;s=1,... ,r,
\]
($v_{js}$ is determined by formula (27.19)).

\textit{One-product problem with transportation costs.} Let there be a homogeneous product placed in $m$ warehouses. Permissible
amount of product in the warehouse $i$ is $a_i\;(i=1, ..., m)$. In the considered period, the product will be required in $n$ markets, and the demand of market $j$ cannot be determined in advance, and is considered as a random variable $\theta_j$. Denote by
\[
y_j=\sum^m_{i=1}x_{ij},\quad\quad j=1,...,n,
\]
the amount of product transported to the $j$-th market. Then the loss under  demand  $\theta_j$, equals to
\[
f_j(y_j,\theta_J) =
\left\{
\begin{array}{ll}
\alpha_j\big(\sum^m_{j=1}x_{ij}-\theta_j\big), & \quad \sum^m_{j=1}x_{ij}\geq\theta_j, \\
\beta_j\big(\theta_j-\sum^m_{j=1}x_{ij}\big), & \quad \sum^m_{j=1}x_{ij}<\theta_j.
\end{array}
\right.
\]
It is necessary to determine the transportation plan $x=(x_{ij},i = 1,...m, j = 1, ...,n)$, which minimizes the sum of transport costs and expected losses
\begin{equation}
F(x)=\sum_{i,j}a_{ij}x_{ij}+\sum^n_{s=1}\mathbb{E}\max\Big\{\alpha_j\Big(\sum^m_{i=1 }x_{ij}-\theta_j\Big),\;
\beta_j\Big(\theta_j-\sum^m_{i=1}x_{ij}\Big)\Big\}\quad\label{27.24}
\end{equation}
subject to constraints
\begin{equation}
\sum^n_{i=1}x_{ij}\leq a_j,\quad i=1,...,m, \quad x_{ij}\geq 0. \quad \label{27.25}
\end{equation}

Due to the specificity of constraints (27.25), the linear programming problem   in the conditional gradient method is subdivided into $m$ subproblems,
each of which has a trivial solution: for each $i$ only one value of $x_{ij} = a_j$ is greater than zero, the others $x_{ij}$ are equal to zero.

\textit{A stochastic model of allocation of agricultural crops.} Let's consider one model of the optimal distribution of sowing plots taking into account the fact that the yield on different plots is random.

We use the following notations: let $n$ be number of cultures, $m$ be
the number of plots of land of different types, $S_j$ be the $j$-th area
plot, $a_{ij}$ $(\theta)$ be the yield of the $i$-th crop on the $j$-th plot, $\gamma_i$ be the specific weight of the $i$-th crop in the total production under consideration of
agricultural production, $x_{ij}$ is the area occupied by
$i$-th culture on the $j$-th site.

The objective is to determine a crop allocation plan that satisfies the area constraints for each land type and maximizes the expected value of the quantity of agricultural products in specified proportions, i.e.,
\[
f(x)\equiv \mathbb{E}\Big(\min_i \frac{1}{\gamma_i}\sum^m_{j=1}a_{ij}(\theta)x_{ij}\Big)\xrightarrow{}\max_x ,
\]
\[
\sum^n_{i=1}x_{ij}\leq S_j, \quad\quad j=1,...,m,
\]
\[
x_{ij}\geq 0, \quad i=1,...,n.
\]
Due to the properties of the constraints of the described model, the linear programming problem in the stochastic conditional gradient method is divided into $t$ subproblems,  each of which has an obvious solution.

\textbf{5. The stochastic method of the reduced gradient.} 
\label{Sec.27.5}
Now we investigate the stochastic variant of the method described in $\S$ 16 for solving the minimization problem
\[
F(x)=\mathbb{E}f(x,\theta)
\]
subject to the conditions
\[
x\in X = \{x: Ax = b,\quad x\geq 0\},
\]
where $A$ is an $m \times n$-matrix, $b$ is a vector from $E_m$. For the set
$X$ the prerequisites formulated in $\S$ 16 are fulfilled.

A stochastic variant of the reduced gradient method for solution
problems (27.1), (27.2), when $f(x,\theta)$ is a Lipschitz function, is defined as follows.

\textit{I n i t i a l \;\;s t a g e}. Choose an admissible point $x^1$ satisfying the conditions $Ax = b$, $x\geq 0$. Set $k = 1$ and go
to Step 1.

\textit{S t e p 1.} Choose  indices $I_k$ of $m$ largest components of the vector $x^k$,
\[
x_B = \{x_i,i\in I_k\},\quad B=\{a_j, j\in I_k\},\quad N = \{a_j, j\notin I_k\}.
\]
Calculate the reduced gradient
\begin{equation}
r^T_N = (v^k_N)^T-(v^k_B)^TB^{-1}N.\quad\quad \label{27.26}
\end{equation}
The descent direction $d^k$ is determined:
\begin{equation}
d^k_j =\left\{
\begin{array}{lll}
     -r_j, \quad & j\notin I_k, \quad &r_j\leq 0, \\
     -x_jr_j, \quad & j\notin I_k, \quad& r_j > 0,
\end{array}\right.
\label{eqn:27.27}
\end{equation}
\[
d^k_B = -B^{-1}Nd_N.
\]

\textit{S t e p 2.} Sequences of points $x^k, v^k$ are recalculated using the formulas
\[
x^{k+1} = x^k + \sigma_kd^k,
\]
\[
v^{k+1} = v^k + \alpha_k(\xi(x^k,\alpha_k)-v^k),
\]
\begin{equation}
\sigma_k =\min(\rho_k, \lambda),\label{eqn:27.28}
\end{equation}
\[
\lambda = \left\{
\begin{array}{ll}
    \min\Big\{-\frac{x^k_j}{d^k_j}\Big|d^k_j<0\Big\},\quad & \exists d^k_j < 0, \\
    \quad \infty, \quad & d^k\geq 0.
\end{array}\right.
\]
Set $k = k+1$ and proceed to Step 1. The vectors $\xi(x^k,\alpha_k)$ are calculated according to formulas (25.2) -- (25.4). If $f(x,\theta)$ is a convex function with respect to $x$, then
\[
\xi(x^k,\alpha_k) = g(\tilde{x}^k,\theta^k).
\]
\begin{theorem}
\label{th:27.7}
{Let the conditions of Theorem 27.5 be fulfilled. Then the sequence $\{F(x^k)\}$ converges with probability 1 and the limits
points of the sequence $\{x^k\}$ with probability 1 satisfy the Kuhn-Tucker conditions.}
\end{theorem}

The proof is carried out similarly to the proof of Theorem 16.1.

Let us now consider the situation when the function $F(x)$ is continuously
differentiable and there exists a random vector $\xi(k)$ such that
\[
\mathbb{E}\{\xi(k)|B_k\} = \nabla F(x^k) + b_k.
\]
\noindent
In particular, if $f(x, \theta)$ is continuously differentiable with respect to $x$ at each $\theta$, then we can take $\xi(k)=\nabla f\left(x^{k}, \theta^{k}\right)$. In this case, the method (27.26) -- (27.28) is simplified. No averaging of $v^{k}$ is necessary, and the reduced gradient $r_{N}$ in (27.26) is determined by the formula
\begin{equation}\tag{27.26}
r_{N}^{T}=\xi_{N}^{T}(k)-\xi_{B}^{T}(k) B^{-1} N .
\end{equation}
\begin{theorem}
\label{th:27.8}
Let the conditions are fulfilled:
$$
\begin{gathered}
\sum_{k=0}^{\infty} \rho_{k}=\infty, \sum_{k=0}^{\infty} \rho_{k}^{2}<\infty, \quad \mathbb{E}\left\|b_{k}\right\| \le r_{k}, \quad r_{k} \rightarrow 0, \\
\mathbb{E}\left\|b_{k}\right\|^{2}<\infty, \quad \mathbb{E}\|\xi(k)\|^{2}<\infty,
\end{gathered}
$$
\noindent
{$\nabla{F}(x)$ satisfies the Lipschitz condition in $X$. Then $\{F\left(x^{k}\right)\}$ converges with probability 1 and the limit points of the sequence $\{x^{k}\}$, defined by formulas (27.26) -- (27.28), satisfy the Kuhn-Tucker conditions with probability 1}.
\end{theorem}

\textbf{6. Stochastic methods of feasible directions.} 
\label{Sec.27.6}
In the previous paragraphs of this section, we considered the problem of minimizing $F(x) \equiv \mathbb{E} f(x, \theta)$ under the constraint $x \in X$. In the case where the set $X$ is of simple structure, the projection operation on $X$ was used. In the conditional gradient method, the projection operation was replaced by minimizing in $X$  some linear form. In the reduced gradient method, the values of the basis variables were recalculated by changing the independent out-of-basis variables according to the requirement of admissibility preservation.

Let us now discuss a more complicated problem when the region $X$ is defined by nonlinear constraints $f_{i}(x) \le 0\;(i=1, \ldots, m)$. We assume that $X$ is bounded and the regularity condition is satisfied.

Let $F(x), f_{i}(x)\;(i=1, \ldots, m)$ be continuously differentiable functions. To solve the problem (27.1), (27.2), consider the following stochastic variant of the method of possible directions given in $\S$ 17:
\begin{equation}\tag{27.29}
x^{k+1}=x^{k}-\gamma_{k} d^{k},
\end{equation}
\noindent
where the possible direction $d^{k}$ at a given admissible point $x^{k}$ is found from the solution of the following linear programming problem:
\begin{equation}\tag{27.30}
\min_{d,\sigma} \sigma
\end{equation}
\noindent
subject to conditions
\begin{equation}\tag{27.31}
\begin{gathered}
\left<y^{k}, d\right>-\sigma \le 0, \\
f_{i}(x^{k})+\left<\nabla f_{i}(x^{k}), d\right>-\sigma \le 0, \quad i=1, \ldots, m, \\
-1 \le d_{j} \le 1, \quad j=1, \ldots, n .
\end{gathered}
\end{equation}
\noindent
The sequence $y^{k}$ is found using the averaging operation
\begin{equation}\tag{27.32}
y^{k+1}=y^{k}+a_{h}\left(\xi(k)-y^{k}\right),
\end{equation}
\noindent
where $\xi(k)$ is a random vector such that
$$
\mathbb{E}\{ \xi(k) | B_{k}\}=\nabla F(x^{k})+b_{k}.
$$
Let us set
\begin{equation}\tag{27.33}
\begin{aligned}
& \rho_{k}^{\prime}=\max \left\{\lambda \mid x^{k}+\lambda d^{k} \in X\right\}, \\
& \gamma_{k}=\min \left\{\rho_{k}^\prime, \rho_{k}\right\} .
\end{aligned}
\end{equation}
\begin{theorem}
\label{th:27.9}
{Let the conditions of Theorem 27.4 be satisfied, then with probability 1 the limit points of the sequence $\{x^{k}\}$ defined by formulas (27.29) -- (27.33) satisfy the Kuhn-Tucker conditions, the sequence $\{F\left(x^{k}\right)\}$ converges with probability 1.}
\end{theorem}

{\it P r o o f.} From the mean value theorem it follows
\begin{equation}\tag{27.34}
F(x^{k+1}) \le F(x^{k})+C \gamma_{k}^{2}+\gamma_{k}\left<\nabla F(x^{k}), d^{k}\right>.
\end{equation}
\noindent
In inequality (27.34), we need to evaluate $\left<\nabla F(x^{k}), d^{k}\right>$. This can be done by considering the auxiliary problem:
$$
\min_{d,\sigma} \sigma
$$
\noindent
subject to conditions
\begin{equation}\tag{27.35}
\begin{gathered}
\left<\nabla F(z), d\right> \le \sigma, \\
f_{i}(z)+\left<\nabla f_{i}(z), d\right> \le \sigma, \quad i=1, \ldots, m, \\
-1 \le d_{j} \le 1 .
\end{gathered}
\end{equation}

Let $(\bar{\sigma}(z), \bar{d}(z))$ be a solution of this problem. Then the following property holds: for the points $z$ belonging to a sufficiently small neighborhood of the point $z^{\prime} \notin X^{*}$ (the set of Kuhn-Tucker points), the inequality $\bar{\sigma}(z)<\sigma<0$ is satisfied. It follows that
$$
\left<\nabla F(z), \bar{d}(z)\right> \le \sigma<0.
$$

Obviously, to prove the theorem it is sufficient to check conditions (4.12), (4.13). Let there exists a subsequence
$$
x^{s} \rightarrow x^{\prime} \notin X^{*}, \quad s \in S \subset\{1,2, \ldots, k, \ldots\} .
$$
\noindent
There exists such an integer $\bar{s}$ and a positive number $\bar{\delta}$ such that $U_{2 \bar{\delta}}\left(x^{s}\right) \cap X^{*}=\emptyset$ $(s \ge \bar{s})$. Therefore,
$$
\begin{gathered}
\left<\nabla F\left(x^{s}\right), \bar{d}\left(x^{s}\right)\right> \le \sigma<0, \\
f_{i}\left(x^{s}\right)+\left<\nabla f_{i}\left(x^{s}\right), \bar{d}\left(x^{s}\right)\right> \le \sigma, \;\;\;i=1,\ldots,m;\\
-1 \le \bar{d}_{j}\left(x^{s}\right) \le 1, \quad j=1, \ldots, n,
\end{gathered}
$$
\noindent
where $\bar{d}\left(x^{s}\right)$ is the solution of problem (27.35).

Suppose that condition $(4.12)$ is not satisfied, i.e., all points $x^{k}$ $(k \ge s)$ belong to $U_{\delta}\left(x^{s}\right)(\delta \le \bar{\delta} / 2)$. Then for sufficiently small $\delta$ it follows from the previous inequalities that
$$
\begin{gathered}
\left<\nabla F(x^{k}), \bar{d}\left(x^{s}\right)\right> \le \sigma / 2, \\
f_{i}(x^{k})+\left<\nabla f_{i}(x^{k}),  \bar{d}\left(x^{s}\right)\right> \le \sigma / 2, 
\;\;\;i=1,\ldots,m;\\
-1 \le \bar{d}_{j}\left(x^{s}\right) \le 1, \quad j=1, \ldots, n .
\end{gathered}
$$

Since with probability 1 $\left(y^{k}-\nabla F\left(x^{k}\right)\right) \rightarrow 0$, then starting at some number $k \ge s$,
$$
\left(y^{k}, \bar{d}\left(x^{s}\right)\right) \le \sigma / 4
$$
\noindent
Therefore, for the optimal solution $d^{k}$ of the problem $(27.30)-(27.31)$ the following relations hold:
\begin{equation}\tag{27.36}
\begin{gathered}
\left<y^{k}, d^{k}\right> \le \sigma / 4, \\
f_{i}(x^{k})+\left<\nabla f_{i}(x^{k}), d^{k}\right> \le \sigma / 4, 
\;\;\;i=1,\ldots,m;\\
-1 \le d^{k} \le 1, \quad j=1, \ldots, n .
\end{gathered}
\end{equation}
\noindent
From inequalities (27.36) it is easy to show that the step size $\gamma_{k}$ is equal to $\rho_{k}$. Hence, from inequality (27.34) for sufficiently large $k$, we have
\begin{equation}\tag{27.37}
F(x^{k}) \le F(x^{s})+\frac{\sigma}{8} \sum_{r=s}^{k-1} \rho_{r} .
\end{equation}

Coming to the limit on $k \rightarrow \infty$, we obtain a contradiction. Therefore
$$
k(s)=\min _{r>s}\left\{r \mid\left\|x^{r}-x^{s}\right\|>\delta\right\}<\infty \text {. }
$$
\noindent
In turn, from the relation
$$
x^{k(s)}=x^{s}+\sum_{r=s}^{k(s)-1} \rho_{r} d^{r}
$$
\noindent
it follows that
$$
\sum_{r=s}^{k(s)-1} \rho_{r}>\frac{\delta}{C} \text {. }
$$
\noindent
Hence, from inequality (27.37) we conclude that condition (4.13) is satisfied:
$$
\varlimsup_{s \rightarrow \infty} F\left(x^{k(s)}\right)<\lim _{s \rightarrow \infty} F\left(x^{s}\right)
$$
\noindent
The theorem is proved.

Let us now consider the problem of minimization of the function
$$
F(x) \equiv \mathbb{E} f(x, \theta)
$$
\noindent
assuming that $f(x, \theta)$ satisfies the Lipschitz condition in $x$ on $D \supset X$ with constant $L(\theta)$ ; $f_{t}(x)$ $(i=1, \ldots, m)$ are continuously differentiable functions.

Let us define the sequences:
$$
\begin{aligned}
& x^{k+1}=x^{k}+\gamma_{k} d^{k}, \\
& v^{k+1}=v^{k}+a_{k}\left(\xi(x^{k}, \alpha_{k})-v^{k}\right),
\end{aligned}
$$
\noindent
where $d^{k}$ is the solution of the problem
$$
\min_{\sigma,d} \sigma
$$
\noindent
subject to
$$
\begin{gathered}
\left<v^{k}, d\right>-\sigma \le 0, \\
f_{i}(x^{k})+\left<\nabla f_{i}(x^{k}), d\right>-\sigma \le 0, \quad i=1, \ldots, m_{,} \\
-1 \le d_{j} \le 1, \quad \ddot{j}=1, \ldots, n .
\end{gathered}
$$

The step size $\gamma_{k}$ is chosen according to formula (27.33). In finite-difference methods for minimizing $F(x)$, the vectors $\xi\left(x^{k}, \alpha_{k}\right)$ are computed using the formulas $(25.2)-(25.4)$. If $f(x, \theta)$ is a convex function, then in the stochastic quasi-gradient method
$$
\xi(x^{k}, \alpha_{k})=g(\tilde{x}^{k}, \theta^{k}) \in \partial f(\tilde{x}^{k}, \theta^{k}) .
$$
\begin{theorem}
\label{th:27.10}
{Suppose that the conditions of Theorem 27.5 are satisfied, then $\left\{F\left(x^{k}\right)\right\}$ converges with probability 1 and the limit points of the sequence $\left\{x^{k}\right\}$ satisfy the Kuhn-Tucker condition with probability 1}.
\end{theorem}

The proof is similar to the proofs of Theorems 17.1 and 27.9.

\smallskip
\textbf{7. Stochastic methods of minimization of Lipschitz functions}.
\label{Sec.27.7}
Let it be required to solve the problem:
\begin{equation}\tag{27.38}
F(x) \equiv \mathbb{E} f(x, \theta) \rightarrow \min_x
\end{equation}
\noindent
subject to
\begin{equation}\tag{27.39}
h(x) \le 0,
\end{equation}
\noindent
where $f(x, \theta), h(x)$ are Lipschitz functions. The method for solving the problem (27.38), (27.39) is defined by the relations:
$$
x^{k+1}=x^{k}+\rho_{k} d^{k},
$$
\noindent
where 
$$
\begin{array}{ll}
d^{k}=-z_{F}^{k}, & h\left(x^{k}\right)<0, \\
d^{k}=-z_{F}^{k}-u_{k} z_{h}^{k}, & h\left(x^{k}\right) \ge 0, \\
u_{k}  =\max \left\{0, \frac{h\left(x^{k}\right)-\left<z_{F}^{k}, z_{h}^{k}\right>}{\left\|z_{h}^{k}\right\|^{2}}\right\} .&
\end{array}
$$
\noindent
Vectors $z_{F}^{k}, z_{b}^{k}$ are calculated by the formulas:
$$
\begin{aligned}
& z_{F}^{k+1}=z_{F}^{k}+a_{k}\left(\xi_{F}(x^{b}, \alpha_{k})-z_{F}^{k}\right), \\
& z_{h}^{k+1}=z_{h}^{k}+a_{k}\left(\xi_{h}(x^{k}, \alpha_{k})-z_{h}^{k}\right) .
\end{aligned}
$$

\noindent
In finite-difference methods $\xi_{F}, \xi_{h}$ are determined by relations (25.2) to (25.4). If $f,$ $h$ are convex functions, then in the stochastic quasi-gradient method
$$
\begin{aligned}
& \xi_{F}(x^{k}, \alpha_{k})=g(\tilde{x}^{k}, \theta_{k}) \in \partial f(\tilde{x}^{k}, \theta^{k}), \\
& \xi_{h}(x^{k}, \alpha_{h})=g(\tilde{x}^{k}) \in \partial h(\tilde{x}^{k}) .
\end{aligned}
$$

Let us denote by $X^{*}$ the set of points satisfying the necessary  condition of extremum for  problem (27.38), (27.39): $x^{*} \in X^{*}$ if there are such vectors $g_{F}\left(x^{*}\right) \in \partial F\left(x^{*}\right),\; g_{h}\left(x^{*}\right) \in \partial h\left(x^{*}\right)$ that
$$
\begin{gathered}
g_{F}\left(x^{*}\right)+\lambda g_{h}\left(x^{*}\right)=0, \quad \lambda \ge 0, \\
\lambda h\left(x^{*}\right)=0, \quad h\left(x^{*}\right) \le 0 .
\end{gathered}
$$

The following result holds.
\begin{theorem}
\label{th:27.11}
{Let the conditions of Theorem 27.5 be satisfied. Then $\{F\left(x^{k}\right)\}$ converges with probability 1 and the limit points of the sequence $\{x^{k}\}$ belong to $X^{*}$ with probability 1.}
\end{theorem}

\section*{$\S$ 28. The stochastic finite-difference Arrow-Hurwicz method}
\label{Sec.28}
\setcounter{section}{28}
\setcounter{definition}{0}
\setcounter{equation}{0}
\setcounter{theorem}{0}
\setcounter{lemma}{0}
\setcounter{remark}{0}
\setcounter{corollary}{0}
Let us now consider the solution of a general stochastic programming problem:
\begin{equation}\tag{28.1}
F_{0}(x) \equiv \mathbb{E} f_{0}(x, \theta) \rightarrow \min_x
\end{equation}
subject to constraints
\begin{equation}\tag{28.2}
F_{i}(x) \equiv \mathbb{E} f_{i}(x, \theta) \le 0, \quad i=1, \ldots, m,
\end{equation}
\begin{equation}\tag{28.3}
x \in X,
\end{equation}
where $f_{i}(x, \theta)$ $(i=0,1, \ldots, m)$ are convex Lipschitz in $x$ (at each $\theta$) on $D \supset X$ functions, $X$ is a convex bounded closed set in $E_{n}$. The method under study is defined by the relations:
\begin{equation}\tag{28.4}
x^{k+1}=\pi_{X}\left(x^{k}-\rho_{k}\left[\xi^{0}(x^{k}, d^{k})+\sum_{i=1}^{m} u_{i}^{k} \xi^{i}(x^{k}, \alpha_{k})\right]\right),
\end{equation}
\begin{equation}\tag{28.5}
u^{k+1}=\pi_{U}\left(u^{k}+\rho_{k} L_{u}(x^{k}, u^{k}, \theta^{k})\right),
\end{equation}
where $L_{u}(x, u, \theta)$ is the gradient of the function
$$
L(x, u, \theta)=f_{0}(x, \theta)+\sum_{i=1}^{m} u_{i} f_{i}(x, \theta)
$$
over variables $u$.

The vectors $\xi^{\nu}\left(x^{k}, \alpha_{k}\right)\;(\nu=0,1, \ldots, m)$ are calculated using one of the formulas:
\begin{equation}\tag{28.6}
\begin{aligned}
\xi^{\nu}(x^{k}, \alpha_{k})=\frac{1}{2 \alpha_{k}} \sum_{i=1}^{n}\left[f_{\nu}\left(\tilde{x_{1}^{k}}, \ldots, x_{i}+\alpha_{k}, \ldots, \tilde{x}_{n}^{k}, \theta^{k}\right)-\right. \\
\left.-f_{\nu}\left(\tilde{x}_{1}^{k}, \ldots, x_{i}-\alpha_{k}, \ldots, \tilde{x}_{n}^{k}, \theta^{k}\right)\right] e_{i},
\end{aligned}
\end{equation}
\begin{equation}\tag{28.7}
\xi^{\nu}(x^{k}, \alpha_{k})=\sum_{i=1}^{n} \frac{f_{\nu}\left(x^{k}+\Delta_{k} e_{i}, \theta^{k}\right)-f_{\nu} \tilde{\left(x^{k}, \theta^{k}\right)}}{\Delta_{k}} e_{i} \text {. }
\end{equation}
\noindent
Convergence of trajectories $\left\{x^{k}\right\},\left\{u^{k}\right\}$ to the set of saddle points of the Lagrangian function
$$
L(x, u)=F_{0}(x)+\sum_{i=1}^{m} u_{i} F_{i}(x)
$$
\noindent
is explored in the Cesaro sense.

Let us denote
$$
\begin{aligned}
& \hat{x}^{k}=\sum_{s=0}^{k} \rho_{s} x^{s} / \sum_{s=0}^{k} \rho_{s}, \\
& \hat{u}^{k}=\sum_{s=0}^{k} \rho_{s} u^{s} / \sum_{s=0}^{k} \rho_{s},
\end{aligned}
$$
\noindent
$K_{\nu}(\theta)$ is the Lipschitz constant of the function $f_{\nu}(x, \theta)$ $(\nu=0,1, \ldots, m)$ on $D$.
\begin{theorem}
\label{th:28.1}
{Let the conditions are fulfilled:}
$$
\begin{gathered}
\sum_{k=0}^{\infty} \rho_{k}=\infty, \quad \sum_{k=0}^{\infty} \rho_{k}^{2}<\infty, \quad \frac{\Delta_{k}}{\alpha_{k}} \rightarrow 0, \quad \sum_{k=0}^{\infty} \rho_{k} \alpha_{k}<\infty, \quad \alpha_{k} \rightarrow 0, \\
\mathbb{E} K_{\nu}^{2}(\theta)<\infty, \quad \mathbb{E} f_{\nu}^{2}(x, \theta)<\infty, \quad \nu=0,1, \ldots, m .
\end{gathered}
$$
{Then with probability 1 the limit points of the sequences $\left\{\hat{x}^{k}\right\},\left\{\hat{u}^{k}\right\}$ belong to the set of saddle points of $L(x, u)$.}
\end{theorem}

We give a proof for the case where the vectors $\xi^{\nu}\left(x^{k}, \alpha_{k}\right)$ $(\nu=0,1, \ldots, m)$ are defined by formulas (28.6). We denote the bounded quantities by the symbol $C$.

Let
$$
\begin{gathered}
v^{k}=\xi_{0}(x^{k}, \alpha_{k})+\sum_{i=1}^{m} u_{i}^{k} \xi_{i}(x^{k}, \alpha_{k}), \\
\gamma_{k}=\xi_{0}(x^{k}, \alpha_{k})-\nabla F_{0}(x^{k}, \alpha_{k})+\sum_{i=1}^{m} u_{l}^{k}\left(\xi_{i}(x^{k}, \alpha_{k})-\nabla F_{i}(x^{k}, \alpha_{k})\right),
\end{gathered}
$$
$\left(X^{*}, U^{*}\right)$ is the set of saddle points of the function 
$L(x, u),\; x^{*} \in X^{*}$.

{ The following inequalities hold true:} 
\begin{equation}\tag{28.8}
\begin{aligned}
& \left\|x^{k+1}-x^{*}\right\|^{2} \le\left\|x^{k}-\rho_{k} v^{k}-x^{*}\right\|^{2} \le \\
& \le\left\|x^{k}-x^{*}\right\|^{2}+2 \rho_{k}\left<v^{k}, x^{k}-x^{n}\right>+\rho_{k}^{2}\left\|v^{k}\right\|^{2} \le \\
& \le\left\|x^{k}-x^{*}\right\|^{2}+2 \rho_{k}\left<\nabla F_{0}\left(x^{k}, \alpha_{k}\right)+\sum_{i=1}^{m} u_{i}^{k}\nabla F_{i}(x^{k}, \alpha_{k}), x^{*}-x^{k}\right>+ \\
& +2 \rho_{k}\left<\gamma_{k}, x^{*}-x^{k}\right>+C \rho_{k}^{2}+\rho_{k}^{2}\left(\left\|v^{k}\right\|^{2}-\mbox{E}\left\|v^{k}\right\|^{2}\right) \\
& \le\left\|x^{k}-x^{*}\right\|^{2}+2 \rho_{k}\left(\bar{L}\left(x^{*}\right)-L\left(x^{k}, u^{k}\right)\right)+C \rho_{k} \alpha_{k}+ \\
& +C \rho_{k}^{2}+2 \rho_{k}\left<\gamma_{k}, x^{*}-x^{k}\right>+\rho_{k}^{2}\left(\left\|v^{k}\right\|^{2}-\mbox{E}\left\|v^{k}\right\|^{2}\right),
\end{aligned}
\end{equation}
\noindent
where $\bar{L}(x)=\max _{u \in U} L(x, u)$.

For dual variables $u$, the inequalities are satisfied:
\begin{equation}\tag{28.9}
\begin{gathered}
\left\|u^{k+1}-u\right\|^{2} \le\left\|u^{k}+\rho_{k} L_{u}\left(x^{k}, u^{k}, \theta^{k}\right)-u\right\|^{2} \le\left\|u^{k}-u\right\|^{2}+ \\
+2 \rho_{k}\left(L\left(x^{k}, u^{k}, \theta^{k}\right)-L\left(x^{k}, u, \theta^{k}\right)\right)+\rho_{k}^{2}\left\|L_{u}\left(x^{k}, u^{k}, \theta^{k}\right)\right\|^{2} \le \\
\le\left\|u^{k}-u\right\|^{2}+2 \rho_{k}\left(L\left(x^{k}, u^{k}, \theta^{k}\right)-L\left(x^{k}, u, \theta^{k}\right)\right)+ \\
+C \rho_{k}^{2}+\rho_{k}^{2}\left(\left\|L_{u}\left(x^{k}, u^{k}, \theta^{k}\right)\right\|^{2}-\mathbb{E}\left\|L_{u}\left(x^{k}, u^{k}, \theta^{k}\right)\right\|^{2}\right) .
\end{gathered}
\end{equation}

Adding inequalities $(28.8)$ and (28.9), we obtain
$$
\begin{aligned}
\left\|x^{k}-x^{*}\right\|^{2} & +\left\|u^{k+1}-u\right\|^{2} \le\left\|x^{k}-x^{*}\right\|^{2}+\left\|u^{k}-u\right\|^{2}+ \\
+ & 2 \rho_{k}\left(\bar{L}\left(x^{k}\right)-L\left(x^{k}, u\right)\right)+2 \rho_{k}\left(L\left(x^{k}, u^{k}, \theta^{k}\right)-L\left(x^{k}, u^{k}\right)\right)+ \\
& +2 \rho_{k}\left(L\left(x^{k}, u\right)-L\left(x^{k}, u, \theta^{k}\right)\right)+2 \rho_{k}\left<\gamma_{k}, x^{*}-x^{k}\right>+ \\
& +\rho_{k}^{2}\left(\left\|v^{k}\right\|^{2}-\mathbb{E}\left\|v^{k}\right\|^{2}\right)+\rho_{k}^{2}\left(\left\|L_{u}\left(x^{k}, u^{k}, \theta^{k}\right)\right\|^{2}-\right. \\
& \left.-\mathbb{E}\left\|L_{u}\left(x^{k}, u^{k}, \theta^{k}\right)\right\|^{2}\right)+C \rho_{k}^{2}+C \rho_{k} \alpha_{k} .
\end{aligned}
$$
From here it follows that
\begin{equation}\tag{28.10}
\begin{aligned}
&\sum_{s=0}^{k} \rho_{s}\left(L(\hat{x^{k}}, u)-\bar{L}\left(x^{s}\right)\right) \le C+\frac{1}{2} \sum_{s=0}^{k} C\left(\rho_{s}^{2}+\rho_{s} \alpha_{s}\right)+ \\
&+\frac{1}{2} \sum_{s=0}^{k}\Big\{2 \rho_{s}  \left(L\left(x^{s}, u^{s}, \theta^{s}\right)-L\left(x^{s}, u^{s}\right)\right)+2 \rho_{s}\left(L\left(x^{s}, u\right)-L\left(x^{s}, u, \theta^{s}\right)\right)+ \\
& +2 \rho_{s}\left<\gamma_{s}, x^{*}-x^{s}\right>+\rho_{s}^{2}\left(\left\|v^{s}\right\|^{2}-\mathbb{E}\left\|v^{s}\right\|^{2}\right)+ \\
& +\rho_{s}^{2}\left(\left\|L_{u}\left(x^{s}, u^{s}, \theta^{s}\right)\right\|^{2}-\mathbb{E}\left\|L_{u}\left(x^{s}, u^{s}, \theta^{s}\right)\right\|^{2}\right)\Big\}.
\end{aligned}
\end{equation}
Hence, $\bar{L}\left(\hat{x}^{k}\right) \rightarrow \bar{L}\left(x^{*}\right)$ with probability 1, since the series in relation (28.10) converge. Similarly, we obtain that with probability 1,
$$
\bar{L}\left(\hat{u}^{k}\right) \rightarrow \max _{u} L(u), \quad \bar{L}(u)=\min _{x \in X} L(x, u).
$$
The theorem is proved.

\newpage
\section*{Bibliographical Notes}
\label{BibNotes}

\small

The following literature references are by no means intended to be exhaustive: they are mainly those works that are directly related to the material presented, as well as the main monographs; monographs are given preference over articles.

\textbf{Preface}

General optimization problems are considered in the works of V. M. Glushkov [16], N. N. Krasovsky [54], V. S. Mikhalevnch and V. L. Volkovich [62], N. N. Moiseev [67], L. S. Pontryagin, V. G. Boltyansky, R. V. Gamkrelidze and E. F. Mishchenko [101], A. N. Tikhonov and V. Y. Arsenin [118].

The theory of linear programming and linear algebra can be found in the monographs by Dantzig [32], D. K. Faddeev and V. N. Faddeeva [120], and Murtagh [71]. The theory and methods of nonlinear programming are presented in the books by Bazaraa and Shetty [6], F. P. Vasiliev [11], Y. G. Evtushenko [36], I. I. Eremin and N. N. Astafiev [38], Zangwill [46], V. G. Karmanov [51], N. N. Moiseev, Y. P. Ivanilov, and E. M. Stolyarova [68], Ortega and Reinboldt [92], Polak [98], B. T. Polyak [100], B. N. Pshenichny and Y. M. Danilin [105], Himmelblau [125], and the collection [131].

{\textbf{Chapter 1}

Elements of the theory of convex analysis are contained in almost every textbook on optimization, exhaustive material is available in the monograph by Rockafellar [109], a great deal of space is devoted in the books by V. F. Demyanov and L. V. Vasiliev [34], B. N. Pshenichny $[104]$.

$\S$ 1. Various classes of nonsmooth functions have been studied in the works of Danskin [31], V. F. Demyanov and V. N. Malozemov [33], V. F. Demyanov and L. V. Vasiliev [34], A. D. Ioffe and V. M. Tikhomirov [49], V. I. Norkin [80, 83], E. A. Nurminskii [88], B. N. Pshenichny [103, 104], N. Z. Shor [136], Clarke [145], Leuenberger [162], Mangasarian [163], Mirica [169], Mifflin [167], Rockafellar [176-178], and Bihain [143].

Generalized gradients have been introduced for non-convex locally Lipschitz [145], quasi-differentiable [34, 103], admitting upper convex approximation [104], generalized differentiable [80, 83], weakly convex [88] functions. Various properties of generalized gradients have been studied for locally Lipschitz [146], qusi-differentiable [34], and generalized differentiable [83] functions. Numerical methods for minimizing non-convex nonsmooth functions using generalized gradients have been proposed for Lipschitz [22], generalized differentiable [80], weakly convex [88] and semi-smooth [166] functions.

$\S$ 2. Differential properties of locally Lipschitz functions were systematically studied by Clarke [145]; he introduced the notions of "generalized directional derivative" and "subdifferential". The exposition of $\S 2$ is based on [146, 22]. The operation of smoothing or averaging functions has long been used in mathematics. For the purposes of nonsmooth optimization it began to be used by A. M. Goupal [23] and V. J. Katkovnik [52]. Properties of locally Lipschitz functions were also studied in the works of Rockafellar [176-178], Hiriart-Urruti [153]. A theorem on the mean was established in [159]. A nontraditional definition of locally Lipschitz functions was given by Mirica [169]. $\S$ 3. The problematics of the theory of necessary extremum conditions for nonsmooth functions can be found in the monographs of B. N. Pshenichny [103, 104]. Necessary extremum conditions for problems with locally Lipschitz functions were studied by Clarke [146], Warga [10], A. M. Gupal [22], Hiriart-Urruti [152, 153], and Mifflin [167]. Specific necessary extremum conditions are established for quasi-differentiable functions [34]. The question of boundedness of Kuhn-Tucker multipliers under more general conditions was studied in [150].

\textbf{Chapter 2}

$\S$ 4. Finite-difference methods for minimizing locally Lipschitz functions were proposed by A. M. Gupal in [23] and developed in [22]. Numerical minimization methods  based on the smoothing operation were studied in the book by V. J. Katkovnik [52]. Numerical methods of optimization without computing gradients were studied for smooth problems in the book by B. N. Pshenichny and Y. M. Danilin [105], and for nonsmooth convex problems in the book by A. S. Nemirovsky and D. B. Yudin [73].

$\S$ 5. Methods for solving minimax problems were considered in the books by Danskin [31], V. F. Demyanov and V. N. Malozemov [33], V. V. Fedorov [122], N. Z. Shor [136], and E. A. Nurminskii [91]. The construction of finite differences for maximum functions is done in [17].

$\S$ 6. Random search methods can be found in the book by L. A. Rastrigin [106]; the convergence and convergence rate of these methods were investigated by V. G. Karmanov [51]. The result of $\S$ 6 can be found in the book [22].

$\S$ 7. The complexity of problems and the efficiency of convex optimization methods were investigated by A. S. Nemirovsky and D. B. Yudin [73]. The results of $\S$ 7 are contained in [27]. The convergence rate of finite-difference methods was investigated by G. G. Murauskas [70].

$\S$ 8. The method of Lyapunov functions for investigating the convergence of relaxation optimization methods is given in Bazaraa and Shetty [6], Y. G. Evtushenko [36], and B. T. Polyak [100]. For non-relaxation algorithms this method is generalized in [5, 46, 87, 91].

The first necessary and sufficient conditions for convergence of nonlinear programming algorithms were formulated by Zangwill [46]. Based on these results, E. A. Nurminskii [87, 91] proposed very general sufficient conditions for convergence of nonlinear programming algorithms, which found wide application in the study of non-relaxation algorithms of nonlinear and stochastic programming. These results were then generalized by many authors. In this paragraph, necessary and sufficient convergence conditions are formulated in Theorems 8.1-8.3; Theorem 8.1 is the most general. L. G. Bazhenov [5] transformed the necessary and sufficient conditions of Zangwill [46] to a form convenient for application to non-relaxation algorithms of nonlinear programming and got rid of formal assumptions to the contrary. S. P. Uryasiev obtained necessary and sufficient convergence conditions close to Theorem 8.3.

Sufficient conditions for convergence in this paragraph are formulated in Theorems 8.4-8.9. Condition 1) of Theorem 8.8 is a relaxation of a similar requirement of E. A. Nurminskii [87]; it was used by V. I. Norkin [80]. Its further weakening (Condition 1) of Theorem 8.6) was found by S. P. Uryasiev. Condition 2) of Theorem 8.8 was introduced by E. A. Nurminskii [89] and weakened (Condition 2) of Theorem 8.6) by S. I. Lyashko [59]. Condition 3) of Theorem 8.8 is equivalent to the fact that $W^{*}$ is not dense anywhere in 
$E_{1}$ [80], and is a generalization of the requirement of E. A. Nurminskii [87] that $W^{*}$ is countable at most. Originally, E. A. Nurminskii [87] proved approximately in the conditions of Theorem 8.8 that all limit points of the sequence $\left\{x^{k}\right\}$ belong to $X^{*}$. He also showed [91] that without condition 3) of Theorem 8.8, there exist limit points $\left\{x^{k}\right\}$ belonging to $X^{*}$, and V. I. Norkin [84] pointed out that these are the minimal in $W(x)$ limit points of $\left\{x^{k}\right\}$. Y. M. Ermoliev [41] and Y. M. Ermoliev and P. I. Verchenko [42] noticed that under the conditions of E. A. Nurminskii [87], $\left\{W\left(x^{k}\right)\right\}$ has a limit. M. V. Mikhalevich [65] showed that without condition 3) in Theorem 8.8, $\overline{\lim} _{k \rightarrow \infty}W\left(x^{k}\right) \le \sup \left\{w \in E_{1} \mid w \in W^{*}\right\}$.
The next step was taken by P. A. Dorofeev [35]; he showed that under the conditions of E. A. Nurminskii [87] without assumptions on the set $W^{*}$ there is convergence of $\left\{x^{k}\right\}$ by one limit point and by a functional. The connectivity of the set of limit points $\left\{x^{k}\right\}$ under the conditions of Lemma 8.3 is noted in [92, p. 459]. In the work [44] S. K. Zavriev  proposed other general sufficient conditions for convergence of non-convex optimization algorithms. In one part these conditions are more general than the sufficient convergence conditions of this paragraph, and in another part they are more specific: there are no formal assumptions to the contrary. However, the character of the statements about convergence of algorithms obtained in [44] is the same; in particular, he obtained the convergence conditions for non-convex optimization algorithms on a functional, proved a theorem similar to Theorem 8.8, and proved Lemma 8.3.

$\S$ 9, 10. The generalized gradient method was proposed by N. Z. Shor [134]; its convergence conditions for convex functions are studied [40, 99, 136]. The method has been generalized and improved in [18, 34, 38,  73,136]; it has also been extended to solve some non-convex nonsmooth optimization problems [4, 80, 88, 136, 166]. Generalized gradient methods with space stretching are studied in detail in [136]. On the basis of the ellipsoid method [141], which is a special case of the method [135], L. G. Khachnyan [124] constructed a polynomial algorithm in linear programming.

Problems with importance-ordered criteria (lexicographic optimization problems) are considered in the monographs by V. V. Podinovskii and V. M. Gavrilov [97], V. V. Fedorov [122], and in [121].

Improper problems are studied in the monograph by I. I. Eremin, V. D. Mazurov, and H. N. Astafiev [39].

The books [13,136] describe numerous applications of the generalized gradient method and its generalizations.

$\S$ 11. Relaxation methods for minimizing a nonsmooth convex maximum function are described in the monograph by V. F. Demyanov and V. N. Malozemov [33].  Relaxation methods for minimizing general convex functions are proposed by Lemareschal [160] and Wolfe [185]; relaxation methods for convex optimnization are also considered in the papers [107, 108, 130,157, 157, 161, 166, 172, 179]. Relaxation methods for minimizing semi-smooth functions under constraints are proposed by Mifflin [166, 168] and developed by Kiwil [156, 157]. The method of this paragraph fits into the general scheme given in [172].

$\S$ 12. Reviews of global optimization methods were made by McCormick [165], P. G. Strongin [114], and A. G. Sukharev [117] (see also the collection [181]). A notable place is devoted to global optimization problems in the books by D. I. Batischev [7], V. P. Bulatov [8], F. P. Vasiliev [11], I. B. Motskus [69], L. A. Rastrigin [106], and A. G. Sukharev [116].

Combined global optimization algorithms have been studied in [1, 11, 12, 69, 82, 86].

The method of nonconvex approximations was proposed and justified in the works of S. A. Piyavskii [95, 96], Yu. M. Danilin and S. A. Piyavskii [29], Yu. M. Danilin [30] and has been widely used. The generalization to problems with non-convex constraints is new.

The smoothing method for global optimization purposes was proposed by D. B. Yudin [138] and studied in [52, 81, 94, 106]. The method of one-dimensional global optimization under non-convex constraints is proposed in [115].

Methods for global minimization of a concave function under linear constraints  are proposed in [3, 8, 126, 154, 188]. The application of the heavy ball and the gully step methods to  global optimization is discussed in [100], and random inertia methods are discussed in [106].

\textbf{Chapter 4}

Stochastic gradient averaging to construct convergent stochastic optimization methods was applied in [20, 22, 41, 52, 91, 91, 128].

Methods for minimizing nonsmooth functions with generalized gradient averaging were proposed in [22, 84, 85, 129, 171].

$\S$ 13. Nonmonotone methods with averaging of finite-difference analogs of the gradient were studied in [22].

$\S$ 14. Stability issues of numerical optimization methods are discussed in the books by V. G. Karmanov [51], B. T. Polyak [100], and the stability of difference schemes, in the monograph A. A. Samarsky [111].

   The stability of convex optimization methods with respect to random errors in the calculation of the gradient was studied in the paper by M. V. Mikhalevich [66].

The stability of the generalized gradient method with respect to deterministic errors in determining the gradient was studied in the monographs of N. Z. Shor [136], A. S. Nemirovskii and D. B. Yudin [73].

The stability of the subgradient method with respect to errors in the position of the points where the gradients are taken was established in [50, 132].

The method of averaged gradients for minimizing generalized differentiable functions was proposed in [84]. Adaptive step adjustments in nonrelaxation nonsmooth optimization methods are discussed in [37, 90, 91, 100, 129, 136].

The heavy ball method is described in [100]. The ravine (gully) step method was proposed by I. M. Gelfand and M. L. Tsetlin [14], and studied and generalized by Yu. E. Nesterov [74,76]. The justification of these methods for minimizing nonconvex nonsmooth functions is new.

\textbf{Chapter 5}

$\S$ 15. The conditional gradient method for quadratic problems was developed by Frank and Wolfe [148]. An analogue of the conditional gradient method for minimizing convex functions was studied in [78], and for minimizing Lipschitz functions, in [23].

$\S$ 16. The reduced gradient method was studied by Wolfe [184]; it is a generalization of the simplex method to non-linear programming problems and therefore has become widely used [6, 71, 131, 170].

$\S$ 17. Methods of possible directions were proposed by Zeutendijk [47]; they are also studied in the book [6].

$\S$ 18. The method of penalty functions was introduced by R. Courant, the monograph by Fiacco and McCormick [123] is devoted to it, and it occupies a central place in the monograph by V. V. Fedorov [122].

The existence of exact penalty functions was established by I. I. Eremin [37] and Zangwill [187]; they were also studied in [164]. Global properties of exact penalty functions were studied in [64, 142, 146].

$\S$ 19. A finite-difference method for minimizing Lipschitz functions under constraints was constructed as an analogue of B. N. Pshenichny's linearization method [105] and studied in [25].

$\S$ 20. The Arrow-Hurwicz method was proposed in the book [137], and the Cesaro averaging was studied in [73]. The finite difference method from $\S$ 20 was studied in [26].

\textbf{Chapter 6}

$\S$ 21. Information on multivalued mappings can be found in the books [34, 50, 60, 77, 104].

Measurable set-valued mappings are the  generalization of measurable functions. It suffices to get acquainted with measurable functions from the book by A. N. Kolmogorov and S. V. Fomin [53]; much useful information on measurability can be found in Warga's monograph [10]. Information on measurable set-valued mappings can be found in the survey [182] and in the papers [10, 48, 49, 57, 110, 144, 175].

We mainly follow the works of Castaing [144], Rockafellar [175], Himelberg [151]. Kasten's results are also presented in [10]. Lemma 21.7 is established in [79, 86].

$\S$ 22, 23. Random Lipschitz functions were studied in [79] (see Remark 22.1), while generalized differentiable functions were studied in [86].

\textbf{Chapter 7}

Problems and methods of stochastic programming can be found in the monographs of Yu. M. Ermoliev [41], D. B. Yudin [139, 140]. The behavior of optimization methods in the presence of random noise is discussed in the book by B.T. Polyak [100].

Stochastic programming methods originate from the paper by Robbins and Monroe [174]. The finite-difference method of stochastic approximation was proposed by Kiefer and Wolfowitz [155]. The monographs of Vazan [9], M. B. Nevelson and R. Z. Khasminsky [72] are devoted to the method of stochastic approximation. Techniques for accelerating the convergence of the stochastic approximation method are given in [112], and adaptive step adjustment was considered in [119].

$\S$ 24. Problems of convex stochastic programming were studied in [41, 73,50], a generalization to the weakly convex case was made in [91], and to the generalized differentiable case in [86].

The interpretation of the convergence of stochastic methods as a consequence of the stability of their deterministic variants is accepted in the book [50]. Conditions for the convergence of stochastic analogues of deterministic methods were studied in [89].

Information on the theory of probability is contained in the books of Loeve [58], B. N. Prokhorov and Yu. A. Rozanov [102], and A. N. Shiryaev [133].

Stochastic gradient averaging methods were studied in [20-22, 41,52, 85, 86, 100, 128].

The behavior of the stochastic heavy ball method in the quadratic case is discussed in [100]. The justification of the stochastic heavy ball and ravine step methods for solving nonconvex stochastic optimization problems is new.

$\S$ 25. Finite-difference methods of Lipschitz stochastic programming generalizing the classical method of Kiefer and Wolfowitz are constructed in the book [22]. The rate of convergence of such methods was studied in [70].

$\S$ 26, 27. The averaging operation was first used by Ya. Z. Tsypkin [128] and Yu. M. Ermoliev [41]; it is connected to tracking the minimum of a nonstationary quadratic function. The presentation of the material $\S$ 26,27 is based on [22]. The rate of convergence of the stochastic conditional gradient method was studied by S. V. Pashko [93]. Theorems 26.1–26.3 can also be obtained as a consequence of the results in [45].

$\S$ 28. A stochastic variant of the Arrow-Hurwicz finite-difference method described in $\S$ 20 is considered.

\newpage
\section*{BIBLIOGRAPHY}
\label{biblio}
 

\begin{enumerate}
\item Algorithms for optimization of design solutions /Edited by A. I. Polovinkin. -- M.: Sov. radio,
\item A l e x a n d r o v, P. S. Introduction to set theory and general topology. -- M. Nauka, 1977.
\item A s h e p k o v, A. G., Belov, B. I., Bulatov, V. P. et al. Methods of solving problems of mathematical programming and optimal control.- Novosibirsk: Nauka, -- 1984.
\item B a z h e n o v,  L. G. Convergence conditions of a minimization method of almost-differentiable functions. Cybern Syst Anal 8 (4), 607–609 (1972). https://doi.org/10.1007/BF01068282
\item B a zh e n o v,  L. G. Necessary and sufficient conditions of convergence of iterative procedures // Dokl. of the Academy of Sciences of the Ukrainian SSR. Ser. A.-- 1985.-- No. 2.-- P. 64--65.
\item B a z a r a a,  M., Shetty, K. Nonlinear Programming. Theory and algorithms. -- Moscow: Mir, 1982. 
\item B a t i s c h e v,  D. I. Methods of optimal design. -- Moscow: Radio i Svyas, 1984.
\item B u l a t o v,  V. P. Immersion methods in optimization problems. -- Novosibirsk: Nauka, 1977.
\item 
V a z a n,  M. Stochastic approximation. -- Moscow: Mir, 1972.
\item 
W a r g a,  J. Optimal control of differential and functional equations.- Moscow: Nauka, 1977.
\item 
V a s i l i e v,  F. P. Numerical methods for solving extreme problems. -- Moscow: Nauka, 1980.
\item 
V i l k o v,  A. V., Zhidkov, N. N., Shchedrin, B. M. Methods of finding the global minimum of a function of one variable // Zhurnal Vychislitelnoi Matematiki i Matematicheskoi Fiziki. -- 1975.- Vol. 15, No. 4.- P. 1040-1042.
\item 
Computational methods for selecting optimal design solutions / Edited by V. S. Mikhalevich. -- Kiev: Naukova Dumka, 1977.
\item 
G e l f a n d, I. M., Tsetlin, M. L. Principles of nonlocal search in automatic optimization systems // Dokl. of the USSR Academy of Sciences.- 1961.- Vol. 137, No. 2.- P. 295-298.
\item 
G i h m a n, I. I., Skorokhod, A. V. Introduction to the theory of random processes. -- Moscow: Nauka, 1977.
\item 
G l u s h k o v,  V. M. Fundamentals of paperless informatics. -- Moscow: Nauka, 1982.
\item 
G l u sh k o v a,  O. V., G u p a l, A. M. Numerical methods of minimization of maximum functions without gradient calculation // Kibernetika.- 1980.- No. 5.- P. 141-143.
\item 
G o l s h t e i n,  E. G. Generalized gradient method for finding saddle points // Ekonomika i matematicheskie metody. -- 1972. -- Vol. 8, No. 4. -- P. 569-579.
\item 
G o r b a ch u k,  V. I. General conditions of convergence of algorithms to arbitrarily specified sets// Kibernetika.- 1983.- No. 6.- P. 112-113.
\item 
G u p a l,  A. M., Bazhenov, L.G. Stochastic analog of the conjugant-gradient method. Cybern Syst Anal 8 (1), 138–140 (1972). \\https://doi.org/10.1007/BF01069146
\item
G u p a l, A.M., Bazhenov, L.G. A stochastic method of linearization. Cybern Syst Anal 8, 482–484 (1972). https://doi.org/10.1007/BF01069005 
\item 
G u p a l,  A. M. Stochastic methods of solving nonsmooth extreme problems.-- Kyiv: Naukova Dumka, 1979. 
\item 
G u p a l,  A. M. Algorithms for searching the extremum of nondifferentiable functions with constraints / Institute of Cybernetics of the Academy of Sciences of Ukrainian SSR. -- Kyiv, 1976, 18 p. -- Preprint 76-8.
\item 
G u p a l,  A. M. A method for the minimization of almost-differentiable functions. Cybern Syst Anal 13 (1), 115–117 (1977). \\https://doi.org/10.1007/BF01071397
\item 
G u p a l,  A. M. Minimizing method for functions that satisfy the Lipschitz condition. Cybern Syst Anal 16(5), 733–737 (1980). \\https://doi.org/10.1007/BF01078505
\item 
G u p a l,  A. M., Dubrovsky, V. B. B. Finite-difference method of Arrow-Hurwicz with averaging // Kibernetika. -- 1985. -- No. 4. -- P.120-121.
\item 
G u p a l  A. M., Loskutov V. G. Efficiency of finite-difference methods of minimization of nondifferentiable functions // Kibernetika. -- 1982. -- No. 2. --P. 115-117.
\item 
G u r i n,  L. P., Lobach, V.P. Combination of the Monte Carlo method with the steepest descent method when solving some extremal problems // Zhurnal Vychislitelnoi Matematiki i Matematicheskoi Fiziki.-- 1962.-- T. 2, No. 3.-- P. 499--502.
\item 
D a n i l i n,  Y. M., Piyavskiy, S. A. About one algorithm for finding the absolute minimum // Teoria Optimalnyh Resheniy. Vol. 2.- Kiev : Institute of Cybernetics of the Academy of Sciences of the Ukrainian SSR, 1967. -- P. 25-37.
\item 
D a n i l i n,  Y. M. Evaluation of the efficiency of one algorithm for finding the absolute minimum // Zhurnal Vychislitelnoi Matematiki i Matematicheskoi Fiziki.- 1971.- Vol. 11, No. 4.- P. 1026-1031.
\item 
D a n s k i n,  J. M. Theory of maximin. -- Moscow: Sovetskoe radio, 1970.
\item 
D a n z i g,  J. Linear programming, its generalizations and applications. -- Moscow: Progress, 1966.
\item 
D e m y a n o v,  V. F., M a l o z e m o v,  V. N. Introduction to minimax. -- Moscow: Nauka, 1972.
\item 
D e m y a n o v,  V. F., Vasiliev, L. V. Nondifferentiable optimization. -- Moscow: Nauka, 1981.
\item 
D o r o f e e v,  P. A. About some properties of the generalized gradient method // Zhurnal Vychislitelnoi Matematiki i Matematicheskoi Fiziki.- 1985.- Vol. 25, No. 2.- P. 181-189.
\item 
E v t u s h e n k o,  Y. G. Methods of solving extreme problems and their application in optimization systems. -- Moscow: Nauka, 1982.
\item 
E r e m i n, I. I. The penalty method in convex programming. Cybern Syst Anal 3(4), 53–56 (1967). https://doi.org/10.1007/BF01071708
\item 
E r e m i n,  I. I., Astafiev, N. N.  Introduction to the theory of linear and convex programming. -- Moscow: Nauka, 1976.
\item 
E r e m i n,  I. I., Mazurov, V. D., Astafiev, N. N.  Nonsimple problems of linear and convex programming. -- Moscow: Nauka, 1983.
\item 
E r m o l' e v, Y.M. Methods of solution of nonlinear extremal problems. Cybern Syst Anal 2(4), 1–14 (1966). https://doi.org/10.1007/BF01071403
\item 
E r m o l i e v,  Y. M. Methods of stochastic programming. -- Moscow: Nauka, 1976.
\item 
E r m o l' e v, Y. M., Verchenko, P.I. A linearization method in limiting extremal problems. Cybern Syst Anal 12(2), 240–245 (1976). \\https://doi.org/10.1007/BF01069893
\item 
E r m o l i e v,  Yu. S., Yastremsky, A. I. Stochastic models and methods in economic planning. -- Moscow: Nauka, 1979.
\item 
Z a v r i e v,  S. K. Stochastic gradient methods for solving minimax problems. -- Moscow: Publishing House of Moscow State University, 1984.
\item 
Z a v r i e v, S. K. About one analog of the amplified law of large numbers for sums of dependent random variables // Software and models of operations research / Edited by L. N. Korolev and P. S. Krasnoshchekov. -- Moscow: Publishing House of Moscow State University, 1985.
\item 
Z a n g w i l l,  W. Nonlinear programming: a unified approach. -- Moscow: Sovetskoe radio, 1973. 
\item 
Z o y t e n d i j k,  G. Methods of feasible directions. -- Moscow: Inostrannaya Literatura, 1963. 
\item 
I o f f e,  A. D., Levin, V. L. Subdifferentials of convex functions // Proceedings of Moscow Mathematical Society, vol. 26.- Moscow: Nauka, 1972.- P. 3-73.
\item 
I o f f e,  A. D., Tikhomirov, V. M. Theory of Extreme Problems. -- Moscow: Nauka, 1974.
\item 
Iterative Methods in Game Theory and Programming / Edited by V. V. Kuznetsov. 3. Belenkiy, and V. A. Volkonsky. -- Moscow: Nauka, 1979.
\item 
K a r m a n o v,  V. G. Mathematical Programming. -- Moscow: Nauka, 1979.
\item 
K a t k o v n i k,  V. Ya. Linear estimation and stochastic optimization problems. -- Moscow: Nauka, 1976.
\item 
K o l m o g o r o v,  A. N., Fomin, S. V. Elements of the theory of functions and functional analysis. -- Moscow: Nauka, 1968.
\item 
K r a s s o v s k s k i y,  N. N. Theory of motion control. -- Moscow: Nauka, 1968.
\item 
K r u z h k o v,  S. N. To the question of differentiability of functions of many variables almost everywhere // Bulletin of Moscow State University. -- 1976.- No. 6.- P. 67-70.
\item 
L e v i n,  A. Yu. About one algorithm of minimization of convex functions// Dokl. of the USSR Academy of Sciences.- 1965.- Vol. 160, No. 6.- P. 1244-1247.
\item 
L e v i n,  V. L. Measurable sections of multivalued mappings into topological spaces and upper envelopes of Karateodori integrants // Reports of the USSR Academy of Sciences. -- Vol. 252, No. 3. -- P. 535-539.
\item 
L o e v,  M. Probability theory. -- Moscow: Inostrannaya Literatura, 1962.
\item 
L y a s h k o,  S. I. Linearization method in solving extremal problems for functions with unknown parameter // Mathematical methods of operations research and reliability theory.- Kyiv: Institute of Cybernetics of Ukrainian SSR, 1978. -- P. 13-20.
\item 
M a k a r o v,  V. L., Rubinov A. M. Mathematical theory of economic dynamics and equilibrium. -- Moscow: Nauka, 1973.
\item 
M i k h a l e v i c h, V. S. Sequential optimization algorithms and their application. Part I. Cybern Syst Anal 1(1), 44–55 (1965). \\https://doi.org/10.1007/BF01071444
\item 
M i k h a l e v i c h,  V. S., V o l k o v i c h,  V. L. Computational methods of research and design of complex systems. -- Moscow: Nauka, 1982.
\item 
M i k h a l e v i c h,  V. S., Gupal, A.M. \& Bolgarin, I.V. Numerical methods of solving nonstationary and limit extremal problems without evaluating gradients. Cybern Syst Anal 20, 684–687 (1984). https://doi.org/10.1007/BF01071613
\item 
M i k h a l e v i c h,  V. S., G u p a l,  A. M., Dubrovskii, V. B. Global character of exact nonsmooth penalties// Kibernetika.- 1988.- No. 4. -- P.111-112.
\item 
M i k h a l e v i c h, M.V. Generalized stochastic method of centers. Cybern Syst Anal 16(2), 292–296 (1980). https://doi.org/10.1007/BF01069122
\item 
M i k h a l e v i c h, M. V. Stability of stochastic programming methods to stochastic quasigradient computing errors. Cybern Syst Anal 16(3), 434–438 (1980). https://doi.org/10.1007/BF01078265
\item 
M o i s e e v,  N. N. Elements of the theory of optimal systems. -- Moscow: Nauka, 1975.
\item 
M o i s e e v N. N., Ivanilov, Y. P., Stolyarova, E. M. Optimization methods.  -- Moscow: Nauka, 1978.
\item 
M o t s k u s,  I. B. Multiextremal tasks in design.  -- Moscow: Nauka, 1967.
\item 
M u r a u s k a s s,  G. G. Investigation of the convergence rate of iterative stochastic algorithms of optimization and estimation: Author's thesis .... Cand. of Phys.-Math. sciences.  -- Kiev: Institute of CYbernetics of the Academy of Sciences of the Ukranian SSR, 1983.
\item 
M u r t a g h,  B. Modern Linear Programming.  -- Moscow: Mir, 1984.
\item 
N e v e l s o n,  M. B., Khasminsky R. 3. Stochastic approximation and recurrence estimation.  -- Moscow: Nauka, 1972.
\item 
N e m i r o v s k y,  A. S., Yudin, D. B. Complexity of problems and efficiency of optimization methods.  -- Moscow: Nauka, 1979.
\item 
N e s t e r o v, Y. E. Method for solving the problem of convex programming with convergence rate O (1/$k^2$) // Dokl. of the USSR Academy of Sciences. -- 1983. -- Vol. 269, No. 3. -- P. 543-547.
\item 
N e s t e r o v,  Yu. E. Methods of minimization of nonsmooth convex and quasi-convex functions // Ekonomika i matematicheskie metody. -- 1984.-T. 20, No. 3.-P. 519-531.
\item 
N e s t e r o v,  Yu. E. About one class of methods of unconditional minimization of a convex function with high convergence rate // Zhurnal Vychislitelnoi Matematiki i Matematicheskoi Fiziki. -- 1984. -- Vol. 24, No. 7. -- P. 1090-1093.
\item 
N i k a i d o,  X. Convex structures and mathematical economics.  -- Moscow: Mir, 1972.
\item 
N i k o l a e v a,  N. D. About one algorithm for solving problems of convex programming//Ekonomika i matematicheskie metody. -- 1974. -- Vol. 10, No. 5. -- P. 941-946.
\item 
N o r k i n, V. I. Stochastic Lipschitz functions. Cybern Syst Anal 22(2), 226–233 (1986). https://doi.org/10.1007/BF01074785
\item 
N o r k i n, V. I. Nonlocal minimization algorithms of nondifferentiable functions. Cybern Syst Anal 14(5), 704–707 (1978). https://doi.org/10.1007/BF01069307
\item 
N o r k i n,  V. I. Two remarks: on the smoothing method in multiextremal optimization and on the finite-difference method in nondifferentiable optimization // Abstracts of 3-rd All-union seminar "Numerical Methods of Nonlinear Programming".  -- Kharkov: Kharkov State University. -- 1979.
\item 
N o r k i n, V. I. Conditions of convergence of combined algorithms of nonlinear programming. Cybern Syst Anal 15(6), 854–860 (1979).\\ https://doi.org/10.1007/BF01069396
\item 
N o r k i n,  V. I. Generalized-differentiable functions. Cybern Syst Anal 16, 10–12 (1980). https://doi.org/10.1007/BF01099354
\item 
N o r k i n,  V. I. A method of minimizing an undifferentiable function with generalized-gradient averaging. Cybern Syst Anal 16, 890–892 (1980). \\https://doi.org/10.1007/BF01069064
\item 
N o r k i n, V. I. Method of generalized gradient descent. Cybern Syst Anal 21(4), 495–505 (1985). https://doi.org/10.1007/BF01070609
\item 
N o r k i n,  V. I. Generalized gradient method in the problem of non-convex nonsmooth optimization: Author's thesis ... Cand. of Phys.-Math. sciences. -- Kyiv: Institute of Cybernetics of the Academy of Sciences of the Ukrainian SSR. -- 1983.
\item 
N u r m i n s k i i, E. A. Convergence conditions for nonlinear programming algorithms. Cybern Syst Anal 8(6), 959–962 (1972). \\https://doi.org/10.1007/BF01068520
\item 
N u r m i n s k i i, E. A. The quasigradient method for the solving of the nonlinear programming problems. Cybern Syst Anal 9(1), 145–150 (1973). https://doi.org/10.1007/BF01068677
\item 
N u r m i n s k i i, E. A. Convergence of stochastic analogues of deterministic nonlinear programming procedures // Teoria Optimalnyh Resheniy.  -- Kyiv: Institute of Cybernetics of the Academy of Sciences of the Ukrainian SSR, 1973.  -- pp. 60-73.
\item 
N u r m i n s k i i, E. A., Zhelikhovskii, A.A. Investigation of one regulating step in a quasigradient method of minimizing weakly convex functions. Cybern Syst Anal 10(6), 1027–1031 (1974). https://doi.org/10.1007/BF01069448
\item 
N u r m i n s k i i, E. A. Numerical methods for solving deterministic and stochastic minimax problems.  -- Kyiv: Naukova Dumka, 1979.
\item 
O r t e g a, J., Reinboldt, V. Iterative methods for solving nonlinear systems of equations with many unknowns.  -- Moscow: Mir, 1975. *
\item 
P a s h k o, S.V. On the rate of convergence of the stochastic linearization method // Kibernetika.  -- 1984.  -- No. 4. -- P. 118-119.
\item 
P e t e r s e n, I. Statistical optimization through smoothing // Izvestia Academii Nauk of the USSR. Technique, cybernetics.  -- 1969.  -- No. 2.
\item 
P i y a v s k i, S. A. Algorithm for finding the absolute minimum of functions and Theory of optimal solutions. Vol. 2.-- Kyiv: Institute of Cybernetics of the Academy of Sciences of Sciences of the Ukrainian SSR, 1967. -- P. 13-24.
\item 
P i y a v s k i, S. A. One algorithm for finding the absolute extremum of a function // Zhurnal Vychislitelnoi Matematiki i Matematicheskoi Fiziki.  -- 1972.  -- Vol. 12, No. 4.-- P. 888--896.
\item 
P o d i n o v s k y, V. V., Gavrilov, V.M. Optimization according to consequently applied criteria.-- Moscow: Sovetskoe radio, 1975.
\item 
P o l a k, E. Numerical optimization methods. A unified approach.-- Moscow: Mir, 1974.
\item 
P o l y a k, B. T. One general method for solving extremal problems // Dokl. of the Academy of Sciences of USSR .  -- 1967.  -- Vol. 174, No. 1.  -- P. 33-36.
\item 
P o l y a k, B. T. Introduction to optimization.  -- Moscow: Nauka, 1983.
\item 
P o n t r y a g i n, L.S., Boltyansky, V.G., Gamkrelidze, R.V., Mishchenko, E.F. Mathematical theory of optimal processes.  -- Moscow: Fizmatgiz, 1961.
\item 
P r o k h o r o v, B. N., Rozanov, Yu. A. Theory of Probability.  -- Moscow: Nauka, 1967.
\item 
P s h e n i c h n y, B. N. Necessary conditions for extremum.  -- Moscow: Nauka, 1969.
\item 
P s h e n i c h n y, B. N. Convex analysis and extremal problems.  -- Moscow: Nauka, 1980.
\item 
P s h e n i c h n y, B. N., Danilin Yu. M. Numerical methods in extremal problems.  -- Moscow: Nauka, 1975.
\item 
R a s t r i g i n, L. A. Extreme control systems.  -- Moscow: Nauka, 1974.
\item 
R z h e v s k y, S.V. $\epsilon$-subgradient method for solving the problem of convex programming // Zhurnal Vychislitelnoi Matematiki i Matematicheskoi Fiziki.  -- 1981.  -- Vol. 21, No. 5.  -- P. 1126-1132.
\item 
R z h e v s k i i, S. V. Monotone algorithm for seeking saddle point of unsmooth function. Cybern Syst Anal 18(1), 111–115 (1982). \\https://doi.org/10.1007/BF01078057
\item 
R o c k a f e l l a r, R. Convex analysis.  -- Moscow: Mir, 1973.
\item 
R o c k a f e l l a r, R. Convex integral functionals and duality // Mathematical Economics.  -- Moscow: Mir, 1974.  -- P. 222-237.
\item 
S a m a r s k y, A. A. Theory of difference schemes.  -- Moscow: Nauka, 1977.
\item 
S a r i d i s, J. Self-organizing stochastic control systems.  -- Moscow: Mir, 1980. 
\item S t e r n b e r g, S. Lectures on differential geometry.  -- Moscow: Mir, 1970.
\item 
S t r o n g i n, R. G. Numerical methods in multiextremal problems.  -- Moscow: Nauka, 1978.
\item 
S t r o n g i n, R. G., Markin D. L. Algorithm for solving multiextremal problems with nonlinear non-convex constraints // Methods of mathematical programming and software.  -- Sverdlovsk: IMM UC USSR Academy of Sciences, 1984, p. 103-104.
\item 
S u k h a r e v, A. G. Optimal search for extremum.  -- Moscow: Moscow State University Publishing House, 1975.
\item 
S u k h a r e v, A. G. Global extremum and methods for finding it and Mathematical methods for operations research / Ed. N. N. Moiseeva and P. S. Krasnoshchekova.  -- Moscow: Moscow State University Publishing House, 1981.  -- P. 4-37.
\item 
T i k h o n o v, A. N., Arsenin V. Ya. Methods for solving ill-posed problems.  -- Moscow: Nauka, 1979.
\item 
U r y a s e v, S. P., Mirzoakhmedov, F. Step adjustments in stochastic programming methods // Zhurnal Vychislitelnoi Matematiki i Matematicheskoi Fiziki.  -- 1983.  -- No. 6.  -- P. 1314-1325.
\item 
F a d d e e v, D. K., Faddeeva, V.N. Computational methods of linear algebra.  -- Moscow: Fizmatgiz, 1960.
\item 
F e d o r o v, V. V. On optimization problems with an ordered set of constraints // Zhurnal Vychislitelnoi Matematiki i Matematicheskoi Fiziki.  -- 1975.  -- Vol. 15, No. 5.  -- P. 1126-1137.
\item 
F e d o r o v, V. V. Numerical methods of maximin.  -- Moscow: Nauka, 1979.
\item 
F i a c c o, A., McCormick, D. P. Nonlinear programming. Methods of sequential unconditional minimization.  -- Moacow: Mir, 1972.
\item 
X a c h i y a n, L. G. Polynomial algorithm in linear programming // Dokl. Academy of Sciences of the USSR.  -- 1979.  -- Vol. 244, No. 5.  -- P. 1093-1096.
\item 
H i m m e l b l a u, D. Applied nonlinear programming.  -- Moscow: Mir, 1975.
\item 
H o a n g T u y. Concave programming under linear constraints // Reports of the Academy of Sciences of the USSR.  -- 1964.  -- Vol. 159, No. 1.  -- P. 32-35.
\item 
T s y p k i n, Ya. Z. Adaptation and training in automated systems.  -- M.: Nauka, 1968.\item 
Tsypkin Ya. 3. Generalized learning algorithms // Automation and telemechanics.  -- 1976.  -- No. 1.  -- P. 97-103.
\item 
C h e p u r n o y, N. D. Method of averaged quasigradients with stepwise step adjustment for minimizing weakly convex functions // Kibernetika.  -- 1981.  -- No. 6.  -- P. 131.
\item 
C h e p u r n o y, N. D. On a relaxation algorithm for minimizing non-smooth functions // Kibernetika.  -- 1984.  -- No. 2.  -- P. 110-111.
\item 
Numerical methods of conditional optimization / Eds. F. Gill and W. Murray  -- Moscow: Mir, 1977.
\item
S h e p i l o v, M. A. Method of the generalized gradient for finding the absolute minimum of a convex function. Cybern Syst Anal 12, 547–553 (1976). https://doi.org/10.1007/BF01070389 
\item 
S h i r y a e v, A. N. Probability.  -- Moscow: Nauka, 1980.
\item 
S h o r, N. Z. Application of the gradient descent method for solving a network transport problem // Materials of a scientific seminar on theor. and adj. question cybernetics and research operations.--Kyiv: PK  of Academy of Sciences of the Ukrainian SSR, 1962, issue. 1. -- P. 9--17.
\item 
S h o r, N. Z. Cut-off method with space extension in convex programming problems. Cybern Syst Anal 13, 94–96 (1977). https://doi.org/10.1007/BF01071394
\item 
S h o r, N. Z. Methods for minimizing non-differentiable functions and their applications.  -- Kyiv: Naukova Dumka, 1979.
\item 
A r r o u, K. J., Hurwicz, L., Uzawa, X. Research on linear and nonlinear programming.  -- Moscow: Inostrannaya Literatura, 1962.
\item 
Y u d i n, D. B. Methods of quantitative analysis of complex systems // Izv. Academy of Sciences of the USSR. Technique, cybernetics.  -- 1965.  -- No. 1.  -- P. 3-13.
\item 
Y u d i n, D. B. Mathematical methods of control under conditions of incomplete information. -- Moscow: Sovetskoe radio, 1974.
\item 
Y u d i n, D.B. Problems and methods of stochastic programming.-- Moscow: Sovetskoe radio, 1979.
\item 
Y u d i n, D. B., Nemirovsky A. S. Information complexity and effective methods for solving convex extremal problems // Ekonomika i matematicheskie metody.-- 1976.-- Vol. 12, No. 2.-- P. 357--369.
\item 
B a z a r a a,  M. S., Good, J. J. Sufficient conditions for a globally exact penalty function without convexity // Math. Progr. Study 19.-- 1982.-- Vol. 19.-- P. 1--15. 
\item 
B i h a i n,  A. Optimization of upper semidifferentiable functions // J. Optimiz. Theory Appl.-- 1984.-- Vol. 44.-- P. 545--568. 
\item 
C a s t a i n g, C. Sur le multi-applications measurables // Revue Francaise d’Informatique et de Recherche Operationnelle.-- 1967.-- No. 1.-- P. 91--126.
\item 
C l a r k e, F. H. Generalized gradients and applications// Trans. Amer. Math. ' Soc.-- 1975.-- Vol. 205.-- P. 247--262. ’
\item 
C l a r k e, F. H. A new approach to Lagrange multipliers// Math. Oper. Res.-- 1976.-- Vol. 1, No. 2.-- P. 165--174. 
\item D e b r e u,  G. Integration of correspondences // Proc. 5-th Berckeley Symp. on Math. Stat.-- Prob. 11, part 1.-- P. 351-372. 
\item 
F r a n k, M., Wolfe, P. An algorithm for guadratic programming // Naval ' Res. Log. Quart.-- 1956.-- Vol. 3, No. 1--2.-- P. 113-120. 
\item 
F u b i a n, V. Stochastic approximations of constrained minima // Trans. 4-th | Prague Conf. Inf. Theory. Stat. Dec. Funct. Random Process, 1965.-- Prague.-- | 1967.-- P. 277-290. 
\item 
H i e n  N g u y e n, V., Strodiot, J. J., Mifflin, R. On conditions to have bounded multipliers in locally lipschitz programming // Math. Progr.-- 1980.  -- Vol. 18, 1.--P. 100--106. 
\item 
H i m e 1 b e r g, C. J. Measurable relations// Fundamenta Mathematicae. -- 1975.-- Vol. 87, No. 1.-- P. 53-72. :
\item 
H i r i a r t-U r r u t y, J. B. On optimality conditions in nondifferentiable ‘ programming // Math. Progr.-- 1978.-- Vol. 14.-- P. 73-86.
\item 
H i r i a r t-U r r u t y, J. B. Refinements of necessary optimality conditions  in nondifferentiable programming. I // Appl. Math. Optim.-- 1979.-- No. 5.-- P. 63-82; II // Math. Progr. Study 19.-- 1982.-- P. 120-139.
\item 
H o a n g T u y. Concave Minimization under Linear Constraints with Special Structure// Optimization.-- 1985.-- Vol. 16, No. 3.-- P. 335-352.
\item 
K i e f e r, J., Wolfowitz, J. Stochastic estimation of the maximum of a regression function // Ann. Math. Stat.-- 1952.-- Vol. 23, No. 8.--P. 462-466.
\item 
K i w i e l,  K. C. An aggregate subgradient method for nonsmooth convex minimization // Math. Progr.-- 1983.-- Vol. 27.-- P. 320-341.
\item 
K i w i e l,  K. C. Methods of Descent for Nondifferentiable Optimization// Lecture Notes in Mathematics.-- 1985.-- Vol. 1133.-- P. 362.
\item 
K u s h n e r, H. J., Clark, D. S. Stochastic approximation methods for constrained and unconstrained systems.-- N. Y., Heidelberg, Berlin: Springer, . 1978. 
\item 
L e b o u r g, M. G. Valeur moyenne pour gradient generalize // C. R. Acad. Sc Paris.-- 1975.-- Vol. 281, ser. A, No. 19.-- P. 775-779.
\item 
L e m a r e c h a l, C. An extension of Davidon methods to nondifferentiable problems// Math. Progr. Study 3.-- 1975.-- P. 95-109.
\item 
L e m a r e c h a l, C., Strodiot, J. J., Bihain, A. On a bundle algo­ rithm for nonsmooth minimization: Nonlinear Progr. 4 / Eds O. L. Mangasarian, R. R. Meyer, S. M. Robinson.-- N. Y.: Academic Press, 1981.-- P. 245-281.
\item 
L e u e n b e r g e r, D. G. Quasi-convex programming //SIAM J. Appl. Math.-- 1968.-- Vol. 16.-- P. 1090-1095.
\item 
M a n g a s a r i a n, O. L. Pseudo-convex functions // SIAM J. Control.-- 1965.-- Vol. 3, No. 3.-- P. 281-290.
\item 
M a n g a s a r i a n, O. L. Exact penalty functions in nonlinear programming// Math. Progr.-- 1979.-- K2 17.-- P. 251-269.
\item 
M c C o r m i c k, G. P. Attempts to calculate global solutions of problems that may have local minima: Numerical Methods for Nonlinear Optimization / Ed. F. Lootsma.-- Academic Press, London, N. Y. 1972,-- P. 209-221.
\item 
M i f f l i n, R. An algorithm for constrained optimization with semismooth func­ tions // Math. Oper. Res.-- 1977.--No. 2,--P. 191-207.
\item 
M i f f l i n, R. Semismooth and semiconvex functions in constrained optimization // SIAM J. Contr. and Optim.-- 1977.-- Vol. 15, No. 6.-- P. 959-972.
\item 
M i f f l i n, R. A modification and an extentions of Lemareshal’s algorithm for nonsmooth minimization // Math. Progr. Study 17.-- 1982.-- P. 77-90.
\item 
M i r i c \'a, S. A note on the generalized differentiability of mappings // NonL Anal. Theory, Methods, Applications.-- 1980.-- Vol. 4, No. 3.-- P. 567-575. https://doi.org/10.1016/0362-546X(80)90092-9
\item 
M u r t a g h, B. A., Saunders M. Large scale linearly constrained optimization // Math. Progr.-- 1978.-- Vol. 14.-- P. 41-72.
\item 
N u r m i n s k i,  E. A. On $\epsilon$-subgradient methods of non-differentiable optimization // Lect. Notes in Control and Information Sciences.-- 1979.-- No. 14.-- P. 187-195.
\item 
P o l a k, E., Mayne, D. Q. Algorithm models for nondifferentiable optimization // SIAM J. Control and Optimization.-- 1985.-- Vol. 23, No.3.-- P. 477-491.
\item 
R a d e m a c h e r, H. Uber partielle und total differenzierbarkeit von funktionen mehrerer variabeln und uber die transformation der doppelintegrale // Math. Anna- len.-- 1919.-- Vol. 79.- P. 340-359.
\item 
R o b b i n s, H., Monro, S. A stochastic approximation method // Ann. Math. Stat.-- 1951.-- Vol. 22, No. 3.-- P. 400-407.
\item 
R o c k a f e l l a r, R. T. Measurable dependence of convex sets and functions on parametries // J. Math. Anal. Appl.-- 1969.-- Vol. 28.-- P. 4-25.
\item 
R o c k a f e l l a r, R. T. Directionally lipschitzian functions and subdifferen­ tial calculus // Proc. London Math. Soc.-- 1979.-- Vol. 39.-- P. 331-355.
\item 
R o c k a f e l l a r, R. T. Generalized directional derivatives and subgradients ofnonconvex functions // Canad. J. Math.-- 1980.-- Vol. 32.-- P. 257-280.
\item 
R o c k a f e l l a r, R. T. The theory of subgradients and its applications to problems of optimization. Convex and nonconvex functions. Research Notes in Mathematics 1 / Eds. K. H. Hoffman, R. Wille.-- Berlin: Heldermann.-- 1981.
\item 
S t r o d i o t,  J. J., Hien Nguyen, V., Heukemez, N. $\epsilon$-optimal solution in nondifferentiable convex programming and related questions // Math. Prog.-- 1982.-- Vol. 25.-- P. 307-323.
\item 
T h i b a u l t,  Z. Quelques proprietes des sous-differentiels de fonctions reels localment lypchitziennes definies sur un espace de Banach separable // C. R. Acad. Sc. Paris.-- 1976 -- Vol. 282, Ser. A, No. 10.-- P. 507-610.
\item Towards global optimization / Eds L. C. W. Dixon, G. P. Szego.-- Amsterdam: North-Holland, 1975.
\item 
W a g n e r, D. H. Survey of measurable selection theorems // SIAM J. Control and Opt.-- 1977.-- Vol. 15, No. 5.-- P. 859-903.
\item 
W e t s, R. Stochastic programming. Solution techniques and approximation schemes // Math. Program. State Art. 11-th Int. Symp. Bonn., 23--27 Aug. 1982. Berlin e. a. 1983.-- P. 566-603.
\item 
W o l f e, P. Methods for nonlinear constraints.-- Heidelberg: North-Holland, 1967.-- P. 120-131.
\item 
W o l f e, P. A method of conjugate subgradients for minimizing nondifferentiable functions // Math. Progr. Study 3.-- 1975.-- P. 145-173.
\item W o l f e, P. Finding a nearest point in a polytope // Math. Progr.-- 1976.-- Vol. 11.-- P. 128-149.
\item 
Z a n g w i l l, W. I. Nonlinear programming via penalty functions//Management Science.-- 1967.-- Vol. 13.-- P. 344-388.
\item 
Z w a r t, P. B. Global maximization of a convex function with linear inequality constraints // Operat. Res.-- 1974.-- Vol. 22, No. 3.-- P. 602-609.
\end{enumerate}

\end{document}